\documentclass[reqno, 10pt]{amsbook}

\usepackage{color}

\usepackage{amssymb}
\usepackage{mathtools}

\makeatletter
\def\ckech{\mathaccent"\accentclass@014}
\def\tah{\mathaccent"\accentclass@05E}
\makeatother

  \renewcommand\check{\bm\ckech}
\renewcommand\hat{\bm\tah}

\swapnumbers
\makeindex 

\usepackage{bm}

\makeatletter
\renewcommand{\l@subsection}{\@tocline{3}{0pt}{3pc}{2.5pc}{}}
\def\indexchap#1{\global\topskip 2pc\relax
  \twocolumn[\fontsize{\@xivpt}{18}%
    \vskip\topskip\vskip-\baselineskip\hbox{}
    \bfseries\centering #1\par]%
  \global\topskip 34\p@
}
\makeatother

\newtheorem{theorem}{Theorem}[section]
\newtheorem{lemma}[theorem]{Lemma}
\newtheorem{corollary}[theorem]{Corollary}

\theoremstyle{definition}
\newtheorem{assumption}[theorem]{Assumption}
\newtheorem{definition}[theorem]{Definition}

\theoremstyle{remark}
\newtheorem{example}[theorem]{Example}
\newtheorem{remark}[theorem]{Remark}

\newtheorem*{theorem*}{Theorem}

\newcommand\loc{\textnormal{loc}}
\newcommand{\shharp}{=\kern -.5em\|}
\newcommand{\vsharp}{\asymp\kern -.5em\|}

\newcommand\esssup{\operatornamewithlimits{ess\,sup\,}}

 \newcommand{\EO}{\overset{\,\scalebox{0.4}
{$\boldsymbol 0$}}{E}\,\!}

 \makeatletter 
 \def\dashint{\operatorname{\,\,\,\mathclap{\!\int}\! \!\text{\bf--}\!\!}}
 \makeatother

\makeatletter
\def\dashintindex{\operatorname%
{-\kern-.7em\DOTSI\intop\ilimits@}}%
\def\dashint{\operatorname%
{\,\,\text{\bf--}\kern-.98em\DOTSI\intop\ilimits@\!\!}}
\makeatother

\def\sfA{{\sf A}}
\def\sfB{{\sf B}}

\def\sfL{{\sf L}}

\def\lpq{$L_{{p,q}}\,$}

\newcommand\sfb{{\sf b}}
\newcommand\sfbd{{\sf{d}\!\sf{I}}}
\newcommand\sfd{{\sf{d}}}
 
\newcommand\sfp{{\sf p}}

\newcommand\sft{{\sf t}}

\newcommand\gb{\mathfrak{b}}

\newcommand\bB{\mathbb{B}}
\newcommand\bC{\mathbb{C}}
 
\newcommand\bM{\mathbb{M}}
\newcommand\bR{\mathbb{R}}

\newcommand\bS{\mathbb{S}}

\newcommand\bZ{\mathbb{Z}}

\def\sft{{\sf t}}

\newcommand\frA{\mathfrak{A}}
\newcommand\frB{\mathfrak{B}}

\newcommand\frF{\mathfrak{F}}

\newcommand\frM{\mathfrak{M}}

\newcommand\frT{\mathfrak{T}}
\newcommand\gW{\mathfrak{W}}

\newcommand\cA{\mathcal{A}}
\newcommand\cB{\mathcal{B}}

\newcommand\cF{\mathcal{F}}

\newcommand\cL{\mathcal{L}}
\newcommand\cM{\mathcal{M}}
\newcommand\cN{\mathcal{N}}
\newcommand\cO{\mathcal{O}}
\newcommand\cP{\mathcal{P}}

\newcommand\cS{\mathcal{S}}
\newcommand\cT{\mathcal{T}}
\newcommand\cZ{\mathcal{Z}}

\newcommand\usigma{\underline{\sigma}}
\newcommand\ua{\underline{a}}
\newcommand\ub{\underline{b}}

\newcommand\scB{\normalfont\textsc{b}}  

\renewcommand\({{\rm(}}

\renewcommand\){{\rm)}}

\def\+){\tmspace+\thinmuskip{.05em}\)}

\def\dashnorm{\,\,\text{\bf--}\kern-.5em\|}

\def\binLpq{$b\in E_{(p,q),1},
\sfd_{0}/p+1/q\leq 1$}

\newcommand{\nlimsup}{\operatornamewithlimits{\overline{lim}}}
\newcommand{\nliminf}{\operatornamewithlimits{\underline{lim}}}

\newcommand{\sign}{\text{\rm\,sign}\,}

\newcommand{\osc}{\operatornamewithlimits{osc}}

\newcommand{\tr}{{\rm tr}\,}

\renewcommand{\eqref}[1]{\text{\rm(\ref{#1})}}

\newcommand\dist{{\rm dist}\,}

\makeatletter
\renewcommand{\thechapter}{\arabic{chapter}}
\newcommand{\mychapter}[2][1]{\chapter{#2}
\setcounter{equation}{0}
\setcounter{theorem}{0}
 \edef\@currentlabel{.\arabic{chapter}}   
 \markboth{
 \textsc{Chapter \arabic{chapter}. {#1}}%
 }{}}
\newcommand{\mysection}[2][1]{%
\section{#2}\setcounter{equation}{0}
 \edef\@currentlabel{.\thechapter.\arabic{section}}    
 \markright{
 \textsc{Section
 \arabic{section}. {#1}}}}
\newcommand{\mysubsection}[1]{\subsection{#1}
\edef\@currentlabel{.\arabic{chapter}.\arabic{section}:\arabic{subsection}}}
\def\mojtag#1{\hbox{\m@th\normalfont#1}}
\def\tagform@#1{\mojtag{(\ignorespaces\arabic{equation}\unskip)}}
\def\swappedhead#1#2#3{%
\kern-\parindent
{\bf\arabic{theorem}}\thmname{\@ifnotempty{#2}{. }#1}%
  \thmnote{ \textmd{\upshape(#3)}}}
\def\swappedhead@plain#1#2#3{%
{\bf\arabic{theorem}}\thmname{\@ifnotempty{#2}{. }#1}%
  \thmnote{ \textmd{\upshape(#3)}}}
\def\granddadref#1{\expandafter\@setref\csname r@#1\endcsname\@firstoftwo{#1}}
\renewcommand{\ref}[1]{%
$\def\aref{\if.\granddadref{#1}\else-1{\bf??}%
\protect\G@refundefinedtrue
\@warning{Reference `#1' on page \thepage \space is not defined}
\fi}%
\def\refa{\ifnum\arabic{chapter}=\aref\else\aref\fi}%
\def\refaa{\ifnum-1=\thechapter\else-1\fi}%
\def\refab{\if\Alph{chapter}\refaa-1\else-2\fi}
\def\refac{\ifnum\refab=\aref\else\refa\fi}
\def\refad{\if\Alph{chapter}\refac\else\refa\fi}
\def\refb{\if.\refad\else\arabic{section}.\aref\fi}%
\def\refc{\ifnum\arabic{section}=\refb\else.\refb\fi}%
\def\refd{\if.\refc\else\refb\fi}
\textup{\ifnum-1=\refd\else\refd\fi}$}

\makeatother

\begin{document}

\frontmatter
\title{Stochastic It\^o Equations and Parabolic
Second-Order Equations with singular Drift}
\author{N.V. Krylov}
\address{School of Mathematics, University of Minnesota, Minneapolis, MN, 55455}
\email{nkrylov@umn.edu}

\renewcommand{\subjclassname}{\textup{2010} 
Mathematics Subject Classification}
 
\subjclass{Primary  ;\\Secondary  }

\begin{abstract}
The aim of the book is to present some
recent results in the theory of
stochastic It\^o equations
with singular deterministic part (drift) and its
applications to second-order elliptic
and parabolic equations with singular
first-order coefficients. The singularity is characterized by means of Morrey spaces
and this allows for much more singular
coefficients than those from Lebesgue spaces. For instance, first-order
coefficients having behavior like $1/|x|$
near the origin are allowed.

In the first part of the book we are dealing with equations having just measurable coefficients and treat the Markov diffusion 
time-inhomogeneous processes $X$ corresponding
to parabolic operators. In particular,
mixed-norm parabolic Aleksandrov estimates, Harnack inequality and H\"older
continuity of $X$-caloric functions
are investigated. This produces the corresponding results in PDEs such as
extended Aleksandrov maximum principle,
 Harnack inequality and H\"older
continuity of PDE-caloric functions.

In two remaining chapters we concentrate
on weak and strong solutions of It\^o
equations which requires some 
regularity restrictions
on the diffusion matrix (or second-order
coefficients in the PDE language).
We give the best to date conditions
in terms of Morrey spaces
for the existence and uniqueness of weak and strong solutions
of It\^o equations with singular drift.
The majority of our main results are new even if
the drift part is zero.

 For mathematicians working in the areas of
 parabolic second-order equations and It\^o
 stochastic equations with singular ingredients.
\end{abstract}

\maketitle

\setcounter{page}{4}
 \setcounter{tocdepth}{3}
 \tableofcontents

\chapter*{Preface}

Let $\bR^{d}$, $d\geq2$,
\index{$B$@Sets!$\bR^{d}$}%
\index{$B$@Sets!$\bR^{d+1}$}%
be a $d-$dimensional Euclidean space of points
$x=(x^{1},...,x^{d})$, $\bR^{d+1}
=\{(t,x):t\in\bR,x\in\bR^{d}\}$. Let $(\Omega,\cF,P)$ be a 
complete probability space and
let $\{\cF_{t}\}$ be an increasing filtration of 
$\sigma$-fields $\cF_{t}\subset \cF$, that are complete.
Let  
$w_{t}$ be a $d_{1}$-dimensional Wiener process relative to
$\{\cF_{t}\}$, where $d_{1}\geq d$.

Assume that on $\bR^{d+1}$ we are given 
Borel
functions $\bR^{d}$-valued  $b $
and $d\times d_{1}$-matrix valued $\sigma$. One of the main topics of this book
is the investigation of various issues
related to
the It\^o equation 
 \begin{equation}
                         \label{6.15.2}
 x_{s}=x +\int_{0}^{s}\sigma(t+r,x_{r})\,dw_{r}+  
\int_{0}^{s}b(t+r,x_{r})\,dr, 
 \end{equation}
 where 
 $(t,x) \in\bR^{d+1}$ is nonrandom. 
The diffusion matrix $a:=\sigma\sigma^{*}$
\index{$S$@Miscelenea!$a_{\pm}$@$a:=\sigma\sigma^{*}$}%
will be often assumed to take values
in $\bS_{\delta}$,
\index{$B$@Sets!$\bS_{\delta}$}%
\index{$B$@Sets!$\bS_{0}$}%
 that is the subset
of the set of symmetric nonnegative $d\times d$-matrices
$\bS_{0}$  
whose eigenvalues are in $[\delta,
\delta^{-1}]$, where $\delta\in(0,1]$
\index{$S$@Miscelenea!$\delta\in(0,1]$}%
is fixed throughout the book.
The drift coefficient $b$ is allowed
to have some singularities expressed 
trough its belonging to either
  Lebesgue spaces $L_{q,p}$  or to  Morrey spaces $E_{q,p,\beta}$. The latter may be
less familiar to some readers
and it is good to know that $b$ with
$|b(x)|\leq 1/|x|$ is in the Morrey space
$E_{q,p,1}$ with any $p\in(1,d)$.

We are interested in just any solution of
(\ref{6.15.2}) ({\em weak solutions\/})
\index{$D$@Processes!weak solution}%
 and also in the
 so-called {\em strong solutions\/}, 
\index{$D$@Processes!strong solution}%
that are
solutions such that, for each $t\geq0$, $x_{t}$ is $\cF^{w}_{t}$-measurable,
\index{$B$@Sets!$\cF^{w}_{t}$}%
where $\cF^{w}_{t}$ is the completion of $\sigma(w_{s}:s\leq t)$. 

In the recent past the author published
several articles related to equations
like (\ref{6.15.2}) when sometimes $b$ also depended explicitly on  
  $\omega$.
They started with imposing, {\em different\/}
in different articles, conditions on
$b$ in terms of Lebesgue spaces but
the results of previous articles were used
in the next ones. This caused using the
arguments like ``as in the proof of...it is not hard to see''. This created a pyramid
of through references which became very disturbing
at some moment and at the same time the author realized 
that the initial
conditions  on $b$ in terms of Lebesgue
spaces  are  very inconvenient and  are  not satisfied in a number of situations where,
however, the implications of these condition,
the main of which is contained in Assumption
\ref{assumption 8.19.2}, are still true.  Therefore,
there was a need to restructure the
``pyramid'' on a new basis and this was 
one of   motivations of the book.

Another very strong motivation and inspiration came from
remarkable papers by F.~Flandoli, M.~Gubinelli, and M.~Maurelli
(\cite{BFGM_19}),  M.~R\"ockner and Guohuan  Zhao (\cite{RZ_20}, \cite{RZ_25},   and
a quite recent article by D. Kinzebulatov and K.R. Madou (\cite{KM_24},   
in which the authors gave a new 
powerful impulse
to the theory of weak and strong solutions
of It\^o equations with singular $b$, even though
in their papers $\sigma=(\delta^{ij})$.
These papers and also 
\cite{Ki_24}, \cite{Ki_25} and the references 
found in them
contain the most advanced
information about solvability of It\^o equations with singular $b$ if $\sigma=(\delta^{ij})$. However, none of the results
in this book   is covered by the results in the articles mentioned above for the simple reason that
our diffusion coefficient is not constant.
For that matter, the majority of our main results are new even if
the drift part is {\em zero\/}.

Since  the classical work
by K.~It\^o \cite{It_51_1} (1951), who introduced stochastic integral
equations even more general than
(\ref{6.15.2}) and proved their strong
solvability and also proved the famous
It\^o's formula, all kinds of results
from the theory of second-order elliptic and parabolic PDEs became conveniently
available to probabilists dealing with
 stochastic processes. On the other hand,
I.I.~Gikhman \cite{Gi_47} (1947, the proofs are in \cite{Gi_51}, 1951) by using the theory
of stochastic {\em differential\/} equations
proved the {\em first\/} result on the {\em classical
solvability\/} of second order parabolic degenerate equations. This showed the power of probabilistic approach to PDEs.
In this book we will see more of such
interplay between stochastic  
equations  and PDEs.

The book consists of six chapters
and an appendix. In the first 
chapter we start by deriving from 
Aleksandrov parabolic $L_{d+1}$-estimates
and elliptic $L_{d}$-estimates
a mixed-norm estimate in $L_{(q,p)}$ following an idea of A.I.~Nazarov. This leads us to
proving It\^o's formula for functions in $W^{1,2}
_{(p,q)}$ with $d/p+2/q< 2$, and to the uniqueness results for the 
   second- and first-order parabolic equations in Sobolev spaces.
Then we deal with the solvability of uniformly
nondegenerate stochastic equations with
measurable coefficients and drift
of class $L_{(q,p)}$, $d/p+1/q\leq1$. By using
the Skorokhod approach and selecting the solutions, among many, in an appropriate way
we construct time inhomogeneous Markov diffusion processes with
trajectories being solutions of the given stochastic equation.

Chapter 2 deals with some properties of stochastic integrals with uniformly bounded and nondegenerate diffusion which are not necessarily solutions of It\^o's equations. Here we present
the estimates of the time spent by such processes in the space-time sets of small measure, which long time ago were proved by the author together with M.V. Safonov for diffusion
processes. By using an idea of Fabes-Stroock
we show the higher summability of Green's functions of stochastic integrals, which is higher than what the Aleksandrov estimates
guarantee. The Fabes-Stroock paper 
deals with time-homogeneous case and is based on Gehring's lemma. Our argument is based on the parabolic analog of Gehring's lemma.

Starting with Chapter 3 we focus on Markov
diffusion processes corresponding to stochastic
equations. The main emphasis is on studying
the particular conditions on the drift $b$ which
allow  singularities and at the same time
do  not distort the pure diffusion too much
so that such properties as the H\"older
continuity and Harnack inequality of the
caloric functions associated with the process
are still valid. These conditions on $b$
are expressed in terms of the Morrey spaces
and, as we show, they are practical in the sense that the corresponding Markov diffusion
processes do exist.
We also prove the Liouville theorem, that the probabilistic
Green's functions are summable to a small negative power and that the probabilistic solutions
of PDEs are {\em $W^{1,2}_{(q,p)}$-viscosity\/} solutions.

In Chapter 4 we give some applications to the
theory of elliptic and parabolic equations when
$b$ is in the mixed-norm Morrey space with $\sfd_{0}/p+1/q\leq 1$, where $\sfd_{0}\in(d/2,2)$ is the Fabes-Stroock constant.
Extended Aleksandrov estimates, Harnack
inequality, and H\"older continuity are given
in the analytic rather than probabilistic form.

In Chapters 5 and 6  we deal with the weak and strong solutions, respectively. 
We attract the attention of the reader 
to the following disclaimer about the order of
summation while defining mixed-norm
Lebesgue spaces found in the beginning of Chapter \ref{chapter 3.7.1}:

{\em One of the ways to choose the norm
is fixed throughout the rest of the book
unless specifically stated otherwise.
We will be referring to some results that are
proved elsewhere for only one of the norms   \eqref{3.27.3} or \eqref{4.3.2}. In such
situations we mean that the result, we are referring to, actually, holds
for both norms and is proved by 
insignificant changes in the original proof.
This is, for instance, explicitly mentioned 
and underlined in \cite{Kr_26_1}.}

The  main results in Chapters 5 and 6
are based on some analytic facts,  
the exposition of which is done in pure PDE
terms. In chapter 5 this is the theory of Morrey-Sobolev spaces which provides the ground
for applying It\^o's formula and showing the weak uniqueness of solutions. 
We show that the weakly unique solutions
form a strong Markov strong Feller diffusion 
processes.
In Chapter 6 the main analytic fact is
Theorem \ref {theorem 6.21.1} allowing to estimate solutions
of some special parabolic PDEs and prove
the existence, uniqueness, and the differentiability with respect to initial data
of strong solutions of \eqref{6.15.2}.
It is worth mentioning that as a rule our
uniqueness results bear on solutions
(which are shown to exist) such that their
potentials admit certain estimates. However,
there are also unconditional results.
For instance, if $d=d_{1}$, $\sigma=(\delta^{ij})$ and $|b|=cf$, where
$$
f(t,x)=I_{t>0}\frac{1}{|x|^{\gamma}(|x|+\sqrt t)^{1-\gamma}
}, \quad \gamma\in\Big(\frac{d}{d+1},\frac{2d}{2d+1}\Big)
$$
and the constant $c>0$ is sufficiently small,
then  equation \eqref{6.15.2} with these data
has a weak solution and each solution has
the same finite dimensional distributions
(see Remark \ref{remark 1.28.1}).
This solution is shown to be a strong one
in Remark \ref{remark 3.21,3}.

  We finish the introduction with some notation 
and stipulations. Throughout the book
  the summation
convention over repeated indices
(even at the same level) is enforced.

  In the proofs of various results  we use
the symbol $N$  to denote finite 
nonnegative constants
which may change from one occurrence to another and,
 if in the statement of a result there are constants
called $N$ which are claimed to depend only on certain
parameters, then in the proof of the result
the constants $N$ also depend only on the same
parameters unless specifically stated otherwise.
Of course, if we write 
$
N=N(...),
$
 this means that $N$ depends only
on what is inside the parentheses. 
 
Introduce  
$$
B_{R}(x)=\{y:|y-x|<R\},\quad B_{R}=B_{R}(0),
$$
$$
C_{T,R}=[0,T)\times B_{R},\quad C_{T,R}(t,x)=(t,x)+C_{T,R},
\quad C_{R}(t,x)=C_{R^{2},R}(t,x),
  $$
 $C_{R}=C_{R}(0,0)$, and 
\index{$B$@Sets!$B_{R}(x)$}%
\index{$B$@Sets!$B_{R}$}%
\index{$B$@Sets!$C_{T,R}$}%
\index{$B$@Sets!$C_{T,R}(t,x)$}%
\index{$B$@Sets!$C_{R}(t,x)$}%
\index{$B$@Sets!$C_{R}$}%
\index{$B$@Sets!$\bC_{R}$}%
\index{$B$@Sets!$\bC$}%
\index{$B$@Sets!$\bB_{R}$}%
\index{$B$@Sets!$\bB$}%
let $\bB_{R}$ be the collection of $B_{R}(x)$
and $\bC_{R}$ be the collection of $C_{R}(t,x)$. 
Set $\bB=\bigcup_{R>0}\bB_{R}$, $\bC=\bigcup_{R>0}\bC_{R}$.
$$
\quad a_{\pm}=a^{\pm}=(1/2)(|a|\pm a),
$$
$$
D_{i}u=u_{x^{i}}=\frac{\partial}{\partial x^{i}},
\quad Du=(D_{i}u),
\quad D_{ij}u=u_{x^{i} x^{j}}=D_{i}D_{j}u ,
$$
$$
 D^{2}u=(D_{ij}u), \quad \partial_{t}=\frac{\partial}{\partial t}.
$$
By derivatives of functions we always mean the Sobolev derivatives.

We use  
\index{$S$@Miscelenea!$a_{\pm}$}%
\index{$C$@Operators!$D_{i}$}%
\index{$C$@Operators!$u_{x^{i}}$}%
\index{$C$@Operators!$u_{x^{i}x^{j}}$}%
\index{$C$@Operators!$D_{ij}$}%
\index{$C$@Operators!$D$}%
\index{$C$@Operators!$D^{2}$}%
\index{$C$@Operators!$\partial_{t}$}%
\index{$C$@Operators!$u^{(\varepsilon)}$}%
the notation $u^{(\varepsilon)}=u*\zeta_{\varepsilon}$,
where
$\zeta_{\varepsilon}(x)=\varepsilon^{-d}\zeta(x/\varepsilon)$,
$\varepsilon>0$,
and $\zeta$ is a nonnegative $C^{\infty}$-function with
support in $B_{1}$ whose integral is equal to one.
 
In Chapters \ref{chapter 3.24,2} through \ref{chapter 3.24,1},
for $p,q\in[1,\infty]$, $L_{q,p}$ is the space
of functions $f$ on $\bR^{d+1}$ with finite norm
$$
\|f\|_{L_{q,p}}=\Big(\int_{\bR}\Big(\int_{\bR^{d}}|f|^{p}
\,dx\Big)^{q/p}\,dt\Big)^{1/q}.
$$
In subsequent chapters the meaning of $L_{p,q}$
is defined in Chapter \ref{chapter 3.7.1}.
We write $\|u,v,...\|_{L_{q,p}}$ to mean
 the sum of the $L_{q,p}$-norms of what is inside.
\index{$N$@Norms!$"|"|u,v,..."|"|_{L_{q,p}}$}%

If $\cO$ is a Borel subset 
\index{$A$@Sets of functions!$L_{q,p}(\cO)$}%
of $\bR^{d+1}$ we set
$L_{q,p}(\cO)=\{f:f{I_{\cO}}\in L_{q,p}\}$.
Next $L_{p}=L_{p,p}$, $L_{p}(\bR^{d})$ has the usual sense as well as $L_{p}(\cO)$
\index{$A$@Sets of functions!$L_{p}$}%
\index{$A$@Sets of functions!$L_{p}(\bR^{d})$}%
\index{$A$@Sets of functions!$L_{p}(\cO)$}%
for
Borel $\cO$ in $\bR^{d}$ or $\bR^{d+1}$.
 By $L_{(q,p)}$
we mean 
\index{$A$@Sets of functions!$L_{q,p}$}%
\index{$N$@Norms!$"|"|f"|"|_{L_{q,p}}$}%
\index{$A$@Sets of functions!$L_{(q,p)}$}%
\index{$N$@Norms!$"|"|f"|"|_{L_{(q,p)}}$}%
the space of functions on $\bR^{d+1}$
such that
$$
\|f\|_{L_{(q,p)}}:=\begin{cases}
\|f\|_{L_{q,p}}<\infty\quad\text{if}
\quad p>q,\\ 
\Big(\int_{\bR^{d}}\Big(\int_{\bR}|f|^{q}
\,dt\Big)^{p/q}\,dx\Big)^{1/p}<\infty\quad\text{if}
\quad q\geq p.
\end{cases}
$$
Of course, formally this definition makes sense
only for finite $p,q$. We extend it to cover
infinite values in a well-known way.

Accordingly we define $W^{1,2}_{(q,p)}$
as the space
\index{$A$@Sets of functions!$W^{1,2}_{(q,p)}$}%
\index{$A$@Sets of functions!$W^{1,2}_{q,p}$}%
 of functions $u$ such that
$u$ and its
the Sobolev derivatives $\partial_{t}u,
D^{2}u,Du$ are in $L_{(q,p)}$. 
Similarly, $W^{1,2}_{q,p}$ is defined
as the space  of functions $u$ such that
$u$ and its
the Sobolev derivatives $\partial_{t}u,
D^{2}u,Du$ are in $L_{q,p}$.
The norms in
$W^{1,2}_{(q,p)}$ and $W^{1,2}_{q,p}$ are defined in a natural way. If $\cO$ is an open
sunset of $\bR^{d+1}$ the space 
$W^{1,2}_{(q,p)}(\cO)$ is defined
the space  of functions $u$ such 
\index{$A$@Sets of functions!$W^{1,2}_{(q,p)}(\cO$}%
\index{$A$@Sets of functions!$W^{1,2}_{q,p}(\cO)$}%
that
$uI_{\cO}$ and   $I_{\cO}\partial_{t}u,
I_{\cO}D^{2}u,I_{\cO}Du$ are in $L_{q,p}$.
Similarly, $W^{1,2}_{q,p}(\cO)$ is defined.
The norms in those spaces are introduced
in a natural way.
 
If $(\Omega,\cF,P)$ is a probability space
and a $\sigma$-field $\cN\subset\cF$ we often
use the
\index{$C$@Operators!$E_{\cN}$}%
 notation
$$
E_{\cN }\quad{\text{for}}
\quad E\{\cdot\mid
\cN\}.
$$
 
The following notation may look too complex to absorb at once.
However, in most cases of using 
them we remind the reader
the meaning of these notation and the most important is the difference between $\tau$ and
$\tau'$.
In the situation when an underlying
$\bR^{d}$-valued continuous random 
process $x_{t},t\geq0$, is involved
and $\cO$ is an open subset of $\bR^{d}$,
we define $\tau'_{\cO} $
\index{$S$@Miscelenea!$\tau'_{\cO}$}%
 as the first
time $ x_{t}$ exits from $\cO$ ($=\infty$ if $x_{t}$ never leaves $\cO$).
 If $\cO=B_{R}(x)$ we write $\tau'_{R}(x)$
\index{$S$@Miscelenea!$\tau'_{R}(x)$}%
 in place
of $\tau'_{B_{R}(x)}$ and if $x=0$ we drop it
 in these notation.

For $t,R \in(0,\infty)$ and $x\in\bR^{d}$
denote by $\theta_{t}\tau'_{R}(x)$ the first exit time of
\index{$S$@Miscelenea!$\theta_{t}\tau'_{R}(x)$}%
the process
$ x_{t+s} $, $s\geq0$, from $B_{R}(x)$ or,
equivalently, the first exit time of $x_{t+s}-x $
from $B_{R} $.
 Set 
$\theta_{0}\tau'_{R}(x)=\tau'_{R}(x)$ the first time
$ x_{t}-x$ exits from $B_{R}$. 
Next, 
\index{$S$@Miscelenea!$\tau_{R}(x)$}%
\index{$S$@Miscelenea!$\theta_{t}\tau_{R}(x)$}%
we set 
$$
\tau_{R}(x)=R^{2}\wedge
\tau'_{R}(x),\quad \theta_{t}\tau_{R}(x)=R^{2}\wedge \theta_{t}\tau'_{R}(x).
$$

If we are given a stopping time $\tau$,   set $\theta_{\tau}\tau'_{R}(x)(\omega)=I_{\tau<\infty}\theta_{\tau(\omega)}\tau'_{R}(x)$. Similar sense is 
\index{$S$@Miscelenea!$\theta_{\tau}\tau'_{R}(x)$}%
given to $\theta_{\tau}\tau_{R}(x)(\omega)$.
If $x=0$, we drop it in the above notation.
It is useful to note that $t+\theta_{t}\tau_{R}(x)$
are stopping times and
\begin{equation}
                                      \label{10.17.5}
\theta_{t}\tau'_{R }(x)\leq \theta_{t}\tau'_{2R },
\end{equation}
whenever $|x|<R$, because $B_R(x)\subset B_{2R}$.

By $\bar\tau'_{\cO}(x)$
\index{$S$@Miscelenea!$\tau'_{\cO}$@$\bar\tau'_{\cO}(x)$}%
\index{$S$@Miscelenea!$\tau'_{\cO}$@$\bar\tau'_{R}(x)$}%
 we mean the first
exit time of $x+x_{t}$ from $\cO$,
$\bar\tau'_{R}(x)$ is the first exit time
of $x+x_{t}$ from $B_{R}$. Note that
$\bar\tau'_{R}(x)= \tau'_{R}(-x)$,
$\bar\tau'_{R}(0)= \tau'_{R}(0)=\tau'_{R}$.

If $F$ is a closed subset of $\bR^{d}$,
by $\gamma_{F}$
\index{$S$@Miscelenea!$\tau_{\cO}$@$\theta_{t}\gamma_{F}$}%
\index{$S$@Miscelenea!$\theta_{\tau}\tau'_{R}(x)$}%
 we denote the first time
$x_{t}$ hits $F$. By $\theta_{t}\gamma_{F}$ we mean the first
time $x_{t+s},s\geq0, $ hits $F$. We use the same
agreement as above to define $\gamma_{R}(x)$,
$\theta_{t}\gamma_{R}(x)$   if $F=\bar B_{R}(x)$  and if $x=0$.

We denote by $\bar\gamma_{F}(x)$
\index{$S$@Miscelenea!$\tau_{\cO}$@$\bar\gamma_{R}(x)$}%
\index{$S$@Miscelenea!$\tau_{\cO}$@$\bar\gamma_{F}(x)$}%
\index{$S$@Miscelenea!$\tau_{\cO}$@$\bar\gamma_{R}(x)$}%
\index{$S$@Miscelenea!$\tau_{\cO}$@$\gamma_{R}(x)$}%
\index{$S$@Miscelenea!$\tau_{\cO}$@$\theta_{t}\gamma_{R}(x)$}%
the first time $x+x_{t}$ hits $F$,
$\bar\gamma_{R}(x)$ is the first time $x+x_{t}$ hits $\bar B_{R}$.

Other notation are introduced wherever
appropriate and the list of notation 
is found at the end of the book.

 \mainmatter

\mychapter[Preliminaries]{Preliminaries}

                    \label{chapter 3.24,2}
 
\mysection[Aleksandrov's parabolic estimates]{The Aleksandrov estimates
for potentials of stochastic integrals in  
$L_{(q,p)}$, $d/p+1/q\leq1$}
                    \label{section 8.19.1}

The first parabolic Aleksandrov estimates
for processes given by $dx_{t}=\sigma_{t}\,dw_{t}+b_{t}\,dt$ with bounded $b_{t}$ 
and possibly degenerating $\sigma_{t}$ appeared in
\cite{Kr_77}. In \cite{NU_85} they were
extended to the case that $b_{t}=b(t,x_{t})$ with
$b(t,x)$ of class $L_{d+1}$ and
$\sigma$ uniformly nondegenerate. The author
in \cite{Kr_86} developed a new technique  
to achieve somewhat more precise result
allowing again degeneration of $\sigma$
and by adding to it an interpolation argument
 A.I. Nazarov in \cite{Na_87} 
developed mixed-norm estimates. 
Here we follow \cite{Kr_86} where
we somewhat sharpened the probabilistic versions of some arguments in \cite{Na_87}, that by the way,   treats the problem in PDEs terms rather than in the 
probabilistic ones.  This, in particular,
allows us to prove the maximum principle
for the first-order parabolic equations with singular $b$
(see Theorem \ref{theorem 12.20.2}).

Let $(\Omega,\cF,P)$ be a complete probability
space, let $\cF_{t}, t\geq0$, be an increasing family of
complete $\sigma$-fields $\cF_{t}\subset\cF$, $t\geq0$,
let $m_{t}$ be an $\bR^{d}$-valued continuous local martingale
relative to $\cF_{t}$, let  $A_{t}$ be a continuous
$\cF_{t}$-adapted
nondecreasing process, and let $B_{t}$ be a continuous $\bR^{d}$-valued
$\cF_{t}$-adapted process which has finite variation (a.e.)
on each finite time interval. Assume that
$$
A_{0}=0,\quad m_{0}=B_{0}=0,\quad d\langle m
 \rangle_{t} \ll  dA_{t}
$$
and that we are also given progressively measurable relative to
$\cF_{t}$ nonnegative processes $r_{t}$ and $c_{t}$. Finally, take  
  $\cF_{0}$-measurable $\bR^{d}$-valued $x_{0}$ and $\sft_{0}\in\bR$ and 
introduce
$$
x_{t}=x_{0}+m_{t}+B_{t},\quad \sft_{t}=\sft_{0}+\int_{0}^{t}r_{s}\,dA_{s},
\quad \phi_{r,t}=\int_{r}^{t}c_{s}\,dA_{s},\quad \phi_{t}=\phi_{0,t},
$$
$$
a^{ij}_{t}= \frac{d\langle m^{i},m^{j}\rangle
_{t}}
{dA_{t}}.
$$

Here is Lemma 1.4.1 of \cite{Kr_25}.
\begin{lemma}
                                   \label{lemma 4.23.3}
Let $\gamma,\tau$ be   $\{\cF_{t}\}$-stopping times,
$\gamma\geq\tau$, and set   
$$
A=E_{\cF_{\tau}}\int_{\tau}^{\gamma}e^{-\phi_{\tau,t}} \tr a_{t} \,dA_{t},\quad 
B=E_{\cF_{\tau}}\int_{\tau}^{\gamma}e^{-\phi_{\tau,t}}\,|dB_{t}|.
$$
Then for any Borel $f(t,x)\geq0$ we have
$$
E_{\cF_{\tau}}\int_{\tau}^{\gamma}e^{-\phi_{\tau,t}}(r_{t}\det a_{t})^{1/(d+1)}
f(\sft_{t},x_{t})\,dA_{t}
$$
\begin{equation}
                                            \label{4.24.1}
\leq N(d) (B^{2}+A)^{d/(2d+2)}\|f\|_{L_{ d+1} }.  
\end{equation}
 
\end{lemma} 

The expressions like the left-hand side 
of \eqref{4.24.1} are called {\em potentials\/}
of function $f$.
\index{$S$@Miscelenea!potentials}
 Lemma \ref{lemma 4.23.3}  implies the following important result.

\begin{theorem}
                 \label{theorem 6.3,1}
Suppose that,
\begin{equation}
                         \label{5.4,1}
\gamma-\tau\leq \theta_{\tau}\tau'_{ R}(x_{\tau}),
\end{equation}
where $\theta_{\tau}\tau'_{ R}(x_{\tau})$ is  the first exit time of $x_{ \tau +s}$
from $B_{R}(x_{\tau})$, then $A\leq 2B^{2}+RB\leq 3B^{2}+R^{2}$,
so that $A$ in \eqref{4.24.1} can be replaced
with $R$. If, in addition, $B\leq KR$,
where $K$ is a constant, then
$$
E_{\cF_{\tau}}\int_{\tau}^{\gamma}e^{-\phi_{\tau,t}}(r_{t}\det a_{t})^{1/(d+1)}
f(\sft_{t},x_{t})\,dA_{t}
$$
$$
\leq N(d)(1+K)^{d/(d+1)}R^{d/(d+1)}
\|f\|_{L_{ d+1} }.
$$

\end{theorem}

Proof. By It\^o's formula
$$
|x_{t\wedge\gamma}-x_{t\wedge\tau}|^{2} = 2\int_{t\wedge\tau}^{t\wedge\gamma}
e^{-\phi_{t\wedge\tau,s}}\big[\tr a_{s}-c_{s}|x_{s}-x_{t\wedge\tau}|^{2}\big]\,dA_{s}
$$
$$
+2
\int_{t\wedge\tau}^{t\wedge\gamma}
e^{-\phi_{t\wedge\tau,s}}(x^{i}_{s}-x^{i}_{t\wedge\tau})\,dB^{i}_{s}
+2\int_{t\wedge\tau}^{t\wedge\gamma}e^{-\phi_{t\wedge\tau,s}}
(x^{i}_{s}-x^{i}_{t\wedge\tau})\,dm^{i}_{s}.
$$
Here the stochastic integral is a local martingale.
Therefore, by replacing $t$ with $ \tau_{n}$ for an appropriate sequence
of stopping times  $\tau_{n}\to\infty$, then taking expectations and using that
$$
|x_{\tau_{n}\wedge\gamma}-x_{\tau_{n}\wedge\tau}|^{2} \leq  R^{2},\quad c_{s}\geq0.
$$
$$
E_{\cF_{\tau}}\Big|\int_{\tau_{n}\wedge\tau}^{\tau_{n}\wedge\gamma}
e^{-\phi_{\tau_{n}\wedge\tau,s}}(x^{i}_{s}-x^{i}_{\tau_{n}\wedge\tau})\,dB^{i}_{s}\Big|
$$
$$
=I_{\tau_{n}\geq\tau}E_{\cF_{\tau}}\Big|\int_{ \tau}^{\tau_{n}\wedge\gamma}
e^{-\phi_{\tau,s}}(x^{i}_{s}-x^{i}_{ \tau})\,dB^{i}_{s}\Big|
\leq RB,
$$
we find
$$
2E_{\cF_{\tau}}\int_{\tau_{n}\wedge\tau}^{\tau_{n}\wedge\gamma}
e^{-\phi_{\tau,s}} \tr a_{s} \,dA_{s}
\leq 3R^{2}+2RB.
$$
Sending $n\to\infty$ yields 
the desired result. \qed

Next comes Lemma 1.4.2 of \cite{Kr_25}.

\begin{lemma}
                                            \label{lemma 5.6.1}
In the notation of Lemma \ref{lemma 4.23.3}
for any Borel $f( x)\geq0$ we have
\begin{equation}
                                            \label{5.6.1}
E_{\cF_{\tau}}\int_{\tau}^{\gamma}e^{-\phi_{\tau,t}} (\det a_{t}) ^{1/d}
f( x_{t})\,dA_{t}\leq N(d)
 (B^{2}+A)^{1/2}\|f\|_{L_{d }(\bR^{d })}.
\end{equation}
\end{lemma} 

The following corollary provides the case
when $c_{t}$ plays a major role, regardless of how irregular $a$ and $B$
could be.

\begin{corollary}
             \label{corollary 10.5.1}
Suppose that $|dB_{t}|\ll dA_{t}$ and   there is a constant  $\mu
\geq0$ such  that $\mu c_{t}\geq \tr a_{t}$ and  $\sqrt \mu c_{t}\geq |dB_{t}|/dA_{t}$. Then $A\leq \mu, B\leq\sqrt\mu$,
so that $B^{2}+A$ in \eqref{4.24.1} and
\eqref{5.6.1} can be replaced with $2\mu$.
\end{corollary}

 \begin{definition}
            \label{definition 2.3.1}
Let $\mu, q ,p \in[1,\infty]$. 
\index{$S$@Miscelenea!properly tight}%
\index{$S$@Miscelenea!$\nu(\mu,q,p)$}%
We say that $(\mu, q ,p)$ are {\em properly tight\/} if
$$
\nu(\mu,q,p):=1-\frac{\mu}{p}-\frac{1}{q}\geq 0.
$$
\end{definition}

\begin{theorem}
                                         \label{theorem 5.5.1}
Assume  the notation of Lemma \ref{lemma 4.23.3} and let
$(d,q,p)$ be properly tight.  Then for any Borel $f(t, x)\geq0$   we have  
\begin{equation}
                                                   \label{5.6.40}
I(p,q,f):= E_{\cF_{\tau}}\int_{\tau}^{\gamma}e^{-\phi_{\tau,t}}\kappa_{t}
f(\sft_{t},x_{t})\,dA_{t}\leq N(d)  (A+B^{2})^{ d/(2p)}
\|f\|_{L_{(q,p)}},
\end{equation}
where 
$\kappa_{t}=
r_{t}^{1/q}(\det a_{t})^{1/p}c^{\nu}_{t}$ 
and for any $\alpha\geq0$ we set $\alpha^{0}=1$ \(say, if $\nu=0$\).
\end{theorem}

Proof. We use an idea of Nazarov from \cite{Na_87}. 
If $\nu=1$ ($p=q=\infty$), \eqref{5.6.40} is obvious with $N=1$,
By H\"older's inequality, if $1>\nu>0$,  
$$
I(p,q,f)\leq  \Big(I(p (1-\nu ),q (1-\nu),
f^{1/(1-\nu)})
\Big)^{1-\nu}.
$$
It follows that it suffices to concentrate on   $\nu=0$.
Then we observe that if $q=\infty$, then $p=d$ and
$$
\|f\|^{p}_{L_{(q,p)}}=\int_{\bR^{d}}\sup_{t }f^{d}(t,x)\,dx.
$$
In that case \eqref{5.6.40} follows from Lemma \ref{lemma 5.6.1}.  
If $p=\infty$, then $q=1$,
and
$$
I(p,q,f)= E_{\cF_{\tau}}\int_{\tau}^{\gamma}e^{-\phi_{\tau,t}}r_{t} f (\sft_{t},x_{t})
\,dA_{t}\leq  E_{\cF_{\tau}}\int_{\tau}^{\gamma} \sup_{x}f (\sft_{t},x)
\,d\sft_{t}
$$
$$
\leq \int_{\bR}\sup_{x}f (t,x)\,dt= \|f\| _{L_{(q,p)}}.
$$
In the third simple situation when $q=p=d+1$ estimate
\eqref{5.6.4} follows from Lemma \ref{lemma 4.23.3}.
We prove the lemma in the remaining cases of $p,q<\infty$ by interpolating between
the above ones.

If  $p>q$ (and hence $p>d+1$) we take a
  nonnegative
function $h(t )$ such that $(hf)/h=f$ ($0/0:=0$) and
 use
$$
 r_{t}^{1/q}(\det a_{t})^{1/p}f=
 \Big(r_{t}^{1/q-1/p}h^{-1}\Big)\Big((r\det a_{t})^{1/p}fh\Big)
$$
along with H\"older's inequality. By  
performing simple manipulations we find
$$
I(p,q,f)\leq IJ
$$
\begin{equation}
                                                   \label{5.6.50}
:= 
\Big(I(\infty,1,
h^{-p /(p-d-1)})
\Big)^{ (p-d-1)/p}\Big(I(d+1 ,d+1  ,
(hf)^{p/(d+1)})
\Big)^{(d+1)/p}.
\end{equation}
Here  
$$
I\leq \Big(\int_{\bR} h^{-p /(p-d-1)}(t)\,dt\Big)^{ (p-d-1)/p}.
$$
Also
$$
J\leq
N(d)(B^{2}+A)^{d/(2p)}\|(hf)^{p/(d+1)}\|^{(d+1)/p}_{L_{ d+1}(\bR^{d+1})}
$$
$$
=N(d)(B^{2}+A)^{d/(2p)}\Big(\int_{\bR} 
\Big(\int_{\bR^{d}}f^{p}(t,x)\,dx\Big)h^{p}(t)\,dt\Big)^{1/p}.
$$
We now choose $h$ so that
$$
h^{-p /(p-d-1)}(t)=\Big(\int_{\bR^{d}}f^{p}(t,x)\,dx\Big)h^{p}(t).
$$
Then both estimates become
$$
J\leq N(d)(B^{2}+A)^{d/(2p)}\|f\|^{q/p}_{L_{(q,p)}},\quad
I\leq \|f\|^{q(p-d-1)/p}_{L_{(q,p)}}
$$
and coming back to \eqref{5.6.50} we get \eqref{5.6.40}.

In the remaining case $q>p$ (and $q>d+1$)
we use
$$
r_{t}^{1/q}(\det a_{t})^{1/p}f=
 \Big((\det a_{t})^{1/p-1/q}h^{-1}\Big)\Big((r\det a_{t})^{1/q}fh\Big).
$$
This time for $h=h(x)$
$$
I(p,q,f)\leq IJ
$$
\begin{equation}
                                                   \label{5.6.6}
:= 
\Big(I(d ,\infty,
h^{-q /(q-d-1)})
\Big)^{ (q-d-1)/q}\Big(I(d+1 ,d+1  ,
(hf)^{q/(d+1)})
\Big)^{(d+1)/q}.
\end{equation}
 Here
$$
I\leq N(d)(B^{2}+A)^{(q-d-1)/(2q) }\Big(\int_{\bR^{d}}h^{-qd/(q-d-1)}
(x)\,dx\Big)^{(q-d-1)/(qd )},
$$
$$
J\leq N(d)(B^{2}+A)^{d/(2q)) }\Big(\int_{\bR^{d}}h^{q}(x)
\Big(\int_{\bR}  f^{q}(t,x)\,dt\Big)\,dx\Big)^{1/q}.
$$
We choose $h$ so that
$$
h^{-qd/(q-d-1)}
(x)=h^{q}(x)
\Big(\int_{\bR} f^{q}(t,x)\,dt\Big)
$$
and then easily come to \eqref{5.6.4}. The theorem
is proved. \qed

The rigorous statement and the proof of the following
corollary, not used in the subsequent text, is left to the
interested reader.
\begin{corollary}
                                         \label{corollary 5.1.1}
Introduce a measure (Green's measure) on Borel 
subsets $\Gamma$ of $\bR^{d+1}$
\index{$S$@Miscelenea!Green's measure}%
by the formula
$$
G(\Gamma)=
E_{\cF_{\tau}}\int_{\tau}^{\gamma}e^{-\phi_{\tau,t}}\kappa_{t}
I_{\Gamma}(\sft_{t},x_{t})\,dA_{t}.
$$
Assume that $A,B<\infty$ and set $p'=p/(p-1),q'=q/(q-1)$.
 Then $G(\Gamma)$ is absolutely continuous and
its  density $G(t,x)$ is such that, if $p\geq q$,
$$
\Big(\int_{\bR} \Big(\int_{\bR^{d}}G^{p'}(t,x)\,dx\Big)^{q'/p'}
dt\Big)^{1/q'}
$$
and, if $p\leq q$,
$$
\Big(\int_{\bR^{d}}\Big(\int_{\bR} G^{q'}(t,x)\,dt\Big)^{p'/q'}
dx\Big)^{1/q'}
$$
is dominated by
$$
  N(d)(B^{2}+A)^{ d/(2p)}.
$$
\end{corollary}

\begin{assumption}
                    \label{assumption 9.24.1}
We have that $(d,q_{0},p_{0})$ are properly tight,
  $d|B_{t}|\ll dA_{t}$
 and there exists a Borel $h(t,x)$
such that   ($P(d\omega)\times dA_{t}$-a.e.)
$$
|b_{t}|\leq \kappa^{0}_{t}h(\sft_{t}
,x_{t}),
$$
where $b_{t}=dB_{t}/dA_{t}$ and $\kappa^{0}_{t}=
r_{t}^{1/q_{0}}(\det a_{t})^{1/p_{0}}c^{\nu_{0}}_{t}$
 ($\alpha^{0}\equiv1,\alpha\geq0$), $\nu_{0}=\nu(d,q_{0},p_{0})$. 
\end{assumption}

\begin{theorem}
                                         \label{theorem 5.5.2}   
 Under Assumption \ref{assumption 9.24.1}
suppose that $(d,q,p)$ are properly tight and, if $p_{0}=d$, then either

(a) $\|h\|_{L_{(\infty,d)}}< \varepsilon(d)$,
where $\varepsilon(d)\in(0,1)$ is defined in the proof,
or 

(b) $A_{t}=t$ and $\gamma\leq \tau+\theta_{\tau}\tau'_{R}(x)$, for some $R,x$.

Then for any Borel $f(t, x)\geq0$   we have
\begin{equation}
                                            \label{5.6.4}
I(p,q,f):= E_{\cF_{\tau}}\int_{\tau}^{\gamma}e^{-\phi_{\tau,t}}\kappa_{t}
f( \sft _{t},x_{t})\,dA_{t}\leq N(d,p_{0} )C^{ d/(2p)}
\|f\|_{L_{(q,p)}},
\end{equation}
where 
\begin{equation}
                           \label{9.14,1}
\kappa_{t}=r_{t}^{1/q}(\det a_{t})^{1/p}c^{\nu}_{t},  
\end{equation}
  for any number $\alpha\geq0$ we set $\alpha^{0}=1$ (say, if $\nu=0$) and in all cases but  (b)
$$
C=  A+ \|h\|_{L_{(q_{0},p_{0})}}^{ 2p_{0}/(p_{0}-d)  }\quad\text{if}
\quad p_{0}>d,\quad C=A\quad\text{if}
\quad p_{0}=d,
$$
whereas in case (b), 
$C=   NR^{2}$, where $N$ depends only on
$d$ and $\|h\| _{L_{(\infty,d)}}$.

\end{theorem}

Proof.   
Using stopping times we easily reduce the general situation
to the one in which $A,B<\infty$. After that, in light of 
Theorem
\ref{theorem 5.5.1}, we  need only prove that 
in all cases excluding (b)
\begin{equation}
                                                   \label{5.6.5}
B\leq N(d,p_{0} )\Big(A^{1/2}+I_{p_{0}>d} \|h\|_{L_{(q_{0},p_{0})}}
 ^{ p_{0}/(p_{0}-d)  }\Big).
\end{equation}

By Theorem
\ref{theorem 5.5.1}
\begin{equation}
                              \label{5.22,1}
B=E_{\cF_{\tau}}\int_{\tau}^{\gamma}e^{-\phi_{\tau,t}}|dB_{t}|\leq
I(p_{0},q_{0},h)\leq N(d)
 (A+B^{2})^{  d/(2p_{0})}\|h\|_{L_{(q_{0},p_{0})}}.
\end{equation}

Here if $B^{2}\leq A$, estimate \eqref{5.6.5} 
obviously holds.
If $A\leq B^{2}$, then the above inequality yields
$$
B\leq N(d)B^{ d/ p_{0} } \|h\|_{L_{(q_{0},p_{0})}}.
$$
If $p_{0}=d$ and $\|h\|_{L_{(q_{0},p_{0})}}
<\varepsilon=N^{-1}(d)$, this implies $B=0$.
Otherwise we have
 $ 
  B^{(p_{0}- d)/p_{0}}
\leq N(d) \|h\|_{L_{(q_{0},p_{0})}}$  
and we obtain \eqref{5.6.5} again. 
In case $p_{0}=d$ and (b) holds, 
owing to Theorem \ref{theorem 6.3,1},
we need only show that $B\leq NR$. We 
have $q_{0}=\infty$, $\nu_{0}=0$
$$
|b_{t}|\leq (\det a_{t})^{1/d}\sup_{s}h(s
,x_{t}),
$$
and the needed estimate 
of $B$ follows from  Corollary 1.1.13 of \cite{Kr_25}.
The theorem is proved. \qed  

\begin{remark}
                                                   \label{remark 5.7.1}
In Theorem 2.17 of \cite{Kr_19}  estimate 
\eqref{5.6.4} is given if $A_{t}=t$ and   
$c_{t}=\lambda \tr a_{t}$,
where $\lambda>0$ is a number (and $\gamma=\infty$).

\end{remark}

\begin{remark}
                 \label{remark 1.31.2}
The case (a) and the assumption that $A_{t}=t$ in Theorem \ref{theorem 5.5.2} can be actually, eliminated
on account of making a random time
change by using the properly defined inverse function
to $A_{t}$. This will only amount
to replacing $dt$ with $dA_{t}$
everywhere.
 \end{remark}

Here is a corollary of Theorem \ref{theorem 5.5.2} stated in the most common case
when Theorem \ref{theorem 6.3,1} is applicable, so that
$A\leq N (R^{2}+B^{2})$.

\begin{theorem}
         \label{theorem 2.3,1}         
Under Assumption \ref{assumption 9.24.1}
suppose that $(d,q,p)$ are properly tight.
Also   let $A_{t}=t$, take $R>0$ and let $\cO$
be a domain belonging to the cylinder $\bR\times B$ with $B\in \bB_{R}$.
Suppose that $\gamma\leq \tau+\theta_{\tau}\tau _{\cO}$ ($\theta_{\tau}\tau_{\cO}$ is the first exit time of $(\tau+t, x_{\tau+t}) $ 
from  $\cO$).

Then for any Borel $f(t, x)\geq0$   we have
\begin{equation}
                                            \label{2.3.2}
  E_{\cF_{\tau}}\int_{\tau}^{\gamma}e^{-\phi_{\tau,t}} \kappa_{t} 
f( t,x_{t})\,dt \leq N(d,\delta,p_{0},q_{0} )M^{ d/p}
\|f\|_{L_{(q,p)}(\cO)},
\end{equation}
where   $\kappa_{t}$ is from \eqref{9.14,1} and
$$
M=  R + \|h\|_{L_{(q _{0},p _{0})}(\cO)}
^{ p_{0} /(p_{0}-d)   } ,
$$
if $p _{0}>d$,
whereas   $M=N(d,\|h\|_{L_{(\infty,d)}(\cO)})R   $ if $p _{0}=d$ ($q _{0}=\infty$).

\end{theorem}

Proof. On the set $\{\omega:(\tau,x_{\tau})\not\in \cO\}$ we have $\gamma=\tau$ and the left-hand
side of \eqref{2.3.2} is zero. On the complement of this set we replace $b_{t}$
with $\tilde b_{t}:=b_{t}I_{\cO}(t,x_{t})$ and replace $x_{t}$
with
$$
y_{t}=x_{0}+m_{t}+\int_{0}^{t}\tilde b_{s}\,dA_{s}.
$$
Obviously $y_{t}=x_{t}$ on $[\tau,\gamma]$.
Therefore, replacing $x_{t}$ with $y_{t}$
does not
  affect  the left-hand
side of \eqref{2.3.2}. Accordingly, we replace
$h$ with $hI_{\cO}$ preserving 
Assumption \ref{assumption 9.24.1}.
Finally, $f$ in the left-hand
side of \eqref{2.3.2} can obviously be replaced
with $fI_{\cO}$. This shows that without restricting generality we may exclude $\cO$
from the norms above and use that $\gamma
\leq \tau+\theta_{\tau}\tau'_{B}$.

Then the case $p_{0}=d$ is taken care of
by Theorem \ref{theorem 5.5.2}. In case $p_{0}>d$   we only need to prove that
$$
B\leq N(R+\|h\|_{L_{( q_{0},p_{0})}}
^{ p_{0}/(p_{0}-d)}),
$$
which trivially holds if $B\leq R$. If $B\geq R$,
estimate \eqref{5.22,1} shows that
$$
B\leq N B ^{d/p_{0}}\|h\|_{L_{(q_{0},p_{0})}},
\quad B\leq \|h\|_{L_{(q_{0},p_{0})}}^{p_{0}/(p_{0}-d)}.
$$ \qed

The following result will be used much later.
For a Borel set $\Gamma$ in a Euclidean
space set $|\Gamma|$ be its volume and for
a function $f=f(x)$ on $\Gamma$ 
\index{$S$@Miscelenea!$\dashintindex_{\Gamma}f\,$}%
set
$$
\dashint_{\Gamma}f\,dx=\frac{1}{|\Gamma|}
\int_{\Gamma}f\,dx.
$$
If there is a
Banach space $L$ of functions
on $\Gamma$ such that $\|1\|_{L}>0$
and $f\in L$ 
\index{$N$@Norms!$\dashnorm f"|"|$}%
define
$$
\dashnorm f\|_{L}=\|1\|_{L}^{-1}\|f\|_{L}.
$$ 
For instance, 
$$
\dashnorm f\|_{L_{q,p}(C_{\rho})}
=N(d)r^{-d/p-2/q}\|f\|_{L_{q,p}(C_{\rho})}
$$
$$
=\Big(\dashint_{(0,\rho^{2})}
\Big(\dashint_{B_{\rho}}|f(t,x)|^{p}\,dx\Big)^{q/p}\,dt\Big)^{1/q}.
$$

Next lemma gives an analytic condition
for the assumptions of Theorem \ref{theorem 6.3,1}
to be satisfied.

\begin{lemma} 
                      \label{lemma 3.25.1}
Under the assumptions of Theorem \ref{theorem 2.3,1} fix $R\in(0,\infty)$ and suppose that 
$\nu_{0}=0$ and there exists a constant
$\hat h\leq 1$ such that, for any  
$C\in\bC_{R}$, we have
\begin{equation}
                            \label{3,26,1}
\dashnorm h\|_{L_{(q_{0},p_{0})}(C)}\leq
\hat h R^{-1}.
\end{equation} 
Then for any $x\in\bR^{d}$ and $t\geq0$
\begin{equation}
                            \label{3,25,2}
  E_{\cF_{t}}\int_{t}^{t+\theta_{t}\tau_{R}(x)}e^{-\phi_{t,s}} 
|b_{s}|\,ds\leq N(d,\delta,p_{0},q_{0} )
\hat h R.
\end{equation}

\end{lemma}

Proof. By Theorem \ref{theorem 2.3,1} 
the left-hand side of \eqref{3,25,2}
is less than
$$
NM^{d/p_{0}}\|h\|_{L_{(q_{0},p_{0})}(C)}
\leq N\hat h M^{d/p_{0}}R^{d/p_{0}+2/q_{0}-1}
=N\hat h M^{d/p_{0}}R^{1/q_{0} },
$$
where for $p_{0}>d$
$$
M\leq R+N\hat h^{q_{0}} R^{-q_{0}+q_{0}(d/p_{0}+2/q_{0})}\leq N(1+\hat h^{q_{0}})R,
$$
which leads to  \eqref{3,25,2}. For $p_{0}=d$
we have $q_{0}=\infty$ and by Theorem
\ref{theorem 2.3,1} we have 
$$
M\leq N(d,\|h\|_{L_{(\infty,d)}(C)})R
\leq N(d,\hat h)R,
$$
which leads to  \eqref{3,25,2} again. \qed

\begin{remark}
                                  \label{remark 5.7.2}
As in \cite{Na_87} we note that estimate
\eqref{5.6.4} also, obviously, holds 
if 
$$
|b_{t}|\leq \sum_{k=1}^{n}\kappa^{k}_{t}h_{k}(\sft_{t}
,x_{t}),
$$
where $\kappa^{k}_{t}=r_{t}^{1/q_{k}}
(\det a_{t})^{1/p_{k}}c^{\nu_{k}}_{t}$, $p_{k}\in[1,\infty],
q_{k}\in[1,\infty)$, $\nu_{k}=1-d/p_{k}-1/q_{k}\geq0$,
and $h_{k}$ are nonnegative Borel functions. In that
case the constant $C$ depends only on $d,p,q,p_{k},q_{k}$,
  $\|h_{k}\|_{L_{(q_{k},p_{k})}}$, $k=1,..,n$, in a somewhat complicated way.

\end{remark}

\mysection[It\^o's formula]
{It\^o's formula for $W^{1,2}_{(q,p)}$-functions with $d/p+2/q<2$} 

In this section we show the role
of estimates like \eqref{2.3.2} in
establishing It\^o's formula for
functions having generalized rather than continuous derivatives.
Our result is an improvement of previous results in which $u\in W^{1,2}_{d+1} $.
We suppose that we are given a
complete probability space $(\Omega,\cF,P)$ with increasing filtration
of complete $\sigma$-fields $\cF_{t}\subset
\cF$, $t\geq0$, and a $d_{1}$-dimensional
process $w_{t}$, which is
a Wiener process relative to $\{\cF_{t}\}$, where the integer $d_{1}\geq d$.
We also suppose that we are given
progressively measurable
$d\times d_{1}$-matrix valued
process $\sigma_{t}$ and $\bR^{d}$-valued process $b_{t}$ such that
$\sigma_{t}$ is a bounded function
of $(\omega,t)$ and
$$
\int_{0}^{T}|b_{t}|\,dt<\infty
$$
for all $\omega,T$. Under these assumptions the process
$$
x_{t}=\int_{0}^{t}\sigma_{s}\,dw_{s}
+\int_{0}^{t}b_{s}\,ds
$$
is well defined. Our final assumptions
are that

(i) there exists a Borel measurable
function $b$ on $\bR^{d+1}$ such that   
$|b_{t}|\leq b(t,x_{t})$ and
for some
$\rho_{b}\in(0,\infty)$
\index{$S$@Miscelenea!$\bar b_{R}$@$\hat b_{p_{b},\rho_{b}}$}%
 and
$p_{b}\in (1,\infty)$
\begin{equation}
                    \label{6.16.1}
\hat b_{p_{b},\rho_{b}}
:=\sup_{\rho\leq \rho_{b}}
\rho\sup_{C\in\bC_{\rho}}
\dashnorm b\|_{L_{p_{b}}(C)}<\infty;
\end{equation}

(ii)  
for some $p,q\in(1,p_{b})$, satisfying
\begin{equation}
                    \label{6.16.2}
\frac{d}{p}+\frac{2}{q}<2,
\end{equation}
an $R\in(0,\infty)$
and any Borel nonnegative $f$ on 
$\bR^{d+1}$
\begin{equation}
                             \label{12.11.6}
E\int_{0}^{\tau_{R}}f(s,x_{s})\,ds\leq N_{0}\|f\|_{L_{(q,p)}(C_{R})},
\end{equation}
where $N_{0}$ is independent of $f$ and $\tau_{R}$ is the first exit time of $(s,x_{s})$ from $C_{R}$.

\begin{remark}
                \label{remark 6.18,1}
(i) If $b_{t}$ is bounded by a constant, condition \eqref{6.16.1} is obviously
satisfied with any $p_{b}$ and $\rho_{b}$. 
If $b\in L_{p_{b}} $
with $p_{b}\geq d+2$,
condition \eqref{6.16.1} is  
satisfied because
$$
\hat b_{p_{b},\rho_{b}}=N(d)\rho_{b}^{1-(d+2)/p}\sup_{C\in\bC_{\rho_{b}}}\|b\|_{L_{p_{b}}(C)}.
$$
It is an easy exercise to show that
$b(t,x)=1/|x|$ satisfies \eqref{6.16.1} with any $p_{b}\in(1,d)$ and any $\rho_{b}$.

(ii) Theorem \ref{theorem 2.3,1} provides estimates like  
\eqref{12.11.6} under certain conditions, one of which is
$d/p+1/q\leq 1$, which implies
\eqref{6.16.2}.
\end{remark}

We use the following version of Theorem 6.4 of \cite{Kr_26_1}. 
\index{$B$@Sets!$\bR^{d+1}_{t}$}%
Introduce
$$
\bR^{d+1}_{t}=(t,\infty)\times\bR^{d}.
$$

\begin{theorem}
                        \label{theorem 5.25,10}
Let $\hat p,\hat q \in(1,\infty),\hat p_{b}>\max(\hat p,\hat q) $. Then for any   $u\in C^{\infty}_{0}$
$$
\|b |Du|\, \|_{L_{(\hat q,\hat p) }(\bR^{d+1}_{0})}
\leq N\|b\|_{\dot E_{\hat p_{b},1}(\bR^{d+1}_{0})}\|D^{2}u,\partial_{t}u\|_{L_{(\hat q,\hat p)}(\bR^{d+1}_{0})},
$$
where   $N$ depends 
\index{$N$@Norms!$"|"|b"|"|_{\dot E_{p,1}}$}%
only on $d,\hat p,\hat q,\hat p_{b}$  and
$$
\|b\|_{\dot E_{\hat p_{b},1}(\bR^{d+1}_{0})}:=
\sup_{\rho>0}\rho\sup_{\substack{C\in\bC_{\rho}\\C\subset \bR^{d+1}_{0}}}
\dashnorm b\|_{L_{\hat p_{b}}(C)}.  
$$
\end{theorem}

We write that this is a {\em version\/}
of    Theorem 6.4 of \cite{Kr_26_1} because it is proved there
for only one version of $L_{q,p}$
when the interior integration is done 
with respect to $x$. But 
in a few places in \cite{Kr_26_1}
it is emphasized that the order
of integration in $L_{q,p}$ is irrelevant. This is the case of 
Theorem 6.4 of \cite{Kr_26_1} 
as well. Another distinction
is that in \cite{Kr_26_1} we have
$\bR^{d+1}$ in place of $\bR^{d+1}_{0}$. This version is easily shown
to be true if one
\index{$S$@Miscelenea!It\^o's formula}%
 takes $fI_{\bR^{d+1}_{0}}$ in place of $f$ in Theorem 6.4 of \cite{Kr_26_1}.

Here is It\^o's formula.
\begin{theorem}   
                      \label{theorem 12.11.3}
 Under the above assumptions let
  $u\in W^{1,2}_{(q,p)}(C_{R})$.
Then, with probability one for all $t\geq0$,
$$
u(t\wedge\tau_{R},x_{t\wedge\tau _{R}})
=u( 0)+\int_{0}^{t\wedge\tau _{R}}D_{i}u 
 (s,x_{s})\sigma^{ik}_{s}\,dw^{k}_{s}
$$
\begin{equation}
                      \label{6.17,1}
+\int_{0}^{t\wedge\tau _{R}}[
\partial_{t}u(s,x_{s})+(1/2) a^{ij}_{s} D_{ij}u(s,x_{s})
+b^{i}_{s} D_{i}u(s,x_{s})]\,ds
\end{equation}
and the stochastic integral above is a square-integrable
martingale, where $\tau _{R}$ is the first
exit time of $(t,x_{t})$ from $C_{R}$
and $a_{s}=\sigma_{s}\sigma^{*}_{s}$.
\end{theorem}

The proof of the theorem follows the path
suggested in \cite{Kr_69}, \cite{Kr_77} to use the smooth
approximation of $u$ and then use 
\eqref{12.11.6}.  
We prove the theorem after we prove the following.

\begin{lemma}
                \label{lemma 6.16.1}
Under the assumptions of the theorem
we have
\begin{equation}
                             \label{6.16,2}
E\int_{0}^{\tau_{R}}|b|\,|Du|(s,x_{s})\,ds\leq N \|u\|_{W^{1,2}_{(q,p)}(C_{R})},
\end{equation}
where $N$ depends only on $N_{0},R,
\rho_{b}, p, q, p_{b}$, and $\hat b_{p_{b},\rho_{b}}$.
\end{lemma}

Proof. By virtue of \eqref{12.11.6}
it suffices to prove that 
\begin{equation}
                             \label{6.16,3}
I:=\||b|\,|Du|\|_{L_{(q,p)}(C_{R})}\leq N \|u\|_{W^{1,2}_{(q,p)}(C_{R})}.
\end{equation}

There are two very different cases.

{\em Case $p_{b}>d+2$\/}.  By using
H\"older's inequality we obtain that
$$
I\leq\|b\|_{L_{p_{b}}(C_{R})}
\|Du\|_{L_{(r,s)}(C_{R})},
$$
where 
$$
s=\frac{pp_{b}}{p_{b}-p},\quad
r=\frac{qp_{b}}{p_{b}-q}.
$$
After that it only remains to use
the embedding theorems to dominate
the last norm of $Du$ by $\|u\|_{W^{1,2}_{(q,p)}}$.

{\em Case $p_{b}\leq d+2$\/}.
It is not hard to see  that, in light of
\eqref{6.16.1} and the fact that
$p_{b}\leq d+2$,
$$
\sup_{\rho>0}\rho\sup_{C\in\bC_{\rho}}
\dashnorm I_{C_{2R}}b\|_{L_{p_{b}}(C)}
\leq N
<\infty.
$$
Then extend $u$ to $\bR^{d+1}_{0}$
not increasing its $W^{1,2}_{(q,p)}$-norm
by much in such a way that it vanishes
in $\bR^{d+1}_{0}\setminus C_{2R}$. 
Keep the symbol $u$ for such extension.
Then by Theorem \ref{theorem 5.25,10}
$$
I\leq \|I_{C_{2R}}b|Du|\,\|_{L_{(q,p)}
(\bR^{d+1}_{0})}\leq N
\| u\,\|_{W^{1,2}_{(q,p)}
(\bR^{d+1}_{0})}\leq N
\| u\,\|_{W^{1,2}_{(q,p)}
(C_{R})} .
$$
The lemma is proved. \qed   

{\bf Proof of Theorem \ref{theorem 12.11.3}}. The last statement, of course,
follows from \eqref{12.11.6} 
and the fact that
$|Du|^{2}\in L_{(q,p)}(C_{R})$ by embedding theorems in light of \eqref{6.16.2}.
To prove the rest we   approximate  $u$ in $W^{1,2}_{(q,p)}(C_{R})$
by smooth functions $u^{n}, n=1,2,...$ and write for each of them It\^o's
formula like \eqref{6.17,1}. Since $d/p+2/q<2$,
by embedding theorems $u\in C(\bar C_{R})$ and  $|u^{n}-u|\to u$ in $C(\bar C_{R})$ as $n\to \infty$.
This guarantees the convergence of
the terms which are not integrals.

The convergence of integrals is
guaranteed by Lemma \ref{lemma 6.16.1},
estimate \eqref{12.11.6},
and the fact that $Du^{n}\to Du$
in $L_{2q,2p}(C_{R})$ owing to embedding
theorems. The theorem is proved. \qed

\begin{remark}
                        \label{remark 6.17,1}
  Condition \eqref{6.16.2} implies
that $p\vee q>(d+2)/2$. Hence
$p_{b}>(d+2)/2$. If $p_{b}\leq (d+2)/2$, It\^o's formula also holds but
with $u$ from Morrey-Sobolev classes
(see Lemma \ref{lemma 3.16.1}).
\end{remark}

\mysection[Maximum principle]{Application to PDEs. Maximum principle
in $W^{1,2}_{(q,p)}$, $d/p+1/q\leq 1$}

 In this section $\cO\subset\mathbb{R}^{d+1}$ 
is a bounded domain, $R\in(0,\infty)$, and we suppose that
 \smallskip
\begin{equation}
                                                                    \label{6.10.01}
\cO\subset \bR\times B_{R}=\big\{ (t,x):|x|<R\big\} .  
\end{equation} 

Here we apply the previous results
to obtain generalizations of parabolic Aleksandrov's
estimates in the mixed norm setting.
The line of arguments is taken from
\cite{Kr_86},
\cite{Kr_18} and is adapted to
the mixed norm situation as in \cite{Na_87}.

 Assume  that, for any $(t,x)\in \cO$, the following objects
are defined: nonnegative definite $\bS_{0}$-valued
      $a(t,  x)=(a^{ij}(t,  x)) $, $\bR^{d}$-valued  
  $b(t,x)=(b^{i}(t,x)) $
and nonnegative   $c(t,x)$ and $ r(t,x)$. Suppose
that  these functions  are Borel measurable and 
 introduce  
$$
\cL u=r\partial_{t}u+(1/2)a^{ij}D_{ij}u +b^{i}D_{i}u -cu.
 $$

According to the way the operator $\cL$ is defined
we introduce the notion of parabolic boundary 
as follows.

\begin{definition}
                                                   \label{def:3.3.4}   
Given a domain $\cO\subset\mathbb{R}^{d+1}$, the
(right)  {\em parabolic boundary }of $\cO$ is the set of all points
$(t_{0},x_{0})$ belonging to the boundary $\partial \cO$ of $\cO$ (as
a set in $\mathbb{R}^{d+1}$) for each of which one can find a continuous
function $x_{t}\in\bR^{d}$ and a
\index{$S$@Miscelenea!parabolic boundary}%
\index{$S$@Miscelenea!$\partial'$}%
 number $\delta>0$ such that $x_{t_{0}}=x_{0}$
and $(t,  x_{t})\in \cO$ for $t\in[t_{0}-\delta, t_{0})$.
The parabolic boundary of $\cO$ is denoted by $\partial^{\prime}\cO$. 
\end{definition}

\begin{remark}
               \label{remark 6.19,2}
A useful property of $\partial^{\prime}\cO$ is that, if for $\rho>0$
 we denote by $\cO^{ (\rho )}$ the set of
    $(t,x)\in \cO$ for which $\dist((t,x),\partial \cO)>\rho$, then  it turns out that $\dist(\partial'\cO,\partial'\cO^{(\rho)})=
\rho$ if $\cO^{(\rho)}\ne\emptyset$. Indeed, if $(t_{1},x_{1})
\in \partial'\cO^{(\rho)}$ and $(t,x_{t})$, $t\in[t_{1}-\kappa,t_{1}]$,
is a continuous trajectory such that $(t,x_{t})\in \cO^{(\rho)}$
for $t\in[t_{1}-\kappa,t_{1})$
and $ x_{t_{1}}=  x_{1}$, then the point $(t_{1},x_{1})$ is in
$\partial \cO^{(\rho)}$ and its distance to $\partial \cO$ equals $\rho$.
It follows that one can shift the trajectory  $(t,x_{t})$ as a rigid body
by distance $\rho$ so that the shifted trajectory lies in $\cO$
apart from its end point which is on $\partial \cO$ and hence on
$\partial'\cO$. Therefore, $\dist\big(\partial'\cO,(t_{1},x_{1})\big)=\rho$
as claimed.

\end{remark}

\begin{assumption}
                    \label{assumption 9.24.2}   
We have that $(d,q_{0},p_{0})$ are properly tight
 and there exists a Borel $h(t,x)\in L_{(p_{0},q_{0})}$
such that   
\begin{equation}
                                  \label{9.24.3}
|b |\leq r ^{1/q_{0}}(\det a )^{1/p_{0}}c^{\nu_{0}}h 
\quad (\alpha^{0}\equiv1,\alpha\geq0,\nu_{0}=\nu(d,q_{0},p_{0})).
\end{equation}
   
\end{assumption}

In our first result we also use the following. 

\begin{assumption}
\label{assumption 6.2.1}
 
The functions $a,b,c,r$ are bounded and, for a constant $\delta>0$,
we have $a\geq\delta(\delta^{ij})$, 
  $c,r\geq
\delta$.
\end{assumption}

\begin{lemma}
                                              \label{lemma 6.2.2} 
Let 
\begin{equation}
                              \label{9.23.1}
p,q\in(1,\infty),\quad
\nu:=1-\frac{d}{p}-\frac{1}{q}\geq0.
\end{equation}
Suppose that Assumptions \ref{assumption 9.24.2}
and
\ref{assumption 6.2.1}  are satisfied and let   $u\in W^{1,2}_{(q,p)}
(\cO)\cap C(\bar \cO)$.  
Then on $Q$
\begin{equation}
                                                 \label{6.2.1}
u \leq \sup_{\partial'\cO}u_{+}+N(d,p_{0})C
\big\|c^{-\nu}r^{-1/q}(\det a)^{-1/p}(\cL u)_{-}
\big\|_{L_{(q,p)}(\cO )},
\end{equation}
where  
$$
C= \Big(R+ \|h\|_{L_{(q_{0},p_{0})}}
^{  p_{0}/(p_{0}-d)  } 
\Big)^{ d/p}  
 $$
if $p_{0}\ne d$ and $C= N(d,\|h\|_{L_{(\infty,d)}}) R$ if $p_{0}=d$.
\end{lemma}

Proof.  By dividing all coefficients by $r$,
we reduce the general case to the one with $r=1$.
Next step is to observe that Assumption
\ref{assumption 6.2.1} allows us to assume
that $a,b,c$ are smooth.
Then by   using Remark \ref{remark 6.19,2} we see that it suffices
to prove \eqref{6.2.1} with
$\cO^{(\rho)}$ in place of $\cO$.
This shows that we may assume that
$u$ is smooth.

In that case take $(t_{0},x_{0})\in \cO$
and solve the stochastic equation
$$
x_{t}=x_{0}+\int_{0}^{t}\sqrt{a(t_{0}+s,x_{s})}
\,dw_{s}+\int_{0}^{t}b(t_{0}+s,x_{s})
\,ds.
$$
Let $\tau$ be the first exit time of $(t_{0}+s,x_{s})$ from $\cO$. Since $\cO$ is bounded,
$\tau$ is bounded as well. Then  applying It\^o's
formula to
$$
e^{-\phi_{t}}u(t_{0}+s,x_{s}),\quad
\phi_{t}=\int_{0}^{t}c(t_{0}+s,x_{s})
\,ds,
$$
plugging into it $\tau$ in place of $t$,
noticing that $(t_{0}+\tau,x_{\tau})\in
\partial'Q$, we obtain
$$
u(t_{0},x_{0})\leq \sup_{\partial'\cO }u_{+}
+E\int_{0}^{\tau}e^{-\phi_{t}}f(t_{0}+t,x_{t})\,dt,
$$
where $f=-\cL u$. Now introduce
$g=fc^{-\nu}r^{-1/q}(\det a)^{-1/p}$. Then we get
$$
u(t_{0},x_{0})\leq \sup_{\partial'\cO }u_{+}
+E\int_{0}^{\tau}e^{-\phi_{t}}
[c^{ \nu}r^{ 1/q}(\det a)^{ 1/p}
g_{+}](t_{0}+t,x_{t})\,dt,
$$
Now it only remains to apply Theorem
\ref{theorem 5.5.2} along with Theorem 
\ref{theorem 6.3,1}. The lemma is proved. \qed  

In the rest of the section we  
replace Assumption    \ref{assumption 6.2.1}  with the following.

\begin{assumption}
                                            \label{assumption 6.9.2} 
In $\cO$  (a.e.)
\begin{equation}
                                                                   \label{6.10.02}
r+\tr a+c >0.
\end{equation}
\end{assumption}

\begin{theorem}
                                               \label{theorem 6.10.1}

Under Assumptions \ref{assumption 9.24.2}
and \ref{assumption 6.9.2} suppose that
\eqref{9.23.1} is satisfied,
 $u\in W^{1,2}_{(q,p),\loc}
(\cO)\cap C(\bar \cO)$ and  $u\leq0$ on $\partial'\cO$.
Then \eqref{6.2.1} holds with the
same constant  $N$.
 
\end{theorem}

It is important to make precise that in \eqref{6.2.1}  
the norm is taken only over $\{\cL u>0\}\cap \cO $
even if the factor of $(\cL u)_{-}$ is $0^{-1}$
outside this set. Also, if $\nu=0$, we set by definition
$c^{-\nu}\equiv1$ even if $c=0$.

Theorem \ref{theorem 6.10.1} is somewhat
more general than Theorem 4.1 of \cite{Na_87}, where, in particular, if
$d/p+1/q=1$ and $p,q<\infty$, condition
\eqref{6.10.02} is replaced with
$r+\tr a>0$.

To prove the theorem we need three lemmas. 

\begin{lemma}
                                                \label{lem:3.3.1} 
 For $x\ne0$ and   constant $\varepsilon,K\in(0,\infty) $ denote 
$$
\lambda=\frac{x}{|x|},\quad\rho=|x|,
\quad f=\varepsilon\,\tr a+K  (a^{ij}\lambda^{i}\lambda^{j}
+ c  ).
 $$
Take $R<1/\varepsilon$.
Then there exists a function $\psi\in C^{\infty}_{\loc}(\bR^{d})$,
which depends only on $R$, $K $, $\varepsilon$, and $x$, such
that $\psi\ge0$  in $B_{R}$  and  
$$
 a^{ij}D_{ij}\psi +f|D\psi |-c\psi+f+\tr a\leq 0
 $$ 
on  $ \bR\times  B_{R}\setminus\{(t,x):x= 0\}$.
\end{lemma}

This lemma is Lemma 3.1.8 of \cite{Kr_18}
and is proved by  defining  $\psi=\beta-\cosh\big(\alpha|x|\big)$, where
  $\alpha>0$ and $\beta>\cosh(\alpha R)$
 are chosen appropriately.

\begin{lemma}
                                                                       \label{lemma 6.10.5}
Theorem \ref{theorem 6.10.1} holds true if  we additionally assume that
   $u\in W^{1,2}_{(q,p)}(\cO)\cap C(\bar \cO)$ and $h$ is bounded.
 
\end{lemma}

Proof.
We first  
 note that the norm in \eqref{6.2.1}   is not affected
if we multiply 
all coefficients of $\cL$   by the same 
  strictly  positive function.
Therefore,  since by assumption $r+\tr a+ c >0$, we  
 may suppose  that 
\begin{equation}
                                                       \label{eq:3.3.9}
r+\tr a+c=1\quad \mbox{in $\cO$ (a.e.) }
\end{equation}
 without loss of generality. 
\noindent
 Now since $h$ is bounded,
all $r,a,b$, and $c$ are bounded.

For $\delta>0$, we  introduce  
$$
 \cL ^{\delta}=\cL +\delta( \partial_t+\Delta-1)
$$
and denote by  $a(\delta)$, $b(\delta)$, $c(\delta)$, and
$r(\delta)$  the coefficients of $\cL ^{\delta}$. 
 Clearly, $b(\delta)=b$ satisfies \eqref{9.24.3}
with the same $h$ and by Lemma \ref{lemma 6.2.2}
for any $w\in W^{1,2}_{(q,p)}(\cO)\cap C(\bar \cO)$ such
that $w\le0$ on $\partial^{\prime}\cO$,  
\begin{equation}
                                                                  \label{eq:3.3.10}
w\le N\big\| g (\delta)
\big[\cL w+\delta(\partial_{t}w+\Delta w-w)\big]_{-}\big\|
 _{L_{(q,p)}(\cO)}
\end{equation}

\noindent
 in $\bar \cO$, where
$$
 g(\delta)=(c+\delta)^{-\nu}(r+\delta)^{-1/q}(\det(a+\delta
I))^{-1/p}.
$$   

Furthermore, since $h$ is bounded, by the inequality between
geometric an arithmetic means, for any $\kappa>0$ there exists an $N$
such that 
$$
r ^{1/q_{0}}(\det a )^{1/p_{0}}c^{\nu_{0}}h \leq r+\kappa\tr a+ Na^{ij} 
\frac{x^{i}x^{j}}{|x|^{2}}.
$$
Hence owing to \eqref{9.24.3}, there exist constants $\varepsilon_{1}
\in[0,R^{-1})$ and $N\in[0,\infty)$ such that
$$
|b|\leq \varepsilon_{1}\tr a+N \Big[a^{ij}\frac{x^{i}x^{j}}{|x|^{2}}+ c 
+r\Big]=:f_{1}
$$
and $f_{1}>0$. By using Lemma \ref{lem:3.3.1} we find a function $\psi$
such that $\psi>0$ and  
\begin{equation}
                                                                \label{6.10.05}
a^{ij}D_{ij}\psi +f_{1}|D\psi|-(c+r)\psi+f_{1}
+\tr a\leq 0
\end{equation}

\noindent
in $\bR\times B_{R}$.

After that we  take $\gamma>0$ and set 
$$
w=u^{\gamma}:=u+\gamma\psi_{1}
$$  
in (\ref{eq:3.3.10}), where  
$$
\psi_{1}=-e^{-t}(1+ \psi).
$$
 In light   of \eqref{6.10.05}, one easily sees that in $\cO$  
$$
\cL \psi_{1}=e^{-t}\big[r(1+\psi)-\cL (1+\psi)\big]\ge re^{-t}+ce^{-t}+\tr ae^{-t}>0.
$$

 Next, introduce 
$$
g=\lim_{\delta\downarrow0}g(\delta),\quad\Gamma=\big\{
 (t,x)\in \cO\,:\, g(t,x)=\infty\big\} ,
$$
$$
f^{\gamma}=\partial_{t}u^{\gamma}+\Delta u^{\gamma}-u^{\gamma}.
$$

Then  (\ref{eq:3.3.10})  yields  
\begin{align} 
u^{\gamma} & \le N\big\| g 
(\delta) [\cL u^{\gamma}+\delta f^{\gamma} ]_{-}\big\|
_{L_{(q,p)}(\Gamma)}\nonumber 
\\ 
 & +N\big\| g 
(\delta) [\cL u^{\gamma}+\delta f^{\gamma} ]_{-}\big\|
_{L_{(q,p)}(\cO\setminus\Gamma)}.\label{eq:3.3.11}
\end{align}

If  $|\{\cL u<0\}\cap\Gamma|>0$, then the right-hand side of
\eqref{6.2.1}  is infinite and the estimate holds. 
 In the remaining case $\cL u\ge0$  on 
 $\Gamma$ (a.e.)  and  $\cL u^{\gamma}>0$   on $\Gamma$ (a.e.). 
 Observe that 
$g(\delta)\leq\delta^{-1}$  so that  the first  norm on the right 
in (\ref{eq:3.3.11})  is less than 

$$\big\| \, [ \delta^{-1}\cL u^{\gamma}+f^{\gamma} ]_{-}\big\|
_{L_{(q,p)}(\Gamma)}.
$$

 By the Lebesgue dominated convergence theorem
this term tends to zero as $\delta\downarrow 0$ since  
$$
\cL u^{\gamma}>0\quad \mbox{on $\Gamma$ (a.e.)},\quad
  [ \delta^{-1}\cL u^{\gamma}+f^{\gamma} ]_{-}
\le|f^{\gamma}|\in L_{(q,p)}(\cO).
$$

 In the second term  in (\ref{eq:3.3.11}), we  have  $g<\infty$
on $\cO\setminus\Gamma$  implying that  $\det  a>0$ and 
$$
1\ge \delta  g(\delta) =\frac{\delta^{\nu}}
{(c+\delta)^{\nu}}\frac{\delta^{1/q}}{(r+\delta)^{1/q}}\frac{\delta^{d/p}}
{(\det(a+\delta I))^{1/p}}\to0
$$
 on $\cO\setminus\Gamma$ as $\delta\downarrow0$.  
Also observe  that $g(\delta)\le g$ and 
$$
[\cL u^{\gamma}+\delta f^{\gamma}]_{-}\le(\cL u^{\gamma})_{-}
+\delta|f^{\gamma}|\le(\cL u)_{-}+\delta|f^{\gamma}|.
$$

 This and  the monotone convergence theorem and the Lebesgue dominated
convergence theorem  convince us  that, as $\delta\downarrow0$, the limit of the second
 norm  in (\ref{eq:3.3.11})  is less than  
$$
  \lim_{\delta\to0}\big\{ \big\|
g   (\delta) (\cL u )_{-}\big\|
_{L_{(q,p)}(\cO\setminus\Gamma)}+\big\| \delta
g  (\delta)f^{\gamma}\big\|
_{L_{(q,p)}(\cO\setminus\Gamma)}\big\} 
$$
$$
 =  \big\|
g (\cL u)_{-}\big\|
_{L_{(q,p)}(\cO\setminus\Gamma)}\leq\big\|
g (\cL u)_{-}\big\|_{L_{(q,p)}(\cO)}.
$$

 We emphasize that precisely these arguments lead to
our definitions of $c^{0}$. 
Finally,  we  get from (\ref{eq:3.3.11})  that 
$$u^{\gamma} \le N\big\|
 g(\cL u)_{-}\big\| _{L_{(q,p)}(\cO)}.
$$
 After that by
letting $\gamma\to0$, we  arrive at  \eqref{6.2.1}. 
The lemma is proved.   \qed

\begin{lemma}
                                                         \label{lemma 6.10.4}
Lemma \ref{lemma 6.10.5} holds true without the assumption that
$h$ is bounded \(but still   $u\leq0$ on $\partial'\cO$\).
\end{lemma}

Proof.  Introduce
$$
\cL _{n}=I_{h\leq n}\cL +I_{h>n} ( \partial_{t}+\Delta-1\ )
$$
and let $r_{n},a_{n},b_{n},c_{n}$ be the coefficients of $\cL_{n}$.
Obviously, $r_{n},a_{n},b_{n},c_{n}$ satisfy   Assumptions \ref{assumption 9.24.2}
and \ref{assumption 6.9.2}
with   $hI_{h\leq n}$
in place of $h$. They also satisfy \eqref{6.10.02}. Furthermore,  
\begin{align*}
&\big\|c_{n}^{-\theta}r^{-1/q}_{n}(\det a_{n})^{-1/(p)}(\cL _{n}u)_{-}
\big\|_{L_{(q,p)}(\cO )}
\\ 
&\leq \big\|(\partial_{t}u+\Delta u-u)I_{h>n}\big\|_{L_{(q,p)}(\cO )}
\\ 
&\quad+\big\|c^{-\theta}r^{-1/q}(\det a)^{-1/p}(\cL u)_{-}I_{h\leq n}
\big\|_{L_{(q,p)}(\cO )}
\end{align*}

\noindent
and the latter tends to the norm in the right-hand side of \eqref{6.2.1}
owing to the dominated and monotone convergence theorems. Here 
the fact that $\partial_{t}u+\Delta u-u\in L_{(q,p)}(\cO)$ is crucial.
The combination of this argument and Lemma \ref{lemma 6.10.5}
proves the current lemma.     \qed

{\bf Proof of Theorem \ref{theorem 6.10.1}}. 
We may assume that, for a constant $\gamma>0$, we have $u\le-\gamma$ on $\partial^{\prime}\cO$.
Indeed,  generally  we have that   $u-\gamma\le-\gamma$
on $\partial^{\prime}\cO$ and, if estimate 
\eqref{6.2.1}
 holds for $u-\gamma$, then using the fact that 
\begin{equation}
                               \label{9.24.5}
\cL (u-\gamma)=\cL u+\gamma c\ge \cL u,\quad
\big(\cL (u-\gamma)\big)_{-}\le(\cL u)_{-}
\end{equation}
and letting $\gamma\to0$, we obtain \eqref{6.2.1}
 as is.

In light of Remark \ref{remark 6.19,2},
since $u$ is continuous in $\bar \cO$ and
$u\le-\gamma$ on $\partial^{\prime}\cO$, it follows that,
if $\rho$ is small enough,
then $u\le0$ on $\partial'\cO^{(\rho)}$.
 As a result of this argument we have
that $u\in W^{1,2}_{(q,p)}(\cO^{(\rho)})$
and $u\le 0$ on $\partial'\cO^{(\rho)}$. 
After that it only remains to use Lemma
\ref{lemma 6.10.4} with $\cO^{(\rho)}$
in place of $Q$ and let $\rho\downarrow0$.
The theorem is proved. \qed

By substituting $u-\sup_{\partial'\cO}u_{+}$
in place of $u$ in \eqref{6.2.1} and observing
\eqref{9.24.5} we obtain the estimate
in more general form.
\begin{theorem}
                                               \label{theorem 9.24.01}

Under Assumptions \ref{assumption 9.24.2}
and \ref{assumption 6.9.2} suppose that
\eqref{9.23.1} is satisfied and
 $u\in W^{1,2}_{(q,p),\loc}
(\cO)\cap C(\bar \cO)$.
Then 
on $\cO$
\begin{equation}
                                  \label{9.24.6}
u \leq  \sup_{\partial'\cO }u_{+}+N(d,p_{0})C
\big\|c^{-\nu}r^{-1/q}(\det a)^{-1/p}(\cL u)_{-}
\big\|_{L_{(q,p)}(\cO )},
\end{equation}
where  $N$ is the constant in \eqref{6.2.1} and
$$
C= \Big(R+ \|h\|_{L_{(q_{0},p_{0})}}
^{  p_{0}/(p_{0}-d)  } 
\Big)^{ d/p}  
 $$
if $p_{0}>d$ and $C=N(d,\|h\|_{L_{(\infty,d)}})R$ if $p_{0}=d$.
 
\end{theorem}

Just in case, we draw the reader's attention
to the remark after Theorem \ref{theorem 6.10.1}
concerning possible undefined expressions in
\eqref{9.24.6}.

\mysection{An application to
first-order parabolic equations}
                                     \label{section 12.21.1}

Take 
  $p_{0},q_{0},p,q\in[1,\infty)$ such that
$$
 \frac{d}{p_{0}}+\frac{1}{q_{0}}=1,\quad\frac{d}{p}+\frac{1}{q}= 1,\quad p>q_{0}d.
$$

\begin{theorem}

                                  \label{theorem 12.20.2}
Let $Q$ be a bounded domain in $\bR^{d+1}$
 and let $u\in W^{1,2}_{(q,p)}(Q)
\cap C(\bar Q)$ be such that $u\leq 0$
on $\partial'Q$ and
$$
\partial_{t}u+b^{i}D_{i}u \geq  0
$$
in $Q$ (a.e.), where $b\in L_{(q_{0},p_{0})}(Q)$. Then $u \leq0$ in $Q$.
\end{theorem}

This theorem is an immediate consequence of
the following result in which one need only
send $\varepsilon\downarrow 0$ and take into account
that $ q_{0}d/p<1$.

\begin{lemma}
                            \label{lemma 12.20.02}
Let $0\in Q$. Then
under the assumptions of Theorem \ref{theorem 12.20.2}
for any $\varepsilon\in(0,1]$ we have
\begin{equation}
                                      \label{12.20.03}
u(0)\leq N\varepsilon^{-q_{0}d/p }
\|\varepsilon\Delta u\|_{L_{(q,p)}(Q)},
\end{equation}
where $N$ is independent of $\varepsilon$.
\end{lemma}

Proof. Observe that with $I=(\delta^{ij})$,
for any $\varepsilon>0$,
we have
$$
|b|\leq  (\det(\varepsilon I)
 ) )^{1/p_{0}}
\varepsilon^{- d/p_{0}}|b|\qquad (I=(\delta^{ij})).
$$
Hence, by Theorem \ref{theorem 9.24.01} 
$$
u\leq N(1+\varepsilon^{- d/(p_{0}-d)})^{d/p}
\varepsilon^{- d/p}
\big\|\partial_{t}u+\varepsilon  \Delta u
+b^{i}D_{i}u)_{-}\big\|_{L_{(q,p)}(\cO)},
$$
where $N$ is independent of $\varepsilon$.
Here
$$
\partial_{t}u+\varepsilon \Delta u
+b^{i}D_{i}u\geq \varepsilon \Delta u,
\quad 
(\partial_{t}u+\varepsilon \Delta u
+b^{i}D_{i}u)_{-}\leq \varepsilon |\Delta u|.
$$
Therefore,
 $$
u\leq N\varepsilon^{- d^{2}/(pp_{0}-pd)- d/p}
\|\varepsilon\Delta u\|_{L_{(q,p)}(\cO)}.
$$
This yields \eqref{12.20.03} after simple
computations and the lemma is proved. \qed

The result of Theorem \ref{theorem 12.20.2}
is close to be sharp in the following sense.

\begin{example}
                                    \label{example 12.21.1}
Take $\varepsilon\in(0,1)$ and $p_{0},q_{0}\in[1,\infty)$
such that 
$d/p_{0}+1/q_{0}=1+\varepsilon$,
$p_{0}<q_{0}d$ (say $q_{0}\geq 2$),
 and $p_{0}\geq d$. Then it turns out that
there exists $b\in L_{(q_{0},p_{0})}$, $p,q
\in[1,\infty)$ such that $d/p+1/q=1$ and $p>q_{0}d$,
a bounded domain $\cO\subset\bR^{d+1}$ such that
$0\in \cO$, and a function $u\in W^{1,2}_{(q,p)}(\cO)
\cap C(\bar \cO)$   such that $u\leq 0$
on $\partial'\cO$ and
$$
\partial_{t}u+b^{i}D_{i}u \geq  0
$$
in $\cO$ (a.e.), but $u(0)>0$.

To show this set
$$
\alpha=\frac{1-\varepsilon^{2}}{q_{0}},
\quad \beta=\frac{1-\varepsilon^{2}}{p_{0}}d
$$
  and observe that, since $p_{0}<q_{0}d$,
we have $\alpha<\beta$ and $1-\alpha>1-\beta$,
so that there exists $p$ satisfying
\begin{equation}
                                          \label{12.21.03}
\frac{d}{1-\beta}>p>\frac{d}{1-\alpha}.
\end{equation}
Here the left-hand side is strictly bigger than
$q_{0}d$ since 
$$
q_{0}(1-\beta)=q_{0}\Big(1-\frac{d}{p_{0}}+
\frac{\varepsilon^{2}}{p_{0}}d\Big)
=q_{0}\Big( \frac{1}{q_{0}}-\varepsilon+
\frac{\varepsilon^{2}}{p_{0}}d\Big)<1
$$
in light of $p_{0}\geq d$ and $\varepsilon\in(0,1)$.
Therefore, we can choose $p$ satisfying 
\eqref{12.21.03} and such that $p>q_{0}d$
as required. 
 
After that $q$ is also defined and we set
$$
u(t,x)=2-\exp\big(|t|^{1-\alpha}+|x|^{1+\beta}\big),
\quad \cO=\{(t,x):u(t,x)>0\},
$$
$$
  b(t,x)=-\frac{1-\alpha}{1+\beta}\,\frac{1}{|t|^{\alpha}
|x|^{\beta}}\frac{x}{|x|}\,\text{sign}\,\,\! t.
$$
Since  $\alpha q_{0}<1$ and $\beta p_{0}<d$,
 we have $b\in L_{p_{0},q_{0}}(\cO)$.
Also the inequality $\alpha q <1$,
guaranteeing that  $\partial_{t}u
\in L_{q_{0},p_{0}}(\cO)$, is equivalent to 
the right inequality in \eqref{12.21.03},
whereas $p(1-\beta)<d$,
guaranteeing that  $D^{2}u
\in L_{(q_{0},p_{0})}(\cO)$, is equivalent to 
the left inequality in \eqref{12.21.03}.
Hence, $u\in W^{1,2}_{(q_{0},p_{0})}(\cO)$,
$u$ is also continuous, equals zero
on the whole boundary of $\cO$, $u(0)=1$,
and, as is easy to see,
$\partial_{t}u+b^{i}D_{i}u=0$ apart from
the plane $t=0$.

\end{example}

\mysection[Skorokhod's approach]{Passing to the limit in stochastic integrals. Skorokhod's approach}

We will use a few times the
  following
  results  due to A. V. Skorokhod
(see Ch.~1, \S6 and Ch.~2, \S3 in \cite{Sk_61}).  
\begin{lemma}
                           \label{lemma 4.23.1}
  Suppose that $d_1$-dimensional random
 processes $\xi^{(n)}_{t}$  $(t\geq 0, n =  
1,2, . . .)$ are defined
 on some probability spaces equipped with probability measures $P^{n}$. Assume that for each $T> 0$ and
$\varepsilon > 0$
\begin{equation}
                                                  \label{5.10.5}
\lim_{c\to\infty} \sup_{n} \sup_{t\leq T} P^{n} (|\xi^{(n)}_{t}|>c) = 0,
\end{equation}
\begin{equation}
                                                  \label{5.10.6}
\lim_{h\downarrow 0} \sup_{n} \sup_{\substack{t_{1},t_{2}\leq T\\
|t_{1}-t_{2}|\leq h}} P^{n}(|\xi^{(n)}_{t_{1}}-\xi^{(n)}_{t_{2}}|>
\varepsilon)=0.
\end{equation}
Then  one can find a sequence of integers $n'\to\infty$,
 a probability space equipped with a probability measure
 $P$, and
random processes $\tilde \xi_{t},\tilde \xi_{t}^{(n')}$
 defined on this probability space such that all finite-dimensional
distributions of $\tilde \xi_{t}^{(n')}$ coincide with 
the corresponding finite-dimensional
distributions of $\xi_{t}^{(n')}$ and 
$$
P (|\tilde \xi_{t}-\tilde \xi_{t}^{(n')}|>\varepsilon) \to 0   
$$
as $n'\to\infty$ for any $\varepsilon>0$ and $t\geq0$.
\end{lemma}

\begin{lemma}
                                               \label{lemma 4.23.2}
Suppose that on a complete
 probability space $(\Omega,\cF,P)$
we are given random processes $\xi^{(n)}_{t}$,
$w^{(n)}_{t}$, $n=0,1,2,...$. Suppose that
the assumptions of Lemma  \ref{lemma 4.23.1} 
 are satisfied and
\begin{equation}
                                                     \label{5.14.1}
 \xi^{(n)}_{t}\to \xi^{(0)}_{t},\quad w^{(n)}_{t}\to w^{(0)}_{t}
\end{equation}
in probability as $n\to\infty$
for each $t\geq0$. Finally, assume that   $w^{(n)}_{t}$ are $d_{1}$-dimensional
Wiener processes relative to some
increasing families of complete $\sigma$-fields
$\cF_{t}^{n}\subset\cF$,
$t\geq0$, $n=0,1,2,...$, the functions 
$\xi^{(n)}_{t}(\omega)$ are bounded on $[0,\infty)\times\Omega$
uniformly in $n$, and each of them is
progressively measurable relative to $\cF^{n}_{t}$. 

Then the stochastic integrals 
$$
I^{n}_{t}:=\int_{0}^{t}\xi^{(n)}_{s}\,dw^{(n)}_{s}
$$
 are well defined for $t\geq0$, $n=0,1,2,...$ and
   $I^{n}_{t}\to I^{0}_{t}$
in probability as $n\to\infty$
for each $t\geq0$.
\end{lemma}

\begin{remark}
                                                 \label{remark 5.11.1}
As it follows from the proof of Lemma \ref{lemma 4.23.2}
given in \cite{Sk_61}
  we need conditions
\eqref{5.10.5}, \eqref{5.10.6}, and \eqref{5.14.1} to hold
only for $t,t_{1},t_{2}$ restricted to a set
of full measure in order for the assertion of the lemma to be true.
\end{remark}

 In the following Lemma 3.1.4 of \cite{Kr_25} the function $\sigma(t,x)$
is a bounded Borel $d\times d_{1}$-matrix valued
function on $\bR^{d+1}_{0}=(0,\infty)\times\bR^{d}$,
$b(t,x)$ is a Borel $\bR^{d}$-valued function
defined on the same set.

\begin{lemma}
                                             \label{lemma 5.12.1}
Let $\bR^{d+d_{1}}$-valued processes $(x^{(i)}_{t},w^{(i)}_{t})$, $t\geq0$,
$i=1,2$, 
defined on perhaps different complete probability spaces, have the same
finite-dimensional distributions. Define 
$\cF^{i}_{t}$ as the  completion of
  $\sigma(x^{(i)}_{s},\\ w^{(i)}_{s}:s\leq t)$ and assume that $w^{(1)}_{t}$
is a Wiener process with respect to   $\cF^{ 1}_{t}$. Also suppose that $x^{(1)}_{t}$ is continuous and
(a.s.) for all $t\geq 0$
\begin{equation}
                                              \label{5.13.1}
\int_{0}^{t}|b(s,x^{(1)}_{s})|\,ds<\infty,\quad
x^{(1)}_{t}=\int_{0}^{t}\sigma(s,x^{(1)}_{s})\,dw^{(1)}_{s}+
\int_{0}^{t}b(s,x^{(1)}_{s})\,ds.
\end{equation}
Then $x^{(2)}_{t},w^{(2)}_{t}$ have   modifications 
(called again $x^{(2)}_{t},w^{(2)}_{t}$) such that $w^{(2)}_{t}$
is a Wiener process with respect to   $\cF^{2}_{t}$
and (a.s.) for all $t\geq 0$
\begin{equation}
                                              \label{5.13.2}
\int_{0}^{t}|b(s,x^{(2)}_{s})|\,ds<\infty,\quad
x^{(2)}_{t}=\int_{0}^{t}\sigma(s,x^{(2)}_{s})\,dw^{(2)}_{s}+
\int_{0}^{t}b(s,x^{(2)}_{s})\,ds.
\end{equation}
\end{lemma}
 
 The following is similar to Lemma 3.1.5 of \cite{Kr_25}.

 \begin{lemma}
                            \label{lemma 8.13.1}
Let $T\in(0,\infty)$ and let $x^{(n)}_{t}$, $t\in[0,T]$, $n=0,1,2,...$,
be $\bR^{d}$-valued random functions measurable with respect to $(\omega,t)$. Assume that

(i) For any $\varepsilon>0$
$$
\lim_{n\to\infty}\int_{0}^{T}P(|x^{(n)}_{t}-
x^{(0)}_{t}|\geq\varepsilon)\,dt=0;
$$

(ii) For any $R\in(0,\infty)$ there exist $N,p\in[1,\infty)$ such that
for any Borel $f(t,x)\geq0$
vanishing for $|x|\geq R$ and $n\geq1$ we have
\begin{equation}
                                 \label{8.13.2}
E\int_{0}^{T} f(t,x^{(n)}_{t})\,dt
\leq  N\|f\|_{L_{p} }.
\end{equation}

Let a sequence $t^{(n)}\to0$ and let Borel $f^{n}(t,x)$, $n=0,1,2,...$, be given
on $\bR^{d+1}$ such that $f^{n}\to f^{0}$ in  measure
  on any bounded subset
of $\bR^{d+1}$ 
and
\begin{equation}
                                 \label{8.13.3}
\lim_{c\to\infty}\sup_{n\geq0}
E\int_{0}^{T} |f^{n}(t^{(n)}+t,x^{(n)}_{t})|
I_{[c,\infty)}\big(|f^{n}(t^{(n)}+t,x^{(n)}_{t})|\big)\,dt=0.
\end{equation}

Then  
\begin{equation}
                                 \label{8.13.1}
\lim_{n\to \infty}E \int_{0}^{T}
\big|f^{n}(t^{(n)}+t,x^{(n)}_{t})
-f^{0}(t,x^{(0)}_{t})\big|\,dt =0.
\end{equation}

\end{lemma}

Proof. Observe that for $\hat f^{n}(t,x)=f(t^{(n)}+t,x)$ we have
  $\hat f^{n}\to f$ in  measure
 on any bounded subset
of $\bR^{d+1}$. Therefore, we may assume that $t^{(n)}=0$.

By (i) the functions $x^{(n)}_{t}(\omega)$
converge to $x^{(0)}_{t}(\omega)$ with respect to
the product measure.  
Taking $f$ in \eqref{8.13.2} bounded and continuous
and passing to the limit   we conclude
that \eqref{8.13.2} holds with $n=0$. After that
it is extended in a usual way to all Borel
nonnegative $f$ vanishing for $|x|\geq R$.

Upon noting that $f^{n}=(-M)\vee f^{n}\wedge M+
[f^{n}-(-M)\vee f^{n}\wedge M]$ and using \eqref{8.13.3}, we see that it suffices
to prove \eqref{8.13.1} assuming that $|f^{n}|\leq M$
for all $n$. In that case
$$
E\int_{0}^{T}|f^{n}( t,x^{(n)}_{t})|I_{B_{R}^{c}}
(x^{(n)}_{t})\,dt\leq
M\int_{0}^{T}P(|x^{(n)}_{t}|\geq R)\,dt
$$
$$
\to
M\int_{0}^{T}P(|x^{(0)}_{t}|\geq R)\,dt
$$
at all $R$ that are the points of continuity of
the last expression. When $R$ is large this expression
is as small as we like. This shows that we may concentrate on $f^{n}$ that vanish for $|x|\geq R$.

Then note that 
$$
f^{n}( t,x^{(n)}_{t})
-f^{0}( t,x^{(0)}_{t})=\big[f^{n}( t,x^{(n)}_{t})
-f^{0}( t,x^{(n)}_{t})\big]
+\big[
f^{0}(t,x^{(n)}_{t})-f^{0}(t,x^{(0)}_{t})\big]
$$
 and owing to \eqref{8.13.2}
$$
E\int_{0}^{T}\big|f^{n}( t,x^{(n)}_{t})
-f^{0}( t,x^{(n)}_{t})\big|\,dt
\leq N\|I_{B_{R}}(f^{n}-f^{0})\|_{L_{p}(\bR^{d}_{T})}\to 0
$$
as $n\to\infty$. Hence to prove the lemma,
it suffices to prove that
\begin{equation}
                                 \label{8.13.4}
\lim_{n\to \infty}E \int_{0}^{T}
\big|f^{0}(t,x^{(n)}_{t})
-f^{0}(t,x^{(0)}_{t})\big|\,dt =0.
\end{equation}

For any $\varepsilon>0$ one can find a bounded continuous
function $g(t,x)$ vanishing for $|x|\geq R$ and such that
$$
\|f^{0}-g\|_{L_{p}(\bR^{d}_{T})}\leq\varepsilon,
$$
where $\bR^{d}_{T}=(0,T)\times \bR^{d}$.
\index{$B$@Sets!$\bR^{d}_{T}$}%
Then in light of \eqref{8.13.2} it follows that the left hand side of
\eqref{8.13.4} is dominated by
$$
2N\varepsilon+\lim_{n\to \infty}E \int_{0}^{T}
\big|g(t,x^{(n)}_{t})
-g(t,x^{(0)}_{t})\big|\,dt ,
$$
where the latter limit is zero since $x^{(n)}_{t}\to
x^{(0)}_{t}$ with respect to the product measure.
The lemma is proved. \qed

\mysection[An existence theorem
for stochastic equations]{An existence theorem
for stochastic equations with $b\in L_{(q,p)}$, $d/p+1/q\leq 1$}
                                            \label{section 4.23.2}

In this section  we prove a result saying that
in a wide class of cases there exists a 
probability space and a Wiener process
on this space such that a stochastic 
equation having measurable coefficients run by this Wiener process is solvable. In other words,  
 according to the conventional terminology, we are talking
 here about the ``weak'' solutions
of a stochastic equation. The   difference between ``weak'' solutions and
usual (``strong'') solutions consists in
 the fact that the latter can be constructed
on any a priori given probability space on the basis of any given
Wiener process. 
We will use the Skorokhod method, which its creator introduced 
in \cite{Sk_61} to show the solvability of stochastic equations with coefficients
continuous in $x$. Later in \cite{Kr_77}
the author proved the solvability of
uniformly nondegenerate stochastic equations 
with bounded drift.
We follow \cite{Kr_20_2}.

Let $\sigma(t,x)$ be Borel $d\times d$ symmetric
matrix valued, $b(t,x)$ be Borel $\bR^{d}$-valued functions
given on $\bR^{d+1}  =(-\infty,\infty)\times \bR^{d}$. 
We assume that the eigenvalues of $\sigma(t,x)$ are in
$[\delta,\delta^{-1}]$, where $\delta\in(0,1]$ is a fixed number.
Recall that the set of such matrices we denote by $\bS_{\delta}$.

Next, fix   numbers $p,q\in(1,\infty)$, $ \gb\in(0,\infty)$,
let   $b^{(n)}(t,x)$, $n= 1,2,...$, be   $\bR^{d}$-valued 
Borel functions
on $\bR^{d+1} $ and suppose that     
$$
\|b \|_{L_{(q,p)}},
\|b^{(n)}\|_{L_{(q,p)}}\leq  \gb,
\quad n= 1,2,...,\quad \frac{d}{p}+\frac{1}{q}\leq  1,
$$
and $b^{(n)}\to b $ as $n\to\infty$ in  $L_{(q,p)}$.
 Let   $\sigma^{(n)}(t,x)$, $n= 1,2,...$,    
be  Borel functions on $\bR^{d+1}$ with values in $\bS_{\delta}$
  such that
$\sigma^{(n)}\to \sigma $ as $n\to\infty$ ($\bR^{d+1} $-a.e.).

We take $(t^{0},x^{0})\in\bR^{d+1}$ and will be dealing with the equation
\begin{equation}
                                      \label{11.29.2}
x _{t}=x^{(0)}  +\int_{0}^{t}\sigma (t^{(0)}+s,x_{s})\,dw_{s}
+\int_{0}^{t}b (t^{(0)}+s,x_{s}) \,ds,
\end{equation}

\begin{theorem}
              \label{theorem 9.6.4} 
(i) There exists 
a probability space $(\Omega ,\cF ,P )$,
a filtration of $\sigma$-fields $\cF _{t}\subset \cF $, $t\geq0$,
a process $w _{t}$, $t\geq0$, which is a $d$-dimensional Wiener process
relative to $\{\cF _{t}\}$, and an $\cF _{t}$-adapted
process $x_{t}$ such that 
 (a.s.) for all   $t\geq0$ equation \eqref{11.29.2} holds.

(ii) Furthermore,
let $(t^{(n)},x^{(n)})\in\bR^{d+1} $, $n= 1,2,...$, be such that $(t^{(n)},x^{(n)})
\to (t^{(0)},x^{(0)})$  as $n\to\infty$. 
Assume that for each $n=1,2,...$  
 there exists 
a probability space $(\Omega^{n},\cF^{n}, P^{n})$,
a filtration of   $\sigma$-fields $\cF^{n}_{t}\subset \cF^{n}$, $t\geq0$,
a process $w^{(n)}_{t}$, $t\geq0$, which is a $d$-dimensional Wiener process
relative to $\{\cF^{n}_{t}\}$, and an $\cF^{n}_{t}$-adapted
process $x^{(n)}_{t}$ such that (a.s.) for all   $t\geq0$
\begin{equation}
                                                 \label{11.29.1}
x^{(n)}_{t}=x^{(n)} +\int_{0}^{t}\sigma^{(n)}(t^{n}+s,x^{(n)}_{s}) \,dw^{(n)}_{s}
+\int_{0}^{t}b^{(n)}(t^{n}+s,x^{(n)}_{s}) \,ds.
\end{equation}

Then  the finite dimensional distributions of 
a subsequence of
$ x^{(n)}_{\cdot} $ converges weakly to 
the corresponding distributions of one of the solutions of
\eqref{11.29.2}  described in (i). Moreover, if $p\geq q$,  
the set of distributions of $ x^{(n)}_{\cdot} $ on 
$C([0,\infty),\bR^{d })$ is tight.
\end{theorem}

We prove the theorem after some preparations. Define $q_{0}$ from $d/p+1/q_{0}=1$ so that
$q_{0}\leq q$ and introduce  
$$
B(t)=  \|bI_{(0,t)}\|_{L_{(q_{0},p)} }
 ^{ q_{0}}.
$$
In the following lemmas we use an idea from \cite{GM_01}.
\begin{lemma}  
                                 \label{lemma 4.25.1}
Suppose that $p\geq q$ and
let $x_{t}$ be a solution of \eqref{11.29.2}. Then   
for $0\leq s<t<s+1<\infty$ and $n=1,2,...$, we have
\begin{equation}
                                                  \label{4.25.3}
E|x_{t}-x_{s}|^{n}\leq N\big(t-s+B^{2}(t_{0}+t)-B^{2}(t_{0}+s)\big)^{nd/(2p)},
\end{equation}
where $N=N(n,d,\delta,p,\gb )$.
\end{lemma}
 
Proof. We may assume that $t_{0}=0$.
Then observe that for any integer $n=1,2,...$
$$
I_{n+1}:=E\Big(\int_{s}^{t}|b (u,x_{u})| \,du\Big)^{n+1}
$$
$$
=(n+1)!E\int_{s\leq u_{1}\leq...\leq u_{n}}|b (u_{1},x_{u_{1}})|\cdot...\cdot
|b (u_{n},x_{u_{n}})|
$$
\begin{equation}
                                      \label{8.18.1}
\times E_{\cF_{u_{n}}}\Big(\int_{u_{n}}^{t}|b (u,x_{u})| \,du \Big)\,
du_{1}\cdot...\cdot du_{n},
\end{equation}
where by Theorem \ref{theorem 5.5.2}
$$
E_{\cF_{u_{n}}}\int_{u_{n}}^{t}|b (u,x_{u})| \,du 
\leq eE_{\cF_{u_{n}}}\int_{u_{n}}^{t}e^{-(u-u_{n})}|b (u,x_{u})| \,du 
$$
$$
\leq N(d,p)\Big(t-s+\|bI_{(s,t)}\|_{L_{q_{0},p}}^{ 2q_{0}}\Big)^{d/(2p)}
\|b \|_{L_{q,p}}.
$$

Therefore,
$$
I_{n+1}\leq N(n+1)I_{n} 
\Big(t-s+\|bI_{(s,t)}\|_{L_{q_{0},p}}^{ 2q_{0}}\Big)^{d/(2p)}
\|b \|_{L_{q,p}} ,
$$
where $N$ depends only on $d$, $p$, $q$, and $\delta$.  Here
$$
\|bI_{(s,t)}\|_{p,q}^{ 2q}=
\Big(B(t)-B(s)\Big)^{2 }\leq\ 
 B^{2 }(t)-B^{2}(s).
$$
Therefore,
$$
I_{n+1}\leq N(n+1)I_{n}\Big(t-s+B^{2}(t)-B^{2}(s)
\Big) ^{d/(2p)} \|b \|_{L_{q,p}}.
$$

The induction on $n$ yields
$$
I_{n}\leq N^{n}n!\Big(t-s+B^{2}(t)-B^{2}(s)
\Big) ^{nd/(2p)}\|b \|^{n}_{L_{q,p}}.
$$
Also, as is well known,
$$
E\Big|\int_{s}^{t}\sigma (u,x_{u}) \,dw_{u}\Big|^{ n}
\leq N(n,\delta)(t-s)^{n/2}.
$$
It follows that the left-hand side of \eqref{4.25.3}
is less than a constant $N$ times
$$
(t-s)^{n/2}+\Big(t-s+B^{2}(t)-B^{2}(s)
\Big) ^{nd/(2p)},
$$
which less than twice the  factor of $N$
in \eqref{4.25.3} because $p>d$ and $ t-s \leq 1$.
This proves the lemma. \qed

\begin{lemma}
                             \label{lemma 4.25.3}
Under the assumptions in Theorem \ref{theorem 9.6.4} (ii)
the set of distributions of $ x^{(n)}_{\cdot} $ on 
$C([0,\infty),\bR^{d })$ is tight if $p\geq q$.  

\end{lemma}

Proof. Define
$$
B_{n}(t)=\|b^{(n)}I_{(t^{n},t ^{(n)}+t)}\|^{q_{0}}_{L_{  q_{0},p }}
$$
and let $\phi^{n}(s)$ be the inverse function of $ 
\psi^{n}(t):=t ^{(n)}+t+B_{n}^{2}(t ^{(n)}+t)$.
By Lemma \ref{lemma 4.25.1} and Kolomogorov's criterion
the set of distributions of $y^{(n)}_{\cdot}:=
 x^{(n)}_{\phi^{n}(\cdot)} $ on 
$C([0,\infty),\bR^{d })$ is tight.

Observe that, as $n\to\infty$, $\psi^{n}(t)$
converges to $t_{0} +t+B^{2} (t_{0} +t)$ which is continuous and monotone.
By Polya's theorem the convergence is uniform
on any finite time interval, and hence, the
functions $\psi^{n}(t)$ are   
 equi-continuous
on any finite time interval. Now 
define

$$
\Phi(s)=\inf_{n\geq 1}\phi^{n}(s)
$$
and take $S\in(0,\infty)$. By the tightness,
for any  
  $\varepsilon>0$ there is a compact set $K_{\varepsilon}$
in $C([0,S],\bR^{d })$ such that $P^{n}(\{y^{(n)}_{s},s\leq S\}
\in K_{\varepsilon})\geq 1-\varepsilon$ for all $n$.
Due to the uniform continuity of $\psi^{n}$
and of the elements of $K_{\varepsilon}$, the elements of
$$
\hat K_{\varepsilon}
:=\{\{f(\psi^{n}(t)),t\leq \Phi(S)\}: \{f(s),s\leq S\}
\in K_{\varepsilon},n=1,2,...\}
$$
are uniformly continuous and, of course, uniformly bounded,
so that $\hat K_{\varepsilon}$ is a compact set 
in $C([0,\Phi(S)],\bR^{d})$
and
$$
P(\{y^{(n)}_{\psi^{n}(t)},t\leq \Phi(S)\}\in\hat K_{\varepsilon})
\geq 1-\varepsilon.
$$
It only remains to observe that
 $y^{(n)}_{\psi^{n}(t)}=x^{(n)}_{t}$, $S$ is arbitrary,
and $\Phi(S)\to\infty$ as $S\to\infty$. The lemma is proved. \qed

{\bf Proof of Theorem  \ref{theorem 9.6.4}}. 
Due to the possibility to use
 mollifiers we see that assertion (ii)
implies (i).
In the proof of (ii), thanks to   
 Lemma \ref{lemma 4.25.3}, we need only prove
the assertion concerning
the convergence of finite dimensional distributions. 

Having in mind Lemma \ref{lemma 4.23.1} define for $M>0$
$$
\xi^{(n)}_{t}=\int_{0}^{t}b^{(n)}(t^{(n)}+s,x^{(n)}_{s})\,ds,
$$
$$
\xi^{(n)M}_{t}=\int_{0}^{t}b^{(n)}(t^{(n)}+s,x^{(n)}_{s})
I_{|b^{(n)}( t^{(n)}+s, x^{(n)}_{s})|\leq M}\,ds.
$$

Since the derivative of $\xi^{(n)M}_{t}$ is bounded,
both conditions \eqref{5.10.5} and \eqref{5.10.6}
are satisfied for $\xi^{(n)M}_{t}$. Furthermore,
for any $T\in(0,\infty)$ by Aleksandrov's estimates
\begin{equation}
                                  \label{8.14.3}
E^{n}\int_{0}^{T}|b^{(n)}(t^{(n)}+s,x^{(n)}_{s})|
I_{|b^{(n)}(t^{(n)}+s,x^{(n)}_{s})|\geq M}\,ds 
\leq  
N\|b^{(n)}I_{|b^{(n)}|\geq M}\|_{L_{(q,p)}},
\end{equation}
where $N$ is independent of $n$. Since $b^{n}\to b$
in the $\|\cdot\|_{(q,p)}$-norm, the latter quantity can be made
as small as we like on account of choosing $M$
 large enough. Therefore, Lemma \ref{lemma 4.23.1} is applicable
to $\xi^{n}_{t}$. It is, obviously, also applicable to
$$
\eta^{(n)}_{t}=x^{(n)}+\int_{0}^{t}\sigma^{n}(t^{(n)}+s,x^{(n)}_{s}) \,dw^{(n)}_{s}.
$$
Hence, there is a subsequence, which by common abuse of notation
we identify with the original one, a probability space and
random $\bR^{2d}$-valued
processes $(\tilde x^{(n)}_{t},\tilde w^{(n)}_{t})$,
$(\tilde x^{(0)}_{t},\tilde w^{(0)}_{t})$  
defined on this probability space such that all finite-dimensional
distributions of $(\tilde x^{(n)}_{t},\tilde w^{(n)}_{t})$ coincide with 
the corresponding finite-dimensional
distributions of $(x_{t}^{(n)},w^{(n)}_{t})$ and 
\begin{equation}
                                                       \label{5.13.6}
P (|(\tilde x^{(n)}_{t},\tilde w^{(n)}_{t})
-(\tilde x^{(0)}_{t},\tilde w^{(0)}_{t})|\geq \varepsilon) \to 0   
\end{equation}
as $n \to\infty$ for any $\varepsilon>0$ and $t\geq0$.
 
Furthermore (as a result of   
\eqref{5.10.5}), for any $T\in(0,\infty)$ as $R\to
\infty$
\begin{equation}
                                   \label{8.14.1}
P (|\tilde x^{(n)}_{t}|>R)\to0
\end{equation}
uniformly with respect to $t\leq T$ and $n\geq 1$ and, 
as \eqref{5.13.6}   implies, with respect to
$n\geq 0$.

For $n\geq0$ introduce $\tilde \cF^{n}_{t}$ as the completion
of $\sigma(\tilde x^{(n)}_{s}, \tilde w^{(n)}_{s},s\leq t)$. It is easy to see,
using Kolmogorov's continuity criterion, that
$\tilde w^{(0)}_{t}$ admits a continuous modification
$\hat w^{(0)}_{t}$ such that
$\{\hat w^{(0)}_{t},\tilde \cF^{0}_{t}\}$  is a Wiener
process.

By Lemma \ref{lemma 5.12.1}, for each $n\geq1$, the process
$(\tilde x^{(n)}_{t},\tilde w_{t}^{(n)})$
admits a continuous
modification denoted by   $(\hat x^{(n)}_{t},\hat w_{t}^{(n)})$
such  that 
 $(\hat w_{t}^{(n)},
\tilde \cF_{t}^{n})$ is a Wiener process and  (a.s) for all $t\geq0$
\begin{equation}
                                                     \label{5.12.1}
\hat x^{(n)}_{t}=x^{(n)}+\int_{0}^{t}\sigma^{(n)}(t_{n}+s,
\hat x^{(n)}_{s})\,d\hat w^{(n)}_{s}+\int_{0}^{t}b^{(n)}(t_{n}+s,
\hat x^{(n)}_{s})\,ds.
\end{equation}

In light of \eqref{5.13.6}   we have
\begin{equation}
                                                       \label{5.17.1}
P (|(\hat x^{(n)}_{t},\hat w^{(n)}_{t})
-(\tilde x^{(0)}_{t},\tilde w^{(0)}_{t})|\geq \varepsilon) \to 0   
\end{equation}
as $n \to\infty$ for any $\varepsilon>0$ and $t\geq0$.

Now the fact that $\tilde x^{(0)}_{t}$ may be not measurable
in $t$ causes some problems. However, 
set $\phi(x)=x/(1+|x|)$  and observe that, owing to \eqref{5.17.1},
 $\phi(\hat x^{(n)}_{t})$ form a Cauchy sequence
in $L_{1}(\Omega\times[0,T])$ and, hence, converges in that space
to $\phi(\hat x^{(0)}_{t})$, where $\hat x^{(0)}_{t}$ is measurable
with respect to $(\omega,t)$. By Fubini's theorem
there is a set $\cS\subset [0,\infty)$ of full measure
such that, for any $t\in\cS$, $\hat x^{(0)}_{t}=\tilde x^{(0)}_{t}$
(a.s.).   We set
$\hat x^{(0)}_{t}=0$ for $t\not\in \cS$ and 
observe that $\hat x^{(0)}_{t}$
is $\tilde \cF^{0}_{t}$-adapted.

Also note that   \eqref{5.17.1}  remains 
valid if we replace $(\tilde x^{(0)}_{t},\tilde w^{(0)}_{t})$
by $(\hat x^{(0)}_{t},\hat w^{(0)}_{t})$ and 
restrict the ranges of $t,s$ to $t,s\in\cS$.
This is done to accommodate Remark \ref{remark 5.11.1}.
 Then by Lemma
\ref{lemma 4.23.2} for any  $t\geq0$ and 
bounded continuous
$d\times d$ symmetric matrix-valued
 $\alpha(t,x)$ we have 
\begin{equation}
                                                     \label{5.13.7}
 \int_{0}^{t}\alpha(t^{(n)}+s,\hat x_{s}^{(n)})\,d\hat w^{(n)}_{s}
\to \int_{0}^{t}\alpha(t^{(0)}+s,\hat x_{s}^{(0)})\,d\hat w^{(0)}_{s}
\end{equation}
as $n\to\infty$ in probability.
We want to use this  to pass to the limit in the stochastic term
in \eqref{5.12.1}.
But first observe that
by Theorem \ref{theorem 5.5.2}  
for any $T\in(0,\infty)$, Borel $f(t,x)\geq0$, and $n\geq1$  
\begin{equation}
                                                  \label{4.26.1}
E\int_{0}^{T}
f(t,\hat x^{n}_{t})\,dt\leq N\|fI_{(0,T)}\|_{L_{d+1}},
\end{equation}
where $N$ is independent of $f$ and $n$. The convergence
in probability implies that \eqref{4.26.1} holds for
$n=0$ as well with the same constant $N$, first for 
nonnegative $f\in C^{\infty}_{0}
(\bR^{d+1})$ and then, due to general measure-theoretic arguments,
for any Borel nonnegative $f$.

Then take an $\alpha$ as above with values in $\bS_{\delta}$ and write
$$
I_{n}(t):= \int_{0}^{t}\sigma^{(n)}(t^{(n)}+s,\hat x_{s}^{(n)})\,d\hat w^{(n)}_{s}
-\int_{0}^{t}\sigma^{(0)}(t^{(0)}+s,\hat x_{s}^{(0)})\,d\hat w^{(0)}_{s}
$$
$$
=J_{n}(t)+I_{n}(B_{R},t)
+I_{n}(B_{R}^{c},t)+I  (B_{R},t)+I (B_{R}^{c},t),
$$
where
$$
J_{n}(t):=\int_{0}^{t}\alpha(t^{(n)}+s,\hat x_{s}^{(n)})\,d\hat w^{(n)}_{s}
- \int_{0}^{t}\alpha(t^{(0)}+s,\hat x_{s}^{(0)})\,d\hat w^{(0)}_{s},
$$
$$
I _{n}(\Gamma,t):=\int_{0}^{t}\beta^{(n)}_{\Gamma}(
t^{(n)}+s,\hat x_{s}^{(n)}) \,d\hat w^{(n)}_{s}, \quad
\beta^{(n)}_{\Gamma} =(
 \sigma^{(n)} -\alpha)I_{\Gamma}, 
$$
$$
I (\Gamma,t):=\int_{0}^{t}
\beta^{(0)}_{\Gamma}(t^{(0)}+s,\hat x_{s}^{(0)})  \,d\hat w^{(0)}_{s}.
$$
Our goal is to show that
\begin{equation}                                                 
                                \label{8.14.2}
I_{n}(t)\to 0,
\end{equation}
as $n\to\infty$ in probability and we already know
that this holds for $J_{n}(t)$. Therefore,
it suffices to show that on account of choosing
$R$ and
$\alpha$, for any $\varepsilon>0$, we can make the probabilities that
$|I_{n}(B_{R},t)|\geq\varepsilon$, 
$|I _{n}(B_{R}^{c},t)|\geq\varepsilon$,
$|I  (B_{R},t)|\geq\varepsilon$,
$|I (B_{R}^{c},t)|\geq\varepsilon$, as small as we like for all large $n$.

Observe that
$$
P(|I _{n}(B_{R}^{c},t)|\geq\varepsilon)\leq
\varepsilon^{-2}NE\int_{0}^{t}I_{B_{R}^{c}}(\hat x^{(n)}_{s})\,ds,
$$
where $N$ depends only on $d$ and $\delta$.
Here the right-hand side is uniformly in $n$
small if $R$ is large on  account of \eqref{8.14.1}
which obviously holds true if we replace $\tilde x^{(n)}_{t}$ with $\hat x^{(n)}_{t}$.

Then by \eqref{4.26.1}, for $G=[0,t]\times B_{R}$,
$$
P(|I _{n}(B_{R} ,t)|\geq\varepsilon)\leq
\varepsilon^{-2} E\int_{0}^{t}\|\beta^{(n)}_{B_{R}}(
t^{(n)}+s,\hat x_{s}^{(n)})\|^{2}\,ds
$$
$$
\leq
\varepsilon^{-2} NE\int_{0}^{t}\|\beta^{(n)}_{B_{R}}(t^{(n)}+s,\hat x_{s}^{(n)})\| \,ds 
\leq \varepsilon^{-2}N\|(\sigma^{(n)}-\alpha)(t^{(n)}+\cdot,\cdot)\|_{L_{d+1}(G)} .
$$
The last term tends to $\varepsilon^{-2}\|\sigma^{(0)}-\alpha\|_{L_{d+1}(G)}$ as $n\to\infty$
and this shows how to choose $\alpha$.
One deals with $|I  (B_{R},t)|\geq\varepsilon$,
$|I (B_{R}^{c},t)|\geq\varepsilon$ similarly
and arrives at \eqref{8.14.2}.

Finally, in light of \eqref{8.14.3} a direct application of Lemma
\ref{lemma 8.13.1} proves that in probability
$$
\int_{0}^{t}b^{(n)}(t_{n}+s,
\hat x^{(n)}_{s})\,ds\to \int_{0}^{t}b^{(0)}(t_{0}+s,
\hat x^{(0)}_{s})\,ds.
$$

This and \eqref{8.14.2} allow us to pass to the limit
in \eqref{5.12.1} when $t\in \cS$ (and 
$\hat x^{(n)}_{t}\to \hat x^{(0)}_{t}$ in probability) and shows that
\eqref{5.12.1} holds true for $n=0$ and
$t\in \cS$. In turn this implies that
$\hat x^{(n)}_{t}$ is extendible from the 
set of full measure $\cS$ to all $t$ (as the right-hand side of \eqref{5.12.1}) as a continuous
function satisfying \eqref{5.12.1} with $n=0$
for all $t$ at once (a.s.).
The theorem is proved. \qed
 
\mysection[Examples of nonexistence
and nonuniqueness]{Examples of nonexistence
and nonuniqueness}
                      \label{section 4.23.1}

The following example is taken from
\cite{Kr_20_2}.

 \begin{example}
                       \label{example 3.22.1}

Suppose that   numbers $\alpha$ and $\beta$ satisfy 
\begin{equation}
                                                  \label{4.18.1}
0<\alpha\leq \beta <1,\quad \alpha+\beta=1
\end{equation}
 and set
$$
b(t,x)=-\frac{1}{t^{\alpha}|x |^{\beta}}\frac{x }{|x |}
I_{0<|x|,t\leq 1}.
$$
Observe that, if $d/p+1/q=1+\varepsilon$, $\varepsilon>0$, 
one can take $\beta=d/(p+p\varepsilon)$, $\alpha=1/(q+q\varepsilon)$
and then 
$$
\int_{\bR}\Big(\int_{\bR^{d}}|b(t,x)|^{p}\,dx\Big)^{q/p}dt<\infty,
\quad 
\int_{\bR^{d}}\Big(\int_{\bR}|b(t,x)|^{q}\,dt\Big)^{p/q}dx<\infty.
$$
Also note
 that if  $p\leq qd$ (say $p=q$), condition \eqref{4.18.1}  is satisfied.

However, it turns out that no matter which $\alpha,\beta$
we take satisfying \eqref{4.18.1}    there are no solutions of the equation $dx_{t}=dw_{t}
+b(t,x_{t})\,dt$ starting at zero, where $w_{t}$ is a $d$-dimensional
Wiener process.

To prove this assume the contrary. Namely, assume the there is a stopping
time $\tau$ such that $P(\tau>0)>0$ and for $t\leq \tau$
there is $x_{t}$ such that
$$
x_{t}=w_{t}+\int_{0}^{t}b(s,x_{s})\,ds.
$$
We may assume that $\tau\leq 1$ and before $\tau$ the process 
is in $B_{1}$.
Then for $t\leq \tau$ 
\begin{equation}
                                                \label{4.18.4}
dx_{t}=-\frac{1}{t^{\alpha}|x_{t}|^{\beta}}\frac{x_{t}}{|x_{t}|}
I_{ x_{t}\ne 0 }\,dt+dw_{t},
\end{equation}
$$
d|x_{t}|^{2}=-2\frac{|x_{t}|}{t^{\alpha}|x_{t}|^{\beta}}\,dt
+d\,dt+2x_{t}\,dw_{t}.
$$
 
We will be  interested in $|x_{t}|^{1+\beta}=\xi_{t}^{(1+\beta)/2}$,
where $\xi_{t}=|x_{t}|^{2}$. 
By It\^o's formula for any $\varepsilon>0$ we have
$$
d(\xi_{t}+\varepsilon)^{(1+\beta)/2}=
\frac{1+\beta}{2}(\xi_{t}+\varepsilon)^{( \beta-1)/2}\,d\xi_{t}
+\frac{\beta^{2}-1}{8}(\xi_{t}+\varepsilon)^{( \beta-3)/2} 4|x_{t}|^{2}\,dt
$$
\begin{equation}
                                                           \label{3.18.1}
=I_{t}(\varepsilon)\,dt+J_{t}(\varepsilon)\,dt 
+(1+\beta)(\xi_{t}+\varepsilon)^{( \beta-1)/2}x_{t}\,dw_{t},
\end{equation}
where
$$
I_{t}(\varepsilon)=-(1+\beta)(\xi_{t}+\varepsilon)^{( \beta-1)/2}
\frac{|x_{t}|^{\alpha}}{t^{\alpha}} ,
$$

$$
J_{t}(\varepsilon) =
\frac{1+\beta}{2} \Big[ d+(\beta-1)
(\xi_{t}+\varepsilon)^{-1}  |x_{t}|^{2} \Big](\xi_{t}+\varepsilon)^{
(\beta-1)/2}.
$$
Since $(\xi_{t}+\varepsilon)^{-\alpha/2}
 |x_{t}|^{\alpha}\uparrow I_{x_{t}\ne0}$ as $\varepsilon\downarrow0$,
by the dominated convergence theorem
$$
\int_{0}^{t}I_{s}(\varepsilon)\,ds\to
-(1+\beta)\int_{0}^{t}I_{x_{s}\ne0}
\frac{1}{s^{\alpha}}\,ds,
$$
which is finite. 

Furthermore,  since $|x_{s}|^{ \beta-1 }x_{s}$ is bounded
on each trajectory, by the dominated
convergence theorem
$$
\int_{0}^{t}|(\xi_{s}+\varepsilon)^{(\beta-1)/2}x_{s}
-|x_{s}|^{\beta-1 }x_{s}|^{2}\,ds\to0,
$$
and we conclude from \eqref{3.18.1} that for $t\leq\tau$
$$
|x_{t}|^{1+\beta}=-(1+\beta)\int_{0}^{t}I_{x_{s}\ne0}
\frac{1}{s^{\alpha}}\,ds
$$
\begin{equation}
                                                       \label{4.18.2}
+\lim_{\varepsilon\downarrow0}\int_{0}^{t}J_{s}(\varepsilon)\,ds
+(1+\beta)\int_{0}^{t}|x_{s}|^{\beta-1}x_{s}I_{x_{s}\ne0}\,dw_{s}
\end{equation}
and the above limit exists and is finite. Since
$2J_{s}(\varepsilon)\geq (\xi_{s}+\varepsilon)^{(\beta-1)/2}$,
it follows that
$$
\int_{0}^{t}|x_{s}|^{\beta-1 }\,ds=
\lim_{\varepsilon\downarrow0}\int_{0}^{t}
(\xi_{s}+\varepsilon)^{(\beta-1)/2}\,ds
$$
and the left-hand side is finite. 
In particular, 
\begin{equation}
                                                    \label{4.18.5}
\int_{0}^{\tau}I_{x_{s}=0}\,ds=0.
\end{equation}

Now by the dominated convergence theorem
\eqref{4.18.2} implies that  
$$
|x_{t}|^{1+\beta}=-(1+\beta)\int_{0}^{t} 
\frac{1}{s^{\alpha}}\,ds
$$
$$
+(1/2)(1+\beta)\int_{0}^{t}(d+\beta-1)|x_{s}|^{\beta-1 }\,ds
+(1+\beta)\int_{0}^{t}|x_{s}|^{\beta-1}x_{s} \,dw_{s}.
$$

Next, use $\alpha\leq \beta$ and H\"older's inequality to
conclude that
$$
\int_{0}^{t}|x_{s}|^{-\alpha}\,ds=
\int_{0}^{t}\Big(\frac{1}{s^{\alpha}|x_{s}|^{\beta}}
\Big)^{\alpha/\beta}s^{\alpha^{2}/\beta}\,ds
$$
 $$
\leq\Big(\int_{0}^{t}\frac{1}{s^{\alpha}
|x_{s}|^{\beta}}\,ds\Big)^{\alpha/\beta}
\Big(\int_{0}^{t}s^{\alpha^{2}/(\beta-\alpha)}
\,ds\Big)^{(\beta-\alpha)/\beta}.
$$
Since, $\alpha^{2}/(\beta-\alpha)+1=(\alpha^{2}+1-2\alpha)/(\beta-\alpha)
=\beta^{2}/(\beta-\alpha)$
$$
\int_{0}^{t}|x_{s}|^{-\alpha}\,ds\leq N 
\Big(\int_{0}^{t}\frac{1}{s^{\alpha}
|x_{s}|^{\beta}}\,ds\Big)^{\alpha/\beta}t^{\beta},
$$
where $N =N (\alpha,\beta)$
(which is trivial if $\alpha=\beta$).
Thus,
$$
|x_{t}|^{1+\beta}+ct^{\beta}\leq N_{1}
\Big(\int_{0}^{t}\frac{1}{s^{\alpha}
|x_{s}|^{\beta}}\,ds\Big)^{\alpha/\beta}t^{\beta}+
(1+\beta)\int_{0}^{t}|x_{s}|^{\beta-1}x_{s}\,dw_{s},
$$
where $c>0$ is a constant.
For   equation \eqref{4.18.4} to make sense we should have
\begin{equation}
                                            \label{4.18.3}
\int_{0}^{\tau}\frac{1}{s^{\alpha}
|x_{s}|^{\beta}}\,ds<\infty 
\end{equation}
  (a.s.). Therefore 
$$
\gamma:=\tau\wedge\inf\{t\geq0:N_{1}
\Big(\int_{0}^{t}\frac{1}{s ^{\alpha}
|x_{s}|^{\beta}}\,ds\Big)^{\alpha/\beta}\geq c/2\},
$$
is a stopping time such that $P(\gamma>0)=P(\tau>0)$.
It follows that  for any $t>0$
$$
\int_{0}^{t}I_{s<\gamma}|x_{s}|^{\beta-1}x_{s}\,dw_{s}
\geq  0,
$$
which is only possible if $I_{s<\gamma}|x_{s}|^{\beta-1}x_{s}=0$
for almost all $(\omega,s)$. Then $x_{s}=0$ for $s<\gamma$
and \eqref{4.18.5} is only possible if $P(\tau=0)=1$.
\end{example}
 
\begin{remark}
                         \label{remark 12.22.1}
One may ask if the size of $b$ plays a role
in the above argument. More precisely,
take $\varepsilon>0$ and consider the equation
$dx_{t}=dw_{t}+\varepsilon b(t,x_{t})\,dt$
with zero initial condition. It turns out that
this equation does not have solutions either.

To see this, set $c=\varepsilon^{1/\alpha}$ 
and denote $y_{t}=c^{-1}x_{c^{2}t}$, $B_{t}=c^{-1}
w_{c^{2}t}$. Then the equations becomes
$dy_{t}=dB_{t}+b(t,y_{t})\,dt$, and since $B_{t}$
is a Wiener process, it does no have solutions.

\end{remark}

Next, we are going to present an example
where different solutions have different
finite-dimensional distributions
(no weak uniqueness). This
 example  was
brought to the author's attention by
M.~Gerensc\'er (also see  \cite{GG_25}).

\begin{example}
                       \label{example 3.22.2}
Take $1<q<2$, set $b^{1}(t,x)=t^{-1/q}I_{0<t\leq1,|x^{1}|\leq1}\text{sign}\, x^{1}$, $b^{i}=0$, $i\geq 2$,
and consider the equation
\begin{equation}
                            \label{3.21.1}
x_{t}=x+\int_{0}^{t}b(s,x_{s})\,ds+w_{t},
\end{equation}
where $w_{t}$ is a $d$-dimensional
Wiener process. Obviously, $b\in L_{(q*,p*)}$
for $p*$ large enough satisfying
$d/p^*+1/q^*<1$, so that \eqref{3.21.1}
is solvable on appropriate probability spaces.

Next, set $x'=(x^{2},...,x^{d})$ and find a constant $t_{0}$ such that for
$$
A=\{\sup_{t\leq t_{0}}|w'_{t}|\leq 1,
\inf_{t\leq t_{0}}(3t^{1-1/q}+w^{1}_{t})\geq 0\}
$$
($1-1/q<1/2$) we have
$$
P(A) \geq 3/4.
$$

If $x^{1}>0$ and $\omega\in A$, then $x^{1}_{t}\geq0$ for $t\leq t_{0}$. Indeed, otherwise
there is $s\leq t_{0}$ such that $x^{1}_{t}$
becomes $0$ for the first time and then
$$
0=x^{1}+4s^{1-1/q}+w_{s},\quad4s^{1-1/q}+w_{s}=-x^{1},\quad \inf_{t\leq t_{0}}(3t^{1-1/q}+w^{1}_{t})\leq -x^{1}.
$$
 
Now, when $x=x_{n}$ is such that $x^{1}_{n}
\downarrow0$ and $p^*>q^*$ the distributions of solutions of \eqref{3.21.1} on $C[0,\infty),\bR^{d})$
converge in the weak topology to
the distribution of a solution of \eqref{3.21.1} starting at the origin.
Since the set $\{x_{\cdot}:x_{t}\geq0,t\leq t_{0}\}$ is closed in  $C[0,\infty),\bR^{d})$, the limiting probability of this set
is greater that $3/4$. If we approximate $0$
by $x^{1}_{n}$ from below we will have another
solution for which this probability is less
that $1/4$. Hence, nonuniqueness of finite-dimensional distributions.

Observe for the future that, for any
 $q_{b}\in (1,q)$,
$$
r\Big(\dashint_{(0,r^{2})}t^{-q_{b}/q}\,dt
\Big)^{1/q_{b}} =Nr^{1-1/q_{b}-1/q}
$$
which tends to infinity as $r\downarrow 0$
since $p,q_{b}<2$. It would go to zero
in case $q_{b}=2$. However, the above example does not work in that case.
\end{example}

\mysection[Markovian families of random processes]{On Markovian families of random processes}
                             \label{section 8.15.1}

This section is based 
on \cite{Kr_73_1}.
Let $E$ be  a Polish space representable as  the countable union
 of  compact sets, $\bB(E)$
 \index{$B$@Sets!$\bB(E)$}%
be the space of bounded Borel functions on $E$ with the
norm 
$$
\|f\| = \sup_{x \in E}|f(x)|.
$$ 
Define  $\bar C_{0}$ as the closure in the norm $\|\cdot\|$ of the family $C_{0}$ of continuous
functions on $E$ with compact support, $H$ the closure in the norm 
$$
\sup_{ t\geq 0}\|f(t,\cdot)\|
$$
of the family of continuous functions $f(t, x)$ defined on $[0,\infty)\times E$ and having compact
support in $[0,\infty)\times E$. 

Let $\Omega = \{\omega\}$ be a set; let $\frM^{0}$ and $\frM_{t}$, $t\geq0$, be $\sigma$-fields of the subsets
of $\Omega$ such that $\frM^{0}\supset \frM_{t}\supset \frM_{s}$ for all $t \geq s$  . Assume that for each 
$\omega\in \Omega$ a function
$x_{t}(\omega)$ is defined which is continuous from the right with respect to $t$ on $[0, \infty)$ and takes on values
in E. Assume that the collection of functions $x_{t}(\omega)$  is such that for each $\omega\in \Omega$ and
$s \geq0$   an $\omega' \in \Omega$ can be found such that $x_{t}(\omega')= x_{t+s}(\omega)$ for all 
$t\geq 0$. As usual, $\cN_{t}$ denotes the $\sigma$-field of subsets of $\Omega$ generated by the sets of the form $\{\omega: x_{s}(\omega) \in \Gamma\}$
for $s\leq t$  and Borel sets $\Gamma\subset E$;
$\cN_{\infty}=\sigma(x_{s},s<\infty)$.   Assume that $\cN_{t}\subset\frM_{t}$, for all $t\geq0$.

Let $P^{n}$, $n=0,1,...$, be a sequence of probability measures on $(\Omega,\frM^{0})$. We write
$P^{n}\to P^{0}$ provided that for all $f\in H$
$$
E^{n}\exp\int_{0}^{\infty}e^{-t}f(t,x_{t})\,dt
\to E^{ 0}\exp\int_{0}^{\infty}e^{-t}f(t,x_{t})\,dt,
$$
where $E^{ n}$ is the expectation sign with respect to $P^{n}$. Likewise  we will be carrying   to the expectation signs the indices 
the probabilities are supplied with.

If $\Pi$ is a family of probability measures on $(\Omega,\frM^{0})$ we say that $\Pi$ is a compactum
provided that a subsequence $P^{n_{k}}$ can be chosen from any sequence of measures $ P^{n}\in \Pi$
such that $P^{n_{k}}\to P$ for some $ P \in \Pi$.
Let a family $\Pi_{x}$ of probability measures on $(\Omega,\frM^{0})$ be defined for each $x\in E$.

\begin{definition}
                     \label{definition 8.15.1}
We say that the system $\{\Pi_{x}\}$ is a $\bB(E)$-system if for each $x$ the family $\Pi_{x}$ is a compactum
and if for any $n, f_{1},...,f_{n}\in H$ and 
$\alpha_{1},...,\alpha_{n}\geq 0$  we have
\begin{equation}
                                \label{8.10.1}
\sup_{P\in \Pi_{x}}E^{P}\sum_{i=1}^{n}\alpha_{i}\exp
\int_{0}^{\infty}e^{-t}f(t,x_{t})\,dt\in \bB(E)
\end{equation}
\end{definition}
 
It is useful to note for applications that $\{\Pi_{x}\}$ is a $\bB(E)$-system provided $\{\Pi_{x}\}$ is semicontinuous
in $x$; in other words the fact that $x_{n}\to x$ and $P^{n}\in \Pi_{x_{n}}$ implies the existence
of a sequence $P^{n_{k}}$ convergent to some $P\in\Pi_{x}$. The reader can easily verify that in
this case the left side of \eqref{8.10.1} is an upper semicontinuous function and hence belongs
to $\bB(E)$.

In what follows a certain family $\frT$ of bounded  $  \frM_{t} $-stopping   times
is assumed to be given. This family is supposed to include all   constant moments of time.

\begin{definition}
                     \label{definition 8.15.2}
The system $\{\Pi_{x}\}$ of families $\Pi_{x}$ of probability measures is called Markovian relative  to $(\frT,\frM_{t})$ provided the following three conditions are satisfied:

1) $\{\Pi_{x}\}$ is a $\bB(E)$-system and 
$P(x_{0}=x)=1$ for any $x\in E$ and $P\in \Pi_{x}$.

2) For any $x\in E$, $P\in \Pi_{x}$,
$\tau\in\frT$, and $f\in C_{0}$
$$
E^{P}\Big\{ \int_{0}^{\infty}
e^{-t}f(x_{\tau+t})\,dt\mid \frM_{\tau}\Big\}
\leq v(f,\Pi_{x_{\tau} })\quad P-(\rm a.s.),
$$
where
$$
v(f,\Pi_{x }):=\sup_{P\in\Pi_{x}}E^{P}
\int_{0}^{\infty}
e^{-t}f(x_{t})\,dt.
$$

3)   For any $x\in E$, $P\in \Pi_{x}$,
$\tau\in\frT$, and $f\in C_{0}$
\begin{equation}
                                \label{8.10.3}
E^{P}\Big[\int_{0}^{\tau}e^{-t}f(x_{t})\,dt
+e^{-\tau}v(f,\Pi_{x_{\tau }})\Big]\leq v(f,\Pi_{x })
\end{equation}
(as we will see below, $v(f,\Pi_{x })\in\bB(E)$ so that the left-hand side is meaningful).
\end{definition}

This definition has much to do with the problem of
finding $P\in\Pi_{x}$ maximizing
$$
E^{P}
\int_{0}^{\infty}
e^{-t}f(x_{t})\,dt
$$
over {\em controls\/} $P\in\Pi_{x}$. Roughly speaking, property 3)
means that, if  after time $\tau$ we   use any control
$P\in \Pi_{x_{\tau}}$, the combined action
before $\tau$ and after $\tau$ does not give us greater reward than using just a plain $P\in\Pi_{x}$.
Property 2) says that, no matter what we did
before time $\tau$, the conditional gain after $\tau$
cannot be larger than $v(f,\Pi_{x_{\tau} })$, that is the maximum of what we can get just staring from
$x_{\tau}$.

\begin{definition}
            \label{definition 8,15.3}

Let $P_{x}\in \Pi_{x}$ for any $x$ and let $\{P_{x}\}$ be a Markovian system relative  to $(\frT,\frM_{t})$. Then 
we call $X=(x_{t},\frM_{t},P_{x})$
a Markov process. We call it strong Markov if $\frT$ contains all bounded
$\frM_{t}$-stopping times.

\end{definition}

Since the notion of the Markovian property of the system $\{\Pi_{x}\}$ plays a basic role
in the succeeding arguments, it is necessary to present conditions which are sufficient
for the fulfillment of 2) and 3).
Define the operator $\theta_{t}$
\index{$C$@Operators!$\theta_{t}$}%
 acting
on elements of $\cN_{\infty}$ as an operator
preserving all operations on sets and
acting on the generating $\cN_{\infty}$ sets by the formula
$$
\theta_{t}\{\omega:x_{s}(\omega)\in\Gamma\}
=\{\omega:x_{t+s}(\omega)\in\Gamma\}.
$$

Assume that for any $x\in E$, $P\in \Pi_{x}$, and
$\tau\in\frT$ a regular conditional probability
$P(\omega,D)=P(D\mid \frM_{\tau})$, $P$-(a.s.), $D\in\cN_{\infty}$, exists (for appropriate conditions for that see, for instance, Theorem 1.1.6 of \cite{SV_79} or \cite{Kr_73}). Set $P^{\tau}_{\omega}=  
P(\omega,\theta_{\tau}D)$ and assume that the measure
$P^{\tau}_{\omega}$ extends on $\frM^{0}$
in such a way that $P^{\tau}_{\omega}\in
\Pi_{x_{\tau(\omega)}(\omega)}$ $P$-almost all $\omega$. It is clear that in this case
condition 2) is automatically satisfied.

Furthermore, it is shown in Section 3.3 of
\cite{Kr_25}  that it follows from 1) that for any $f\in C_{0}$ and
every $x\in E$ there exists $P_{x}\in\Pi_{x}$ such that $P_{x}(D)\in \bB(E)$ for all $D\in \cN_{\infty}$, and for all $x\in E$
$$
E_{x}\int_{0}^{\infty}e^{-t}f(x_{t})\,dt=v(f,\Pi_{x}).
$$

Therefore condition 3) will be satisfied, if 1)
is satisfied and, each time we have a function
$P_{x}\in\Pi_{x}$ such that $P_{x}(D)\in \bB(E)$ for all $D\in \cN_{\infty}$, we have that for any $y\in E$, $P\in \Pi_{y}$, and $\tau\in\frT$ there exists $P'\in \Pi_{y}$ such that
$$
P'(A\theta_{\tau}D)=EI_{A}P_{x_{\tau}}(D),
\quad \forall A\in\frM_{\tau},D\in\cN_{\infty}.
$$
Indeed in this case
$$
v(f,\Pi_{y})\geq E'\Big[\int_{0}^{\tau}
e^{-t}f(x_{t})\,dt+e^{-\tau}\theta_{\tau}\int_{0}^{\infty}e^{-t}f(x_{t})\,dt\Big]
$$
$$
=E \Big[\int_{0}^{\tau}
e^{-t}f(x_{t})\,dt+e^{-\tau}E_{x_{\tau}}\int_{0}^{\infty}
f(x_{t})\,dt\Big].
$$

Unfortunately, such a ``pasting" of measures $P$ and 
$P_{x}$ is not always possible,
since for $\frM_{t}=\cN_{t+}$ the variable $\tau$ may depend on the ``infinitesimal" future. Therefore,
when checking \eqref{8.10.3} in particular cases it is useful to keep in mind that if $\frM_{t}=\cN_{t+}$ and
$v(f,\Pi_{x})$ is lower semicontinuous and   \eqref{8.10.3} holds for all bounded 
$\cN_{t}$-stopping times, then \eqref{8.10.3} is valid for all bounded $\frM_{t}$-stopping  times.
This follows from Fatou's lemma and from the fact that, in that case, if $\tau$ is an $\frM_{t}$-stopping   time, then $\tau+\varepsilon$ is an
$\cN_{t}$-stopping  time for all $\varepsilon> 0$.  In general, the conditions under which
$v(f,\Pi_{x})$ is lower semicontinuous
are unknown. 

The following theorem gives another sufficient condition for the Markovian property
of a system.

We call a function $\eta[P,\omega]$ $\{\Pi_{x}\}$-{\em
admissible\/}  if $E^{p}\eta[P]$
is finite and continuous on each $\Pi_{x}$ and for any
$n,\alpha_{0},\alpha_{1},...,\alpha_{n}\geq0,f_{1},...,f_{n}\in H$
$$
\sup_{P\in\Pi_{x}}M^{P}\Big[\alpha_{0}\eta[P]
+\sum_{i=1}^{n}\alpha_{i}\exp\int_{0}^{\infty}e^{-t}f(
t,x_{t})\,dt\Big]\in \bB(E).
$$
Observe that
$$
\sum_{i=1}^{n}\alpha_{i}\exp\int_{0}^{\infty}e^{-t}f(
t,x_{t})\,dt
$$
for $f_{i}\in H,\alpha_{i}\geq0$ is $\{\Pi_{x}\}$-admissible
for any $\bB(E)$-system $\{\Pi_{x}\}$.

Let a functional $\xi^{s}_{t}[P,\omega]$ 
$(0 \leq s \leq t \leq\infty)$ be defined for each 
$P \in \bigcup_{x}\Pi_{x}$, $\omega\in\Omega$, and be such that
for any $\tau\in\frT$
$$
\xi^{s}_{\infty}[P,\omega]=\xi^{s}_{\tau}[P,\omega]
+e^{-\tau}\xi^{\tau}_{\infty}[P,\omega].
$$
\begin{theorem}
                         \label{theorem 8.10.1}
Let $\{\Pi_{x}\}$ satisfy condition 1)
and let $\xi^{0}_{\infty}[P,\omega]$
be $\{\Pi_{x}\}$-admissible. Also suppose that
for any $x\in E,P\in\Pi_{x}$, $\tau\in\frT$, $f\in C_{0}$, and $\alpha\geq0$ we have
$$
E^{P}|\xi^{\tau}_{\infty}[P]|<\infty,
$$
$$
E^{P}\Big\{\alpha\xi^{\tau}_{\infty}[P]
+\theta_{\tau}\int_{0}^{\infty}e^{-t}f(x_{t})\,dt
\mid \frM_{\tau}\Big\}\leq w_{\alpha}(f,x_{\tau})
\quad P-(a.s.),
$$
$$
E^{P}\Big\{\alpha\xi^{0}_{\tau}[P]
+ \int_{0}^{\tau}e^{-t}f(x_{t})\,dt
+e^{-\tau}w_{\alpha}(f,x_{\tau})\Big\}\leq w_{\alpha}(f,x),
$$
where
$$
w_{\alpha}(f,x):=\sup_{P\in\Pi_{x}}E^{P}
\Big\{\alpha\xi^{\tau}_{\infty}[P]
+ \int_{0}^{\infty}e^{-t}f(x_{t})\,dt\Big\}
$$
(we will see that $w_{\alpha}(f,\cdot)\in\bB(E)$).

Then the system
$$
\big\{P\in\Pi_{x}: E^{P}\xi^{0}_{\infty}[P]
=\sup_{P\in\Pi_{x}}E^{P}\xi^{0}_{\infty}[P]\big\}
$$
is Markovian relative to $(\frT,\frM_{t})$.

\end{theorem}

This is Theorem 3.3.4 of \cite{Kr_25}
and the following is Theorem 3.3.5 of~\cite{Kr_25}.

\begin{theorem}
                            \label{theorem 8.11.1}
Let $\{\Pi_{x}\}$ be a 
Markovian system relative to $(\frT,\frM_{t})$.
Then there exists a function $P_{x}$ on $E$ such that
$P_{x}\in\Pi_{x}$ for each $x$ and

a) $X=(x_{t}, \frM_{t},P_{x})$ is a Markov process;

b) if, additionally, $\frT$ contains all bounded
$\frM_{t}$-stopping times, then $X$ is a strong Markov
process.
\end{theorem}

\mysection[Markov diffusion processes]{Markov diffusion processes}
                           \label{section 5.14.1}

We are going to use the results of Section
\ref{section 8.15.1}
applied in the case when   $E=\bR^{d+1}$, that is, when the $t$-variable
is considered just as one of the coordinates of points $(t,x)
\in\bR^{d+1}$. Let $\sigma(t,x)$ be Borel 
$\bS_{\delta}$-valued and $b(t,x)$ be Borel $\bR^{d}$-valued functions  
given on $\bR^{d+1}$. 
  Also we suppose that
$b\in L_{(q,p)}$ with $p,q\in(1,\infty)$ satisfying
$$
\frac{d}{p}+\frac{1}{q}\leq 1.
$$

Let $\Omega$ be the set of $\bR^{d+1}$-valued
 continuous function $(t_{0}+t,x_{t})$, $t_{0}\in \bR$,
defined for $t\in[0,\infty)$.
For $\omega=\{(t_{0}+t,x_{t}),t\geq0 \}$, define
$\sft_{t}(\omega)=t_{0}+t$, $x_{t}(\omega)=x_{t}$,
and set $ \cN_{t}=\sigma((\sft_{s},x_{s}),s\leq t)$,
$  \cN_{\infty}= \sigma((\sft_{s},x_{s}),s< \infty)$. Denote by $\frT$ the set of {\em bounded\/} stopping times
relative to $\{\cN_t\}$.  

In the following result we use the terminology from
\cite{Dy_63}. In the time-homogeneous case 
with bounded $b$ Theorem
\ref{theorem 4.27.1}  was obtained in \cite{Kr_73_1}.  
If $b$ is bounded, but there also
jumps, the result is found in \cite{AP_77}
and if $p=q$ in \cite{GM_01}.

\begin{theorem}
                             \label{theorem 4.27.1}
 On $\bR^{d+1}$ there exists a strong Markov process
$$
X=\{(\sft_{t},x_{t}) ,\cN_{t}, P_{t,x})
$$
such that   for any $(t,x)\in\bR^{d+1}$
there exists a $d$-dimensional Wiener process $w_{t}$, $t\geq0$,
which is a Wiener process relative to $\bar \cN_{t}$,
where $\bar \cN_{t}$ is the completion
\index{$B$@Sets!$\bar \cN_{t}$}%
 of $\cN_{t}$
with respect to all $P_{s,y}$, and such that with 
$P_{t,x}$-probability one, for
all $s\geq 0$, $\sft_{s}=t+s$ and
\begin{equation} 
                             \label{4.27.10}
x_{s}=x+\int_{0}^{s}\sigma(t+u,x_{u})\,dw_{u}
+\int_{0}^{s}b(t+u,x_{u})\,du.
\end{equation}
\end{theorem}

Proof. Define $a=\sigma^{2}$,
$$
\cL u(t,x)=\partial_{t}u(t,x)+(1/2)a^{ij}D_{ij}u(t,x)+b^{i}D_{i}u(t,x)
$$
and introduce $\Pi_{t,x}$ as the set of probability
measures on $(\Omega,\cN_{\infty})$ such that $P((\sft_{0},x_{0})=(t,x))=1$, $P$-(a.s.)
\begin{equation}
                                                     \label{4.27.5}
 \int_{0}^{T}|b(\sft_{s},x_{s})|\,ds<\infty,\quad\forall T<\infty,
\end{equation}
and the process
$$
\eta_{t}(u)=u(\sft_{t},x_{t})-\int_{0}^{t}\cL u(\sft_{s},x_{s})\,ds-u(\sft_{0},x_{0})
$$
is a local martingale relative to $\{\cN_{t}\}$ for all $u\in C^{\infty}
_{0}(\bR^{d+1})$. In the terminology of
Stroock-Varadhan $P$ is a solution of a martingale
problem.

Right away observe that owing to 
Lemma 3.4.1 of \cite{Kr_25},  if $P \in\Pi_{t,x}$, then $x_{t}$ is a solution of \eqref{4.27.10}  and in light of
It\^o's formula, $\eta_{t}(u)$ is a square-integrable
martingale and
\begin{equation}
                                                     \label{8.16.1}
E \sup_{t\in[0,T]} \eta_{t}^{2}(u) \,dt\leq N(d,\delta,Du)T,\quad\forall T<\infty.
\end{equation}

By Theorem \ref{theorem 8.11.1}   to prove
the present theorem, it suffices to show that $\Pi_{t,x}\ne\emptyset$
and
$\{\Pi_{t,x}\}$ is a Markovian system relative to $(\frT,\cN_t)$.

That $\Pi_{t,x}\ne\emptyset$ follows from Theorem \ref{theorem 9.6.4} (i)  and It\^o's formula.
Let us prove that $\{\Pi_{t,x}\}$ is a $\bB(E)$-system. To achieve this,
as it is observed after \eqref{8.10.1}, it suffices to show that
if $(t_{n},x_{n})\to (t,x)$ and $P^{n}\in\Pi_{t_{n},x_{n}}$,
then there exists a subsequence $n(k)\to\infty$ and $P^{0}\in\Pi_{t,x}$
such that for any $f\in C^{\infty}_{0}(\bR^{d+2})$
$$
E^{ n(k) }\exp\Big(\int_{0}^{\infty}e^{-t}f(t,\sft_{t},x_{t})\,dt\Big)
\to E^{0}\exp\Big(\int_{0}^{\infty}e^{-t}f(t,\sft_{t},x_{t})\,dt\Big),
$$
where $E^{ n(k) },E^{0}$ are the expectation signs with respect
to $P^{ n(k) },P^{0}$, respectively.
The reader will easily derive this property   
from Theorem \ref{theorem 9.6.4} (ii) by 
using Taylor's series and
observing that
$$
E\Big(\int_{0}^{\infty}e^{-t}f(t,\sft_{t},x_{t})\,dt\Big)^{n}
$$
$$
=E \int_{0}^{\infty}...\int_{0}^{\infty}
e^{-t_{1}}f(t_{1},\sft_{t_{1}},x_{t_{1}})\cdot...
\cdot e^{-t_{n}}f(t_{n},\sft_{t_{n}},x_{t_{n}})\,
dt_{1}\cdot...\cdot dt_{n}.
$$

What remains is to prove that for $(\frT,\cN_{t})$  the conditions 2) and 3)
in Definition \ref{definition 8.15.2} are satisfied.   
In our space-time situation the operators $\theta_{t}$
are defined starting with
$$
\theta_{t}\{\omega:(\sft_{s},x_{s})(\omega) \in\Gamma\}
=\{\omega:(\sft_{t+s},x_{t+s})(\omega)\in\Gamma\},
$$
which is naturally extended to all sets in $\cN_{\infty}$.
Introduce
$$
\zeta[P,f,\omega]=\zeta(f)=\int_{0}^{\infty}e^{-t}
f(t,\sft_{t},x_{t})\,dt.
$$

By Theorem 1.1.6 of \cite{SV_79}
or \cite{Kr_73} for any $P\in \Pi_{t,x}$
and $\tau\in\cT$ there exists a regular conditional probability $P(\omega,A)$ relative to $\cN _{\tau}$.
Define a measure $P^{\tau}_{\omega}$ on $\cN_{\infty}$
by the formula $P^{\tau}_{\omega}(A)=P(\omega,\theta_{\tau}A)$. Then 
$$
E \big\{\theta_{\tau}\zeta(f)\mid \cN_{\tau}\big\}
=E^{\tau}_{\omega}\zeta(f)\quad P\rm-(a.s.).
$$
Therefore to prove that $\{\Pi_{t,x}\}$
possesses property 2), it suffices to show that
$P^{\tau}_{\omega}\in\Pi_{\sft_{\tau},x_{\tau}}$ $P$-(a.s.).

For any $f\in C_{0}(\bR^{d+1})$, $\tau\in\cT$, and $A\in\cN_{\tau}$ we have
$$
\int_{A}E^{\tau}_{\omega}f(\sft_{0},x_{0})\,P(d\omega)=
\int_{A}E \big\{\theta_{\tau}f(\sft_{0},x_{0})
\mid \cN_{\tau}\big\}\,P(d\omega)
$$
$$
=
\int_{A}f(\sft_{\tau},x_{\tau})\,P(d\omega).
$$
It follows that $E^{\tau}_{\omega}f(\sft_{0},x_{0})
=f(\sft_{\tau},x_{\tau})$ $P$-(a.s.) and, due to the arbitrariness
of $f$,
$P^{\tau}_{\omega'}\big(\sft_{0}=\sft_{\tau(\omega')}(\omega'),x_{0}=x_{\tau(\omega')}(\omega')\big)=1$ for $P$-almost all $\omega'$.
Also clearly, \eqref{4.27.5} holds with 
$P^{\tau}_{\omega}$ in place of $P$
for $P$-almost all $\omega$.

Next, take $u\in C^{\infty}_{0}(\bR^{d+1})$, $\tau\in\cT$, 
rational numbers such that $0\leq s_{1}\leq...\leq s_{n}= s<t$ and a box in $Q\subset\bR^{nd}$ with edges parallel to the coordinate axes and vertices with rational coordinates. Note that 
$$
\theta_{\tau}\eta_{t}(u)=u(\sft_{\tau+t},x_{\tau+t})
-\int_{\tau}^{\tau+t}\cL u(\sft_{s},x_{s})\,ds
-u(\sft_{\tau},x_{\tau})
=\eta_{\tau+t}(u)-\eta_{\tau}(u)
$$
so that
$P$-(a.s.)
$$
E^{\tau}_{\omega}I_{Q}((\sft,x)_{s_{1}},...,(\sft,x)_{s _{n}})
\eta_{t}(u)=E\big\{\theta_{\tau}\big(
I_{Q}((\sft,x)_{s_{1}},...,(\sft,x)_{s _{n}})
\eta_{t}(u)\big)\mid \cN_{\tau}\big\}
$$
$$
=E\big\{  
I_{Q}((\sft,x)_{\tau+s_{1}},...,(\sft,x)_{\tau+s _{n}})\big(
\eta_{\tau+t}(u)-\eta_{\tau}(u)\big)\mid \cN_{\tau}\big\}
$$
$$
=E\big\{  
I_{Q}((\sft,x)_{\tau+s_{1}},...,(\sft,x)_{\tau+s _{n}})\big(
\eta_{\tau+s}(u)-\eta_{\tau}(u)\big)\mid \cN_{\tau}\big\}
$$
$$
=E^{\tau}_{\omega}I_{Q}((\sft,x)_{s_{1}},...,(\sft,x)_{s _{n}})
\eta_{s}(u).
$$
In short, for our choice of $n,s_{i},t$ and $Q$
there exists an event $\Omega'$ of full $P$-probability such that for all $\omega\in \Omega'$
\begin{equation}
                                \label{8.16.2}
E^{\tau}_{\omega}I_{Q}((\sft,x)_{s_{1}},...,(\sft,x)_{s _{n}})
\eta_{t}(u)=E^{\tau}_{\omega}I_{Q}((\sft,x)_{s_{1}},...,(\sft,x)_{s _{n}})
\eta_{s}(u).
\end{equation}
Since there are only countably many such
$n,s_{i},t$ and $Q$, there is a smaller $\Omega'$
of full $P$-probability such that \eqref{8.16.2}
holds for all $n,s_{i},t$ and $Q$ with the above
properties. Then the usual measure-theoretic
argument shows that
\begin{equation}
                                \label{8.16.3}
E^{\tau}_{\omega}I_{A} 
\eta_{t}(u)=E^{\tau}_{\omega}I_{A} 
\eta_{s}(u)
\end{equation}
for any $A\in\cN_{s}$ as long as $\omega\in\Omega'$
and $s$ and $t$ are rational. Extending $s$
to be any number $<t$ is trivially possible
since $(\sft,x)_{t}$ is continuous. Extending the range
of $t$ is possible due to \eqref{8.16.1}
and perhaps requires reducing the $\Omega'$
to a different set of $P$-measure 1.
From \eqref{8.16.3} we conclude that
$\eta_{t}(u)$ is  a $P^{\tau}_{\omega}$-martingale
$P$-(a.s.).

This shows that
$P^{\tau}_{\omega}\in\Pi_{\sft_{\tau},x_{\tau}}$ $P$-(a.s.) and condition 2) is satisfied for $\{\Pi_{t,x}\}$.

As we have explained in Section \ref{section 8.15.1},
to verify that condition 3) is satisfied, it suffices to make sure that,  each time we have a function
$P_{t,x}\in\Pi_{t,x}$ such that $P_{t,x}(D) $
is Borel for all $D\in \cN_{\infty}$ and we have  $P\in \Pi_{s,y}$ for some $(s,y)\in \bR^{d+1}$  and $\tau\in\cT$, there exists $P'\in \Pi_{s,y}$ such that
\begin{equation}
                                  \label{8.16.8}
P'(A\theta_{\tau}D)=E I_{A}P_{\sft_{\tau},x_{\tau}}(D),
\quad \forall A\in\cN_{\tau},D\in\cN.
\end{equation}

If $\omega',\omega\in\Omega$, set $\omega'\theta_{\tau}\omega$ to be the function $(\sft,x)'_{t}$
(generally not in $\Omega$) such that $(\sft,x)'_{t}=(\sft,x)_{t}
(\omega')$ for $t\leq \tau(\omega')$ and
$(\sft,x)'_{t}=(\sft,x)_{t-\tau(\omega')}(\omega)$
for $t> \tau(\omega')$.

For any $\cN_{\infty}$-measurable $\xi(\omega)$ and $\omega'
\in\Omega$ set
$$
\xi_{\omega'}(\omega)=\xi(\omega'\theta_{\tau}\omega)
I_{\Omega}(\omega'\theta_{\tau}\omega).
$$
Finally, if $A\in\cN_{\infty}$ define $\xi^{A}(\omega):=I_{A}(\omega)$ and
\begin{equation}
                           \label{8.16.5}
P'(A)=\int_{\Omega}E_{(\sft,x)_{\tau(\omega')}(\omega')}
\xi^{A}_{\omega'}\,P(d\omega').
\end{equation}

We are going to show that $P'$ is well defined and 
is the measure we are after.

Let $f_{1} ,...,f_{n} $ be bounded continuous functions on $\bR^{d+1}$ and $0=t_{0}\leq t_{1}\leq...\leq t_{n}<t_{n+1}=\infty$.
Clearly the function
$$
E_{t,x}\prod_{i=j}^{n}f_{i}((\sft,x)_{t_{i}-s})
$$
is left-continuous in $s$ for $s\leq t_{j}$
and is Borel in $(t,x)$. Hence, it is jointly
measurable. It follows for $\xi=\prod_{i=1}^{n}f_{i}
((\sft,x)_{t_{i}})$ that ($\prod_{1}^{0}=\prod^{n}_{n+1}:=1$)
$$
E_{(\sft,x)_{\tau(\omega')}(\omega')}
\xi _{\omega'}
$$
$$=\sum_{j=1}^{n+1}I_{t_{j-1}\leq
\tau(\omega')<t_{j}}\prod_{i=1}^{j-1}f_{i}((\sft,x)_{t_{i}}
(\omega'))E_{(\sft,x)_{\tau(\omega')}(\omega')}
\prod_{i=j}^{n}f_{i}((\sft,x)_{t_{i}-\tau(\omega')})
$$
is measurable with respect to $\omega'$. Then
usual argument shows that the function  $E_{(\sft,x)_{\tau(\omega')}(\omega')}
\xi _{\omega'}$ is measurable with respect to $\omega'$ for any bounded $\cN_{\infty}$-measurable $\xi$. 
Hence, the right-hand side of \eqref{8.16.5}
is well defined and, obviously, gives a probability
measure on $\Omega$.

Next, we need a property of $\cN_{\tau}$-measurable functions. Let $\omega' $ and $\omega''$ be {\em fixed\/} such that
$(\sft,x)_{t}(\omega'  )=(\sft,x)_{t}(\omega'')$ for
$t\leq \tau(\omega')$. Clearly, the set of functions
$\xi(\omega)$ such that $\xi(\omega')=\xi(\omega'' )$ contains all functions of the type $\prod_{i=1}^{n}f_{i}
((\sft,x)_{t_{i}})$  if $t_{n}\leq \tau(\omega')$. Then this set contains all $\cN_{\tau(\omega')}$-measurable functions ($\tau(\omega')$ is a fixed number). Consequently, if $\xi$
is $\cN_{\tau}$-measurable, then
\begin{equation}
                                 \label{8.16.7}
\xi(\omega'')I_{\tau(\omega'')=\tau(\omega')}
=\xi(\omega')I_{\tau(\omega')=\tau(\omega')}=\xi(\omega')
\end{equation}
since $\xi(\omega)I_{\tau(\omega)=\tau(\omega')}$
is $\cN_{\tau(\omega')}$-measurable. For $\xi\equiv1$
this yields $\tau(\omega'')=\tau(\omega')$ and coming back to \eqref{8.16.7} we get that $\xi(\omega'')
=\xi(\omega' )$.

Having this in mind take an $\cN_{\tau}$-measurable
$\xi$ and observe that
$$
\Big(\xi\theta_{\tau}\prod_{i=1}^{n}f_{i}((\sft,x)_{t_{i}})
\Big)_{\omega'}(\omega)=\Big(
\xi \prod_{i=1}^{n}f_{i}((\sft,x)_{\tau+t_{i}})
\Big)(\omega'\theta_{\tau}\omega)
$$
$$
=\xi(\omega'\theta_{\tau}\omega)
\prod_{i=1}^{n}f_{i}\big((\sft,x)_{\tau(\omega'\theta_{\tau}\omega)+t_{i}}(\omega'\theta_{\tau}\omega)\big)
$$
$$
=I_{(\sft,x)_{\tau(\omega')}(\omega')=(\sft,x)_{0}(\omega)}
\xi(\omega')\prod_{i=1}^{n}f_{i}\big((\sft,x)_{t_{i}}
(\omega)\big).
$$
It follows that for any $\cN_{\infty}$-measurable $\eta$
$$
(\xi\theta\eta)_{\omega'}(\omega)
=I_{(\sft,x)_{\tau(\omega')}(\omega')=(\sft,x)_{0}(\omega)}
\xi(\omega')\eta(\omega).
$$
By taking here $\xi$ and $\eta$ as the indicators
of  appropriate sets and using \eqref{8.16.5}
we get \eqref{8.16.8}.

Now it remains to prove that $P'\in \Pi_{s,y}$.
That $P'((\sft_{0},x_{0})=(t,x))=1$ is obvious.
To check \eqref{4.27.5} write
$$
P'\Big(\int_{0}^{T}|b(\sft_{t},x_{t})|\,dt<\infty\Big)
$$
$$
\geq
P'\Big(\int_{0}^{\tau}|b(\sft_{t},x_{t})|\,dt<\infty,
\theta_{\tau}\Big\{
\int_{0}^{T}|b(\sft_{t},x_{t})|\,dt<\infty\Big\}\Big)
$$
$$
=EI_{[0,\infty)}\Big(\int_{0}^{\tau}|b(\sft_{t},x_{t})|\,dt\Big)P_{\sft_{\tau},x_{\tau}}
\Big(\int_{0}^{T}|b(\sft_{t},x_{t})|\,dt<\infty\Big)=1.
$$

To check that $\eta_{t}(u)$ is a   martingale
with respect to $P'$ 
first observe that, as it follows from the definition of $P'$, for any $A\in \cN_{\tau}$ and  
$\cN_{\infty}$ measurable $f((\sft,x)_{\cdot})\geq0$
$$
E'I_{A}f((\sft,x)_{\tau+\cdot})=EI_{A}
 E_{\sft_{\tau},x_{\tau}}f.
$$
By considering the $f$'s which are the products
of $f(t)$ and $f((\sft,x)_{\cdot})$ and then using
well-known techniques one proves that
for any $\cB([0,\infty))\times\cN_{\infty}$ measurable $f_{t}((\sft,x)_{\cdot})\geq0$ we have
$$
E'I_{A}f_{\tau}((\sft,x)_{\tau+\cdot})=EI_{A}
 E_{\sft_{\tau},x_{\tau}}f_{t}\big|_{t=\tau}.
$$

Next,
for any $T\in(0,\infty)$
$$
E'\sup_{t\leq T}|\eta_{t}(u)|
\leq I_{1}+I_{2},
$$
where 
$$
I_{1}=E\sup_{t\leq\tau}|\eta _{t}(u)|\leq
N(d,\delta,Du)\sup_{\omega}\tau^{1/2},
$$
$$
I_{2}=E'I_{\tau\leq T}\big(|\eta_{\tau}(u)|+\theta_{\tau}
\sup_{t\in[0,T]}|\eta_{t}(u)|\big)
$$
$$
\leq N(d,\delta,Du)T^{1/2}
+EE_{\sft_{\tau},x_{\tau}}\sup_{t\in[0,T]}|\eta_{t}(u)|
\leq N(d,\delta,Du)T^{1/2}.
$$

Finally, take $0\leq s \leq t$ and  $A\in\cN_{s}$ and write
$$
E'I_{A}\eta_{t}(u)=J_{1}+J_{2},
$$
where
$$
J_{1}=E'I_{A,\tau\leq t}\big(\eta_{\tau}(u)+
\theta_{\tau}\eta_{t-\kappa}(u)\big|_{\kappa=\tau}\big)
$$
 $$
=EI_{A,\tau\leq t} \eta_{\tau}(u)
+EI_{A,\tau\leq s}E_{\sft_{\tau},x_{\tau}}
\eta_{t-\kappa}(u)\big|_{\kappa=\tau}
=EI_{A,\tau\leq t} \eta_{\tau}(u),
$$
$$
J_{2}=E'I_{A,\tau> t}\eta_{t}(u)
=EI_{A,\tau> t}\eta_{t}(u).
$$
Hence, 
$$
E'I_{A}\eta_{t}(u)=EI_{A}\eta_{t\wedge\tau}=
EI_{A}\eta_{s\wedge\tau}=E'I_{A}\eta_{s}(u).
$$
The theorem is proved.\qed

\mychapter[It\^o
 processes with moderated drift]
{Nondegenerate It\^o
 processes with moderated drift}

                  \label{chapter 10.20.1}

\mysection[Introduction]{Introduction}

Let  $d_{1}$ be an integer $\geq d$,
$(\Omega,\cF,P)$ be a complete probability space,
and let $(w_{t},\cF_{t})$ be a $d_{1}$-dimensional
Wiener process on this space with complete, relative to
$\cF,P$, $\sigma$-fields $\cF_{t}$. Fix $\delta\in(0,1]$. Let $\sigma_{t},t\geq0$,
be a progressively measurable process with values in the set 
of $d\times d_{1}$-matrices such that
$a_{t}:= \sigma_{t}\sigma^{*}_{t}\in \bS_{\delta}$ for all $(\omega,t)$, and let $b_{t},t\geq0$, be an $\bR^{d}$-valued
progressively  measurable process.
Assume that for any $T\in[0,\infty)$ and $\omega$
$$
\int_{0}^{T} |b_{t}| \,dt<\infty.
$$

Under this condition the stochastic process  
$$
x_{t}=\int_{0}^{t}\sigma_{s}\,dw_{s}
+\int_{0}^{t}b_{s}\,ds
$$
is well defined.

Define one of the main
 quantities, we will be using,
by
$$
\bar b_{R}=\sup_{\rho\leq R} b'_{\rho},\quad
  b'_{\rho}:=\frac{1}{\rho}\sup_{x\in\bR^{d}} 
\sup_{ t\geq0}
\esssup E_{\cF_{t}}\int_{0}^{\theta_{t}\tau_{ \rho}(x)}
|b_{t+s} |\,ds  .
$$
One can say that $\bar b_{\rho}$ ``moderates'' 
\index{$S$@Miscelenea!$\bar b_{R}$}%
$b$ on scale $\rho$. Sufficient analytic conditions
for $b'_{R}$ to be finite can be found
by using Lemma \ref{lemma 3.25.1}.
In particular, this happens if $b$
is bounded.

\begin{remark}
                              \label{remark 10.14.1}
For any $\rho>0$ and stopping time $\tau$ 
and $\cF_{\tau}$-measurable $\bR^{d}$-valued $y$ we have
$$
E_{\cF_{\tau}}\int_{0}^{\theta_{\tau}\tau_{ \rho}(y)}
|b_{\tau+s} |\,ds
\leq \bar b_{\rho}\rho
$$
or, in other words, for any $A\in \cF_{\tau}$
\begin{equation}
                                   \label{10.14.1}
EI_{A}\int_{\tau}^{\tau+\theta_{\tau}\tau_{ \rho}(y)}
|b_{ s} |\,ds
\leq \bar b_{\rho}\rho P(A).
\end{equation}

Indeed, if $\tau$ and $y$ take only countably many values
(including $\infty$ for $\tau$),
estimate \eqref{10.14.1} immediately follows from the definition of $\bar b_{\rho}$. 
In the case of general $\tau$, one knows that there exists
a sequence of stopping times $\tau_{n}$
with values in the set of dyadic rationals such that
$\tau_{n}\downarrow \tau$. Furthermore, as is easy
to see, for any $\varepsilon\in(0,\rho)$
$$
\nliminf_{n\to\infty}\theta_{\tau_{n}}\tau_{ \rho}(y) 
\geq \theta_{\tau }\tau_{ \rho-\varepsilon}(y) 
, 
$$
$$ \int_{\tau}^{\tau+\theta_{\tau}\tau_{ \rho-\varepsilon}(y)} 
|b_{ s} |\,ds\leq \nliminf_{n\to\infty}
\int_{\tau_{n}}^{\tau_{n}+\theta_{\tau_{n}}\tau_{ \rho }(y)} 
|b_{ s} |\,ds
$$
and Fatou's lemma shows that \eqref{10.14.1} holds
if we replace $\rho$ in its left hand side with $\rho-\varepsilon$. After this replacement, still in the case of discrete
$y$, it will only remain
to use the monotone convergence theorem sending $\varepsilon\downarrow0$. To pass to general $y$ we approximate
it with the discrete ones $y_{n}$ and use the fact that
$$
\nliminf_{n\to\infty}\theta_{\tau}\tau_{\rho}(y_{n})
\geq \theta_{\tau }\tau_{ \rho-\varepsilon}(y) .
$$
 
\end{remark}

Next assumption, in which $\sfb_{0}=\sfb_{0}(d,\delta)\in(0,1]$ is a   number
to be specified later in Theorem \ref{theorem 8.2.1}, is  {\em supposed to hold throughout
this chapter after Theorem \ref{theorem 8.2.1}.}
\begin{assumption} 
            \label{assumption 8.19.2}
  We have a $\rho_{b}\in(0,\infty)$ such that $\bar b_{\rho_{b}}
\leq \sfb_{0}$.
\end{assumption}

\begin{remark}
                         \label{remark 8.21.01}
A very important feature of this assumption is
that it is preserved under self-similar {\em dilations\/}.
To be more precise, take a constant $c\in (0,1] $
and introduce $ \hat x_{t} =
 c^{-1} x_{c^{2} t} $. Then for
$$
\hat\sigma_{s} =\sigma_{c^{2}s} ,\quad \hat b_{s} =
cb_{c^{2}s} ,\quad \hat w_{s}=c^{-1}w_{c^{2}s}.
$$
we have
$$
\int_{0}^{T}|\hat b_{t} |\,dt<\infty,
$$
$$
\hat x_{t}= \int_{0}^{t}\hat\sigma_{s} \,
d\hat w_{s}+\int_{0}^{t}\hat b_{s} \,ds,
$$
and $\hat w_{s}$ is a Wiener process relative to
$\{\cF_{c^{2}t}\}$.
 
In addition, 
for $\theta_{t}\hat \tau_{ R}$ being the 
minimum of $R^{2}$ and the
first exit time of $(t+s,\hat x_{t+s})$ from $C_{R}(t,\hat x_{t})$ we have $\theta_{t}\hat \tau_{R}=c^{-2}\theta_{c^{2}t}\tau_{C_{cR}}$ so that
$$
E_{\cF_{c^{2}t}}\int_{0}^{\theta_{t}\hat \tau_{ \rho}}
|\hat b_{t+s} |\,ds  
=cE_{\cF_{c^{2}t}}\int_{0}^{c^{-2}\theta_{c^{2}t}\tau_{ c\rho}}
|  b_{c^{2}t+c^{2}s} |\,ds  
$$
$$
=c^{-1}E_{\cF_{c^{2}t}}\int_{0}^{ \theta_{c^{2}t}\tau_{ c\rho}} 
|  b_{c^{2}t+ s} |\,ds  
\leq c^{-1}\bar b_{c\rho}(c\rho)=
\bar b_{c\rho} \rho \leq \sfb_{0}  \rho .
$$
 Here for simplicity we considered $\theta_{t}\hat \tau_{ R}(x)$ only for $x=0$. The general case is not hard either.

\end{remark}

\begin{remark}
                                       \label{remark 3.27.1}
Usual way to deal with additive functionals shows 
(see, for instance,   \eqref{8.18.1}) that
for any $n=1,2,...$, $\rho,t\in[0,\infty)$, $x\in\bR^{d}$,
$$
E_{\cF_{t}} \Big(\int_{0}^{\theta_{t}\tau_{ \rho }(x)}
|b_{t+s}|\,ds\Big)^{n} \leq n!\,\bar b_{\rho}^{n}\rho^{n}.
$$
Furthermore, by taking into account that
for any random variable $\xi\geq0$ and $\alpha
\in[1,2]$ we have $\big(E\xi^{\alpha}\big)^{1/\alpha}
\leq \big(E\xi \big)^{(2-\alpha)/\alpha}
\big(E\xi^{2}\big)^{(\alpha-1)/\alpha}$, we find that
for any $\varepsilon>0$ there exists $\alpha=\alpha(\varepsilon)
>1$ such that  
$$
\Big(E_{\cF_{t}} \Big(\int_{0}^{\theta_{t}\tau_{ \rho}(x)} 
|b_{t+s}|\,ds\Big)^{\alpha} 
 \Big)^{1/\alpha}\leq 2^{(\alpha-1)/\alpha}
\bar b_{\rho}\rho \leq (1+\varepsilon)\bar b_{\rho}\rho .
$$
 
\end{remark}

Recall that
\begin{equation}
                                  \label{7.1.1}
\dashnorm h\|_{L_{(q,p)}(C)}=\|h\|_{L_{(q,p)}(C)}\|
1\|_{L_{(q,p)}(C)}^{-1}.
\end{equation}

\begin{definition}
            \label{definition 2.3.2}
Let $\kappa,q,p$ be properly tight,
$\hat b\leq 1,\rho>0$ be some constants,
\index{$S$@Miscelenea!properly tight}%
 $b $ be a function
on $\bR^{d+1}$, and $\rho\in(0,\infty)$.  
\index{$S$@Miscelenea!$\bar b_{R}$@$\hat b_{(q,p),\rho}$}%
Introduce
\begin{equation}
                            \label{8.19.30}
 \hat b_{(q,p),\rho}=\sup_{r\leq \rho}r\sup_{C\in \bC_{r}}
\dashnorm b\|_{L_{(q ,p) }(C)} .
\end{equation}   
 
\end{definition}

In the next lemma we give a simple analytic condition
for Assumption     \ref{assumption 8.19.2} to be satisfied.
\begin{lemma}
                       \label{lemma 12.11.1}
 Let $d,q,p$ be properly tight, $\nu(d,q,p)=0$,
and let a number $\rho_{b}\in (0,\infty)$.
 Suppose that $\hat b_{(q,p),\rho_{b}}
\leq 1$ and  $N 
\hat b_{(q,p),\rho_{b}}
\leq \sfb_{0}$, where $N=N(d,\delta,p ,q)$ is 
taken from Lemma \ref{lemma 3.25.1}. 
Then Assumption     \ref{assumption 8.19.2}    
  is satisfied.

\end{lemma}

This lemma is a trivial consequence of
Lemma \ref{lemma 3.25.1}.

\begin{remark}
                              \label{remark 12.11.1}

It turns out that \eqref{8.19.30} with $q/2$ in place of  $q$
can be made
 as small  
as we like on  account of taking $\rho$
small if $\|b\|_{L_{(q,p)} }<\infty$ for
some $p,q\in(1,\infty)$
satisfying $d/p+2/q=1$, which implies that
$d,q/2,p$ are properly tight with $\nu=0$. Indeed, by H\"older's
inequality, if $C\in\bC_{\rho}$, then
$$
\dashnorm b\|_{L_{(q/2,p)}(C)}\leq 
\dashnorm b\|_{L_{(q,p)}(C)}=
N(d)\rho^{-1}\|b\|_{L_{(q,p)}
(C) },
$$
where the last norm tends to zero as $\rho\to0$.
\end{remark}

However, there are many situations when
\eqref{8.19.30} is not finite but 
Assumption    \ref{assumption 8.19.2}
  is still satisfied.

\begin{example}
                    \label{example 5.23.1}

In $\bR^{d}$ (with $d\geq 2$) take a $d$-dimensional
Wiener process $w_{t}$ and consider the system
$dx^{1}_{t}=dw^{1}_{t}+b(x_{t}) \,dt$,
$dx^{i}_{t}=dw^{i}_{t}$, $i\geq2$, where
$$
b(x )=\beta(x^{1} ),\quad \beta(r)=-|r|^{-\alpha}I_{(-1,1)}(r) \sign r,
$$ 
  and $\alpha<1$ is as close
to $1$ as we wish. The solutions of our system
form a strong Markov time-homogeneous process for which
Assumption  \ref{assumption 8.19.2}
  is rewritten
as 
\begin{equation}
                                 \label{10.10.1}
\bar b_{\rho_{b}}:=\sup_{\substack{\rho\leq \rho_{b} 
\\ \,C\in\bC_{\rho}}}
\frac{1}{\rho}
\sup_{ x \in\bR^{d }}E_{ x}\int_{0}^{\tau_{C}}
|b( x_{s})|\,ds\leq \sfb_{0},
\end{equation}
where $\tau_{C}$ is the 
\index{$S$@Miscelenea!$\bar b_{R}$@$\bar b_{\rho_{b}}$}%
first exit time of $(t,x_{t})$
from $C$.

Here $\hat b_{(q,p),\rho_{b}}$
is definitely infinite if
$\kappa>1$ and $\alpha$ is too close to $1$. 
However, if
$C=C_{\rho}(s,y)$ and $|y^{1} |\leq 2\rho$,
then with $P_{ x}$-probability one
$\tau_{C}\leq \tau_{(-3\rho,3\rho)}$, where by
$\tau_{(a,b)}$ we denote the first exit time of $x^{1}_{t}$ from $(a,b)$. In that case by using It\^o's formula one gets that, for $(0,x)\in C$
$$
 E_{ x}\int_{0}^{\tau_{C}}|b(x_{s})|\,ds\leq E_{ x}\int_{0}^{\tau_{(-3\rho,3\rho)}}
|b(x_{s})|\,ds=:\phi(|x^{1}|),
$$
where $\phi(r)=0$ if $|x^{1}|\geq 3\rho$ and  otherwise
$$
\phi(r)=r- 3\rho+ \int_{r}^{3\rho}
\exp\Big(\frac{2}{1-\alpha}(t\wedge 1)^{1-\alpha}\Big)
\,dt\leq N\rho.
$$

In case $ y^{1} \geq 2\rho$ observe that,
for $r_{\pm}:=y^{1}\pm \rho$,
$$
E_{ x}\int_{0}^{\tau_{C}}|b(x_{s})|\,ds\leq E_{ x}\int_{0}^{\tau_{( r_{-},r_{+})}}
|b(x_{s})|\,ds=: \psi(x^{1}) ,
$$
which is zero if $x^{1}\not \in ( r_{-},r_{+}) $
and if $x^{1}  \in ( r_{-},r_{+}) $
 by It\^o's formula equals (observe that $|\beta(r)|=
-\beta(r)$ on $( r_{-},r_{+}) $)
$$
\psi(x^{1})=-E_{t,x}(x^{1}_{\tau_{( r_{-},r_{+})}}-x^{1}_{0})
$$
which  is less than $2\rho$.
Therefore, $\psi(x^{1})\leq 2\rho$ and,
since similar situation occurs if $y^{1}\leq -2\rho$,
$\bar b_{\infty}<\infty$.

To show that \eqref{10.10.1} is satisfied 
we show that $\bar b_{\rho}\to0$ as $\rho\to0$.

Observe that $\phi(r)$ is a decreasing function for
$r\geq0$ and if $3\rho\leq1$
$$
\phi(r)\leq\phi(0)=\int_{0}^{3\rho}\Big[
\exp\Big(\frac{2}{1-\alpha}(t\wedge 1)^{1-\alpha}\Big)
-1\Big]\,dt
$$
$$
\leq 3\rho\Big[
\exp\Big(\frac{2}{1-\alpha}(3\rho)^{1-\alpha}\Big)
-1\Big]\leq N\rho^{2-\alpha}.
$$

To estimate $\psi$ for $r \in ( r_{-},r_{+}) $
introduce
$$
 \xi(r)=r-r_{+}+\eta(r),\quad \eta(r)=\frac{2\rho}{  \gamma -1}
\big(e^{2\hat b(r-r_{+})}-1 \big),
$$
where
$$
\hat b= ( 2/|y^{1}| )^{\alpha},\quad \gamma=e^{-4\hat b\rho}.
$$
By observing that on $( r_{-},r_{+}) $ we have $\beta(r)\geq -\bar b$ and $\eta'\leq0$, we obtain
$$
(1/2)\eta''+b\eta'\leq (1/2)\eta''-\hat b \eta'=0,
\quad (1/2)\xi''+b\xi'\leq b=-|b|.
$$
Furthermore, $\xi(r_{\pm})=0$ and a simple application 
of It\^o's formula shows that
$\psi\leq \xi$ on $ ( r_{-},r_{+}) $.

To estimate $\xi$ use that $e^{t}-1\geq t$
implying that $\eta(r)\leq 4\rho\hat b(r_{+}-r)(1-\gamma)^{-1}$, so that
$$
\xi(r)\leq (r_{+}-r)\big(4\rho\hat b(1-\gamma)^{-1}-1\big)\leq 2\rho \big(4\rho\hat b(1-\gamma)^{-1}-1\big).
$$
Note that $\hat b\rho\leq \rho^{1-\alpha}$ since $y^{1}
\geq 2\rho$ and $\psi(r)\leq\xi(r)\leq 10\rho^{2-\alpha}$ if $\rho$
is small enough.
\end{example}

\mysection[Preliminary results]{Preliminary results} 

                     \label{section 10.25.2}

We use the following which 
combines particular cases of   Lemmas
\ref{lemma 4.23.3} and \ref{lemma 5.6.1}. 
We have two stopping times $\gamma\geq\tau$.

\begin{theorem}
                                        \label{theorem 9.27.1}

For any $\lambda\geq 0$ and  Borel $f(s,y),g(y)\geq0$  
\begin{equation}
                                   \label{5.6.401}
  E_{\cF_{\tau}}\int_{\tau}^{\gamma} e^{-\lambda(s-\tau)}  
f(  s,x_{ s})\,ds\leq N(d, \delta  )
\big(A_{\lambda} + B ^{2}_{\lambda}
\big)^{d/(2d+2 )}
\|f\|_{L_{d+1 }},
\end{equation}
\begin{equation}                      \label{3.7.1} 
  E_{\cF_{\tau}}\int_{\tau}^{\gamma}e^{-\lambda(s-\tau)}   
g(  x_{ s})\,ds\leq N(d, \delta  )
\big(A_{\lambda} + B_{\lambda}^{2}
\big)^{1/2}
\|g\|_{L_{d}(\bR^{d})},
\end{equation}
where  
$$
A_{\lambda}=E_{\cF_{\tau}}\int_{\tau}^{\gamma}e^{-\lambda(s-\tau)}\,ds,\quad
B_{\lambda}=E_{\cF_{\tau}}\int_{\tau}^{\gamma}e^{-\lambda(s-\tau)}|b_{ s}|\,ds.
$$

\end{theorem}
 
Observe that if $\gamma=\tau+\theta_{\tau}\tau_{ R}(y) $ 
in Theorem \ref{theorem 9.27.1}, 
then obviously $\gamma-\tau
\leq R^{2}$ and $A_{0} \leq R^{2}$. In that case also
$B_{0} \leq \bar b_{R}R$ by definition. Hence
we have the following.

\begin{lemma}
                                      \label{lemma 8.16.1}
 
For any Borel $f,g \geq0$, 
$\cF_{\tau}$-measurable $\bR^{d}$-valued $y$, and $R>0$ 
we have
\begin{equation}
                                          \label{9.29.2}
 E_{\cF_{\tau}}\int_{0}^{\theta_{\tau}\tau_{ R}(y)  }  
f( \tau+ s,x_{\tau+ s})\,ds\leq 
 N(d,  \delta)(1+\bar b_{R}) ^{d/  (d+1)  }
 R 
 ^{d/ (d+1)  }
\|f\|_{L_{d+1 }},
\end{equation}
\begin{equation}
                                            \label{8.22.2}
  E_{\cF_{\tau}}\int_{0}^{\theta_{\tau}\tau_{ R}(y)  }   
g(  x_{\tau+s})\,ds\leq N(d, \delta  )
(1+\bar b_{R})  
 R 
\|g\|_{L_{d}(\bR^{d})},
\end{equation}
 
\end{lemma}
 
Since $\theta_{\tau}\tau_{ R }(y)\leq R^{2}$,
we have
$E\theta_{\tau}\tau_{ R }(y)\leq R^{2}$, and this shows   that on average $\theta_{\tau}\tau_{ R }$ is of order
not more than $R^{2}$ for small $R$.
A very important fact which is implied by Corollary
\ref{corollary 7.29.1} is that 
$\theta_{\tau}\tau_{ R }$ is of order
not less than $R^{2}$.
To show this we need  the following result,
in which
\begin{equation}
                                                  \label{1.6.1}
\theta_{t}\gamma'_{ R}(x)  =\inf\{s\geq0:x_{t+s}  \in \bar B_{R}(x)  \}
\end{equation}
and which introduces one of the most important
conditions \eqref{12.18.3} in the book
under which $x_{t}$ behaves,
roughly speaking, as is there were no drift.

\begin{theorem}
                           \label{theorem 8.2.1}
There is a
   constant 
 $\sfb_{0}=\sfb_{0}(d,  \delta)>0$ such
\index{$S$@Miscelenea!$\bar b_{R}$@$\sfb_{0}$}%
that if, for a $\rho_{b}\in(0,\infty]$, 
we have
\begin{equation}
                                     \label{12.18.3}
  \bar b_{\rho_{b}}\leq \sfb_{0},
\end{equation}
then there is a constant 
 $\sfp_{0}=\sfp_{0}(d,\delta)\in (0,1) $ such 
\index{$S$@Miscelenea!$\sfp_{0}$}%
that for any $\rho\leq\rho_{b}$, any stopping time $\tau$
and $\cF_{\tau}$-measurable $\bR^{d}$-valued $y$, on the set $\{\tau<\infty\}$  we have
\begin{equation}
                                          \label{8.2.2} 
  P_{\cF_{\tau}}( \theta_{\tau}\tau'_{\rho}(y)   \geq   \rho^{2} )\leq 1-\sfp_{0},\quad
   P_{\cF_{\tau}}( \theta_{\tau}\tau'_{\rho}(x_{\tau})  \geq   \rho^{2} )\geq 
 \sfp_{0}  ;   
\end{equation}
moreover for $n=1,2,...$  
\begin{equation}
                                          \label{1.3.1} 
P_{\cF_{\tau}}( \theta_{\tau}\tau'_{\rho}(y)  > n\rho^{2})
 \leq (1-\sfp_{0})^{n},   
\end{equation}
so that 
\begin{equation}
                                          \label{3.7.2}
E _{\cF_{\tau}}\theta_{\tau}\tau'_{\rho}(y)  \leq N(d,\delta)\rho^{2},
\end{equation}
and
\begin{equation}
                                          \label{1.3.3}  
I:=E_{\cF_{\tau}}\int_{0}^{\theta_{\tau}\tau'_{\rho}(y)}|b_{\tau+s}|\,ds
\leq N(d ,\delta) \bar b_{\rho}\rho.
\end{equation}

Furthermore,   on the set
 $\{|y-x_{\tau}|\leq 9\rho/16,\tau<\infty\}$  
\begin{equation}
                                          \label{1.2.1} 
P_{\cF_{\tau}}\big(\theta_{\tau}\tau'_{\rho}(y) >\theta_{\tau}\gamma'_{\rho/16}(y)    \big)\geq\sfp_{0}.   
\end{equation}

\end{theorem}

 To prove the theorem we need an auxiliary result, in which
$$
m_{t,s}=-\int_{t}^{t+s}\sigma_{r}\,dw_{r},\quad a_{s}= 
\sigma_{s}\sigma_{s}^{*}.
$$

\begin{lemma}
                                       \label{lemma 1.2.1}
(i) There exists $\kappa=\kappa(d)>0$ such that,
for any $\rho>0$,
$$
\psi_{t}(s,y)=\rho^{-4}\big(\rho^{2}-4|y|^{2}\big)^{2}\phi_{t,s},\quad
\phi_{t,s}=\exp\int_{t}^{t+s}
\kappa \rho^{-2}\,\tr a_{r}\,dr,
$$
the process $\{\psi_{t}(s, m_{t,s}),\cF_{t+s}\}$ is a local
 submartingale.

(ii) Take a $\zeta\in C^{\infty}_{0}(\bR)$
such that it is even, nonnegative, and decreasing
on $(0,\infty)$.
For  $T\in(0,\infty)$ and $x\in \bR$ and $t\leq T$define
$u(t,x)=E\zeta(x+w^{1}_{T-t}) )$. Also 
take $t\geq 0$, $x\in\bR^{d}$ and set
$$
\xi_{t,s}=\frac{(x+ m_{t,s},a_{t+s}(x+ m_{t,s}) )}
{|x+m_{s}|^{2}}\quad (0/0:=1),\quad
\eta_{t,s}= \int_{0}^{s}\xi_{t,r}\,dr.
$$
Then   the process $\{u(\eta_{t,s},
|x+m_{t,s}|),\cF_{t+s}\}$
is a  supermartingale before $\eta_{t,s}$ reaches $T$,
in particular, on $[0,\delta^{2} T]$.

(iii) There exists $\alpha=\alpha(d,\gamma)>1$
such that for $u(x)=|x|^{-\alpha}$
and any $a\in\bS_{\gamma}$ we have
$$
a^{ij}D_{ij}u(x)\geq0,\quad x\ne0.
$$

\end{lemma}

Proof. (i) It is easy to see that for a $\kappa=\kappa(d)>0$
we have $\kappa\mu^{2}-8\mu+16\delta^{2}d^{-1}(1-\mu)\geq0
$ for all $\mu$, which implies that
for all $\lambda$
$$
\kappa(1-4\lambda^{2})^{2}-8(1-4\lambda^{2})
+64\delta^{2} d^{-1} \lambda^{2}\geq0.
$$

It follows that  (dropping $t$)
$$
\rho^{4}\phi_{s}^{-1}d\psi(s, m_{s} )=\kappa 
\big(\rho^{2}-4| m_{s}|^{2}\big)^{2}\rho^{-2}\tr a_{s}\,ds
$$
$$
-8\big(\rho^{2}-4|  m_{s}|^{2}\big)\big(2  m_{s}\,dm_{s}
+ \tr a_{s}\,ds\big)
+64( m_{s},a_{s}  m_{s})\,ds\geq dM_{s},
$$
where $M_{s}$ is a local martingale. This proves (i).

(ii) Observe that $u$ is smooth, even in $x$, and satisfies
$\partial_{t}u+(1/2)u''=0$.
Furthermore, as is easy to see $u'(t,x)\leq 0$
for $x\geq0$.
It follows by It\^o's formula that 
before $\eta_{t,s}$ reaches $T$ we have
(dropping obvious values of some arguments)
$$
du(\eta_{t,s},|x+m_{t,s}|)= \xi_{s} ( \partial_{t}u+(1/2)u'')
\,ds 
$$
$$
+\frac{u'}{2|x+m_{t,s}|} (\tr a_{t+s}-
 \xi_{s} )\,ds
+dM_{s},  
$$
where $M_{s}$ is a stochastic integral. Here the second
term with $ds$ is negative since $u'\leq0$, and this proves 
that $u(\eta_{t,s},|x+m_{t,s}|)$ is a local supermartingale.
Since it is nonnegative, it is a supermartingale.

Assertion (iii) is proved by simple computations.
The lemma is proved.\qed

{\bf Proof of Theorem \ref{theorem 8.2.1}}. 
Notice that for 
$$
\gamma:=\rho^{2}\wedge \inf\{s\geq 0:|m_{\tau,s} |\geq \rho/2\}
$$
we have
$
\phi_{\tau,\gamma}\leq e^{\kappa d/\delta}.    
$
Hence, by Lemma \ref{lemma 1.2.1} (i) 
$$
1=\psi(0,0) \leq  E_{\cF_{\tau}}\psi_{\tau}(\gamma,m_{\tau,\gamma})
\leq e^{\kappa d/\delta^{2}}
P_{\cF_{\tau}}(\sup_{s\leq \rho^{2}}|m_{\tau,s} |< \rho/2),
$$
$$
P_{\cF_{\tau}}(\sup_{s\leq \rho^{2}}|m_{\tau,s} |< \rho/2)
\geq 2\sfp_{0}(d,\delta)>0
$$
(the latter is used to predefine $\sfp_{0}$).
Also note that by Remark \ref{remark 10.14.1}
$$
P_{\cF_{\tau}}( \theta_{\tau}\tau'_{\rho}(x_{\tau})  < \rho^{2} )\leq
P_{\cF_{\tau}}\big(\int_{0}^{\theta_{\tau}\tau_{ \rho }(x_{\tau})}|b_{\tau+s}| \,ds
\geq \rho/2\big)
$$
$$+
P_{\cF_{\tau}}(\sup_{s\leq \rho^{2}}|m_{\tau,s} |\geq \rho/2)
$$
$$
\leq 2\bar b_{\rho}+1-P_{\cF_{\tau}}(\sup_{s\leq \rho^{2}}|m_{\tau,s} |< \rho/2)\leq 2\bar b_{\rho}+1-2\sfp_{0} 
$$
and we get the second relation  in \eqref{8.2.2}  for 
 $2\sfb_{0}\leq 
\sfp_{0}$.

To prove the first relation take $\zeta$ such that
$\zeta(z)=\eta(z/\rho)$, where $\eta(z)=1$
 for $|z|\leq 2 $ and take
$T=\delta^{2}\rho^{2} $, in which case $u(0,x)\leq u(0,0)<1$
and $u(0,0)$ depends only on $\delta$ (and $\eta$).
Also define $\mu$ as the first time $\eta_{\tau,s}$
reaches $T$, which is certainly less than or
equal to $\rho^{2}$. Now observe that
$u(\eta_{\tau,\mu},|y+m_{\tau,\mu}|)=u(T,|y+m_{\tau,\mu}|)=\zeta
(|y+m_{\tau,\mu}|)$. It follows that
 $$
 P_{\cF_{\tau}}(\sup_{s\leq R^{2}}|y+m_{\tau,s}|< 2\rho)
\leq P_{\cF_{t}}( |y+m_{\tau,\mu}|< 2\rho)
$$
$$
\leq E_{\cF_{t}}u(\eta_{\tau,\mu},|y+m_{\tau,\mu}|)\leq u(0,|y|)\leq u(0,0).
$$
Hence,
$$
P_{\cF_{\tau}}( \theta_{\tau}\tau'_{\rho} (y)  < \rho^{2} )\geq
 P_{\cF_{\tau}}\big(\int_{\tau}^{\tau+\theta_{\tau}\tau_{\rho} (y) }|b_{s}| \,ds
\leq \rho/2,\sup_{s\leq \rho^{2}}|y+m_{\tau,s}|\geq 2\rho \big)
$$
$$
\geq 1-P_{\cF_{\tau}}\big(\int_{\tau}^{\tau+\theta_{\tau}\tau_{\rho} (y)}|b_{s}| \,ds
\geq \rho/2\big)-P_{\cF_{\tau}}(\sup_{s\leq \rho^{2}}|y+m_{\tau,s}|\leq 2\rho)
$$
and it is clear how to adjust \eqref{12.18.3}
to get both inequalities in \eqref{8.2.2} with perhaps
$\sfp_{0}$ different from the above one.

To prove \eqref{1.3.1}   observe that
in light of \eqref{8.2.2}  for any $i=0,1,2...$
$$
P_{\cF_{\tau+i\rho^{2}}}(\max_{s\leq \rho^{2}}|x_{s+ \tau+i\rho^{2}}
-x_{ \tau+i\rho^{2}}+\xi_{i}|<\rho^{2})\leq 1-\sfp_{0},
$$
where $\xi_{i}=x_{_{\tau+i\rho^{2}}}-x_{\tau}-y $. In other words,
$$
P_{\cF_{\tau+i\rho^{2}}}(\max_{s\leq \rho^{2}}|x_{s+ \tau+i\rho^{2}}
-x_{\tau}-y|<\rho^{2})\leq 1-\sfp_{0}.
$$
Now \eqref{1.3.1} follows since its left hand side is
the conditional expectation given $\cF_{t}$
of the product of the above probabilities over $i=0,..,n-1$.

To prove \eqref{1.3.3} note that  
$$
I=\sum_{n=1}^{\infty}E_{\cF_{\tau}}I_{\theta_{\tau}\tau'_{\rho} (y)
>(n-1)\rho^{2} }E\Big\{
\int_{(n-1)\rho^{2} }^{(n\rho^{2})\wedge\theta_{\tau}\tau'_{\rho} (y)} |
b_{\tau+s}|\,ds\mid \cF_{\tau+(n-1)\rho^{2}  }\Big\}
$$
$$
\leq \bar b_{\rho} \rho\sum_{n=1}^{\infty}P_{\cF_{\tau}}( \theta_{\tau}\tau'_{\rho} (y)
>(n-1)\rho^{2} )
 \leq \bar b_{\rho}\rho \sum_{n=1}^{\infty}
(1-\sfp_{0})^{n-1}.
$$
This yields \eqref{1.3.3}.  

To prove \eqref{1.2.1} we may assume that $|y-x_{\tau}|>\rho/16$ and
using assertion (iii) of Lemma \ref{lemma 1.2.1}
 conclude that
$$
du(| x_{\tau+s}-x_{\tau}+y|)\geq b^{i}_{\tau+s}D_{i}u(|x_{\tau+s}-x_{\tau}+y|)\,ds+dM_{s},
$$
where $M_{s}$ is a local martingale
before $x_{\tau+s}-x_{\tau}+y$ hits the origin. 
For our $x$, on the time interval, which we denote
 $(0,\nu)$, when $ x_{\tau+s}-x_{\tau}+y
\in B_{\rho}\setminus \bar B_{\rho/16}$ we have
$|D u(|x_{\tau+s}-x_{\tau}+y|)| \leq N(d,\alpha)\rho^{-\alpha-1}$.  
Furthermore, at starting point $u(y)\geq (9\rho/16)^{-\alpha}$.
Consequently and by \eqref{1.3.3}
$$
(9\rho/16)^{-\alpha}\leq N\rho^{-\alpha-1}
E_{\cF_{\tau}}\int_{0}^{\theta_{\tau}\tau'_{ \rho} (y) }|b_{\tau+s}|\,ds
+P_{\cF_{\tau}}\big(\nu=\theta_{\tau}\tau'_{  \rho} (y)  \big)\rho^{-\alpha}
$$  
$$
+
P_{\cF_{\tau}}\big(\nu=\theta_{\tau}\gamma'_{ \rho/16}(x)  \big)(\rho/16)^{-\alpha},
$$
$$
(16/9)^{\alpha}\leq N \bar b_{\rho}
+1+(16^{\alpha}-1)
P_{\cF_{\tau}}\big(\theta_{\tau}\tau'_{\rho} (y) >\theta_{\tau}\gamma'_{\rho/16}(y)  \big) .
$$
It follows easily that \eqref{1.2.1} holds
with $\sfp_{0}$ perhaps different from the above ones,
once    \eqref{12.18.3} holds
with $2\sfb_{0}\leq\sfp_{0}$.
The theorem is proved.  \qed

 {\em We remind the reader that from this point on
throughout the chapter
we suppose that Assumption     \ref{assumption 8.19.2}  
 is satisfied. \/}

In light of \eqref{3.7.2} and \eqref{1.3.3} 
estimate \eqref{3.7.1} implies the following.

\begin{corollary}
                                     \label{corollary 3.7.1}
For any Borel $g\geq0$, $\rho\leq\rho_{b}$ and $\cF_{\tau}$-measurable $\bR^{d}$-valued $y$
we have
\begin{equation}
                                          \label{9.29.20}
 E_{\cF_{\tau}}\int_{0}^{\theta_{\tau}\tau'_{\rho}(y)  }  
g( x_{\tau+s})\,ds\leq 
 N(d,  \delta)(1+\bar b_{\rho})  
\rho 
\|g\|_{L_{d}(\bR^{d})}.
\end{equation}
\end{corollary}

From \eqref{1.3.1} we immediately obtain the following
  \begin{corollary}
                                  \label{corollary 2.3.1}
Let   $\rho\in(0,\rho_{b}] $. Then there
exists a constant $N$, depending only on
$\sfp_{0} $, such that, for any $\cF_{\tau}$-measurable $\bR^{d}$-valued $y$,
$T \geq 0$,
$$
P_{\cF_{\tau}}(\theta_{\tau}\tau'_{\rho}(y)   >T)\leq Ne^{-T/(N\rho^{2})}.
$$
\end{corollary}
\begin{corollary}
                                          \label{corollary 8.18.1}
For  $\mu \in[0,1]$,   stopping time $\tau$,  and any $\rho\in(0,\rho_{b}]$  we have
\begin{equation}
                                          \label{8.18.3}
E_{\cF_{\tau}}e^{-\mu \rho^{-2} \theta_{\tau}\tau_{ \rho}(x_{\tau})} \leq e^{- \mu\sfp_{0}/2}I_{\tau<\infty}.
\end{equation}
\end{corollary} 

 Indeed, for $\gamma:=\theta_{\tau}\tau_{\rho} 
(x_{\tau})$
on the set $\{\tau<\infty\}$
the derivative with respect to $\mu$ of the left-hand
side of \eqref{8.18.3} is
$$
-\rho^{-2} E_{\cF_{\tau}} \gamma e^{-\rho^{-2}\mu \gamma} \leq -
 e^{-\mu } 
P_{\cF_{\tau}}( \gamma\geq \rho^{2})\leq-e^{-\mu}\sfp_{0},
$$
where the last inequality follows from \eqref{8.2.2}.
By integrating we find
$$
E_{\cF_{\tau}}e^{-\mu \rho^{-2} \gamma} -1\leq  
(e^{-\mu   }-1)
\sfp_{0},
$$
which after using
$$
e^{-\mu  }-1\leq- \mu /2,\quad
1-\mu
\sfp_{0}/2\leq e^{-\mu 
\sfp_{0}/2}
$$
leads to \eqref{8.18.3}.
 
\begin{theorem}
                                     \label{theorem 8.20.1}
For any $\lambda\geq \rho_{b}^{-2} $, stopping 
time $\tau$, and $\rho\in(0,\infty)$ on 
the set $\{\tau<\infty\}$ we have
\begin{equation}
                           \label{8.20.1}
E_{\cF_{\tau}}e^{-\lambda  \theta_{\tau}
\tau_{ \rho}(x_{\tau}) }\leq
e^{\sfp_{0}/2}e^{- \sqrt{\lambda}  \rho
\sfp_{0}/2} .
\end{equation}
   
In particular, if $s\leq \rho \rho_{b}
\sfp_{0}/4 $, we have
\begin{equation}
                            \label{10.2.2}
P_{\cF_{\tau}}(  \theta_{\tau}\tau_{\rho}(x_{\tau}) \leq  s )\leq 
 e^{\sfp_{0}/2}\exp\Big(-\frac{{\sfp_{0}}^{2}\rho^{2}}{16 s}\Big).
\end{equation}
\end{theorem}

Proof. In case $\tau<\infty$ take an integer $n\geq 1$ and introduce $\tau^{k}$, $k=1,...,n$,
as the first exit time of $(\tau+s,x_{\tau+s} )$, $s\geq0$,
from $C_{\rho/n}( \tau+\tau^{k-1 } ,x_{\tau+\tau^{k-1}} )$ after $\tau^{k-1}$
($\tau^{0}:=0$). 
If
 $
\lambda\leq n^{2}/\rho^{2}$ and $ \rho/n\leq \rho_{b}$ 
then by \eqref{8.18.3} with $\mu=(\rho/n)^{2}
\lambda$ we have
$$
E_{\cF_{\tau^{k-1}}} e^{-\lambda( \tau^{k} - \tau^{k-1} )} \leq e^{- (\rho/n)^{2}
\lambda \sfp_{0}/2}.
$$
Hence,
\begin{equation}
                                           \label{12.27.1}
E_{\cF_{\tau}}e^{-\lambda \theta_{\tau}\tau_{\rho}(x_{\tau}) }\leq E_{\cF_{\tau}}\prod_{k=1}^{n}
e^{-\lambda( \tau^{k} - \tau^{k-1}  )}\leq  
 e^{- \rho^{2}n^{-1}
\lambda \sfp_{0}/2}.
\end{equation}
By taking
$n=\lceil \rho\sqrt\lambda\rceil$ and observing that that
$ \rho/n\leq  \rho_{b}$ and
$\rho^{2}n^{-1}
\lambda\geq \rho\sqrt\lambda-1$, we come to \eqref{8.20.1}.  

To prove \eqref{10.2.2} again consider the case that $\tau<\infty$ and
note that  for $\lambda\geq\rho_{b}^{-2}$ we have
$$
P_{\cF_{\tau}}( \theta_{\tau}\tau_{\rho}(x_{\tau})   \leq  s  )=P_{\cF_{\tau}}\big(
 \exp(-\lambda 
 \theta_{\tau}\tau_{\rho}(x_{\tau})  ) \geq \exp(-\lambda s)\big)
$$
$$
\leq  
 \exp(\sfp_{0}/2+ \lambda s-\sqrt{\lambda}\rho \sfp_{0}/2).
$$
For $\sqrt{\lambda}=\rho\sfp_{0}/(4s)$
 we get \eqref{10.2.2} provided that $\rho\sfp_{0}/(4s)\geq\rho_{b}^{-1}$.
The theorem is proved. \qed

\begin{corollary}
                                       \label{corollary 7.29.1}
Let   $\lambda>0$, $\rho\in (0,\rho_{b}]$.
Then there are constants $N=N(\sfp_{0}),\nu=\nu(\sfp_{0})>0$ such that on the set $\{\tau<\infty\}$ we have
\begin{equation}
                                                   \label{8.21.1}
N E_{\cF_{\tau}}\int_{0}^{\theta_{\tau}\tau_{\rho}(x_{\tau})}
e^{-\lambda t} \,dt  \geq \lambda^{-1}
(1-e^{-\lambda\nu \rho^{2}}) .
\end{equation}
In particular (as $\lambda\downarrow0$),
$NE_{\cF_{\tau}}\theta_{\tau}\tau_{\rho}(x_{\tau})\geq \nu \rho^{2}$.

\end{corollary}

Indeed, for any $\nu\leq   \sfp_{0}/4$
we have
$$
E_{\cF_{\tau}}\int_{0}^{\theta_{\tau}\tau_{\rho}(x_{\tau})}
e^{-\lambda t} \,dt=\lambda^{-1}
E_{\cF_{\tau}}(1-e^{-\lambda \theta_{\tau}\tau_{ \rho}(x_{\tau})})
$$
$$
\geq\lambda^{-1}
 E_{\cF_{\tau}}I_{ \theta_{\tau}\tau_{\rho}(x_{\tau}) >\nu \rho^{2}} (1-e^{-\lambda\nu \rho^{2}})
$$
$$
=\lambda^{-1} P_{\cF_{\tau}}( \theta_{\tau}\tau_{\rho}(x_{\tau}) >\nu \rho^{2})(1-e^{-
\lambda\nu \rho^{2}})
$$
$$
\geq \lambda^{-1}\Big(1-e^{\sfp_{0}/2} \exp
\Big(-\frac{\sfp_{0}^{2} }{16\nu }\Big)\Big)
(1-e^{-\lambda\nu \rho^{2}}),
$$
which yields \eqref{8.21.1} for an appropriate 
small 
$\nu =\nu( \sfp_{0} )>0$.

This result will be used in proving a higher summability
of the Green's functions of $x_{\cdot}$.
The next one is aimed at proving the precompactness
of distributions of various processes like $x_{\cdot}$.

\begin{corollary}
                           \label{corollary 10.26.1}
For any $n>0$ and
  $t\geq 0$ on the set $\{\tau<\infty\}$ we have
\begin{equation}
                                  \label{10.28.2}
E_{\cF_{\tau}}\sup_{r\in[0,t]}|x_{\tau+r}-x_{\tau}|^{ n}
\leq N(  t ^{ n/2}+ t ^{ n}),
\end{equation}
where $N=N(n, \rho_{b},\sfp_{0})$.
\end{corollary}

Indeed,  
for   $t\leq \mu \rho_{b}\sfp_{0}/4$
on the set $\{\tau<\infty\}$
we have
$$
P_{\cF_{\tau}}(\sup_{r\leq t }|x_{\tau+r}-x_{\tau}|\geq \mu)
\leq P_{\cF_{\tau}}(\theta_{\tau}\tau_{ \mu}\leq  t)
\leq e^{\sfp_{0}/2}\exp\Big(-\frac{  \mu^{2}\sfp_{0}^{2}}
{16 t}\Big).
$$
Consequently,
$$
E_{\cF_{\tau}}\sup_{r\leq t}|x_{\tau+r}-x_{\tau} |^{ n}
=n\int_{0}^{\infty}\mu^{n-1}
P_{\cF_{\tau}}(\sup_{r\leq t }|x_{\tau+r}-x_{\tau}|\geq \mu)\,d\mu
$$
$$
\leq n\int_{0}^{4t/(\rho_{b}\sfp_{0})}\mu^{n-1}\,d\mu
+ne^{\sfp_{0}/2}\int_{0}^{\infty}\mu^{n-1}
\exp\Big(-\frac{  \mu^{2}\sfp_{0}^{2}}
{16 t}\Big)
\,d\mu,
$$
and   the result follows. 

\begin{remark}
                          \label{remark 7.7.1}
If
Assumption \ref{assumption 8.19.2}  holds with any $\rho_{b}>0$, then 
the right-hand side of \eqref{10.28.2} becomes $Nt^{n/2}$.
\end{remark}

A few more general results are related to going through
  long ``tubes".
\begin{theorem} 
                                        \label{theorem 1.24.1}
Let $\rho\in(0,\rho_{b}]$, $\tau$ be a  stopping time, $\cF_{\tau}$-measurable $y_{\tau}\in\bR^{d}$ be such that $16|x_{\tau}-y_{\tau}|\geq 3\rho$ on the set $\{\tau<\infty\}$.  On the same set
for $r>0$ denote by $S_{r}(x_{\tau},y_{\tau})$   the open convex hull
of $B_{r}(x_{\tau})\cup B_{r}(y_{\tau})$. Then there exist
nonrandom
$T_{0},T_{1}$, depending only on $\sfp_{0}$,
such that $0<T_{0}<T_{1}<\infty$ and on the set $\{\tau<\infty\}$ the $P_{\cF_{\tau}}$-probability $\pi$
that $x_{\tau+s}$, $s\geq0$, will reach $\bar B_{\rho/16}(y_{\tau})$ before exiting
from $S_{\rho}(x_{\tau},y_{\tau})$ and this will happen
on the time interval $[nT_{0}\rho^{2},nT_{1}\rho^{2}]$
is greater than $\pi_{0}^{n}$, where
$$
n= \Big\lfloor \frac{16|x_{\tau}-y_{\tau}|+\rho}{4\rho}\Big\rfloor ,\quad \pi_{0}=\sfp_{0}/3.
$$

\end{theorem}

Proof.  We argue in case $\tau(\omega)<\infty$.
Introduce $\nu=\nu(x_{\tau},y_{\tau}) $ as the first time $ x_{\tau+s}$
reaches $\bar B_{\rho/16}(y_{\tau})$ and $\gamma=\gamma(x_{\tau},y_{\tau}) $ as the first time
it exits from $S_{\rho}(x_{\tau} ,y_{\tau})$. Owing to
$16|x_{\tau}-y_{\tau} |\geq 3\rho$, we have $n\geq1$ and we are going to use the induction
on $n$ with the induction hypothesis that,
for all $\rho\in(0,\rho_{b}]$,
$$
\Big\lfloor \frac{16|x_{\tau}-y_{\tau} |+\rho}{4\rho}\Big\rfloor=n
\Longrightarrow 
P_{\cF_{\tau}}(\gamma >\nu \in[nT_{0}\rho^{2},n
T_{1}\rho^{2}])\geq \pi^{n}_{0}.
$$

   If $n=1$, then $3\rho/16\leq |x_{\tau}-y_{\tau}|< 7\rho/16$ and by Theorem
\ref{theorem 8.2.1} (see \eqref{1.2.1}) we have $P_{\cF_{\tau}}(\theta_{\tau}\tau'_{\rho} >\nu )\geq\sfp_{0}$.
Furthermore, in light of  Theorem
\ref{theorem 8.2.1}, there is $T_{1}=T_{1}(\sfp_{0})$
such that $P_{\cF_{\tau}}(\theta_{\tau}\tau'_{\rho} >T_{1}\rho^{2})\leq \sfp_{0}/3$.
Using \eqref{10.2.2} we also see that there is
$T_{0}=T_{0}(\sfp_{0})<T_{1}$ such that
$P_{\cF_{\tau}}(\nu \leq T_{0}\rho^{2})\leq \sfp_{0}/3$.
It follows that $P_{\cF_{\tau}}(\gamma >\nu \in[T_{0}\rho^{2},
T_{1}\rho^{2}])\geq\sfp_{0}/3=\pi_{0}$.
This justifies the start of the induction.

Assuming that our hypothesis is true for some $n\geq 1$
 suppose that
$(n+2)\rho/4>|x_{\tau}-y_{\tau}|+\rho/16\geq  (n+1)\rho/4$. In that case, let 
$$
z_{\tau}=n\rho(x_{\tau}-y_{\tau})/(4|x_{\tau}-y_{\tau}|)
$$
and let $\nu'$ be the first time $ x_{\tau+s}$ reaches $\bar B_{\rho/16}(z_{\tau})$,
and let $\gamma'$ be the first time it exits
from $S_{\rho}(x_{\tau} ,z_{\tau})$. As is easy to see,
$$
P_{\cF_{\tau}}(\gamma >\nu\in[(n+1)T_{0}\rho^{2},(n+1)
T_{1}\rho^{2}])
$$
$$
\geq P_{\cF_{\tau}}\big(\gamma'>\nu'\in[T_{0}\rho^{2},T_{1}\rho^{2}],  
\gamma(x_{\nu'},y_{\tau}) >\nu(x_{\nu'},y_{\tau}) \in[nT_{0}\rho^{2},n
T_{1}\rho^{2}] \big)
$$
$$
=E_{\cF_{\tau}}I_{\gamma'>\nu'\in[T_{0}\rho^{2},T_{1}\rho^{2}]}
P_{\cF_{\nu'}} \big(\gamma(x_{\nu'},y_{\tau}) >\nu(x_{\nu'},y_{\tau})\in[nT_{0}\rho^{2},n
T_{1}\rho^{2}] \big).
$$
Observe that on the set $\nu'<\infty$ we have
$n\rho/4\leq |x_{\nu'}-y_{\tau}|+\rho/16<(n+1)\rho/4$, so that, by  our induction hypothesis, the last conditional probability
above is greater than $\pi_{0}^{n}$. Then   using the first part of the proof
we obtain our result for $n+1$ in place of $n$.
The theorem is proved. \qed

\begin{remark}
                                      \label{remark 1.25.1}
Notice that, for any fixed $x_{\tau},y_{\tau}$,  the time interval
$[nT_{0}\rho^{2},nT_{1}\rho^{2}]$ is as close to the origin
as we wish if we choose $\rho$ small enough.
Then, of course, the corresponding probability will be
quite small but $>0$.

\end{remark} 

\begin{corollary}
                                    \label{corollary 10.17.1}
Let $\rho\in(0,\infty)$,  $\tau$ be a  stopping time
and let $y$ be $\cF_{\tau}$-measurable $\bR^{d}$-valued. Then
there is a constant $N=N(d,\delta,\rho_{b},\rho)$ such that
on $\{\tau<\infty\}$ for all $T\in(0,\infty)$ we have
\begin{equation}
                                            \label{10.17.3}
P_{\cF_{\tau}}(\max_{t\leq T}|x_{\tau+t}-y|<\rho)
\leq Ne^{-T/N}.
\end{equation}

\end{corollary}

Indeed, on the set where $\tau<\infty$ and $|x_{\tau}-y|
\geq \rho$ estimate \eqref{10.17.3} is obvious. So we may
concentrate on $y\in B_{\rho}(x_{\tau})$ Then
by using 
Theorem \ref{theorem 1.24.1} we see that, given that
$\tau<\infty$ with    
 probability not less than some $\beta>0$,
  depending only on $d,\delta$, and $\rho_{b}$,
the process   $ x_{\tau+t}$
 will reach $\bar B_{\rho_{b}/16}(x_{\tau} +\rho_{b}e_{1}/4)$,
where $e_{1}$ is the first basis vector before time $T_{1}=T_{1}
(d,\delta)$. Therefore, its first coordinate will increase
by at least $3\rho_{b}/16$.   Repeating this argument
and taking into account that $\rho<\infty$, we see that
with probability
$\pi>0$ depending only on $d,\delta$, $\rho_{b}$, and $\rho $,
the process $ x_{\tau+t} $ will leave $ B_{\rho}(y)$ before time
$S$, where $S$ depends only on $d,\delta$, $\rho_{b}$, and $\rho $,  
that is
$$
P_{\cF_{\tau}}(\max_{t\leq S}|x_{\tau+t}-y|<\rho)\leq 1-\pi.
$$
Iterating this inequality, which is also true
if $|x_{\tau}-y|\geq \rho$, we obtain $
P_{\cF_{\tau}}(\max_{t\leq nS}|x_{\tau+t}-y|<\rho)\leq (1-\pi)^{n}$ for $n=1,2,...$ and this yields the result.

The following complements Corollaries \ref{corollary 2.3.1}
and \ref{corollary 10.17.1}.
 
\begin{corollary}
                                     \label{corollary 1.25.1}
Let   $\kappa\in[0,1)$, $\rho\in(0,\rho_{b}]$,  and let $\tau$ be a stopping time  $y$ be $\cF_{\tau}$-measurable $\bR^{d}$-valued.
Then for any $T>0$ on the set $\{\tau<\infty,
 |x_{\tau}-y|\leq\kappa \rho\}$
\begin{equation}
                         \label{1.25.2}
NP_{\cF_{\tau}}(\theta_{\tau}\tau' _{\rho} (y )  > T)\geq e^{-\nu T/[(1-\kappa)\rho]^{2}},
\end{equation}
where $N$ and $\nu>0$ depend only on $\sfp_{0}$.
\end{corollary}

Indeed, 
$\theta_{\tau}\tau' _{\rho} (y )\geq
\theta_{\tau}\tau' _{(1-\kappa)\rho} (x_{\tau} )$, which shows that  we may assume that $y=x_{\tau}$ and $\kappa=0$.
In that case, consider meandering of $x_{\tau+s}$
between $\bar B_{\rho/16}(x_{\tau})$ and $\partial B_{\rho/16}(y)$, where
$|y-x_{\tau}|=\rho/4$, without exiting from $B_{\rho}(x_{\tau})$. As is easy to deduce
from Theorem \ref{theorem 1.24.1},
given that the $n$th loop happened, with probability
$\pi^{4}_{0}$ the next loop will occur and take at least
$4\rho^{2}T_{0}$ of time. Thus  the $n$th loop
will happen and will take   at least $4n\rho^{2}T_{0}$ of time
with probability at least $\pi_{0}^{4n}$.
It follows that, for any $n$,
$$P_{\cF_{\tau}}(\theta_{\tau}\tau' _{\rho} 
(x_{\tau}) \geq 4n\rho^{2}T_{0}))\geq \pi_{0}^{4n},
$$
and this yields \eqref{1.25.2} for $y=x_{\tau}$ and $\kappa=0$.

Roughly speaking Corollary   
\ref{corollary 7.29.1} implies that the average time
  $x_{\tau+t}-x_{\tau}$ spends in
  $B_{\rho}$ before exiting from it is larger than a constant times $\rho^{2}$. Here the process $x_{\tau+t}-x_{\tau}$ starts from the center of $B_{\rho}$. An important fact,
leading to the so-called doubling property of
Green's functions, is that one can start not
from the center but even outside but not too far.

\begin{theorem}
                          \label{theorem 7.16.1}
Let   $\lambda\geq 0$, $\rho\in (0,\rho_{b}]$, and let $y$ be 
an $\cF_{\tau}$-measurable $\bR^{d}$-valued variable.
Then there are constants $N=N(\lambda, \sfp_{0}),\nu=\nu( \sfp_{0})>0$ such that on the set $A:=\{\tau<\infty,|x_{\tau}-y|\leq 2\rho\}$ we have
\begin{equation}
                                                   \label{7.16.2}
N E_{\cF_{\tau}}\int_{\tau}^{\infty}I_{B_{\rho}(y)}(x_{ t})
e^{-\lambda  (t-\tau) } \,dt  \geq \lambda^{-1}
(1-e^{-\lambda\nu \rho^{2}}) ,
\end{equation}
 where the right-hand side is understood
as $\nu\rho^{2}$ if $\lambda=0$.
\end{theorem}

Proof. By Theorem \ref{theorem 1.24.1} on $A$ with
$P_{\cF_{\tau}}$-probability $\pi=\pi(\sfp_{0})>0$ after time $\tau$ the process
$x_{t}$ will reach $B_{\rho/2}(y)$ and, if $\gamma$ is the first time this occurred, then $P_{\cF_{\tau}}(\gamma-\tau
\leq T )\geq \pi$, where $T=T(\sfp_{0})<\infty$. Then
$$
E_{\cF_{\tau}}\int_{\tau}^{\infty}I_{B_{\rho}(y)}(x_{ t})
e^{-\lambda  (t-\tau) } \,dt
$$
$$
\geq E_{\cF_{\tau}}e^{\lambda(\tau-\gamma)}\int_{\gamma}^{\infty}I_{B_{\rho/2} }(x_{ t}-x_{\gamma})
e^{-\lambda  (t-\gamma) } \,dt
$$
$$
\geq e^{-\lambda T }
E_{\cF_{\tau}}I_{\gamma-\tau
\leq T}E_{\cF_{\gamma}}
\int_{0}^{\infty}I_{B_{\rho/2} }(x_{\gamma+ t}-x_{\gamma})
e^{-\lambda t } \,dt.
$$
By Corollary   
\ref{corollary 7.29.1} the interior conditional expectation is larger than the right-hand side 
of \eqref{7.16.2} divided by $N(\sfp_{0})$ and then it only remains to
recall that $P_{\cF_{\tau}}(\gamma-\tau
\leq T)\geq \pi$. The theorem is proved. \qed

\mysection[On oblique cylinder]{An analog of the ``lemma on  oblique cylinder''}

In the previous section we concentrated on the behavior of the increment of $x_{t}$ after
time $\tau$. Here the main emphasis is
on the sets in $\bR^{d+1}$, the probability
$(t,x_{t})$ reaches these sets after $\tau$, and the time spent in these sets after $\tau$. 

 Our first big project is to prove a version of
Theorem 4.17 of 
\cite{Kr_19},
which provides an important step
toward establishing Harnack's inequality
for caloric functions. It is worth saying that
in the case of bounded $b$ Theorem \ref{theorem 12.7.2}
is proved by constructing a rather simple barrier,
see the PDE argument in the proof of Lemma 9.2.1
(``lemma on an oblique cylinder'')
of \cite{Kr_18} or the probabilistic argument
in the proof of Lemma 2.3 of \cite{Ya_20}.
In our case for the same purpose,
 we need a rather tedious argument like in
Theorem 4.17 of 
\cite{Kr_19}
just to get a good control of the {\em spatial\/}
process $x_{t}$.

Below $\sfp_{0}=\sfp_{0}(d,\delta)\in(0,1)$
is taken from Theorem \ref{theorem 8.2.1} (recall that Assumption
\ref{assumption 8.19.2}
 is supposed to hold throughout
this chapter).

 \begin{theorem}
             \label{theorem 12.7.2}
Let   $R\leq \rho_{b}$, $ \kappa,\eta\in(0,1)$,
$\tau$ be a stopping time, $(s,y)\in\bR^{d+1}$, $s>0$,
$z\in \bar B_{\kappa R}(y)$.
Then there exist $N,\nu>0$, depending only on 
$\kappa,\eta,\sfp_{0} $,  
 such that, for any
$\rho\in(0,1]$, on the set 
$$
A:=\{\eta R^{2}\leq s-\tau\leq  \eta^{-1}R^{2},
|x_{\tau}-y|< \kappa R\}
$$
we have
\begin{equation}
                                                 \label{12.9.1}
NP_{\cF_{\tau}}( |x_{s}-z|< \rho R  ,
\max_{[\tau,s]}|x_{t}-y| <R )\geq \rho^{\nu }.
\end{equation}
\end{theorem}

The proof of this theorem, given below after 
appropriate preparations,
 follows that of
Theorem 4.17 of 
\cite{Kr_19} and, roughly speaking, consists
of splitting the interval $[0,t]$ into several 
parts, estimating the probability that on the first part
the process will reach a neighborhood of $z$ without 
exiting from $B_{R}(y)$,
and then on the consecutive time intervals shrink
the neighborhood with constant coefficient
in such a way as to arrive at time $t$ in $B_{\rho R}(z)$
without exiting from $B_{R}(y)$.

Set $s'=s-(1/2)\eta R^{2} $.
The following lemma says that on $A$ on the time interval
$[\tau,s']$ with positive probability the process
$ x_{t}$ will reach in $B_{\rho_{0}\kappa R}(z)$
without exiting from $B_{R}(y)$,
where $\rho_{0}\in(0,1)$ is fixed. The idea of the proof is that, if $x_{\tau}$ is  close to $z$
in terms of $\rho_{b}$, then
with positive probability it will stay close on any fixed finite time interval. However, if it is far
from $z$, then the process can go through an 
appropriate tube in almost no time and reach
the neighborhood of $z$, to which the first case scenario applies.

\begin{lemma} 
                    \label{lemma 12.15.1}
Let  $R\leq \rho_{b}$ and
let
  $ \rho_{0}\in(0,1)$.
Then there exists
$\mu=\mu( \sfp_{0},  \kappa, \rho_{0},\eta)>0$ such that
on $A$
\begin{equation}
                                                 \label{12.7.3}
 P_{\cF_{\tau}}( |x_{s'}-z|<\rho_{0} \kappa  R ,
\max_{[\tau,s']}|x_{t}-y| <R  )\geq \mu.
\end{equation}

\end{lemma}

Proof. Observe that \eqref{12.7.3} becomes stronger if
$\rho_{0}$ becomes smaller. Therefore we may assume that
\begin{equation}
                               \label{12.7.4}
  \rho_{0}\leq \min\big(1/16,  \sqrt{\eta / 
(288T_{1} )},\kappa^{-1}-1  \big)  ,
\end{equation}
where $T_{1}=T_{1}(\sfp_{0})$ is taken from 
Theorem \ref{theorem 1.24.1}.
Then we split the proof into two cases.  

{\em Case 1: $|x_{\tau}-z|\leq 3\rho_{0}^{2}\kappa R$\/}.
By Corollary \ref{corollary 1.25.1},
applied to the ball of radius $R_{1}:=\rho_{0}\kappa R$
($\leq\rho_{b}$) and $x=z-x_{\tau} $,
after noticing that $|x_{\tau}-z|\leq (1/2)R_{1}$ we obtain
$$
NP_{\cF_{\tau}}\big(\max_{t\leq \eta^{-1} R^{2}}| x_{\tau+t}-z|<\rho_{0}\kappa R        \big )
$$
$$
\geq \exp(-\nu \eta^{-1} R^{2}/R_{1}^{2})
=\exp(- \nu \eta^{-1}  /[\rho_{0}\kappa ]^{2}).
$$
The probability here is obviously less
($\rho_{0}\kappa\leq 1-\kappa$) than the probability in
\eqref{12.7.3}  
and this proves \eqref{12.7.3} in the first case.

{\em Case 2: $|x_{\tau}-z|\geq 3\rho_{0}^{2}\kappa R$\/}.
Set
$R_{0}=16\rho_{0}^{2}\kappa R$ and note that $|x_{\tau}-y|+R_{0}
<\kappa R+(1-\kappa)16\rho_{0}R<R$.   Therefore, the sausage $S_{R_{0}}(x_{\tau},z )$,
defined as the open convex hull of $B_{R_{0}}(x_{\tau})\cup B_{R_{0}}(z)$,
 belongs
to $B_{R}(y)$. By Theorem \ref{theorem 1.24.1} with probability
not less than $\pi_{0}^{n}$ before time
$nT_{1}R_{0}^{2}$ the process $ x_{\tau+s}$
will hit $\bar B_{R_{0}/16}(z)$ without exiting
from $S_{R_{0}}(x_{\tau},z)$, where
$$
n\leq\frac{4|x_{\tau}-z|}{R_{0}}+\frac{1}{4}
\leq \frac{1}{2\rho^{2}_{0}\kappa}+\frac{1}{4}=:n_{0}.
$$
Since $R_{0}<  R$ and $|x_{\tau}-z|<  2R$ and also thanks to
$144T_{1}\rho_{0}^{2}\leq \eta/2$, we  have
$$  
nT_{1}R_{0}^{2}\leq T_{1}R_{0}(4|x_{\tau}-z|+ R_{0}/4)
\leq T_{1}R_{0}9   R=144 T_{1}\rho_{0}^{2}\kappa
R^{2}  \leq (\eta/2) R^{2}.
$$
By introducing $\gamma$ as the first time $ x_{\tau+t}$
hits $\bar B_{R_{0}/16}(y)$ we conclude that
\begin{equation}
                                               \label{1.27.1}
P_{\cF_{\tau}}(\max_{[0,\gamma]}|x_{\tau+t}-y|<R, \gamma\leq (\eta/2) R^{2})\geq \pi_{0}^{n}.
\end{equation}

Observe also that $R_{0}/16=\rho_{0}^{2}\kappa  R $  and at time $\gamma \leq (\eta/2) R^{2}$
the point $ x_{\tau+\gamma}$ is in $\bar B_{\rho_{0}^{2}\kappa  R}(z)$ and $\tau+\gamma\leq  s '$.
It follows from Case 1
that, given that $\tau+\gamma<\infty,\gamma \leq (\eta/2) R^{2}$,
 with probability $\pi_{1}>0$ depending only on $\sfp_{0}$,
$\rho_{0}$, $\kappa$, and $\eta$
the process $ x_{\tau+t}$ will stay in $B_{\rho_{0}\kappa R}(y)$
on the time interval $[\gamma,\gamma+\eta^{-1} R ^{2}]$.
Notice that $\gamma+\eta^{-1} R ^{2}\geq s'$.
Along with \eqref{1.27.1} this imply   \eqref{12.7.3} with $\mu=\pi_{0}^{n_{0}}\pi_{1}$.
The lemma is proved.  \qed

\begin{lemma} 
                    \label{lemma 12.15.10}
There are  constants $N,\nu>0$, depending only on
$\kappa,d,\delta $,  such that,  
for any $R\leq \rho_{b},\rho\in (0,1)  $, 
on the set
$$
 \{|x_{s'}-z|< \kappa (1-\kappa)R   \}
$$
we have
\begin{equation}
                            \label{12.14.2}
NP_{\cF_{s'}}\big( \max_{[s',s ]}|x_{ t}-y|<R ,
 |x_{s}-z|<\rho R  
\big) \geq \rho^{\nu }.
\end{equation}
\end{lemma}

Proof.  
  Set $\rho_{0}=1/2$,
observe that it suffices to prove
\eqref{12.14.2} for $\rho\leq\kappa $,
and define
$$
 n(\rho)=\Big\lfloor
\frac{\ln(\rho/\kappa)}{\ln \rho_{0} }\Big\rfloor+1\quad(\geq1),\quad \bar\eta
=\eta\frac{1-\rho_{0}}{1-\rho_{0}^{2n(\rho)}}.
$$

Note that $\bar\eta\in [\eta(1-\rho_{0}),\eta
(1+\rho_{0})^{-1}]$, 
so that by Lemma \ref{lemma 12.15.1} 
 estimate \eqref{12.7.3},
with constant $t$ in place of $\tau$
with $t\geq \bar \eta R^{2}$ and $z=y$, 
is valid and implies that  on the set 
$|x_{t}-z|< \kappa R $
\begin{equation} 
                         \label{12.16.10}
P_{\cF_{t}}( |x_{t+\bar\eta R^{2}/2}-z|<\rho_{0}\kappa R ,
\max_{[0,\bar\eta R^{2}/2]}|x_{ t+r}-z| <R  )\geq \mu,
\end{equation}     
whenever $R\leq \rho_{b}$. For $n=1,2,...$ 
 introduce  $R_{0}=(1-\kappa) R$
so that $|x_{s'}-z|\leq\kappa R_{0}$ and set
$$
R_{n}=\rho_{0}^{n-1}R_{0} ,\quad s_{n}= 
 \bar\eta  R_{n}^{2}/2=
 \bar\eta R_{0}^{2} \rho_{0}^{2(n-1)}/2,\quad t_{n}=\sum_{k=1}^{n}s_{k}.
$$
Then for each $n$  we get from \eqref{12.16.10} that 
on the set $\{|x_{s'+t_{n-1}}-z|< \kappa R_{n}\} $
($t_{0}:=0$) we have
\begin{equation}
                                \label{8.31.01}
 P_{\cF_{s'+t_{n-1}}}\big( |x_{s'+ t_{n}}-z|<  \kappa R_{n+1} , \max_{s\leq s_{n}}| x_{s'+t_{n-1}+t}-z|< R_{0} \big)\geq 
\mu.
\end{equation}
The conditional expectation given $\cF_{s'}$
of the product of the left-hand sides of \eqref{8.31.01}
over $n=1,...,n(\rho)$
is certainly less than
$$
P_{\cF_{s'}}\big( \max_{[0,t_{n(\rho)}]}|x_{s'+t}-y|<R ,
 |x_{s'+t_{n(\rho)}}-z|<\kappa R_{n(\rho)+1}  
\big) .
$$

Here  
$$
t_{n(\rho)}=\bar\eta R_{0}^{2}\frac{1-\rho_{0}^{2n(\rho)}}
{2(1-\rho_{0}) }=\eta R_{0}^{2}/2=s-s',\quad
\kappa R_{n(\rho)+1}= \kappa\rho_{0}
^{n(\rho)}R_{0}\leq \rho R_{0}.
$$
Therefore, the the probability in the left-hand side of \eqref{12.14.2}
is larger than $\mu^{n(\rho)}$.
Now to finish the proof,
it only remains to note that
$$
\mu^{n(\rho) }\geq
\mu \exp\Big(
\frac{\ln(\rho/\kappa)}{\ln \rho_{0} }\ln\mu\Big)
=N^{-1}\rho^{\nu}.
$$
The lemma is proved.  \qed     

{\bf Proof of Theorem \ref{theorem 12.7.2}}. 
Set $\rho_{0}=1-\kappa$.
In light of Lemmas \ref{lemma 12.15.1} and \ref{lemma 12.15.10} we have
on $A$ that
$$
P_{\cF_{\tau}}( |x_{s}-z|< \rho R  ,
\max_{[\tau,s]}|x_{t}-y| <R )\ 
$$
$$
\geq E_{\cF_{\tau}}I_{  |x_{s'}-z|<(1-\kappa)\kappa R ,
\max_{[\tau,s']}|x_{t}-y| <R  }
$$
$$
\times 
P_{\cF_{s'}}\big( \max_{[s',s ]}|x_{ t}-y|<R ,
 |x_{s}-z|<\rho R  
\big)\geq \mu N^{-1}\rho^{\nu}.
$$
  The theorem is proved. \qed

The following improvement of
Theorem \ref{theorem 12.7.2} consists
of allowing $z$ to be close
to the boundary of $B_{R}(y)$.

 \begin{theorem}
             \label{theorem 6.7,1}
Let $R\leq \rho_{b}$, $ \kappa,\eta,
  \in(0,1)$, $s>0,
z\in \bar B_{(1-\rho) R}(y)$.
Then the assertion of Theorem \ref{theorem 12.7.2} still holds. 

\end{theorem} 

Proof. We may assume that  $|z-y|>  \kappa R$, since the case of $|z-y|\leq\kappa R$ is taken care of by Theorem
\ref{theorem 12.7.2}. 
Also we may suppose that $\kappa\geq \rho$, so that $|z-y|>\rho R$.
Introduce 
$$
z_{0}=z+\rho R\frac{z-y}{|z-y|}.
$$
Since $ |z_{0}-y|>\kappa R$ and it
suffices to concentrate on $\rho
\leq \kappa/4$,
  we may assume that   for an integer
$m\geq2$,  we have 
\begin{equation}
                         \label{5.6,1}
|z_{0}-y|=2^{m}\rho R.
\end{equation}
Indeed, if \eqref{12.9.1} is true
for any $\rho=|z_{0}-y| 2^{-m}R^{-1} $
and actually $|z_{0}-y| 2^{-m-1}R^{-1}<\rho<|z_{0}-y| 2^{-m}R^{-1}$, then the left-hand side
of \eqref{12.9.1} will become
smaller if we replace there $\rho$
with $|z_{0}-y| 2^{-m-1}R^{-1}$ and after this replacement, by assumption,
it will be greater that $(|z_{0}-y| 2^{-m-1}R^{-1})^{\nu}=2^{-\nu}(|z_{0}-y| 2^{-m}R^{-1})^{\nu}\geq 2^{-\nu}\rho^{\nu} $.
This justifies \eqref{5.6,1}.

Then introduce the points $z_{k}$,
$k=1,...,m+1$, lying on the straight segment connecting the points
$y+(z_{0}-y) $ and $y-(z_{0}-y)$, by 
$$ 
|z_{0}- z_{k}|=2^{k  }\rho R=:2R_{k}.
$$
Also introduce   open balls $B^{k}$,
$k=1,...,m+1$,
by requiring the straight open
segment connecting $z_{k}$ and 
$z_{0}$ to be its diameter,
so that $B^{1}=B_{\rho}(z)$ and
$B^{m+1}=B_{|z_{0}-y|}(y)$.

Finally,  represent the interval $[s_{m+1},s_{1})$, where $s_{m+1}=s-(1/2)\eta R^{2}$, $
s_{1}=s$, as the union
of 
$$
 [s_{k+1},s_{ k }) , \quad k=1,...,m ,
$$
 such that
$s_{i}$ decreases and $s_{k}-s_{k+1}=(1/4) (s_{k+1}-s_{k+2})$ for $k= 1,...,m-2 $.
In fact $s_{k }-s_{k+1}=\alpha 4^{ k}$, where $\alpha$ is found from
$$
s_{1}-s_{m+1}=\alpha\sum_{k=1}^{m }4^{ k}\quad\Big(4^{m}<\sum_{k=1}^{m }4^{ k}<4\cdot4^{m}\Big) .
$$   

Observe that, since $2R>|z-y|>\kappa R$,
we have $2R>2^{m}\rho R>\kappa R$,
$2>2^{m}\rho>\kappa$. Also $s_{1}-s_{m+1}=(1/2)\eta R^{2}$ which implies certain estimates on $\alpha$ and shows that
$$
 (1/8)\eta R^{2}_{k+1} \rho^{-2}4^{-m}\leq s_{k}-s_{k+1}\leq (1/2)\eta R^{2}_{k+1}\rho^{-2}4^{-m}.
$$
Here
$$
 1/4<\rho^{-2}4^{-m}<\kappa^{-2}.
$$
It follows that
$$
\beta R_{k+1}^{2}\leq s_{k}-s_{k+1}\leq \gamma R_{k+1}^{2},
$$
where $\beta,\gamma$ depend only on
$\eta,\kappa$.

By using an argument very similar
to the one from the proof of
Lemma \ref{lemma 12.15.1} one shows that, if $r\leq \rho_{b}$, balls
$B^{'}\in \bB_{r},B^{''}\in\bB_{r/2}$,
$B^{''}\subset B^{'}$ and $B^{'},B^{''}$
have a common boundary point, then,
given $x_{t_{1}}\in (1/2)B'$ and $\beta r^{2}\leq t_{2}-t_{1}\leq \gamma r^{2}$, the $P_{\cF_{t_{1}}}$-probability that,   $x_{t_{2}}\in (1/2)B''$ and $x_{t}\in B'$ for all
$t\in[t_{1},t_{2}]$, is bigger that
$\nu=\nu(d,\delta,\sfp_{0},\eta)>0$.

By applying this fact to $t_{1}=s_{k+1},
t_{2}=s_{k},r=R_{k+1}, B'=B^{k+1},
B''=B^{k}$ and doing so for $k=m,m-1,...,2$ we conclude that,
given that $|x_{s_{m+1}}-y|<(1/2)|z_{0}-y|$,
the $P_{\cF_{s_{m+1}}}$-probability
that $x_{s}\in B_{\rho}(z)$ and is not exiting before $s$ from $B_{R}(y)$ is bigger
than $\nu^{m-1}\geq \varepsilon \rho
^{-\ln_{2}\nu}$, where $\varepsilon
=\varepsilon(\kappa)>0$.
 
After that it only remains to recall
that $ |z_{0}-y|>\kappa R$ and to add
that
by  Lemma \ref{lemma 12.15.1}
$$
P_{\cF_{\tau}}( |x_{s_{m+1}}-y|<(1/2)
|z_{0}-y|  ,
\max_{[\tau,s_{m+1}]}|x_{t}-y| <R  )\geq \mu.
$$
The theorem is proved. \qed

\mysection[Time spent in space-time sets 
of small measure]{Estimating time spent in space-time sets 
of small measure}
                                  \label{section 12.29.2}

The central result of this section is Theorem \ref{theorem 12.21.1}
which needs some auxiliary constructions and results.

We present 
extensions to the case that $b\in L_{d+1}$  of
 probabilistic versions of some
PDE results found in \cite{KS_80}, \cite{Sa_80},  
\cite{Kr_18}. Recall the notation introduced in
the Preface and also
introduce
$$
C^{o}_{T,R}=(0,T)\times B_{R},\quad C^{o}_{T,R}(t,x)=(t,x)+C^{o}_{T,R},
\quad C^{o}_{R}(t,x)=C^{o}_{R^{2},R}(t,x),
  $$
 $C^{o}_{R}=C^{o}_{R}(0,0)$, which are open sets.
Fix 
 $$
q, \eta ,  \kappa\in(0,1).
$$

 For cylinders $Q= C^{o}_{ \rho}(t ,x ) $
 define       
$$
Q'=(t ,x )-
C^{o}_{\eta^{-1} \rho^{2},\rho},\quad Q''=\big
(t -\eta^{-1}\rho^{2},x \big)
+C^{o}_{\eta^{-1}\rho^{2}\kappa^{2},\rho\kappa},
 $$
$$
Q'_{+}=Q\cup Q'\cup(\{t\}\times B_{\rho}(x)).
$$

Imagine that the $t$-axis is pointed up vertically. Then 
 $Q'$ is   adjacent to $Q$ from below, the 
two cylinders have a common base, and along the $t$-axis
$Q'$ is $\eta^{-1}$ times longer than $Q$.  
 The cylinder
$Q''$ is obtained by  contracting $Q'$  
to    the
center  of  its  lower base   with the contraction factor  
$\kappa^{-2}$ for the $t$-axis and $\kappa^{-1}$ for the spatial
axes. 

\begin{remark}
                                                   \label{remark 2,14,1}
If  $Q= C^{o}_{ \rho}(t ,x )$, then
  the   distance
between $Q$ and $Q''$ along the $t$ axis is     
\begin{equation}
                                                          \label{2,14,1}
\eta^{-1}\rho^{2}-\eta^{-1}\rho^{2}\kappa^{2}=\eta^{-1}\rho^{2}
(1-\kappa^{2}),
\end{equation}
 which is $\rho^{2}$ if $\eta=1-\kappa^{2}$.

\end{remark}

 Let $\Gamma$ be a measurable subset of $C_{1}$
and
 introduce   $\mathcal{B}=\cB(\Gamma,q)$ 
 as the family  of   
{\em open\/}
cylinders $Q$ of  type 
$ C^{o}_{ \rho}(t_{0},x_{0})$ such that   
$$
Q\subset  C_{1 } \quad\text{and}\quad
|Q\cap\Gamma|\ge q|Q|.
$$ 

Finally,  define 
$$
\Gamma''=\bigcup_{Q\in\mathcal{B}}Q'',\quad
\Gamma''_{\varepsilon}=\bigcup_{Q\in\mathcal{B}:|Q|\geq
\varepsilon}Q''.
 $$

 Observe that for $Q\in\mathcal{B}$ the set 
   $Q''$ is open. Hence, 
$\Gamma''$ and $\Gamma''_{\varepsilon}$  are open and measurable. 
\begin{lemma}
                                                     \label{lem:4.1.6} 
If $|\Gamma|\le q|C_{1 }|$, then   (i)
$$|\Gamma''|\ge\Big(1-\frac{1-q}{3^{d+1}}\Big)^{-1}
(1+\eta)^{-1}\kappa^{d+2}|\Gamma|
 $$
and (ii) there exist $\eta=\eta_{ 0}(d,q)\in(0,1)$, $\kappa=\kappa_{ 0}(d,q)\in(0,1)$ and $\vartheta=\vartheta_{ 0}(d,q )>1$
such that for any sufficiently small $\varepsilon>0$
there exists a closed $\Gamma_{\varepsilon}
\subset \Gamma''_{\varepsilon}$ such that
\begin{equation}
                                                 \label{12.21.3}
|\Gamma_{\varepsilon}|\geq\vartheta |\Gamma|.
\end{equation}
\end{lemma}

The first assertion of the lemma originated
in \cite{KS_80}, \cite{Sa_80},  is presented, for instance
as Lemma 9.3.6 in
\cite{Kr_18}. The second one is proved
in the same way as the second assertion
of Lemma 4.8 of \cite{Kr_21}.

\begin{lemma}
                            \label{lemma 12.20.2}
Let $\kappa\in(0,1)$. Then
there is a constant 
$q_{0}=q_{0}(d,\delta,\kappa )
\in(0,1)$ such that for any $\rho\leq  \rho_{b}$,
  $y\in\bR^{d},s\geq0$, 
on the set 
$$
A:=\{ x_{s}\in \bar B_{\kappa \rho }(y) \}
$$
 for any Borel set $\Gamma\subset C_{\rho}(s,y) $
satisfying $|\Gamma|\geq q_{0}|C_{\rho}|$, 
we have
\begin{equation}
                                    \label{12.20.3}
E_{\cF_{s}}\int_{s}^{s+\theta_{s}\tau_{\rho}
( y)   }I_{\Gamma}(t , x_{t})
\,dt\geq\mu_{1}  E _{\cF_{s}}\theta_{s}\tau_{\rho}
(y) \geq \mu_{0} \rho^{2} ,
\end{equation}
where   $\mu_{i}=\mu_{i}(d,\delta,\kappa )\in(0,1)$.

\end{lemma}

Proof.  Note that
in light of
 Corollary \ref{corollary 7.29.1} on $A$ we have
$$
  E _{\cF_{s}}\theta_{s}\tau_{\rho}
(y)   \geq 
E _{\cF_{s}}\theta_{s}\tau_{(1-\kappa)\rho}(x_{s})  \geq\nu(\sfp_{0},\kappa)  \rho^{2}.
$$
Next, by using Theorem \ref{theorem 6.3,1} we get that 
$$
E _{\cF_{s}}\theta_{s}\tau_{\rho}
(y) -
E_{\cF_{s}}\int_{s}^{s+\theta_{s}\tau_{\rho}
(y)   }I_{\Gamma}(t , x_{t})
\,dt    
$$
$$
=E_{\cF_{s}} \int_{s}^{s+\theta_{s}\tau_{\rho}
(y)   }I_{C_{\rho}(s,y)
\setminus \Gamma}(t, x_{t})
\,dt
$$
$$
\leq N\rho^{2}\dashnorm I_{C_{\rho}(s,y)\setminus\Gamma}(\cdot+s,\cdot +y)\|_{L_{d +1}(C_{\rho})}
=NR^{2}\dashnorm I_{C_{\rho}(s,y)\setminus\Gamma} \|_{L_{d +1}(C_{\rho}(s,y))}
$$
$$
\leq N\rho^{2}(1-q_{0})^{1/(d +1)}\leq N(1-q_{0})^{1/(d +1)}
 E _{\cF_{s}}\theta_{s}\tau_{\rho}
(y),
$$
where the constants $N$ depend only on $ d,\delta,\kappa$. We see how to choose
$q_{0}$ to get the desired result.
The lemma is proved.  \qed

In Lemma \ref{lemma 12.20.30} by $q_{0}$ we mean the one from Lemma \ref{lemma 12.20.2}.

\begin{lemma}
                               \label{lemma 12.20.30}
Take   $Q= C^{o}_{ \rho}(s,y)$
with $\rho\leq  \rho_{b}$,
use the notation $Q',Q'',Q'_{+}$ introduced above, set $\cO=Q'_{+}$,
assume that $\eta=1-\kappa^{2}$,
and suppose that Borel $\Gamma\subset Q$ is such that
$|\Gamma|\geq q_{0}|Q|$. Then there is a constant
$\nu_{0}>0$, depending only on $ \kappa$,
$ d,\delta$,  
such that on  $A:=\{\tau<\infty,(\tau ,x_{\tau} )\in Q''\}$ we have
\begin{equation}
                                                    \label{12.20.4}
E_{\cF_{\tau}} \int_{\tau}^{ \tau
+\theta_{\tau}\tau_{ \cO}} I_{\Gamma}(t, x_{t})
\,ds\geq \nu_{0}E _{\cF_{\tau}}\theta_{\tau}\tau_{ \cO}
\end{equation}
($\theta_{\tau}\tau_{ \cO}$ on $A$ is the   time spent by $(t, x_{t})$
in $Q'_{+}$ after $\tau$).
\end{lemma}

Proof. Thanks to   Remark \ref{remark 2,14,1}
we have $s-\tau \in(\rho^{2}, \eta^{-1}\rho^{2})$.
Also $|y-x_{\tau} |<\kappa\rho$ on $A$. It follows by
 Theorem \ref{theorem 12.7.2} that
$$
P _{\cF_{\tau} }\big(\sup_{r\in[\tau,s ]}| x_{r}-y|< \rho,
| x_{s}-y|<  \rho/2\big)\geq \nu ,
$$
where $\nu=\nu(\kappa,\sfp_{0})>0$.

Next,   in light of Lemma \ref{lemma 12.20.2}, we have on $A$
$$
E_{\cF_{\tau}} \int_{\tau}^{ \tau
+\theta_{\tau}\tau_{ \cO}}I_{\Gamma}(t, x_{t})
\,dt 
$$ 
$$
\geq E_{\cF_{\tau}}I_{\sup_{r\in[\tau,s ]}| x_{s}-y|< \rho,
| x_{s}-y|<  \rho/2}
E_{\cF_{s}}\int_{s}^{ s+\theta_{s}\tau_{\rho}
( y) }I_{\Gamma}(t, x_{t})
\,dt 
$$
$$
\geq \mu_{0}\rho^{2}
P _{\cF_{\tau}}\big(\sup_{r\in[\tau,s ]}| x_{s}-y|< \rho,
| x_{s}-y|<  \rho/2\big)\geq \mu_{0}\nu\rho^{2}.
$$
On the other hand, the height of $Q'_{+}$ is $(1+\eta^{-1})
\rho^{2}$, so that $(t, x_{t})$ cannot spend
in $Q'_{+}$ more time than $(1+\eta^{-1})
\rho^{2}$. This proves the lemma.  \qed

\begin{theorem}
                                             \label{theorem 12.21.1}
For
 any $\kappa,\eta\in(0,1)$ there exist  $\gamma\in(0,1)$
and $N$,
depending only on $\kappa,d,\delta$ with $N$
also depending on $\eta$,  such that
for any $R\leq \rho_{b}$, $q\in(0,1)$, $y\in\bR^{d},s\geq R^{2}$,  Borel
$\Gamma\subset C_{R}(s ,y)$ satisfying
$|\Gamma|\geq q |C_{R} |$, 
on the set $A:=\{s-\eta^{-1}R^{2}\leq\tau\leq s-R^{2},x_{\tau}\in B_{\kappa R}(y)\}$
we have

\begin{equation}   
                                                   \label{10.1.10}
  E_{\cF_{\tau}}\int_{\tau }^{\tau
+\theta_{\tau}\tau'_{  R}(y) }
I_{\Gamma}(t, x_{t})\,dt\geq N^{-1}q^{ \gamma}R^{2}.
 \end{equation}

\end{theorem}

Proof. Set $s'=s-R^{2}$ and observe that on $A$
$$
E_{\cF_{\tau}}\int_{\tau }^{ \tau
+\theta_{\tau}\tau'_{  R}(y)}
I_{\Gamma}(t, x_{t})\,dt
$$
$$
\geq E_{\cF_{\tau}}I_{|x_{s' }-y|< \kappa R  ,
\max_{[\tau,s']}|x_{t}-y| <R }E_{\cF_{s'}}
\int_{s'}^{ s'+\theta_{s'}\tau'_{ R}(y)}I_{\Gamma}(t, x_{t})\,dt
$$
and in light of Theorem \ref{theorem 12.7.2} (with $\rho=\kappa$)  
$$
N(d,\delta,\kappa,\eta)P_{\cF_{\tau}}(|x_{s'}-y|< \kappa R  ,
\max_{[\tau,s']}|x_{t}-y| <R )\geq 1.
$$
It follows that to prove the theorem, it suffices to show that \eqref{10.1.10} holds for $\tau=s'$ with
$\gamma\in(0,1)$ and $N$ depending only on 
$\kappa,d,\delta$.

Fix $R\leq  \rho_{b}$, $y\in\bR^{d}$
and denote
$$
G _{s}(\Gamma ):= E_{\cF_{s'}}
\int_{s'}^{ s'+\theta_{s'}\tau'_{ R}(y)}I_{\Gamma}(t, x_{t})\,dt.
$$
Also fix 
$\kappa_{0}\in[\kappa,1)$,    such that for $\eta_{0}=1-\kappa_{0}^{2}$ 
and $q=q_{0}$ the factor of $|\Gamma|$ in 
Lemma \ref{lem:4.1.6} is strictly bigger than one and
take $\theta=\theta(d,q_{0},1-\kappa_{0}^{2},\kappa_{0})>1$
from that lemma. 
Now we set ourselves the problem of finding
the largest (nonrandom) $\mu_{R}(q)$ such that
$$
R\leq  \rho_{b},\,\,
R^{2}\leq s ,x_{s-R^{2}}\in B_{\kappa R}(y),\,\,
|\Gamma|\geq q |C_{R}|,\,\,
\Gamma\subset C_{R}(s,y)
$$
\begin{equation}    
                           \label{2.12.1}
 \Longrightarrow G_{s}(\Gamma)\geq 
\mu_{R}(q)R^{2}.
\end{equation}
Observe that such $\mu_{R}(q)$ obviously exists  and the assertion of the theorem (with $\tau=s$) now says that
$$
\mu_{R}(q)\geq N^{-1}q^ {\gamma} .
$$

Notice that a combination of Lemma \ref{lemma 12.20.2}
and Theorem \ref{theorem 12.7.2}  leads to the conclusion that
there exist  $ \mu_{0},q_{0}\in(0,1)$, depending only on
$ \kappa,d,\delta $,  
 such that for $q\in [q_{0} ,1]$   
$$    
\mu_{R}(q)\geq \mu_{0} .
$$

We will be comparing $\mu_{R}(q')$ and $\mu_{R}(q'')$
for $ 0 <q'<q''<1$ such that
\begin{equation}
                                      \label{1.5.1}
(1+\theta) q'\geq 2q''.
\end{equation}

In case $R^{2}\leq s ,x_{s'}\in B_{\kappa R}(y)$ take a Borel
$\Gamma\subset C_{R}(s,y)$ satisfying
$|\Gamma|\geq q' |C_{R} |$ and
in the construction before Lemma \ref{lem:4.1.6}  replace $C_{1}$
by $C_{R}(s,y)$, keep all other notation, and
from the chosen $\Gamma,\kappa,\eta=1-\kappa_{0}^{2}$, and $q_{0}$ (not $q'$)   
 build up  the closed
set  $\Gamma _{\varepsilon}$ and take
$\varepsilon$ so small that \eqref{12.21.3} holds.
  There are  two cases:  

(i) $\big|\Gamma_{\varepsilon}\setminus 
C_{R}(s,y)\big|\leq (q''-q')|C_{R} |$,

(ii) $\big|\Gamma_{\varepsilon}\setminus
C_{R}(s,y)\big|> (q''-q')|C_{R} |$.

{\em Case  (i )\/}. Our  goal is to show that
\begin{equation}
                                                 \label{7,26,5}
G_{\tau} (\Gamma )
\geq\min\big(\mu_{0}, \nu_{0}\mu_{R}(q'')\big)R^{2},
\end{equation}
where $\nu_{0}$
  depends only on  $\kappa, d, \delta $.

Observe that, if $|\Gamma|\geq q_{0}|C_{R }|$, by definition
$G_{s} (\Gamma )\geq \mu_{R}(q_{0})R^{2}\geq\mu_{0}R^{2}$.
 Hence, we may assume that
$$
|\Gamma|< q_{0}|C_{R }|.
$$
In that case define
$$
\hat\Gamma_{\varepsilon}=\Gamma_{\varepsilon}\cap C_{R }(s,y).
 $$

Notice that by definition and  Lemma \ref{lem:4.1.6}   
$$
q'|C_{R }|\leq|\Gamma|\leq \theta^{-1}|\Gamma_{\varepsilon}|.
 $$
Moreover, by assumption   
$$
|\Gamma_{\varepsilon}|=\big|\Gamma_{\varepsilon}
\setminus C_{R }(s,y)\big|+|\hat  \Gamma_{\varepsilon}|
\leq (q''-q')|C_{R }|+|\hat  \Gamma_{\varepsilon}|.
 $$
Due to \eqref{1.5.1}, it follows that
$$
|\hat \Gamma_{\varepsilon}|\geq q''|C_{R }|,
$$
so that      
$$
G_{s} (\hat \Gamma_{\varepsilon})\geq\mu_{R}(q'')R^{2} .   
 $$

We now estimate $G_{s} ( \Gamma )$ from
below by means of $G_{s} (\hat \Gamma_{\varepsilon})$
using Lemma \ref{lemma 12.20.30}. Since $\Gamma_{\varepsilon}
\subset \Gamma''_{\varepsilon}$,
the closed set $\Gamma_{\varepsilon}$ is 
covered by the family $\{Q'':Q\in\mathcal{B},|Q|\geq
\varepsilon\}$. Then there is
finitely many $Q(1),...,Q(n)\in \mathcal{B}$ such that
$|Q(i)|\geq\varepsilon$, $i=1,...,n$, and
$$
\Gamma_{\varepsilon}
\subset \bigcup_{i=1}^{n}Q''(i)
=:\Pi_{\varepsilon}.
$$

Now for $(t,x)\in \Pi_{\varepsilon}$ define  
$i (t,x) $
as the first $i\in\{1,...,n\}$ for which
$(t,x)\in Q''(i)$. 
Also set $Q '_{+} (0):=C_{2R^{2},R}(s' ,y)$
and $i(t,x)=0$ if $(t,x)\in\partial' C_{2R^{2},R}(s,y)$.
Then note that $(s',x_{s'})\in Q '_{+} (0)$ and define recursively $\gamma^{0}=s'$,
$\tau^{1}$ as the first time after $\gamma^{0}$ when $(t, x_{t})$ exits
from $Q '_{+} (0)   \setminus \Gamma _{\varepsilon}$,
$\gamma^{1}$ as the first  time after $\tau^{1}$
when $(t, x_{t})$ exits from $Q '_{+}  (i(
 \tau^{1} , x_{\tau^{1}}))$,
and generally, for $k=2,3,...$ define
$\tau^{k}$ as the first time after $\gamma^{k-1}$ when $(t, x_{t})$  exits
from $Q'_{+}(0) \setminus \Gamma _{\varepsilon}$,
$\gamma^{k}$ as the first  time after $\tau^{k}$
when $(t, x_{t})$ exits from $Q'_{+}(i( \tau^{k}, x_{\tau^{k}}))$.
It is easy to check that so defined
$\tau^{k}$ and $\gamma^{k}$ are stopping times
and, since $|Q(i)|\geq\varepsilon$ and the trajectories 
of $(t, x_{t})$  are continuous,  
$\tau^{k}\uparrow s'+(2R^{2})\wedge\theta_{s'}\tau'_{ R}(y)$ as $k\to\infty$.
Furthermore,
($A$-a.s.) all the $\tau^{k}$'s equal $s'+(2R^{2})\wedge\theta_{s'}\tau'_{ R}(y)$
 for all large $k$.

On $A$ we have   
$$
G_{s} (\Gamma)\geq\sum_{k=1}^{\infty}E_{\cF_{s'}}\int_{\tau^{k}}
^{\gamma^{k}}I_{\Gamma}(t, x_{t})\,dt
$$
$$
=\sum_{k=1}^{\infty} E_{\cF_{s'}} E_{\cF_{\tau_{k}}} \int_{\tau^{k}}
^{\gamma^{k}}I_{\Gamma\cap Q'_{+}(i( \tau_{k} ,x_{\tau_{k}}))}(t,x_{t})\,dt.
$$
We estimate the interior expectation from below
by Lemma \ref{lemma 12.20.30} and get that $G_{s} (\Gamma )/\nu_{0}$ ($\nu_{0}=\nu_{0}(\kappa,d,\delta)$)
is greater than or equal to  
$$
\sum_{k=1}^{\infty} E_{\cF_{s'}}(
\gamma^{k}-\tau^{k} )
\geq \sum_{k=1}^{\infty} E_{\cF_{s'}}
\int_{\tau^{k}}^{\gamma^{k}}I_{\Gamma_{\varepsilon}}
(t,x_{t})\,dt
$$
$$=E_{\cF_{s'}}\int_{s'}^{ s'+\theta_{s'}\tau'_{ R}(y)}
I_{\Gamma_{\varepsilon}}
(t,x_{t})\,dt
$$
$$
=G_{s}(\Gamma_{\varepsilon})\geq
G_{s}(\hat \Gamma_{\varepsilon}) \geq \mu (q'')R^{2}.
$$

This proves \eqref{7,26,5}.

{\em Case (ii)\/}. Here the goal is to prove that  
\begin{equation}
                                                   \label{7,27,1}
G_{s} (\Gamma)\geq \mu_{0}\xi\eta^{n}(q''-q')^{n}R^{2},
 \end{equation}
  where   $\xi>0,n\geq1$ depend 
only on $d,\delta $,   and $\kappa$.

First we make a simple observation that
for some $(t,x)\in\Gamma_{\varepsilon}$
it holds that $t<(q'-q'')R^{2}+s$. It follows that there is a cylinder
$$
Q=C^{o}_{ \rho}(r,z)\in\cB
$$
 such that $Q'$
contains points in the half-space $t<(q'-q'')R^{2}+s$. 
Since $q'<q''$ and   $Q'$ is adjacent
to $Q\subset C_{R }(s,y)$, this implies that the height of $Q'$
is at least $(q''-q')R^{2}$, that is,  
\begin{equation}
                                            \label{12.22.7}
\rho^{2}\eta^{-1}\geq (q''-q')R^{2},\quad\rho^{2}\geq
\eta(q''-q')R^{2}.
\end{equation}
On the other hand, $Q\subset C_{R}(s,y)$, 
  $\rho< R$, and $s\leq r< s+R^{2} $    

Moreover, 
 by construction,  $|\Gamma\cap Q|\geq q_{0}|Q|$ and
by Lemma \ref{lemma 12.20.2} on the set where
$|x_{r}-z|\leq\rho/2$        
$$
I :=E_{\cF_{r}} \int_{r}^{r+\theta_{r}
\tau_{\rho}(z)}I_{\Gamma}(t, x_{t})
\,dt 
\geq \mu_{0}\rho^{2}\geq \mu_{0}\eta(q''-q')R^{2}.
$$
   Now by Theorem \ref{theorem 6.7,1}
$$
 E_{\cF_{s'}}\int_{s'}^{s'+\theta_{s'}
\tau'_{ R}(y)}
I_{\Gamma}(t, x_{t})\,dt
\geq E_{\cF_{s'}} I_{|x_{r}-z|< \rho /2  ,
\max_{[s',r]}|x_{t}-y| <R}I
$$
$$
\geq \mu_{0}\eta(q''-q')R^{2}P_{\cF_{s'}}\big
(|x_{r}-z|< \rho /2  ,
\max_{[s',r]}|x_{t}-y| <R\big)
$$
$$
\geq N^{-1}
(\rho/R)^{\nu}
\mu_{0}\eta(q''-q')R^{2}.
$$
This proves \eqref{7,27,1}.

By combining the two cases (i) and (ii) we conclude that
$$
G_{\tau} (\Gamma)\geq \min\big(\mu_{0}, \nu_{0}\mu_{R}(q''),
\mu_{0}\xi \eta^{n}(q''-q')^{n}\big)R^{2},
$$
and the arbitrariness of $\Gamma$ allows us to conclude that
\begin{equation}
                                               \label{7,27,5}
\mu (q')\geq \min\big(\mu_{0}, \nu_{0}\mu (q''),
\mu_{0}\xi\eta^{n}(q''-q')^{n}\big),
\end{equation}
whenever \eqref{1.5.1} holds.

Next, introduce   
$$
\hat{\mu}(q)=\min\big(\mu_{0},\mu (q)\big).
$$
Observe that in light of
\eqref{7,27,5} 
there exists $\varepsilon_{0}\in(0,1)$,
depending only on  $\delta, d,\kappa$,   such that
for any $0<q'<q''<1$ 
  satisfying $(1+\xi)q'\geq2q''$  we have   
\begin{equation}
                                                 \label{12.23.2}
\hat{\mu}(q')\geq\varepsilon_{0}\min\big( 
(q''-q')^{n },\hat{\mu}(q'')\big).
\smallskip\end{equation}
We also know that $\hat{\mu}(q)\geq \mu_{0} $ for $q\geq q_{0}$.

 We may certainly assume that $\varepsilon_{0}
\leq\bar{\varepsilon}:=2/(1+\theta)$ (recall that
$\theta>1$) and  
 we claim that for $q_{k}=
 \bar{\varepsilon}^{k}  q_{0}$,
$k=0,1,2,...$, we have   
\begin{equation}
                                                 \label{12.23.1}
\hat{\mu}(q_{k})\geq\varepsilon^{kn }_{0}\chi,\quad
\chi:=
\min\big(\mu_{0},q_{0}^{n }(1-
 \bar{\varepsilon})^{n }
\big).
 \end{equation}

To prove the claim we use induction. If $k=0$,
\eqref{12.23.1} is obvious. If it is true for a $k$,
then $q_{k}-q_{k+1}= \bar{\varepsilon} ^{k}
q_{0}(1- \bar{\varepsilon })$,   
$$
(q_{k}-q_{k+1})^{n }= \bar{\varepsilon }^{kn }
q_{0}^{n }(1- \bar{\varepsilon })^{n }
\geq \varepsilon_{0}^{kn }\chi,
$$
so that by \eqref{12.23.2}  
and the fact that $(1+\theta)q_{k+1}=2q_{k}$   \smallskip
$$
\hat{\mu}(q_{k+1})\geq\varepsilon_{0}\min\big(\varepsilon_{0}^{kn }
\chi,\hat{\mu}(q_{k})\big)\geq\varepsilon_{0}\varepsilon_{0}^{kn }
\chi\geq \varepsilon_{0}^{(k+1)n }
\chi.
$$

This proves \eqref{12.23.1} and shows that,
  if we define $r>1$ so that $\varepsilon_{0}^{n }
=\bar{\varepsilon}^{r}$, then 
$
\hat{\mu}(q_{k})\geq Nq_{k}^{ r }$ with 
$ r ,N>0$ depending only on
 $\delta, d,\kappa$. By observing that $\hat{\mu}$ 
is an increasing function we obtain that $\hat{\mu}(q)
\geq Nq^{ r }
$, $\mu(q)\geq Nq^{ r }$ for $q\leq1$. 
This yields \eqref{10.1.10} with $\gamma=r$
and  proves the theorem.  \qed

\begin{corollary}
                               \label{corollary 10.11.1}
For any 
$\kappa\in(0,1)$ there exists  
$N$, depending only on $\kappa,d,\delta $,  such that, for any $R\leq  \rho_{b}$, $y\in\bR^{d}$, $s\geq R^{2}$
on the set 
$$
A=\{ x_{s-R^{2}}\in \bar B_{\kappa R}(y)\}
$$
for any closed set
$\Gamma\subset  C_{R}(s,y)$, the (conditional) probability that 
after time $s'=s-R^{2}$ the process 
$(t, x_{t})$  
reaches $\Gamma$ before exiting from $[s-R^{2},s+R^{2})\times B_{R}(y) $
is greater than or equal to $N^{-1} (|\Gamma|/|C_{R}|)^{\mu-1/
(d +1)}$:
$$
P_{\cF_{s'}} (\tau_{\Gamma}  <
s'+(2R^{2})\wedge\theta_{s'}\tau'_{ R}(y) )
\geq N^{-1} (|\Gamma|/|C_{R}|)^{\gamma-1/(d +1)},
$$
where   $\tau_{\Gamma}  $ is the first time $(t, x_{t})$
hits $\Gamma$  after $s'$
  and $\gamma$ is taken from Theorem \ref{theorem 12.21.1}.
\end{corollary}

Indeed, set $q=|\Gamma|/|B_{R}|$ and  observe
that  on  $A$ owing to Theorem \ref{theorem 6.3,1}
$$
R^{2}q^{\gamma}\leq  NE_{\cF_{s'}}
\int_{s'}^{s'+\theta_{s'}\tau'_{ R}(y)}I_{\Gamma}(t,x_{t})\,dt
$$
$$
\leq NE_{\cF_{s'}}
I_{\tau_{\Gamma}  <s'+(2R^{2})\wedge\theta_{s'}\tau'_{ R}(y)}
E_{\cF_{\tau_{\Gamma}}}\int_{\tau_{\Gamma}}^{\tau_{\Gamma}+\theta_{\tau_{\Gamma}}\tau_{R}(y)}I_{\Gamma}(t,x_{t})\,dt
$$
$$
\leq NP_{\cF_{\tau}} (\tau_{\Gamma}  <s'+(2R^{2})\wedge\theta_{s'}\tau'_{ R}(y) )R^{2}q^{1/(d +1)}.
$$
 
The following corollary may look of having dubious
value. However, its PDE version is one of the main ingredients
in the Sobolev space theory of fully nonlinear
elliptic and parabolic equations (see \cite{Kr_18}).
 
\begin{corollary}
                 \label{corollary 10.1,1}    
For any 
$\kappa\in(0,1)$ there exists  
$N$, depending only on $\kappa,d,\delta $,  such that, for any $R\leq  \rho_{b}$, $y\in\bR^{d}$, $s\geq R^{2}$
and any Borel nonnegative $f$ 
on the set $A=\{ x_{s-R^{2}}\in B_{\kappa R}(y)\}$  
$$
\int_{C_{R}(s,y)}f^{\mu}(t,y)\,dydt
\leq NR^{d+2-2\mu}
\Big(E _{\cF_{s'}}\int_{s'}^{s'+\theta_{s'}\tau'_{ R}(y)}f (t,x_{t})\,dt\Big)^{\mu},
$$
where $N$ depends only on $\kappa,d,\delta $,
$\mu=1/(2\gamma)$.

\end{corollary}

Indeed, without losing generality assuming that
$f=0$ outside $C_{R}(s,y)$ and
setting
$$
u:=E _{\cF_{s-R^{2}}}\int_{s'}^{s'+\theta_{s'}\tau'_{ R}(y)}f (t,x_{t})\,dt,
$$
we have that 
for any $\lambda>0$  
$$  
u\geq \lambda E _{\cF_{s'}}\int_{s'}^{s'+\theta_{s'}\tau'_{ R}(y)}I_{f(t,x_{t})\geq\lambda}\,dt 
$$
$$
\geq\lambda N^{-1}R^{2}\big(|\{f\geq\lambda\} |/|C_{R}|\big)^{\gamma},\smallskip
$$
$$
|\{f\geq\lambda\} |\leq N R^{-2/\gamma}\lambda^{-1/\gamma}
|C_{R}|u^{1/\gamma}.
$$
It follows that for any $c>0$
$$
\int_{C_{R}(s,y)}f^{1/(2\gamma)}(x)I_{f>c}\,dx=
\big(1/(2\gamma)\big)\int_{c}^{\infty}\lambda^{1/(2\gamma)-1}
|\{f(y)>\lambda\}|\,d\lambda
$$
$$
\leq N R^{-2/\gamma} 
|C_{R}|u^{1/\gamma}c^{-1/(2\gamma)}.
$$
Also
$$
\int_{C_{R}(s,y)}f^{1/(2\gamma)}(x)I_{f\leq c}\,dx\leq
c^{1/(2\gamma)}|C_{R}|.
$$
For $c=uR^{-2}$ we have
$$
R^{-2/\gamma} 
|C_{R}|u^{1/\gamma}c^{-1/(2\gamma)}=
c^{1/(2\gamma)}|C_{R}|,
$$
$$
\int_{C_{R}(s,y)}f^{1/(2\gamma)}(x) \,dx
\leq N u^{1/(2\gamma)}R^{d+2-1/\gamma  }.
$$
This is what is claimed.

The following helps investigate the boundary
behavior of the probabilistic solutions
of parabolic equations.
\begin{corollary}
                      \label{corollary 10.21,1}
 
Let $R\leq  \rho_{b}$, $\nu,\kappa\in(0,1)$,
$y\in\bR^{d}$
 and assume that a closed set $\Gamma\subset
B_{R}(y)$ is such that, for any $r\in(0,R)$, 
$$
|B_{r}(y)\cap\Gamma|\geq
\nu |B_{r}|.
$$
Let $\tau$ be a stopping time.
 Then there exist    constants $\alpha\in(0,1)$
and $N$, depending only on $\kappa,d,\delta $,
   and $\nu$, such that,  
on the set $A:=\{\tau<\infty,x_{\tau}\in \bar B_{R/2}(y)\}$,
\begin{equation}
                        \label{10.21.1}
 P_{\cF_{\tau}}(\gamma_{\tau,R}(y) <\tau_{\Gamma} )\leq N( |x_{\tau}-y| /R)^{\alpha},
\end{equation}
where $\tau_{\Gamma} $ is the first time $  x_{ t}$
hits $\Gamma$ after $\tau$ and $\gamma_{\tau,R}(y)$ is the first time it exits from $B_{R}(y)$ after $\tau$
$(\gamma_{\tau,R}(y)=\tau+\theta_{\tau}\tau'_{R}(y))$.

\end{corollary}

Indeed, if $x_{\tau}\in\Gamma$,
\eqref{10.21.1} is obvious. Otherwise, let $R_{n}=R2^{-n}$, $\Gamma_{n}=\Gamma\cap  \bar B _{R_{n}}(y)$, 
($n=0,1,...$),
and let $n_{0}\geq 1$ satisfy  $x_{\tau}\in \bar B_{R_{n_{0}}}(y)\setminus\Gamma
 $. Then  define $\gamma_{n_{0}}=\tau$ and
$\gamma_{n_{0}-1}$ as the first  time    $x_{t}$   exits from
$B_{R_{n_{0}-1}}(y)$ after time $\tau$.
On $\{\gamma_{\tau,R}(y) <\tau_{\Gamma}\}$ we have $\gamma_{n_{0}-1}<\tau_{\Gamma}$ and, if $n_{0}\ge 2$, so that $x_{\gamma_{n_{0}-1}}\in\partial B_{R_{n_{0}-1}}(y)\setminus \Gamma$, define
$\gamma_{n_{0}-2}$ as the first exit time of $x_{t}$
from $B_{R_{n_{0}-2}}(y)$ after $\tau$, and so on.
In this way on $\{ 
x_{\tau}\in \bar B_{R_{n_{0}}}(y)\setminus\Gamma\}$ we define $\gamma_{n_{0}},
...,\gamma_{0}$ such that on $\{ 
x_{\tau}\in \bar B_{R_{n_{0}}}(y)\setminus\Gamma \}$ 
$$
\{\gamma_{\tau,R}(y) <\tau_{\Gamma} 
 \}= \bigcap_{i=1}^{n_{0} }\{\gamma_{n_{0}-i}
 <\tau_{\Gamma} \}.
$$
    
By Corollary \ref{corollary 10.11.1}
with 
$$
\gamma_{n_{0}-i},\quad 
[\gamma_{n_{0}-i}+R^{2}_ {n_{0}-i-1} ,\gamma_{n_{0}-i}+2R^{2}_ {n_{0}-i-1})\times \Gamma_{n_{0}-i-1}
$$
in place of $s-R^{2}$ and $\Gamma$, respectively, on the set $\{\gamma_{n_{0}-i}<\tau_{\Gamma}\}$
for $i= 1,...,n_{0}  $ we have
$$
P_{\cF_{\gamma_{n_{0}-i}}}(\gamma_{n_{0}-i-1 }
 <\tau_{\Gamma}) \leq q = q(d,\delta, \nu)<1
$$
It follows that on $A\cap \{ 
x_{\tau}\in \bar B_{R_{n_{0}}}(y)\setminus\Gamma \}$
$$
P_{\cF_{\tau}}(\gamma_{\tau,R}(y) <\tau_{\Gamma} )\leq q^{n_{0} }\leq q^{-1}(|x_{\tau}-y|/R)^{\ln_{2}(1/q)},
$$
where the second inequality is achieved
by using the largest possible $n_{0}$.
What we got is just a different form of \eqref{10.21.1}.

 The following makes one of crucial steps
in the proof of Harnack's inequality. Observe that in this theorem
we do not claim that
$q(\xi)\ne 0$ for $\xi$ not close to one. This fact
will be proved next.   
The main feature that distinguishes 
Theorem \ref{theorem 11.8.1} from
Corollary \ref{corollary 10.11.1} is that now the time
after which  the event in question
might happen is not separated from the time
at which $(t,x_{t})$ may reach  $\Gamma$:
$\tau_{\Gamma}\geq s'+R^{2}$ in 
Corollary \ref{corollary 10.11.1}.

\begin{theorem}
                                     \label{theorem 11.8.1}
For any $\kappa\in(0,1)$
there is a 
function $q(\xi)$, $\xi\in(0,1)$,
depending only on $\kappa,\delta,d $, 
and, naturally, on $\xi$,
such that for any $R\leq  \rho_{b}$, 
$y\in \bR^{d}$, $s\geq0$ on the set $A=\{ x_{s}\in
\bar B_{\kappa R}(y)\}$, for any
  closed $\Gamma\subset   C_{ R}(s,y)$ satisfying
$|\Gamma|\geq \xi|C_{ R}|$ we have
\begin{equation}
                               \label{9.17,1}
P _{\cF_{s}}(\tau_{\Gamma}  < s+\theta_{s}
\tau_{R}(y))\geq q(\xi),
\end{equation}
where $\tau_{\Gamma}  $ is the first time after $s$ the process 
$(t,  x_{t})$
hits $\Gamma$. Furthermore, 
$q(\xi)\to 1$ as $\xi\uparrow 1$.
Finally, for any closed $\Gamma'\subset B_{ R}$ satisfying
$|\Gamma'|\geq \xi|B_{ R}|$ on $A$ we have
\begin{equation}
                               \label{9.17,2}
I:=P _{\cF_{s}} (\tau'_{\Gamma}  <  s+\theta_{s}
\tau'_{R}(y) )\geq q(\xi),
\end{equation}
where $\tau'_{\Gamma}  $ is the first time after $s$ the process 
$ x_{t} $
hits $\Gamma$  $(\!$and  $s+\theta_{s}
\tau'_{R}(y)$
  is its first exit time  from
$B_{ R}(y)$  after $s$$)$.
\end{theorem}

Proof.
For any $\varepsilon\in(0,(1-\kappa)\sfp_{0}/4)$ we have   
$$
1-I \leq
P_{\cF_{s}}\Big(\theta_{s}
\tau _{R}(y) =\int_{s}
^{s+\theta_{s}
\tau _{R}(y)}  I_{C_{R}(s,y)
\setminus\Gamma}(t,x_{t})
\,dt\Big)
$$
$$
\leq P_{\cF_{s}}(\theta_{s}
\tau _{R}(y)\leq\varepsilon R^{2} )+\varepsilon^{-1}R^{-2} 
E_{\cF_{s}}\int_{s}
^{s+\theta_{s}
\tau _{R}(y) }  I_{C_{R}(s,y)
\setminus\Gamma}(t,x_{t})
\,dt.
$$
In light of Theorem  \ref{theorem 8.20.1}   and 
Lemma \ref{lemma 8.16.1}
we can estimate the right-hand side and then obtain
that on $A$
$$
P_{\cF_{s}}(\theta_{s}
\tau _{R}(y)\leq\varepsilon R^{2} )
\leq P_{\cF_{s}}\big( \theta_{s}
\tau _{(1-\kappa)R}(x_{s})
\leq\varepsilon R^{2} )\big)\leq Ne^{-1/(N\varepsilon)},
$$
$$
I\leq Ne^{-1/(N\varepsilon)}
+N\varepsilon^{-1}R^{d/(d+1)-2}|C_{R}(s,y\setminus\Gamma|^{1/(d+1)} 
$$
$$
\leq Ne^{-1/(N\varepsilon)}   
+N\varepsilon^{-1}(1-\xi)^{1/(d+1)},
$$
where the constants $N$  depend only on $d,\delta,\kappa $.
By denoting
$$
q(\xi)=1-\inf_{\varepsilon\in(0,(1-\kappa)\sfp_{0}/4)}\big(
Ne^{-1/(N\varepsilon)}
+N\varepsilon^{-1}(1-\xi)^{1/(d+1)}\big),
$$
we get what we claimed about \eqref{9.17,1}.

Estimate \eqref{9.17,2} follows from \eqref{9.17,1}
if one takes in the latter $\Gamma=[0,R^{2}]\times
\Gamma'$ and observes that
$$
\{\tau_{\Gamma}  < s+\theta_{s}
\tau_{R}(y)\}
\subset \{\tau'_{\Gamma}  <  s+\theta_{s}
\tau'_{R}(y) \}.
$$
The theorem is proved. \qed

In the next section we will need the following
fact of crucial importance, the origin
of which lies in \cite{KS_80} and \cite{Sa_80}.
A few other related results below also have their origins
 in \cite{KS_80} and
\cite{Sa_80} where the drift is bounded. 

In the following theorem we prove that $q(\xi)>0$,
$\xi\in(0,1)$,
in Theorem~\ref{theorem 11.8.1}.

\begin{theorem}
                                      \label{thm:4.1.10} 
For any 
$\kappa,\xi\in(0,1)$,  any $R\leq  \rho_{b}$, $y\in\bR^{d}$, $s\geq  0$
on the set $A=\{ x_{s }\in B_{\kappa R}(y)\}$
for any closed set
$\Gamma\subset  C_{R}(s,y)$ satisfying $|\Gamma|
\geq \xi |C_{R}|$, the (conditional) probability that 
after time $s $ the process 
$(t, x_{t})$  
reaches $\Gamma$ before exiting from $C_{R}(s,y) $
is greater than or equal to $q 
=q(\kappa,d,\delta ,\xi)>0$:
$$
P_{\cF_{s }} (\tau_{\Gamma}  <s+\theta_{s}
\tau_{R}(y))
\geq q ,
$$
where $\tau_{\Gamma}  $ is the first time $(t, x_{t})$
hits $\Gamma$. 
\end{theorem}

Proof.  It is convenient for the sake of
simplicity of the proof to assume that $\xi=1/n$
where $n$ is an integer. Obviously this does not restrict generality.
Observe that for $s'=s+\xi R^{2}/2$ and $\Gamma'=\Gamma\cap C_{R}(s',y)$ we have 
$$
|\Gamma'|\geq |\Gamma|
-(\xi/2 )|C_{R}|\geq ( \xi/2 )|C_{R}|.
$$
Then consider the cylinders $(r,r+\xi R^{2}/4)\times
B_{R}(y)$, $r\in [s+\xi R^{2}/2,s+R^{2}]$. It is easy to see that for at least one of them
$$
|\Gamma'\cap (r,r+\xi R^{2}/4)\times
B_{R}(y) |\geq \frac{1}{ 2n (4n-2)}
|(r,r+\xi R^{2}/4)\times
B_{R}(y) |.
$$
Let $(r_{0},r_{0}+\xi R^{2}/4)\times
B_{R}(y)$ be one of them and set $\rho^{2}=\xi R^{2}/4$, $\rho>0$. Then by representing
$B_{R}(y)$ as the union of $B_{\rho}(z)$ we can find
$y_{0}$ such that
$$
|\Gamma'\cap C_{\rho}(r_{0},y_{0}) |\geq \eta
|C_{\rho}(r_{0},y_{0}) |,
$$
where $\eta>0$ depends only on $d,\xi$.

  Then
by Corollary \ref{corollary 10.11.1}, for 
$x_{r_{0}-\rho^{2}}\in B_{  \kappa \rho}(y_{0})$
 the  probability that the process
$(t, x_{t})$ will hit $\Gamma$ before exiting
from $C_{2\rho_{0}^{2} ,\rho_{0} }(r_{0}-\rho^{2} ,y_{0})$ is estimated from 
below by a strictly positive constant depending only
on $\kappa,  d,\delta,\xi$.
After that it only remains to invoke Theorem
\ref{theorem 12.7.2} observing that that $r_{0}-\rho^{2}\geq s+\xi R^{2}/4$.
The theorem is proved.  \qed

\mysection[Estimates in \protect\lpq  
of potentials]{Further estimates of potentials
in $L_{(q,p)}$}

Some arguments in the future have to be repeated
at least twice in slightly different
situations. In order to avoid this
we consider the following setting.

 Take some   $p,q\in[1,\infty]$, $\alpha \in[0,1]$
and fix $ \varkappa_{0}=\varkappa_{0}(\sfp_{0}) >1 $ such that
\index{$S$@Miscelenea!$\varkappa_{0}$}%
the right-hand side of \eqref{8.20.1} equals $1/2$ when
 $\rho=(\varkappa_{0}-1)\lambda^{-1/2}$.

 \begin{assumption} 
                            \label{assumption 8.24.1}
There exists a constant $\ell>0$ such that
 for
any $\lambda\geq \varkappa_{0}^{2}\rho_{b}^{-2}$,   stopping time $\tau$, and Borel $f\geq0$ given on $\bR^{d+1}$, it holds that
\begin{equation}
                              \label{8.24.1}
 E_{\cF_{\tau}}\int_{0}^{\theta_{\tau}\tau_{\rho(\lambda)}(x_{\tau})}e^{-\lambda s}f(\tau+s,x_{\tau+s})\,ds\leq \ell\lambda ^{-\alpha}\|f\|_{L_{(q,p)}},
\end{equation}
 where $\rho(\lambda)=
 \varkappa_{0}\lambda^{-1/2} $ $(\leq\rho_{b})$).    
\end{assumption}

\begin{remark}
                        \label{remark 8.24.1}
As it follows from Theorem \ref{theorem 5.5.1} with
$c_{t}=\lambda,r_{t}=1, A_{t}=t,\gamma=\tau+\theta_{\tau}\tau_{\rho(\lambda)}$, Assumption \ref{assumption 8.24.1}
is satisfied, for instance, if $p,q\in[1,\infty]$,  $d/p+1/q\leq 1$  
and $\alpha =1-d/(2p)-1/q$  with 
$$
\ell=N(d,\delta)(1+\bar b^{2}_{\rho(\lambda)})^{d/(2p)}
\leq N(d,\delta)(1+\sfb_{0}^{2})^{d/(2p)}=N(d,\delta).
$$
 
 Below we will see that we can take a wider
range of parameters $p,q$.
 \end{remark}

\begin{remark}
                        \label{remark 10.15.1}
Since the $L_{(q,p)}$-norm is translation invariant,
one gets an equivalent assumption if 
$(\tau+s,x_{\tau+s})$ is replaced with  
$(\nu+s ,x_{\tau+s}+y)$ for any $\cF_{\tau}$-measurable
real-valued
$\nu$ and $\bR^{d}$-valued $y$. In short, in Assumption \ref{assumption 8.24.1} one can replace \eqref{8.24.1} with
\begin{equation}
                              \label{10.15.1}
 E_{\cF_{\tau}}\int_{0}^{\theta_{\tau}\tau_{\rho(\lambda)}(x_{\tau})}e^{-\lambda s}f( \nu+s,x_{\tau+s}+y)\,ds\leq \ell\lambda ^{-\alpha}\|f\|_{L_{(q,p)}}.
\end{equation}

\end{remark}

\begin{lemma}
                                     \label{lemma 8.22.1}

 Under Assumption \ref{assumption 8.24.1} 
for any  stopping time $\tau$, $s_{0}\geq0$, 
$y_{0}\in \bR^{d}$, $\lambda\geq \varkappa_{0}^{2}\rho_{b}^{-2} $,
and Borel nonnegative $f$ vanishing outside 
$C_{\lambda^{-1/2}}(s_{0},y_{0})$
we have on the set $\{\tau<\infty\}$ that
\begin{equation}
                                             \label{8.22.10}
 E_{\cF_{\tau}}\int_{0}^{\infty}e^{- \lambda s}
f(s,x_{\tau+s}-x_{\tau})  \,ds\leq 
N(\sfp_{0})\ell\lambda^{-\alpha } 
\Phi_{\lambda}(s_{0},y_{0})\|  f\|_{L_{(p ,q) }},
\end{equation}
where   $\Phi_{\lambda}
(t,x)=e^{-\sqrt{\lambda}  (\sqrt t+|x|)\sfp_{0}/4}$.
 
\end{lemma}

Proof. We assume that the event $\{\tau<\infty\}$
occurred. Introduce
  $\gamma^{0}$ as the first time $(s,x_{\tau+s}-x_{\tau})$, $s\geq0$,
hits $\bar C_{\lambda^{-1/2}}(s_{0},y_{0})$ and set $\tau^{0}$
as the first time after $\gamma^{0}$ this process
 exits from $C_{ 
\varkappa_{0}\lambda^{-1/2}  }(s_{0},y_{0})$ after $\gamma^{0}$.
We define recursively $\gamma^{k}$, $k=1,2,...$, as the first
time after $\tau^{k-1}$ the process
$(s ,x_{s}-x_{\tau+s}-x_{\tau})$
hits $\bar C_{\lambda^{-1/2}}(s_{0},y_{0})$ and $\tau^{k}$ as 
the first time after $\gamma^{k}$ this process
 exits from 
$C_{ \varkappa_{0}\lambda^{-1/2} }(s_{0},y_{0})$.
As is easy to see $\tau+\tau^{k}$ and $\tau+\gamma^{k}$
are stopping times and they  are either
infinite or lie between $\tau+s_{0}$ and $\tau+s_{0}+
\varkappa_{0}^{2}\lambda^{-1}$.

The left-hand side of
\eqref{8.22.10} on $\{\tau<\infty\}$ equals
\begin{equation}
                                                \label{8.22.1}
 E_{\cF_{\tau}}
\sum_{k=0}^{\infty} e^{- \lambda \gamma^{k}  }I_{k},
\end{equation}
where
$$
I_{k}=I_{\gamma^{k}<\infty }
E_{\cF_{\tau+\gamma^{k}}} \int_{\gamma^{k} }^{ \tau^{k} 
 }
e^{-\lambda  (s-\gamma^{k}) }
f( s,x_{\tau +s}-x_{\tau}) \,ds  
$$
$$
=I_{\gamma^{k}<\infty }E_{\cF_{\tau+\gamma^{k}}} \int_{0 }^{\tau^{k}-\gamma^{k}
 }
e^{-\lambda r }
f(  \gamma^{k} +r,x_{\tau+\gamma^{k}+r }-x_{\tau}) \,dr .
$$
By Remark \ref{remark 10.15.1} the last expression
is less than $\ell\lambda^{-\alpha}\|f\|_{L_{  (p , q) }}$,
that is $I_{k}\leq \ell\lambda^{-\alpha}\|f\|_{L_{  p , q }}$.

Next,   observe that, if $\sqrt{s_{0}}>|y_{0}| $, then  $\gamma^{0} $
is bigger than the first exit time of $( s,x_{\tau+s}-x_{\tau})$
from $C_{\sqrt{s_{0}}}$, that is $\gamma^{0} \geq
\theta_{\tau}\tau_{ \sqrt{s_{0} }}$ and
by Theorem~\ref{theorem 8.20.1} 
$$
E_{\cF_{\tau}}e^{- 
\lambda  \gamma^{0}  }\leq N(\sfp_{0})e^{-\sqrt{\lambda} 
\sqrt{s_{0}}\sfp_{0}/2}.
$$
In case $\sqrt{s_{0}}\leq |y_{0}| $ and $|y_{0}|> \lambda^{-1/2}$
our $\gamma^{0}$
is bigger than $\theta_{\tau}\tau_{ |y_{0}|-\lambda^{-1/2}} $ and
$$
E_{\cF_{\tau}}e^{- 
\lambda  \gamma^{0}  }\leq Ne^{-\sqrt{\lambda}(|y_{0}|-\lambda^{-1/2})\sfp_{0}/2}.
$$
The last estimate (with   $N=1$) also holds if 
$|y_{0}|\leq\lambda^{-1/2}$, so that
in case $\sqrt{s_{0}}\leq |y_{0}| $
$$
E_{\cF_{\tau}}e^{- 
\lambda  \gamma^{0}  }
\leq Ne^{-\sqrt{\lambda} |y_{0}|\sfp_{0}/2}
$$
and we conclude that in all cases
$$
E_{\cF_{\tau}}e^{- 
\lambda  \gamma^{0} } \leq 
N(\sfp_{0})e^{-\sqrt{\lambda} (\sqrt{t_{0}}+|y_{0}|)\sfp_{0}/4}.
$$

Furthermore,  by the choice of $\varkappa_{0}$ and
 Theorem \ref{theorem 8.20.1}
$$
E_{\cF_{\tau+\tau^{k-1}}}e^{-\lambda ( \gamma^{k}
- \tau^{k-1} )} \leq \frac{1}{2},   
$$
$$
E_{\cF_{\tau}}e^{-\lambda  \gamma^{k }}=
E_{\cF_{\tau}}e^{-\lambda  \tau^{k-1}  }E_{\cF_{\tau+\tau^{k-1}}}
 e^{-\lambda ( \gamma^{k} 
- \tau^{k-1} )}  
\leq \frac{1}{2}E_{\cF_{\tau}}e^{- \lambda   \tau^{k-1}  },
$$
so that  
$$
E_{\cF_{\tau}}e^{-\lambda   \gamma^{k} }\leq\frac{1}{2}
 E_{\cF_{\tau}}e^{- \lambda  \gamma^{k-1} } ,\quad
 E_{\cF_{\tau}}e^{- \lambda  \gamma^{k} } \leq 2^{-k}
E_{\cF_{\tau}}e^{- \lambda  \gamma^{0}  }.
$$

Recalling \eqref{8.22.1} we see that the left-hand side 
of \eqref{8.22.10}
is indeed dominated by
$$
N\ell\lambda^{-\alpha}\Phi_{\lambda}(s_{0},y_{0})\|f\|_{L_{(p ,q) }}
$$  
and the lemma is proved. \qed

 The following theorem shows that the time spent
by $(s,x_{s})$ in cylinders $C_{1}(0,x)$
decays very fast as $|x|\to\infty$. Recall that
$$
\bR^{d+1}_{0}=(0,\infty)\times\bR^{d}.
$$
\begin{theorem}
                               \label{theorem 8.22.1}
 Under Assumption \ref{assumption 8.24.1}
 for any  stopping time $\tau$, any $\lambda\geq \varkappa_{0}^{2}\rho_{b}^{-2}$, and   Borel nonnegative $f$  
on the set $\{\tau<\infty\}$
we have 
\begin{equation}
                                     \label{8.22.4}
I:= E_{\cF_{\tau}}\int_{0}^{\infty}e^{- \lambda s}
f(s,x_{\tau+s}-x_{\tau}) \,ds\leq 
N (\sfp_{0})\ell\lambda^{-\alpha }\|\Psi_{\lambda}  f\|_{L_{(p ,q) }(\bR^{d+1}_{0})},
\end{equation}
where $\Psi _{\lambda}(t,x)=\exp(- \sqrt{\lambda} 
(|x|+ \sqrt t)\sfp_{0}/16)$.
\end{theorem}

Proof.  
Set $\zeta(t,x)=\lambda^{ (d+2)/2 }\eta(\lambda t,\sqrt{\lambda} x)$, where $\eta$ has unit integral
and is proportional to the indicator of $C_{1}$,  and for $(t,x),(r,y)\in\bR^{d+1}$
set
$$
f_{r,y}(t,x)=f(t,x)\zeta(t-r,x-y).
$$
Clearly, due to Lemma \ref{lemma 8.22.1},
$$
I
=\int_{\bR}\int_{\bR^{d }}E_{\cF_{\tau}}\int_{0}^{\infty}
e^{-\lambda s}
f_{r,y}(s ,x_{\tau+s}-x_{\tau}) \,ds\,dydr
$$
$$
\leq N(\sfp_{0})\ell\lambda^{-\alpha } \int_{-1/\lambda}^{\infty}\int_{\bR^{d }}
\Phi_{\lambda}(r_{+},y)\|  f_{r,y}\|_{L_{(p ,q) }}\,dydr.
$$
{\em Case $\infty>p\geq q$}. Introduce
$$
M_{1}^{q/(q-1)}= 
\int_{-1/\lambda}^{\infty}\int_{\bR^{d }}\Phi^{q/(2q-2)}_{\lambda}(r_{+},y)
\,dydr,
$$
$$
M_{2}^{p/(p-q)}=
\int_{-1/\lambda}^{\infty}\int_{\bR^{d }}\Phi^{pq/(4p-4q)}_{\lambda}(r_{+},y)
\,dydr,\quad p\ne q,\quad M_{2}=1,\quad p=q.
$$
It follows by H\"older's 
inequality  that
$$
\lambda^{\alpha } I
\leq N\ell M_{1}\Big(
\int_{ -1/\lambda}^{\infty}\int_{\bR^{d }}
\Phi_{\lambda}^{q /2}(r_{+},y)\int_{0}^{\infty}\Big(\int_{\bR^{d}}
f_{r,y}^{p }(t,x)\,dx\Big)^{q /p }dt
\,dydr\Big)^{1/q }
$$
$$
=N\ell M_{1}\Big(
\int_{0}^{\infty}dt \Big(\int_{-1/\lambda}^{\infty}
\int_{\bR^{d }}\Phi_{\lambda}^{q /2}(r_{+},y)
 \Big(\int_{\bR^{d}}
f_{r,y}^{p }(t,x)\,dx\Big)^{q /p } 
\,dydr\Big)\Big)^{1/q } 
$$

$$
\leq N\ell M'\Big(
\int_{0}^{\infty}\Big(\int_{-1/\lambda}^{\infty}\int_{\bR^{d }}\int_{\bR^{d }}
\Phi_{\lambda}^{p /4}(r_{+},y)f_{r,y}^{p }(t,x)\,dydrdx 
\Big)^{q /p }dt\Big)^{1/q },
$$
where $M'=M_{1}M^{1/q}_{2}$.

We replace $\Phi_{\lambda}^{p /4}(r_{+},y)$ by $\Phi_{\lambda}^{p /4}(t,x)$
taking into account that these values are comparable
as long as $\zeta(t-r,x-y)\ne0$. After that
integrating over $dydr$ and computing $M_{1},M_{2}$
lead immediately to
\eqref{8.22.4}. 

{\em Case $p< q<\infty$}. 
It follows by H\"older's inequality that  
$$
\lambda^{ \alpha }I\leq N\ell M_{3}\Big(
\int_{-1/\lambda}^{\infty}\int_{\bR^{d }}
\Phi_{\lambda}^{p /2}(r_{+},y)\int_{\bR^{d}}\Big(\int_{0}^{\infty}
f_{r,y}^{q}(t,x)\,dt\Big)^{p/q  }dx
\,dydr\Big)^{1/p }
$$
$$
\leq N\ell M''\Big(\int_{\bR^{d}}dx\Big(
\int_{-1/\lambda}^{\infty}\int_{\bR^{d }}\int_{0}^{\infty}
\Phi_{\lambda}^{q /4}(r_{+},y)f_{r,y}^{q}(t,x)\,dtdydr\Big)^{p/q}
\Big)^{1/p},
$$
where $M''=M_{3}M^{1/p}_{4}$,
$$
M_{3}^{p/(p-1)}= 
\int_{-1/\lambda}^{\infty}\int_{\bR^{d }}\Phi^{p/(2p-2)}_{\lambda}(r_{+},y)
\,dydr,
$$
$$
M_{4}^{q/(q-p)}=
\int_{-1/\lambda}^{\infty}\int_{\bR^{d }}\Phi^{pq/(4q-4p)}_{\lambda}(r_{+},y)
\,dydr .
$$
This leads   to
\eqref{8.22.4} as above.

{\em Case $p\leq q=\infty$}. If $p<\infty$
it suffices to use H\"older's inequality only once.
If $p=\infty$ it suffices to observe that
$$
\Phi _{\lambda} (r_{+},y)\sup_{t,x}[f(t,x)\zeta(t-r,x-y)]
$$
$$
\leq N(\sfp_{0})\Psi _{\lambda} (r_{+},y)
 \sup_{t,x} [\Psi_{\lambda}  (t,x)
f(t,x) ].
$$
Similarly one treats the remaining case 
$p=\infty>q$.
 The theorem is
proved. \qed

If  $f(t,x)=f(x)$ we come to the following.
\begin{corollary}
                           \label{corollary 10.4.1}
Under Assumption \ref{assumption 8.24.1}
for any $\lambda\geq \varkappa_{0}^{2}\rho_{b}^{-2}$ and   Borel nonnegative $f(x)$  
we have on $\{\tau<\infty\}$
\begin{equation}
                                          \label{10.4.2}
 E_{\cF_{\tau}}\int_{0}^{\infty}e^{- \lambda s}
f( x_{\tau+s}-x_{\tau}) \,dt\leq 
N(\sfp_{0})\ell\lambda^{-\alpha-1/q} \|\Psi_{\lambda}  f\|_{L_{p}(\bR^{d})},
\end{equation}
where $\Psi _{\lambda}( x)=\exp(- \sqrt{\lambda}  |x|\sfp_{0}/16)$.
\end{corollary}    
 
Next result  is dealing with the exit times of the process
$x_{s}$ rather than $(s,x_{s})$. We will need it
while showing an improved integrability of Green's functions. Assumption \ref{assumption 8.24.1}
is no longer needed.

Estimate \eqref{9.29.5} below in case $b$ is bounded
was the starting point for the theory
of {\em time homogeneous\/} controlled diffusion processes
about fifty five years ago.
 \begin{lemma}
                                      \label{lemma 9.29.1}
Let $p \in [d,\infty]$ and let $\tau$ be a stopping time. Then for any Borel 
nonnegative $f(x)$  and $\rho\leq \rho_{b}$, 
 
\begin{equation}
                                       \label{9.29.5}
 E_{\cF_{\tau}}\int_{0}^{\theta_{\tau}\tau' _{ \rho}(x_{\tau}) }   
f( x_{\tau+s})\,ds\leq N(\delta,d )
\rho 
 ^{2-d/  p  } 
\|f\|_{L_{p}(\bR^{d})}.
\end{equation}

\end{lemma}

Proof. If $p=d$, the result follows from 
Lemma \ref{lemma 5.6.1}. Indeed, as
in Theorem \ref{theorem 6.3,1} we have $A\leq 3\rho^{2}+\rho^{2}$, and by definition, $B/\rho\leq \bar b_{\rho} \leq \bar b_{\rho_{b}}\leq \sfb_{0} $.

In case $p>d$ observe that in
$$
 E_{\cF_{\tau}}\int_{0}^{\theta_{\tau}\tau' _{ \rho}(x_{\tau}) }   
f( x_{\tau+s} )\,ds\leq N(\delta,d )
 \rho 
\|f\|_{L_{p}(\bR^{d})} 
$$
 the   norm can be taken only over $B_{\rho}(x_{\tau})$ because
$x_{\tau+s}\in B_{\rho}(x_{\tau})$ before $\theta_{\tau}\tau' _{\rho}(x_{\tau})$. After that we replace $\rho 
\|f\|_{L_{d}(B_{\rho}(x_{\tau}))}$ by $N\rho ^{2-d/p}
\|f\|_{L_{p}(B_{\rho}(x_{\tau}))}$ by using H\"older's inequality
and then come  to \eqref{9.29.5}. The lemma is proved. \qed

\mysection[Green's functions]{Green's functions. After Fabes--Stroock}  
                                 \label{section 10.26.1}
 It could be a good time to remind the reader that
the assumptions which hold throughout this chapter are stated
in the Introduction to the chapter including Assumption 
\ref{assumption 8.19.2}.
In this section we  
take a stopping time $\tau$ and $\lambda\geq \varkappa_{0}^{2}\rho_{b}^{-2}$, and introduce a measure on $\cF_{\tau}\otimes
\frB(\bR^{d+1})$ by
$$
G(\Gamma):= E I_{\tau<\infty}\int_{0}^{\infty}e^{- \lambda s}
I_{\Gamma}(\omega,s,x_{\tau+s}-x_{\tau}) \,ds.
$$
We do not include $\tau,\lambda$ in the notation $G(\Gamma)$ because they are assumed to be fixed
up until Theorem \ref{theorem 8.27.1}.

Take  $\cA\in\cF_{\tau}$ such that
$\cA\subset \{\tau<\infty \}$ and define
a measure on $\bR^{d+1}$ by
$$
G_{\cA}(\Lambda)= G(\cA\times\Lambda) ,
$$
Then for any Borel $f\geq 0$
on $\bR^{d+1}$ by
Remark \ref{remark 8.24.1}  and  
Theorem \ref{theorem 8.22.1} we have
$$
\int_{\bR^{d+1}}f\,G_{\cA}(dtdx)=
\int_{\Omega\times \bR^{d+1}}I_{\cA}f(t,x)\,G(d\omega dt dx)
$$
$$
=EI_{\cA }\int_{0}^{\infty}e^{- \lambda s}
f(s,x_{\tau+s}-x_{\tau}) \,ds
$$
$$\leq 
N (d,\delta) \lambda^{-d/(2d+2) }\|\Psi_{\lambda}  f\|_{L_{d+1}(\bR^{d+1}_{0})}P(\cA ).
$$
This shows that the measure $G_{\cA}(\Lambda)$
has a density, which we denote by $G_{\cA}(t,x)$, and the following result holds true.
 
\begin{theorem}
                         \label{theorem 9.3.1}
The function $G_{\cA}$ is Borel measurable,
is such that $G_{\cA}(t,x)=0$ for $t\leq0$,
\begin{equation}
                                         \label{9.3.5}
 \|\Psi^{-1}_{\lambda} G_{\cA}\|_{L_{(d+1)/d } } 
\le N (d,\delta)\lambda^{ -d/(2d+2 )}P(\cA) ,
\end{equation}
and for any Borel nonnegative $f$ given on $\bR^{d+1}$   
we have
\begin{equation}
                                         \label{10.16.1}
EI_{\cA}\int_{0}^{\infty}e^{- \lambda s}
f(s,x_{\tau+s}-x_{\tau}) \,ds =\int_{\bR^{d+1}}f(t,x)G_{\cA}(t,x)\,dxdt.
\end{equation}
\end{theorem}

In light of \eqref{10.16.1} it is natural to call $G_{\cA}$
the Green's function of the process $(s,x_{\tau+s}-x_{\tau})$
on $\cA$.

Here is a remarkable property of
$G_{\cA}$.  

\begin{theorem}
                 \label{theorem 9.28,1}
Let $\tau=0$, $\cA=\Omega$, so that
\eqref{10.16.1} becomes
\begin{equation}
                      \label{9.28,2}
E \int_{0}^{\infty}e^{- \lambda s}
f(s,x_{ s} ) \,ds =\int_{\bR^{d+1}}f(t,x)G_{\Omega}(t,x)\,dxdt,
\end{equation}
where $G_{\Omega}(t,x)$ is nonrandom. 
Then for any $\varepsilon\in(0,\rho^{2}_{b})$
\begin{equation}
                      \label{9.28,3}
\int_{C_{\rho_{b}}\setminus\{t<\varepsilon\}}G^{-\mu}_{\Omega}(t,x)\,
dxdt\leq N(\varepsilon,d,\delta,\rho_{b},\lambda)<\infty,
\end{equation}
where $\mu>0$ is taken from
Corollary \ref{corollary 10.1,1}.
\end{theorem}

Proof. By Corollary \ref{corollary 10.1,1}
the left-hand side of \eqref{9.28,3}
is dominated by a constant times
$$
\Big(E\int_{0}^{\tau_{R}}G^{-1}_{\Omega}(t,x_{t})\,dt\Big)^{\mu}\leq N
\Big(E\int_{0}^{\infty}e^{-\lambda t}I_{C_{\rho_{b}}}G^{-1}_{\Omega}(t,x_{t})\,dt\Big)^{\mu}
$$
where in light of \eqref{9.28,2} the expectation is less than $|C_{\rho_{b}}|$.
\qed

Estimate \eqref{9.3.5} shows that
$G_{\cA}$ is summable to the power $(d+1)/d$.
 It turns out that,
actually,
 is summable to a higher  power. The proof of this
is based on the parabolic version of
Gehring's lemma from \cite{GS_82} (see our Appendix).

Introduce $\bC_{+}$ as the set of cylinders $C_{R}(t,x)$, $R>0$,  
$t\geq0$, $x\in\bR^{d}$. 
\index{$B$@Sets!$\bC_{+}$}%
For $C=C_{R}(t,x)\in \bC_{+}$ let
 $2C=C_{2R}(t,x)$.
 \index{$B$@Sets!$2C$}%
  If $C\in \bC_{+}$ and $C=C_{R}(t,x)$,
we call $R$ the radius of $C$.  

\begin{theorem}
                                  \label{theorem 9.3.2}
There exist  $d_{0}\in(1,d)$ and    $N $, depending only
on $\delta,d $,  
such that for any     $C\in \bC_{+}$ of radius 
$\rho\leq \varkappa_{0} /(2\sqrt\lambda)$ $(\leq\rho_{b}/2)$ and $p\geq d_{0}+1 $, we have
\begin{equation}
                                \label{10.14.01}
\| G_{ \cA}\|_{L_{p/(p-1)}(C)}\leq N \rho^{-(d+2) /p }
\| G_{\cA}\|_{L_{1}(2C )} , 
\end{equation}
which is equivalently rewritten as
$$
\Big(\dashint_{C}G^{p/(p-1)}_{\cA}\,dxdt\Big)^{(p-1)/p}
\leq N\dashint_{2C}G_{\cA}\,dxdt.
$$
 
\end{theorem}

Proof. We basically follow the idea in \cite{FS_84}.
Take $C\in \bC_{+}$ of radius $\rho\leq\varkappa_{0}/(2\sqrt\lambda) $ and on the set $\{\tau<\infty\}$
define 
recursively  $\gamma^{0}$ as the first  time
  when the process $(t ,x_{\tau+t}-x_{\tau})$, $t\geq0$,
hits $\bar C$, $\tau^{0}$ as the first time after $\gamma^{0}$
when this process leaves $2C$, $\gamma^{n }$ as the first  time
after $\tau^{n-1}$ when the process $(t ,x_{\tau+t}-x_{\tau})$
hits $\bar C$, $\tau^{n }$ as the first time after $\gamma^{n }$
when this process leaves $2C$.
 
Then for any nonnegative Borel $f$ vanishing outside $C$
such that $\|f\|_{L_{d+1}(C)}=1$ 
we have
$$ 
I:=
\int_{C}f (t,x)  G_{\cA}(t,x)\,dxdt
$$  
$$
=\sum_{n=0}^{\infty}EI_{\cA}e^{- \lambda \gamma^{n} } 
E_{\cF_{\tau+\gamma^{n}}} \int_{0}^{\tau^{n}-\gamma^{n}}e^{-\lambda t }f(t+\gamma^{n} ,x_{\tau+t+\gamma^{n} }-x_{\tau})\,dt.
$$

Next, we use  
\eqref{9.29.2}  and what was said about the relation       
of \eqref{8.24.1} to \eqref{10.15.1}  and take into account 
that $\bar b_{\varkappa_{0}/\sqrt\lambda}\leq\bar b_{\rho_{b}}\leq 1$ 
 to see that the conditional expectations
above are less than $N(d,\delta)\rho^{d/(d+1)}$. After that
we use   Corollary
\ref{corollary 7.29.1}   
to get that
$$
\rho^{2}\leq N(\sfp_{0}) E_{\cF_{\tau+\gamma^{n}}} \int_{\gamma^{n}}^{\tau^{n}}
e^{-\lambda(t-\gamma^{n})} \,dt. 
$$
Then we obtain   
$$
I
\leq N  \rho ^{-(d+2)/(d+1)}
\sum_{n=0}^{\infty}EI_{\cA} e^{-\lambda \gamma^{n} }
\int_{\gamma^{n}}^{\tau^{n}}
e^{-\lambda(t-\gamma^{n})} \,dt
$$
$$
=N  \rho ^{-(d+2)/(d+1)}
\sum_{n=0}^{\infty}EI_{\cA}  
 \int_{\gamma^{n} }^{\tau^{n} }e^{-\lambda t } \,dt
$$
$$
\leq N  \rho ^{-(d+2)/(d+1)}
EI_{\cA} \int_{0}^{\infty}e^{-\lambda t}I_{2C}(t,x_{\tau+t}-x_{\tau}) \,dt
$$
$$
=N  \rho ^{-(d+2)/(d+1)}\int_{2C}G_{\cA}(t,x)\,dxdt.
$$

The arbitrariness of $f$ 
implies that
$$
\Big(\dashint_{C}  G_{\cA} ^{(d+1)/d}(t,x)\,dxdt
\Big)^{d/(d+1)}\leq N \dashint_{2C}  G_{\cA}(t,x)\,dxdt.
$$

Now the assertion of the theorem for $p=d_{0}+1$ follows directly from
the parabolic version of the famous Gehring's lemma
stated as Proposition 1.3  in \cite{GS_82}
(see our Appendix). 
For larger $p$ it suffices to use H\"older's inequality.
The theorem is proved. \qed

\begin{theorem}
                                           \label{theorem 2.3.1}
For any $p\geq d_{0}+1$  
\begin{equation}
                              \label{8.26.2}
\|G_{\cA}\|_{L_{p/(p-1)}(\bR^{d+1}_{0})}\leq N(\delta,d  )\lambda^{(d+2)/(2p)-1}P(\cA)
\end{equation}
In particular, for any Borel $f\geq0$ given on $\bR^{d+1}$ (and $\lambda\geq \varkappa_{0}^{2}\rho_{b}^{-2}$) on $\{\tau<\infty\}$ we have
\begin{equation}
                               \label{8.26.5} 
E_{\cF_{\tau}}\int_{0}^{\infty}e^{-\lambda t}f(t,x_{\tau+t})\,dt
\leq N(\delta,d )
\lambda^{(d+2)/(2d_{0}+2)-1}\|f\|_{L_{d_{0}+1}}.
\end{equation}
 
\end{theorem}

Proof.   Represent $ [0,\infty)\times \bR^{d}$
as the union of countably many $C_{1},C_{2},...\subset \bC_{+}$
of radius $\varkappa_{0}/(2\sqrt\lambda)$
so that each point in $\bR^{d+1}_{0}$ belongs  
to no more than $m(d)$  of the $2C_{i}$'s. Then 
$$
\|G_{\cA}\|_{L_{p/(p-1)}(\bR^{d+1}_{0})}
\leq \big\|\sum_{i}I_{C_{i}}G_{ \cA}\big\|_{L_{p/(p-1)}(\bR^{d+1}_{0})}
$$
$$
\leq\sum_{i} \| G_{\cA}\|_{L_{p/(p-1)}(C_{i})}
\leq N(\delta,d)\lambda^{(d+2)/(2p)}\sum_{i} \| G_{\cA}\|_{L_{1}(2C_{i})}
$$
$$
\leq N_{1}\lambda^{(d+2)/(2p)}
\| G_{\cA}\|_{L_{1}(\bR^{d+1}_{0})}= 
N_{1}\lambda^{(d+2)/(2p)-1}P( \cA ).
$$
This proves \eqref{8.26.2} and the fact that
$$
EI_{\cA}\int_{0}^{\infty}e^{- \lambda t}
f(t,x_{\tau+t}-x_{\tau}) \,dt
$$
$$
\leq 
N(\delta,d )
\lambda^{(d+2)/(2d_{0}+2)-1}\|f\|_{L_{d_{0}+1}}P(\cA).
$$
The arbitrariness of $\cA$ shows that \eqref{8.26.5}
holds with $f(t,x_{\tau+t}-x_{\tau})$ in place of 
$f(t,x_{\tau+t} )$. One then eliminates 
$x_{\tau}$ as in Remark \ref{remark 10.15.1}.
The theorem is proved. \qed

\begin{remark}
                             \label{remark 8.26.1}

If $\lambda\in(0,\varkappa_{0}^{2}\rho_{b}^{-2})$, one can also give
an estimate of the left-hand side $J$ of \eqref{8.26.5} 
by taking nonnegative $f\in L_{p}(\bR^{d+1}_{0})$
and observing that, for $\lambda_{0}=\varkappa_{0}^{2}\rho_{b}^{-2}$,
$$
J 
=\sum_{n=0}^{\infty}e^{-\lambda n}
E_{\cF_{\tau}} \int_{n}^{n+1}e^{-\lambda (t-n)}f(t,x_{\tau+t})\,dt
$$
$$
\leq \sum_{n=0}^{\infty}e^{\lambda_{0}-\lambda}e^{-\lambda n}
E_{\cF_{\tau}}\int_{n}^{n+1}e^{- \lambda_{0} (t-n)}f(t,x_{\tau+t})\,dt
$$
$$
=\sum_{n=0}^{\infty}e^{\lambda_{0}-\lambda}e^{-\lambda n}
E_{\cF_{\tau}}E_{\cF_{\tau+n}}\int_{0}^{ 1}e^{- \lambda_{0} t}f(n+t,x_{\tau+n+t})\,dt
$$
where  each conditional expectation in the sum is dominated
by 
$$
N\|fI_{[n,n+1)}\|_{L_{p }(\bR^{d+1}_{0})}
$$
 in light of \eqref{8.26.5} . Therefore
$$
J\leq N\sum_{n=0}^{\infty} e^{-\lambda n}\|fI_{[n,n+1)}\|_{L_{p }(\bR^{d+1}_{0})}
\leq N(1-e^{-\lambda})^{-(p-1)/p}\|f\|_{L_{p}(\bR^{d+1}_{0})},
$$
where the second inequality follows
from H\"older's inequality. 
\end{remark}

Similar improvement of integrability occurs for
the Green's function of $x_{t}$ rather than $(t,x_{t})$.
Observe that 
$$
g_{\cA}(x):=\int_{0}^{\infty}G_{\cA}(t,x)\,dt
$$
satisfies
$$
EI_{\cA}\int_{0}^{\infty}e^{- \lambda s}
f( x_{\tau+s}-x_{\tau}) \,ds =\int_{\bR^{d }}f( x)g_{\cA}( x)\,dx 
$$
for any Borel nonnegative $f$ on $\bR^{d}$. For this reason
we call $g_{\cA}$ the Green's function of $x_{\tau+s}-x_{\tau}$ on
$\cA$.

If $\mu>0$ by Jensen's inequality
$$
g^{-\mu}_{\cA}(x)\leq\Big(\int_{\rho_{b}/2}^{\rho_{b}}
G_{\cA}(t,x)\,dt\Big)^{-\mu}\leq N
\int_{\rho_{b}/2}^{\rho_{b}}G^{-\mu}_{\cA}(t,x)\,dt,
$$
which along with Theorem \ref{theorem 9.28,1} leads
to the following.
\begin{theorem}
                 \label{theorem 9.28,2}
Let $\tau=0$, $\cA=\Omega$.
Then  
\begin{equation}
                      \label{9.28,30}
\int_{B_{\rho_{b}} }g^{-\mu}_{\Omega}(x)\,
dx \leq N(  d,\delta,\rho_{b},\lambda)<\infty,
\end{equation}
where $\mu>0$ is taken from
Corollary \ref{corollary 10.1,1}.
\end{theorem}

By using Remark \ref{remark 8.24.1} and Corollary \ref{corollary 10.4.1}
with $p=d$, $q=\infty$ we come to the following.
\begin{theorem}
                                       \label{theorem 10.4.1}
 We have
\begin{equation}
                                         \label{10.4.5}
 \|\Psi^{-1}_{\lambda} g_{\cA}\|_{L_{d/(d-1)}(\bR^{d}) } 
\le N(d,\delta)\lambda^{-1/2}P(\cA ),
\end{equation}
where   $\Psi _{\lambda}( x)=\exp(- \sqrt{\lambda}
  |x| \sfp_{0}/16)$.
 
\end{theorem}

According to this theorem 
$g_{\cA}$ is summable to the power $d/(d-1)$.
Again it turns out that this power can be increased. 
If $B$ is an open ball in $\bR^{d}$ by $2B$ we 
\index{$B$@Sets!$2B$}%
denote
the concentric open ball of twice the radius of $B$.  

\begin{theorem}
                                  \label{theorem 10.4.3}
There exist  $d_{0}\in(1,d)$ and a constant    $N $, depending only
on $d,\delta$, 
 such that for any ($\lambda\geq\varkappa_{0}^{2}\rho_{b}^{-2}$)  ball  $B$ of radius $\rho\leq \varkappa_{0}/(2\sqrt\lambda)$ 
 and $p\geq d_{0}  $, we have
\begin{equation}
                                          \label{10.4.6}
\| g_{\cA}\|_{L_{p/(p-1)}(B)}\leq N \rho^{-d /p }
\| g_{\cA}\|_{L_{1}(2B )} , 
\end{equation}
which is equivalently rewritten as
$$
\Big(\dashint_{B}g^{p/(p-1)}_{\cA}\,dx \Big)^{(p-1)/p}
\leq N\dashint_{2B }g_{\cA}\,dx.
$$
 
\end{theorem}

Proof. We again follow the idea in \cite{FS_84}.
Take a ball  $B$ of radius $\rho\leq \varkappa_{0}/(2\sqrt\lambda)$ and
on the set $\{\tau<\infty\}$
define 
recursively  $\gamma^{0}$ as the first  time
  when the process $ x_{\tau+t}-x_{\tau} $
hits $\bar B$, $\tau^{0}$ as the first time after $\gamma^{0}$
when this process leaves $2B$, $\gamma^{n }$ as the first  time
after $\tau^{n-1}$ when the process $ x_{\tau+t}-x_{\tau} $
hits $\bar B$, $\tau^{n }$ as the first time after $\gamma^{n }$
when this process leaves $2B$.
 
Then for any nonnegative Borel $f$ vanishing outside $B$
with $\|f\|_{L_{d }(B)}=1$ 
we have
$$ 
I:=
\int_{B}f ( x)  g_{\cA}( x)\,dx
$$  
\begin{equation}
                                 \label{7.17.1}
=\sum_{n=0}^{\infty}EI_{\cA}e^{- \lambda \gamma^{n} } 
E_{\cF_{\tau+\gamma^{n}}} \int_{0}^{\tau^{n}-\gamma^{n}}e^{-\lambda t }f( x_{\tau+\gamma^{n}+t}-x_{\tau})\,dt.
\end{equation}

Next we use   
\eqref{9.29.5} with $p=d$
 to see that the conditional expectations
above are less than $N\rho $. After that
we use   Corollary   
\ref{corollary 7.29.1} 
to get that   
$$
\rho^{2}I_{\tau+\gamma^{n}<\infty}\leq N E_{\cF_{\tau+\gamma^{n}}} \int_{\gamma^{n}}^{\tau^{n}}
e^{-\lambda(t-\gamma^{n})} \,dt.
$$
Then we obtain  
$$
\int_{B}f   g_{\cA}( x)\,dx 
\leq N  \rho ^{-1}
\sum_{n=1}^{\infty}EI_{\cA} e^{-\lambda \gamma^{n} }
\int_{\gamma^{n}}^{\tau^{n}}
e^{-\lambda(t-\gamma^{n})} \,dt
$$
$$
=N  R ^{-1}
\sum_{n=1}^{\infty}EI_{\cA}  
 \int_{\gamma^{n} }^{\tau^{n} }e^{-\lambda t } \,dt
$$
$$
\leq N  \rho ^{-1}
EI_{\cA} \int_{0}^{\infty}e^{-\lambda t}I_{2B}( x_{\tau+t}-x_{\tau}) \,dt
=N  \rho ^{-1}\int_{2B}g_{\cA}( x)\,dx .
$$

The arbitrariness of $f$ 
implies that
$$
\Big(\dashint_{B}  g_{\cA} ^{d/(d-1)}(x)\,dx 
\Big)^{(d-1)/d}\leq N \dashint_{2B}  g_{\cA}(x)\,dx,
$$
and again it only remains to use Gehring's lemma
in case $p=d$.
For larger $p$ it suffices to use H\"older's inequality.
The theorem is proved.  \qed

By mimicking the proof of Theorem \ref{theorem 2.3.1}
one gets its ``elliptic'' counterpart.

\begin{theorem}
                                      \label{theorem 2.3.2}
For any $p\geq d_{0}$  we have
$$
\|g_{\cA}\|_{L_{p/(p-1)}(\bR^{d})}
 \leq N(\delta,d )
\lambda ^{d/(2p)-1}.
$$
In particular, for any Borel $f\geq0$ given on $\bR^{d}$ \(and $\lambda\geq \varkappa_{0}^{2}\rho_{b}^{-2}$\),
on $\{\tau<\infty\}$ we have
\begin{equation}
                               \label{8.26.4}
E_{\cF_{\tau}}\int_{0}^{\infty}e^{-\lambda t}f(x_{\tau+t})\,dt
\leq N(\delta,d )
\lambda ^{d/(2d_{0})-1}\|f\|_{L_{d_{0}}(\bR^{d})}.
\end{equation}
\end{theorem}

\begin{remark}
                                           \label{remark 2.7.1}
 Denote by   $\sfd_{0}= \sfd_{0}(d,\delta)<d$,
  the maximum of the $d_{0}$'s 
  \index{$S$@Miscelenea!$\sfd_{0}$}%
from Theorems \ref{theorem 9.3.2} and 
 \ref{theorem 10.4.3}.

Observe that, as the simple example of $a^{ij}=\delta^{ij}$
and $b\equiv0$
shows, $\sfd_{0}(d,\delta)>d/2$ and $\sfd_{0}(d,1)$
can be taken to be as close to $d/2$  as we wish. We call $ \sfd_{0}(d,\delta)$ {\em the Fabes-Stroock constant\/} because these authors discovered and proved
in \cite{FS_84}
its existence in terms of PDEs  for elliptic equations.  
 
\end{remark}

The function $g_{\cA}$ also possesses additional
properties which make it interesting from the point of view
of Real Analysis. Similar results for the Green's
functions in domains can be found in \cite{Kr_21}.

\begin{theorem}[doubling property]
                                          \label{theorem 10.28.2}
Let $\rho\leq \rho_{b}/2$, $B=B_{\rho}(y)$ and 
for Borel $\Gamma$ define
$$
g_{\cA}(\Gamma)=\int_{\Gamma}g_{\cA}(x)\,dx.
$$  
Then
$g_{\cA}(2B )\leq Ng( B   )$, where $N=N(d,\delta,\lambda)$.
 \end{theorem}

Proof. We follow part of the proof of Theorem
\ref{theorem 10.4.3} but take $f=I_{2B}$
and again use \eqref{9.29.5} with $p=d$. Then we see
that
$$
I\leq N\sum_{n=0}^{\infty}EI_{\cA}e^{- \lambda \tau^{n} }\rho^{2}I_{\tau+\tau^{n}<\infty}.
$$
After that repeating the manipulations from the 
proof of Theorem
\ref{theorem 10.4.3} but using Theorem \ref{theorem 7.16.1}
instead of Corollary \ref{corollary 7.29.1} leads to the desired conclusion.
The theorem is proved.
  \qed

\begin{corollary}[$A_{\infty}$-property of $g_{\cA}$]
                                  \label{corollary 10.28.100}
There are constants $\mu\geq1$ and $N$,
depending only on $d,\delta,\lambda$, such that for any
ball of radius $\rho\leq \varkappa_{0}/(2\sqrt\lambda)$  and Borel $\Gamma\subset B$
we have
\begin{equation}
                                          \label{10.28.5}
N\frac{g_{\cA}(\Gamma)}{g_{\cA}(B)}\geq\Big(\frac{|\Gamma|}{|B|}\Big
)^{\mu}.
\end{equation}
\end{corollary}

Proof. Take the same $\gamma^{n},\tau^{n}$
as in the proof of Theorem \ref{theorem 10.4.3}
and observe that  $g_{A}(\Gamma)$   
is the $E_{\cF_{\tau}}$-expectation
of the sum  of
$$
 E_{\cF_{\sigma^{n}}}
\int_{\sigma^{n}}^{\tau +\tau^{n }}
e^{-\lambda (t-\tau)}I_{\Gamma+x_{\tau}}(x_{t})\,dt
$$
$$
=e^{-\lambda\gamma^{n}}E_{\cF_{\sigma^{n}}}
\int_{\sigma^{n}}^{\tau +\tau^{n }}
e^{-\lambda (t-\sigma^{n})}I_{\Gamma+x_{\tau}}(x_{t})\,dt=:e^{-\lambda\gamma^{n}}I_{n} .
$$
over $n=0,1,...$, where $\sigma^{n}=\tau+
\gamma^{n}$. Define
$$
\hat \Gamma=(\sigma^{n}+\rho^{2},
\sigma^{n}+5\rho^{2})\times(\Gamma+x_{\tau}),
$$
$$
C_{2\rho}(s,y)=[\sigma^{n}+\rho^{2},
\sigma^{n}+5\rho^{2})\times 2( B+x_{\tau}).
$$
and notice that $\hat \Gamma \subset C_{2\rho}(s,y)$ and 
$$
\frac{|\hat\Gamma|}{|C_{2\rho}|}=2^{-d}\frac{|\Gamma|}{|B|}=:q.
$$
Clearly,
$$
I_{n}\geq e^{-5\lambda\rho^{2}}
E_{\cF_{\sigma^{n}}}\int_{\sigma^{n}}
^{\sigma^{n}+\theta_{\sigma^{n}}\tau'_{2\rho}(y)}I_{\hat \Gamma}(t,x_{t})\,dt.
$$

In light of Theorem \ref{theorem 12.21.1}
and Corollary \ref{corollary 3.7.1} this allows us to conclude that
$$
Nq^{-\gamma}e^{-\lambda\gamma^{n}}I_{n}\geq  e^{-\lambda\gamma^{n}}\rho^{2}\geq
 e^{-\lambda\gamma^{n}}
E_{\cF_{\sigma^{n}}}\int_{\sigma^{n}}
^{\sigma^{n}+\theta_{\sigma^{n}}\tau'_{2\rho}(y)}I_{B}( x_{ t}-x_{\tau})\,dt
$$
$$
\geq 
E_{\cF_{\sigma^{n}}}\int_{\sigma^{n}}
^{\sigma^{n}+\theta_{\sigma^{n}}\tau'_{2\rho}(y)}e^{-\lambda(t-\tau)}I_{B}( x_{ t}-x_{\tau})\,dt=:E_{\cF_{\sigma^{n}}}J_{n}.
$$

We note that
$$
\sum_{n=0}^{\infty}J_{n}=\int_{0}^{\infty}
e^{-\lambda t}I_{B}(x_{\tau+t}-x_{\tau})\,dt
$$
and get
$
g_{\cA}(\Gamma)\geq N^{-1}q^{\gamma}g_{\cA}(B)$. \qed

Corollary \ref{corollary 10.28.100} is almost identical 
to Corollary 2.3 in \cite{FS_84}. However, there are no lower order terms
in \cite{FS_84} and the comparable situations would be only when
$x_{t}$ was a solution of \eqref{11.29.2}. 
 
\begin{remark}
                              \label{remark 12.4.1}
Once we know that $g_{\cA}$ is an $A_{\infty}$-weight, it is also
an $A_{p}$-weight for certain large $p$. In particular,
on any bounded closed $\Gamma $,
$g_{\cA}^{-\alpha}$ is  summable  for some $\alpha>0$. This we 
already know from
  Theorem~\ref{theorem 9.28,2}.

\end{remark}

Above $\tau,\lambda$ were fixed in order to
make the notation shorter.
Now we allow them to change.

The estimates \eqref{8.26.5} and \eqref{8.26.4}
have the same spirit as \eqref{4.24.1} and \eqref{5.6.1}
and, by virtually repeating the proof of
Theorem \ref{theorem 5.5.1},
  we come to the following. Recall that
we say that $(\sfd_{0},q,p)$ are properly tight  if  
\begin{equation}  
                            \label{9.27.3}
p ,q  \in[1,\infty],\quad
 \nu:=1-\frac{\sfd_{0}}{p }-\frac{1}{q }\geq 0 .
\end{equation}

\begin{theorem}
                 \label{theorem 8.27.1}
Suppose   that $(\sfd_{0},q,p)$ are properly tight.
Then for any Borel $f\geq0$ given on $\bR^{d+1}$ 
\(recall that $\lambda\geq \varkappa_{0}^{2}\rho_{b}^{-2}$\)
we have on $\{\tau<\infty\}$ that
\begin{equation}
                                 \label{8.27.01}
E_{\cF_{\tau}}\int_{0}^{\infty}e^{-\lambda t}f(t,x_{\tau+t})\,dt
\leq N(d,\delta )
\lambda^{ (1/2)(d/p+2/q)-1}\|f\|_{L_{(q,p)}}.
\end{equation}

\end{theorem}

By using the same argument as in Remark \ref{remark 10.15.1}
we can replace $(t,x_{\tau+t})$ with $(\tau+t,x_{\tau+t})$ in
\eqref{8.27.01} and then we see that Assumption
\ref{assumption 8.24.1} is satisfied for $p,q$ as in
   \eqref{9.27.3}, $\ell=N(d,\delta )$ and $\alpha=
1-(1/2)(d/p+2/q)$. Then Theorem \ref{theorem 8.22.1}
is valid, which yields the following.

\begin{theorem}
                                  \label{theorem 8.30.1}
 
Suppose   that $(\sfd_{0},q,p)$ are properly tight.
Then for any Borel $f\geq0$ given on $\bR^{d+1}$ 
and $\lambda\geq \varkappa_{0}^{2}\rho_{b}^{-2}$ 
we have on the set $\{\tau<\infty\}$ that
\begin{equation}
                                 \label{8.27.1}
E_{\cF_{\tau}}\int_{0}^{\infty}e^{- \lambda s}
f(s,x_{\tau+s}-x_{\tau}) \,ds
\leq N(\delta,d )
\lambda^{ -\chi}\|\Psi_{\lambda}f\|_{L_{(q,p)}(\bR^{d+1}_{0})}.
\end{equation}
where $\Psi _{\lambda}(t,x)=\exp(- 
\sqrt{\lambda} (|x|+ \sqrt t)\sfp_{0}/16)$, $\chi=1-(1/2)(d/p+2/q)$.
In particular, if $f$ is independent of $t$, $p\geq \sfd_{0}$,
 and $q=\infty$
$$
E_{\cF_{\tau}}\int_{0}^{\infty}e^{- \lambda t}
f( x_{\tau+t}-x_{\tau}) \,dt\leq 
N\lambda ^{-1+d/(2p)}\|\bar \Psi_{\lambda}^{d_{0}/p} f\|_
{L_{p}(\bR^{d})},
$$
where $\bar \Psi _{\lambda}( x)=\exp(- 
\sqrt{\lambda}  |x| \sfp_{0}/16)$.  
\end{theorem}

\begin{theorem}
                           \label{theorem 9.7.1}
Suppose   that $(\sfd_{0},q,p)$ are properly tight.
Then

(i)
  for any
$n=1,2,...$, nonnegative Borel $f$ on $\bR^{d+1}_{0}$, 
  and
 $T\leq \varkappa_{0}^{-2}\rho_{b}^{2}$  we have
on $\{\tau<\infty\}$ that
\begin{equation}
                                          \label{9.7.1}
E_{\cF_{\tau}}\Big[\int_{0}^{T}  
f(t,x_{\tau+t}-x_{\tau})\,dt\Big]^{n}\leq n!N^{n} (d,\delta)
T^{n\chi }\| \Psi^{(1-\nu)/n} _{1/T}
f\|^{n}_{L_{(q,p)}(\bR^{d+1}_{0}) },
\end{equation}

(ii)  for any
  nonnegative Borel $f$ on $\bR^{d+1}_{0}$,
$\rho\leq 1$,  and
 $T\geq \varkappa_{0}^{-2}\rho_{b}^{2}$  we have
on $\{\tau<\infty\}$ that
\begin{equation}
                               \label{9.7.10}
I:=E_{\cF_{\tau}} \int_{0}^{T}  
f(t,x_{\tau+t} )\,dt \leq N(d,\delta,\rho_{b}) T\rho^{-2-d}\sup_{C\in\bC_{\rho}}
 \|  
f\| _{L_{(q,p)}(C) } .  
\end{equation}

\end{theorem}

Proof. To prove   (i) we proceed  by induction
on $n$. The induction hypothesis is that for   $\kappa\in[0,1/n]$, any $\tau$, and any $\bR^{d+1}_{0}$-valued $\cF_{\tau}$-measurable $(\gamma,\xi)$
$$
J_{n}(\tau,\gamma,\xi,T):=
E_{\cF_{\tau}}\Big[\int_{0}^{T}  
f(\gamma+t, x_{\tau+t}-x_{\tau}+\xi )\,dt\Big]^{n}
$$
\begin{equation}
                                   \label{10.28.1}
\leq n!N^{n} 
T^{n\chi } \Psi_{1/T}^{( \nu-1)\kappa n}
(\gamma,\xi)\| \Psi^{( 1-\nu)\kappa} _{1/T}
f\|^{n}_{L_{(q,p)}(\bR^{d+1}_{0}) }.
\end{equation}

Denote $\theta_{\tau}x_{t}=x_{\tau+t}-x_{\tau}$ and observe that for $s\geq \tau+t$ we have $\theta_{\tau}x_{s}
=\theta_{\tau+t}x_{s-t}+\theta_{\tau}x_{t}$.
If the hypothesis  holds true for some $n\geq1$,
then  by observing that
$$
 E_{\cF_{\tau}}\Big[\int_{0}^{T}  
f(\gamma+t,\theta_{\tau}x_{t}+\xi)\,dt\Big]^{n+1}
$$
$$
=(n+1)E_{\cF_{\tau}}\int_{0}^{T}f(\gamma+t,\theta_{t}x_{\tau}+\xi) J_{n}(\tau+t,\gamma+t,\theta_{\tau}x_{t}+\xi,T-t)\,dt,
$$
we see that, for any $\kappa\in[0,1/n]$,
\begin{equation}
                                        \label{8.31.1}
J_{n}(\tau,\gamma,\xi,T)
 \leq (n+1)!N^{n}T^{n\chi }\|\Psi^{(1-\nu)\kappa}_{1/T}
f\|^{n}_{L_{(q,p)}(\bR^{d+1}_{0}) }I,
\end{equation}
where
$$
I=E_{\cF_{\tau}}\int_{0}^{T}\Psi^{(\nu-1)\kappa n}_{1/T}f(\gamma+t,\theta_{\tau}x_{t}+\xi)\,dt.
$$
  We have for any $\lambda>0$
$$
I\leq e^{\lambda T}
E_{\cF_{\tau}}\int_{0}^{\infty}e^{-\lambda t}\Psi^{(\nu-1)\kappa n}_{1/T}f(\gamma+t,\theta_{\tau}x_{t}+\xi)\,dt,
$$
 where the last term,
owing to Theorem \ref{theorem 8.30.1}, for $\lambda=1/T$ 
and $\mu\in[0,1]$ is
dominated by
$$
N(\delta,d ) T^{\chi} \|\Psi^{(\nu-1)\kappa n}_{1/T}f(\gamma+\cdot,\xi+\cdot)\Psi^{\mu}_{1/T} \|_{L_{(q,p)}(\bR^{d+1}_{0})}
$$
$$
\leq N(\delta,d ) T^{\chi}\Psi^{-\mu}(\gamma,\xi) \|\Psi^{(\nu-1)\kappa n+\mu}_{1/T}f(\gamma+\cdot,\xi+\cdot)  \|_{L_{(q,p)}(\bR^{d+1}_{0})},
$$
where the last inequality is due to the fact that
$\Psi _{\lambda}(s,y)\leq \Psi  _{\lambda}(t+s,x+y)
\Psi^{-1} _{\lambda}(t, x)$.

The estimate of $I$ for $n=1$, $\kappa\in[0,1]$,
and $\mu=(1-\nu)\kappa$ yields \eqref{10.28.1}
with $n=1$  after replacing 
$\Psi^{(\nu-1)\kappa }_{1/T}f$ by $f$,
which justifies the start of the induction.

For $\mu=(1-\nu)\kappa(n+1)$, $\kappa\in[0,1/(n+1)]$, we have
$\Psi^{(\nu-1)\kappa n+\mu}_{1/T}=
\Psi^{(1-\nu)\kappa   }_{1/T}$ and this along with 
\eqref{8.31.1} show  that our hypothesis holds true
also for $n+1$.
This proves \eqref{9.7.1}.

While proving \eqref{9.7.10} we may assume that  
$\rho=1$ (see Remark \ref{remark 2.29.1}) and that
  $T=k \beta$, where $k\geq 1$ is an integer
and $\beta=\varkappa_{0}^{-2}\rho_{b}^{2}$.  Then 
first consider the case of $\nu=0$. Note that
owing to \eqref{9.7.1}  
\begin{equation}
                                    \label{12.18.03}
E_{\cF_{\tau}}\int_{0}^{\beta}f(t,x_{\tau+t})
\,dt\leq N(d,\delta)\rho_{b}^{(2d_{0}-d)/p}\|f(\cdot,\cdot+x_{\tau})\Psi_{1/\beta}
 \|_{L_{(q,p)}(\bR^{d+1}_{0})} 
\end{equation}

Let $\cZ=\{0,1,2,...\}\times \bZ^{d}$ and for
$z=(z_{1},z_{2})\in\cZ$ let $C^{z}=C_{ d}(z)$. Observe that
on $C^{z}$ we have  
$$
\Psi_{1/\beta} \leq  \exp(- 
2\mu(|z_{2}|+ \sqrt z_{1} ) ),
$$
where $2\mu=\beta^{-1/2} \sfp_{0}/16$.
Furthermore, for each $z\in\cZ$
$$
\|f(\cdot,\cdot+x_{\tau})\|_{L_{(q,p)}(C^{z})}\leq \sup_{C\in\bC_{1}}
 \|  
f\| _{L_{(q,p)}(C) }.
$$
Therefore, by noting that
$ f(\cdot,\cdot+x_{\tau})\Psi_{1/\beta}  \leq \sum_{\cZ}  f(\cdot,\cdot+x_{\tau} )\Psi_{1/\beta} 
I_{C^{z}}$ and using Minkowski's inequality we get
that the norm in \eqref{12.18.03} is dominated by
$$
N(d)\sup_{C\in\bC_{1}}
 \|  
f\| _{L_{(q,p)}(C) }\sum_{\cZ}
 \exp(- 2
\mu(|z_{2}|+ \sqrt z_{1} ) ).
$$
By majorating the last sum by an integral
  we obtain
that it is dominated by
$$
\int_{0}^{\infty}\int_{\bR^{d}}e^{-2\mu (|x|+ \sqrt t-d-1)_{+} }\,dxdt
$$
$$
\leq N+
\int_{0}^{\infty}\int_{\bR^{d}}e^{-2\mu (|x|+ \sqrt t-d-1)_{+} }I_{|x|+\sqrt t>2(d+1)}\,dxdt
$$
$$
\leq N+
\int_{0}^{\infty}\int_{\bR^{d}}e^{- \mu (|x|+ \sqrt t )  }I_{|x|+\sqrt t>2(d+1)}\,dxdt
$$
$$
\leq N+
\int_{0}^{\infty}\int_{\bR^{d}}e^{- \mu (|x|+ \sqrt t )  } \,dxdt=N+N\mu^{-d-2}.
$$
Hence, for $n=0$
$$
E_{\cF_{\tau+n\beta}}\int_{n\beta}^{(n+1)\beta}f(t,x_{\tau+t})
\,dt\leq N(d,\delta)\rho_{b}^{(2d_{0}-d)/p}
(1+\rho_{b}^{d+2})
\sup_{C\in\bC_{1}}\|
f\| _{L_{(q,p)}(C) }.
$$
Clearly, this also holds for any $n=1,2,...$ and since
$T=k\beta=k\varkappa_{0}^{-2}\rho_{b}^{2}$,
$$
I \leq N(d,\delta)T\rho_{b}^{-2}\rho_{b}^{(2d_{0}-d)/p}
(1+\rho_{b}^{d+2})
\sup_{C\in\bC_{1}}\|
f\| _{L_{(q,p)}(C) } 
$$
and this  
 proves \eqref{9.7.10} if $\nu=0$.

If $\nu=1$ ($p=q=\infty$), \eqref{9.7.10} is obvious,
and if $\nu<1$, by the above, \eqref{9.7.10}
holds with $(1-\nu)(p,q)$ in place of $(p,q)$,
which yields \eqref{9.7.10} as is after using 
H\"older's inequality.
The theorem is 
proved. \qed

\begin{remark}
                     \label{remark 2.29.1}
The term $\sup_{C\in\bC_{1}}
 \|  
f\| _{L_{(q,p)}(C) }$ in \eqref{9.7.10}
can be replaced
with $\rho^{-2-d}\sup_{C\in\bC_{\rho}}
 \|  
f\| _{L_{(q,p)}(C) }$ if $\rho\leq 1$.

Indeed, by simple inspection one proves that
for any $\rho\geq 1$, $C\in\bC_{\rho}$ and $C'$,
defined as the union of $2C$ and its reflection
in its lower base,
$$
 I_{ C}( t,  x)\leq N\int_{C'}I_{C_{1}}( t- s,  x-  y)\,
dsdy,
$$
where $N=N(d)$. Dilations show that, for any
$\rho\leq 1$  and $C\in \bC_{1}$ we have
$$
 I_{ C}( t,  x)\leq N\rho^{-2-d}\int_{C'}I_{C_{\rho}}( t- s,  x-  y)\,
dsdy.
$$
 It follows
for $C\in \bC_{1}$ and $\rho\leq1$ that
$$
 |f|I_{ C}\leq N\rho^{-2-d}\int_{C'}|f|
I_{C_{\rho}(s,y)}\,
dsdy, 
$$
$$
\rho^{2+d} \|f\|_{L_{q,p}(C)}
\leq N\int_{C'}
\sup_{C\in\bC_{\rho}}\|
f\| _{L_{(q,p)}(C) }\,dsdy=
N 
\sup_{C\in\bC_{\rho}}\|
f\| _{L_{(q,p)}(C) }.
$$
\end{remark}

Next theorem improves Theorem \ref{theorem 5.5.1}  in what concerns
the restrictions on $p,q$.
\begin{theorem}
                          \label{theorem 9.5.1}  
Suppose   that $(\sfd_{0},q,p)$ are properly tight.
 Then 
for any $\rho\leq \rho_{b} $, $\cF_{\tau}$-measurable $\bR^{d}$-valued $y$,  and Borel nonnegative $f$ 
given on $\bR^{d+1}$,
we have on $\{\tau<\infty\}$ that  
\begin{equation}
                                \label{9.5.4}
E_{\cF_{\tau}}\int_{0}^{\theta_{\tau}\tau_{\rho}(y) }f( t, x_{\tau+t} -y)\,dt\leq
N(d,\delta,p,q)\rho^{2}\dashnorm f\|_{L_{(q,p)}(C_{\rho})},
\end{equation}
\begin{equation}
                             \label{9.3.1}
E_{\cF_{\tau}}\int_{0}^{\theta_{\tau}\tau'_{\rho}(y)  }f( t, x_{\tau+t}-y)\,dt\leq
N(d,\delta,p,q)\rho^{2}
\sup_{C\in\bC_{\rho}}
 \dashnorm  
f\| _{L_{(q,p)}(C) }.
\end{equation}

\end{theorem}

Proof.  Since $\theta_{\tau}\tau_{\rho}(y) \leq \rho^{2}$, the left-hand side of
\eqref{9.5.4} is smaller than
$$
e^{\lambda \rho^{2}\varkappa_{0}^{2}}E_{\cF_{\tau}}\int_{0}^{\infty}e^{-\lambda \varkappa_{0}^{2}t}
I_{C_{\rho}}f(t,x_{\tau+t} -y)\,dt
$$
for any $\lambda>0$. For $\lambda=\rho^{-2}$ we have
$\lambda \varkappa_{0}^{2}\geq \varkappa_{0}^{2}\rho_{b}^{-2}$
and \eqref{9.5.4} follows from \eqref{8.27.1}.

To prove \eqref{9.3.1}, it suffices to note that \eqref{9.5.4} remains valid if its right-hand side is replaced with that of \eqref{9.3.1}, and then repeat the same argument
as in the proof of \eqref{1.3.3}.
 The theorem
is proved. \qed

Here is a key to finding analytic
conditions ensuring  that $\bar b_{\rho_{b}}\leq\sfb_{0}$.

\begin{corollary}
                                       \label{corollary 3.26.1}
Assume that there exists 
  functions $b_{i}(t,x)\geq0$,
$i=0,1,...,k$, on $\bR^{d+1}$ such that
$|b_{t}|\leq (b_{0}+...+b_{k})(t,x_{t})$ for all $(\omega,t)$. Take some
$p_{i},q_{i}$ such that  that $(\sfd_{0},q_{i},p_{i})$  are  properly tight. Suppose that
there is
a constant $\hat b\in(0,\infty)$
such that, for any $\rho\leq  \rho_{b}$ and $C\in\bC_{\rho}$
\begin{equation}
                                                \label{3.26.2}
\sum_{i=0}^{k}\dashnorm b_{i}\|_{L_{(q_{i},p_{i})}(C)}\leq\hat b \rho^{ -1}.
\end{equation}
Then $\bar b_{  \rho_{b}}\leq  N
(d,\delta)\hat b$. 

\end{corollary}

Indeed, \eqref{9.5.4} implies that
$$
E_{\cF_{t}}\int_{0}^{\theta_{t}\tau_{\rho}(x) }(b_{0}+...+b_{k})( t+s, x_{t+s}  )\,ds
$$
$$
\leq N(d,\delta)\rho^{2}\sum_{i=1}^{k}\dashnorm b_{i}I_{C_{\rho}(t,x_{t}+x) }\|_{L_{(q_{i},p_{i})} }
\leq N(d,\delta)\hat b \rho .
$$

 \begin{remark}
                                 \label{remark 10.18.1}
In light of Corollary \ref{corollary 3.26.1} it is tempting
to claim that if the analytic condition \eqref{3.26.2} holds and $N(d,\delta)\hat b
\leq \sfb_{0}$, then our main Assumption    \ref{assumption 8.19.2} 
is satisfied.  However, there is a vicious circle:
Corollary \ref{corollary 3.26.1} was obtained on the basis
of Assumption     \ref{assumption 8.19.2}. Nevertheless,
we will see in the case of stochastic equations that
if condition \eqref{3.26.2} is satisfied and 
$N(d,\delta)\hat b
< \sfb_{0}$, then there {\em exists\/} at least one solution
of the equation, for which Assumption      \ref{assumption 8.19.2} 
is satisfied. Recall that we already know 
 from Lemma \ref{lemma 12.11.1} some sufficient analytic conditions for Assumption      \ref{assumption 8.19.2}  to be satisfied.

\end{remark}

One also has an estimate similar to 
\eqref{9.3.1} for $\rho>\varkappa_{0}^{-1}\rho_{b}$, albeit, with not so sharp control of the constants.

\begin{theorem} 
                      \label{theorem 3.27.20}

Suppose   that $(\sfd_{0},q,p)$ are properly tight.
Then
for   any $\rho\in(0,\infty)$,   $\cF_{\tau}$-measurable $\bR^{d}$-valued $y$,
and Borel $f\geq0$ on $\{\tau<\infty\}$ we have
\begin{equation}
                                              \label{3.21.80}
E_{\cF_{\tau}}\int_{0}^{\theta_{\tau}\tau'_{\rho}(y)  }f( t, x_{\tau+t} )\,dt\leq
\hat N   \sup_{C\in\bC_{1}}
 \|  
f\| _{L_{(q,p)}(C) } , 
\end{equation}
where $\hat N$ depends only on 
$d,\delta ,\rho_{b}$, and $\rho $.
\end{theorem}

Proof.   By Corollary \ref{corollary 10.17.1} we have
$P_{\cF_{\tau}}(\theta_{\tau}\tau'_{\rho}(y) >T)\leq Ne^{-T/N}$
for all $T$ with $N=N(d,\delta,\rho_{b},\rho)$ and there exists
$T=T(d,\delta,\rho_{b},\rho)$ such that   $
Ne^{-T/N}\leq 1/2$.
This shows, by the same argument
as in the proof of \eqref{1.3.3},     
that to prove the current theorem 
it suffices to prove that
\begin{equation}
                                              \label{3.1.4}
E_{\cF_{\tau}}\int_{0}^{T\wedge\theta_{\tau}\tau'_{\rho}(y) }f( t, x_{\tau+t}  )\,dt\leq
\hat N  \sup_{C\in\bC_{1}}
 \|  
f\| _{L_{(q,p)}(C) } ,
\end{equation}
where $\hat N$ depends only on 
$d,\delta ,\rho_{b}$, and $\rho$. Here the left-hand side
is less than
$$
E_{\cF_{\tau}}\int_{0}^{T    }f( t, x_{\tau+t}  )\,dt,
$$
so that \eqref{3.1.4} follows from \eqref{9.7.10}.
The theorem is proved. \qed

Theorem \ref{theorem 9.7.1} allows us to prove It\^o's formula
for functions $u\in W^{1,2}_{(q,p)}(\cO)$, where $\cO$
is a 
\index{$B$@Sets!$W^{1,2}_{(q,p)}(\cO)$}%
domain in $\bR^{d+1}$ and 
$$
W^{1,2}_{(q,p) }(\cO)=\{v: v, \partial_{t}v,
 Dv,  D^{2}v\in L_{(q,p) }(\cO) \}
$$
with norm introduced in a natural way.       
Before, the formula was known only for
(smooth, It\^o, and) $W^{1,2}_{d+1}$-functions
and processes with bounded drifts or
for $W^{2}_{\sfd_{0}}$-functions in case the drift
of the process is dominated by $h(x_{t})$
with $h\in L_{d}(\bR^{d})$ (see \cite{Kr_19}).  

The following extends Theorem 2.10.1  of \cite{Kr_77}
to functions with lower summability of the derivatives
and to spaces with mixed norms.

\begin{theorem}[It\^o's formula]
                   \label{theorem 10.15.1}
Suppose  
\index{$S$@Miscelenea!It\^o's formula}%
 that $(\sfd_{0},q,p)$ are properly tight and    $p<\infty$,
$q<\infty$. Let
$\cO$ be a bounded domain in $\bR^{d+1}$, $0\in \cO$,
$b$ be {\em bounded\/}, and  $u\in W^{1,2}_{(q,p)}(\cO)\cap C(\bar \cO)$. Then,
for $\tau$ defined as the first exit time of $(t,x_{t})$
from $\cO$ with  probability one for all $t\geq0$,
$$
u(t\wedge\tau,x_{t\wedge\tau})
=u(0,0)+\int_{0}^{t\wedge\tau}D_{i}u(s,x_{s})\sigma^{ik}_{s}\,dw^{k}_{s}
$$
\begin{equation}
                                      \label{10.15.01}
+\int_{0}^{t\wedge\tau}[
\partial_{t}u(s,x_{s})+ a^{ij}_{s}D_{ij}u(s,x_{s})
+b^{i}_{s}D_{i}u(s,x_{s})]\,ds
\end{equation}
and the stochastic integral above is a square-integrable
martingale.
\end{theorem}

Proof. First assume that $u$ is smooth and its derivatives
 are bounded. Then  
\eqref{10.15.01} holds by It\^o's formula and, moreover,
by denoting $\tau^{n}= n\wedge\tau$
for any $n\geq0$ we have
$$
E \int_{ \tau^{n}}^{\tau^{n+1}}|Du(s,x_{s})|^{2}\,ds
\leq NE\Big(\int_{\tau^{n}}^{ \tau^{n+1}
}D_{i}u(s,x_{s})\sigma^{ik}_{s}\,dw^{k}_{s}\Big)^{2}
$$
$$
=NE\Big(u( \tau^{n+1},x_{ \tau^{n+1}})-
u( \tau^{n },x_{ \tau^{n }})
$$
$$-
\int_{ \tau^{n}}^{\tau^{n+1}}[
\partial_{t}u(s,x_{s})+ a^{ij}_{s}D_{ij}u(s,x_{s})
+b^{i}_{s}D_{i}u(s,x_{s})]\,ds\Big)^{2}
$$
$$
\leq N\sup_{\bar \cO}|u|^{2}
+NE\Big(\int_{ \tau^{n}}^{\tau^{n+1}}I_{\cO}
\big(|\partial_{t}u|+|Du|+|D^{2}u|\big)(s,x_{s})\,ds\Big)^{2}.
$$
Since $\cO$ is bounded, $\tau$ is bounded as well
and
in light of Theorem \ref{theorem 9.7.1} we conclude that
\begin{equation}
                                      \label{10.15.2}
E \int_{0}^{\tau}|Du(s,x_{s})|^{2}\,ds\leq N
\sup_{\bar \cO}|u|^{2}+N \| \partial_{t}u,
 Du, D^{2}u\|^{2}_{L_{(q,p)}(\cO)},
\end{equation}
where $N$ are independent of $u$ and $\cO$
as long as the size of $\cO$ in the $t$-direction
is under control.
Owing to Fatou's theorem,
this estimate is also true for those $u\in W^{1,2}_{p,q}(\cO)
\cap C(\bar \cO)$
that can be approximated uniformly and in the 
$W^{1,2}_{p,q}(\cO)$-norm by smooth functions with bounded
derivatives (recall that $p<\infty$,
$q<\infty$). For our $u$
there is no guarantee that such approximation is possible.
However, mollifiers do such approximations
 in any subdomain $\cO'\subset \bar \cO'\subset
\cO$ since $d/p+2/q<2$ ($\sfd_{0}>d/2$), so that by embedding theorems $u\in C_{\loc}(\cO)$. Hence, \eqref{10.15.2} holds for our $u$ if we replace $\cO$
by $\cO'$
(containing $(0,0)$). Setting $\cO' \uparrow \cO$  proves \eqref{10.15.2}
in the general  case and proves the last assertion
of the theorem.

After that \eqref{10.15.01} with $\cO'$ in place of $\cO$
is proved by routine approximation of $u$ by smooth 
functions. Setting $\cO' \uparrow \cO$ finally proves \eqref{10.15.01}.
The theorem is proved.\qed

\begin{remark}
                   \label{remark 9.29,1}
We remind the reader that Assumption
\ref{assumption 8.19.2} is supposed
to hold in this chapter, in particular,
in Theorem \ref{theorem 10.15.1}.
However, due to 
Lemma \ref{lemma 3.25.1}, it is 
automatically satisfied if $b$ is bounded
(as in Theorem \ref{theorem 10.15.1})
on account of taking $\rho_{b}$
small enough.
\end{remark}

\begin{remark}
                                \label{remark 3.9.1}
If $b\equiv0$, it turns out
that for any properly tight  $(d_{0},p,q)$, $\rho\in(0,\infty)$,
$x\in \bR^{d}$ and Borel $f(t,x)\geq 0$
\begin{equation}
                             \label{3.9.2}
E\int_{0}^{\tau}f(s,x_{s})\,ds
\leq N(d,\delta)\rho^{2}\dashnorm f\|_{L_{(q,p)}(C_{\rho}(0,x))},
\end{equation}
where $\tau$ is the first exit time of $(s,x_{s})$
from $C_{\rho}(0,x)$.

Indeed, if $\rho=1$, this follows from 
Theorem  \ref{theorem 9.7.1} where we take   $T=1$,
any appropriate $\rho_{b}$ and observe that
$\tau\leq 1$ and we may assume that $f=0$
outside $C_{1}(0,x)$. The case of general $\rho$
is treated by parabolic scaling of $\bR^{d+1}$.
\end{remark}

This simple observation has the following implication
in which
$$
\cL_{0} u(t,x)=\partial_{t}u+(1/2)a^{ij}(t,x)D_{ij}u(t,x) ,
$$
where $a(t,x)$ is a Borel $\bS_{\delta}$-valued
function on $\bR^{d+1}$.
\begin{lemma}
                          \label{lemma 3.9.2}
Suppose   that $(\sfd_{0},q,p)$ are properly tight, $x\in\bR^{d}$,
$\rho\in(0,\infty)$,   $u\in W^{1,2}_{(q,p)}(C_{\rho}(0,x))$
and $u=0 $ on $\partial'C_{\rho}(0,x)$.
Then
\begin{equation}
                             \label{3.9.3}
|u(0,0)| 
\leq N(d,\delta)\rho^{2}\dashnorm (\cL_{0}u)_{-}\|_{L_{(q,p)}(C_{\rho}(0,x))}.
\end{equation}
\end{lemma}

Proof. First recall that, since $\sfd_{0}>d/2$, we have
$d/p+2/q<2$ and $u$ is continuous in 
$\bar C_{\rho}(0,x)$ by embedding theorems. Then
approximate $u$ in $W^{1,2}_{(q,p)}$-norm by smooth
functions $u^{n}$ vanishing on $\partial'C_{\rho}(0,x)$.
By It\^o's formula
$$
u^{n}(0,0)=-E\int_{0}^{\tau}\cL_{0}u^{n}
(s,x_{s})\,ds,
$$
where $x_{s}$ is a solution of $dx_{s}=\sqrt{a(s,x_{s})}\,dw_{s}$ 
with $x_{0}=0$ and $\tau$ is the first exit time of $(s,x_{s})$
from $C_{\rho}(0,x)$.
In light of \eqref{3.9.2} estimate \eqref{3.9.3}
holds with $u^{n}$ in place of $u$. Sending
$n\to\infty$ yields \eqref{3.9.3} as is and
proves the lemma. \qed

Here is an ``elliptic'' version of Theorem
\ref{theorem 10.15.1} proved in the same way
on the basis of the same Theorem  \ref{theorem 9.7.1} with $q=\infty$. For $p\geq d$ 
Theorem \ref{theorem 12.23.1} can be found in \cite{Kr_77}.

\begin{theorem}[It\^o's formula]
                          \label{theorem 12.23.1}
Assume that  $p\in[\sfd_{0},\infty)$. 
\index{$S$@Miscelenea!It\^o's formula}%
  Let
$\cO$ be a bounded domain in $\bR^{d}$, $0\in \cO$,
$b$ be {\em bounded\/}, and  $u\in W^{ 2}_{p }(\cO)\cap C(\bar \cO)$ ($u$ is independent of $t$). Then,
for $\tau$ defined as the first exit time of $ x_{t} $
from $\cO$ with  probability one for all $t\geq0$,
$$
u( x_{t\wedge\tau})
=u(0 )+\int_{0}^{t\wedge\tau}D_{i}u( x_{s})\sigma^{ik}_{s} \,dw^{k}_{s}
$$
$$
+\int_{0}^{t\wedge\tau}[(1/2)
 a^{ij}_{s}D_{ij}u( x_{s})
+b^{i}_{s}D_{i}u( x_{s})]\,ds 
$$
and the stochastic integral above is a square-integrable
martingale.  
\end{theorem}

\mychapter[Regular diffusion
processes]{Regular diffusion
processes}
                      \label{chapter 3}

\def\ineq{$\bar b_{\rho_{b}}\leq \sfb_{0}$\,}

\mysection[Analytic criterion for \protect\ineq]
{Analytic criterion for
\ineq (Assumption \protect\ref {assumption 8.19.2})}   
                \label{section 1.20.1}

  In Lemma \ref{lemma 12.11.1}    
we have already pointed out
an  analytic conditions sufficient
for the inequality $\bar b_{\rho_{b}}\leq \sfb_{0}$
to hold. In this section we present its generalization
 in terms of Morrey spaces. 
We suppose that on $\bR^{d+1}$ we are given Borel   $\bS_{\delta}$-valued $a$ ($\delta\in(0,1])$ and $\bR^{d}$-valued $b$. Set $\sigma=\sqrt a$.
 We follow Section 3 of \cite{Kr_21_1}.

Let $(\Omega,\cF,P)$ be a complete probability
space, let $\cF_{t}, t\geq0$, be an increasing family of
complete $\sigma$-fields $\cF_{t}\subset\cF$,  
and let $w_{t}$ be an $\bR^{d }$-valued Wiener process
relative to $\cF_{t}$.  We will be dealing with the equation (system)
\begin{equation}
                                        \label{11.29.20}
x _{s}=x  +\int_{0}^{s}\sigma (\sft_{r},x_{r})\,dw_{r}
+\int_{0}^{s}b (\sft_{r},x_{r}) \,dr,\quad \sft_{s}=t+s,
\end{equation}
where $(t,x)\in\bR^{d+1}$ are given initial conditions.

We assume that  $\sigma$ and $b$ are smooth
and $b$ is bounded. Then it is well known that the solutions of  system
\eqref{11.29.20}
form a strong Markov process $X$ with trajectories $(\sft_{s},x_{s})$. 

Set   
$$
  \bar b_{R}=\sup_{r\leq R}r^{-1}
\sup_{(t,x)\in\bR^{d+1}}\sup_{C\in\bC_{r}}
E_{t,x}\int_{0}^{\tau_{C}}|b(\sft_{s},x_{s})|\,ds,
$$
where $\tau_{C}$ is 
\index{$S$@Miscelenea!$\bar b_{R}$}%
the first exit time of 
$(\sft_{s},x_{s})$ from $C$.

Also recall that the Fabes-Stroock constant $\sfd_{0}=\sfd_{0}(d,\delta)\in(d/2,d)$
is introduced in Remark \ref{remark 2.7.1} and $\sfb_{0}=\sfb_{0}(d,\delta)$
is introduced in Theorem \ref{theorem 8.2.1}.

\begin{theorem}
                     \label{theorem 9.27.10} 
Assume that 

(i) $\sigma$ and $b$ are smooth
and $b$ is bounded;

(ii)  there is a nonnegative integer $k$ and
there are Borel functions $b_{i}(t,x)$,
$0\leq i\leq k$, 
such that $b =\sum_{i=0}^{k}b_{i} $,  and 
we are given properly tight $(d_{0},q_{i},p_{i})$,
$i\leq k$. Furthermore, for
$$
\hat{b}_{\rho}=\sup_{r\leq\rho}r
\sup_{C\in \bC_{r}}\sum_{i=0}^{k}
\dashnorm b_{i}\|_{L_{(q_{i},p_{i})}(C)},
$$  
and the constant $  N_{0}=N_{0}(d,\delta)$, which is the constant from   Corollary \ref{corollary 3.26.1},
 we have  
\begin{equation}
                 \label{9.27.4}   
N_{0}\hat{b}_{\rho_{b}}< \sfb_{0}   
\end{equation}
for some $\rho_{b}\in(0,\infty)$
(note strict inequality). 
Then
\begin{equation} 
                          \label{9.27.5}
 \bar b_{ \rho_{b}}\leq \sfb_{0}.
\end{equation}
\end{theorem} 
\begin{example}
                \label{example 1.27.1}
One of situations when $\hat{b}_{\rho}$ is finite presents when
$k=1$,
$|b_{0}(t,x)|\leq h_{0}(x)$, $|b_{1}(t,x)|
\leq h_{1}(t)$ and, say $h_{0}(x)\leq c |x|^{-1}$,
where $c$ is sufficiently small, and $h_{1}\in L_{2}(\bR)$.
In that case one can take $p_{0}=d_{0},
q_{0}=\infty$, $p_{1}=\infty,q_{1}=1$.

Indeed, if $|x_{0}|\leq 2r$, then
$$
\dashint_{B_{r}(x_{0})}|x|^{-d_{0}}\,dx
\leq 2^{d}\dashint_{B_{2r} }|x|^{-d_{0}}\,dx=N(d,d_{0})r^{-d_{0}},
$$
and if $|x_{0}|\geq 2r$, then 
$|x|^{-1}\leq r^{-1}$ on $B_{r}(x_{0})$
and $\dashnorm |\cdot|^{-1}\|_{L_{d_{0}}
(B_{r}(x_{0}))}\leq r^{-1}$.

Also    
$$
\dashint_{\!s}^{\,\,\,s+r^{2}}h_{1}(t)\,dt
\leq r^{-1}\Big(\int_{s}^{s+r^{2}}
h_{1}^{2}(t)\,dt\Big)^{1/2}
$$
and the integral here tends to zero
as $r\downarrow 0$ uniformly with respect to $s$.
Therefore, by taking $c$ small enough
 and taking appropriately small $\rho_{b}$ we can satisfy
\eqref{9.27.4} with any given $\hat b>0$.

Bounded $b$ also satisfy \eqref{9.27.4}.
\end{example}

{\bf
Proof of Theorem \ref{theorem 9.27.10}}.  
For $\mu\in[0,\infty)$ denote by $x^{\mu}_{t}$
 the diffusion process corresponding to $\mu b$
in place of $b$ and use the superscript $\mu$ for other
objects related to $x^{\mu}_{t}$. Call a $\mu$
``good'' if  
$$
\bar b^{\mu}_{\rho_{b}}\leq \hat \sfb_{0}=(1/2)(\sfb_{0}+N_{0}\hat b_{\rho_{b}})\quad (<\sfb_{0})
$$
and define $\cM$ as the set of good $\mu$.
 
Our claim is that $1\in \cM$. Observe that $0\in\cM$.
We are going to use the method of continuity proving,
first, that $\cM\cap[0,1]$ is closed and, second,
that $\cM \cap[0,1]$ is open to the right.  Below $\rho\in(0,\rho_{b}]$,  $C\in\bC_{\rho}$,
and $(t,x)\in C$, are arbitrary. 

If $\mu_{n}\in \cM\cap[0,1]$,  $n=1,2,...$,
 converge to $\mu_{0}$,  
then  
  we have 
\begin{equation}   
                                                \label{3.26.41}
E^{\mu_{n}}_{t,x}\int_{0}^{\tau_{C}}
\mu_{n} |
b(\sft_{s},x _{s})|\,ds
\leq   \rho\hat\sfb_{0} ,
\end{equation}
where by writing $E^{\mu_{n}}$ we mean
that the symbol $\mu_{n}$ should
be placed inside the expectation sign
in appropriate positions and
$\tau _{C}$ is the first exit time
of $(\sft_{s},x _{s})$ from $C$.
By using Girsanov's theorem and Fatou's lemma
one easily shows that \eqref{3.26.41} is also true for $n=0$.
But in that case, $  \bar b^{\mu_{0}}_{\rho_{b}}\leq
\hat\sfb_{0}  $ so that, indeed,
$\cM\cap[0,1]$ is closed.

To prove that $\cM$ is open to the right,
 first take $\mu=0$ and
$\varepsilon>0$ and
observe that since $b$ is bounded  
and $\tau^{ \varepsilon}_{C}\leq \rho^{2}$, there is a constant $K$ such that
$$
E^{ \varepsilon}_{t,x}\int_{0}^{\tau_{C}}|  \varepsilon  
b(\sft_{s},x _{s})|\,ds
\leq \varepsilon K\rho.
$$
Hence, for $\varepsilon$ small enough
we have $   \bar b^{\varepsilon}_{\rho_{b}}<\hat\sfb_{0}$,
so that all small $\varepsilon$'s are good.
Next,  take a
$\mu\in \cM \cap(0,1] $,   $\varepsilon>0$,
 and use Girsanov's theorem
to see that
\begin{equation}
                                                \label{3.26.4}
E^{\mu+\varepsilon}_{t,x}\int_{0}^{\tau_{C}}| (\mu+\varepsilon) 
b(\sft_{s},x^{\mu+\varepsilon}_{s})|\,ds
=E^{\mu }_{t,x}e^{\phi(\varepsilon)}\int_{0}^{\tau _{C}}| (\mu+\varepsilon) 
b(\sft_{s},x _{s})|\,ds,
\end{equation}
where  for $\check b=\sigma^{*}(\sigma\sigma^{*})^{-1}b$
$$
\phi(\varepsilon)=\varepsilon\int_{0}^{\rho^{2}_{b}}
 \check b(\sft_{s},x_{s})\,dw_{s}
-(\varepsilon^{2}/2)
\int_{0}^{\rho^{2}_{b}}|\check b(\sft_{s},x_{s})|^{2}\,ds.
$$
 
Recall that $E^{\mu}_{t,x}e^{\phi(\beta\varepsilon)}=1$ for any $\beta$
and observe that for any $\beta>1$
$$
E^{\mu}_{t,x}e^{\beta\phi(\varepsilon)}=
E^{\mu}_{t,x}e^{\phi(\beta\varepsilon)}
\exp\Big((\varepsilon^{2}/2)( \beta^{2} -1)
\int_{0}^{\rho_{b}^{2}}|\check b(\sft_{s},x_{s})|^{2}\,ds\Big)
\leq e^{\varepsilon^{2}\beta^{2}K},
$$ 
where $K$ is a constant independent of $t,x$.  We use this
and H\"older's inequality to obtain from
\eqref{3.26.4} that
$$
E^{\mu+\varepsilon}_{t,x}\int_{0}^{\tau_{C}}| (\mu+\varepsilon) 
b(\sft_{s},x _{s})|\,ds
$$
\begin{equation}
                                                 \label{3.27.4}
\leq e^{\varepsilon^{2}\beta K } 
\Big(E^{\mu}_{t,x}\Big(\int_{0}^{\tau _{C}}| (\mu+\varepsilon) 
b(\sft_{s},x _{s})|\,ds\Big)^{\alpha}\Big)^{1/\alpha},
\end{equation}
where $\alpha=\beta/(\beta-1)$. 
 
 For any $\varepsilon_{1}
>0$
according to  Remark \ref{remark 3.27.1},
for an appropriate choice of $\beta$, the second factor
on the right in \eqref{3.27.4} is less than
$(1+\varepsilon_{1})
 (1+\varepsilon/\mu)   \bar b^{\mu}_{\rho_{b}}\rho$.
Since $\mu$ is good, 
$$
(1+\varepsilon_{1})
 (1+\varepsilon/\mu)   \bar b^{\mu}_{\rho_{b}}\rho\leq (1+\varepsilon_{1})
 (1+\varepsilon/\mu)   \hat \sfb_{0}\rho.
$$

We can choose
$\varepsilon$ and $\varepsilon_{1}$ arbitrarily and make
the left-hand side of \eqref{3.27.4} 
  less than $ \sfb_{0}\rho $ (at this point we use that $\hat \sfb_{0}<\sfb_{0}$).
This shows that $  \bar b^{\mu+\varepsilon}_{\rho_{b}}
<  \sfb_{0} $,
and once we have this, 
$$
  \bar b^{\mu+\varepsilon}_{\rho_{b}}
\leq N_{0}(d,\delta)(\mu+\varepsilon)\hat b_{\rho_{b}}
$$
by Corollary \ref{corollary 3.26.1}.  
It follows that  $\mu+\varepsilon$ is good
for all small enough $\varepsilon>0$ and this
 brings the proof of the theorem to an end. 
\qed

\mysection[Regular diffusion
processes]{Regular diffusion
processes. H\"older
continuity and Harnack inequality for caloric functions} 
                \label{section 11.17,1}

We suppose that on $\bR^{d+1}$ we are given Borel   $\bS_{\delta}$-valued $a$ ($\delta\in(0,1])$ and $\bR^{d}$-valued $b$. Set $\sigma=\sqrt a$.
Define
$$
\cL=\partial_{t}+(1/2)a^{ij}D_{ij}
+b^{i}D_{i}.
$$

Let $\Omega$ be the set of $\bR^{d+1}$-valued
 continuous function $(t_{0}+t,x_{t})$, $t_{0}\in \bR$,
defined for $t\in[0,\infty)$.
For $\omega=\{(t_{0}+t,x_{t}),t\geq0 \}$, define
$\sft_{t}(\omega)=t_{0}+t$, $x_{t}(\omega)=x_{t}$,
and set $ \cN_{t}=\sigma((\sft_{s},x_{s}),s\leq t)$,
$  \cN_{\infty}= \sigma((\sft_{s},x_{s}),s< \infty)$.  

 Let
$$
X=\{(\sft_{t},x_{t}), \cN_{t}, P_{t,x})
$$
be a strong Markov $\bR^{d+1}$-valued process on
$(\Omega,\cN_{\infty})$.

\begin{definition}
           \label{definition 6.24,1}

We say that $X$ is a {\em regular
diffusion process 
\index{$D$@Processes!regular
diffusion process}%
corresponding to
$\cL$ or to $a,b$\/} if

(i)  for any $(t,x)\in\bR^{d+1}$
there exists a $d$-dimensional Wiener process $w_{t}$, $t\geq0$,
which is a Wiener process relative to $\bar \cN_{t}$,
where $\bar \cN_{t}$ is the completion of $\cN_{t}$
with respect to all $P_{s,y}$, and such that with    
$P_{t,x}$-probability one, for
all $s\geq 0$, $\sft_{s}=t+s$ and
\begin{equation} 
                             \label{4.27.010}
x_{s}=x+\int_{0}^{s}\sigma(t+u,x_{u})\,dw_{u}
+\int_{0}^{s}b(t+u,x_{u})\,du.
\end{equation}

(ii) there exists $\rho_{b}\in(0,\infty)$ such that
\begin{equation}
                     \label{7.14,1}
  \bar b_{\rho_{b}}=\sup_{r\leq\rho_{b}}r^{-1}
\sup_{(t,x)\in\bR^{d+1}}\sup_{C\in\bC_{r}}
E_{t,x}\int_{0}^{\tau_{C}}|b(\sft_{s},x_{s})|\,ds \leq \sfb_{0},
\end{equation}
where $\tau_{C}$ is the first exit time of 
$(\sft_{s},x_{s})$ from $C$.

\end{definition}

According to Theorem \ref{theorem 9.27.10} such processes exist and
in this section we assume that we are
given a regular
diffusion process $X$ corresponding to
$\cL$. 

By Theorem \ref{theorem 6.3,1} we have
the following.

\begin{theorem}
                \label{theorem 6.26,5}
For any $\rho\leq\rho_{b}$ and nonnegative Borel $f$ we have
$$
E_{0,0}\int_{0}^{\tau_{R}}
f(\sft_{s},x_{s})\,ds\leq N(d,\delta)
\rho^{2}\dashnorm f\|_{L_{d+1}(C_{\rho})}.
$$
\end{theorem}

The requirement (ii) of Definition
\ref{definition 6.24,1} implies that
Assumption \ref{assumption 8.19.2}, the standing assumption in Chapter \ref{chapter 10.20.1} after Theorem \ref{theorem 8.2.1} is proved,
is satisfied with respect to
any measure $P_{t,x}$. Therefore,
all results in Chapter \ref{chapter 10.20.1} after Theorem \ref{theorem 8.2.1}
are applicable to the solutions
of \eqref{4.27.010} as trajectories of $X$
(there could be other solutions,
but we are talking only about those that
are part of $X$). Of course, one should
make appropriate adjustments in these results since the processes in 
Chapter \ref{chapter 10.20.1} starts
from the origin unlike solutions
of \eqref{4.27.010}.

Corollary \ref{corollary 10.26.1} implies the following.

\begin{theorem}
               \label{theorem 6.24,1}
For any $n>0$, $(t,x)$, and
  $s\geq 0$  
\begin{equation}
                                  \label{6.24,3}
E_{t,x}\sup_{r\in[0,t]}|x_{ r}-x_{0}|^{ n}
\leq N(  t ^{ n/2}+ t ^{ n}),
\end{equation}
where $N=N(n, \rho_{b},d,\delta)$.

\end{theorem}

 \begin{definition}
                                \label{definition 12.27.1}
If $Q$ is a set in $\bR^{d+1}=\{(t,x)
:t\in\bR,x\in\bR^{d}\}$
and $u$ is a bounded Borel function on $Q$,
we call it {\em caloric\/} 
\index{$S$@Miscelenea!caloric function}%
(in $Q$ relative to the process $X$) if
for any $(s,y)$ and $T,R\in(0,\infty)$
such that   $\bar C_{T,R}(s,y)\subset Q$
and any $(t_{0},x_{0})\in C:=C_{T,R}(s,y)$ we have
$$
u(t_{0},x_{0})=E_{t_{0},x_{0}} u(\sft_{\tau_{C}},x_{\tau_{C }}),
$$
where $\tau_{C}$ is the first exit time of $(t_{0}+
t,x_{t})$
from $C$.
  
\end{definition}

\begin{example}
                \label{example 6.24,1}
Let $f$ be a bounded Borel function
on $\bR^{d}$, $T\in\bR$, and for $t\leq T$ introduce
$$
u(t,x)=E_{t,x}f(x_{T-t}).
$$

Then $u$ is a caloric function
in $Q:=(-\infty,T]\times\bR^{d}$.
To show this it suffices to concentrate
on continuous $f$. Then  take $(t,x)\in \bR^{d+1}$
such that $t\leq T$ and, for $\varepsilon>0$, define $g_{\varepsilon}=\varepsilon^{-1}I_{(T,T+\varepsilon)}$. Observe that
for any $(t,x)\in Q$ owing to the continuity of $x_{s}$ and~$f$
$$
v_{\varepsilon}(t,x):=
E_{t,x}\int_{0}^{\infty}g_{\varepsilon}(\sft_{s},x_{s})f(x_{s})
e^{-s}\,ds\to e^{-(T-t)}u(t,x)
$$
as $\varepsilon\downarrow0$. On the other hand, by the strong Markov property, for $t_{0},x_{0}, C$ as 
in Definition \ref{definition 12.27.1}
$$
v_{\varepsilon}(t_{0},x_{0})=
E_{t_{0},x_{0}}\int_{0}^{\tau_{C}}g_{\varepsilon}(\sft_{s},x_{s})f(x_{s})
e^{-s}\,ds
+E_{t_{0},x_{0}}e^{-\tau_{C}}
v_{\varepsilon}(\sft_{\tau_{C}},x_{\tau_{C}}).
$$
Here the first term on the right is zero and letting $\varepsilon\downarrow0$ we get the desired result.
\end{example}
 
Here is a version of Theorem 
\ref{theorem 9.5.1} improving Theorem
\ref{theorem 6.26,5}.
 
\begin{theorem}
                \label{theorem 6.26,6}
Suppose   that $(\sfd_{0},q,p)$ are properly tight.
 Then 
for any $\rho\leq \rho_{b}$,  Borel nonnegative $f$ 
given on $\bR^{d+1}$ 
we have  that
$$
E_{0,0}\int_{0}^{ \tau_{\rho}  }f(\sft_{s},  x_{s} )\,ds\leq
N(d,\delta,p,q)\rho^{2}\dashnorm f\|_{L_{(q,p)}(C_{\rho})},
$$
$$
E_{0,0}\int_{0}^{ \tau'_{\rho}  }f(\sft_{s},  x_{s} )\,ds\leq
N(d,\delta,p,q)\rho^{2}
\sup_{C\in\bC_{\rho}}
 \dashnorm  
f\| _{L_{(q,p)}(C) }.
$$

\end{theorem}

Next,
we deal with the H\"older norm estimates for 
caloric functions and potentials.
If $z_{1}=(t_{1},  x_{1})$ and $z_{2}=(t_{2}, x_{2})$,
we define    
\begin{equation}
                                                   \label{eq:4.2.5}
\rho(z_{1}, z_{2})=|x_{1}-x_{2}|+|t_{1}-t_{2}|^{1/2}
\end{equation}
 and call  $\rho(z_{1},z_{2})$   the parabolic
distance between $z_{1}$ and $z_{2}$.
The PDEs versions of Lemma 
\ref{lem:4.2.2} and Theorem \ref{thm:4.2.1} below belong
to Krylov-Safonov (\cite{KS_80}).
The proofs below are based on a
probabilistic adaptation of the PDE arguments
from \cite{KS_80}.

\begin{lemma} 
                                        \label{lem:4.2.2} 
Let $ R\leq  \rho_{b} $ and let $u$ be a caloric
function in $\bar C_{ R}$. Then there
exist constants $N$ 
and   
$$
\alpha_{0}\in(0,1),
$$ 
depending only on $\delta,d $,
  such that, for any $\alpha\in(0,\alpha_{0}]$ and 
$z_{1},z_{2}\in C_{ R/2 }$, we have   
\begin{equation}
                                       \label{eq:4.2.6}
\big|u(z_{1})-u(z_{2})\big|\le 
NR^{-\alpha}\rho^{\alpha}(z_{1},z_{2})
\sup\big(|u|,\bar C_{  R}\big).
\end{equation}

Furthermore, $\sup(|u|,\bar C_{ R})$ in 
\eqref{eq:4.2.6} can be replaced by $\osc(u,\bar 
C_{ R})$, where we use the 
notation
$$
 \osc(g,\Gamma) =\osc_{\Gamma}g=\sup_{\Gamma}g-
\inf_{\Gamma}g.
$$

\end{lemma}

Proof.  We use the classical arguments
of E. De Giorgi. For $r\leq R$  set
$$
w(r)=\osc(u,\bar C_{r}),\quad 
m(r)=\inf_{\bar C_{r}}u,\quad M(r)=\sup_{\bar C_{r}}u,
$$
$$
\mu(r)=(1/2)
\big(m(r)+M(r)\big).
$$

  Take  $r\le R/2$ and  suppose that
$$
\big|C_{ 2r}\cap\big\{ u\le
\mu(r)\big\} 
\big|\ge (1/2)|C_{ 2r}|.  
$$
Then there is a closed $\Gamma\subset
C_{ 2r}\cap\big\{ u\le
\mu(r)\big\}$ such that
\begin{equation}
                                   \label{eq:4.2.7}
\big|C_{3 r^{2},2r}(r^{2},0)\cap\Gamma 
\big|\ge (1/4)|C_{3 r^{2},2r}|
\end{equation}

By Theorem \ref{thm:4.1.10} (with $s=0$) for any $(t ,x )  
\in \bar C_{r}$ we have
$$
P _{t,x}(\tau_{\Gamma} <\tau_{2r} )\geq \pi_{0},
$$
where $\pi_{0}>0$ depends only on $\delta$ and $d $,  
$\tau_{\Gamma} $ is the first time $(\sft_{s},x_{s})$
hits $\Gamma$, $\tau_{2r}$ is its first exit
time from $C_{2r}$.
Then by definition and the strong Markov property
for $\tau=\tau_{\Gamma}\wedge\tau_{2r}$ we have
$$
u(t ,x )= 
E_{t,x }  u(t +\tau_{2r} ,x_{\tau_{2r} })
$$
$$
=E_{t,x}  u(t+\tau _{2r},x_{\tau_{2r} })I_{\tau_{\Gamma}<\tau_{2r}}
+E_{t,x}  u(t+\tau _{2r},x_{\tau _{2r}})I_{\tau_{\Gamma}\geq\tau_{2r}}
$$
$$
=E_{t,x}  u(t+\tau _{\Gamma},x_{\tau_{\Gamma} })I_{\tau_{\Gamma}<\tau_{2r}}
+E_{t,x}  u(t+\tau_{2r} ,x_{\tau_{2r} })I_{\tau_{\Gamma}\geq\tau_{2r}}
$$
$$
\leq \mu(r)\pi_{0}+M(2r)(1-\pi_{0})
$$
(we used that $\mu(r)\leq M(2r)$).
It follows  that     
$$
M(r)\le \pi_{0}\frac{1}{2}\big(m(r)
+M(r)\big)+(1-\pi_{0})M(2r),
$$
$$
\big(1-\frac{\pi_{0}}{2})M(r)\leq \frac{\pi_{0}}{2} m(r)+(1-\pi_{0})M(2r).
$$

Adding  to this   the obvious inequality 
$$
\big(\frac{\pi_{0}}{2}-1)m(r)\leq -\frac{\pi_{0}}{2} m(r)
+(\pi_{0}-1)m(2r),
$$
we  get  
\begin{equation} 
                                             \label{eq:4.2.8}
\big(1-\frac{\pi_{0}}{2}\big)w(r)\le(1-\pi_{0})w(2r),\quad
w(r)\le\varepsilon w (2r),
\end{equation} 
where $\varepsilon<1$,  $\varepsilon=
\varepsilon(\pi_{0})$.  We may, certainly, assume that 
 $\varepsilon>1/2$.

 We have proved (\ref{eq:4.2.8}) assuming 
that  (\ref{eq:4.2.7}) is true. However if (\ref{eq:4.2.7}) 
is false, then  
$-u$ satisfies an inequality similar to \eqref{eq:4.2.7}
and  this leads to
 (\ref{eq:4.2.8}) again.

\label{iteration}
 Therefore,  $w(r)\le\varepsilon w(2r)$ for all $r\le R/2$.  
Iterations then yield    \smallskip
$$
w(r)\le\varepsilon^{2}w(4r)\quad\text{for}\quad r\le R/4,..., 
w(r)\le\varepsilon^{n}w(2^{n}r)  \quad\text{for}\quad r\le2 ^{-n}R.\smallskip
$$
 If $r\le R/2$ and  we take  $n:=\lfloor\log_{2}(R/r)\rfloor$,
 then  $r\le2^{-n}R$ and   \smallskip
$$
w(r)\le\varepsilon^{n}w(2^{n}r)\le
\varepsilon^{-1}(r/R)^{\alpha}w(R)\le2
\varepsilon^{-1}(r/R)^{\alpha}
\sup\big(|u|,\bar C_{R }\big),\smallskip
$$
where $\alpha=-\log_{2}\varepsilon\in(0,1)$.  
This provides an estimate of 
the oscillation of $u$ in any   $C_{r}$
with $r\le R/2$. The same  estimate
obviously holds for the oscillation of $u$ in  any   
  $ C_{ r}(t,x)\subset C_{  R}$ as long as  $r\le R/2$
and $(t,x)\in C_{R/2}$.

Now  take  $z_{1}=(t_{1},x_{1}),z_{2}=(t_{2},x_{2})\in C_{R/2}$
 such that
$r:=\rho(z_{1},z_{2})\le R/2$ and define 
$$
t=t_{1}\wedge t_{2},\quad x= (x_{1}+x_{2})/2.
  $$
 Then  we have $z_{i}\in\bar C_{ r}(t,x)$, $i=1,2$,
and     
$$
\big|u(z_{1})-u(z_{2})\big|  \le    2\varepsilon^{-1}
\rho^{\alpha}(z_{1},z_{2})R^{-\alpha}\sup\big(|u|,\bar C_{  R}\big).
$$

 In  the case that $\rho(z_{1},z_{2})\geq R/2$  
 we have 
\begin{align*}
\big|u(z_{1})-u(z_{2})\big|   \le &\, 2
 \sup\big(|u|, \bar C_{2} \big)\\ 
  \le &\, 2^{1+\alpha}\rho^{\alpha}(z_{1},z_{2})R^{-\alpha}
\sup\big(|u|,\bar C_{ 2}\big).
\end{align*}

Thus,   $N=2^{1+\alpha}+2 
\varepsilon^{-1} $ in (\ref{eq:4.2.6})  is always
 a good choice
with   $\alpha $ found above. One can take any smaller
$\alpha$ as well since $\rho(z_{1},z_{2})\leq N(d)R$.
The lemma is proved.   \qed

\begin{remark}[Liouville theorem]
                      \label{remark 3.9,1}
Letting $R\to\infty$ in \eqref{eq:4.2.6} we see that if \eqref{7.14,1} holds for any $\rho_{b}>0$
and $u$ is caloric and bounded in $\bR^{d+1}_{0}$,
then $u$ is constant. 

\end{remark}

\begin{corollary}
                       \label{corollary 9.5.1}
The process $X$ is strong Feller in the sense that
for any Borel bounded $f(x)$ and $T\in\bR$ the function
$$
u(t,x)=E_{t,x}f(x_{T-t})
$$
is a (H\"older) continuous function of $(t,x)\in(-\infty, T)
\times\bR^{d}$. As a further standard
consequence of this,
the process $((\sft_{t},x_{t}),\bar\cN_{t+},P_{t,x})$ is strong Markov.

\end{corollary}

The importance of the fact that
$((\sft_{t},x_{t}),  \cN _{t+},P_{t,x})$ is (even just)
  a Markov process is well
seen from the following $0-1$ law
\index{$S$@Miscelenea!Blumenthal's $0-1$ law}%
of Blumenthal.

\begin{theorem}
             \label{theorem 8,16.1}
Let $((\sft_{t},x_{t}),  \cN _{t+},P_{t,x})$
be a Markov process. Then for any
$A\in  \cN_{0+}$ and $(t,x)\in\bR^{d+1}$ we have $P_{t,x}(A)=
P^{2}_{t,x}(A)$, that is $P_{t,x}(A)=0$
or $1$.
\end{theorem}

Indeed, by definition   
$$
P_{t,x}(A\cap A)
=E_{t,x}I_{A}P_{\cN_{0+}}(A)=E_{t,x}I_{A}P_{(\sft_{0},x_{0})}(A)
$$
$$
=E_{t,x}I_{A}P_{t,x}(A)=P^{2}_{t,x}(A).
$$
Here is a surprising albeit very well-known
corollary of Theorem \ref{theorem 8,16.1}.

\begin{corollary}
           \label{corollary 8,16.4}
Let $w_{t}$ be  a one-dimensional
Wiener process. Then
\begin{equation}
                                                           \label{11}
\nlimsup_{t \downarrow 0} \frac{ w_t}{\sqrt t} = \infty
\quad \text{(a.s.)}.  
\end{equation}
\end{corollary}

Indeed
$$
\xi := \nlimsup_{t \downarrow 0} \frac{w_t}{\sqrt t} =
\lim_{t \downarrow 0} \sup_{s \in (0,t)} \frac{w_s}{\sqrt s},
 $$
where the supremum can be confined to rational $s \in (0,t)$.
We see that this supremum is $\cF^{w}_t $-measurable ($\cF^{w}_t $ is the 
completion of the $\sigma$-field $\sigma(w_{s},s\leq t)$), it decreases
with decreasing $t$ and the limit can be taken over the sequence 
$t = \frac 1n$. Furthermore, $\xi$ is $\cF_t^w$-measurable for any 
$t > 0$, that is, $\cF_{0+}^w$-measurable. In particular, 
$P(\xi\in B)$ is zero or one for any Borel set $B$.
Next use that $I_{(-\infty, n)}(x)$ is a right-continuous function
of $x$ for fixed $n$ and $\frac 1{\sqrt t}w_t $ is a normal $(0,1)$
variable.
Then 
\begin{multline}
P ( \xi < n ) = E \, I_{(-\infty, n)} (\xi) 
= \lim_{t \downarrow 0} P \Big\{ \sup_{s \in (0,t)} \frac 1{\sqrt s} w_s
< n \Big\} \\ 
\le \nliminf_{t \downarrow 0} P \Big\{ \frac 1{\sqrt t} w_t < n \Big\}
= \frac 1{\sqrt{2 \pi}} \int\limits_{-\infty}^n e^{-\frac 12 x^2} \, dx < 1.
\end{multline}

By the $0-1$ law, $P(\xi < n ) = 0$, $\xi \ge n$
(a.s.) for any constant $n$ and $\xi = \infty$ (a.s.) indeed.

  Observe that \eqref{11} is also true
for $-w_{t}$ in place of $w_{t}$ and this  shows that, in an 
arbitrarily
small time interval $[0,t]$, the sample path of the Wiener process
passes through the origin infinitely many times.

Here is the statement of the Harnack inequality.

\begin{theorem}
                     \label{thm:4.2.1} 
Let   $R\leq  \rho_{b}$, and let 
$u$ be a nonnegative 
caloric function in $\bar C_{2 R^{2},R}$.
 Then there exists a constant $N$, which depends 
only on $ \delta,d $, such that   
$$
u(R^{2},0)\le Nu(0,x) 
$$
whenever $|x|\le R/2$. 
\end{theorem}

Proof. We basically repeat the proof of  
Theorem 6.1 in  \cite{Kr_20} based on an idea
of E.M. Landis and techniques from \cite{KS_80} 
 and, to exclude a trivial situation, additionally assume that 
$$
u(R^{2},0)>0.
$$
 
  For  
$\kappa=1/2,\eta=1/2$, we  take  $N$ and $\nu$ from
Theorem \ref{theorem 12.7.2}, call  $N_{1}$ this $N$,
 and,  having in mind
    Theorem \ref{theorem 11.8.1},
find $\gamma\in(0,1)$    close to 1 and $\varepsilon>0$
close to zero that
\begin{equation}
                                       \label{eq:4.2.3}
1-\varepsilon\geq q(\gamma)2^{-1}+\big[1-q(\gamma)  \big]2^{\nu}.
\end{equation}

  Next, for $r\in[0,R)$, introduce  
$$\mu(r)=u(R^{2},\,0)(1-r/R)^{-\nu},\quad 
n(r)=\sup\{ u,\bar C_{ r}(R^{2},0)\},
$$
$(n(0)=u(R^{2},0))$
  and define
  $r_{0}$ as the greatest number in $r\in [0,R)$ satisfying 
$$
n(r)= \mu(r).
$$ 
 Such a number does exist because
  $n(0)=\mu(0)$, $\mu(r)\to\infty$
as $r\uparrow R$, and $n(r)$ is bounded, increasing,
and (H\"older) continuous.  
Choose $(t^{0},x^{ 0 })
\in \bar C_{ r_{0}}(R^{2},0)$
such that $n(r_{0})=
u(t^{ 0},x^{0})$ and consider  
the cylinder  
$$
Q:=\Big\{ (t,x)\,:\,0\leq t-t^{ 0}
<\frac{(R-r_{0})^{2}}{4},
\quad|x-x^{ 0}|<\frac{R-r_{0}}{2}\Big\} .
$$

  As is easy to see   $\bar Q\subset \bar C_{ r_{1}}(R^{2},0)$,
where $r_{1}=(R+r_{0})/2 $.  By the definition of $r_{0}$,
this implies   that    
$$
\sup_{\bar Q}u<\mu(r_{1})=u(R^{2},0)\Big(\frac{R-r_{0}}{2R}\Big)^{-\nu}
\leq 2^{\nu}n(r_{0}).
$$

We claim that  owing to this and 
 (\ref{eq:4.2.3}), 
\begin{equation}
                                                       \label{eq:4.2.4}
\big|Q\cap\big\{ u>n(r_{0})/2\big\}\big|\ge  (1- \gamma ) |Q|.
\end{equation}

 To argue by contradiction, assume \eqref{eq:4.2.4}
is false. Then
$$
\big|Q\cap\big\{ u\leq n(r_{0})/2\big\}\big|
>    \gamma   |Q|
$$
and there is a closed set $\Gamma\subset
Q\cap\big\{ u\leq n(r_{0})/2\big\}$ such that
$|\Gamma|> \gamma  |Q|$. 
Introduce $\tau_{\Gamma}$ as the first time the process
$(\sft_{s},x_{s})$ hits $\Gamma$
and $\tau_{Q}$ as the first time it exits from $Q$.
It follows by definition, the strong Markov property
as in the proof of Lemma \ref{lem:4.2.2}, and from 
Theorem \ref{theorem 11.8.1}  that
(note that $n(r_{0})/2\leq\sup_{\bar Q}u$)
$$
u(t^{0},x^{0})=E_{t^{0},x^{0}}
I_{\tau_{\Gamma}<
\tau_{Q}}u(t^{0}+\tau_{\Gamma},
x_{\tau_{\Gamma}})
+E_{t^{0},x^{0}}I_{\tau_{\Gamma}
\geq\tau_{Q}}u(t^{0}+\tau_{Q},
x_{\tau_{Q}})
$$
$$
\leq P_{t^{0},x^{0}}
(\tau_{\Gamma}
<\tau_{Q})n(r_{0})/2+
(1-P_{t^{0},x^{0}}
(\tau_{\Gamma}
<\tau_{Q}))\sup_{\bar Q}u
$$
$$
\leq q(\gamma)n(r_{0})/2+
(1-q(\gamma))\sup_{\bar Q}u
$$
$$
\leq q(\gamma)n(r_{0})/2+
(1-q(\gamma))2^{\nu}n(r_{0}).
$$
Owing to \eqref{eq:4.2.3}
we  now have
$$
n(r_{0})\leq(1+\varepsilon)n(r_{0})
\big[q(\gamma)2^{-1}+(1-q(\gamma))2^{\nu}\big]
\leq (1-\varepsilon^{2})  n(r_{0}),
$$
which is impossible. This proves 
 (\ref{eq:4.2.4}).

  Next we  apply Theorem \ref{thm:4.1.10}
and get that 
$$
u(t^{0},x)\ge \pi_{0}n(r_{0})2^{-1}
$$
if $|x-x^{0}|\le (R-r_{0})4^{-1}$, where  
$\pi_{0}=\pi_{0}( d,\delta )>0$.   
 After that it only remains to apply  
Theorem \ref{theorem 12.7.2}
to conclude that for $|x|\le R/2$ we have    
$$
u(0,x)\ge\frac{1}{2}\pi_{0}n(r_{0})N_{1}^{-1}
\Big(\frac{R-r_{0}}{4}\Big)^{\nu}= 2^{-2\nu-1}\pi_{0}
N_{1}^{-1} u(4,0).   
$$
 The theorem is proved.   \qed

\begin{remark}
                             \label{remark 9.4.1}
If $ \bar b_{\rho_{b}}   
< \sfb_{0}$ for any $\rho_{b}$, we can take $\rho_{b}$
as large as we wish and then 
the one sided Liouville theorem is available: 
If $u\geq 0$ is  caloric in $\bR^{d+1}  $ and {\em independent of $t$ (harmonic)\/},
then $u$ is constant. Indeed, in this case
$u( x)\leq Nu( 0)$ for any $x\in\bR^{d } $,
so that $u$ is bounded and Remark \ref{remark 3.9,1}
is applicable.
\end{remark}

\mysection[Further results]{Further results.
Viscosity solutions} 

We work in the setting of Section
\ref{section 11.17,1} and as there suppose that we are
given a regular
diffusion process $X$ corresponding to
$\cL$.

By using   Lemma \ref{lem:4.2.2} and Theorem \ref{theorem 6.26,6}
one derives in three lines the following analog
of Theorem 6.5 of \cite{Kr_20}.
 \begin{theorem}
                                     \label{theorem 10.8.1}
Let $(\sfd_{0},q,p)$ be properly tight.
Let $ R\leq  \rho_{b} $ and
 let $g$ be a Borel bounded
function on $\bar C_{ R}$ and $f\in L_{(q,p)}(C_{ R})$.
For $(t ,x )\in C_{ R}$ define
\begin{equation}  
                                                  \label{10.19.2}
u(t ,x )=E_{t ,x }
\int_{0}^{\tau_{ R}}f(t +s,x_{s})\,ds+
E_{t  ,x }g(t +\tau_{ R},x_{\tau_{ R}})
\end{equation}
\($\tau_{ R} $ is the first exit time of 
$(\sft_{s},x_{s})$ from $C_{ R}$\).
  Then there exists a constant $N$, which depends
only on $\delta $  and $d$, such that   
\begin{equation}
                                                  \label{eq:4.2.9}
\big|u(z_{1})-u(z_{2})\big|\le N\big(R^{-\alpha}
\rho^{\alpha}(z_{1},z_{2})\sup_{\bar C_{ R}}|g|
+R^{2}\dashnorm f\|_{L_{(q,p)}(C_{ R})}\big)
 \end{equation}
for $z_{1}$, $z_{2}\in  C_{ R/2}$,
$\alpha\in(0,\alpha_{0}]$, where $\alpha_{0}$
is taken from Lemma \ref{lem:4.2.2}.
 
\end{theorem}

Proof.   Observe that $h(t ,x ):
=E_{t,x }g(t +\tau_{ R},x_{\tau_{R}})$
is a caloric function, to which
Lemma \ref{lem:4.2.2} is applicable, and
$u(t ,x )-h(t ,x )$ admits the remaining estimate 
in light of Theorem  \ref{theorem 6.26,6}.
The theorem is proved.  \qed

Here is a  version
of Theorem \ref{theorem 10.8.1} which sometimes
is slightly more convenient.
\begin{theorem}
                                    \label{theorem 12,14.2}
Under the conditions and notation from
 Theorem \ref{theorem 10.8.1}
there exists a constant $N$, which depends
only on $\delta,d $,  
 such that    
\begin{equation}
                                           \label{12,14.6}
\big|u(z_{1})-u(z_{2})\big|\le NR^{-\beta}
\rho^{\beta}(z_{1},z_{2})
\big(\sup_{\bar C_{  R}}|u|
+R^{2}\dashnorm f\|_{L_{(q,p)}(C_{ R})}\big)
\end{equation}
for $z_{1}$, $z_{2}\in  C_{R/2}$, where
$$
\beta=\frac{\alpha_{0}\nu }{\alpha_{0} +\nu},
\quad \nu:=2-\frac{d}{p}-\frac{2}{q}.
$$

\end{theorem}

Proof. Fix $z_{1}$, $z_{2}\in  C_{R/2}$.
Since  there is the sup norm of $u$
on the right, it suffices to prove \eqref{12,14.6}  
assuming that 
$$
\xi:=\Big(\frac{R}{\rho(z_{1},z_{2})}\Big)^{\beta/\alpha_{0}}>4.  
$$
Then set
$$
\bar R=\xi \rho(z_{1},z_{2}).
$$

If $z_{i}=(t_{i},x_{i})$, $i=1,2$, without losing generality
we may assume that $t_{1}\leq t_{2}$. Then  for
\begin{equation}
                                         \label{6,4,1}
|x_{1}|+\bar R/2\leq  R\quad \text{and}\quad
 t_{1}+\bar R^{2}/4\leq  R^{2}
\end{equation}
we have
\begin{equation}
                                         \label{12,14.7}
z_{1},z_{2}\in\bar C_{\bar R/4}(z_{1}) \subset  
C_{\bar R/2}(z_{1})
\subset C_{ R}.
\end{equation}
 Since $z_{1}\in\bar C_{R/2}$,
we have $|x_{1}|\leq R/2$ and $t_{1}\leq R^{2}/4$
and, for any of the inequalities \eqref{6,4,1} to go wrong,
we have to have $\bar R> R $, that is,    
$$   
\Big(\frac{R}{\rho(z_{1},z_{2})}\Big)^{\beta/\alpha_{0}-1}>1,\quad \rho(z_{1},z_{2})> R.
$$
The latter is impossible for $z_{1},z_{2}\in C_{R/2}$. Therefore,
  we assume \eqref{12,14.7} and that 
$\bar R\leq R $.

 Then by Theorem \ref{theorem 10.8.1}
applied to $C_{\bar R/2}(z_{1})$ in place of $C_{R}$
we obtain   
$$
\big|u(z_{1})-u(z_{2})\big|\le N\big(\bar R^{-\alpha_{0}}
\rho^{\alpha_{0}}(z_{1},z_{2})\sup_{\bar C_{ R}}|u|
+\bar R^{\nu}\|f\|_{L_{(q,p)}(C_{ R})}\big),
$$
where the right-hand side is transformed 
to that of \eqref{12,14.6} by simple arithmetics.
 The theorem is proved.  
\qed

\begin{theorem}
                 \label{theorem 11.17,2}
Let $\cO$ be a bounded domain in $\bR^{d+1}$, $f\in L_{(q,p)}(\cO)$
where  $p,q$ are finite and $(\sfd_{0},q,p)$
is properly tight,
and let $g$ be Borel bounded.
Assume that $b$ is {\em bounded\/}.
Then
$$
u(t,x):=E_{t,x}\int_{0}^{\tau_{\cO}}
f(t+s,x_{s})\,ds+E_{t,x}g(t+\tau_{\cO},
x_{\tau_{\cO}}),
$$
which is called a {\em probabilistic solution\/}
\index{$S$@Miscelenea!probabilistic solution}%
of
\begin{equation}
                        \label{11.18,1}
\cL u+f=0\quad\text{in}\quad \cO
\end{equation}
with the boundary data $u=g$ on $\partial'\cO$,
is a $W^{1,2}_{(q,p)}$-viscosity solution 
\index{$S$@Miscelenea!viscosity solution}%
of \eqref{11.18,1}
(the definition of viscosity solutions
will be clear from the proof).
\end{theorem}

Proof. Let $\phi\in W^{1,2}_{(p,q),\loc}(\cO)$
and $(t_{0},x_{0})\in\cO$ be such that
$u-\phi$ has a local maximum at 
$(t_{0},x_{0})$. Set $M=(u-\phi)
(t_{0},x_{0})$ and  introduce $\psi_{\varepsilon}(t,x)=\varepsilon(t-t_{0}+|x-x_{0}|^{2})$. Then for  
$\varepsilon>0$  and sufficiently
small $\rho$, by the strong Markov
property and formula It\^o we have
$$
M=
E_{t_{0},x_{0}}\big(u-\phi-\psi)
(\tau_{\rho}(t_{0},x_{0}),x_{\tau_{\rho}(t_{0},x_{0})}) 
$$
$$
+E_{t_{0},x_{0}}\int_{0}^{\tau_{\rho}(t_{0},x_{0}}\big(f+\cL(\phi+\varepsilon\psi)\big)(
t_{0}+t,x_{t})\,dt.
$$
By taking into account that $u-\phi-\varepsilon\psi\leq M
-\varepsilon\rho^{2}$ on the parabolic boundary
of $C_{\rho} (t_{0},x_{0})$
and applying Theorem \ref{theorem 6.26,6}
we conclude
$$
\varepsilon\leq N\dashnorm \big
(f+\cL(\phi+\varepsilon\psi)\big)_{+}
\|_{L_{q,p}(C_{\rho} (t_{0},x_{0}))},
$$
which implies that
$$
\esssup_{C_{\rho} (t_{0},x_{0})}
(f+\cL(\phi+\varepsilon\psi))>0,\quad
\esssup_{C_{\rho} (t_{0},x_{0})}
(f+\cL\phi)>-
\varepsilon \esssup_{C_{\rho} (t_{0},x_{0})}\cL \psi.
$$
Finally, the arbitrariness of $\rho,\varepsilon$ leads to
$$
\lim_{\rho\downarrow0}\esssup_{C_{\rho} (t_{0},x_{0})}
(f+\cL\phi)\geq0,
$$
and we have thus proved that $u$ is 
a $W^{1,2}_{(q,p)}$-viscosity subsolution of \eqref{11.18,1}. Similarly
one proves that $u$ is 
a $W^{1,2}_{(q,p)}$-viscosity supersolution of \eqref{11.18,1} and this proves the theorem. \qed

 Similarly to Theorem \ref{theorem 10.8.1}
one proves the following
on the basis of 
Theorem~\ref{thm:4.2.1}.

\begin{theorem}
                  \label{theorem 1.18.13}
Under the assumptions of Theorem
\ref{theorem 10.8.1} suppose that
$g\geq 0$. Then there exists
a constant $N$, depending only on $d,\delta$, such that 
$$
u(R^{2}/2,0)\le Nu(0,x) +N
R^{2}\dashnorm f\|_{L_{(q,p)}(C_{ R})}
$$
whenever $|x|\le R/2$. 

\end{theorem}

Next, we mention only a few corollaries
of the results proved in Chapter \ref{chapter 10.20.1}.
All other results admit their versions
in our situation of regular diffusion processes
as well, but it would be probably unnecessary and, certainly, boring to formulate them.

The following is a reformulation
of Corollary \ref{corollary 10.11.1}
for our case of regular diffusion processes.

\begin{theorem}
                 \label{theorem 6.26,1}

For any 
$\kappa\in(0,1)$ there exists  
$N$, depending only on $\kappa,d,\delta$,
   such that, for any $R\leq  \rho_{b}$,
$x\in   B_{\kappa R}$, and closed set
$\Gamma\subset  C_{R}(R^{2},0)$, the $P_{0,x}$-probability that the process 
$(\sft_{s}, x_{s})$  
reaches $\Gamma$ before exiting from $C_{2R^{2},R}$
is greater than or equal to $N^{-1} (|\Gamma|/|C_{R}|)^{\gamma-1/
(d +1)}$:
\begin{equation}
                               \label{10.2.10}
P_{0,x} (\tau_{\Gamma}  <\tau_{2R^{2},R}  )
\geq N^{-1} (|\Gamma|/|C_{R}|)^{\gamma-1/(d +1)},
\end{equation}
where $\tau_{\Gamma} $ is the first time $(\sft_{s}, x_{s})$
hits $\Gamma$, $\tau_{2R^{2},R} $ is the first exit time of
$(\sft_{s}, x_{s})$ from $C_{2R^{2},R}$,
   and $\gamma$ is taken from Theorem \ref{theorem 12.21.1}.
\end{theorem}

Here is a version of Corollary
\ref{corollary 10.1,1}.
 
\begin{theorem}
               \label{theorem 6.26,3}    
For any $R\leq  \rho_{b}$,
$\kappa\in(0,1)$, Borel nonnegative $f$ 
vanishing outside $C_{R}(R^{2},0)$, and $x\in B_{\kappa R}$
$$
\int_{C_{R}(R^{2},0)}f^{1/(2\gamma)}(t,y)\,dydt\leq NR^{d+2-1/\gamma}
\Big(E_{0,x} \int_{0}^{\tau_{2R^{2},R}  }f(t, x_{t})\,dt\Big)^{1/(2\gamma)},
$$
where $N$ depends only on $\kappa,d,\delta$.

\end{theorem}

We know that under the condition of 
Theorem \ref{theorem 6.26,3} there exists
a Borel function $g(x,t,y)\geq0$, $(t,y)\in C_{2R^{2},R}$, such that for any Borel nonnegative~$f$
$$
E_{0,x} \int_{0}^{\tau_{2R^{2},R}  }f(t, x_{t})\,dt=\int_{C_{2R^{2},R}}f(t,y)g(x,t,y)\,dydt.
$$
This $g(x,t,y)$ is the Green's function
\index{$S$@Miscelenea!Green's function}%
of $X$ in $C_{2R^{2},R}$. By using
Theorems \ref{theorem 6.26,3} with
$f(t,y)=g^{-1}(x,t,y)I_{C_{R}(R^{2},0)}(t,y)$
we arrive at the following.

\begin{corollary}
                     \label{corollary 3.6,1}
Under the condition of 
Theorem \ref{theorem 6.26,3}
$$
\int_{C_{R}(R^{2},0)}g^{-1/(2\gamma)}(x,t,y)\,dydt
$$
$$
\leq NR^{d+2-1/\gamma}
\Big(E_{0,x} \int_{0}^{\tau_{2R^{2},R}  }g^{-1}(x,t, x_{t})I_{C_{R}(R^{2},0)}(t,x_{t})\,dt\Big)^{1/(2\gamma)}
$$
$$
=NR^{d+2-1/\gamma}|C_{R} |^{1/(2\gamma)},
$$
where $N$ is from Theorem \ref{theorem 6.26,3}.
\end{corollary}

 The following reformulation
of Corollary \ref{corollary 10.21,1} is used
to investigate the boundary behavior
of probabilistic solutions of
parabolic PDEs.

\begin{theorem}
               \label{theorem 6.26,4}
 
Let $R\leq  \rho_{b}$, $\xi\in(0,1)$,
 and assume that a closed set $\Gamma\subset
B_{R}$ is such that, for any $r\in(0,R)$, $|B_{r}\cap\Gamma|\geq
\xi |B_{r}|$. Then there exist    constants $\alpha\in(0,1)$
and $N$, depending only on $ d,\delta$,
   and $\xi$, such that, for any $ x\in B_{R/2}$,
\begin{equation}
                         \label{10.21,1}
 P_{0,x}(\tau'_{R} <\tau_{\Gamma} )\leq N( |x| /R)^{\alpha},
\end{equation}
where $\tau_{\Gamma} $ is the first time $  x_{t}$
hits $\Gamma$.
\end{theorem}

The next result has the same spirit
as Theorem 4.11 of \cite{Kr_21_2}
and shows the way Theorem \ref{theorem 6.26,4} can be applied
 investigating the boundary
behavior of solutions
of parabolic equations with drift in $L_{(q,p)}$.

 \begin{theorem}  
                     \label{theorem 10.20.3}
Let  $(\sfd_{0},p,q)$ be properly tight, $p<\infty$, $r\leq 1$, 
$T\in(0,\infty]$,
 and
let $D$ be a  domain in $\bR^{d}$ with $0\in\partial D$. 
Assume that for some constants $\rho,\rho_{1},\xi>0$, $\rho_{1}\in[0,\rho/2)$, and 
any $r\in [\rho_{1},\rho)$ we have 
$$
|B_{r}\cap D^{c}|\geq \xi
|B_{r}|.
$$
 Then there exist  $\beta >0$
and $N$, depending only on $ d,\delta $,
    $\xi$, $p$, with 
 $N$ also depending on $\rho$, $\rho_{b}$ and
either $T$, if $T<\infty$
and $D$ is unbounded, or the diameter of $D$, if $D$
is bounded,
such that, for any nonnegative $f\in L_{(q,p)}(Q)$, and $x\in D$, such that $|x|\geq \rho_{1}$,
\begin{equation}
                                                 \label{10.20.4}
u(x):=E_{0,x}\int_{0}^{\tau }f(t, x_{t})\,dt
\leq N |x|^{\beta}\sup_{C\in\bC_{1}}\|I_{Q}f\|_{L_{(q,p)}(C)},
\end{equation}
where $\tau $ is the first exit time of $(t, x_{t})$
from   $Q:=[0,T)\times D$.
\end{theorem}

Proof. We may assume that $\rho\leq  \rho_{b} $.
In light of Theorems \ref{theorem 9.7.1}
and  \ref{theorem 3.27.20}
we also may concentrate on $x\in B_{\rho/2}$. For $2\rho_{1}\leq 2|x|\leq r\leq \rho$  and $\tau^{r} $ being the first
exit time of $ x_{t}$ from $B_{r}\cap D$, we have
thanks to Theorems \ref{theorem 9.5.1} and  \ref{theorem 3.27.20} that  
\begin{equation}
                             \label{12.29.1}
u(x)=E_{0,x}\int_{0}^{T\wedge\tau^{r} }f(t, x_{t})\,dt
+E_{0,x}I_{\tau^{r} <\tau }E_{\tau^{r},x_{\tau^{r}}} \int_{0}^{\tau }f(\sft_{s},x_{s})\,ds 
\end{equation}
$$
\leq N r^{2-(d/p+2/q)}\sup_{C\in\bC_{1}}\|I_{Q}f\|_{L_{(q,p)}(C)}
+N\sup_{C\in\bC_{1}}\|I_{Q}f\|_{L_{(q,p)}(C)}P_{0,x}(\tau^{r} <\tau ). 
$$
Observe that $\{\tau^{r} <\tau \}\subset\{\tau^{r} <\tau_{\Gamma_{r}} \}$,
where $\Gamma_{r}=\bar B_{r}\cap D^{c}$, and by Theorem
\ref{theorem 6.26,4} we have
$P_{0.x}(\tau^{r} <\tau )\leq N( |x| /r)^{\alpha}$. Thus,
$$
u(x)\leq N\sup_{C\in\bC_{1}}\|I_{Q}f\|_{L_{(q,p)}(C)}\big( r^{2-(d/p+2/q)}
+ ( |x| /r)^{\alpha}\big),
$$
whenever $2\rho_{1}\leq 2|x|\leq r\leq \rho$.
For $\rho/2\geq |x|\geq \rho_{1}$, by choosing $r=\sqrt{2|x|\rho}$,
we get 
$$
u(x)\leq N\|f\|_{L_{q,p}(Q)}\big( |x|^{1-(1/2)(d/p+2/q)}
+ |x|^{\alpha/2}\big).
$$
Here
$$
1-\frac{d}{2p}-\frac{1}{q}\geq \frac{d_{0}}{p}
-\frac{d}{2p}>0,
$$
so the we get
the result with  $\beta=
\min(\alpha/2,(2d_{0}-d)/(2p))$. 
\qed

\begin{remark}
                        \label{remark 10.20.3}
Theorem \ref{theorem 10.20.3} is applicable
in case of time independent coefficients and $f$,
 $T=\infty$ when $b\in L_{d}(\bR^{d})$
(cf.~Lemma \ref{lemma 12.11.1}). Observe that
if $b$ and $f$ are bounded and a part of $\partial D$ near the origin
is flat, then one can take $\beta=1$ in \eqref{10.20.4}.
However, even in the case of flat boundary and bounded $f$,
if $b\in L_{d}(\bR^{d})$, then in the general case 
certainly $\beta<1$ (see Example 4.1  in \cite{Sa_10}) and most likely $\beta\to0$
as $\delta\to0$.

\end{remark}

One more result we need in the future
is the following.

\begin{theorem}
              \label{theorem 12.31.1}
Let $f(t,x)$ be a bounded Borel 
function vanishing for $|t|\geq T$
for some $T$. Introduce
$$
u(t,x)=E_{t,x}\exp\int_{0}^{\infty}
f(\sft_{s},x_{s})\,ds .
$$
Then there exists a constant $N$, which depends
only on $ d,\delta$,  $\rho_{b}$, such that    for any $R\leq \rho_{b} $
and $z_{1}$, $z_{2}\in  C_{R/2}$
\begin{equation}
                                           \label{12,14.60}
\big|u(z_{1})-u(z_{2})\big|\le NR^{-\gamma}
\rho^{\gamma}(z_{1},z_{2})\exp(T\sup|f|)
\big(1
+R^{2 }\sup|f|\big),
\end{equation}
  where
$\gamma=\gamma(d,\delta)\in(0,1)$.
\end{theorem}

Proof. By using the Markov property
we get
$$
u(t,x)=E_{t,x}\int_{0}^{\infty}
f(\sft_{r},x_{r})
\Big(\exp\int_{r}^{\infty}
f(\sft_{s},x_{s})\,ds\Big)\,dr
$$
$$
=E_{t,x}\int_{0}^{\infty}
f(\sft_{r},x_{r})u(\sft_{r},x_{r})\,dr.
$$
It follows that for $(t_{0},x_{0})\in C_{ R}$ by the strong Markov property
$$
u(t_{0},x_{0})=E_{t_{0},x_{0}}
\int_{0}^{\tau_{ R}}(fu)(\sft_{t},x_{t})\,dt+
E_{t_{0},x_{0}}u(\sft_{\tau_{ R}},x_{\tau_{ R}}),
$$
where $\tau_{ R}$ is the first exit time of 
$(\sft_{t},x_{t})$ from $C_{ R}$.
After that it only remains to use
Theorem \ref{theorem 12,14.2} along with the observation that
$$
u\leq \exp(T\sup|f|),\quad
\dashnorm fu
\| _{L_{q,p}(C_{ 2R})}\leq \sup|f|\exp(T\sup|f|).
$$
The theorem is proved. \qed 

We finish the section with couple
of estimates of the {\em resolvent operator\/},
\index{$S$@Miscelenea!resolvent operator}%
which is defined by
$$
 R_{\lambda}f(t,x):=E_{t,x}\int_{0}^{\infty}
e^{-\lambda s}f(t+s,x_{s})\,ds.
$$

\begin{lemma}
                \label{lemma 1.27.1}

If $f\geq0$ is independent of $t$, $p\in [d_{0},\infty]$, then for   $\lambda\geq \varkappa_{0}^{2}\rho_{b}^{-2}$
$$
\|\sup_{t}R_{\lambda}f(t,\cdot)\|_{L_{p}(\bR^{d})}
\leq N(d,\delta)\lambda^{-1}\|f\|_{L_{p}(\bR^{d})}.
$$

\end{lemma}

Proof. If $p=\infty$, the result is obvious. In the remaining case by Theorem \ref{theorem 8.30.1}
with $\tau\equiv 0$
we have
$$
R_{\lambda}f (0,0)\leq 
N \lambda^{-1+d/(2p)}
\| f  \bar\Psi_{\lambda}^{d_{0}/p}\|_
{L_{p   }(\bR^{d } )},
$$
where $\bar \Psi _{\lambda}( x)=\exp(- 
\sqrt{\lambda}  |x| \sfp_{0}/16)$, which by shifting the origin yields
$$
\sup_{t}
R_{\lambda}f (t,x)\leq 
N \lambda^{-1+d/(2p)}
\| f ( x+\cdot)\bar\Psi_{\lambda}^{d_{0}/p}\|_
{L_{p   }(\bR^{d } )}  .
$$

By observing that
$$
\int_{\bR^{d}}\| f ( x+\cdot)\bar\Psi_{\lambda}^{d_{0}/p}\|^{p}_
{L_{p   }(\bR^{d } )}\,dx=
\int_{\bR^{d}}f^{p}\,dx
\int_{\bR^{d}}\bar\Psi_{\lambda}^{d_{0}}\,dx,
$$
we easily finish the proof. The lemma is proved. \qed
\begin{theorem}
                                 \label{theorem 9.24.1}
Let  $(\sfd_{0},q,p)$ be properly tight.
Then
for any $\lambda\geq \varkappa_{0}^{2}\rho_{b}^{-2}$,   and Borel nonnegative
$f(t,x)$ 
we have
\begin{equation}
                                         \label{10.11.1}
\| R_{\lambda}f \|_{L_{(q,p)}(\bR^{d+1}_{0})}\leq N(d,\delta)
\lambda^{-1} \| f \|_{L_{(q,p)}(\bR^{d+1}_{0})}.
\end{equation}

\end{theorem}

Proof. The result is obvious if $\nu:=\nu(\sfd_{0},q,p)=1$.
Therefore, we assume that $\nu<1$,
so that at least one of $p,q$ is finite.  
If $p=\infty$, we have
$$
\| R_{\lambda}f \|^{q}_{L_{(q,p)}(\bR^{d+1}_{0})}=\int_{0}^{\infty}\sup_{x}
R^{q}_{\lambda}f(t,x)\,dt
$$
$$
\leq\int_{0}^{\infty}\Big(
\int_{0}^{\infty}e^{-\lambda s}\sup_{x}
f(t+s,x)\,ds\Big)^{q}\,dt
\leq\lambda^{-q}\int_{0}^{\infty}\sup_{x}
f^{q}(t ,x)\,dt,
$$
where the last inequality follows from the Minkowski inequality. If $q=\infty$,  the result follows from Lemma \ref{lemma 1.27.1}.
Therefore, we may concentrate on $p,q<\infty$.

By Theorem \ref{theorem 8.30.1}
we have
\begin{equation}
                      \label{10.23.1}
R_{\lambda}f (t,x)\leq 
N \lambda^{-\chi}
\| f (t+\cdot,x+\cdot)\Psi_{\lambda}^{1-\nu}\|_{L_{(q,p)}
(\bR^{d+1}_{0})}=:NJ(t,x),
\end{equation}
where $\chi=1-(1/2)(d/p+2/q)$ and $\Psi_{\lambda}(t,x)
=\exp(- \sqrt{\lambda} 
(|x|+ \sqrt t)\sfp_{0}/16)$.

If $ p\geq q$, we have
\begin{equation}
                                         \label{11.3.010}
I(t):=\int_{\bR^{d}}|\lambda^{\chi}  J(t,x)|^{p}\,dx
= 
\int_{\bR^{d}}\Big(\int_{0}^{\infty}
F^{q/p}(t,s,x)\,ds\Big)^{p/q}dx ,
\end{equation}
where
$$
F(t,s,x)=\int_{\bR^{d}}\Psi^{(1-\nu)p}_{\lambda}(s,y)f^{p} 
(t+s,x+y)\,dy.
$$
By Minkowski's inequality the   integral
on the right in \eqref{11.3.010} is dominated by
$$
\Big(\int_{0}^{\infty}\Big(\int_{\bR^{d}}F(t,s,x)\,dx\Big)^{q/p}
\,ds\Big)^{p/q},
$$
where
$$
\int_{\bR^{d}}F(t,s,x)\,dx=\int_{\bR^{d}}f^{p} (t+s,y)\,dy
\int_{\bR^{d}}\Psi^{(1-\nu)p}_{\lambda}(s,y)\,dy
$$
$$
\leq N(d,\sfp_{0}) [\mu(p)\sqrt\lambda]^{-d }e^{-\mu (p)\sqrt{\lambda s}}
\int_{\bR^{d}}f^{p} (t+s,y)\,dy,
$$
with $\mu(p)=(1-\nu)p\sfp_{0}/16$. 
It follows that
$$
I(t) 
\leq  N [\mu(p)\sqrt\lambda]^{-d } \Big(\int_{0}^{\infty}
e^{-\mu(q) \sqrt{\lambda s}}\Big(
\int_{\bR^{d}}f^{p}(t+s,y)\,dy\Big)^{q/p}ds\Big)^{p/q},
$$
$$
\|\lambda^{\chi}  J\| _{L_{(q,p)}(\bR^{d+1}_{0})}
\leq N    \mu^{-d /p}(p)  \mu^{-2/q}(q)
\lambda^{- d/(2p)-1/q}\|f \| _{L_{(q,p)}(\bR^{d+1}_{0})},
$$
which along with \eqref{10.23.1} and the fact that 
$\mu^{-d /p}(p)$ and  $\mu^{-2/q}(q)$ are bounded
from above ($1-\nu\geq 1/p,1/q$) yield \eqref{10.11.1}.  

If $q\geq p$,
$$ J(x):=
\int_{0}^{\infty}|\lambda^{\chi} J(t,x)|^{q}\,dt
= \int_{0}^{\infty}\Big(\int_{\bR^{d}}
 F^{p/q}\,(t,x,y) dy\Big)^{q/p}dt
$$
where
$$
F(t,x,y)=\int_{0}^{\infty}
\Psi^{(1-\nu)q}_{\lambda}(s,y)f^{q} 
(t+s,x+y)\,ds.
$$
By Minkowski's inequality
$$
\Big(\int_{0}^{\infty}|\lambda^{\chi}
  J(t,x)|^{q}\,dt\Big)^{p/q}
\leq  \int_{\bR^{d}}\Big(\int_{0}^{\infty}
F(t,x,y)\,dt\Big)^{p/q}dy,
$$
where
$$
\int_{0}^{\infty}
F(t,x,y)\,dt\leq
\int_{0}^{\infty}
 f^{q} 
( s,x+y)\,ds\int_{0}^{\infty}
\Psi^{(1-\nu)q}_{\lambda}(s,y) \,ds
$$
$$
= N\mu^{-2}(q)\lambda^{-1}e^{-\mu(q)\sqrt\lambda|y|}\int_{0}^{\infty}
 f^{q} 
( s,x+y)\,ds.
$$
Hence,
$$
J^{p/q}(x)
\leq N\mu^{-2p/q}(q)\lambda^{-p/q}\int_{\bR^{d}}e^{-\mu(p)\sqrt\lambda |y|}\Big(
\int_{0}^{\infty}
 f^{q} 
( s,x+y)\,ds\Big)^{p/q}\,dy,
$$
$$
\|\lambda^{\eta}  J\|_{L_{(q,p)}(\bR^{d+1}_{0})}^{p}
\leq N\mu^{-2p/q}(q)\mu^{-d}(p)\lambda^{-p/q-d/2}\|f \|_{L_{(q,p)}(\bR^{d+1}_{0})}^{p}
$$
and we again come to \eqref{10.11.1}. The theorem is proved. \qed

\mysection[Existence of regular
diffusion processes]{Existence of regular
diffusion processes}

Here we come back to processes
from Section \ref{section 5.14.1}
but instead of requiring
$b\in L_{(q,p)}$ with $d/p+1/q\leq 1$
we impose a different assumption
which results in additional
nice properties of the processes.

We suppose that on $\bR^{d+1}$ we are given Borel   $\bS_{\delta}$-valued $a,a(n),n=1,2,...$   and $\bR^{d}$-valued $b,
b(n),n=1,2,...$. Set $\sigma=\sqrt a$.

\begin{theorem}
                  \label{theorem 7.4,1}
Suppose that we are given $p_{i},q_{i}$,
$i=1,...,k$, such that   $(\sfd_{0},q_{i},p_{i})$ are properly tight and $\bR^{d}$-valued
$b_{i},b_{i}(n)\in L_{(q_{i},p_{i})}$,
 such that each $b_{i}(n)$ is bounded and
$$
b=\sum_{i=1}^{k}b_{i},\quad
b(n)=\sum_{i=1}^{k}b_{i}(n),\quad
b_{i}(n)\to b_{i}\quad\text{in}\quad
L_{(q_{i},p_{i})}(\cO)
$$
in any ball    $\cO\subset \bR^{d+1}$, and
$a(n)\to a$ (a.e.) as $n\to\infty$.
Also assume that, for each $n$, there exists a regular diffusion process $X^{n}=\{(\sft_{s},x_{s}),\cN_{s},P^{n}_{t,x}\}$
corresponding to $a(n),b(n)$
with $\rho_{b}\in(0,\infty)$
(in Definition \ref{definition 6.24,1}) independent of $n$.

Then there exists a regular diffusion process $X =\{(\sft_{s},x_{s}),\cN_{s},P _{t,x}\}$
corresponding to $a ,b $ and the same $\rho_{b}$.

\end{theorem}

Proof. In light of Theorem
\ref{theorem 6.24,1}, for 
any $R<\infty$, the family of distributions $\{P^{n}_{t,x},|t|,|x|\leq R,
n=1,2,...\}$
is tight. It follows that
for each $(t,x)$ there is a sequence
$P^{ n(k)}_{t,x}$ which
weakly converges to a probability
 distribution $P_{t,x}$ on
$C([0,\infty),\bR^{d+1})$. Cantor's
diagonal procedure shows that one can find such sequence which suits all $t,x$
with rational coordinated. We have thus defined $P_{t,x}$ for all $(t,x)$
with rational coordinated. 
Obviously, the statement of Theorem
\ref{theorem 6.24,1} is valid for $P_{t,x}$.
For any other
$(t,x)$ take a sequence with rational coordinates  $(t_{k},x_{k})\to (t,x)$
such that
$P_{t_{k},x_{k}}$ converges weakly and call
the limit $P_{ t,x }$.  

It turns out that
the limit is independent of the approximating sequence. Indeed,
the probability distribution $P_{t,x}$
on $C([0,\infty),\bR^{d+1})$ is uniquely
characterized by the set of numbers
\begin{equation}
                            \label{2.12.2}
\Phi(P_{t,x},f):=E_{t,x}\exp\int_{0}^{\infty}
f(\sft_{t},x_{t})\,dt
\end{equation}
when $f$ runs through the set of
continuous functions with compact support on $ \bR^{d+1}$.
By Theorem \ref{theorem 12.31.1}
the family $\{\Phi(P^{n}_{t,x},f),n=1,2,...
 \}$ is equicontinuous in $t,x$,
which implies that $\Phi(P_{t,x},f)$
and $P_{x}$ are indeed independent of the approximating sequence. It also
implies that for each $f$ 
and $R<\infty$
$$
\Phi(P^{ n(k)}_{t,x},f)\to 
\Phi(P_{t,x},f)
$$
uniformly for $|t|,|x|\leq R$. 

Next,  for each $r\geq0$ by the   Markov property of $X^{n}$
$$
E _{t,x} \exp\int_{0}^{r}
f(\sft_{s},x_{s})\,ds\,\Phi(P _{ \sft_{r},x_{r} },f )
$$   
$$
=\lim_{k\to\infty}E_{t,x}^{ n(k)}
\exp\int_{0}^{r}
f(\sft_{s},x_{s})\,ds\,\Phi(P _{ \sft_{r},x_{r} },f )
$$
$$
=\lim_{k\to\infty}E_{t,x}^{ n(k)}
\exp\int_{0}^{r}
f(\sft_{s},x_{s})\,ds\,\Phi(P^{ n(k)}_{\sft_{r},x_{r}},f )
$$ 
$$
=\lim_{k\to\infty}E_{t,x}^{ n(k)}
\exp\int_{0}^{r}
f(\sft_{s},x_{s})\,ds\,E_{\sft_{r},x_{r}}^{ n(k)}\exp\int_{0}^{\infty}
f(\sft_{s},x_{s})\,ds
$$
$$
=\lim_{k\to\infty}E_{t,x}^{ n(k)}
\exp\int_{0}^{\infty}
f(\sft_{s},x_{s})\,ds =E_{t,x} 
\exp\int_{0}^{\infty}
f(\sft_{s},x_{s})\,ds.
$$

The arbitrariness of $f$ implies that
$$
E_{t,x}\Big(\exp\int_{r}^{\infty}
f(\sft_{s},x_{s})\,ds\mid \cN_{r}\Big)=
\Phi(P_{\sft_{r},x_{r}},f )
$$
and that $X=((\sft_{s},x_{s}),\cN_{s},P_{t,x})$ is a Markov process. 

Since for bounded continuous $f$ on
$\bR^{d+1}$ we have $E^{n(k)}_{t,x}f(\sft_{s},x_{s})\to E_{t,x} f(\sft_{s},x_{s})$ and the convergence
is locally uniform, $E_{t,x} f(\sft_{s},x_{s})$ is continuous in $(t,x)$,
$X$ is   Feller and strong Markov.

By Theorem \ref{theorem 8.30.1}
 if $(\sfd_{0},p,q)$ is  properly tight,
$T\in (0,\infty)$, and
  Borel $f\geq0$ is given on $\bR^{d+1}$, then
\begin{equation}
                       \label{7.4,3}
E^{n}_{0,0}\int_{0}^{T} 
f(s,x_{s}) \,ds
\leq N\|\Psi f\|_{L_{(q,p)} }.
\end{equation}
where $\Psi _{\lambda}( x)=e^{- 
 |x|} $, and $N$ depend only  on
$d,\delta,p,q,T$.

A simple consequence of this  estimate 
and the weak convergence of $P_{x}^{ n(k)}$ is that for any
Borel $f$ such that $\Psi^{-1}f\in L_{(q,p)}$ we have
\begin{equation}
                       \label{2.12.3}
E_{0,0} \int_{0}^{T}  
f(\sft_{t},x_{t})\,dt=\lim_{k\to\infty}
E_{0}^{ n(k)} \int_{0}^{T}   
f(\sft_{t},x_{t})\,dt.
\end{equation}

Next we prove that the requirements   
in Definition \ref{definition 6.24,1} are satisfied for $X$,
so that $X$ is a regular diffusion process
corresponding to $\cL$. To start fix $(t,x)$ and observe that, in light of
\eqref{2.12.3}, \eqref{7.4,3}, and the assumption that $b_{i}\in L_{(q_{i},p_{i})}$,
$P_{(t,x)}$-(a.s.) we have
$$
\int_{0}^{T}|b(\sft_{s},x_{s})|\,ds<\infty,
\quad \forall T<\infty.  
$$
 Then take 
a twice continuously differentiable
function $u$ with compact support
and observe that in light of \eqref{2.12.3}
$$
I:= E_{0,0}u(\sft_{t},x_{t})-E_{0,0}\int_{0}^{t}\cL u(\sft_{s},x_{s})\,ds
$$
$$
= \lim_{k\to\infty}\Big[
E^{ n(k)}_{0,0}u(\sft_{t},x_{t})-E^{ n(k)}_{0,0}\int_{0}^{t}\cL u(\sft_{s},x_{s})\,ds\Big].
$$
Here for    $\cO$ such that $u=0$ in $\cO^{c} $ we have
$$
E^{ n(k)}_{0,0}\int_{0}^{t}|\big(\cL-\cL(n(k))\big) u(\sft_{s},x_{s})|\,ds
$$
$$
\leq NE^{ n(k)}_{0}\int_{0}^{t}
I_{\cO}\big(|a-a(n(k))|+
|b-b(n(k))|\big)(\sft_{s},x_{s})\,ds
$$
$$
\leq N\big(\|a-a(n(k))\|_{L_{\sfbd_{0}+1}(\cO)}+\sum_{i=0}^{k}\|b_{i}-b_{i}(n(k))\|_{L_{(q_{i},p_{i})}(\cO)}\big)\to 0
$$
as $k\to\infty$. Furthermore, for any $k$
$$
E^{ n(k)}_{0,0}u(\sft_{t},x_{t})-E^{ n(k)}_{0,0}\int_{0}^{t}\cL(n(k)) u(\sft_{s},x_{s})\,ds=u(0,0).
$$
Therefore, $I=u(0,0)$. Similar relation
holds at any other point $(t,x)\in\bR^{d+1}$. The strong Markov property
now implies that 
$$
u(\sft_{t},x_{t})-\int_{0}^{t}\cL u(\sft_{s},x_{s})\,ds
$$
is a $P_{t,x}$-martingale for any $(t,x)$
and smooth $u$ with compact support.
This fact and, for instance, an easy
adaptation of Lemma 3.4.1  of \cite{Kr_25}, proved there for time homogeneous
situation, prove that the requirement
(i) of Definition \ref{definition 6.24,1}
is fulfilled.

To finish proving the theorem it remains to check \eqref{7.14,1}.
For that we only need the following.
\begin{lemma}
                   \label{lemma 1.6.1}
Let $\rho\in(0,\infty)$, $C\in\bC_{\rho}$. Let $(\sfd_{0},q,p)$ be properly tight, $f,f^{n(k)}\in L_{(q,p)}(C)$, $f^{n(k)}\to f$ in $L_{(q,p)}(C)$.
Then for any $(t,x)\in\bR^{d+1}$
\begin{equation}
                  \label{2.12.6}
E_{t,x}\int_{0}^{\tau_{C}} f(\sft_{s},x_{s})\,ds
=\lim_{n\to \infty}
E^{ n(k)}_{t,x}\int_{0}^{\tau_{C}} f^{n(k)}(\sft_{s},x_{s})\,ds,
\end{equation}
where $\tau_{C}$ is the first exit time
of $(\sft_{s},x_{s})$ from $C$.
In particular, for any $\rho>0$,
$C\in\bC_{\rho}$, and $(t,x)\in\bR^{d+1}$
$$
E_{t,x}\int_{0}^{\tau_{C}}|b(\sft_{s},x_{s})|\,ds
=\lim_{n\to \infty}
E^{ n(k)}_{t,x}\int_{0}^{\tau_{C}}|b^{ n(k)}(\sft_{s},x_{s})|\,ds
\leq \sfb_{0} \rho.
$$

\end{lemma}

Proof.  First take a nonnegative and continuous  $f$ on $\bR^{d+1}$.
Then, as is easy to see, $\tau_{C}(x_{\cdot})$ and the function
$$
\Phi[f](\sft_{\cdot},x_{\cdot}):=\int_{0}^{\tau_{C}(\sft_{\cdot},x_{\cdot})} f(\sft_{s}, x_{s})\,ds
$$
on $C([0,\infty),\bR^{d+1})\cap\{(\sft_{\cdot},x_{\cdot}):(\sft_{0},x_{0})\in C \}$ are lower semicontinuous, and hence \eqref{2.12.6} with $f^{n(k)}=f$
holds with $\leq \nliminf$ in place of $=\lim$. In that case Theorem \ref{theorem 9.7.1}, 
valid for $X^{n}$, implies that
\begin{equation}
                                                     \label{1.7.1}
E_{t,x}\Phi[f]\leq
N\|f\|_{L_{(q,p)}(C)},  
\end{equation}
where $N$ depends only on $d,\delta,\rho_{b},\rho $. By standard means estimate \eqref{1.7.1} extends over to any Borel $f\geq0$.
In particular, this implies that,
to prove the lemma,
  it suffices to concentrate on 
$f^{n(k)}=f$ with
continuous $f\geq0$.

Furthermore, if $f\leq M$, for a constant $M$, then
$$
E_{t,x}\Phi[f]=ME_{t,x}\Phi[1]-
E_{t,x}\Phi[M-f]\geq ME_{t,x}\Phi[1]
-\nliminf_{k\to\infty}E^{ n(k)}_{t,x}\Phi[M-f]
$$
$$
\geq \nliminf_{k\to\infty}E^{ n(k)}_{t,x}\Phi[f]+ME_{t,x}\Phi[1]-M
\nliminf_{k\to\infty}E^{ n(k)}_{t,x}\Phi[1].
$$
This and what was said above about $\leq \nliminf$, implies that it only remains
to prove that
\begin{equation}
                                                     \label{1.7.2}
E_{t,x}\Phi[1]\geq
\nliminf_{k\to\infty}E^{ n(k)}_{t,x}\Phi[1].
\end{equation}

First, define $\rho(t,x)$ as the distance
between $(t,x)$ and $C^{c}$ and observe that
for any $t,\gamma>0$
$$
\{\tau_{C}>s\}=\{\inf_{r\leq s}\rho(
\sft_{r},x_{r})
>0\}\supset
\{\inf_{r\leq s}\rho(\sft_{r},x_{r})
\geq \gamma\}
\supset
\{\tau_{C^{\gamma}}>s\},
$$
where $C^{\gamma}=\{x:\rho(x)> \gamma\}$
and $\tau_{C^{\gamma}}$ is the first exit
time of $(\sft_{s},x_{s})$ from $C^{\gamma}$. Since
$\inf_{r\leq s}\rho(\sft_{r},x_{r})$ is a continuous function on $C([0,\infty),\bR^{d+1})$ we conclude
$$
P_{t,x}(\tau_{C}>s)
\geq \nlimsup_{k\to\infty}P_{t,x}^{ n(k)}(\inf_{r\leq s}\rho(\sft_{r},x_{r})
\geq \gamma)\geq 
\nlimsup_{k\to\infty}P_{t,x}^{ n(k)}
(\tau_{C^{\gamma}}>s).
$$
 
It follows that
$$
E_{t,x}\Phi[1]=\int_{0}^{\rho^{2}}
 P_{t,x}(\tau_{C}>s)\,ds
$$
$$
\geq \nlimsup_{k\to \infty}\int_{0}^{\rho^{2}}
 P^{ n(k)}_{t,x}(\tau_{C}>s)\,ds=
\nlimsup_{k\to \infty}E_{t,x}^{ n(k)}\tau_{C^{\gamma}}.
$$
Now clearly, to finish the proof it suffices to show that for any $\varepsilon,\gamma$, and $x$ we have
$|u^{n}_{\gamma}(t,x)-u^{n}(t,x)|\leq q(\gamma)$,
where $q(\gamma)\to 0$ as $\gamma
\downarrow 0$ and
$$
u^{n}_{\gamma}(t,x)
=E^{n}_{t,x}\tau_{C^{\gamma}},\quad
u^{n} (t,x)
=E^{n}_{t,x}\tau_{C }.
$$

By the strong Markov property
$$
u^{n} (t,x)-u^{n}_{\gamma}(t,x)=
E^{n}_{t,x} u^{n} (\sft_{\tau_{C^{\gamma}}},x
_{\tau_{C^{\gamma}}})
\leq \sup_{\rho(s,y)\leq \gamma} u^{n}(s,y)
\leq \sup_{\rho(s,y)\leq \gamma} E^{n}_{s,y}\tau_{C}.
$$
The last quantity goes to zero as
$\gamma\downarrow0$ according to
Theorem \ref{theorem 10.20.3}.
The lemma is proved. \qed

If $g\in C^{2}(\bar C)$,  by It\^o's formula,
$$
E_{t,x}e^{-\lambda \tau_{C}}g(\sft_{\tau_{C}},x_{\tau_{C}})=g(t,x)+E_{t,x}\int
_{0}^{\tau_{C}} 
 \cL g(\sft_{s},x_{s})\,ds.
$$
It follows that
\begin{equation}
                       \label{1.26.3}
E_{t,x}e^{-\lambda \tau_{C}}g(\sft_{\tau_{C}},x_{\tau_{C}})
=\lim_{k\to \infty}
E^{n(k)}_{t,x}e^{-\lambda \tau_{C}}g(\sft_{\tau_{C}},x_{\tau_{C}})
\end{equation}
for smooth $g$ and, by approximation,
for all continuous $g $ on $\bar C$. \qed
 
By using the arbitrariness
of   $g$ in \eqref{1.26.3}
we obtain the following.
\begin{corollary}
              \label{corollary 1.27.1}
For any $t,x$, the $P^{n(k)}_{t,x}$-distributions of $\tau_{C},x_{\tau_{C}}$ weakly converge to their $P_{t,x}$-distributions.
\end{corollary}
 
 The split
in the condition on $p_{i},q_{i}$ below
is needed because we are going to approximate $b$ by its mollifications. 

\begin{theorem}
               \label{theorem 10.21.1} 

Suppose that Assumption (ii)   
of Theorem \ref{theorem 9.27.10} 
is satisfied.
Also suppose that for each $i$,  either
$p_{i}+q_{i}<\infty$, or $p_{i}<\infty$,
$q_{i}=\infty$ and $b_{i}$ is independent of $t$.
Then there exists a regular strong Markov
process corresponding to $a,b$.
\end{theorem} 

Proof.   Approximate $a,b$ by smooth   
$a(n),b_{i}(n)$, $n=1,2,...$, 
by using mollifying kernel $n^{d+1}\zeta
(nt,nx)$, where nonnegative
$\zeta\in C^{\infty}_{0}(\bR^{d+1})$ has unit integral. Observe that $b_{i}(n)$
satisfy Assumption (ii) of Theorem \ref{theorem 9.27.10}
with the same $\rho_{b} $.
Therefore, by that theorem, for each $n$ there exists a regular diffusion process $X^{n}=((\sft_{s},x_{s}),\cN_{s},
P^{n})$
corresponding to $a(n),b(n)$.
Since the convergencies $b_{i}(n)\to b_{i}$
locally in $L_{(p_{i},q_{i})}$ are well known, to finish the proof, it only remains to refer to
Theorem~\ref{theorem 7.4,1}.  \qed

\begin{remark}
                 \label{remark 1.26.1}
It may look like  Theorem \ref{theorem 10.21.1}
is a generalization of Theorem 
\ref{theorem 9.6.4} (i) about the  solvability of
\eqref{11.29.20} with $b\in L_{(q,p)}$
and $d/p+1/q\leq 1$. However,
in the typical case of $k=0$,
along with $b \in L_{(q_{0},p_{0}),\loc}$,
$d_{0}/p_{0}+1/q_{0}\leq1$,
we require \eqref{9.27.4} to hold
and, if we ask ourselves what
$p,q$ should be in order for the
inclusion $b\in L_{(q,p)}$ to imply 
\eqref{9.27.4}, the answer is 
$d/p+2/q\leq 1$, somewhat disappointing.
At the same time in the next example 
we show that, in turn, Theorem 
\ref{theorem 9.6.4} does not cover
all applications of Theorem \ref{theorem 10.21.1}.
\end{remark}

In assumption  
\eqref{9.27.4}   the size
of $\hat b$ could not be too large.

\begin{example}
                           \label{example 12.21.01}
Let
$$
b(t,x)=b(x)=-\frac{d}{|x|}\frac{x}{|x|}I_{x\ne0},\quad 
a^{ij}=\delta^{ij},\quad \sigma=\sqrt2 (\delta^{ij}).
$$
Then as is easy to see, for any $p\in(1,d)$
and any $q$ the quantity $\rho\dashnorm b\|_{L_{q,p}(C)}$, $\rho>0,C\in\bC_{\rho}$, is bounded.
However, the equation $dx_{t}=\sigma\,dw_{t}
+b(x_{t})\,dt$ with initial condition $x_{0}=0$
does not have any solution.

Indeed, if it does, then by It\^o's formula
\begin{equation}
                              \label{12.21.5}
|x_{t}|^{2}=2d\int_{0}^{t}I_{x_{s}=0}\,ds
+2\sqrt2\int_{0}^{t}x_{t}\,dw_{t}.
\end{equation}
Here the first integral is the time spent at the
origin by $x_{s}$ up to time $t$. This integral is zero, because by using It\^o's formula for
$|x^{1}_{t}|$, one sees that the local time
of $x^{1}_{t}$ at zero exists and is finite,
implying that the real time spend at zero is zero.

Then \eqref{12.21.5} says that the local martingale
starting at zero which stands on the right is
nonnegative. But then it is identically zero,
implying the same for $x_{t}$. However,
$x_{t}\equiv0$, obviously, does not satisfy our equation.

At the same time according to Theorem
\ref{theorem 10.21.1},
the equation $dx_{t}=\sigma\,dw_{t}
+\varepsilon b(x_{t})\,dt$ with initial condition $x_{0}=0$ {\em does\/} have solutions if $\varepsilon$ is sufficiently small.
Observe that $b\not \in L_{(q,p),\loc}$
for any $p,q\in(1,\infty)$ satisfying
$d/p+1/q\leq1$, so this example is not covered by Theorem 
\ref{theorem 9.6.4}.
\end{example}

It turns out that in Theorem \ref{theorem 10.21.1} in 
the definition of $\hat b_{\rho_{b}}$   one
cannot replace $r $ with $r^{ 1+\alpha}$,
no matter how small $\alpha>0$ is.

\begin{example}
                         \label{example 12.21.4}

As in Section \ref{section 4.23.1}
take   numbers $\alpha$ and $\beta$ satisfying
$$
0<\alpha\leq \beta <1,\quad \alpha+\beta=1
$$
 and set
$$
b(t,x)=-\frac{1}{t^{\alpha}|x |^{\beta}}\frac{x }{|x |}
I_{0<|x|\leq 1,t\leq 1}.
$$
Using that $d_{0}<d$, it is not hard to find $p,q$
such that  $d_{0}/p+1/q<1$ and the quantity
$\rho^{1+\alpha}\dashnorm b\|_{L_{(q,p)}(C)}$,
$\rho>0, C\in\bC_{\rho}$, is bounded.
However, as we know from Section \ref{section 4.23.1}, the equation $dx_{t}=dw_{t}+\varepsilon
b(t,x_{t})\,dt$ with zero initial condition
does not have solutions no matter how small $\varepsilon>0$ is.
\end{example}

\mychapter[Applications to elliptic and parabolic equations]
{Applications to elliptic and parabolic  equations. Case \binLpq}

                    \label{chapter 3.24,1}

In the whole chapter we suppose that we are given
on $\bR^{d+1}$ a
Borel $\bS_{\delta}$-valued $a$ and a Borel
$\bR^{d}$-valued $b$. We set 
$$
\cL=\partial_{t}+(1/2)a^{ij}D_{ij}+b^{i}D_{i} 
$$

\mysection[Aleksandrov, Harnack, H\"older]{Aleksandrov's estimates,
Harnack inequality, H\"older continuity
of caloric functions}

The following has the same flavor as Nazarov's
Theorem 4.1 of  \cite{Na_87}
or Theorem 4.3 of  \cite{Kr_20_2}. We get a wider, than before,
range of $p,q$ on  account
of restricting $b$.
Here is a qualitative form
of the maximum principle (notice $L_{(q,p)},W^{1,2}_{(q,p)}$
and not $ L_{q,p}, W^{1,2}_{q,p}$). Its elliptic counterpart is found in
\cite{DK_22} with a proof completely different
from what is below.

The following result, in case $\cL u=0$,
$u$ is smooth,
and $b$ is bounded, is usually referred to as a
Krylov-Safonov type result.

\begin{theorem}
               \label{theorem 1.20.1}  
(a) Suppose that   we are given $p_{b},q_{b},p,q$ such that
$$
p_{b},q_{b}\in(1,\infty],\quad
\frac{\sfd_{0}}{p_{b}}+\frac{1}{q_{b}}\leq 1,
\quad p ,q \in(1,\infty),
\quad
\frac{\sfd_{0}}{p }+\frac{1}{q }\leq 1,   
$$
Take a $\rho_{b}\in(0,\infty)$ and set
$$
\hat{b}_{q_{b},p_{b},\rho_{b}}=\sup_{r\leq\rho_{b}}r
\sup_{C\in \bC_{r}} 
\dashnorm b \|_{L_{(q_{b},p_{b})}(C)}.
$$
(b) Suppose that
\index{$S$@Miscelenea!$\bar b_{R}$@$\hat{b}_{q_{b},p_{b},\rho_{b}}$}%
$$
N_{0}\hat{b}_{q_{b},p_{b},\rho_{b}}< \sfb_{0} ,
$$
where $N_{0}$ is taken from \eqref{9.27.4}.

Then for any
  $R\leq \varkappa_{0}^{-1}\rho_{b}$
and $ 
u\in W^{1,2}_{(q,p)}(C_{R})$

(i) there
exist constants $N$ 
and   
 $
\beta\in(0,1),
 $ 
depending only on $d,\delta $,
  such that, for any  
$z_{1},z_{2}\in C_{ R/2 }$, we have   
$$
\big|u(z_{1})-u(z_{2})\big|\le 
NR^{-\beta}\rho^{\beta}(z_{1},z_{2})\big(
\sup\big(|u|,\bar C_{  R}\big)
+NR^{2}\dashnorm \cL u\|_{L_{(q,p)}(C_{R})}\big);
$$

(ii) if $u\geq 0$, there exists a constant $N$, which depends 
only on $ \delta,d $, such that   
$$
u(R^{2}/2,0)\le Nu(0,x) +NR^{2}\dashnorm \cL u\|_{L_{(q,p)}(C_{R})}
$$
whenever $|x|\le R/2$. In particular, if 
$\cL u=0$ in $ C_{R} $, then $u(R^{2}/2,0)\le Nu(0,x)$ whenever $|x|\le R/2$.
\end{theorem}

Proof. One can replace $\cL$ in  (i), (ii)
with $\cL_{n}:=I_{|b|\geq n}(\partial_{t}+\Delta)+I_{|b|< n}\cL$
and then pass to the limit
by the dominated convergence and monotone convergence theorems.
Then we see that we may assume that $b$ is bounded.
After that, having in mind mollifications
and the dominated convergence theorem, we may assume
that $a$ and $b$ are smooth.

In this situation we take
the process $X$ corresponding to  $\sigma=\sqrt{2a}$
from Theorem \ref{theorem 10.21.1}, which  
by Theorem \ref{theorem 9.27.10}  satisfies
Assumption   \ref{assumption 8.19.2}
and makes all results of Chapters 
\ref{chapter 10.20.1}  and \ref{chapter 3} available.  By using 
 It\^o's formula from Theorem \ref{theorem 10.15.1}
we arrive at
\begin{equation}
                                          \label{2.8.1}
u(t,x)=E_{t,x} u(t+\tau,x_{\tau})+E_{t,x} \int_{0}^{\tau}f(t+s,x_{s})\,ds,
\end{equation}
where $f=-\cL u $  and $\tau$ is the first exit time
of $(t+s,x_{s})$ from $Q$.
Then assertion (i) follows from Theorem \ref{theorem 12,14.2}
and assertion (ii) follows from Theorem \ref{theorem 1.18.13}. 
The theorem is proved. \qed

We will see later that assertions (i) and (ii)  hold  true
also under different assumptions on $u,b$.

\begin{remark}
                         \label{remark 12.21.1}
Theorem \ref{theorem 1.20.1} is applicable
to elliptic equations. It suffices to suppose that
$a,b$ are independent of $t$, take $q_{b}=\infty$
and consider $u$ independent of $t$. Then, for $p_{b}
=\sfd_{0}$, the condition on $b$ becomes
\begin{equation}
                          \label{12.15.1}
\hat b_{\sfd_{0},\rho_{b}}:=\sup_{\rho\leq
\rho_{b}}\rho\sup_{B\in\bB_{\rho}}\dashnorm b\|_{L_{\sfbd_{0}}(B) }<\sfb_{0}/N_{0}, 
\end{equation}
for some $\rho_{b}\in(0,\infty)$. Since $\sfd_{0}<d$,
by H\"older's inequality this condition
is satisfied if $\|b\|_{L_{d}(B)}< \sfb_{0}/N_{0}$,
which is true if $|B|$ is small, for instance, if $b\in L_{d}(\bR^{d})$.
In this particular case the assertions of
Theorem \ref{theorem 1.20.1}
(for smooth $u$) are obtained in
Safonov \cite{Sa_10} {\em with no restriction on\/}
$R$ and with constants depending only on $d,\delta$,
and $\|b\|_{L_{d}(\bR^{d})}$ and no other characteristics of $b$
 are involved.  
On the other hand, as we know
from Example 5.3.5 of \cite{Kr_25}, our condition may be satisfied 
with $p<d$ but $b\not\in L_{p+,\loc}(\bR^{d})$.
\end{remark}

One more result  is the parabolic
Aleksandrov estimate with mixed norms
and ``supercritical'' $b$.    

\begin{theorem}
                  \label{theorem 10.14.1}
Suppose that the   assumptions of
Theorem \ref{theorem 1.20.1} are satisfied. 
Let $R\in(0,\infty)$, domain $Q\subset C_{R}$, 
and let 
$u\in W^{1,2}_{(q,p),\loc}(Q)\cap C(\bar Q)$.

  Take a function
$c\geq 0$ on $Q$. Then on $ Q$
\begin{equation}
                                  \label{10.14.10}
u \leq  N
\|I_{Q,u>0}( \cL u-cu)_{-}\|_{L_{(q,p)} }
+\sup_{\partial'Q}u_{+},
\end{equation}
where $N=N(\delta,d,\rho_{b},     R)$.
In addition, if $R\leq \varkappa_{0}^{-1}\rho_{b}$, we have
$N=N(d,\delta)R^{2-(d/p+2/q)}$. In particular (the maximum principle),
if $ \cL u-cu \geq0$ in $Q$ and $u\leq0$
on $\partial'Q$, then $u\leq 0$ in $Q$.
\end{theorem}

Proof. Obviously the right-hand side of
\eqref{10.14.10} decreases if we replace $c$ with zero.
Hence, we may assume that $c=0$. 
Also, we need to prove \eqref{10.14.10}
only in $Q\cap\{u>0\}$ on the parabolic boundary of
which either $u=0$ or $u\leq \sup_{\partial'Q}u_{+}$.
Therefore, we may assume that $u>0$ in $Q$.

Then for $\varepsilon>0$ define $Q^{\varepsilon}$
as the collection of $(t,x)\in Q$ such that 
the closed ball in $\bR^{d+1}$ centered at $(t,x)$
with radius $\varepsilon$ lies in $Q$.
Obviously $Q^{\varepsilon}$ is   open.
As we have seen in the proof of Theorem \ref{theorem 6.10.1},  $\dist(\partial' Q,\partial'Q^{\varepsilon})=\varepsilon$.
It follows, owing to the continuity of $u$
and the monotone convergence theorem, that
it suffices to prove \eqref{10.14.10}
with $Q^{\varepsilon}$ in place of $Q$.
As a consequence of that we may assume that
$u\in W^{1,2}_{(q,p)}(Q)$.

This  gives us the opportunity  
to replace $\cL$ in  \eqref{10.14.10}
with $\cL_{n}:=I_{|b|\geq n}(\partial_{t}+\Delta)+I_{|b|< n}\cL$
and then pass to the limit
by the dominated convergence and monotone convergence theorems.
Hence, we may assume that $b$ is bounded.
After that, having in mind mollifications
and the dominated convergence theorem, we may assume
that $a$ and $b$ are smooth.

In this situation
we can write \eqref{2.8.1} for $u$
and then, to prove \eqref{10.14.10},
 it only remains
to use Theorem \ref{theorem 3.27.20}     
  and the fact that   
$(t+\tau
,x_{\tau})\in 
\partial'Q$. The second assertion follows from 
Theorem \ref{theorem 9.5.1}.
The theorem is proved.  \qed

Needless to say that applying
Theorem \ref{theorem 10.14.1} to $-u$
in place of $u$ one also gets the lower estimate of $u$. 

\begin{remark}
The result of Theorem \ref{theorem 10.14.1}  
in case   $b\equiv0$ and $p=q=\sfd_{0}+1$
can be found   in Remark 1
of \cite{Es_93}, where it is given without proof.
Complete proof for bounded $b(t,x)$
can be extracted from \cite{CKKS_98} dealing with viscosity  
solutions of fully nonlinear parabolic equations.
 We have mixed norms and our $b$ is ``supercritical''.
\end{remark}

\begin{example}[Cf. Example 1.3.1 \cite{LSU_67}]
                       \label{example 1.5.1}
The condition on $\hat b_{q_{b},p_{b},\rho_{b}}$
in assumption (b) of Theorem  \ref{theorem 1.20.1}, basically,
reduces to the requirement of it to be  sufficiently small.
It turns out that this smallness assumption
is essential.   

For instance,  
the function
$$
u (t,x)=\int_{0}^{1-t}\exp(-|x|^{2}/(4s))\,ds,
\quad (t,x)\in C_{1},
$$
has bounded derivatives $\partial_{t}u$ and $Du$,
its second-order derivative $D^{2}u$ is unbounded
only when $x$ is close to $0$ and is of order
$\ln |x|$ as $x\to 0$, and for $x\ne0$ it
 satisfies
\begin{equation}
                             \label{1.9.10}
 \partial_ t u (t,x)
+\Delta u (t,x)-\frac{d}{|x|^{2}}x^{i}D_{i}u (t,x)=0,
\end{equation}
where $|b (t,x)| =d/|x|$ whose constant 
  $\hat b_{q_{b},p_{b},\rho_{b}}$ is finite and even independent of $\rho_{b}$
as long as $p_{b}<d$. In this situation according to   Theorem \ref{theorem 1.20.1}, if $\hat b_{p_{b},q_{b},\rho_{b}}$
were small enough then we would have
$$
1=u(0,0)\le \sup_{t\leq 1}
\int_{0}^{1-t}\exp(-1/(4s))\,ds
= \int_{0}^{1}\exp(-1/(4s))\,ds,
$$
which is wrong.
 
\end{example}

 The next theorem
is taken from \cite{Kr_21_1}. It generalizes
the corresponding result of \cite{FS_84} on  account
of having nonzero $b$.
With bounded $b$ the result is found in Cabr\'e
\cite{Ca_95} 
and in Fok  \cite{Fo_98} for
  $b\in L_{d+\varepsilon}$ adapting his result
to the linear case.

\begin{theorem}
                                         \label{theorem 4.1.1}
Let $ G $ be a bounded domain in $\bR^{d}$,
$u\in W^{2}_{\sfd_{0},\loc}( G )\cap C(\bar  G )$.
Also assume that $a$ and $b$ are independent of $t$
 and for some $\rho_{b}\in(0,\infty)$
condition \eqref{12.15.1} is satisfied.
Take a function $c\geq0$. Then on $ G $
\begin{equation}
                                               \label{4.1.5}
u \leq N 
\|I_{ G ,u>0}(  \cL u-cu)_{-}\|_{L_{\sfbd_{0}}(\bR^{d}) }
+\sup_{\partial G }u_{+},
\end{equation}
where $N$ depends only on $d,\delta, \rho_{b}$, and
the diameter of
$  G $.

\end{theorem}

The proof of this theorem is
obtained by mimicking that of Theorem
\ref{theorem 10.14.1} and is  again based on It\^o's formula and Theorem \ref{theorem 3.27.20} with
$q=\infty$.
 
\begin{remark}
                          \label{remark 12.24.1}
The constant $\hat b_{\sfd_{0},\rho_{b}}$ in condition \eqref{12.15.1}
should not be too large. For instance,
$u(x)=1-|x|^{2}$ satisfies $\Delta u+b^{i}D_{i}u=0$
in $G=B_{1}$, where $b(x)=-cx/|x|^{2}$, $c=d$, and for
this $b$ the let-hand side of \eqref{12.15.1} is finite. However,
\eqref{4.1.5} fails.     

\end{remark}

Results like the next one were used
in the theory of fully nonlinear elliptic and parabolic equations (cf.~\cite{Kr_18}).

\begin{theorem}
                    \label{theorem 10.14.2}
Suppose that the   assumptions of
Theorem \ref{theorem 1.20.1} are satisfied.
Then there exists
 a constant  $N $, depending only on
$ d,\delta  $,
such that, for any $\lambda\geq  \varkappa_{0}^{2}\rho_{b}^{-2} $, $R\in(0,\infty]$ and $u\in W^{1,2}_{(q,p),\loc}(C_{R})\cap C(\bar C_{R})$ ($C_{\infty}=\bR^{d+1}_{0}$, $C(\bR^{d+1}_{0})$
is the set of bounded continuous functions on $\bR^{d+1}_{0}$),  
 we have
$$
\lambda\|u_{+}\|_{L_{(q,p)}(C_{R/2})}
\leq 
 N\| (\lambda u-\cL u )_{+}\|_{L_{(q,p)}(C_{R}) } 
$$
\begin{equation}
                                           \label{10.14.40}  
+ N\lambda R^{d/p+2/q}e^{-  R\sqrt{\lambda}\sfp_{0}/2 }
\sup_{\partial' C_{R}}u_{+} ,
\end{equation}
 where the last term should be dropped if $R=\infty$.   
\end{theorem}

Proof.   By having in mind 
the possibility to approximate
$C_{R}$ from inside by similar domains,
we see that we may assume that $R<\infty$ and 
$u\in W^{1,2}_{(q,p)}(C_{R})$.
Then as in the proof of Theorem \ref{theorem 10.14.1}
we reduce the general case to the one in which
$b$ is bounded. After that   we see that we may assume that $a,b,u$
are smooth. Then we take the process $X$
 as in the proof of Theorem \ref{theorem 1.20.1}  
 for $(t,x)\in C_{R/2}$, 
similarly to  \eqref{2.8.1}, write
$$
u(t,x)=E_{t,x}  e^{-\lambda\tau_{R}}u(t+\tau_{R} ,x_{\tau_{R}})
+E_{t,x}  \int_{0}^{\tau_{R}}e^{-\lambda t}f(t+s,x_{s})\,ds
$$
$$
=:I(t,x)+J(t,x),
$$
where $f=\lambda u-\cL u $,
$\tau_{R}$ is the first exit time of $(t+s,x_{s})$
from $C_{R}$ and $x_{s}$ is a solution of \eqref{11.29.20} with $\sigma=\sqrt{2a}$. 

Here, thanks to \eqref{8.20.1}  
$$
I(t,x)\leq Ne^{-  R\sqrt{\lambda}\sfp_{0}/2 }\sup_{\partial'C_{R}}u_{+},
$$
\begin{equation}
                                               \label{10.4.30}
 \|I_{C_{R/2}}I_{+}\|_{L_{(q,p)}  }\leq
NR^{d/p+2/q}e^{-  R
\sqrt{\lambda}\sfp_{0}/2 }\sup_{\partial'C_{R}}u_{+},
\end{equation}
where the $N  $'s depend  only on    
$ d,\delta $. We get the estimate of $J$ from Theorem
\ref{theorem 9.24.1} and by combining it with \eqref{10.4.30} arrive at \eqref{10.14.40}. 
The theorem is proved. \qed

The full strength of the results like Theorem \ref{theorem 10.14.2}
is seen in the theory of fully nonlinear equations.
But even for linear ones one gets a nontrivial information
as, for instance, in the following theorem which, in particular,
implies that, if $R<\infty$,
the operator $\cL: D\to
L_{(q,p)}(C_{R}) $ with the domain
$$
D:=W^{1,2}_{(q,p)}(C_{R})\cap\{u: \cL u\in
L_{(q,p)}(C_{R}),
u_{\big|\partial'C_{R} }=0\}
$$
is a closed operator in $L_{(q,p)}(C_{R})$.

\begin{theorem}
                                               \label{theorem 2.8.1}

Suppose that the   assumptions of
Theorem \ref{theorem 1.20.1} are satisfied
and
take $R\in(0,\infty)$. Suppose we are given
 $u_{0},u_{1},...\in D$ and $f\in L_{(q,p)}(C_{R})$ such that, for 
$f_{n}:=\cL u_{n} $, we have
$$
\sup_{n\geq1}\sup_{\partial'C_{R}}|u_{n}|<\infty,
\quad \|f_{n}-f\|_{L_{(q,p)}(C_{R})}+\|u_{n}-u_{0}\|_{L_{(q,p)}(C_{R})}
\to 0
$$
as $n\to\infty$. Then $\cL u_{0} =f$ in $C_{R}$.
\end{theorem}

Proof. Take a smooth $\psi$ on $C_{R}$ and apply
\eqref{10.14.40} to $u_{n}-u+\psi/\lambda$
in place of $u$. Then pass to the limit as $n\to \infty$
to find
$$
\|\psi_{+}\|_{L_{(q,p)}(C_{R/2})}
\leq N_{1}\lambda R^{d/p+2/q}
e^{-\kappa R\sqrt{\lambda} }
$$
$$
+N_{2}\|\psi-f+\cL u_{0} 
- L\psi  /\lambda
\|_{L_{(q,p)}(C_{R})},
$$
where $N_{2}$ is independent of $\lambda$ and $\psi$
and $N_{1}$ is independent of $\lambda$. By setting
$\lambda\to\infty$ we get
$$
\|\psi_{+}\|_{L_{(q,p)}(C_{R/2})}
\leq N_{2}\|\psi-f+\cL u_{0} \|_{L_{(q,p)}(C_{R})}.
$$
This is true if $\psi$ is smooth enough and by approximation
is true for any $\psi\in L_{(q,p)}(C_{R})$. For
$\psi=f-\cL u_{0} $ we get that
 $f-\cL u_{0} \leq 0$
in $C_{R/2}$. The reader understands that here as well as in  
\eqref{10.14.40} one can take any number $<R$ in place of $R/2$.
Hence, $f-\cL u_{0} \leq 0$ in $C_{R}$. Passing to
$-u_{n}$, $-f$ yields $f-\cL u_{0} \geq 0$
and proves the theorem. \qed

To investigate the boundary behavior of solutions
we need the following.

\begin{lemma}
                        \label{lemma 21.27.1}
Suppose that the  assumptions of
Theorem \ref{theorem 1.20.1} are satisfied
and $p_{b}<\infty$,
let    $T\in(0,\infty]$, 
 and
let $G$ be a  bounded domain in $\bR^{d}$ with $0\in\partial G$. Set $Q=(0,T)\times G$.
Assume that for some constants $\rho,\rho_{1},\xi>0$, $\rho_{1}\in[0,\rho/2)$, and 
any $r\in [\rho_{1},\rho)$ we have $|B_{r}\cap G^{c}|\geq \xi
|B_{r}|$. Let $u\in W^{1,2}_{(q,p)}(Q)\cap C(\bar Q)$, and let $\omega(r)$, $r\geq0$, be a concave function such that
$$
|u(t,x)-u(s,y)|\leq \gamma+\omega(\gamma+|x-y|+\sqrt{ |t-s| }),
$$ 
whenever  $(t,x), (s,y)\in   \partial' Q$, where the constant $\gamma\geq 0$. Then there exist  $\beta >0$
and $N$, depending only on $ d,\delta $, 
    $\xi$,  $p$, with 
 $N$ also depending on $\rho$, $\rho_{b}$ and
 the diameter of $G$, 
such that  for $x\in G$, satisfying $|x|\geq \rho_{1}$,
\begin{equation}
                           \label{12.27.4}
|u(t,x)-u(t,0)|
\leq \gamma+N|x|^{\beta}\sup_{C\in\bC_{1}}\|I_{Q}\cL u\|_{L_{(q,p)}(C)}+
\omega\big(\gamma+N_{0} |x|^{\beta_{0}}\big),
\end{equation}
where $\beta_{0} >0$
and $N_{0}$, depend  only on $ d,\delta $,
    $\xi$,  $p$,  $p_{b}$ with 
 $N_{0}$ also depending on $\rho$, $\rho_{b}$ and
 the diameter of $G$, 
\end{lemma}

Proof. As a few times before we may concentrate
on the case in which $a,b$ are smooth and $b$
is bounded. In that case the strong
Markov process $X$ corresponding to $a,b$
possesses all properties from Chapter \ref{chapter 3}. In particular, by It\^o's
formula for $(t,x)\in Q$
$$
u(t,x)=E_{t,x}\int_{0}^{\tau}f(t+s,x_{s})\,ds
+E_{t,x}u(t+\tau,x_{\tau}),
$$
where $\tau$ is the first  exit time of
$(\sft_{s},x_{s})$ from $Q$ and $f=-I_{Q}\cL u$. Here by Theorem
\ref{theorem 10.20.3} for $|x|\geq \rho_{1}$  
$$
\Big|E_{t,x}\int_{0}^{\tau}f(t+s,x_{s})\,ds\Big|
\leq N|x|^{\beta}\sup_{C\in\bC_{1}}\|f\|_{L_{(q,p)}(C)}.
$$
Also by using Jensen's inequality we get  
$$
\big|u(t,x)-E_{t,x}u(t+\tau,x_{\tau})\big|
\leq E_{t,x}|u(t+\tau,x_{\tau})-u(t,x)|
$$
$$
\leq \gamma+\omega\big(\gamma+E_{t,x}(
|x_{\tau}-x|+\tau^{1/2})\big),
$$
where by Theorem
\ref{theorem 10.20.3} 
$$
E_{t,x}\tau^{1/2}\leq \Big(E_{t,x}\int_{0}^{\tau}\,ds
\Big)^{1/2}\leq N|x|^{\beta/2}.
$$
Also
$$
E_{t,x}|x_{\tau}-x|\leq E_{t,x}\Big|\int_{0}^{\tau}
\sigma(t+s,x_{s})\,dw_{s}\Big|+
E_{t,x}\int_{0}^{\tau}|b(t+s,x_{s})|\,ds,
$$
where the first term is less than $N(d,\delta)
E_{t,x}\tau^{1/2}$. By Theorem
\ref{theorem 10.20.3} with  $p=p_{b},q=q_{b}$,
the second term is dominated by
$$
N_{1}|x|^{\beta_{1}}\sup_{C\in\bC_{1}}\|b\|_{L_{q_{b},p_{b}}(C)}\leq N_{2}
|x|^{\beta_{1}}\sup_{C\in\bC_{\rho_{b}}}\|b\|_{L_{q_{b},p_{b}}(C)}\leq N_{3}
|x|^{\beta_{1}},
$$
where $\beta_{1}>0$ and $N_{1}$  depend only on $ d,\delta $,
    $\xi$, $p_{b}$, with 
 $N_{1}$ also depending on $\rho$, $\rho_{b}$ and the diameter of $D$ and
the first inequality is valid owing
to Remark \ref{remark 2.29.1}, and the second one
follows from assumption (b) in Theorem
\ref{theorem 1.20.1}. The lemma 
is proved.\qed   

 \begin{theorem}
                     \label{theorem 12.27.3}
In Lemma \ref{lemma 21.27.1} let $\gamma=0$
and assume that $u\in W^{1,2}_{q,p,\loc}(Q)$ instead of $u\in W^{1,2}_{q,p}(Q)\cap C(\bar Q)$.
 Then the assertion of Lemma \ref{lemma 21.27.1} still holds (with $\gamma=0$).

\end{theorem}

Proof.  For $\varepsilon>0$ define $G_{\varepsilon}$
as the set of point in $G$ whose distance to
$\partial G$ is strictly greater than $\varepsilon$.
Set $Q^{\varepsilon}=(0,T)\times G_{\varepsilon}$.
Let $x^{\varepsilon}\in\partial G_{\varepsilon}$
be the closest point on $\partial G_{\varepsilon}$
to the origin. Obviously, $x^{\varepsilon}\to0$
as $\varepsilon\downarrow 0$. Owing to the
uniform continuity
of $u$ in $\bar Q$ there exist  $\gamma^{\varepsilon}
\to 0$ as $\varepsilon\downarrow 0$ such that
$$
|u(t,x)-u(s,y)|\leq \gamma^{\varepsilon}+\omega(\gamma^{\varepsilon}+|x-y|+\sqrt{ |t-s| }),
$$ 
whenever  $(t,x),( s,y)\in \partial' Q_{\varepsilon}$.
Also for $\rho\geq r>|x^{\varepsilon}|$ we have $|B_{r}(x^{\varepsilon})\cap G_{\varepsilon}^{c}|\geq |B_{r-|x^{\varepsilon}|}\cap G^{c}|\geq \xi |B_{r-|x^{\varepsilon}|}|=\xi(1-|x^{\varepsilon}|/r)^{d}|B_{r}|
\geq (\xi/2)|B_{r}|$ if $|x^{\varepsilon}|$
is sufficiently close to zero. It follows by Lemma \ref{lemma 21.27.1}
that for $x\in G_{\varepsilon}$, $|x-x^{\varepsilon}|\geq |x^{\varepsilon}|$, we have
$$
|u(t,x)-u(t,x^{\varepsilon})|
\leq \gamma^{\varepsilon}+N|x|^{\beta}\sup_{C\in\bC_{1}}\|I_{Q}\cL u\|_{L_{q,p}(C)}+
\omega\big(\gamma^{\varepsilon}+N_{0} |x|^{\beta_{0}}\big).
$$
By sending $\varepsilon\downarrow 0$, we get the result.
The theorem is proved.  \qed

\mysection[Fanghua Lin estimate]{Fanghua Lin estimate} 

Here
we suppose that for a constant $\varepsilon
\in(\delta,1)$
\begin{equation}
                          \label{11.21.1}
N_{0}(d,\delta)\hat{b}_{q_{b},p_{b},\rho_{b}}<
\varepsilon \sfb_{0},\quad a\in \bS_{\delta/\varepsilon},
\end{equation}
where $N_{0}=N_{0}(d,\delta)$ is the constant in \eqref{9.27.4}.

\begin{theorem}
               \label{theorem 1.26.1}  
Under the above assumption
suppose also that the first  assumption of
Theorem \ref{theorem 1.20.1} is satisfied, $R\leq \varkappa_{0}^{-1}\rho_{b}$,
$u\in W^{1,2}_{(q,p)}(C_{R})$. Then
there exists an $\alpha=\alpha(d,\delta)\in(0,1)$ such that 
\begin{equation}
                          \label{12.26.1}
\dashnorm D^{2}u,Du\|_{L_{\alpha}(C_{R})}\leq
N\dashnorm \cL u\|_{L_{(q,p)}(C_{R})}+NR^{-2}\sup_{\partial'C_{R}}
|u|,
\end{equation}
where $N=N(d,\delta,  
\rho_{b} )$.
 
\end{theorem}

Proof. As usual,
we may assume, first, that $b$ is bounded
and then that $a,b,u$  are smooth. Then set $r=p\vee q$,
$\bar u(t,x)=u(|t|,x)$  and  introduce $v$ as a   $W^{1,2}_{r}(C_{R}\cup(-C_{R}))$-solution of $\cL v=I_{C_{R}}\cL u$ with boundary condition $v=\bar u$. By classical theory (see, for instance,
\cite{LSU_67}) such solution exists and is unique. By Theorem \ref{theorem 10.14.1}
\begin{equation}
                            \label{12.27.2}
|v(-R^{2},0)|\leq N(d,\delta)R^{2}
\dashnorm \cL u\|_{L_{(q,p)}(C_{R})}+
\sup_{\partial' C_{R}}|u|.
\end{equation}

Next, due to  assumptions \eqref{11.21.1} and the fact that, say,
$$
|Du|=D_{i}u\lim_{\tau\downarrow0}
\frac{D_{i}u}{\sqrt{|Du|^{2}+\tau}},
$$ 
 there exist
 $\kappa=\kappa(d,\delta,\varepsilon,
 \rho_{b} )>0$ and, for any $\gamma>0$, there
 exist  smooth $\bS_{\delta}$-valued $\check a$ 
 and smooth $\bR^{d}$-valued $\check b$ on $\bR^{d+1}$ such that 
$$
N_{0}\sup_{r\leq\rho_{b}}r
\sup_{C\in \bC_{r}}
\dashnorm \check b\|_{L_{(q_{b},p_{b})}(C) }<\sfb_{0},
$$
$$
\|(\partial_{t}+\check a^{ij}D_{ij}+\check
b^{i}D_{i}) v-I_{C_{R}}(\cL u-\kappa|D^{2}u|-
\kappa|Du|)\|_{L_{(q,p)}
(C_{R})}\leq \gamma.
$$
  Then we take the
Markov process $X=\{(\sft_{\cdot}, x_{\cdot}), \cN_{t},P_{t,x}\}$ 
 corresponding to $\check \cL=\partial_{t}+\check a^{ij}D_{ij}+\check
b^{i}D_{i}$
and by It\^o's formula conclude
\begin{equation}
                     \label{12.27.01}
v(-R^{2},0)\geq E_{-R^{2},0}\int_{0}^{\tau}
I_{C_{R}}(-\cL u+\kappa|D^{2}u|+\kappa|Du|)(\sft_{s},x_{s})\,ds
-N_{1}\gamma,
\end{equation}
where $\tau$ is the first exit time of
$(\sft_{s},x_{s})$ from $C_{R}\cup(-  C_{R})$ 
and $N_{1}$ is independent of $\gamma$.
By Theorem \ref{theorem 9.5.1}
$$
E_{-R^{2},0}\int_{0}^{\tau}
I_{C_{R}}(-\cL u)(\sft_{s},x_{s})\,ds
\geq -N(d,\delta)R^{2}
\dashnorm \cL u\|_{L_{(q,p)}(C_{R})}.
$$
This and \eqref{12.27.2}, \eqref{12.27.01} yield  
$$
\kappa E_{-R^{2},0}\int_{0}^{\tau}
I_{C_{R}}(|D^{2}u|+|Du|)(\sft_{s},x_{s})\,ds
$$
$$
\leq
N(d,\delta)R^{2}
\dashnorm \cL u\|_{L_{(q,p)}(C_{R})}+
\sup_{\partial' C_{R}}|u|+N_{1}\gamma,
$$
and to finish the proof it only remains to use
Theorem \ref{theorem 6.26,3} and let $\gamma
\downarrow0$. The theorem is proved. \qed

\mychapter[Weak uniqueness]{Weak uniqueness}

                      \label{chapter 3.7.1}

After the classical work by K. It\^o showing that there exists
a unique (strong) solution of \eqref{6.15.2} if $\sigma$ and $b$
are Lipschitz continuous in $x$ (may also depend on   $\omega$ and the nondegeneracy of $\sigma$ is not required),
  much effort was  applied to relax these   conditions.   The first author who achieved a considerable progress was A.V. Skorokhod
\cite{Sk_61} who proved the solvability assuming
only the continuity of $\sigma$ and $b$ 
in $x$ (which may depend
on $t$ and again without nondegeneracy).
Then by using the Skorokhod method and Aleksandrov
estimates the author proved in \cite{Kr_69}
and \cite{Kr_77} the solvability
for the case of {\em just measurable $\sigma$ and  
bounded $b$}  under the nondegeneracy assumption.
Stroock and Varadhan \cite{SV_79} among
many other things
not only proved
the solvability for the coefficients uniformly 
continuous in $x$, but also proved the uniqueness
of their distributions.

It is worth saying
that
restricting
the situation to the one when $\sigma$ and $b$
are independent of time allows one to
relax the above conditions significantly
further, see, for instance, \cite{KS_19}
and the 
references therein. 

The main results of this chapter
are close to  \cite{Ki_25}, 
which contains the most powerful results
in case $\sigma=(\delta^{ij})$ and $b$ is
in a Morrey class (in $(t,x)$).
Still this paper is not completely
covering the results in \cite{RZ_20} or our results
in case $b\in L_{2,\infty}$, 
 $\sigma=(\delta^{ij})$.   Our uniqueness theorem and uniqueness
theorems in \cite{Ki_25} and \cite{RZ_20} are   conditional.
We prove uniqueness only in the class of solutions
({\em which is proved to be nonempty\/}) admitting certain
estimates.   However,
in Theorem \ref{theorem 12.12.3} we 
mention a sufficient analytic
condition on $b$ when the unconditional weak uniqueness 
holds. 

In Remarks 
\ref{remark 10.31.1}  and 
\ref{remark 1.28.1}  we compare our results with some previous ones  and we
 refer the reader to \cite{BFGM_19}, 
\cite{Ki_24}, \cite{Ki_25}, \cite{RZ_20}   for  
 very good reviews  of the recent history 
 of the problem.    
  
   Recall that according to  
Example \ref{example 3.22.2}  assuming that $b\in L_{q,p}$ with $d/p +1/q \leq 1$
alone does not guarantee weak uniqueness even with unit diffusion
(the existence is known).
In Example \ref{example 3.22.1}  for any $\varepsilon>0$ we have
$b\in L_{q,p}$ with $d/p +1/q \geq 1+\varepsilon$ and there are no solution
of \eqref{6.15.2} with unit diffusion and $(t,x)=0$ at all.
In Example 3 of \cite{KZ_75} it is given an equation
$dx_{t}=\sigma(x_{t})\,dw_{t}$ in $d=2$ with $\sigma\sigma^{*}
=(\delta^{ij})$ such that it has unique and strong solutions
for any starting point apart from the origin. If the
starting point is the origin, only weak solutions exist.
All these examples show that we are dealing with quite
delicate problems, many of which are to date far
from being settled in the most satisfactory way.

Here is an example in which we prove existence 
(and conditional uniqueness) of weak solutions:
$$
\sigma=2(\delta^{ij})+I_{x \ne0}\zeta(x )\sin(\ln|\ln |x |)
$$
 (quite discontinuous),
where $\zeta$ is any smooth symmetric $d\times d$-matrix valued 
function vanishing for $|x|>1/2$ and
satisfying $|\zeta|\leq 1$, and $|b|=c/|x|$
with $c$ sufficiently small.
Another example of $b$ is when
$$
|b|= c|x|^{-\gamma}(|x|+\sqrt{|t|})^{\gamma-1},
\quad \gamma\in (d/(d+1),2d/(2d+1))
$$
with $c$ sufficiently small.
Both examples of $b$ are admissible in
\cite{Ki_25} and inadmissible in 
  \cite{RZ_20}. In both articles $\sigma$ is constant.
By the way, as we know
  the equation $dx_{t}=
  dw_{t}-|x_{t}|^{-1}b(x_{t})\,dt$, where $b(x)
  =(d/2)x/|x|$, with initial data $x_{0}=0$
  does not have solutions, so that in the above examples
  $c$ indeed should be sufficiently small.

In particular, we prove a generalization of the Stroock-Varadhan theorem in \cite{SV_79} obtained
for $\sigma$ that is uniformly continuous in $x$
uniformly in $t$
and bounded $b$.
  We need an additional assumption on $a$
and can relax conditions imposed on $b$ in Section
\ref{section 1.20.1} when $k=0$.
Since $a$ will have some regularity, the range of $p_{b},q_{b}$ can be substantially extended.
Indeed, observe that if $ \sfd_{0}/p_{b}
+1/q_{b}= 1$, then $ 1<d/p_{b}+2/q_{b}< 2$
since $d>\sfd_{0}>d/2 $.

\mysection[Morrey-Sobolev spaces]{Morrey-Sobolev spaces}   

An important distinction of the rest of the book from previous chapters
 is that
here,  
for $p,q\in(1,\infty)$ and domain $Q\subset\bR^{d+1}$, by $L_{q,p}(Q)$
we mean the space of Borel (real-, vector- or matrix-valued)
 functions on $Q$   with finite norm given 
either by
\begin{equation}
                                 \label{3.27.3}
\|f\|_{L_{q,p}(Q)}^{q}=\|fI_{Q}\|_{L_{q,p}}^{q}
=\int_{\bR}\Big(\int_{\bR^{d}}|fI_{Q}(t,x)|^{p}\,
dx\Big)^{q/p}\,dt 
\end{equation}
or
\begin{equation}
                                 \label{4.3.2}
\|f\|_{L_{q,p}(Q)}^{p}=\|fI_{Q}\|_{L_{q,p}}^{p}
=\int_{\bR^{d}}\Big(\int_{\bR}|fI_{Q}(t,x)|^{q}\,dt
 \Big)^{p/q}\,dx.
\end{equation}
 
{\em One of the ways to choose the norm
is fixed throughout the rest of the book
unless specifically stated otherwise.
We will be referring to some results that are
proved elsewhere for only one of the norms   \eqref{3.27.3} or \eqref{4.3.2}. In such
situations we mean that the result, we are referring to, actually, holds
for both norms and is proved by 
insignificant changes in the original proof.
This is, for instance, explicitly mentioned and 
underlined in \cite{Kr_26_1}.}

Take $p,q\in(1,\infty)$, domain $\cO\subset
\bR^{d}$, and
   $\beta> 0$ and introduce the
Morrey space $\dot E_{q,p,\beta}(\cO) $
as the set of $g\in  L_{p,q,\loc}(\cO)$ 
\index{$A$@Sets of functions!$\dot E_{q,p,\beta}(\cO)$}%
\index{$N$@Norms!$"|"|f"|"|_{\dot E_{q,p,\beta}(\cO)}$}%
such that  
\begin{equation}
                             \label{2.16.5}
\|g\|_{\dot E_{q,p,\beta} (\cO)}:=
\sup_{\rho >0,C\in\bC_{\rho}}\rho^{\beta}
\dashnorm g I_{\cO} \|_{ L_{q,p}(C)} <\infty .
\end{equation}  
Define
$$
\dot E^{1,2}_{q,p,\beta}(\cO) =\{u:u,Du,D^{2}u, 
\partial_{t}u\in \dot E_{q,p,\beta}(\cO) \},
$$
where $Du,D^{2}u,
\partial_{t}u$ are Sobolev
\index{$A$@Sets of functions!$\dot E^{1,2}_{q,p,\beta}$}%
 derivatives,
and 
provide $\dot E^{1,2}_{q,p,\beta} (\cO)$ with an obvious norm. If $\cO=\bR^{d+1}$ we drop
$\cO$ in $\dot E _{p,q,\beta} (\cO)$ and
$\dot E^{1,2}_{q,p,\beta} (\cO)$.
The subsets of these spaces
consisting of functions independent of
$t$ are denoted by $\dot E_{p,\beta}$
and $\dot E^{2}_{p,\beta}$, respectively.
\index{$A$@Sets of functions!$\dot E_{p,\beta}$}%
\index{$A$@Sets of functions!$\dot E^{2}_{p,\beta}$}%
It is a good idea for the
 reader to keep in mind that
if $\beta>d/p+2/q $ and $
f\in \dot E_{q,p,\beta} $, then
$f=0$ (a.e.) and if $\beta=d/p+2/q $,
then $E_{q,p,\beta}=L_{q,p}$.

Another useful observation is that if
$C\in \bC_{\rho}$, and $1\leq d/p+2/q$, then $I_{C}\in \dot E_{q,p,1}$ and $\|I_{C}\|
_{\dot E_{q,p,1}}=\rho$.

For functions $f(t,x)$ and $\varepsilon>0$ we define
\begin{equation}
                        \label{2.16.6}
f^{(\varepsilon)}=f*\zeta_{\varepsilon},
\end{equation}
 where
  $\zeta_{\varepsilon}=\varepsilon^{-d-2}\zeta
(t/\varepsilon^{2},x/\varepsilon)$, with a nonnegative
$\zeta\in C^{\infty}_{0}$ which has unit integral and support
in $C_{1}(-1,0)$. Observe that
owing to Minkowski's inequality  
$$
\|f^{(\varepsilon)}\|_{\dot E_{q,p,\beta}}
\leq \|f \|_{\dot E_{q,p,\beta}} 
$$
for any $f\in \dot E_{q,p,\beta}$.
Maximal function boundedness is
another notable property of Morrey spaces. 

For $\beta\geq0$ define
\index{$S$@Miscelenea!$\beta$-maximal function}%
 the parabolic 
$\beta$-maximal function of $f$
by
$$
\bM_{\beta} f(t,x)=\sup_{C\in\bC,C\ni(t,x)}
\dashint_{C}f(s,y)\,dyds,\quad \bM=\bM_{0}.
$$

\begin{theorem}[Theorem 6.7 \cite{Kr_26_1}]
              \label{theorem 2.15.1}
Let   $p,q\in(1,\infty)$, $\beta>0$. Then there is a constant
$N $ such that, for any $f\geq0$

 \begin{equation}
                   \label{2.15.1}
\|\bM f\|_{\dot E_{q,p,\beta} }\leq
N   \| f\|_{\dot E_{q,p,\beta} }.
\end{equation}
In particular,
$$
\|\sup_{\varepsilon>0}f^{(\varepsilon)}
\|_{\dot E_{q,p,\beta} }\leq
N   \| f\|_{\dot E_{q,p,\beta} }.
$$
\end{theorem}   
 
Here are   useful approximation results. Set
$$
\|f\|_{E_{q,p,\beta}}:=\sup_{\rho \leq 1,C\in\bC_{\rho}}\rho^{\beta}
\dashnorm f  \|_{ L_{q,p}(C)},\quad \|f\|_{E_{q,p,\beta}(C)}:=
\|fI_{C}\|_{E_{q,p,\beta}},
$$   
(note $\rho\leq1$) and define $E_{q,p,\beta}$ as the collection of $f$ with finite $E_{q,p,\beta}$-norm. The space
\index{$A$@Sets of functions!$E_{q,p,\beta}$}%
\index{$N$@Norms!$"|"|f"|"|_{E_{q,p,\beta}}$}%
\index{$A$@Sets of functions!$E^{1,2}_{q,p,\beta}$}%
$E^{1,2}_{q,p,\beta}$ is defined as the set
of $u\in E_{q,p,\beta}$ such that
$\partial_{t}u,D^{2}u,Du\in  E_{q,p,\beta}$. We provide this space with a natural norm.
\begin{lemma}[Lemma 6.4  \cite{Kr_26_1}]
                               \label{lemma 2.16.2}
Let $q,p\in(1,\infty)$, $0\leq\beta'<\beta$.
If   $\|f_{n}\|_{E_{q,p,\beta'}}$,
$n=0,1,...$, is a {\em bounded\/} sequence
and $f_{n}\to f_{0}$ in $L_{q,p}(C)$
for any
$C\in \bC$, then for any
$C\in \bC$
$$
\lim_{n\to\infty}\|f_{n}-f_{0}\|_{E_{q,p,\beta}(C)}=0. 
$$
In particular, if
  $f\in E_{q,p,\beta'}$, then for any $C\in \bC$
\begin{equation} 
                            \label{2.16.3}
\lim_{\varepsilon\downarrow0}
\|f^{(\varepsilon)}-f \|_{E_{q,p,\beta}(C) }=0.
\end{equation}
\end{lemma}
 
The proof of this fact is obtained by observing that for $r\leq 1$
$$
r^{\beta}\dashnorm f_{n}-f\|_{L_{q,p}(C_{r})}
\leq 2\varepsilon^{\beta-\beta'}\sup_{n}
\|f_{n}\|_{E_{q,p,\beta'}}+N(\varepsilon,C)
\|f_{n}-f\|_{L_{q,p}(C)}.
$$
  
\begin{lemma}[Lemma 6.5  \cite{Kr_26_1}]
                \label{lemma 2.16.3}
Let $p,q\in(1,\infty)$,  $g(t,x)\geq 0$ be a Borel function such that for any smooth bounded $f(t,x)$ we have
(bounded linear functional  of special type)
\begin{equation}
                       \label{2.16.4}
\int_{\bR^{d+1}}g|f|\,dxdt\leq
\|f\|_{E_{q,p,\beta}}. 
\end{equation}
Then, for any $f\in E_{q,p,\beta}$,
\eqref{2.16.4} holds and, moreover,
$$
\lim_{\varepsilon\downarrow 0}
\int_{\bR^{d+1}}
g|f-f^{(\varepsilon)}|\,dxdt=0.
$$
\end{lemma}

Here are parts of Lemmas 2.5 and 2.8 of \cite{Kr_23_1}.
\begin{lemma}
                   \label{lemma 11.15,1}
(i) Let   $0< \beta <2 $. Then any $u\in E^{1,2}_{q,p,\beta}$
is bounded and continuous  and for any $\varepsilon
\in(0,1]$  
\begin{equation}
                     \label{3.20.06}
|u|
\leq \varepsilon^{2-\beta} \|\partial_{t}u,D^{2}u\|_{E_{q,p,\beta}}+N(d, \beta)\varepsilon^{-\beta} \| u\|_{E _{p,q,\beta}}.
\end{equation} 

(ii) Let $1<\beta\leq d/p+2/q$, $\beta<2$. Then for any
$u\in E^{1,2}_{q,p,\beta}$, $\rho\leq 1$,
$(t_{i},x_{i})\in C_{\rho}$, $i=1,2$, we have
\begin{equation}
                           \label{4.6.1}
|u(t_{1},x_{1})-u(t_{2},x_{2})|\leq
N(d,p,q,\beta)\rho^{2-\beta}\|u\|_{E^{1,2}_{q,p,\beta}}.
\end{equation}
\end{lemma}

The following is a corollary of H\"older's
inequality.

\begin{lemma}[Lemma 6.6  \cite{Kr_26_1}]
                \label{lemma 2.16.5}
If $p,q\in(1,\infty)$, $\beta>1$, 
and
$(p_{0},q_{0})=(p,q)\beta=(r,s)(\beta-1) $, then for any
$f,g$
$$
\|fg\|_{\dot E_{q,p,\beta}}
\leq \|f\|_{\dot E_{q_{0},p_{0},1}}\|g\|_{\dot E_{s,r,\beta-1}}.
$$

\end{lemma}

For $k,s,r>0,\alpha\in \bR $, and appropriate $f(t,x)$'s
on $\bR^{d+1}$
\index{$S$@Miscelenea!$p_{\alpha,k}(s,r)$}%
\index{$C$@Operators!$P_{\alpha,k}f(t,x)$}%
 define
$$
p_{\alpha,k}(s,r)=\frac{1}{s^{(d+2-\alpha)/2}}e^{-r^{2}/(ks)}I_{s>0}, 
$$
$$
P_{\alpha,k}f(t,x)=\int_{\bR^{d+1} }p_{\alpha,k}(s,|y|)f(t+s,x+y)\,dyds.
$$
$$
=\int_{t}^{\infty}\int_{\bR^{d} }p_{\alpha,k}(s-t,|y-x|)
f(s,y)\,dsdy.
$$  

\begin{theorem}
                     \label{theorem 10.5,1}
(i) There is a constant $c(d)>0$ such that
$u=c(d)P_{2,4}(\partial_{t}u+\Delta u)$
if $u\in C^{\infty}_{0}$.

(ii) For $\alpha,\beta,k>0$ we have
$P_{\alpha,k}P_{\beta,k}=c(\alpha,\beta,k)P_{\alpha+\beta,k}$.

(iii) For any integer $n\geq1$, $\alpha\geq d+2+n$, and bounded $f$
with compact support we have $|D^{n}P_{\alpha,k}f|\leq N(d,\alpha,n)P_{\alpha-n,2\kappa}|f|$.

\end{theorem}

Proof. Assertion (i) follows from It\^o's
formula applies to $u(t,\sqrt 2 w_{t})$, where
$w_{t}$ is the $d$-dimensional Wiener process.
Assertion (ii) follows after direct computations. Assertion (iii) is also
proved by direct computations augmented
by the fact that $r^{m}e^{-r^2/\kappa}\leq
N(m,\kappa)e^{-r^2/(2\kappa)}$. \qed

The following fact is for the information only.
It will not be used in the future.

\begin{lemma}
                      \label{lemma 12.14,1}
For any $k>0$, integer $n\geq1$ and $f\in L_{2}(\bR^{d+1})$ we have
$$
\|D^{n}P_{n,k}f\|_{L_{2}}\leq N(d,k,n)\|f\|_{L_{2}}.
$$
\end{lemma}

Proof. We will give the proof for $n=1$
only. The proof in the general case is similar. Denote by $\tilde g(t,\xi)$ the Fourier
transform of $g(t,x)$ with respect to $x$.
Then  
$$
\widetilde{D_{j}P_{1,k}f}(t,\xi)
= \int_{0}^{\infty}\Big(\int_{\bR^{d}} p_{1,k}(s,|y|)e^{i(\xi,y)}\,dy\Big)
i\xi^{j}\tilde f(t+s,\xi)\,ds
$$
$$
=N(d)\int_{0}^{\infty}s^{-1/2}e^{-k|\xi|^{2}s/4}i\xi^{j}\tilde f(t+s,\xi)\,ds.
$$
By  the Young inequality
 $$
\int_{\bR}\big|\widetilde{D_{j}P_{1,k}f}(t,\xi)\big|^{2}\,dt\leq \Big(\int_{0}^{\infty}
s^{-1/2}e^{-k|\xi|^{2}s/4}|\xi|\,ds\Big)^{2}
\int_{\bR}|\tilde f(t ,\xi)|^{2}\,dt.
$$
After that it only remains to note that the first factor on the right is independent of $\xi$ and integrate with respect to $\xi\in\bR^{d}$. \qed

  The following
is nontrivial
only if $\beta\leq d/p+2/q $.
\begin{theorem}[Theorem 6.8 \cite{Kr_26_1}]
                       \label{theorem 2.11.2}
Let    
$q_{1},q_{2}\in(1,\infty]$, $k>0$,
$0<\alpha<\beta $. Then 
there is a constant $N$ such that for any $f\geq0$
we have
\begin{equation}
                          \label{10.7.4}
 \|P_{\alpha,k}f\|_{\dot E_{r_{1},r_{2},\beta-\alpha}}
\leq N \|f\|_{\dot E_{q_{1},q_{2},\beta}},
\end{equation}
where $r_{i}(\beta-\alpha)=q_{i}\beta$,
$i=1,2$.
\end{theorem}
 This theorem and Theorem \ref{theorem 10.5,1}
lead to the following
 
\begin{corollary}[Corollary 6.4 \cite{Kr_26_1}]
                      \label{corollary 10.8.1}
Under the assumptions of Theorem \ref{theorem 2.11.2}, if $\beta>1$,  for any $u\in C^{\infty}_{0}$  we have
$$
\|Du\|_{\dot E_{r_{1},r_{2},\beta-1}}
\leq N\|\partial_{t}u+\Delta u\|_{\dot E_{q_{1},q_{2},\beta }},
$$
where $r_{i}(\beta-1)=q_{i}\beta$,
$i=1,2$.  
\end{corollary}

\begin{corollary} 
                      \label{corollary 10.8.2}
Under the assumptions of Theorem \ref{theorem 2.11.2}, 
if $1<\beta\leq d/p+2/q$,  for any $u\in C^{\infty}_{0}$  we have
$$
\|Du\|_{  E_{r_{1},r_{2},\beta-1}}
\leq N\| u\|_{  E^{1,2}_{q_{1},q_{2},\beta }},
\quad\text{(no dots)}
$$
where $r_{i}(\beta-1)=q_{i}\beta$,
$i=1,2$.  
\end{corollary} 

This follows immediately from Corollary \ref{corollary 10.8.1} after taking there
$\zeta u$ in place of $u$, where $\zeta
\in C^{\infty}_{0}$ and $\zeta=1$ in a ball
of radius 1.

As an easy consequence of H\"older's inequality
 we have the following.
\begin{lemma}
                \label{lemma 6,11.1}
If  
\begin{equation}
                       \label{9,19.1}
\frac{d}{p}+\frac{2}{q}\geq \beta>1 
\end{equation}
and
$(q_{0},p_{0})=(q,p)\beta=(r,s)(\beta-1) $, then for any
$f,g$
$$
\|fg\|_{E_{q,p,\beta}}
\leq \|f\|_{E_{q_{0},p_{0},1}}\|g\|_{E_{r,s,\beta-1}}.
$$

\end{lemma}

Here is a mixed-norm analog of the parabolic Adams  theorem.
\begin{theorem}[Theorem 6.4 \cite{Kr_26_1}]
                        \label{theorem 5.25,1}
Let $ \alpha>0,q_{1},q_{2}\in(1,\infty),q>\max(q_{1},q_{2})$, $k>0$, $b(t,x)\geq0$. Then for any $f(t,x)\geq0$
\begin{equation}
                            \label{5.25,1}
\|bP_{\alpha,k}f\|_{L_{q_{1},q_{2}}}
\leq N\|b\|_{\dot E_{q,q,\alpha} }
\|f\|_{L_{q_{1},q_{2}}},
\end{equation}
where $N$ depends only on $d,q_{i},q,\alpha,k$.
In particular, for any $u\in C^{\infty}_{0}$
\begin{equation}
                        \label{11.22,3}
\|b Du \|_{L_{q_{1},q_{2}}}
\leq N\|b\|_{\dot E_{q,q,1} }K,
\quad \| b u\|_{L_{q_{1},q_{2}}}
\leq N\|b\|_{\dot E_{q,q,2} }K,
\end{equation}
where $K=\|D^{2}u,\partial_{t}u\|_{L_{q_{1},q_{2}}}$
and $N$ depends only $d,q_{i},q$.
\end{theorem}

\begin{remark}
                     \label{remark 11.2.1}
The first estimate in \eqref{11.22,3} follows from \eqref{5.25,1} with $\alpha=1$   and the fact that
for $f=\partial_{t}u+\Delta u$ we have
$$
Du(t,x)=c\int_{\bR^{d +1 }_{0}}
\frac{y}{s^{(d+2)/2}}e^{-|y|^{2}/(4s)}
f(t+s,x+y)\,dyds,
$$
where $c$ is a constant  and
$(|y|/s^{1/2})e^{-|y|^{2}/(4s)}\leq
Ne^{-|y|^{2}/(8s)}$. The second estimate
follows when $\alpha=2$ since
$$
u(t,x)=c\int_{\bR^{d +1 }_{0}}
\frac{1}{s^{d/2}}e^{-|y|^{2}/(4s)}
f(t+s,x+y)\,dyds.
$$

\end{remark}

\begin{corollary}
                 \label{corollary 10.5,1}
Estimate \eqref{5.25,1} says that the operator
$f\to bP_{\alpha,k}f$ is bounded in $L_{q_{1},q_{2}}$. Its conjugate (with time reversed)
is then also bounded as an operator in 
$L_{q'_{1},q'_{2}}$, where $q_{i}'=q_{i}
/(q_{i}-1)$, that is
$$
\|P_{\alpha,k}(bf)\|_{L_{q'_{1},q'_{2}}}
\leq N\|b\|_{\dot E_{q,q,\alpha} }
\|f\|_{L_{q'_{1},q'_{2}}}.
$$
In case $q_{1}=q_{2}=2$ and $\alpha=1$ we have that, if $q>2$, then
$$
\|P_{1,k}(bf)\|_{L_{2} }
\leq N\|b\|_{\dot E_{q,q,1} }
\|f\|_{L_{2} }. 
$$
\end{corollary}

A useful addition to the above properties
of multiplication by $b$ is the following.

\begin{lemma}
                          \label{lemma 11.29,1}
Let $q>2$, $f(x)\geq0$, $b\geq 0$, $T\in \bR$, $c=(4\pi)^{-d/2}$
,$$
u(t,x)=\frac{c}{(T-t)^{d/2}}\int_{\bR^{d}}e^{-|x-y|^{2}/(4T-4t )}f(y)\,dy\, I_{t<T}.
$$
Then
\begin{equation}
                                \label{11.29,2}
\int_{(-\infty,T)\times\bR^{d}}
b^{2}u^{2}\,dxdt\leq N(d,q)\|b\|^{2}_{\dot E_{q,q,1} }\int_{\bR^{d}}f^{2}\,dx.
\end{equation}
\end{lemma}

Proof. We may assume that $T=0$ and $f$ is smooth
and bounded. In that case
set $v(t,x)=u(-t,x)$, $w(t,x)=v(t,x)\zeta(t)$, $t\geq0$,
where $\zeta$ is infinitely differentiable $\zeta=1$ near zero, $\zeta(t)=0$ for $t\geq 1$,
$\zeta\geq0,\zeta'\leq0$.
Observe that $\partial_{t}u+\Delta u=0$ for $t<0$,
$\partial_{t}v=\Delta v$ for $t>0$,
$\partial_{t}w=\Delta w+v\zeta'$ for $t>0$,
which after being multiplied by $w$ and integrating by parts 
yields
\begin{equation}
                             \label{11.29,3}
(1/2)\int_{\bR^{d+1}_{0}}|Dw|^{2}\,dxdt=
\int_{\bR^{d}}f^{2}\,dx
+\int_{0}^{\infty}\Big(\int_{\bR^{d}_{0}}v^{2}\,dx
\Big)
\zeta\zeta'\,dt\leq \int_{\bR^{d}}f^{2}\,dx.
\end{equation}

By It\^o's formula for $t\geq0$ we have
$$
w(t,x)=-cP_{2,4}(\partial_{t}w+\Delta w)(t,x)=-
cP_{2,4}(2\Delta w+v\zeta')(t,x).
$$
For us the most important is that this holds with
$t=0$ when $w(0,x)=f(x)$. Then note that because of
  the semigroup property of the heat semigroup for $t<0$ we have
$$
\frac{c}{(-t)^{d/2}}\int_{\bR^{d}}
e^{-|x-y|^{2}/(-4t)}cP_{2,4}(hI_{0,\infty)})(0,y)\,dy
=cP_{2,4}(hI_{0,\infty)})(t,x).
$$
It follows that for $t\leq0$ we have
$$
u(t,x)=-cP_{2,4}\big((2\Delta w+v\zeta')I_{0,\infty)}\big)(t,x).
$$
Next, we use that $|P_{2,4}\Delta w|=|(D_{i}
P_{2,4}D_{i} w|\leq NP_{1,8}|Dw|$ and
$P_{2,4}=NP_{1,4}P_{1,4}$ combined with the fact 
that, obviously, $\zeta'\in \dot E_{q,q,1}$.
Then by applying Theorem \ref{theorem 5.25,1} and Corollary
\ref{corollary 10.5,1} we arrive at
$$
\int_{(-\infty,T)\times\bR^{d}}
b^{2}\big(P_{2,4}(I_{0,\infty)}\Delta w)\big)^{2}\,dxdt
\leq N\|b P_{1,8}| I_{0,\infty)}Dw|\,\|^{2}_{L_{2}R^{d+1}}
$$
$$
\leq N\|b\|^{2}_{\dot E_{q,q,1} }
\|Dw\|^{2}_{L_{2}(\bR_{0}^{d+1})},
$$
$$
\int_{(-\infty,T)\times\bR^{d}}
b^{2}\big(P_{2,4}|I_{0,\infty)}v\zeta'|\big)^{2}\,dxdt
\leq N\|b\|^{2}_{\dot E_{q,q,1} }
\|P_{1,8}|I_{0,\infty)}v\zeta'|\,\|^{2}_{L_{2}(\bR ^{d+1})}
$$
$$
\leq N\|b\|^{2}_{\dot E_{q,q,1} }\|vI_{(0,1)}\|
^{2}_{L_{2}(\bR ^{d+1})}.
$$
After that it only remains to use \eqref{11.29,3}
and   that, for any $t>0$,
$$
\int_{\bR^{d}}v^{2}(t,x)\,dx\leq \int_{\bR^{d}}f^{2}\,dx.
$$
The lemma is proved. \qed

If $u$ and $b$ are independent of $t$
the first estimate in \eqref{11.22,3}
follows from the well-known Chiarenza-Frasca
result that $\|bu\|_{L_{p}}\leq N\|b\|_{\dot E_{q,1}}\|Du\|_{L_{p}}$, which has the following
``parabolic'' analog.

\begin{lemma}
                        \label{lemma 1.1.1}

 Let $u,f,f_{i}$,
$i=1,...,d$, be  
function on $\bR^{d+1}$ of class $L_{1,\loc}$
such that (in the sense of generalized functions)
\begin{equation}
                            \label{1.5.10}
\text{either}\quad
\partial_{t}u=f+D_{i}f_{i},\quad\text{or}
\quad u\geq0\quad\text{and}\quad 
\partial_{t}u\leq f+D_{i}f_{i}.
\end{equation}
Let $q\in(2,d+2]$, $\rho\in(0,\infty)$,
\index{$S$@Miscelenea!$\bar b_{R}$@$\hat b_{q,\rho}$}%
and
$$
\hat b_{q,\rho}:=\sup_{r\leq\rho}r\sup_{C\in\bC_{r}}\dashnorm b\|_{L_{q}(C)}<\infty.
$$
Let $C\in \bC_{\rho }$ and let $\zeta\in C^{\infty}_{0}(C)$
be a nonnegative functions with the integral
of its square equal to one.
Then  
$$
\int_{\bR^{d+1}}b^{2}\zeta^{2}u^{2}\,dxdt
\leq N\hat b_{q,\rho }
 \int_{\bR^{d+1}}  \zeta^{2} |Du|^{2}
\,dxdt
$$
\begin{equation}
                            \label{1.5.2}
+N\hat b_{q,\rho }\rho ^{-d-4}
\int_{\bR^{d+1}}u^{2}I_{C}\,dxdt
+N\hat b_{q,\rho }\rho ^{-d-2}
\int_{\bR^{d+1}}I_{C} \big[\rho ^{2} f^{2}  +\sum_{i}f_{i}^{2}\big]\,dxdt,
\end{equation}
where the constants $N$ depend only on $d,q$.

\end{lemma}  

Proof. We may assume that the right-hand side of \eqref{1.5.2} is finite and then by using mollifiers we reduce the general case to the one
in which $u,f,f_{i}$ are smooth and bounded.
Then we have that either
$$
\partial_{t}(\zeta u)+\Delta(\zeta u)+F=0,
$$
or $\zeta u\geq0$ and above we have $\leq0$
in place of $=0$,
where
$$
F=-u\partial_{t}\zeta-\zeta f-\zeta D_{i}f_{i}-\Delta(\zeta u)=:g_{1}+...+g_{4} .
$$
It follows (for instance, by It\^o's formula) that the left-hand side of \eqref{1.5.2} is dominated by a constant times the sum $G_{1}+...+G_{4}$, where
$$
G_{i}=\int_{\bR^{d+1}}b^{2}I_{C  }P^{2}_{2,4}g_{i}.
$$

Before proceeding further we note that we may
look at $\zeta$ as a scaled and translated
function with support in $C_{1}$. Then it is seen that
$$
\rho ^{d+2} |\zeta |^{2}
+
\rho ^{d+4} |D\zeta |^{2}
+
\rho ^{d+6} |\partial_{t}\zeta |^{2}
 \leq N (d).
$$
Furthermore, it is easy to see that $bI_{C}
\in \dot E_{q,q,1}$ and 
$
\|I_{C}b\|_{\dot E_{q,q,1}}\leq
\hat b_{q ,\rho  } .
$

Next, owing to Theorem \ref{theorem 5.25,1}  and Corollary \ref{corollary 10.5,1} and the fact that
$P_{2,4}=NP_{1,4}P_{1,4}$
$$
G_{1}\leq N\hat b_{q,\rho }
\int_{\bR^{d+1}}P^{2}_{1,4}(I_{C}u\partial_{t}
\zeta)\,dxdt
$$
$$
\leq
 N\hat b_{q,\rho }\rho ^{2}
\int_{\bR^{d+1}}u^{2}|\partial_{t}
\zeta|^{2}\,dxdt
\leq N\hat b_{q,\rho }\rho ^{-d-4}
\int_{\bR^{d+1}}u^{2}I_{C}\,dxdt.
$$
Similarly,
$$
G_{2}\leq N\hat b_{q,\rho }
\int_{\bR^{d+1}}P^{2}_{1,4}(I_{C} 
\zeta f)\,dxdt\leq N\hat b_{q,\rho }
\rho ^{2}\int_{\bR^{d+1}}\zeta^{2}f^{2}\,dxdt
$$
$$
\leq N\hat b_{q,\rho }
\rho ^{-d}\int_{\bR^{d+1}}I_{C}f^{2}\,dxdt.
$$
While estimating $G_{3}$ we use that
$P_{2,4}(\zeta D_{i}f_{i})=D_{i}P_{2,4}(\zeta f_{i})-P_{2,4}(f_{i}D_{i}\zeta_{i})$ and that
$|DP_{2,4}h|\leq NP_{1,8}|h|$. This yields
$$
G_{3}\leq N\hat b_{q,\rho }\rho ^{-d-2}
 \int_{\bR^{d+1}} \sum_{i}f_{i}^{2}\,dxdt .
$$
Finally,
$$
G_{4}\leq N\hat b_{p_{0},\rho_{0}}
 \int_{\bR^{d+1}} |uD\zeta+\zeta Du|^{2}
\,dxdt\leq N\hat b_{p_{0},\rho_{0}}
 \int_{\bR^{d+1}}  \zeta^{2} |Du|^{2}
\,dxdt
$$
$$
+N\hat b_{p_{0},\rho_{0}}\rho_{0}^{-d-4}
 \int_{\bR^{d+1}} u ^{2}I_{C}
\,dxdt.
$$
The lemma is proved. \qed

Taking $\zeta$ in the form $\rho^{-(d+2)/2}
\xi(t/\rho^{2},x/\rho)$ and sending
$\rho\to\infty$ we arrive at the following result.

\begin{corollary}
                    \label{corollary 1.12,1}
Under the assumptions of Lemma
\ref{lemma 1.1.1} suppose that $f=0$ and
$
 \rho^{-2}u^{2} I_{C_{\rho}}\to 0
$
in $L_{1} $ as $\rho\to\infty$
(for instance,  $u\in L_{2} $).
Then 
$$
\int_{\bR^{d+1}}b^{2} u^{2}\,dxdt
\leq N(d,q)\|b\|_{\dot E_{q,q,1} }
 \int_{\bR^{d+1}} \big[  |Du|^{2}  +\sum_{i}f_{i}^{2}\big]\,dxdt.
$$

\end{corollary}

We will also need some results about $t$-traces
of $P_{\alpha,k}f$. We show that the $t$-traces
of functions in $E^{1,2}_{q,p,\beta}$
possess some regularity as $L_{p}$-functions. For $\gamma=0$ or $1$ set
$$
D^{\gamma}=D\quad\text{if}\quad
\gamma=1\quad\text{and}\quad D^{\gamma}=1\quad
\text{if}\quad \gamma=0.
$$
Below by
$D^{\gamma}u(0,\cdot)$ we mean the limit
in $L_{r}(B)$ for any ball $B$
of $D^{\gamma}u^{(\varepsilon)}(0,\cdot)$
as $\varepsilon\downarrow 0$.
The existence of this limit easily follows from Lemma \ref{lemma 6,17.1} and Corollary \ref{corollary 6,19.1} below.
By the way, note that according to
Lemma \ref{lemma 11.15,1}  the functions in $E^{1,2}_{q,p,\beta}$ are bounded and continuous.

\begin{theorem}
              \label{theorem 6,6,1}
Take  
$r\in[ p,\infty)$, $\mu>0$    and assume
that  
$$
2< \beta+\gamma \leq \gamma+\frac{d}{p}+\frac{2}{q}
<2 +\frac{d}{r},
\quad \kappa:=\gamma+\frac{d}{p}
+\frac{2}{q}-\frac{d}{r}\leq \mu
<2 .
$$

 Then for any
$u\in E_{q,p,\beta}$ the trace $D^{\gamma}u(0,\cdot)$ is uniquely defined
and for any $\varepsilon>0$
\begin{equation}
                     \label{6,6.1}
\|D^{\gamma}u(0,\cdot)\|_{E_{r,\beta+\gamma-\mu}(\bR^{d})}\leq N\varepsilon \|\partial_{t}u,
D^{2}u\|
_{ E _{p,q,\beta}}
+N\varepsilon^{-\mu/(2-\mu)}
\|u\|
_{ E _{p,q,\beta}},
\end{equation} 
\begin{equation}
                     \label{8,27.1}
\|D^{\gamma}u(0,\cdot)\|_{E_{r,\beta+\gamma-2}(\bR^{d})}\leq N  \|u\|
_{ E^{1,2} _{p,q,\beta}},
\end{equation} 
where  the constants $N$ depend only
on $d,p,q,\beta,\mu$.

\end{theorem}

\begin{remark}
                 \label{remark 11.21,1}
If $\gamma=0$ we are dealing with $u(0,x)$
and $\beta>2$. We know from Lemma 
\ref{lemma 11.15,1} that, if $\beta<2$, then $u$ is bounded.
\end{remark}

In the proof of this theorem we may and will
assume without losing generality that 
$u(-t,x)=u(t,x)$. This theorem is proved in
\cite{Kr_25} in case of norms defined
as usual in \eqref{3.27.3}, so until the 
end of this section we will be dealing
with the norms defined in \eqref{4.3.2}.

\begin{remark}
               \label{remark 11.29.1}
Obviously, $E^{1,2}_{q,p,\beta}\subset
W^{1,2}_{q,p,\loc}$ and one can show that
$E^{1,2}_{q,p,\beta}\not\subset
W^{2,1}_{q+\varepsilon,p+\varepsilon,\loc}$ no matter how small $\varepsilon>0$ is.
Therefore, in terms of the local summability of derivatives, general functions in 
$E^{1,2}_{q,p,\beta}$ are not much better
than those in $W^{1,2}_{q,p,\loc}$. For the latter
class the trace theorems for $\gamma=1$ (see Lemma 
\ref{lemma 6,17.1}) only, basically, can guarantee that $Du(0,\cdot)\in L_{r,\loc}$
and, if $r<d$. This does not yield even the boundedness of $u(0,\cdot)$. At the same time \eqref{8,27.1} and the Morrey
theorem (see, for instance, Theorem 10.2.1
in \cite{Kr_08}) imply that $u(0,\cdot)$
is $2-\beta$ H\"older continuous, provided
$1<\beta<2$, that is, almost Lipschitz continuous. Of course, this is at the expense of $u\in E^{1,2}_{q,p,\beta}$.

\end{remark}

\begin{remark}
               \label{remark 8,25.1}
For the probability part in this book 
 the most important particular case
of \eqref{6,6.1} is when $\gamma=1$,
$r=p$ (and $q>2$), $\mu=\kappa=1+2/q$,
and $1<\beta\leq d/p+2/q$.
 
 In that case
for any
$u\in E^{1,2}_{q,p,\beta}$  and any $\varepsilon>0$ (observe that $\beta>\beta+1-\kappa$)
\begin{equation}
                     \label{6,6.10}
\|D u(0,\cdot)\|_{E_{p,\beta } }\leq N\varepsilon \|\partial_{t}u,
D^{2}u\|
_{ E _{p,q,\beta}}
+N\varepsilon^{-(q+2)/(q-2)}
\|u\|
_{ E _{p,q,\beta}}.
\end{equation} 

Another case, we are going to use, is when
$\gamma=1$, 
$q>2$, $1<\beta\leq d/p+2/q$,
$$
r=\frac{p^{2}}{d}\Big(\frac{d}{p}-\frac{2}{q}+1\Big),\quad \mu=\kappa=1+\frac{2}{q}+
\frac{d}{p}\Big(1-\frac{2}{q}\Big)
\Big(\frac{d}{p}+1-\frac{2}{q}\Big)^{-1}.
$$
In that case as is easy to see $r>p$ and
\eqref{8,27.1} says that
$$
\|D^{\gamma}u(0,\cdot)\|_{E_{r,\beta-1}(\bR^{d})}\leq N  \|u\|
_{ E^{1,2} _{q,p,\beta}}.
$$
 
\end{remark}

To prove Theorem \ref{theorem 6,6,1},
first, we need the following corollary of
Theorem 10.2 of \cite{BIN_75} which
we give with a different proof for completeness.

\begin{lemma}
                 \label{lemma 6,17.1}
Let  
$r\geq p$,
$$
 \gamma+\frac{d}{p}+\frac{2}{q}
<2 +\frac{d}{r},
\quad \kappa:=\gamma+\frac{d}{p}
+\frac{2}{q}-\frac{d}{r}\quad (
<2 ).
$$    
 Then for any
$u\in W^{1,2}_{q,p}$ and $\varepsilon>0$ we 
have
\begin{equation}
                            \label{6,17.2}
\|D^{\gamma}u(0,\cdot)\|_{L_{r}}\leq
N\varepsilon \|\partial_{t}u,D^{2}u
\|_{L_{q,p} }+N\varepsilon^{-\kappa/(2-\kappa)}
\|u
\|_{L_{q,p} }.
\end{equation}
\end{lemma}

Proof. The case of arbitrary $\varepsilon>0$
is reduced to that of $\varepsilon=1$
by using   self-similarity.
To treat $\varepsilon=1$ take $\zeta\in
C^{\infty}_{0}(\bR)$ such that $\zeta(t)=1$
for $|t|\in[0,1]$, $\zeta(t)=0$ for $|t|\geq 2$, and define $-f=
\partial_{t}(\zeta u)+\Delta(\zeta u)$.  
We know (from It\^o's formula or
from PDEs) that $\zeta u=P_{2,4}f$. It follows that
$
|D^{\gamma}u(0,x)|\leq NP_{2-\gamma,8}|f|(0,x)$. By H\"older's
inequality we get that
$$
 P_{2-\gamma,8}|f|(0,x) \leq \int_{\bR^{d}}I_{1}(y)I_{2}(x-y)\,dy,
$$
where
$$
I_{1}^{q'}(y)=\int_{0}^{2}p_{2-\gamma,8}^{q'} (t,y)\,dt,\quad I_{2}^{q}=\int_{0}^{2} 
|f(t,y)|^{q}\,dt,\quad q'=q/(q-1).
$$
We see that we are dealing with the $L_{r}$-norm of a convolution.
By Young's
inequality  
$$
\|F\|_{L_{r}}\leq
\|f \|_{L_{q,p}}\|I_{1}
\|_{L_{s}},
$$
where $1/s=1+1/r-1/p$ ($ \leq 1$ since $r\geq p$). 

To estimate the last norm observe that
for certain constants $N_{1},N_{2}$
$$
I_{1}^{q'}(y) 
=N_{1}|y|^{2-(d+\gamma)q'}
\int_{0}^{N_{2}/|y|^{2}}
t^{-(d+\gamma)q'/2}e^{- 1/t }\,
dt.
$$
Since as easy to see $(d+\gamma)q'/2>1$
($d\geq2$), the integral is a bounded function
of $|y|$ which tends to zero as $|y|\to\infty$ faster than $1/|y|^{n}$ for any $n$. Furthermore,
$$
\int_{B_{1}}|y|^{s(2/q'-d-\gamma)}\,dy
<\infty
$$
because $s(2/q'-d-\gamma)>-d $, which is equivalent to $\gamma+d/p+2/q<2+d/r$.
Hence, $\|I_{1}
\|_{L_{s}}<\infty$ and this immediately leads to
\eqref{6,17.2}. \qed 

\begin{corollary}
           \label{corollary 6,19.1}
For any $\rho>0,\varepsilon>0$ and
$u\in W^{1,2}_{q,p}(C_{2\rho})$ we have
$$
\dashnorm D^{\gamma} u(0,\cdot)\|_{L_{r}(B_{\rho})}
\leq N\varepsilon  
\dashnorm \partial_{t}u,D^{2}u
\|_{L_{q,p}(C_{2\rho})}
$$
$$
+N
(\varepsilon \rho^{-2}+\varepsilon^{-\kappa/(2-\kappa)}\rho^{(2\kappa-2\gamma)/(2-\kappa)})
\dashnorm u
\|_{L_{q,p}(C_{2\rho})}.
$$
\end{corollary}

Indeed, the case of arbitrary $\rho>0$ is reduced
to $\rho=1$ by means of parabolic dilation. In the latter case
take $\zeta\in C^{\infty}_{0}
(\bR^{d+1})$ such that $\zeta=1 $ on $C_{1}$ and 
$\zeta=0$ in $\bR^{d+1}_{0}
 \setminus  C_{2}$. Then  use
\eqref{6,17.2} to see that
$$
\|D^{\gamma}u(0,\cdot)\|_{L_{r}(B_{1})}
\leq N\varepsilon  
\|\partial_{t}(\zeta u),D^{2}(\zeta u)
\|_{L_{q,p} }+N\varepsilon^{-\kappa/(2-\kappa)}
\|u
\|_{L_{q,p}(C_{2}) }
$$
$$
\leq N\varepsilon  
\|\partial_{t}u,D^{2}u
\|_{L_{q,p}(C_{2}) }+N\varepsilon \|u,Du\|_{L_{q,p}(C_{2}) }
+N\varepsilon^{-\kappa/(2-\kappa)}
\|u
\|_{L_{q,p}(C_{2}) }.
$$
After that it only remains to use
the interpolation inequality
$$
\|Du\|_{L_{q,p}(C_{2}) }\leq
\|\partial_{t}u,D^{2}u
\|_{L_{q,p}(C_{2}) }+
N\|u
\|_{L_{q,p}(C_{2}) }.
$$

In the following Lemma 4.1.21
of \cite{Kr_25} $\gamma$ can be any
number in $[0,d+2)$.  

\begin{lemma}
                  \label{lemma 11.28.1}

Let $0\leq\gamma<\beta\leq d+2$, $k>0$. Then there exists 
 a constant  $N$
  such that for any $f\geq0$
and $\rho\in(0,\infty)$ we have
\begin{equation}
                             \label{1.17.2}
P_{\gamma,k}(I_{C^{c} _{ \rho}}f)(0)
\leq N\rho^{\gamma-\beta}\bM_{\beta} 
f(0) .
\end{equation}
\end{lemma}

To prove  Theorem \ref{theorem 6,6,1} we also need its homogeneous 
version for the homogeneous Morrey 
space $\dot E^{1,2}_{q,p,\beta} $.

\begin{lemma}
               \label{lemma 6,18.1}
Let   
$r\geq p$  and let
$$
2< \beta +\gamma\leq \gamma+\frac{d}{p}+\frac{2}{q}
<2 +\frac{d}{r}.
$$ 
Then for any
$u\in \dot E^{1,2}_{q,p,\beta}$ its trace $u(0,\cdot)$ is uniquely defined
and  
\begin{equation}
                     \label{6.6.1}
\|D^{\gamma}u(0,\cdot)\|_{\dot E_{r,\beta+\gamma-2} }\leq N \|\partial_{t}u,
D^{2}u\|
_{ \dot E _{p,q,\beta}},
\end{equation} 
where the constant  $N$ depends only
on $d,p,q,r,\beta$.

\end{lemma}

Proof. Take $\zeta\in C^{\infty}_{0} (\bR^{d+1})$, such that $\zeta(0)=1\geq \zeta\geq0$, define $\zeta_{n}(t,x )=
\zeta(t/n^{2},x/n )$ and observe that,
as $n\to\infty$,
$$
 \big|\|
\partial_{t}(\zeta_{n}u)\|_{\dot E_{q,p,\beta}} -\|\zeta_{n}
\partial_{t}u\|_{\dot E_{q,p,\beta} }\big|\leq n^{-2}\sup |\partial_{t}\zeta|
\|
u\|_{\dot E_{q,p,\beta}} \to 0.
$$
Also
$$
 \big|\|
D(\zeta_{n}u)\|_{\dot E_{q,p,\beta}} -\|\zeta_{n}
Du\|_{\dot E_{q,p,\beta} }\big|\leq n^{-1}\sup |D\zeta|
\|
u\|_{\dot E_{q,p,\beta}} \to 0,
$$
$$
 \big|\|
D^{2}(\zeta_{n}u)\|_{\dot E_{q,p,\beta}} -\|\zeta_{n}
D^{2}u\|_{\dot E_{q,p,\beta} }\big|\leq n^{-2}\sup |D^{2}\zeta|
\|
u\|_{\dot E_{q,p,\beta}}
$$
$$
+2n^{-1}\sup |D\zeta|\|
Du\|_{\dot E_{q,p,\beta} } \to 0.
$$

It follows that it suffices to concentrate on $u$ that vanish
for large $|t|+|x|^{2}$. In that case set $-f=\partial_{t}u+\Delta u$.
To further reduce our
problem observe that using translations show that it suffices
to prove that for any
  $\rho>0$,
$$
\rho^{\beta+\gamma-2}\dashnorm D^{\gamma}u(0,\cdot)\|_{L_{r}( B_{\rho})}\leq N\sup_{\rho_{1}\geq \rho}
\rho^{\beta }_{1}\dashnorm f\|
_{L_{q,p} (C_{\rho_{1}})   )}
$$
\begin{equation}
                     \label{6.19.1}
=  N\sup_{\rho_{1}\in[\rho,\rho+\rho_{2}]}
\rho^{\beta }_{1}\dashnorm f\|
_{L_{q,p} (C_{\rho_{1}})   },
\end{equation}
where $\rho_{2}$ is such that
$u(t,x)=0$ for $|t|+|x|^{2}\geq \rho_{2}^{2}$  
and the last equality is due to
$\beta \leq d/p+2/q$.

It is easy to pass to the limit
in \eqref{6.19.1} from smooth functions to arbitrary ones in
$W^{1,2}_{q,p}(C_{\rho_{2}})\supset
E^{1,2}_{q,p,\beta}(C_{\rho_{2}})$.
Therefore, we may assume that $u$ is smooth.
We thus reduced the general case to the task of proving the first estimate
in  \eqref{6.19.1} for smooth $u$ with compact support. One more reduction
is achieved by using the self-similarity
which shows that we only need to concentrate on $\rho=1$, that is,
we only need to prove
\begin{equation}
                    \label{6,19,1}
 \dashnorm D^{\gamma}u(0,\cdot)\|_{L_{r}( B_{1})}\leq N\sup_{\rho \geq 1}
\rho^{\beta} \dashnorm f\|
_{L_{q,p} (C_{\rho })   )}
\end{equation}
for smooth $u$ with compact  support.

  Now   define
  $
g=|f|I_{C_{2}},h=|f|I_{C_{2}^{c}}$.
As it follows from the proof of Lemma~
\ref{lemma 6,17.1},
$$
|D^{\gamma}u(0,x)|\leq NG_{\gamma}(x)+NH_{\gamma}(x),
$$
where
$$
 (G_{\gamma},H_{\gamma})(x) =
\int_{0}^{\infty} 
\int_{\bR^{d}}P_{2-\gamma,8}(t,x-y)(g,h)(t,y)\,dydt.
$$

Estimate  \eqref{1.17.2}
implies that
that for $|x|\leq 1$
$$
H_{\gamma}(x)\leq N\sup_{\rho>1}\rho^{\beta}
\dashint_{(0,x)+C_{\rho}}h\,dydt
\leq N\sup_{\rho>1}\rho^{\beta}
\dashint_{ C_{2\rho}}h\,dydt
$$
$$
\leq N\sup_{\rho \geq1}
\rho^{\beta} \dashnorm f\|
_{L_{q,p}  (C_{\rho} )},
$$
where the last inequality is
due to H\"older's inequality.
Hence,
\begin{equation}
                     \label{6,8,3}
 \| H_{\gamma}\|_{L_{r}(  B_{1})}\leq N\sup_{\rho \geq1}
\rho^{\beta} \dashnorm f\|
_{L_{q,p}  (C_{\rho} )}. 
\end{equation}

Then we get the estimates
$$
 \| G_{\gamma}\|_{L_{r}(\bR^{d}) }\leq N
  \|g\|
_{L_{q,p}}\leq   N\sup_{\rho \geq1}\rho^{\beta}
  \|f\|
_{L_{q,p}(C_{\rho})}
$$
as in the proof of Lemma \ref{lemma 6,17.1}.
This and \eqref{6,8,3} prove
\eqref{6,19,1} and the lemma. \qed

{\bf Proof of Theorem \ref{theorem 6,6,1}}.  To prove \eqref{6,6.1}, it suffices to show that
for any $\rho\in (0, 1]$,
$\varepsilon>0$ 
\begin{equation}
                     \label{6,20.1}
I_{\rho}:=\rho^{\beta+\gamma-\mu}\dashnorm Du(0,\cdot)\|_{L_{r}(B_{\rho})}\leq N\varepsilon \|\partial_{t}u,
D^{2}u\|
_{ E _{p,q,\beta}}
+N\varepsilon^{-\mu/(2-\mu)} 
\|u\|
_{ E _{p,q,\beta}}.
\end{equation} 

By Corollary
\ref{corollary 6,19.1} with
$\epsilon=\varepsilon \rho^{ \gamma-\mu }$
in place of $\varepsilon$  we get
$$
 I_{\rho}\leq N\epsilon \|\partial_{t}u,
D^{2}u\|
_{ E _{p,q,\beta}}+N\big(\epsilon  \rho^{-2} +\epsilon^{-\kappa/(2-\kappa)}\rho^{(2\kappa-2\mu)/(2-\kappa)
}\big)\| u\|_{E_{q,p,\beta}}.
$$
For $\epsilon< \rho^{2-\mu}$ this yields (here we use that $ \kappa\leq \mu<2$)
$$
 I_{\rho}\leq N\epsilon \|\partial_{t}u,
D^{2}u\|
_{ E _{p,q,\beta}}+N\epsilon^{-\mu/(2-\mu)}\| u\|_{E_{q,p,\beta}}.
$$

In the remaining case $\rho^{2-\mu}\leq \epsilon$. In that case for $\zeta\in
C^{\infty}_{0}((-1,1)\times B_{2})$   such that $\zeta
=1$ on $C_{1}$ we have by Lemma \ref{lemma 6,18.1} that
\begin{equation}
                     \label{8,27.2}
I_{\rho}\leq \rho^{2-\mu} \rho^{\beta+\gamma-2}
\dashnorm D(\zeta u)(0,\cdot)\|_{L_{r}(B_{\rho})}\leq N\epsilon
\|\partial_{t}(\zeta u),D^{2}(\zeta u)
\|_{\dot E_{q,p,\beta}}.
\end{equation}
Owing to $\beta\leq d/p+2/q$,
 the last norm here is easily shown to be less than
$$
N\|\partial_{t} u,D^{2}u
\|_{  E_{q,p,\beta} }
+N\| u\|_{  E_{q,p,\beta}}.
$$ 
Therefore, \eqref{6,20.1} holds in this case as well and this proves
\eqref{6,6.1}. Estimate \eqref{8,27.1}
follows from \eqref{8,27.2} with
$\mu=2$ and $\epsilon=1$.
The theorem is proved. \qed

For $\rho>0$, $p ,q  \in(1,\infty)$ introduce 
\index{$S$@Miscelenea!$a_{\pm}$@$a^{\shharp}_{\rho}$}%
\begin{equation}
                         \label{6.3.01}
 a^{\shharp}_{\rho} = \sup
_{r\leq\rho}
 \sup
_{ C\in\bC_{r} }\dashint_{C}|a(t,x)-\tilde a_{C}(t)|\,dxdt ,
\end{equation} 
\index{$S$@Miscelenea!$\bar b_{\rho_{b}}$@$\hat b_{p,q,\rho}$}%
\begin{equation}
                           \label{3.14.2}
\hat b_{p ,q ,\rho }=\sup_{r\leq\rho }r
\sup_{C\in \bC_{r}} 
\dashnorm b \|_{L_{p ,q  }(C)}. 
\end{equation}

Fix  $p ,q ,\beta    $ such that
\begin{equation}
                        \label{3.21.010}
p ,q  \in(1,\infty),\quad 1<\beta\leq\frac{d}{p }+\frac{2}{q } .
\end{equation}

Fix some
$$
\rho_{a},\rho_{b}\in(0,1].
$$

Here is Theorem 2.1 of \cite{Kr_27},
which is a close restatement of Theorem 3.5 of \cite{Kr_23_1} where $\cL$ contains
also singular zeroth-order term.
The objects
$$
\hat a=\hat  a(d,\delta,q,p,\beta)>0,\quad\hat b=
\hat  b(d,\delta,q,p,\beta,\rho_{a} )>0,
$$ 
$$\lambda_{0}=  \lambda_{0}(d,\delta,q,p,\beta,\rho_{a} ) >0,\quad
N_{1}=N_{1}(d,\delta,q,p,\beta,\rho_{a}  )
$$
below are taken from Theorem 2.1 of \cite{Kr_27}.

 \begin{theorem}
                     \label{theorem 5.8,20}
Suppose that
\index{$S$@Miscelenea!$a_{\pm}$@$\hat a$}%
\index{$S$@Miscelenea!$\bar b_{R}$@$\hat b$}%
\begin{equation}
                              \label{12.5,1}
a^{\shharp}_{\rho_{a}}\leq 
\hat  a,\quad\hat b_{q\beta,p\beta,\rho_{b}}\leq \hat  b.
\end{equation}
 Then for any $u\in E^{1,2}_{q,p,\beta}$ 
and  $\lambda\geq \lambda_{0}\rho_{b}^{-2}$ 
\begin{equation}
                        \label{5.10.2}
\|\lambda u,\sqrt\lambda Du, \partial_{t}u, D^{2}u\|_{  E_{q,p,\beta}}
\leq N _{1}\rho_{b}^{-\alpha}\|f\|_{  E_{q,p,\beta}}, 
\end{equation}
where   
$$
f=\cL  u-\lambda u:=\partial_{t}u+(1/2)a^{ij}D_{ij}u+b^{i}D_{i}u-\lambda u,\quad \alpha=
d+2+\beta-\frac{d}{p}-\frac{2}{q}.
$$
Furthermore, for any
$f\in E_{q,p,\beta}$ and $\lambda\geq\lambda_{0}\rho_{b}^{-2}$ there exists a unique
$u\in E^{1,2}_{q,p,\beta}$ such that in
$\bR^{d+1}$
\begin{equation}
                          \label{10.14,1}
\cL u -\lambda u=-f. 
\end{equation}
 
\end{theorem}

It is important to have in mind that if
$\beta<2$ (our main case) and $u\in E^{1,2}_{q,p,\beta}$,
then according to Lemma \ref{lemma 11.15,1},
$u$ is bounded and continuous. 

\begin{remark}
                          \label{remark 12.2,1}
We need a few additions to this theorem.
First, note that, if $f(t,x)=0$ for $t\geq T$
and all $x$, then $u(t,x)=0$ for $t\geq T$
and all $x$. This is because this is true
if $\cL=\partial_{t}+\Delta$ and for general $\cL$
is obtained either by perturbation method or the method of continuity preserving this property.

Next, if real valued Borel $c(t,x)$ is such that
$2N_{1}|c|\leq \lambda_{0}\rho_{b}^{-2}$,
then the assertions of Theorem \ref{theorem 5.8,20} remain true for 
$$
\cL u=\partial_{t}u+(1/2)a^{ij}D_{ij}u+b^{i}D_{i}u
+cu.
$$
One need  only replace $N_{1}$ with $2N_{1}$
in \eqref{5.10.2}. This is proved by perturbation method.

Finally, if $f(t,x)=0$ for $t\geq T$
and all $x$, then $u$ is the solution
of the Cauchy problem for $t<T$ with terminal data
$u(T,\cdot)=0$. In the future we will need
more general data, say $g\in C^{\infty}_{0}(\bR^{d})$. In that case redefine $a$, $b$, and $c$ for
$t\geq T$ by setting $a=(\delta^{ij}),b=0,c=0$ and
for $f\in E_{q,p,\beta}$ consider the equation
$$
\cL u-\lambda u=-fI_{t<T}+I_{t\geq T}(  \Delta -
\lambda) g ,
$$
where $\cL$ is as above (containing $c$). The solution
$u$ will be of class $E^{1,2}_{q,p,\beta}$
and will be equal $g(x) $ for $t\geq T$
(consider $u-g $), so that  $u(t,x)=g(x)$
for $t\geq T$. Such a $u$ is unique.

\end{remark}

\begin{remark}
                       \label{remark 3.16,1}
One more comment to make is about
uniqueness of solutions when $p,q,\beta$
vary. Denote by $\sfA$ 
\index{$B$@Sets!$\sfA$}%
the collection
of $(a,b,p,q,\beta)$, where $a$ is $\bS_{\delta}$-valued as usual,  satisfying  
\eqref{3.21.010} and \eqref{12.5,1} and
suppose that
$$
 (a,b,p',q',\beta') ,(a,b,p'',q'',\beta'')\in \sfA.
$$
Then by Theorem \ref{theorem 5.8,20}
for any $f\in E^{1,2}_{q',p',\beta'}\cap
E^{1,2}_{q'',p'',\beta''}$ and
$$
\lambda\geq \rho_{b}^{-2}\big
(\lambda(d,\delta,p',q',\beta',\rho_{a} )\vee
\lambda(d,\delta,p'',q'',\beta,''\rho_{a} )\big)
$$
equation \eqref{10.14,1} has a unique
solution $u'\in E^{1,2}_{q',p',\beta'}$ and
a unique
solution $u''\in E^{1,2}_{q'',p'',\beta''}$.
An important fact is that $u'=u''$.

This follows from the fact that the explicit
formulas for solutions of the heat equation
shows that indeed $u'=u''$ in that case and in the general case this is seen from the method
of continuity applied in $E^{1,2}_{q',p',\beta'}$ or $E^{1,2}_{q'',p'',\beta''}$.

\end{remark}

\begin{remark}
                          \label{remark 12.5,3}
If $\beta<2$ and $q>2$, then by combining \eqref{5.10.2}
with Lemma \ref{lemma 11.15,1} and Remark \ref{remark 8,25.1} we get
$$
\sup_{\bR^{d+1}}|u(t,x)|
+\sup_{t\in\bR}\|Du(t,x)\|_{E_{r,\beta-1}}\leq
N\|f\|_{E_{q,p,\beta}},
$$
where $N$ depends only on $d,\delta,p,q,\beta,\lambda,\rho_{a},\rho_{b} $ and $r>p$ is from
Remark \ref{remark 8,25.1}.

\end{remark}

A useful addition to Theorem \ref{theorem 5.8,20} is the following result before which
we introduce new spaces. Let a domain $\cO\subset
\bR^{d+1}$. Fix a $\chi_{x}\in C^{\infty}_{0}(\bR^{d}),\chi_{t}\in C^{\infty}_{0}(\bR)$, such that $\chi_{x}=1$ on  $B_{1}$, $\chi_{t}=1$
on $(-1,1)$ and $0\leq \chi_{x},\chi_{t}\leq 1$, and set $\chi_{tn}(t)=\chi_{t}(t/n)$, $\chi_{xn}(x)=\chi_{x}(x/n)$, $\chi_{n} 
=\chi_{tn} \chi_{xn} $, 
$$
\EO_{q,p,\beta}(\cO)=\{u\in E_{q,p,\beta}(\cO):\lim_{n\to\infty}\|u\chi_{n}-u\|_{E_{q,p,\beta}(\cO)}=0\},
$$    
\begin{equation}
                              \label{12.5,2}
\EO^{1,2}_{q,p,\beta}(\cO)=\{u\in E _{q,p,\beta}(\cO):\partial_{t}u,D^{2}u,Du,u
\in \EO_{q,p,\beta}(\cO)\}.
\end{equation}
\index{$A$@Sets of functions!$\EO_{q,p,\beta}$}%
\index{$A$@Sets of functions!$\EO^{1,2}_{q,p,\beta}$}

Observe that
$\EO^{1,2}_{q,p,\beta}(\cO)$ can be equivalently defined as the subset of $E^{1,2}_{q,p,\beta}(\cO)$
of functions $u$ such that $\|u\chi_{n}-u\|_{E^{1,2}_{q,p,\beta}}\to0$ as $n\to \infty$.
This equivalence is an easy consequence of the
formulas
$$
\chi_{n}Du=D(u\chi_{n})-uD\chi_{n},
\quad \chi_{n}\partial_{t}u=\partial_{t}(u\chi_{n})-u\partial_{t}\chi_{n},
$$
$$
\chi_{n}D_{ij}u=D_{ij}(u\chi_{n})
-D_{i}\chi_{n}D_{j}u-D_{j}\chi_{n}D_{i}u
-uD_{ij}\chi_{n}
$$
and the fact that $|D\chi_{n}|+n|D^{2}\chi_{n}|+n|\partial_{t}\chi_{n}|\leq N(d)n^{-1}$.
If $\cO=\bR^{d+1}$ we drop $\cO$ in 
$\EO _{q,p,\beta}(\cO)$ and $\EO^{1,2}_{q,p,\beta}(\cO)$.

 The space $\EO_{p,\beta}$ is defined as the
 \index{$A$@Sets of functions!$\EO_{p,\beta}$}%
  subspace of $\EO_{q,p,\beta}$ consisting of
   functions independent of $t$.

\begin{theorem}
                       \label{theorem 7,22.1}
If in Theorem \ref{theorem 5.8,20}
we have $f\in \EO_{q,p,\beta}$, then the unique
$E^{1,2}_{q,p,\beta}$-solution $u$ of
\eqref{10.14,1} belongs to
$\EO^{1,2}_{q,p,\beta}$.
\end{theorem}
 
Proof. Observe that $u_{n}:=u\chi_{n}$ satisfies
$$
\partial_{t}u_{n}+\cL u_{n}- \lambda  u_{n}=-f\chi_{n}+u\partial_{t}\chi_{n}+a^{ij}D_{i}uD_{j}\chi_{n}
+u\cL \chi_{n}.
$$
Here $|D\chi_{n}|\leq N(d)/n$, so that
$$
\|a^{ij}D_{i}uD_{j}\chi_{n}\|_{E_{q,p,\beta}}
\leq Nn^{-1}\|u\|_{E^{1,2}_{q,p,\beta}}\to0
$$
as $n\to\infty$. Similarly,
$$
\|u\partial_{t}\chi_{n}+ua^{ij}D_{ij}\chi_{n}\|_{E_{q,p,\beta}}
\leq Nn^{-2}\|u\|_{E^{1,2}_{q,p,\beta}}\to0.
$$
Also 
$$
\|ub^{i}D_{i}\chi_{n}\|_{E_{q,p,\beta}}
\leq Nn^{-1}\|u|b|\, \|_{E_{q,p,\beta}}
\leq Nn^{-1}\|u,Du\|_{E_{q,p,\beta}}\| b\, \|_{E_{\beta q,\beta p,1}}
\to 0.
$$

Hence, in light of \eqref{5.10.2}
$$
(\partial_{t},D^{2},D,1)u_{n}\to
(\partial_{t},D^{2},D,1)u 
$$
in $E_{q,p,\beta}$ as $n\to\infty$. 
The theorem is proved. \qed

\mysection[Uniqueness of weak solutions]{Uniqueness of weak solutions}
                     \label{section 3.22,1}

Our basic assumption here are the same as in 
\eqref{3.21.010} and Theorem \ref{theorem 5.8,20}
 with addition that $\beta<2$, that is
 ($a$ is $\bS_{\delta}$-valued and) $(a,b,q,p,\beta)\in\sfA\cap\{\beta<2\}$:
\begin{equation}
                             \label{4.11,1}
p ,q  \in(1,\infty),\quad 1<\beta\leq\frac{d}{p }+\frac{2}{q } ,
\quad \beta<2,
\end{equation}
\begin{equation}
                             \label{4.11,3} 
a^{\shharp}_{\rho_{a}}\leq 
\hat  a(d,\delta,q,p,\beta),\quad\hat b_{q\beta,p\beta,\rho_{b}}
\leq \hat  b(d,\delta,q,p,\beta,\rho_{a} ).  
\end{equation}

 We   set $\sigma=\sqrt a$ and consider
the equation 
\begin{equation}
                            \label{3.24,1}
x _{s}=x  +\int_{0}^{s}\sigma (t+r,x_{r})\,dw_{r}
+\int_{0}^{s}b (t+r,x_{r}) \,dr.
\end{equation}

This and the next section are based on \cite{13}.
\begin{definition}
                    \label{definition 12.15,1}
Let $q',p',\beta'\in(1,\infty)$ and let $x_{\cdot}$  
 be a solution  of
 \eqref{3.24,1}. We call it {\em
$E_{q',p',\beta'}$-admissible \/} 
\index{$E_{q,p,\beta}$-admissible solutions}%
if 
for any $ T\in(0,\infty)$ there exists
a constant $N\in(0,\infty)$ such that for any
  nonnegative Borel  $f$ on $\bR^{d+1} $  
    we have
\begin{equation}
                       \label{12.15,2}
E   \int_{0}^{T}  
f(s,x _{ s} )\,ds \leq N  
 \|  f\| _{E_{q',p',\beta'}} .
\end{equation}

\end{definition}

The following is very important. 

\begin{remark}
                            \label{remark 4.2.1}
Consider equation \eqref{3.24,1} with zero
initial data and $t=0$ and make the change of variables
$ x_{t} = \rho _{b}y_{\rho_{b}^{-2}t}  $, 
$B_{t}=\rho_{b} w_{\rho_{b}^{-2}t}$. Then
\begin{equation}
                                    \label{4.4.1}
dy_{t} = \tilde  b(t,y_{t} )\,dt+\tilde 
\sigma( t,y_{t} )\,dB_{t},
\end{equation}
where $\tilde b(t,x)=\rho _{b}b(\rho_{b}^{ 2}t,\rho_{b}x)$, $\tilde\sigma(t,x)
=\sigma(\rho_{b}^{ 2}t,\rho_{b}x)$, and $B_{t}$
is a Wiener process.

Taking into account that $\rho_{b}\leq1$, it is easy to check that $\tilde \sigma$
and $\tilde b$ satisfy the assumptions of Theorem \ref{theorem 5.8,20} with the {\em same\/}
$\rho_{a},\hat a,\hat b $ and 1
in place of $\rho_{b}$. 
At the same time the issues of existence and uniqueness of solutions of \eqref{4.4.1} and
\eqref{3.24,1} are equivalent.  
Also note that $E_{q,p,\beta}$-admissible 
solutions are still $E_{q,p,\beta}$-admissible
after this transformation.

\end{remark}

{\em This remark shows that without loosing generality
in the rest of the chapter
we impose\/}

\begin{assumption}
                           \label{assumption 4.3.1}
We have $\rho_{b}=1$.
\end{assumption}

\begin{lemma}
                 \label{lemma 3.16.1}
 Let $x_{\cdot}$   be a solution  of
 \eqref{3.24,1} with $(t,x)=(0,0)$
and let $\tau$ be a stopping time
such that $(t,x_{t})\in C_{R}$ for $t\leq \tau$ and some
$R\in(0,\infty)$. Assume that

(a) for any  
  Borel nonnegative $f$ on $\bR^{d+1}$
\begin{equation}
                             \label{12.11.06}
E\int_{0}^{\tau }f(s,x_{s})\,ds\leq N\|f\|_{E_{q,p,\beta }},
\end{equation}
where $N$ is independent of $f$. Then

(b) (It\^o's formula)  for
\index{$S$@Miscelenea!It\^o's formula}%
 any    $u\in E^{1,2}_{q,p,\beta  } $, 
 with probability one for all $t\geq0$,
\begin{equation}
                                \label{3.16.2}
u(t\wedge\tau ,x_{t\wedge\tau })
=u( 0)+\int_{0}^{t\wedge\tau }D_{i}u 
\sigma^{ik} (s,x_{s})\,dw^{k}_{s}
+\int_{0}^{t\wedge\tau }\cL u (s,x_{s})\,ds
\end{equation}
and the stochastic integral above is a 
square-integrable
martingale.
\end{lemma}

Proof. By Corollary \ref{corollary 10.8.2} we have
$|Du|^{2}\in E_{s/2,r/2,2(\beta  -1) }$, where $r=
p\beta  /(\beta  -1)$,  $s=
q\beta  /(\beta  -1)$. Note that  
$$
  2>\beta >1,\quad \beta /(\beta -1) >2,
\quad 2(\beta -1)  < \beta ,
\quad
r /2\geq p, \quad s /2\geq q.
$$
This implies that the last statement of the lemma 
follows from  \eqref{12.11.06}.

Then we apply It\^o's formula to $u^{(\varepsilon)}$.
Since $u$ is bounded and continuous,  we have
the convergence  of the terms 
$u^{(\varepsilon)}(t\wedge\tau ,x_{t\wedge\tau })
,u^{(\varepsilon)}( 0)$ to
$u (t\wedge\tau ,x_{t\wedge\tau })
,u ( 0)$.

The inequality $ \beta  -1   < \beta/2$
and Lemma \ref{lemma 2.16.2}   imply that
$$
E\Big|\int_{0}^{t\wedge\tau }D_{i}\big(u^{(\varepsilon)}-u\big) 
\sigma^{ik} (s,x_{s})\,dw^{k}_{s}\Big|^{2}
$$
$$
\leq N\int_{0}^{t\wedge\tau }|D  u^{(\varepsilon)}-Du|^{2}(s,x_{s})\,ds
\leq N\|D  u^{(\varepsilon)}-Du\|^{2}_{E_{2q,2p, \beta/2 }(C_{R })}\to0
$$
as $\varepsilon\downarrow 0$.
This shows that after we apply It\^o's
formula to $u^{(\varepsilon)}$, we will be able to pass to the limit in the stochastic integral term.

In what concerns the usual integral,
observe that   estimate \eqref{12.11.06} implies 
the existence of  a Borel
function $g(t,x)\geq0$ such that 
$$
\int_{C_{R}}gf\,dxdt=E\int_{0}^{ \tau }f(s,x_{s})\,ds\leq N\|f\|_{E_{q,p ,\beta  }} 
$$
for any $f\geq0$.
Then Lemma \ref{lemma 2.16.3} shows that
$$
\int_{C_{R}}g|a^{ij}D_{ij}(u^{(\varepsilon)}-u)| \,dxdt\leq
N\int_{C_{R}}g| (D^{2} u)^{(\varepsilon)}-D^{2}u| \,dxdt\to 0
$$
as $\varepsilon\downarrow 0$.
This means that we can pass to the limit in the usual integral containing
$a^{ij}D_{ij}u^{(\varepsilon)}$.
Furthermore,  there is a Borel $h\geq0$ such that
$$
\int_{C_{R}}h f\,dxdt=E\int_{0}^{ \tau }|b|f(s,x_{s})\,ds\leq N\|bf\|_{E_{q,p,\beta }} 
\leq N\|f\|_{E_{s,r,\beta-1}}
$$
for any $f\geq0$, where
the last inequality follows from
Lemma  \ref{lemma 6,11.1}.
Then Lemma \ref{lemma 2.16.3} implies that
$$
E\int_{0}^{ \tau }
|b|  \,|D u^{(\varepsilon )}-Du |(s,x_{s})\,ds
=\int_{C_{R}}h |D u^{(\varepsilon )}-Du |\,dxdt\to 0
$$
as $\varepsilon\downarrow 0$.
   The theorem is proved.  \qed

\begin{definition}
                   \label{definition 12.16,1}

Let $x_{\cdot}$   be a solution  of
 \eqref{3.24,1} with $(t,x)=(0,0)$. We call it 
 $(q,p,\beta)$-{\em reasonable\/}
\index{$D$@Processes!$(q,p,\beta)$-reasonable solution}%
 if
there exists a sequence of stopping
times $\tau^{n}\uparrow\infty$ 
such that each $\tau^{n}$ satisfies the condition of Lemma \ref{lemma 3.16.1}
(with $R $ depending on $n$)
and for each $\tau=\tau^{n}$
either
  (a) (with $ N$ depending on $n$) or (b) 
  of Lemma \ref{lemma 3.16.1} holds. 
\end{definition}

In an obvious way one defines 
$(q,p,\beta)$-reasonable solutions of 
\eqref{3.24,1} with arbitrary starting point $(t,x)$.
Notice that $E_{q,p,\beta}$-admissible solutions
are $(q,p,\beta)$-reasonable. In the proof of
Theorem \ref{theorem 12.16,6} we will see an advantage of using
the notion of $(q,p,\beta)$-reasonable solutions.

\begin{theorem}[Unconditional and conditional  weak uniqueness]
           \label{theorem 12.12.3}
 (i)~Suppose that there exist $q',p'$ such that 
 $(d,q',p')$ are properly tight
  and $b\in L_{(q',p'),\loc}$. Also suppose that
\begin{equation}
                        \label{11.14,8}
\frac{d}{p}+\frac{1}{q}\leq 1
\end{equation}
and, in case $p\geq q$,  in the whole chapter 
$L_{q,p}$ is
  defined by using \eqref{3.27.3}, however,
in case $p\leq q$, in the whole chapter 
$L_{q,p}$ is
 defined by using \eqref{4.3.2}. 
Then   all solutions of \eqref{3.24,1} with fixed $(t,x)$ (provided they exist)  
are $(q,p,\beta)$-reasonable  and
have the same finite-dimensional distributions.

(ii) Generally, let 
$$
 (a,b,q_{i},p_{i},\beta_{i}) \in \sfA\cap\{\beta<2\}\quad i=1,2 ,
$$
and let
$y^{(i)}_{\cdot}$   be   $(q_{i},p_{i},\beta_{i})$-reasonable solutions of
 \eqref{11.29.20} with $(t,x)=(0,0)$ perhaps on different probability spaces. 
Then $y^{(1)}_{\cdot}$ and $y^{(2)}_{\cdot}$ have the same finite-dimensional distributions.

\end{theorem}

Proof. First, we prove (ii). Since, by Lemma
\ref{lemma 3.16.1}, (a) implies (b), we need  only show that the fulfillment of (b) for each $\tau=\tau_{n}$ implies weak uniqueness.
Denote by $x_{\cdot}$ one of $y^{(i)}_{\cdot}$.

Take     $ f_{0},...,f_{m}\in C^{\infty}_{0}(\bR^{d}) $ and $0=t_{0}< t_{1}<...<t_{m}<\infty$.
By Theorem \ref{theorem 5.8,20}
and Remark \ref{remark 3.16,1} for each $\lambda$ large enough there is
a bounded function $u$,
such that $u\in E^{1,2}_{q_{1}, p_{1},\beta_{1}}
\cap E^{1,2}_{q_{2}, p_{2} ,\beta_{2} }$ and 
\eqref{10.14,1} holds with $f=f_{m}$. 
In light of (b)
by   It\^o's formula  
applied to $
u(t,x_{t})e^{-\lambda t}$
we obtain
$$
Ef_{0}(x_{t_{0}\wedge \tau^{n}}),...,
f_{m-1}(x_{t_{m-1}\wedge \tau^{n}})
u (t_{m-1}\wedge \tau^{n},x_{t_{m-1}\wedge \tau^{n}})e^{-\lambda
(t_{m-1}\wedge \tau^{n})}
$$
$$
=Ef_{0}(x_{t_{0}\wedge \tau^{n}}),...,
f_{m-1}(x_{t_{m-1}\wedge \tau^{n}})
\int_{t_{m-1}\wedge \tau^{n}}
^{  \tau^{n}}e^{-\lambda t}
f_{m}(t,x_{t})\,dt
$$
$$
+Ef_{0}(x_{t_{0}\wedge \tau^{n}}),...,
f_{m-1}(x_{t_{m-1}\wedge \tau^{n}})
u (\tau^{n},x_{ \tau^{n}})e^{-\lambda
\tau^{n}}.
$$
By letting $n\to\infty$ we get
$$
Ef_{0}(x_{t_{0} }),...,
f_{m-1}(x_{t_{m-1} })
u (t_{m-1},x_{t_{m-1} })e^{-\lambda
 t_{m-1} }
$$
\begin{equation}
                        \label{11.14,3}
=\int_{t_{m-1} }
^{\infty}e^{-\lambda t}Ef_{0}(x_{t_{0} }),...,
f_{m-1}(x_{t_{m-1} })
f_{m}(t,x_{t})\,dt.
\end{equation}
On the right we have the Laplace transform
of a function, knowing which uniquely
defines this function up to almost everywhere, but because the function is continuous in $t$, it defines it uniquely
for all $t\geq t_{m-1}$. Therefore, if
we suppose that the distribution of
$(x_{t_{0} } ,...,
 x_{t_{m-1} })$ is uniquely defined
(independent of $i=1,2$), then
the left-hand side of \eqref{11.14,3} 
is uniquely defined implying that
the distribution of
$(x_{t_{0} } ,...,
 x_{t_{m} })$ is uniquely defined. For $m=1$ \eqref{11.14,3} implies that
the distribution of $x_{t_{1}}$
is uniquely defined and then the induction
on $m$ proves assertion (ii).

Assertion (i) follows from assertion (ii)
and Theorem \ref{theorem 2.3,1},  which guarantees that
any solution of \eqref{3.24,1} admits estimate
\eqref{12.11.06} with $L_{(q,p)}(C_{R})$ in place of 
$E_{q,p,\beta}$ 
and $E_{q,p,\beta}\subset L_{(q,p)}(C_{R})$, if the mixed-norms
are understood as in    (i).
The theorem is proved. \qed 
 
\begin{remark}
                     \label{remark 11.14,3}
It is worth finding out which $(q,p,\beta)$
satisfy  \eqref{4.11,1} and \eqref{11.14,8}. It is not hard to check that
both conditions are satisfied iff
$$
\infty>p>d,\quad 1<\beta\leq 2-\frac{d}{p},
\quad 2-\frac{2d}{p}\geq\frac{2}{q}\geq \beta-\frac{d}{p},
$$
or iff
\begin{equation}
                                       \label{4.11,2}
\infty>p>d, 
\quad \frac{d}{p}+\frac{2}{q}\geq \beta>1,\quad
\frac{d}{p}+\frac{1}{q}\leq1.
\end{equation}
For instance, if $p=d+1$ and $q=d+2$
the above inequalities are satisfied and
$q>p$, so that to apply assertion (i)
of Theorem \ref{theorem 12.12.3} we should
use the norm defined by \eqref{4.3.2}. 

\end{remark}

\begin{remark}
As is shown in Example \ref{example 3.22.2}, assuming $b\in L_{q,p}$ with $d/p +1/q \leq 1$
alone does not guarantee weak uniqueness even with unit diffusion.   
\end{remark}

An interesting situation in Theorem \ref{theorem 12.12.3} (ii) occurs when one of $x_{\cdot}$
or $y_{\cdot}$ is a strong solution. It turns out that then the other one is also strong
and  there could be only one strong solution on a given probability space for which either (a)   or (b) 
  of Lemma \ref{lemma 3.16.1}
holds. It is worth proving this for equations    
more general than \eqref{11.29.20}.
Consider the equation
\begin{equation}
                         \label{2.28.2}
x_{t}=\int_{0}^{t}\sigma(s,x_{s})\,dw_{s}
+\int_{0}^{t}b(s,x_{s})\,ds,
\end{equation}
where   $\sigma =(\sigma^{ik})$ is Borel with values in the set of $d\times d_{1}$-matrices ($d_{1}\geq d$), and $w_{t}
=(w^{1}_{t},...,w^{d_{1}}_{t})$ is a Wiener process
on a probability space.   
We suppose that  $a =(a^{ij}):=\sigma\sigma^{*}$ is $\bS_{\delta}$-valued  and $b$ satisfy the assumptions
stated at the beginning of the section,
that is we have $q,p,\beta$ such that
$(a,b,q,p,\beta)\in\sfA\cap\{\beta<2\}$. Since
\begin{equation}
                             \label{3.16,2}
x_{t}=\int_{0}^{t}\sqrt{a}(s,x_{s})\,dB_{s}
+\int_{0}^{t}b(s,x_{s})\,ds,
\end{equation}
where $B_{s}$, defined by $dB_{s}=a^{-1/2}\sigma(s,x_{s})\,dw_{s}$ with $B_{0}=0$, is a $d$-dimensional Wiener process (compute its bracket),
the above results are applicable to equation
\eqref{2.28.2}. Let us call $x_{t}$ $(q,p,\beta)$-reasonable 
solution of \eqref{2.28.2} if it is 
$(q,p,\beta)$-reasonable solution
\index{$D$@Processes!$(q,p,\beta)$-reasonable solution}%
 of \eqref{3.16,2}.
\begin{theorem}
                         \label{theorem 2.28.1}

Suppose that 
$$
 (a,b,q_{i},p_{i},\beta_{i}) \in \sfA\cap\{\beta<2\}\quad i=1,2 ,
$$
and let
$y^{(i)}_{\cdot}$   be   $(q_{i},p_{i},\beta_{i})$-reasonable solutions of
 \eqref{2.28.2} with the same Wiener process.
Assume that $y^{(1)}_{\cdot}$ is a strong solution. Then  $y^{(1)}_{\cdot}=y^{(2)}_{\cdot}$ (a.s).

\end{theorem}

Proof. First notice that by Theorem \ref{theorem 12.12.3} the processes $y^{(1)}_{\cdot}$
and $y^{(2)}_{\cdot}$ have the same finite-dimensional distributions and, hence, both
are $(p_{1},q_{1},\beta_{1})$-reasonable.
After that
one is tempted to refer
to the result of A. Cherny \cite{Ch_02} saying that weak uniqueness and strong existence imply
the uniqueness of solutions. However,
in his result one needs unconditional
weak uniqueness which we do not know how to prove
in the general case. Therefore, we proceed
differently still using the idea from \cite{Ch_02}.

Define     $\tau=\sigma^{*}\sigma$. This is a symmetric  nonnegative definite matrix and the following
is well defined
$$
\Sigma=\lim_{\varepsilon\downarrow 0}
\tau(\tau+\varepsilon I)^{-1},
$$
where $I$ is the $d_{1}\times d_{1}$ identity matrix. As is easy to see by using the diagonal forms,
$\Sigma^{2}=\Sigma$, $\Sigma\tau=\tau$, and ($\text{tr}\,AB=
\text{tr}\,BA$)
$$
\text{tr}\,( \Sigma \sigma^{*}-\sigma^{*})
(\sigma \Sigma-\sigma)=
\text{tr}\,(\Sigma \sigma^{*}\sigma\Sigma
-\sigma^{*}\sigma\Sigma-\Sigma\sigma^{*}\sigma
+\sigma^{*}\sigma)
=\text{tr}\,(-\Sigma\tau+\tau)=0,
$$
so that $\sigma\Sigma=\sigma$.

Then let $(q,p,\beta)=(q_{1},p_{1},\beta_{1})$
and let $x_{\cdot}$ be any $(q,p,\beta)$-reasonable solution of
\eqref{2.28.2},
so that its distribution coincide with that of
$y^{(1)}_{\cdot}$.

 By extending our probability
space, if necessary, we suppose that we are also
given a $d_{1}$-dimensional Wiener process
$\bar w_{s}$ independent of $w_{t}$. Define
$$
\xi_{s}=\int_{0}^{s}\Sigma(u,x_{u})\,d\bar w_{u}+\int_{0}^{s}\big(I-\Sigma(u,x_{u})\big)\,
dw_{u}.
$$
An easy application of the L\'evy theorem shows that $\xi_{s}$ is a $d_{1}$-dimensional Wiener process.

The crucial step is to prove that
the processes $x_{\cdot}$ and $\xi_{\cdot}$
are independent   because
(dropping arguments $(s,x_{s})$) 
$$
dx^{i}_{s}d\xi^{k}_{s}=\sigma^{ir}\,dw^{r}_{s}(
\delta^{kn}-\Sigma^{kn})\,dw^{n}_{s}
=\sigma^{ir} (
\delta^{kr}-\Sigma^{kr})\,ds
=(\sigma^{ik} 
 -\sigma^{ir}\Sigma^{rk})\,ds=0.
$$

To do that, fix $T\in(0,\infty)$, take two bounded Borel
functions $c'$ and $c''$ with compact support on $\bR^{d+1}$ and $\bR^{d_{1}+1}$. Then by Remark \ref{remark 12.2,1} for sufficiently large $\lambda>0$
there exist   $v',v''\in E^{1,2}_{q,p,\beta}$
such that for $t\geq T$ we have $v'(t,\cdot)=1$,
$v''(t,\cdot)=1$ and for $t<T$
$$
\partial_{t}v'+(1/2)a^{ij}D_{ij}v'+(c'-\lambda)v'=0,\quad \partial_{t}v''+(1/2)\Delta v''+(c''-\lambda)v''=0.
$$

Set
$$
\phi_{t}= \int_{0}^{t}[c'(s,x_{s})+c''(s,\xi_{s})]\,ds-2\lambda t .
$$
 
By It\^o's formula applied to 
$$
v' (t,x_{t})v''(t,\xi_{t})e^{\phi_{t}},
$$
we get
$$
  e^{\phi_{ T}}=
v' (0,0)v''(0,0)
$$
$$
+\int_{0}^{T}e^{\phi_{s}}\big[v''(s,\xi_{s}) \sigma^{ik}D_{i}v' (s,x_{s})\,
dw^{k}_{s}+v'_{n}(u,x_{u})D_{\xi^{i}}
v''(s,\xi_{s})\,d\xi^{i}_{s}\big].
$$
By taking expectations
we see that  
\begin{equation}
                   \label{7.7.2}
E\exp\Big(\int_{0}^{T}c'(s,x_{s})\,ds\Big)\exp\Big(\int_{0}^{T}c''(s,\xi_{s})\,ds\Big)=v' (0,0)v''(0,0)e^{2T}.
\end{equation}
The arbitrariness of $c',c''$  implies first that
$$
E\exp\Big(\int_{0}^{T}c'(s,x_{s})\,ds\Big)\exp\Big(\int_{0}^{T}c''(s,\xi_{s})\,ds\Big)
$$
$$
=
E\exp\Big(\int_{0}^{T}c'(s,x_{s})\,ds\Big)E\exp\Big(\int_{0}^{T}c''(s,\xi_{s})\,ds\Big),
$$
and second that $x_{\cdot}$ and $\xi_{\cdot}$
are independent indeed.

Then we observe that
$$
w_{s}=I_{s} 
+\int_{0}^{s}\big(I-\Sigma(u,x_{u})\big)\,d\xi_{u},
$$
where
$$
I_{s}=\int_{0}^{s}\Sigma(u,x_{u})\,dw_{u}=
\lim_{\varepsilon\downarrow 0}
\int_{0}^{s}(\tau(u,x_{u})+\varepsilon I)^{-1}
\sigma^{*}(u,x_{u})\,dm_{u},
$$
$$
m_{s}=\int_{0}^{s}\sigma (u,x_{u})\,dw_{u}
=x_{s}-\int_{0}^{s}b(u,x_{u})\,du.
$$
We see that $I_{s}$ is a functional of $x_{\cdot}$, so that the distribution
of $I_{\cdot}$   is independent of which
$(q,p,\beta)$-reasonable solution we take. Since
the Wiener process
$\xi_{\cdot}$ is independent of $x_{\cdot}$,
the conditional distribution of $w_{\cdot}$
given $x_{\cdot}$ and
the joint distribution of $(w_{\cdot},x_{\cdot})$ are   independent of which
$(q,p,\beta)$-reasonable solution we take.

Now assume  that   $x_{\cdot}$ is a strong
solution: $x_{\cdot}=x_{\cdot}(w_{\cdot})$. Then
the joint distribution of $(w_{\cdot},x_{\cdot})$
is concentrated on the set $\Gamma:=\{(w_{\cdot},x_{\cdot}
(w_{\cdot}))\}$, and the joint distribution
of $(w_{\cdot},y_{\cdot})$ for any other
$(q,p,\beta)$-reasonable solution $y_{\cdot}$  is also concentrated on this
set. Since for any $w_{\cdot}$ there is only
one point $(w_{\cdot},x_{\cdot}(w_{\cdot}))$
in $\Gamma$, $y_{\cdot}$ should be equal to
$x_{\cdot}(w_{\cdot})$. This proves the theorem.
\qed
 
\mysection[Existence of (weak) solutions]
{Existence of (weak) solutions} 

The setting in this section is the same as 
in Section \ref{section 3.22,1}, that is,
Assumption \ref{assumption 4.3.1} that $\rho_{b}=1$ is supposed to be 
satisfied and
$$
p ,q  \in(1,\infty),\quad 1<\beta\leq\frac{d}{p }+\frac{2}{q } ,
\quad \beta<2,
$$
\begin{equation}
                                            \label{4.25.1}
a^{\shharp}_{\rho_{a}}\leq 
\hat  a(d,\delta,q,p,\beta),\quad\hat b_{q\beta,p\beta,\rho_{b}}
\leq \hat  b(d,\delta,q,p,\beta,\rho_{a} ).  
\end{equation}

We start by drawing
consequences from Theorem
\ref{theorem 5.8,20}.  

\begin{corollary}
                            \label{corollary 3.14.1}
Assume that $a,b$ are smooth and bounded. Take
$R\leq 1$, smooth $f$, and let $u$ be the classical solution of
\begin{equation} 
                            \label{3.14.5}
\cL u+f=0
\end{equation}
in $C_{R}$ with zero boundary condition
on $\partial' C_{R}$.   Then
\begin{equation} 
                            \label{3.14.6}
|u|\leq NR^{2-\beta }\|I_{C_{R}}f\|_{E_{q,p,\beta }},
\end{equation}
where $N$ depends only on $d,\delta,q ,p ,\beta ,\rho_{a} $.
\end{corollary}

Indeed, the case $R<1$ is reduced to $R=1$ by using
parabolic dilations. If $R=1$, 
the maximum principle allows us to concentrate on $f\geq0$ and also shows that
$u(t,x)e^{\lambda _{0}t}$  ($\lambda_{0}$ is from Theorem
\ref{theorem 5.8,20}) is smaller  in $C_{1}$
than the solution $v$ of
$$
\cL v-\lambda _{0} v+I_{C_{1}}fe^{\lambda _{0}t}=0    
$$
in $\bR^{d+1}$. Since $\beta <2$
by embedding theorems we have on $C_{1}$
$$
u\leq v 
\leq N\|v\|_{E^{1,2}_{q,p,\beta }}
\leq N\|I_{C_{1}}f\|_{E_{q,p,\beta }}.  
$$

\begin{corollary}
                            \label{corollary 3.14.2}
Assume that $a,b$ are smooth and bounded and let
$(\sft_{s},x_{s})$ be the corresponding Markov
diffusion process. Then for any $(t,x)\in\bR^{d+1}$,
 $\rho\leq1$, $C\in\bC_{\rho}$,  and Borel $f\geq0$
\begin{equation} 
                            \label{3.14.7}
I(t,x):=E_{t,x}\int_{0}^{\tau_{C}}f(t,x_{t})\,dt\leq N\rho^{2-\beta }\|I_{C }f\|_{E_{q,p,\beta }},
\end{equation}
where $\tau_{C}$ is the first exit time of 
$(\sft_{s},x_{s})$ from $C$. In particular,
for any $(t,x)\in\bR^{d+1}$,
 $\rho\leq1$, $C\in\bC_{\rho}$,  
\begin{equation} 
                            \label{3.14.8}  
 E_{t,x}\int_{0}^{\tau_{C}}|b(t,x_{t})|\,dt\leq  N_{1} 
 \rho\hat b_{p\beta,q\beta,1} ,
\end{equation}
and in both estimates  $N$ and $N_{1}$ depend  only on $d,\delta,q ,p ,\beta ,\rho_{a} $.
\end{corollary}

Indeed, if $f$ is smooth, by It\^o's formula,
$I$ coincides with the solution of \eqref{3.14.5}
in a shifted $C_{\rho}$ and \eqref{3.14.7}
follows from \eqref{3.14.6}. For bounded Borel
$f$ we use the notation $f^{(\varepsilon)}$ from
the proof of Theorem 
\ref{theorem 10.21.1} and observe that $f^{(\varepsilon)}\to f$
almost everywhere, and the corresponding left-hand sides of \eqref{3.14.7} converge because they
are expressed in terms of the Green's function
of $\cL$. As far as the right-hand sides are concerned,
observe that by Minkowski's inequality
$\|f^{(\varepsilon)}\|_{E_{q,p ,\beta }}\leq \|f
 \|_{E_{q,p ,\beta }}$
and this yields \eqref{3.14.7} with
 $\|f \|_{E_{q,p ,\beta }}$
in place of $\|fI_{C} \|_{E_{q,p ,\beta }}$. Plugging 
$fI_{C}$ in such relation  in place of $f$ leads to \eqref{3.14.7}
as is. The passage to arbitrary $f\geq0$ is achieved
by taking $f\wedge n$ and letting $n\to\infty$.

To prove \eqref{3.14.8} while estimating
$$
\rho^{2-\beta }\|I_{C }b\|_{E_{q,p,\beta }}=N
 \sup_{r\leq 1}\sup_{C'\in\bC_{r}}\rho^{2-\beta }r^{\beta-d/p-2/q}
\|I_{C\cap C'}b\|_{L_{q,p}}
$$
consider two possibilities 1) $\rho_{b}=1\geq \rho
\geq r$ and 2) $1\geq r>\rho$. 

In case 1)
we have $r^{\beta}\leq \rho^{\beta-1}r$ and
$$
I:=\rho^{2-\beta }r^{\beta-d/p-2/q}
\|I_{C\cap C'}b\|_{L_{q,p}}\leq N\rho r\dashnorm
b\|_{L_{q,p}(C')}\leq N \rho \hat b_{p\beta,q\beta,1}.
$$

In case 2) ($\beta\leq d/p+2/q$)
$$
I\leq \rho^{2-d/p-2/q}\|I_{C }b\|_{L_{q,p}}=
N\rho(\rho \dashnorm
b\|_{L_{q,p}(C)})\leq N  \rho\hat b_{p\beta,q\beta,1}.
$$

Once $N_{1}$ is specified, we have the following.  

\begin{corollary}
                    \label{corollary 3.14.6}
Suppose that $a,b$ are smooth and bounded and let
$(\sft_{s},x_{s})$ be the corresponding Markov
diffusion process.  Suppose that
\begin{equation}
                      \label{11.15,4}
N_{1} \hat b_{q\beta,p\beta,1}< \sfb_{0}
\end{equation}

Then 
$ \bar b_{1 }\leq \sfb_{0}$ and all results
from Chapter \ref{chapter 10.20.1}
after Theorem \ref{theorem 8.2.1} are applicable
(with $\rho_{b}=1$). In particular, by Corollary
\ref{corollary 10.26.1}
for any $n>0$, $r,s\geq0$, and
$(t,x)$  we have
\begin{equation}
                  \label{10.28.20}
E_{t,x}\sup_{\tau\in[0,r]}|x_{s+\tau}-x_{s}|^{ n}
\leq N(  r ^{ n/2}+ r ^{ n}),
\end{equation}
where $N=N(n,  \sfp_{0}(d,\delta))$.
\end{corollary}

\begin{remark}
                   \label{remark 11.15,9}
One could have obtained \eqref{10.28.20}
on the basis of Theorem \ref{theorem 9.27.10} (with no regularity assumption on $a$). However,
 in condition \eqref{11.15,4}
the definition of $\hat b_{q\beta,p\beta,1}$ should be then modified according to the cases $p\geq q$ and $q\geq p$, and not fixed as in the whole chapter regardless
of these cases.     

\end{remark}

\begin{corollary}
                   \label{corollary 3.14.7}
Suppose that $a,b$ are smooth and bounded and let
$(\sft_{s},x_{s})$ be the corresponding Markov
diffusion process. Then for any $(t,x)\in\bR^{d+1}$,
Borel $f\geq0$, $T\in(0,\infty)$, there exists
$N$ depending only on $d,\delta$, $q $, $p $, $\beta $, $\rho_{a}$,   $T$, such that
\begin{equation}
                            \label{3.14.9}
 E_{t,x}\int_{0}^{T}f(t,x_{t}) \,dt\leq N
\|f\|_{E_{q,p,\beta }}.
\end{equation}
\end{corollary}

The proof of this is almost identical to the proof
of  \eqref{3.14.7} when $\rho=1$.

  Now we abandon the assumption
that   $a$  and $b$ are   smooth and come back
to our   assumptions
(that are supposed to hold throughout the rest of the
 section)
 stated in the beginning of  Section \ref{section 3.22,1}
(including Assumption \ref{assumption 4.3.1}). We also suppose that \eqref{11.15,4} holds
until the end of this section.
 Here is a counterpart of Theorem \ref{theorem 10.21.1} (that we cannot use,
for instance,
because of the discrepancy between $L_{(q,p)}$ and $L_{q,p}$)
\begin{theorem}
               \label{theorem 3.15.1} 
Under the assumptions stated before the theorem:

(i) There is a probability space 
$(\Omega ,\cF ,P )$,
a filtration of $\sigma$-fields $\cF _{s}\subset \cF $, $s\geq0$,
a process $w _{s}$, $s\geq0$, which is a $d$-dimensional Wiener process
relative to $\{\cF _{s}\}$, and an $\cF _{s}$-adapted
process $x_{s}$ such that 
 (a.s.) for all   $s\geq0$ equation \eqref{3.24,1} holds with $(t,x)=(0,0)$. 

(ii) Furthermore, for any 
  nonnegative Borel  $f$ on $\bR^{d+1} $ and
 $ T\in(0,\infty)$   we have
\begin{equation}
                                    \label{3.15.3}
E   \int_{0}^{T}  
f(s,x _{s} )\,ds \leq N  
 \|  f\| _{E_{q,p,\beta  }} ,
\end{equation}
where $N$ is the constant from \eqref{3.14.9},
so that $x_{\cdot}$ is an $E_{q,p,\beta}$-admissible solution.

\end{theorem} 

Proof. As in the proof of Theorem 
\ref{theorem 10.21.1},
approximate $\sigma,b$ by smooth
$\sigma^{(\varepsilon)},b^{(\varepsilon)}$
and take the corresponding Markov processes   
$(\sft_{t},x^{\varepsilon}_{t})$.
By Corollary \ref{corollary 3.14.6}  the  $P_{0,0}$-distributions of $x^{\varepsilon}_{\cdot}$
are precompact on the space $C([0,\infty),\bR^{d})$
and a subsequence   $\varepsilon=\varepsilon_{n}
\downarrow 0$ of them converges
to the distribution of a process $x_{\cdot}=x^{0}_{\cdot}$ defined
on a probability space (the coordinate process on $\Omega=C([0,\infty),\bR^{d})$ 
with cylindrical $\sigma$-field $\cF$ 
completed with respect to
 $P$, which is
 the limiting distribution of $x^{\varepsilon}_{\cdot}$). Furthermore, by Corollary \ref{corollary 3.14.7}     for any 
  nonnegative Borel  $f$ on $\bR^{d+1} $ and
 $\varepsilon,T\in(0,\infty)$   we have
\begin{equation}
                                    \label{3.15.2} 
E _{0,0} \int_{0}^{T}  
f(s,x^{\varepsilon}_{s} )\,ds \leq N 
\|f\|_{E_{q,p,\beta }} \quad (\leq N'
\|f\|_{L_{r}(\bR^{d+1})},r\gg 1)
\end{equation}
where $N$ is the constant from \eqref{3.14.9},
which by continuity is extended to $\varepsilon=0$
for  continuos $f$ with compact support and then by standard arguments for all Borel $f\geq0$.
This proves (ii).
 
Observe that estimate \eqref{3.15.2} also shows
that for any bounded Borel $f$ with compact support
\begin{equation}
                                    \label{3.15.5} 
\lim_{\varepsilon\downarrow 0}E _{0,0} \int_{0}^{T}  
f(s,x^{\varepsilon}_{s} )\,ds =
E _{0,0} \int_{0}^{T}  
f(s,x^{0}_{s} )\,ds .
\end{equation}

 Now  we prove that assertions (i)  holds for
$x_{\cdot}$.
Estimate \eqref{10.28.20} implies that for any finite~$T$
$$
\lim_{c\to\infty}P(\sup_{s\leq T}|x^{0}_{s}|>c)=0,
$$
and estimate   \eqref{3.15.3}
shows that for any finite $c$
$$
E   \int_{0}^{T}  I_{|x^{0}_{s}|\leq c}
|b (s,x^{0}_{s} )|\,dt <\infty.
$$
Hence, with probability one
$$
\int_{0}^{T}  
|b(s,x^{0}_{s} )|\,dt<\infty.
$$

Next, for $0\leq t_{1}\leq...\leq t_{n}\leq t\leq s$, bounded continuos function $\phi(x(1),...,x(n))$,
and smooth bounded $u(t,x)$ with compact support by It\^o's formula we have 
$$
E_{0,0}\phi(x^{\varepsilon}_{t_{1}},...,x^{\varepsilon}_{t_{n}})
\Big[u(s,x^{\varepsilon}_{s})-u(t,x^{\varepsilon}_{t})-
\int_{t}^{s}\cL^{\varepsilon}u(r,x^{\varepsilon}_{r})\,dr\Big]=0,
$$
where
$$
\cL^{\varepsilon}u=\partial_{t}u+(1/2)a^{\varepsilon ij}D_{ij}u+
b^{\varepsilon i}D_{i}u,\quad a^{\varepsilon}= (\sigma^{(\varepsilon)})^{2}.
$$

Using \eqref{3.15.3}, \eqref{3.15.2},  Lemma \ref{lemma 2.16.2} with $\beta'=1$, and the fact that $u$ has compact support show   that
$$
\lim_{\varepsilon_{1}\downarrow0}
\lim_{\varepsilon \downarrow0}E_{0,0}\int_{t}^{s}\big| b^{(\varepsilon)}-b^{(\varepsilon_{1})}\big|(r,x^{\varepsilon}_{r})|Du(r,x^{\varepsilon}_{r})|
\,dr = 0,
$$
$$
\lim_{\varepsilon \downarrow0}E_{0,0}\int_{t}^{s}  b^{(\varepsilon_{1})i} (r,x^{\varepsilon}_{r}) D_{i}u(r,x^{\varepsilon}_{r})
\,dr = E \int_{t}^{s}  b^{(\varepsilon_{1})i} (r,x^{0}_{r}) D_{i}u(r,x^{0}_{r})
\,dr,
$$
$$
\lim_{\varepsilon_{1} \downarrow0}E
\int_{t}^{s}\big| b -b^{(\varepsilon_{1})}\big|(r,x^{0}_{r})|Du(r,x^{0}_{r})|
\,dr = 0.
$$
 After that we easily conclude that
$$
E \phi(x^{0}_{t_{1}},...,x^{0}_{t_{n}})
\Big[u(s,x^{0}_{s})-u(t,x^{0}_{t})-
\int_{t}^{s}\cL u(r,x^{0}_{r})\,dr\Big]=0.
$$
It follows that the process
$$
u(s,x^{0}_{s}) -
\int_{0}^{s}\cL u(r,x^{0}_{r})\,dr
$$
is a martingale with respect to the completion
 of $\sigma\{x^{0}_{t}: t\leq s\}$. Referring to Lemma 3.4.1 of \cite{Kr_25} proves assertion (i).
The theorem is proved. \qed  

\begin{remark}
                          \label{remark 3.7,2}

G. Zhao (\cite{Zh_20_1}) gave an example showing
that, if in condition \eqref{3.14.2}   we replace $r$ with $r^{ \alpha}$, $\alpha>1$, the weak uniqueness
may fail even in the time homogeneous case
and unit diffusion.
\end{remark}  

{\em In the following two remarks the $L_{q,p}$-norm is understood as
in \eqref{3.27.3} and as in \cite{RZ_20} and if the norm
is understood according to \eqref{4.3.2} we denote 
\index{$A$@Sets of functions!$\sfL_{q,p}$}%
\index{$N$@Norms!norm in $\sfL_{q,p}$\,}%
the corresponding space
by $\sfL_{q,p}$\/}.  

\begin{remark}
                           \label{remark 10.31.1} 

In   \cite{RZ_20} the weak uniqueness is proved in the class
of solutions admitting, as they call it, Krylov type estimate
when $\sigma=(\delta^{ij})$   and we have
$p,q\in[1,\infty]$ such that 
\begin{equation}
                           \label{10.29.5}
\frac{d}{p}+\frac{2}{q}=1,\quad\Big(
\int_{\bR}\Big(\int_{\bR^{d}} |b|^{p}\,dx\Big)^{q/p}
\,dt\Big)^{1/q}<\infty 
\end{equation}
(the Ladyzhenskaya-Prodi-Serrin condition).  
\index{$S$@Miscelenea!Ladyzhenskaya-Prodi-Serrin condition}%

Actually, $p=\infty$, $q=2$ is not allowed in \cite{RZ_20}, 
this case fits in \cite{Kr_25_1}, 
\cite{KM_22}, and \cite{Ki_25} where
weak existence and conditional weak uniqueness is obtained. In case $p=d, q=\infty$ the comparison of the results in
\cite{Kr_25_1} and    \cite{RZ_20} can be found
in \cite{Kr_25_1}.  Of course, one has to say that
apart from this result \cite{RZ_20}
and \cite{Ki_25} contain much more nontrivial
information about the solutions. In the further discussion
we assume that $b(t,x)=0$ for $t\not\in(0,T)$, where
$T\in(0,\infty)$.

In case $p\in(d,\infty)$ (we use the norms as in \eqref{3.27.3})
take small $\varepsilon\in(0,1)$ and set $p'=(1-\varepsilon)p$
with the first smallness requirement that $p'>d$. Then define
$q'$ from $1/q'=1-d/p'$. If $\varepsilon$ is small $1-d/p'$
is close to $1-d/p=2/q$, so that $2q'$ will be close to $q$
and then $q'<q$. After that we see
that  there is $\beta\in(1,2)$ close to $1$ such that
\begin{equation}
                                      \label{4.25.03}
\frac{d}{p'}+\frac{1}{q'}=1,\quad
1<\beta\leq \frac{d}{p'}+\frac{2}{q'},\quad \beta p'\leq p,\beta q'<q,
\end{equation}
which, in particular, implies that
for any $C\in \bC_{r}$
\begin{equation}
                                    \label{10.29.6}
\dashnorm b\|_{L_{\beta q',\beta p'}(C)}\leq
\dashnorm b\|_{L_{  q , p }(C)}=Nr^{-1}\|b\|_{L_{  q , p }(C)}.
\end{equation}
The last norm tends to zero as $r\downarrow 0$ if
 the norms are understood as in \eqref{3.27.3} because of
 \eqref{10.29.5}. This guarantees that \eqref{4.25.1}
 is satisfied with some $\rho_{b}>0$ ($a=(\delta^{ij})$)
 and that \eqref{11.15,4} scaled back to $\rho_{b}$
 from $\rho_{b}=1$ is also satisfied for some $\rho_{b}>0$.

Therefore $E_{q',p',\beta}$-admissible  solutions  
of \eqref{3.24,1} exist according to Theorem \ref{theorem 3.15.1}
and by Theorem
\ref{theorem 12.12.3} all $E_{q',p',\beta}$-admissible  solutions 
 have the same finite-dimensional distributions.

In addition, $d/p+1/(q/2)=1$ and $b\in L_{(q/2,p ),\loc}$ if 
$p\geq q/2$ that is $p\geq d+1$.
Therefore, for $p\geq d+1$, Theorem
\ref{theorem 12.12.3} (i) is applicable
and yields unconditional weak uniqueness
of the solutions  
of \eqref{3.24,1}.

 In case $ p\leq d+1$    ($q\geq 2(d+1)\geq 2p$)
 let us use the norm \eqref{4.3.2}, assuming
 that $b\in\sfL_{q,p}$, and set $q'=(1-\varepsilon)q$.
 Then define
$p'$ from $d/p'=1-1/q'$. If $\varepsilon$ is small $1-1/q'$
is close to $1/q +d/p$, so that $d/p'$ will be close to $1/q +d/p$
and then $p'<p$. After that we have \eqref{4.25.03} again
which, in particular, implies that
for any $C\in \bC_{r}$
 \eqref{10.29.6} holds
 with $\sfL$ in place of $L$, where the last norm still
  tends to zero as $r\downarrow 0$ if $b\in \sfL_{q,p}$,
  which case was never addressed before. 
      Therefore, in this situation again solutions
of \eqref{3.24,1} exist according to Theorem \ref{theorem 3.15.1},
weakly unique among $E_{q',p',\beta}$-admissible  solutions.

In addition, $b\in L_{(q/2,p ),\loc}$. Thus, for $p\leq d+1$, Theorem
\ref{theorem 12.12.3} (i) is applicable as well
and yields unconditional weak uniqueness.
Observe that,
  because of Minkowski's
  inequality,
  for $q>p$ the condition: $b\in \sfL_{q,p}$, is stronger than
  $b\in L_{q,p}$.

Also observe that in our setting we can include much more
irregular $b$ (and irregular $\sigma$)  along
 with $b$ from \cite{RZ_20}
and have weak solvability along with conditional 
weak uniqueness. The simplest   $b$ in $d\geq4$ is 
  given by $f(x)(|x^{1}|^{2}+|x^{2}|^{2}+|x^{3}|^{2})^{-1/2}$,
   where $f$ is a bounded $\bR^{d}$-valued function, is way
    away from satisfying  the conditions
in \cite{RZ_20} and does satisfy our conditions.

We continue discussing this example in Remark \ref{remark 6.2.1}.
\end{remark}

\begin{remark}
                  \label{remark 1.28.1}   
There are examples showing that the assumptions 
 of   Theorem \ref{theorem 12.12.3}
  (i) 
concerning $b$ are satisfied with the norm in $\sfL_{q,p}$
but not $L_{q,p}$.  For instance,  
  set  $p_{0}=d+1/2,q_{0}=2d+2$, then
  $$
  \frac{d}{p_{0}}+\frac{2}{q_{0}}>1>\frac{d}{p_{0}}+\frac{1}{q_{0}}.
  $$
  It follows that there exists $\beta\in(1,2)$ such that 
  \eqref{4.11,2} holds with
  $$
  p=p_{0}/\beta,\quad q=q_{0}/\beta,
  $$
and given that \eqref{4.11,3} is satisfied with the norm
in $\sfL_{q,p}$,
according to Remark \ref{remark 11.14,3},
equation \eqref{3.24,1} has a solution and
all solutions have the same finite-dimensional
distributions.

As an example of such $b(t,x)$ let $|b|=cf$, where 
the constant $c>0$ is sufficiently small and $f$ is constructed in the following way.
Take
$$
\gamma\in\Big(\frac{d}{d+1},\frac{2d}{2d+1}\Big)
$$
and set
$$
f(t,x)= 
\frac{1}{|x|^{\gamma}(|x|+\sqrt t)^{1-\gamma}}
I_{t>0}.
$$
Observe that  for $t_{0}\geq0$  
$$
I(r,x):=\int_{t_{0}}^{t_{0}+r^{2}}  f^{q_{0}}(t,x)\,dt=|x|^{2-q_{0}}
\int_{t_{0}/x^{2}}^{(t_{0}+r^{2})/|x|^{2}}  
\frac{1}{(1+\sqrt s)^{q_{0}(1-\gamma)}}\,dt.
$$
Here $\gamma<1$ and the derivative with respect to $t_{0}$
of the integral is negative. Moreover, $q_{0}(1-\gamma)<2$,
so that 
$$
I(r,x)\leq |x|^{2-q_{0}}\int_{0}^{r^{2}/|x|^{2}}s^{-q_{0}(1-\gamma)/2}
\,ds
= Nr^{2-q_{0}(1-\gamma)}|x|^{-q_{0}\gamma}.
$$
Since $p_{0}\gamma<d$, it follows that
$$
\dashnorm f\|^{p_{0}}_{\sfL_{q_{0},p_{0}}(C_{r}(t_{0},x_{0}))}
\leq Nr^{-p_{0}(1-\gamma)}
\dashint_{B_{r}(x_{0})}|x|^{-p_{0}\gamma}\,dx.
$$
Note that, if $x_{0}=0$ the last integral equals $Nr^{-p\gamma}$.
If $2r\geq |x_{0}|$, it is dominated by
$$
N\dashint_{B_{3r } }|x|^{-p_{0}\gamma}\,dx\leq Nr^{-p_{0}\gamma},
$$
and if $2r\leq |x_{0}|$, then the integrand is less than
$r^{-p_{0}\gamma}$. Thus,
\begin{equation}
                                          \label{4.12,1}
\dashnorm f\| _{\sfL_{q_{0},p_{0}}(C_{r}(t_{0},x_{0}))}
\leq Nr^{-1}.
\end{equation}

Hence, for small enough $c$ the assumptions 
 of   Theorem \ref{theorem 12.12.3}
  (i) 
concerning $b$ are satisfied, and equation \eqref{3.24,1}
has a solution and all solutions have the same
finite-dimensional distributions.

Here $\sfL_{q,p}$ is defined with unusual order of integration,
and the results of \cite{RZ_20} are not applicable
since 
$$
\int_{B_{1}}f^{p}(t,x)\,dx=Nt^{(1/2)(d-p) }\int_{0}^{1/\sqrt t}
\frac{\rho^{d-1-p\gamma}}
{ (\rho+1)^{p(1-\gamma)}}\,d\rho,
$$
$$
\int_{0}^{1}\Big(\int_{B_{1}}f^{p}(t,x)\,dx\Big)^{q/p}\,dt
\geq N\int_{0}^{1}t^{(1/2)(d/p-1)q }\,dt=\infty,
$$
where the equality follows from $(1/2)(d/p-1)q=-1$.
Hence, $f\not\in L_{q,p,\loc}$ no matter what $q,p$ are
 satisfying the Ladyzhenskaya-Prodi-Serrin
condition. 
 The results in \cite{Ki_25}
cover this example, with the existence and uniqueness  
proved for a {\em restricted\/} class of solutions
(in which solutions are shown to indeed exist). It maybe worth pointing out a minor difference that 
in \cite{Ki_25} the condition on the {\em singular\/} part $b$ is global: it should belong to
$\dot E_{p,p,1}$, and we do not split $b$
and our condition is local.
Recall also that our $\sigma$ is not necessarily 
constant or even continuous.

  \end{remark}
By changing the origin we can apply Theorem \ref{theorem 12.12.3}
to prove the solvability of \eqref{3.24,1}
with any initial data $(t,x)$ and get  solutions with the
properties as in Theorems \ref{theorem 3.15.1} (ii)  
  weakly unique by Theorem
\ref{theorem 12.12.3}.
For such a solution denote by $P_{t,x}$ the distribution of $(\sft_{s},x_{s}),s\geq0$, ($\sft_{s}=t+s$)
on the Borel $\sigma$-field $\cN_{\infty}$ of $\Omega=C([0,\infty),\bR^{d+1})$. For $\omega=(\sft_{\cdot},x_{\cdot})\in \Omega$
set $(\sft_{s},x_{s})(\omega)=(\sft_{s},x_{s})$.
Also set $\cN_{s}=\sigma \{(\sft_{t},x_{t}),t\leq s\}$.

\begin{theorem} 
                          \label{theorem 12.6.02}
Under the assumptions stated before Theorem 
\ref{theorem 3.15.1} 
the process
$$
X=\{(\sft_{\cdot},x_{\cdot}), \cN_{t},P_{t,x}\} 
$$
is strong Markov 
 regular
diffusion process corresponding to
  $a,b$ 
  with strong Feller resolvent. Furthermore,
for any $(t,x)\in \bR^{d+1}$ and Borel $f\geq0$
\begin{equation}
                               \label{11.22,4}
E _{t,x}  \int_{0}^{T}  
f(s,x _{s} )\,ds \leq N  
 \|  f\| _{E_{q ,p,\beta  }} ,
\end{equation}
where $N$ is the constant from \eqref{3.14.9}.

\end{theorem}

Proof.  
Take $u$ from Theorem \ref{theorem 5.8,20} with 
$
\lambda\geq\lambda_{0}$ and Borel bounded $f$. 
By It\^o's formula for any $(t,x) $
and $ 0\leq r\leq s$ we obtain that with
$P_{t,x}$-probability one  
$$
u(\sft_{s} ,x_{s})e^{-\lambda(s\wedge\tau_{R})} =u(\sft_{r },x_{r })e^{-\lambda(r\wedge\tau_{R})}
+\int_{r\wedge\tau_{R}}^{s\wedge\tau_{R}}e^{-\lambda v}
 \sigma^{ik}D_{i}u(\sft_{v},x_{v})
 \,dw^{k}_{v}
$$
\begin{equation}
                                    \label{12.14.02}
-\int_{r\wedge\tau_{R}}^{s\wedge\tau_{R}}
e^{-\lambda v}f(\sft_{v},x_{v})\,dv,
\end{equation}
where $\tau_{R}$ is the first exit time of 
$(\sft_{v},x_{v})$ from $C_{R}$

From \eqref{12.14.02} with $r=0$ as in the proof of 
Theorem \ref{theorem 12.12.3} we obtain
\begin{equation}
                                    \label{12.14.3}
R_{\lambda}f(t,x):=E_{t,x}\int_{0}^{\infty}e^{-\lambda v}
f(\sft_{v},x_{v})\,dv=u(t,x).
\end{equation}
If $f$ is continuous, this implies that 
the Laplace transform of the continuous in $v$ function
$E_{t,x}f(\sft_{v},x_{v})$
is a  Borel function of $(t,x)$. Then the function
$E_{t,x}f(\sft_{v},x_{v})$ itself
is a  Borel function of $(t,x)$.  Since it is continuous
in $v$, it is Borel with respect to all its arguments.
 This 
fact is obtained for
bounded continuous $f$, but by usual measure-theoretic
arguments  carries it over to all Borel bounded $f$.

Then take $0\leq r_{1}\leq...\leq r_{m}=r$
and continuous $f$
and a  bounded Borel  function $\zeta\big(x(1),...,x(m)\big)$
on $\bR^{md}$  and conclude from \eqref{12.14.02} that
$$
E_{t,x}\zeta(x_{r_{1} },...,x_{r_{m} })u(\sft_{r} ,x_{r})e^{-\lambda r}
$$
$$
=E_{t,x}\zeta(x_{r_{1} },...,x_{r_{m} })\int_{r}^{\infty}e^{-\lambda v}
f(\sft_{v},x_{v})\,dv.
$$
In light of \eqref{12.14.3} this means that
$$
\int_{r}^{\infty}E_{t,x}\zeta(x_{r_{1} },...,x_{r_{m} })
e^{-\lambda v}E_{\sft_{r} ,x_{r}}f(\sft_{v-r},x_{v-r})\,dv
$$
$$
=\int_{r}^{\infty}E_{t,x}\zeta(x_{r_{1} },...,x_{r_{m} })e^{-\lambda v}
f(\sft_{v},x_{v})\,dv.
$$
We have the equality of two Laplace's transforms 
of functions continuous in $v$. It follows that
for $v\geq r$
$$
E_{t,x}\zeta(x_{r_{1} },...,x_{r_{m} })
 E_{\sft_{r} ,x_{r}}f(\sft_{v-r},x_{v-r})
=E_{t,x}\zeta(x_{r_{1} },...,x_{r_{m} }) 
f(\sft_{v},x_{v}).
$$
Again a measure-theoretic argument shows that
this equality holds for any Borel bounded $f$ and then
  the arbitrariness of $\zeta$ yields the Markov property
of $X$.

To prove that it is strong Markov it suffices to 
observe that, owing to \eqref{12.14.3}
 its resolvent $R_{\lambda}$ is strong Feller, that is, maps bounded Borel functions into bounded continuous ones.

To deal with \eqref{7.14,1}, take, for instance,
$(t,x)=(0,0)$ and approximate our (conditionally weakly unique) solution
as in the proof of Theorem \ref{theorem 10.21.1}
by $x^{\varepsilon}_{\cdot}$.
For  
$R\in(0,\infty),y\in\bR^{d}$, 
introduce 
the functional $\gamma_{y,R}(x_{\cdot})$ 
 on $C([0,\infty),\bR^{d+1})$
 as the first exit time of $(s,x_{s})$
from $C_{R}(0, y)$. As is easy to see,
$\gamma_{y,R}(x_{\cdot})$ is lower semi-continuous.
It follows that the same is true
for
$$
 \int_{0}^{\gamma_{y,R}(x_{\cdot})}f( r,x_{r})\,dt,
$$
as  long as
a bounded continuous $f(t,x)\geq0$.
 Therefore, 
\begin{equation}
                                \label{11.26.1}
\nliminf_{n\to\infty} E_{0,0} \int_{0}^{\gamma_{y,R}(x^{\varepsilon_{m}}_{\cdot})}f( r,x^{\varepsilon_{m}}_{r})\,dt
\geq
E_{0,0} \int_{0}^{\gamma_{y,R}(x^{0}_{\cdot})}f( r,x^{0}_{r})\,dt.
\end{equation}
In light of \eqref{3.15.3}, 
inequality
\eqref{11.26.1} holds   for $f=|b |
  $. If $f=|b|$
and $R\leq \rho_{b} $, as   it follows
from \eqref{3.14.8},
the  left-hand side of \eqref{11.26.1}
is smaller that $\sfb_{0}R$. But then  
$$
E_{0,0}\int_{0}^{ \tau_{R}(y)}
|b( s, x^{0}_{ s}) |\,ds  \leq \sfb_{0}R,
$$
and this with the possibility to change the origin leads
 to  \eqref{7.14,1} and according to Definition \ref{definition 6.24,1} means that $X$
is regular.
The theorem is proved. \qed

\begin{remark}
                        \label{remark 11.28,1}
The regularity of $X$ implies that
all results of Section \ref{section 11.17,1}
are applicable. In particular, $X$ is strong Markov and strong Feller (see Corollary
\ref{corollary 9.5.1}). Also Harnack inequality and H\"older
continuity of the caloric functions are valid.

\end{remark}

\begin{corollary}
                  \label{corollary 12.15,3}
Suppose that on a probability space equation
\eqref{11.29.20} has an $E_{q,p,\beta}$-admissible solution $x_{s}$. Then \eqref{12.15,2} holds with the constant $N=N(T)$
from  \eqref{3.14.9} (depending only on
$d,\delta$, $q $, $p$, $\beta $, $\rho_{a}$,   $T$). We can choose $N=N(T)$
to be an increasing function of $T$. Then for any integer $n\geq 1$
and nonnegative $f$ we have
\begin{equation}
                             \label{12.15,1}
E\big(\int_{0}^{T}f(s,x_{s})\,ds\Big)^{n}
\leq n!N^{n}(T)\|  f\|^{n} _{E_{q ,p,\beta  }} .
\end{equation}
 
\end{corollary}
  
Indeed, by Theorem \ref{theorem 12.12.3} the
distribution of $x_{\cdot}$ coincides
with the distribution of the trajectories
of the Markov process $X$ under measure $P_{t,x}$. Then \eqref{12.15,1} is identical
with \eqref{11.22,4} for $n=1$. However, 
by Khasminskii's lemma due to the Markov
property of $X$ \eqref{11.22,4} implies that \eqref{12.15,1}
holds for any $n$ if we replace $E$ with $E_{t,x}$. Obviously, there is no need in doing the replacement.

Here is our most general existence and weak uniqueness theorem. We consider the equation
\begin{equation}
                            \label{12.16,2}
x _{s}= \int_{0}^{s}\sigma (r,x_{r})\,dw_{r}
+\int_{0}^{s}(b+\scB) (r,x_{r}) \,dr,
\end{equation}
where $\sigma, b$ are the same as at the beginning of the section, satisfying
\eqref{11.15,4} (that was declared to hold
until the end of this section before Theorem
\ref{theorem 3.15.1})  and $\scB=\scB(t,x)$
is a Borel $\bR^{d}$-valued function on
$\bR^{d+1}$, such that
\begin{equation}
                        \label{12.16,4}
\sfB^{2}:=\int_{\bR}\sup_{\bR^{d}}|\scB(t,x)|^{2}\,dt<\infty.
\end{equation}
(Just in case, observe
\index{$S$@Miscelenea!$\bar b_{R}$@$\sfB$}%
 that since $\scB$ is Borel,
$\sup_{\bR^{d}}|\scB(t,x)|$ is universally measurable, in particular, Lebesgue measurable,
so that \eqref{12.16,4} makes sense).
\begin{theorem}
                     \label{theorem 12.16,6}
In the above setting equation \eqref{12.16,2}
has an $E_{q,p,\beta}$-admis\-sible solution
on a probability space (weak existence),
and for any $E_{q,p,\beta}$-admissible solution on arbitrary probability space its finite-dimensional distributions
are independent of the solution
(weak uniqueness). Furthermore, for any $T\in(0,\infty)$ and integer $n\geq 1$
there exists a constant $N $
depending only on $d,\delta$, $q $, $p$, $\beta $, $\rho_{a}$,  $n$, $T$, and $\sfB$, such that for any
$E_{q,p,\beta}$-admissible solution
of \eqref{12.16,2} we have
\begin{equation}
                        \label{12.21,3}
E\big(\int_{0}^{T}f(s,x_{s})\,ds\Big)^{n}
\leq N\|  f\|^{n} _{E_{q ,p,\beta  }} .
\end{equation}

\end{theorem}

Proof. {\em Existence\/}. Take the Markov process $X$ from Theorem \ref{theorem 12.6.02}
and concentrate on the measure $P_{0,0}$
on $(\Omega,\cN)$, where $\Omega=C([0,\infty),
\bR^{d+1})$, and $\cN$ is the Borel $\sigma$-field on $\Omega$
completed with respect to $P_{0,0}$.
 We know that there exists
a $d$-dimensional Wiener process $w_{t}$
on  the probability space $(\Omega,\cN,P_{0,0})$
such that $x_{t}(\omega)=\omega(t)$ satisfies
\begin{equation}
                                \label{12.17,1}
x _{s}= \int_{0}^{s}\sigma (r,x_{r})\,dw_{r}
+\int_{0}^{s} b  (r,x_{r}) \,dr.
\end{equation}
Introduce a new probability measure $P$ on $(\Omega,\cN)$ by $P(d\omega)=e^{\phi}P_{0,0}
(d\omega)$, where
$$
\phi=\int_{0}^{\infty}\sigma^{-1}(s,x_{s})
\scB(s,x_{s})\,dw_{s}-(1/2)\int_{0}^{\infty}
|\sigma^{-1}(s,x_{s})
\scB(s,x_{s})|^{2}\,ds.
$$
It is well known that since $\sfB<\infty$, $P$
is indeed a probability measure equivalent to $P_{0,0}$.
Then 
equation \eqref{12.17,1} is rewritten as
\begin{equation}
                             \label{12.17,2}
x _{s}= \int_{0}^{s}\sigma (r,x_{r})\,d\bar w_{r}
+\int_{0}^{s}(b+\scB) (r,x_{r}) \,dr,
\end{equation}
where by Girsanov's theorem 
\begin{equation}
                             \label{12.17,3}
\bar w_{s}:=w_{s}-\int_{0}^{s}
\sigma^{-1}(r,x_{r})\scB(r,x_{r})\,dr
\end{equation}
is a Wiener process on $(\Omega,\cN,P)$. Furthermore, by Corollary \ref{corollary 12.15,3} for each $T\in(0,\infty)$, $n\geq1$, and Borel
$f\geq 0$
$$
E_{0,0}\Big(\int_{0}^{T}f(s,x_{s})\,ds\Big)^{2n}
\leq N\|f\|^{2n}_{E_{q,p,\beta}}.
$$
  It follows that
$$
E  \Big(\int_{0}^{T}f(s,x_{s})\,ds\Big)^{n} 
=E_{0,0}e^{\phi}\Big(\int_{0}^{T}f(s,x_{s})\,ds\Big)^{n}
$$
$$
\leq  \Big(E_{0,0}e^{2\phi}\Big)^{1/2}
\Big(E_{0,0}\Big(\int_{0}^{T}f(s,x_{s})\,ds\Big)^{2n}\Big)^{1/2}\leq N\Big(E_{0,0}e^{2\phi}\Big)^{1/2}\|f\|^{n} _{E_{q,p,\beta}}.
$$
Since, as is well known $E_{0,0}e^{2\phi}
\leq e^{N(\delta)\sfB^{2}}$, this shows that $x_{s}$
is an $E_{q,p,\beta}$-admissible solution
of \eqref{12.17,2} on $(\Omega,\cN,P)$
and \eqref{12.21,3} holds.

{\em Uniqueness\/}. Suppose that on a probability space $(\Omega,\cF,P)$ carrying a $d$-dimensional
Wiener process $\bar w_{s}$ equation
\eqref{12.17,2} has an $E_{q,p,\beta}$-admissible solution $x_{s}$. Then introduce
the process $w_{s}$ by using \eqref{12.17,3}
and a new probability measure $\bar P(d\omega)
=e^{-\phi}P(d\omega)$, where $\phi$ is the same
as above, so that $\phi=\phi_{\infty}$ with
$$
\phi_{t}=\int_{0}^{t}\sigma^{-1}(s,x_{s})
\scB(s,x_{s})\,d\bar w_{s}-(1/2)\int_{0}^{t}
|\sigma^{-1}(s,x_{s})
\scB(s,x_{s})|^{2}\,ds
$$
$$
=\int_{0}^{t}\sigma^{-1}(s,x_{s})
\scB(s,x_{s})\,d w_{s}+(1/2)\int_{0}^{t}
|\sigma^{-1}(s,x_{s})
\scB(s,x_{s})|^{2}\,ds.
$$
By Girsanov's theorem  $w_{s}$ is a Wiener
process on $(\Omega,\cF,\bar P)$ and $x_{s}$
is a solution of
\eqref{11.29.20} with $(t,x)=(0,0)$.

Now we prove that $x_{s}$ is a $(q,p,\beta)$-reasonable solution. To that effect define stopping times
$$
\tau^{n}=  \tau_{n}\wedge
\inf\{t\geq0:e^{-\phi_{t}}\geq n\},
$$
where $\tau_{n}$ is the first exit time
of $(s,x_{s})$ from $C_{n}$. Since 
 the exponential $\bar P$-martingale
$e^{-\phi_{t}}$ has bounded trajectories,   $\tau^{n}\to\infty$
as $n\to\infty$. In addition, for each $n$
and Borel $f\geq0$
$$
\bar E\int_{0}^{\tau_{n}}f(s,x_{s})\,ds
= Ee^{-\phi_{\tau^{n}}}\int_{0}^{\tau_{n}}f(s,x_{s})\,ds
$$  
$$
\leq e^{n}E \int_{0}^{n^{2}}f(s,x_{s})\,ds
\leq N\|f\|_{E_{q,p,\beta}},
$$
where $N$ is independent of $f$. It follows
that $x_{s}$ is indeed a $(p,q,\beta)$-reasonable solution on $(\Omega,\cF,\bar P)$
and by Theorem \ref{theorem 12.12.3} its
$\bar P$-finite-dimen\-sional distributions
are uniquely determined by $a,b$.

Finally, take a Borel bounded $F\geq0$ on
$C([0,\infty),\bR^{d})$ and observe that the formulas
$$
E F(x_{\cdot})=\bar Ee^{\phi}F(x_{\cdot}),
$$
$$
\phi=\int_{0}^{\infty}\sigma^{-1}(s,x_{s})
\scB(s,x_{s})\sigma^{-1}(s,x_{s})\,(dx_{s}-
b(r,x_{r})\,dr)
$$
$$
-(1/2)\int_{0}^{\infty}
|\sigma^{-1}(s,x_{s})
\scB(s,x_{s})|^{2}\,ds
$$
allow us to express $e^{\phi}F(x_{\cdot})$
as a function of $x_{\cdot}$ (involving
$\scB$) and show
that $E F(x_{\cdot})$ is uniquely determined by $a,b,\scB$ and $F$. This proves weak uniqueness. \qed

 The following result is
extended to all $\lambda>0$ in
Theorem \ref{theorem 11.26,1}. Recall that
$\rho_{b}=1$.
\begin{theorem}  
              \label{theorem 7,3.1}
Let 
$\lambda\geq   \lambda_{0}\rho_{b}^{-2}$
and let $f\in E_{q,p,\beta}$. Then
\begin{equation}
                   \label{7,3.7}
u(t,x)=E_{t,x}\int_{0}^{\infty}
e^{-\lambda s}f(t+s,x_{s})\,ds 
\end{equation}
belongs to $E^{1,2}_{q,p,\beta}$
and is a unique solution of
class $E^{1,2}_{q,p,\beta}$ of
equation \eqref{10.14,1}. Furthermore, $u$ is a uniformly (H\"older) continuous function.

\end{theorem}  

Proof. If $f$ is smooth and bounded, then $f\in E_{q,p,\beta}$,
  \eqref{10.14,1} 
has a unique solution $u\in  E^{1,2}_{q,p,\beta}$ and It\^o's formula easily shows that \eqref{7,3.7} holds.

For general $f$ apply the above argument to the standard mollification
$f^{(\varepsilon)}$ of $f$ and call
$u^{\varepsilon}$ the right-hand side of \eqref{7,3.7} with $f^{(\varepsilon)}$ in place of $f$.
Due to \eqref{11.22,4}  
and Lemma \ref{lemma 2.16.3} we have
$u^{\varepsilon}\to u$ as $\varepsilon
\downarrow 0$ and, since $u^{\varepsilon}$
are uniformly continuous, $u$
is continuous.  
Then $D^{2}u^{\varepsilon}$ are
uniformly bounded in $L_{q,p}(C)$
for any cylinder $C$ and, since
$u^{\varepsilon}\to u$, these derivatives converge weakly
 to $D^{2}u$, which thus exist.
The same argument applies to $Du,\partial_{t}u$. It follows, in particular, that $u$
satisfies \eqref{10.14,1}.

Then, for a given cylinder $C\in\bC_{\rho}$ with $\rho\leq 1$ in light of weak convergence we have
$$
\rho^{\beta}\dashnorm D^{2}u\|_{L_{q,p}(C)}
\leq \nliminf_{\varepsilon\downarrow 0}
\rho^{\beta}\dashnorm D^{2}u^{\varepsilon}\|_{L_{q,p}(C)}
$$
$$
\leq N
  \nliminf_{\varepsilon\downarrow 0}
\|f^{(\varepsilon)}\|_{E_{q,p,\beta}}
\leq N\|f \|_{E_{q,p,\beta}}.
$$
Similar relations hold for $u$ and $Du$
and this shows that $u\in E^{1,2}_{q,p,\beta}$. This and the uniqueness in Theorem \ref{theorem 5.8,20} finish
the proof of the current theorem. \qed
 
\begin{theorem}
                       \label{theorem 11.26,1}
  Assertion of Theorem \ref{theorem 7,3.1}
holds true for any $\lambda>0$.
\end{theorem}

Proof. Fix $\lambda>0$. 
We may assume that $\lambda< \lambda_{0}\rho_{b}^{-2}$.
Our first goal is to prove
that
\begin{equation}
                           \label{11.26,1}
E_{t,x}\int_{0}^{\infty}e^{-\lambda s}
|g(t+s,x_{s})|\,ds \leq N\|g\|_{E_{q,p,\beta}},
\end{equation}
where $N$ is independent of $g,t,x$. To do that
observe that we may concentrate on $g\geq0$,
take a constant $K>0$ and set
$$
v(t,x)=E_{t,x}\int_{0}^{\infty}e^{-\lambda s}
K\wedge g(t+s,x_{s})|\,ds. 
$$
By the Markov property, Theorem \ref{theorem 7,3.1},   Theorem \ref{theorem 5.8,20}, and 
Lemma {lemma 11.15,1}
$$
v(t,x)=E_{t,x}\int_{0}^{1}e^{-\lambda s}
K\wedge g(t+s,x_{s})\,ds+e^{-1}E_{t,x}u(t+1,x_{1})
$$
$$
\leq e^{\lambda_{0}\rho_{b}^{-2}}
E_{t,x}\int_{0}^{\infty}e^{-\lambda_{0}\rho_{b}^{-2}s}
K\wedge g(t+s,x_{s})\,ds+e^{-1}\sup u
$$
$$
\leq N\|g\|_{E_{q,p,\beta}}+e^{-1}\sup u.
$$
It follows that $\sup u \leq N\|g\|_{E_{q,p,\beta}}+e^{-1}\sup u$, $\sup u \leq N\|g\|_{E_{q,p,\beta}}$, and to get \eqref{11.26,1} it only remains to let $K\to\infty$.

Now   by referring to the Markov property again, which yields
$$
u(t,x)=E_{t,x}\int_{0}^{\infty}
e^{-(\lambda+\lambda_{0}\rho_{b}^{-2}) s}
\big(f+\lambda_{0}\rho_{b}^{-2}u\big)(t+s,x_{s})\,ds,
$$  
we conclude that indeed $u$ belongs to $E^{1,2}_{q,p,\beta}$
and is a   solution of
class $E^{1,2}_{q,p,\beta}$ of
equation \eqref{10.14,1}. If there are two such
solutions, then their difference $w$ satisfies
$\cL w-\lambda_{0}\rho_{b}^{-2}w=(\lambda-\lambda_{0}\rho_{b}^{-2})w $. By Theorem \ref{theorem 7,3.1} then
$$
|w(t,x)|=\Big|E_{t,x}\int_{0}^{\infty}
e^{-\lambda_{0}\rho_{b}^{-2} s}(\lambda_{0}\rho_{b}^{-2}-\lambda)w(t+s,x_{s})\,ds\Big|
$$
$$
\leq \frac{\lambda_{0}\rho_{b}^{-2}-\lambda}{ \lambda_{0}\rho_{b}^{-2}}\sup|w|.
$$
This implies $w=0$ and the theorem is proved. \qed

\mychapter[Strong solutions]{Strong solutions}

 In this chapter as everywhere we suppose that $d\geq 2$ and,
 in addition, we suppose that
for an integer $d_{1}\geq d$ on $\bR^{d+1}$
we are given a Borel $d\times d_{1}$-valued function  $\sigma$  and $\bR^{d}$-valued functions $b,\scB$. Assume that $a:=\sigma\sigma^*$ is $\bS_{\delta}$-valued.
The main object in this chapter is the equation
 \begin{equation}
                         \label{3.15,1}
 x_{s}=x_{0} +\int_{0}^{s}\sigma(t+r,x_{r})\,dw_{r}+  
\int_{0}^{s}(b+\scB)(t+r,x_{r})\,dr,
 \end{equation}
where $w_{t}$ is a $d_{1}$-dimensional
Wiener process given on a probability 
space and nonrandom $x_{0}\in\bR^{d}$.
The role of $\scB$ will become clear toward
the end of the chapter and
in the following discussion we assume that $\scB\equiv0$.

After the classical work by K. It\^o showing that there exists
a unique strong solution of \eqref{3.15,1} if $\sigma $ and $b$
are Lipschitz continuous in $x$ (may also depend on   $\omega$),
much effort was applied to relax these   conditions.
In   case $d=d_{1}=1$ T. Yamada and S. Watanabe \cite{YW_71} relaxed
the Lipschitz condition on $\sigma$ to the H\"older $(1/2)$-condition
(and even slightly weaker condition) and kept $b$ Lipschitz
(slightly less restrictive). Much attention was paid to equations
with continuous coefficients
satisfying the so-called monotonicity conditions
(see, for instance, \cite{Kr_84} and the references therein).

T. Yamada and S. Watanabe \cite{YW_71} also put forward
a very strong theorem, basically, saying that the existence of weak solutions and strong uniqueness
implies the existence of strong
solutions. Unlike the present book, 
the majority of papers on the subject after that time
are using their theorem.
 S. Nakao (\cite{Na_72}) proved
the strong solvability in time homogeneous case
 if $d=d_{1}=1$ and $\sigma$ is bounded away from zero and infinity
 and is
locally of bounded variation. He also assumed that $b$ is bounded,
but from his arguments it is clear that the summability of $|b|$
suffices. In this respect his result basically shows that
our results are also true if $d=1$
and the coefficients are independent of time. 
Our results are applicable for $d\geq2$. However, in the 
 case that $d=2$ they do not look very satisfactory because $\sigma$ turns out to be H\"older continuous and $b$ locally summable to the power $>2=d$ excluding the singularities
like $1/|x|$.

A. Veretennikov was the first author who in \cite{Ve_80} 
not only proved the existence of strong solutions
in the time inhomogeneous {\em multidimensional\/} case  when $b$ is bounded,
but also considered the case of $\sigma $ in Sobolev class,
namely,   $\sigma _{x}\in L_{2d,\loc}$. He used A. Zvonkin's method (see
\cite{Zv_74}) of transforming the equation in such a way that the
drift term disappears.
In \cite{Zh_11}, \cite{Zh_16}, and \cite{XXZZ_20} (also see the references
there) the result of Veretennikov is
extended to the case of $\sigma$ uniformly continuous in $x$ and $\sigma_{x},b\in L_{q,p}$
with, perhaps, different $p,q$ for $\sigma_{x}$ and $b$ satisfying
\begin{equation}
                           \label{6.12.1}
\frac{d}{p}+\frac{2}{q}< 1
\end{equation}
(the so-called subcritical Ladyzhenskaya-Prodi-Serrin
condition).
In that case much information is available,
we refer the reader
to \cite{Zh_11},  \cite{XZ_20}, \cite{XXZZ_20}, and the references therein.  

Even the case when   $\sigma $  is constant and the process is nondegenerate
attracted
very much attention.   
M. R\"ockner and the author in \cite{KR_05} proved,
among other things, the existence of strong solutions 
when $b\in L_{q,p}$ under condition \eqref{6.12.1}.
We  refer to \cite{PFPR_13}, \cite{MNF_15},  \cite{Zh_16} and the references therein for further results in  this direction.
If $b$ is bounded  A.~Shaposhnikov (\cite{Sh_14}, \cite{Sh_17})
proved the so called path-by-path uniqueness, which, basically,
means that for almost any trajectory $w_{t}$ there is only one
solution (adapted or not). This result was already announced
by A. Davie before with a very entangled proof which left many doubtful.

In the fundamental work by L.~Beck, F.~Flandoli, M.~Gubinelli, and M.~Maurelli
(\cite{BFGM_19}) the authors investigate such equations from
 the points of view of It\^o stochastic equations, stochastic transport 
equations, and stochastic continuity equations. Their article
contains an enormous amount of information and a vast references list. 
In what concerns our situation they require ($\sigma=(\delta^{ij}) $  
and) what they call LPS-condition (slightly imprecise): $b\in L_{q,p,\loc}$,
$q<\infty$, with equality in \eqref{6.12.1} in place of $<$,
or $p=d$ but $\|b\|_{L_{\infty},p}$ to be sufficiently small, or else that $b(t,\cdot)$ to
be continuous as an $L_{d}(\bR^{d})$-function,
and they prove strong solvability and strong uniqueness
(actually, path-by-path-uniqueness which is stronger) but only for {\em almost
all\/} starting points $x$. 

Concerning the strong solutions
starting from any point $x$
in the  time dependent  case
with singular $b$  and constant $\sigma$ probably the best well elaborated results   belong  to R\"ockner and Zhao \cite{RZ_25}, where, among  many other things, they prove existence and uniqueness of strong solutions of
equations like \eqref{3.15,1} with $b\in L_{q,p}$ and $q,p<\infty$, 
with equality in \eqref{6.12.1} in place of $<$, or when $b(t,\cdot)$
is continuous as an $L_{d}(\bR^{d})$-function.
In what concerns the {\em existence and uniqueness
of strong solutions\/} the results in \cite{RZ_25}
are covered by more general results in
\cite{Kr_25_2} contained also in the present book. These more general results are proved
for $b$ in Morrey classes.

In the paper by D. Kinzebulatov and K.R. Madou
\cite{KM_24}  conditions on $b$
are different from  \cite{RZ_25} and \cite{Kr_25_2}. They are stated in terms of
form-boundedness and allow the authors to
  prove  strong
solvability when $b(t,x)$ is form-bounded
for each $t$ with bound uniform in $t$. 
This class of $b$   contains the set of functions
$b$ such that $b(t,\cdot)$ is in a Morrey class
with the norm uniformly bounded in $t$
as required for one part of the drift in \cite{Kr_25_2}, but does not contain the set of functions
$b$ which are in a Morrey class {\em with respect 
to $(t,x)$\/}, for which we prove the
strong solvability. Also it does not contain
major part of the functions $b$ from \cite{RZ_25}.

We refer the reader to
\cite{BFGM_19} and \cite{RZ_25} also for
a very good review  of the motivation
related to the Navier-Stokes equation and history of the problem.

Our approach is absolutely different from
all articles mentioned above and all articles
which one can find in their references.
 We do not use
Yamada-Watanabe theorem or transformations of the noise or   a compactness criterion for random
fields in Wiener-Sobolev spaces as in \cite{RZ_25} and \cite{KM_24}.
Instead, our method is inspired by
an analytic criterion for the existence
of strong solutions which first appeared in \cite{VK_76},
 some 50 years ago and was first used only
in \cite{Kr_21}.
   To make this method work we use  ideas
from many papers, most relevant of which are   
\cite{DK_18}, \cite{Kr_21}, \cite{1}, \cite{Kr_25_1}, \cite{13}, \cite{VK_76}.

Here is an example in which we prove existence
and conditional uniqueness of strong solutions. Take $d=3,d_{1}=12$, and
for some numbers $\alpha,\beta,\gamma\geq0$
let $\sigma^{k}$ be the $k$th column in
 ($0/0:=3^{-1/2}$)
\begin{equation}
                                    \label{6.3.4}
\begin{pmatrix}\alpha & 0 & 0\\
0 & \alpha & 0\\
0 & 0 & \alpha
\end{pmatrix},
\frac{\beta}{|x|}
\begin{pmatrix}
x^{1} &  x^{2} & x^{3} & 0 &  0   & 0   & 0  & 0  & 0 \\
0 & 0 & 0 & x^{1}  & x^{2}   & x^{3}   & 0  & 0  & 0 \\
0 & 0 & 0 & 0  & 0   & 0  & x^{1}& x^{2} & x^{3} 
\end{pmatrix} ,
\end{equation}
$$
b(x)=-\frac{\gamma}{|x|}\,\frac{x}{|x|}I_{0<|x|\leq 1} +\xi(t)\eta(t,x),
$$
 where $\eta$ is bounded 
$\bR^{3}$-valued  and $\xi$ is real-valued of class $
 L_{2}(\bR)$.
Our result shows that if $\alpha=1$ and
$\beta$ and $\gamma$ are sufficiently small,
then \eqref{3.15,1} has a strong solution
which is conditionally unique, however, if
$\xi\equiv0$, then any solution is strong and unique.
By the way, if $\xi\equiv0$, $\alpha=\gamma=0$ and $\beta=1$, there exist
strong solutions of \eqref{3.15,1}
only if the starting point $x\ne0$ (see \cite{Kr_21}). In case $\alpha=1$ and $\beta=0$
strong solutions exist only if $\gamma$ is sufficiently small.  In case   $\alpha=1,\beta=\gamma=0$,
absent in \cite{RZ_25} and \cite{KM_24}, the authors of 
\cite{BFGM_19} prove the unique strong solvability
only for almost all starting points. 
We prove the unique strong solvability for any starting
point.

Observe that for $\beta\ne 0$
and $\gamma\ne0$ we have $D\sigma ,b(t,\cdot)\in L_{d-\varepsilon,\loc}(\bR^{d})$ for any $\varepsilon\in(0,1)$
but not for $\varepsilon=0$. Recall  that the case
of time independent $\sigma,b$ with $D\sigma ,b\in L_{d,\loc}$ is investigated
in \cite{Kr_21} the main idea of which
is used here as well.

Other examples can be found in Remarks  
\ref{remark 6.2.1} 
 and \ref{remark 3.21,3}. There as above we compare
our results with the ones obtained when $\sigma$
is the unit matrix. In this connection note that
our results are new even if $b\equiv0$.

There is another active direction in the investigation of the
strong solutions when $b$ has some singularity
in time but also somewhat regular in space,
see, for instance, \cite{GG_25}, \cite{WHY_25}, and the references therein. This area is out of the scope of the book.

Set
$$
\cL=\partial_{t}+(1/2)a^{ij}D_{ij}+(b^{i}+\scB^{i})D_{i}.
$$
Fix some  
$$
\rho_{0} \in(0,\infty),\quad p_{0}   \in(2,2+d] .  
$$
Suppose that $Da\in L_{1,\loc}(\bR^{d+1})$ and 
\index{$S$@Miscelenea!$a_{\pm}$@$\widehat {Da}_{s,\rho}$}%
\index{$S$@Miscelenea!$\bar b_{R}$@$\hat b_{s,\rho}$}%
introduce
$$
\widehat {Da}_{s,\rho}=\sup_{r\leq\rho  }r
\sup_{C\in \bC_{r}} 
\dashnorm Da \|_{L_{s}(C)},\quad
\hat b_{s,\rho}=\sup_{r\leq\rho  }r
\sup_{C\in \bC_{r}} 
\dashnorm b \|_{L_{s}(C)},
$$
$$
\scB(t)=\sup_{\bR^{d}}|\scB(t,x)|.
$$

It is easy to
\index{$S$@Miscelenea!$\bar b_{R}$@$\scB(t)$}%
 see that 
there is $q_{0}\in(2,p_{0}]$ such that $1\leq d/p_{0}+2/q_{0}$ and for some $\beta\in(1,2)$ and $p= p_{0}/\beta,q=q_{0}/\beta$
we have $p>2,q>2$ and
$$
\beta\leq \frac{d}{p}+\frac{2}{q}.  
$$
Note that $p_{0},q_{0}$ can be taken arbitrarily close to $2$.
Also fix $\rho_{a},\rho_{b}\in[\rho_{0},\infty)$.

Some additional assumptions on the above objects
are stated in Section \ref{section 4.16,1} and are supposed
to hold ever after until the end of Section~\ref{section 7.3.1}.

\mysection[A preliminary estimate. Case $\scB=0$]
{A preliminary estimate. Case $\scB=0$}

\begin{remark}
                         \label{remark 11.30,1}
If $C\in\bC_{\rho_{0}}$, then
as is not hard to prove
$$
\|I_{C}Da\|_{\dot E_{p_{0},p_{0},1}}\leq
\widehat { Da }_{p_{0} ,\rho_{0} },\quad 
\|I_{C}b\|_{\dot E_{p_{0},p_{0},1}}\leq
\hat b_{p_{0} ,\rho_{0} },\quad \|I_{C}\|
_{\dot E_{p_{0},p_{0},1}}\leq \rho_{0}.
$$
Also if $\zeta\in C^{\infty}_{0}(C_{1})$ 
with the integral
of its square equal to one and
$\zeta_{\rho_{0}}(t,x)=\rho_{0}^{-(d+2)/2}
\zeta (t/\rho_{0}^{2},x/\rho_{0})$, then 
$\zeta_{\rho_{0}}$ is in
$C^{\infty}_{0}(C_{\rho_{0}})$, the integral
of its square equals one, and
$$
\rho_{0}^{2}\int_{C_{\rho_{0}}}|D\zeta_{\rho_{0}}|^{2}
\,dxdt=\int_{C_{1}}|D\zeta |^{2}
\,dxdt.
$$
\end{remark}
\begin{lemma}
                        \label{lemma 11.26,3}
Suppose that $a,b$ are infinitely differentiable
with each derivative bounded. Let $f\in C^{\infty}_{0}(\bR^{d})$, $T>0$, and let $u(t,x)$ be the classical solution of
\begin{equation}
                          \label{11.26,5}
\cL u=0\quad\text{in}\quad [0,T]\times \bR^{d}
\end{equation}
with boundary condition $u(T,x)=f(x)$.  
Let $n\in\{iI_{i=1}+2iI_{i\geq2},i=1,2,...\}, \lambda\geq 0$.
Then there
are constants $\widehat{Da},\hat b\in(0,1)$,
depending only on $d,\delta,p_{0}$,  
$n $,
such that if $\widehat { Da }_{p_{0},\rho_{0}}\leq 
e^{- \lambda\rho_{0}}\widehat{Da}$ and $\hat b_{p_{0} , \rho_{0} }
\leq e^{- \lambda\rho_{0}}\hat b$, then
$$
\int_{\bR^{d}}|u(0,x)|^{2n}e^{-\lambda|x|}\,dx
\leq Ne^{ \alpha T}\int_{\bR^{d}}|f(x)|^{2n}e^{-\lambda|x|}\,dx,
$$
$$
\int_{[0,T]\times \bR^{d}} u  ^{2n-2}
|Du |^{2}e^{-\lambda|x|}\,dxdt\leq N 
e^{\lambda
\rho_{0}+\alpha T}\int_{\bR^{d}}|f(x)|^{2n}e^{-\lambda|x|}\,dx,
$$
where 
$$
\alpha=N\rho_{0}^{-2}e^{ \lambda\rho_{0}}
$$
and the constants called $N$ depend  only on $d,\delta,p_{0}$,  
$n $.
\end{lemma}

Proof. Without restricting   generality
we assume that $\widehat { Da }_{p_{0},\rho_{0}}
\leq 1$, $\hat b_{p_{0} , \rho_{0} }
\leq 1$. Take a $C\in\bC_{\rho_{0}}$ and 
a nonnegative $\zeta
\in C^{\infty}_{0}(C)$  with the integral
of its square equal to one.
We    multiply  \eqref{11.26,5}
by $\zeta^{2} u^{2n-1}$ and integrate by parts. Then noting that $a^{ij}D_{i}uD_{j}u\geq\delta|Du|^{2}$ we find for $s\leq T$ that
$$
\int_{\bR^{d}}\zeta^{2}(s,x)u^{2n}(s,x)\,dx+
\int_{[s,T]\times \bR^{d}}\zeta^{2} u  ^{2n-2}
|Du |^{2}\,dxdt
$$
$$
\leq  N\int_{\bR^{d}}\zeta^{2}(T,\cdot)f ^{2n }   \,dx+N\int_{[s,T]\times \bR^{d}}
u^{2n}|\partial_{t}\zeta^{2}|\,dxdt 
$$
$$+N\int_{[s,T]\times \bR^{d}} |u^{n} D \zeta   |\zeta(u^{n-1}|Du|)\,dxdt
$$
$$
+N\int_{[0,T]\times \bR^{d}}\zeta^{2}\big(\big(|Da|+|b| )u^{n}\big)\big(
u^{n-1}|Du|\big)\,dxdt.
$$
The last term is dominated by
$$
(1/2)\int_{[s,T]\times \bR^{d}}\zeta^{2} u  ^{2n-2}
|Du |^{2}\,dxdt+N\int_{[s,T]\times \bR^{d}} \zeta^{2}(|Da|+|b| )^{2}u^{2n} \,dxdt
$$
and the previous one is dominated by
$$
(1/4)\int_{[s,T]\times \bR^{d}}\zeta^{2}
u^{2n-2}|Du|^{2}\,dxdt+N \int_{[s,T]\times \bR^{d}}| D \zeta| ^{2}u^{2n}\,dxdt.
$$ 
 
It follows that
$$
\int_{\bR^{d}}\zeta^{2}(s,x)u^{2n}(s,x)\,dx+
\int_{[s,T]\times \bR^{d}} \zeta^{2} u  ^{2n-2}
|Du |^{2}\,dxdt
$$
$$
\leq N\int_{\bR^{d}}\zeta^{2}(T,\cdot)f ^{2n }   \,dx
+N \int_{[s,T]\times \bR^{d}}\big(|\partial_{t}\zeta^{2}|+|D\zeta |^{2}\big)
  u^{2n}\,dxdt
$$
\begin{equation}
                               \label{11.28,4}
+N\int_{[s,T]\times \bR^{d}} \zeta^{2}(|Da|+|b|)^{2}u^{2n} \,dxdt.
\end{equation}

Before proceeding further we note that we may
look at $\zeta$ as a scaled and translated
function with support in $C_{1}$. Then it is seen that 
\begin{equation}
                              \label{1.21.3}
\rho_{0}^{d+2} |\zeta |^{2}
+
\rho_{0}^{d+4} |D\zeta |^{2}
+
\rho_{0}^{d+6} |\partial_{t}\zeta |^{2}
 \leq N (d),
\end{equation}
and we infer from \eqref{11.28,4} that
$$
\int_{\bR^{d}}\zeta^{2}(s,x)u^{2n}(s,x)\,dx+
\int_{[s,T]\times \bR^{d}} \zeta^{2} u  ^{2n-2}
|Du |^{2}\,dxdt
$$
$$
\leq N\int_{\bR^{d}}\zeta^{2}(T,\cdot)f ^{2n }   \,dx
+N \rho_{0}^{-d-4}\int_{[s,T]\times \bR^{d}}I_{C}
  u^{2n}\,dxdt
$$
\begin{equation}
                               \label{2.10,1}
+N\int_{[s,T]\times \bR^{d}} \zeta^{2}(|Da|+|b|)^{2}u^{2n} \,dxdt.
\end{equation}

To estimate the last term it is convenient
to transform \eqref{11.26,5}.  
For $v:=u^{n}$ we have
$$
\partial_{t}(\zeta v)+ \Delta(\zeta v)+2\zeta\Big(\frac{1}{n}-1\Big)a^{ij}(D_{i}(u^{n/2}))D_{j}(u^{n/2})
$$
$$
+(1/2)\zeta a^{ij}D_{ij}v  
+\zeta  b^{i} D_{i}v -v\partial_{t}\zeta- \Delta(\zeta v)=0 
$$
for $t\leq T$ with boundary value $\zeta v|_{t=T}
=\zeta(T,\cdot)f^{n}$.

Here the second term is either zero if $n=1$
or negative if $n\geq 2$ when $v\geq0$.
Then by It\^o's formula, applied to $(\zeta v)
(t+r,x+\sqrt{2}w_{t})$, where $w_{t}$ is 
a $d$-dimensional Wiener process,
we get  that
$$
\int_{[s,T]\times \bR^{d}}|b|^{2}\zeta^{2}v^{2}\,dxdt
\leq N \int_{\bR^{d+1}_{s}}|b|^{2}I_{C  }P^{2}_{2,4}(F)\,dxdt
$$
\begin{equation}
                                 \label{11.28,5}
+
N \int_{[s,T]\times \bR^{d}}|b|^{2}I_{C }
\hat T^{2}_{T-t}[\zeta(T,\cdot)f^{n}](x)\,dxdt,
\end{equation}
where $\hat T_{r}h(x)=Eh(x+\sqrt2 w_{r})$
and
$F(t,x)=0$ for $t\geq T$ and for $t<T$  
$$
F=(1/2)\zeta a^{ij} D_{ij}v+ \zeta  b^{i} D_{i}v- v\partial_{t}\zeta- \Delta(\zeta v).
$$

Note that 
$$
P_{2,4}(\zeta a^{ij} D_{ij}v)=
 D_{j}\big(P_{2,4}(\zeta a^{ij} D_{i}v)\big)
-NP_{1,4}P_{1,4}\big(D_{j}(\zeta a^{ij})  D_{i}v\big)
$$
and since $|DP_{2,4}h|\leq NP_{1,8}|h|$, we have by Theorem
\ref{theorem 5.25,1}
$$
\int_{\bR^{d+1}_{s}}|b|^{2}I_{C  }
\big| D_{j}P_{2,4}(\zeta a^{ij} D_{i}vI_{(s,T)})\big  |^{2}\,dxdt
$$
$$
\leq N\hat b_{p_{0},\rho_{0}}^{2}
\int_{[s,T]\times \bR^{d}}\zeta^{2}|Dv|^{2}\,dxdt\leq N\hat b_{p_{0},\rho_{0}}^{2}\rho_{0}^{-d-2}
\int_{[s,T]\times \bR^{d}}I_{C}|Dv|^{2}\,dxdt.
$$
Similarly, invoking also Corollary \ref{corollary 10.5,1} and Remark \ref{remark 11.30,1} we get
$$
\int_{\bR^{d+1}_{s}}|b|^{2}I_{C  }
\Big(P_{1,4}P_{1,4}\big(|D(\zeta a)|\, Dv|I_{(s,T)}\big)\Big)^{2}\,dxdt
$$
$$
\leq N\hat b_{p_{0},\rho_{0}}^{2}
\int_{\bR^{d+1}_{s}}\Big( P_{1,4}\big(\{I_{C}
D\zeta|+I_{C}|D a|\zeta\}\, |Dv|I_{(s,T)}\big)\Big)^{2}\,dxdt 
$$
$$
\leq
N\hat b_{p_{0},\rho_{0}}^{2}\int_{[s,T]\times \bR^{d}}(\rho_{0}^{2}|D\zeta|^{2}+\zeta^{2})|Dv|^{2}\,dxdt
$$
$$
\leq N\hat b_{p_{0},\rho_{0}}^{2}\rho_{0}^{-d-2}
\int_{[s,T]\times \bR^{d}}I_{C}|Dv|^{2}\,dxdt.
$$
In the same way applying Theorem
\ref{theorem 5.25,1} and Corollary \ref{corollary 10.5,1} we obtain
$$
\int_{\bR^{d+1}_{s}}|b|^{2}I_{C }
 P^{2}_{2,4} (\zeta b^{i}D_{i}vI_{(s,T)}) \,dxdt\leq N\hat b_{p_{0},\rho_{0}}^{2}\rho_{0}^{-d-2}
\int_{[s,T]\times \bR^{d}}I_{C}|Dv|^{2}\,dxdt,
$$ 
$$
\int_{\bR^{d+1}_{s}}|b|^{2}I_{C }
 P^{2}_{2,4} (I_{C}vI_{(s,T)}\partial_{t}\zeta) \,dxdt\leq N\hat b_{p_{0},\rho_{0}}^{2}\int_{[s,T]\times \bR^{d}}\rho_{0}^{2}|\partial_{t}\zeta|^{2}| v|^{2}\,dxdt
$$
$$
\leq N\hat b_{p_{0},\rho_{0}}^{2}\rho_{0}^{-d-4}
\int_{[s,T]\times \bR^{d}}I_{C}v^{2}\,dxdt,
$$
$$
\int_{\bR^{d+1}_{s}}|b|^{2}I_{C }
 P^{2}_{2,4} (\Delta(\zeta v)I_{(s,T)}) \,dxdt
$$
$$
\leq 
\int_{\bR^{d+1}_{s}}|b|^{2}I_{C }
P^{2}_{1,8}(v|D\zeta|I_{(s,T)}+\zeta|Dv|I_{(s,T)})\,dxdt
$$
$$
\leq
N\hat b_{p_{0},\rho_{0}}^{2}\int_{[s,T]\times \bR^{d}} ( |D\zeta|^{2} v ^{2}+\zeta^{2}|Dv|^{2})\,dxdt
$$
$$
\leq
N\hat b_{p_{0},\rho_{0}}^{2}\int_{[s,T]\times \bR^{d}}I_{C} ( \rho_{0}^{-d-4} v ^{2}+\rho_{0}^{-d-2}|Dv|^{2})\,dxdt.
$$

To finish dealing with \eqref{11.28,5} we apply
Lemma \ref{lemma 11.29,1}
to estimate the last term and get 
$$
\int_{[s,T]\times \bR^{d}}|b|^{2}\zeta^{2}v^{2}\,dxdt
\leq N\hat b_{p_{0},\rho_{0}}^{2}\int_{\bR^{d}}\zeta^{2}(T,\cdot)
f^{2n}\,dx
$$
$$
+
N\hat b_{p_{0},\rho_{0}}^{2}\int_{[s,T]\times \bR^{d}} I_{C}( \rho_{0}^{-d-4} v ^{2}+\rho_{0}^{-d-2}|Dv|^{2})\,dxdt.
$$
Estimating the integral of $|Da|^{2}\zeta^{2}v^{2}
$ is not much different and therefore coming back to \eqref{2.10,1} we conclude that
$$
\int_{\bR^{d}}\zeta^{2}(s,x)u^{2n}(s,x)\,dx+
\int_{[s,T]\times \bR^{d}} \zeta^{2}  u  ^{2n-2}
|Du |^{2}\,dxdt
$$
$$
\leq N\int_{\bR^{d}}\zeta^{2}(T,\cdot)f ^{2n }   \,dx
+N \rho_{0}^{-d-4}\int_{[s,T]\times \bR^{d}}I_{C} u^{2n}\,dxdt
$$
$$
+
N(\hat b_{p_{0},\rho_{0}}^{2}+\widehat{Da}_{p_{0},\rho_{0}}^{2})\rho_{0}^{-d-2}
\int_{[s,T]\times \bR^{d}} 
I_{C}  u  ^{2n-2}
|Du |^{2}\,dxdt.
$$
We substitute here $C=C_{\rho_{0}}(\tau,\xi)$
and $\zeta(t-\tau,x-\xi)$ in place of $\zeta$,
where $(\tau,\xi)\in\bR^{d+1}$. Then we multiply
both parts by $e^{-\lambda |\xi|}$ and integrate
through the resulting inequality with respect to
$(\tau,\xi)\in\bR^{d+1}$. 
At this point it is worth mentioning that since $f\in C^{\infty}_{0}$ and $a$ and $b$
are sufficiently regular, $u$ and its derivatives
go to zero as $|x|\to\infty$ exponentially fast.
Therefore, our manipulations are well justified.
 
Observe that
$$
e^{\lambda \rho_{0}}e^{-\lambda |x|}\geq \int_{\bR^{d+1}}\zeta^{2}(t-\tau,x-\xi)e^{-\lambda |\xi|}\,d\xi d\tau\geq  
e^{-\lambda \rho_{0}}e^{-\lambda |x|},
$$
$$
\int_{\bR^{d+1} }
I_{C_{\rho_{0}}} (t-\tau,x-\xi)e^{-\lambda |\xi|}\,d\xi d\tau
\leq N\rho^{d+2}e^{\lambda \rho_{0}}e^{-\lambda |x|}.
$$
 
Therefore,   we find that
$$
e^{\lambda\rho_{0}}
\int_{\bR^{d}}u^{2n}(s,x)e^{-\lambda|x|}\,dx
+\int_{[s,T]\times \bR^{d}}  u  ^{2n-2}
|Du |^{2}e^{-\lambda|x|}\,dxdt
$$
$$
\leq Ne^{ \lambda \rho_{0}} \int_{\bR^{d}}f^{2n}
e^{-\lambda|x|}\,dx +N e^{2\lambda \rho_{0}}
\rho_{0}^{-2}\int_{[s,T]\times \bR^{d}}  u^{2n}e^{-\lambda|x|}\,dxdt
$$
\begin{equation}
                                 \label{12.9,4}
+
N_{1}e^{2\lambda \rho_{0}}(\hat b_{p_{0},\rho_{0}}^{2}+\widehat{Da}_{p_{0},\rho_{0}}^{2})
\int_{[s,T]\times \bR^{d}} 
   u  ^{2n-2}
|Du |^{2}e^{-\lambda|x|}\,dxdt.
\end{equation}
The last term is absorbed into the left-hand side
if we require
$$
N_{1}e^{2\lambda \rho_{0}}(\hat b_{p_{0},\rho_{0}}^{2}+\widehat{Da}_{p_{0},\rho_{0}}^{2})\leq 1/2
$$
and to deal with the previous one we use 
Gronwall's inequality after throwing away 
the second term on the left. The lemma is proved. \qed

 \mysection[Evolution family $T_{s,t}$. 
 Case $\scB=0$]{Evolution family $T_{s,t}$. Case $\scB=0$}

                                  \label{section 4.16,1}
In this section we suppose that the assumptions 
of Theorem \ref{theorem 5.8,20} are satisfied with 
$\beta,p=p_{0}/\beta,q=q_{0}/\beta,\rho_{a},\rho_{b}$
  specified in the introduction to
the chapter. Also suppose
that 
$$
N_{1}\hat b_{q_{0},p_{0},\rho_{b}}<\sfb_{0},
$$
where $N_{1}=N_{1}(d,\delta,q ,p,\beta ,\rho_{a} )$ 
is taken from \eqref{11.15,4}. More precisely,
we assume that
\begin{equation}
                                         \label{4.16,4}
a^{\shharp}_{\rho_{a}}\leq 
\hat  a,\quad\hat b_{q\beta,p\beta,\rho_{b}}\leq \hat  b,
\quad N_{1}\hat b_{q_{0},p_{0},\rho_{b}}<\sfb_{0},
\end{equation}
where 
$$
\hat a=\hat  a(d,\delta,q,p,\beta)>0,\quad\hat b=\hat  
b(d,\delta,q,p,\beta,\rho_{a} )>0,
$$ 
 are taken from Theorem 2.1 of \cite{Kr_27}.

\begin{theorem}
                  \label{theorem 11.29,1}  
Suppose that
$$
 \widehat {Da}_{p_{0},\rho_{0}}\leq 
e^{-1 }\widehat{Da},\quad \hat b_{p_{0} , \rho_{0} }
\leq e^{-1 }\hat b,
$$
where $\widehat{Da},\hat b$ are from Lemma
\ref{lemma 11.26,3} with $n=1$ there.
Then, as we know from Theorem \ref{theorem 12.6.02} and 
  Remark \ref{remark 11.28,1},
there exists a strong Markov, strong Feller 
 regular
diffusion process $X$ corresponding to
  $a,b$, for which  
estimate \eqref{11.22,4} holds
and its finite-dimensional distributions
are completely determined by $a,b$ (Theorem 
\ref{theorem 12.12.3}). For $s\geq t$ introduce
$$
T_{t,s}f(x)=E_{t,x}f(x_{s-t}).
$$

Then there are constants $N$ depending only 
on $d,\delta,q_{0},p_{0},\rho_{0},\beta$ such that
with $\lambda =\lambda_{0}\rho_{0}^{-2}$ from
Theorem \ref{theorem 5.8,20}

(i) For any $f\in E_{p,\beta }$  
we have
$$
 \| T_{t,s}f\|_{E_{p, \beta} }
\leq N \|f\|_{E_{p,\beta}}.
$$  

(ii) For any $f\in E_{p,\beta }$
and $s_{0}<s$, $e^{ \lambda\cdot}T_{\cdot,s}f\in E^{1,2}_{q,p,\beta}((-\infty,s_{0})\times\bR^{d })$ and
\begin{equation}
                    \label{5,30.1}
\|e^{\lambda  \cdot} T_{\cdot,s}f
\|_{E^{1,2}_{q,p,\beta}((-\infty,s_{0})\times\bR^{d })}\leq
N(s-s_{0}) ^{-1}e^{\lambda s}\|f\|_{E_{p,\beta }}.
\end{equation}
 
(iii) For any $f\in E_{p,\beta }$
and $s-t\in (0,1]$,  
\begin{equation}
                       \label{6,15.3}
 |T_{t,s}f|\leq N(s-t)^{-\beta/2}
\|f\|_{E_{p,\beta }},\quad |T_{t,s}f|\leq N
 \sup_{B\in\bB_{\sqrt{ s-t}}}\dashnorm f\|_{L_{p}(B)},
\end{equation}
\begin{equation}
                       \label{7,7.6}
 \|DT_{t,s}f\|_{E_{p,\beta }}\leq N(s-t)^{-(q+2)/(2q)}
\|f\|_{E_{p,\beta }}.
\end{equation}
Furthermore, $T_{t,s}f(x)$ is a continuous function of $(t,x)$ for $s>t$.

(iv) For any $f\in E_{p,\beta }$  
and $t>0$, (all derivatives are Sobolev derivatives)
\begin{equation}
                     \label{6,13.5}
\partial_{t}T_{t,s}f=
 \cL T_{t,s}f.
\end{equation}

\end{theorem}
 
Proof. First suppose that $a,b$ are infinitely differentiable
with each derivative bounded. Let $f\in C^{\infty}_{0}(\bR^{d})$. Then It\^o's formula shows that
the solution $u$ from Lemma \ref{lemma 11.26,3},
in which we take $\lambda=\rho_{0}^{-1}$, $n=1$, and $T=2$,
admits the representation $u(0,x)=T_{0,2}f(x)$.
Hence,
$$
\int_{\bR^{d}}(T_{0,2}f(x))^{2}e^{-\lambda|x|}\,dx
\leq N\int_{\bR^{d}}|f(x)|^{2 }e^{-\lambda|x|}\,dx.
$$
For $f\geq0$ by Harnack inequality
(see Remark \ref{remark 11.28,1}) $T_{1,2}f(0)\leq NT_{0,2}f(x)$ if $|x|\leq1$. Hence
  by H\"older's inequality
\begin{equation}
                                \label{12.1,1}
|T_{1,2}f(0)|^{p}\leq N 
\int_{\bR^{d}}|f(x)|^{p}e^{-\lambda|x|}\,dx.
\end{equation}
We obtained this for $f\geq0$. The same is true
for $f\leq0$ and then the inequality holds
for any $f\in C^{\infty}_{0}$. As a simple
consequence of this estimate we get that,
as long as $s-t=1,x\in\bR^{d}$,
$$
|T_{t,s}f(x)|^{p}\leq N
\int_{\bR^{d}}|f(x+y)|^{p}e^{-\lambda|y|}\,dy.
$$
Finally, the case that $t-s<1$ reduces to
the one with $t-s=1$ by using self-similarity
(which maps $a,b$ into the new ones  for which
our assumptions hold with the same constants)
and leads to
$$
|T_{t,s}f(x)|^{p}\leq N(s-t)^{-d/2}
\int_{\bR^{d}}|f(x+y)|^{p}e^{-\lambda|y|/\sqrt{s-t}}\,dy.
$$

It follows that for any $B\in \bB_{r}$
$$
\dashnorm T_{t,s}f\|^{p}_{L_{p}(B)}
\leq N\sup_{B'\in \bB_{r}}\dashnorm  f\|^{p}_{L_{p}(B')},
$$
$$
 \sup_{B \in \bB_{r}}\dashnorm T_{t,s}f\| _{L_{p}(B)}
\leq N\sup_{B \in \bB_{r}}\dashnorm  f\| _{L_{p}(B)}
$$
and (i) follows in our particular case if $t-s\leq1$.

Furthermore, Remark \ref{remark 12.4,1} and the inequality
$$
e^{-\lambda|y|/\sqrt{s-t}}\leq N\dashint_{B_{\sqrt{s-t}}}e^{-\lambda|y-z|/\sqrt{s-t}}\,dz
$$
  imply that
$$
|T_{t,s}f(x)|^{p}\leq N (s-t)^{-d/2}
\int_{\bR^{d}}\dashint_{B_{\sqrt{s-t}}}|f(x+y +z)|^{p}\,dz \,e^{-\lambda|y |/\sqrt{s-t}}\,dy
$$
\begin{equation}
                                 \label{12.3,1}
= N (s-t)^{-d/2}
\int_{\bR^{d}}\dashnorm f \|^{p}_{L_{p} (B_{\sqrt{s-t}}(x+y))}  \,e^{-\lambda|y |/\sqrt{s-t}}\,dy
\leq N\sup_{B\in\bB_{\sqrt{s-t}}}\dashnorm f\|^{p}_{L_{p}(B)}.
\end{equation}
This proves \eqref{6,15.3} for $s-t\leq1$.   Since for $s-t>1$ it follows from the evolution
property of $T_{t,s}$ that $|T_{t,s}f|\leq
\sup_{x}T_{s-1,s}|f|\leq N\|f\|_{E_{p,\beta }}$,
we obtain (i) for $s-t>1$ as well.

As it is clear from the beginning of the proof,
for any $s$
the function $v(t,x)=e^{\lambda t}T_{t,s}f(x)$
satisfies the equation
$$
\partial_{t}v+a^{ij}D_{ij}v+b^{i}D_{i}v-\lambda v=0
$$
for $t<s$ with the boundary condition 
$v(s,x)=e^{\lambda s}f(x)$. To prove (ii) we take
infinitely differentiable $\zeta(t)$ such that
$\zeta(t)=1$ for $t\leq s_{0}$ and $\zeta(t)=0$
for $t\geq s$. Then the function $w:=v\zeta$ satisfies
\begin{equation}
                                \label{12.3,2}
\partial_{t}w+a^{ij}D_{ij}w+b^{i}D_{i}w-\lambda w-v\partial_{t}\zeta=0
\end{equation}
in $\bR^{d+1}$ and is smooth bounded
with each derivative bounded. Hence $w\in E^{1,2}_{q,p,\beta}$. Moreover,
$$
\|v\partial_{t}\zeta\|_{E_{q,p,\beta}}
\leq N(s-s_{0})^{-1}\sup_{t\in(s_{0},s)}
\|v(t,\cdot)\|_{E_{p,\beta}}\leq 
N(s-s_{0})^{-1}\|f \|_{E_{p,\beta}}.
$$
Now \eqref{5,30.1} follows by Theorem \ref{theorem 5.8,20} and it implies \eqref{7,7.6} by Theorem
\ref{theorem 6,6,1} and Remark \ref{remark 8,25.1}. The continuity of $T_{t,s}f(x)$
follows from $\beta<2$ and only depends
on the estimate \eqref{5,30.1}. Assertion (iv)
was taken care of in the beginning of the proof.
This finishes the proof in the case
of smooth coefficients.  

In the case of general coefficients still there is a regular diffusion process $X$
corresponding to $a,b$. We use the mollified coefficients  $a^{(\varepsilon)},b^{(\varepsilon)}$ and denote $X^{\varepsilon}$
the corresponding diffusion process. Due to the uniqueness of the finite-dimensional
distributions of $X$ we have $T^{\varepsilon}_{t,s}f
\to T _{t,s}f$ as $\varepsilon\downarrow0$ for any bounded continuous $f$. Then without any trouble,
using the fact of weak convergencies of $DT^{\varepsilon}_{t,s}f, D^{2}T^{\varepsilon}_{t,s}f,
\partial_{t}T^{\varepsilon}_{t,s}f$,
for bounded continuous $f\in E_{p,\beta}$
one proves the assertions (i)-(iii). After that the Fatou lemma allows
us to extend the results to just $f\in E_{p,\beta}$ by using mollifications of $f$
and  the second estimate in \eqref{6,15.3}. Equation \eqref{6,13.5} for $f\in C^{\infty}_{0}$ is obtained as in Theorem \ref{theorem 7,3.1}. Indeed, the fact
that, say $b^{(\varepsilon)}\to b$ strongly
in $L_{2}(C)$ and $DT^{\varepsilon}_{t,s}f
\to DT _{t,s}f$ weakly in $L_{2}(C)$
for any $C\in\bC, \bar C\subset(-\infty,s)
\times \bR^{d} $ implies that the integral of
$  b^{(\varepsilon)i}D_{i}T^{\varepsilon}_{t,s}f$
over $C$ converges to the integral of
$  b^{ i}D_{i}T _{t,s}f$. This shows that
the integral of the difference of the sides of
\eqref{6,13.5} over any such $C$ is zero meaning that the difference itself is zero. Passing from
$f\in C^{\infty}_{0}$ to arbitrary is done similarly on the basis of the second estimate in \eqref{6,15.3}.
The theorem is proved. \qed  
 
\begin{remark}
                        \label{remark 12.4,1}
We derived \eqref{12.1,1} is the ``smooth''
case. Obviously in also holds for $T_{t,s}$
corresponding to $X$. After that self-similar
transformations and shifts of the origin show that for $s-t\leq 1$, $x\in \bR^{d}$
$$                                 
|T_{t,s}f(x)|^{p}\leq N(s-t)^{-d/2}
\int_{\bR^{d}}|f(x+y)|^{p}e^{-\lambda|y|/\sqrt{s-t}}\,dx
$$
with the same $N$ as in \eqref{12.1,1}.

In the future we are also going to use
that for all $t\leq s$
$$
|T_{t,s}f|\leq N[1\wedge(s-t)]^{-\beta/2}
\|f\|_{E_{p,\beta }},\quad |T_{t,s}f|\leq N
 \sup_{B\in\bB_{ 1\wedge \sqrt{ s-t}}}\dashnorm f\|_{L_{p}(B)}.
$$
This estimates follow from \eqref{7,7.6}
and the evolution property of $T_{t,s}$. 
In the last estimate  $p=p_{0}/\beta< d+2$.
Sometimes it is useful to know that it also holds for $p=d+2$ and 
even higher $p$. This follows
from the monotonicity of $\dashnorm f\|_{L_{p}(B)}$ 
in $p$ and yields that for any $r\geq p=p_{0}/\beta$
\begin{equation}
                                      \label{4.18,1}
|T_{t,s}f|\leq N(1\wedge\sqrt{s-t})^{-d/r}
\sup_{B\in B_{1}}\|f\|_{L_{r}(B)},
\end{equation}
where $N$ depends only on $d,\delta,
q_{0},p_{0},\rho_{0},\beta$ and $r$.

\end{remark}

It is  important to note that \eqref{6,15.3}
implies the existence of a density $p(t,x,s,y)$ 
for $t<s$ that is
$$
T_{t,s}f(x)=\int_{\bR^{d}}p(t,x,s,y)f(y)\,dy.
$$

In case $f\in \EO_{p,\beta}$, $T_{t,s}f$
possesses additional properties.

\begin{theorem}
                 \label{theorem 7,24.1}
For any $f\in \EO_{p,\beta}$ and $s_{0}<s$  

(i) the function $ e^{\lambda t}
T_{t,s}f(x)$ belongs to $\EO^{1,2}_{q,p,\beta}((-\infty,s_{0})\times\bR^{d })$, where $\lambda =\lambda_{0}\rho_{0}^{-2}$;

(ii) the functions $T_{s_{0},s}f(x)$
and  $DT_{s_{0},s}f(x)$ belong to $\EO_{p,\beta}$  and, moreover, for any $\varepsilon
\in(0,1]$
$$
\lim_{n\to\infty}\sup_{t\in(s_{0}-\varepsilon^{-1},s_{0} ]}\big(\sup_{\bR^{d}}
|(\chi_{xn}-1)T_{t,s}f|+\|(\chi_{xn}-1)D
T_{t,s}f\|_{E_{p,\beta}})=0
$$
where $\chi_{nx}$ is introduced before
\eqref{12.5,2}.
\end{theorem}

Proof. First, Remark \ref{remark 12.4,1} shows that if
$$
\lim_{|x|\to\infty}\int_{B_{1}(s,x)}|f|^{p}
\,dyds=0,
$$
and $s-t\leq 1$, then $T_{t,s}f(x)\to0$ as $|x|\to\infty$.
It follows that $T_{t,s}f\in \EO_{p,\beta}$
if $f\in \EO_{p,\beta}$.

Then in the notation from the proof of Theorem
\ref{theorem 11.29,1} the function $w=e^{\lambda t}\zeta(t)
T_{t,s}f$ satisfies \eqref{12.3,2} with
$v\partial_{t}\zeta\in \EO_{p,q,\beta}$.
Hence, $w\in \EO_{p,q,\beta}$ by Theorem
\ref{theorem 7,22.1}. This proves assertion (i).

After that assertion (ii)
follows from (i), Lemma \ref{lemma 11.15,1}, and Remark \ref{remark 8,25.1}.    \qed

In the following theorem $T_{t,s}$
are taken from Theorem \ref{theorem 11.29,1}.

 \begin{theorem}
                       \label{theorem 12.4,2}
Suppose that for each $n=1,2,...$ we are given
bounded continuous $\bS_{\delta}$-valued $a^{(n)}(t,x)$
and $\bR^{d}$-valued $b^{(n)}(t,x)$,
which satisfy the assumptions of
Theorem \ref{theorem 11.29,1}
(with fixed $\beta,p,q,\rho_{a},\rho_{b}$).
Then, as we know from Theorem \ref{theorem 12.6.02}
and Remark \ref{remark 11.28,1}, for each $n$
there exists a strong Markov, strong Feller, 
 regular
diffusion process $X^{(n)}$ corresponding to
  $a(n),b(n)$, for which
estimate \eqref{11.22,4} holds
and its finite-dimensional distributions
are completely determined by $a(n),b(n)$ (Theorem 
\ref{theorem 12.12.3}). For $s\geq t$ introduce
$$
T^{n}_{t,s}f(x)=E^{(n)}_{t,x}f(x_{s-t}).
$$
Suppose that $a^{(n)}\to a$ (a.s.)
and $b^{(n)}\to b$ in $L_{p_{0}}(C)$
for any   $C\in\bC$.  

Then  
for any $f\in \EO_{p,\beta}$, $s\in\bR$,  and $\varepsilon\in(0,1]$ we have  
\begin{equation}
                        \label{7,25.10}
\lim_{n\to\infty}\sup_{t\in(s-\varepsilon^{-1},s-\varepsilon)}\sup_{\bR^{d}} |T^{(n)}_{t,s}f-T _{t,s}f
|(x)
=0,
\end{equation}
\begin{equation}
                        \label{7,20.2}
\lim_{n\to\infty}\sup_{t\in(s-\varepsilon^{-1},
s-\varepsilon)}\|DT^{(n)}_{t,s}f-DT _{t,s}f
\|_{E _{p,\beta}(\bR^{d}) }=0.
\end{equation}

\end{theorem}

Proof. First  notice that the arguments in the
proof of Theorem \ref{theorem 3.15.1} and
the uniqueness statement in Theorem 
\ref{theorem 12.12.3} show that
$T^{(n)}_{t,s}f(x)\to T _{t,s}f
(x)$ as $n\to\infty$ for any bounded
continuous $f$, $t\leq s$, $x\in\bR^{d}$.
Since, $T^{(n)}_{t,s}f(x)$ and $ T _{t,s}f
(x)$ are caloric functions in $(-\infty,s)\times\bR^{d}$ and the processes $X^{(n)}$
and $X$ are regular, the functions
$T^{(n)}_{t,s}f(x)\to T _{t,s}f
(x)$ are uniformly continuous on any compact subset of $(-\infty,s)\times\bR^{d}$ and hence
$T^{(n)}_{t,s}f(x)\to T _{t,s}f
(x)$  uniformly on any such subset.

Then the estimate (see \eqref{12.3,1})
$$
|T^{(n)}_{t,s}f(x)|^{p}\leq N (s-t)^{-d/2}
\int_{\bR^{d}}\dashint_{B_{\sqrt{s-t}}}|f(x+y +z)|^{p}\,dz \,e^{-\lambda|y |/\sqrt{s-t}}\,dy
$$
$$
= N (s-t)^{-d/2}
\int_{\bR^{d}}\dashnorm f \|^{p}_{L_{p} (B_{\sqrt{s-t}}(x+y))}  \,e^{-\lambda|y |/\sqrt{s-t}}\,dy
$$
and the fact that $\|f \|_{L_{p} (B_{\rho})}\to0$ as $B$ escapes from any $B_{R}$,
imply that
$$
\sup_{|x|\geq R,1\geq s-t\geq \varepsilon}
(|T^{(n)}_{t,s}f(x)|+|T _{t,s}f(x)|)\to0
$$
as $|x|\to\infty$. For $2\geq s-t\geq 1$
we use that
$$
|T^{(n)}_{t,s}f-T _{t,s}f|\leq |T^{(n)}_{t,s-1+\varepsilon}
T _{s-1+\varepsilon,s}f-T _{t,s-1+\varepsilon}
T _{s-1+\varepsilon,s}f|
$$
$$
+\sup_{x}|T^{(n)} _{s-1+\varepsilon,s}f-T _{s-1+\varepsilon,s}f|(x),
$$
where $T _{s-1+\varepsilon,s}f(x)$ is a continuous function vanishing at infinity. The way to go to higher values of $s-t$ is now obvious and this proves \eqref{7,25.10}
if, additionally, $f$ is bounded and continuous.

In the case of general $f\in \EO_{p,\beta}$
use common mollifiers $f^{(\gamma)}$
which are bounded, continuous, and belong to
$\EO_{p,\beta}$ and according to Remark \ref{remark 12.4,1} are such that
$$
|T^{(n)}_{t,s}f-T^{(n)}_{t,s}f^{(\gamma)}|
\leq N
 \sup_{B\in\bB_{ 1\wedge \sqrt{ s-t}}}\dashnorm f-f^{(\gamma)}\|_{L_{p}(B)}, 
$$
where the right-hand side goes to zero
as $\gamma\downarrow0$ as long as $s-t
\geq \varepsilon$. It is clear how to use these
facts to prove \eqref{7,25.10} in the general case.

  After that \eqref{7,25.10}
and \eqref{7,20.2} follow immediately from
 Theorem \ref{theorem 11.29,1}
(ii) and Remark \ref{remark 8,25.1}.
The theorem is proved. \qed

A useful addition to \eqref{7,20.2} is the following. 

\begin{theorem} 
                       \label{theorem 10,21,2}
Under the assumptions of Theorem \ref{theorem 12.4,2} 
suppose that we are also given real-valued functions
$c( x),c_{n}( x)$ such that $|c|,|c_{n}|\leq 1$ and for any
$B\in\bB_{1}$  
$$
 \int_{B}|c_{n}-c| \,dx \to0
$$
as $n\to\infty$.  
 Then for any $f\in \EO_{p,\beta }$ 
 and $t<s$ we have  
$$
\lim_{n\to\infty} 
\|c_{n}DT^{(n)}_{t,s}f-cDT _{t,s}f
\|_{E _{p,\beta } }
=0.
$$

\end{theorem}

Proof. In light of \eqref{7,20.2}  it suffices to prove that
$$
\lim_{n\to\infty} 
\|(c_{n}-c)DT _{t,s}f 
\|_{E _{p,\beta } }
=0.
$$
Owing to Theorem \ref{theorem 7,24.1} it
suffices to prove that for each $m>0$
$$
\lim_{n\to\infty} 
\|(c_{n}-c)\chi_{xm}DT _{t,s}f 
\|_{E _{p,\beta } }
=0
$$
or that for any ball $B$
$$
\lim_{n\to\infty} 
\|(c_{n}-c) DT _{t,s}f 
\|_{E _{p,\beta }(B) }
=0,
$$  
which is easily obtained by adapting the proof
of Lemma \ref{lemma 2.16.2} and using Remark \ref{remark 12.5,3}  and  H\"older's inequality. \qed

\mysection[A criterion for strong solutions.
 Case $\scB=0$]{A criterion for strong solutions. 
 Case  $\scB=0$ }
                                        \label{section 7.3.1}  
Recall that the assumptions stated in the introduction to the
chapter and at the beginning of Section \ref{section 4.16,1}
are supposed to hold throughout this section.

The following assumptions are basically 
requiring $\widehat {Da}_{p_{0},\rho }$ and $\hat b_{p_{0} , \rho }$ to be sufficiently small.
Suppose that
\begin{equation}
                                \label{12,11,2}
\widehat {Da}_{p_{0},\rho_{0}}\leq 
 e^{-1}\widehat{Da},
\quad\hat b_{p_{0} , \rho_{0} }
\leq  e^{-1}\hat b,
\end{equation}
where $(\widehat{Da},\hat b)=
(\widehat{Da},\hat b)(d,\delta,p_{0})$, are taken from Lemma 
\ref{lemma 11.26,3} when $n=1$ there.

\begin{remark}
                      \label{remark 12.5,1}
In \eqref{12,11,2} the condition on 
$\widehat {Da}_{p_{0},\rho_{0}}$ can be expressed in terms of $D\sigma$ since $a=\sigma
\sigma^{*}$ and in \eqref{4.16,4} the condition on $a^{\shharp}_{\rho_{a}}$ can be expressed
in terms of $\widehat {Da}_{p_{0},\rho_{a}}$
since by the Poincar\'e inequality
$a^{\shharp}_{\rho_{a}}\leq N(d)\widehat {Da}_{p_{0},\rho_{a}}$. We do not do that in order
to be able to  check more easily that
we can use the previous
results. For that matter, observe that,
since $q_{0}\leq p_{0}$, we have $\hat b_{q_{0},p_{0} , \rho} 
 \leq \hat b_{p_{0},p_{0} , \rho }=\hat b_{p_{0},\rho}$.

\end{remark}

Come back to  equation \eqref{3.15,1}.  As we know (see, for instance, Lemma  
3.4.1 of \cite{Kr_25}), any solution of
\eqref{11.29.20} is also a solution 
of \eqref{3.15,1} on an extended probability space (with a different Wiener process) and vice-versa
any solution of \eqref{3.15,1}
 is also a solution 
of \eqref{11.29.20} (on the same probability
space with a different $w_{t}$).
Therefore, from the point of view of
weak solutions there is no difference
which equation to consider. However,
there is a very big difference when it
comes to the strong solutions 
(cf.~Remark \ref{remark 7.10.1}) and this is the reason why in this section we consider
more general equation~\eqref{3.15,1}.

Fix $x_{0}\in\bR^{d}$ and let $(\Omega,\cF,P)$ be a complete probability space with increasing family of complete $\sigma$-fields $\cF_{t}\subset \cF$, $t\in[0,\infty)$. Assume that on this probability space there is a $d_{1}$-dimensional Wiener process $w_{t}$,
which is a Wiener process relative 
to $\{\cF_{t}\}$. Furthermore, assume that
on this probability space with the given Wiener process   equation \eqref{3.15,1} has an
$E_{q ,p,\beta}$-admissible  solution $x_{t}$
(cf. Corollary \ref{corollary 12.15,3}
or Theorem \ref{theorem 12.16,6}).
Theorem  \ref{theorem 3.15.1} 
implies that the objects described  above 
do exist and by the weak uniqueness  Theorem \ref{theorem 12.12.3} (or \ref{theorem 12.16,6}) the finite-dimensional 
distributions of all $E_{q,p,\beta}$-admissible solutions starting from the same point $x_{0}$ 
are the same. In particular,
$Ef(x_{t})=T_{0,t}f(x_{0})$ for any Borel $f\geq0$, and also all other finite-dimensional distributions of $x_{t}$
are the same as for the process $X$
 from Theorem \ref{theorem 12.6.02} under measure  $P_{0,x_{0}}$ with $T_{t,s}$
possessing the properties listed
in Theorem~\ref{theorem 11.29,1}.

The goal of this section is to give an analytical (in pure PDEs terms) criterion
for the solutions of \eqref{3.15,1} to be strong. We follow \cite{VK_76}. 
\begin{theorem} 
                                             \label{theorem 7,27.1}
Let   $t_{0}>0$,
$f\in E_{p,\beta}$, and assume
that $Ef^{2}(x_{t_{0}})<\infty$. Then

(i) With probability one   we
have
\begin{equation}
                 \label{7,26.3}
 f(x_{t_{0}})=   T_{0,t_{0}}f(x_{0})+\int_{0}^{t_{0}}\sigma^{ik}D_{i}
   T_{t_{1},t_{0}}f(x_{t_{1}})\,dw^{k}_{t_{1}},
\end{equation}
where $\sigma^{ik}D_{i}
   T_{t,s}f(x )= \sum_{i}\sigma^{ik}(t,x)D_{i}
   T_{t,s}f (x )$ and similar notation is also used below;

(ii) We have
\begin{equation}
                                                     \label{7,26.4}
   T_{0,t_{0}}f^{2}(x_{0})=(   T_{0,t_{0}}f(x_{0}))^{2}+
\sum_{k}\int_{0}^{t_{0}}   T_{0,t_{1}}\Big[\Big( 
\sigma^{ik}D_{i}   T_{t_{1},t_{0}}f\Big)^{2}\Big](x_{0})\,dt_{1}.
\end{equation}

\end{theorem}

Proof.  First, observe that $T_{0,t_{0}}f^{2}(x_{0})<\infty$.
 As we have noted,
for any $\zeta\in C^{\infty}_{0}(\bR)$, such that it is zero for $t<-1$ and $t>t'$ with $t'<t_{0}$, we have  $\zeta T_{\cdot,t_{0}}f\in E^{1,2}_{q,p,\beta }$. Therefore, we can apply
It\^o's formula (Lemma \ref{lemma 3.16.1})
to $\zeta( t)T_{t,t_{0}}f(x_{t})$ and 
using \eqref{6,13.5} and the arbitrariness of $\zeta$ write that for any $t'<t_{0}$
(a.s.)
\begin{equation}
                       \label{7,26.6}
 T_{t',t_{0} }f(x_{t'})
=T_{0,t_{0}}f(x_{0})
+\int_{0}^{t'} \sigma^{ik}D_{i}
   T_{t_{1},t_{0} }f(x_{t_{1}})\,dw^{k}_{t_{1}}.
\end{equation}

Here owing to \eqref{3.15.3}, embedding
theorem in Corollary \ref{corollary 10.8.2}, 
and \eqref{5,30.1} for certain $\lambda>0$
$$
E\int_{0}^{t'}|DT_{t_{1},t_{0} }f(x_{t_{1}})|^{2}\,dt_{1}\leq N\|e^{\lambda\cdot}(DT_{\cdot, t_{0}  }f)^{2}
\|_{E_{q,p,\beta}((0,t')\times\bR^{d}))}
$$
$$
\leq N\|e^{\lambda\cdot} T_{\cdot, t_{0} }f 
\|^{2}_{E^{1,2}_{q,p,\beta}((-\infty,t')\times\bR^{d}))}\leq N\|f\|^{2}_{E_{p,\beta}}<\infty.
$$

It follows from \eqref{7,26.6} that
$$
\sum_{k}\int_{0}^{t'}   T_{0,t_{1}}\Big[\Big(\sum_{i}
\sigma^{ik}D_{i}   T_{t_{1},t_{0}}f\Big)^{2}\Big](x_{0})\,dt_{1}=E\big|
T_{t',t_{0}}f(x_{t'})
-T_{0,t_{0}}f(x_{0})|^{2}
$$
$$
=E|T_{t',t_{0}}f(x_{t'})|^{2}-\big(
T_{0,t_{0}}f(x_{0})\big)^{2}\leq 
E T_{t',t_{0}}f^{2}(x_{t'})-\big(
T_{0,t_{0}}f(x_{0})\big)^{2}
$$
$$
=T_{0,t_{0}}f^{2}(x_{0})-\big(
T_{0,t_{0}}f(x_{0})\big)^{2}.
$$
Since the last expression is finite
and independent of $t'$, the stochastic
integral in \eqref{7,26.3} is a square-integrable martingale on $[0,t_{0}]$.

Then using the notation associated
with the process $X$
$$
E\big  |T_{t',t_{0} }f(x_{t'})-f(x_{t_{0}})
\big|=E_{x_{0}}|T_{t',t_{0} }f(x_{t'})-f(x_{t_{0}})
\big|
$$
$$
=E_{x_{0}}|E_{x_{0}}\big(f(x_{t_{0}})
\mid \cN_{t'}\big)-f(x_{t_{0}})
\big|=:I_{t'}.
$$
Since $E|f(x_{t_{0}})|<\infty$, by L{\'e}vy's theorem, as $t'\uparrow t_{0}$,
$$
I_{t'}\to E_{x_{0}}|E_{x_{0}}\big(f(x_{t_{0}})
\mid \cN_{t_{0}-}\big)-f(x_{t_{0}})
\big|,
$$
which is zero due to the continuity 
of $x_{t}$. This and \eqref{7,26.6}
prove \eqref{7,26.3}. Formula
\eqref{7,26.4} is obtained by taking the
expectations of squares of the sides of
\eqref{7,26.3}. The theorem is proved.
\qed

For further discussion we need the following result about measurable versions of stochastic integrals
whose integrand depends on a parameter
in measurable way. Denote by  $\cP$ 
the $\sigma$-field of predictable sets 
\index{$B$@Sets!$\cP$}
\index{$B$@Sets!$\cB(0,\infty)$}
and by $\cB(0,\infty)$  the Borel $\sigma$-field in $(0,\infty)$.

\begin{lemma}
                                              \label{lemma 6.17.1}
 Assume that
for $s,r\in (0,\infty)$, $\omega\in\Omega$ we are given
a real-valued function $g(s,r)=g(s,r,\omega)$, $s\in(0,\infty)$, $(r,\omega)
\in(0,\infty)\times\Omega$ which is measurable in $(s,r,\omega)$
with respect to  $\cB(0,\infty)
\otimes \cP$ and such that
 for each $s$
$$
E\int_{0}^{\infty}g^{2}(s,r)\,dr<\infty.
$$
Then there is a function $m_{s,t}=m(s,t,\omega)$ on $[0,\infty)\times\big(
[0,\infty) \times\Omega\big)$
measurable with respect to $\cB[0,\infty)\otimes \cP$,
continuous in $t$ for each $(s,\omega)$ and such that
  for each $s$ it is a martingale starting from zero and, moreover,
for each $s$
(a.s.) for all $t\geq0$
\begin{equation}
                                                     \label{6.17.4}
m_{s,t}=\int_{0}^{t}g(s,r)\,dw_{r}.
\end{equation}
\end{lemma}

Proof. Introduce
$$
\Omega_{s}=\{\omega:\int_{0}^{\infty}g^{2}(s,r)\,dr<\infty\},\quad 
\hat g(s,r)=I_{\Omega_{s}}g(s,r),
$$
$$
  B_{t}(s)=\int_{0}^{t}\hat{g}^{2}(s,r)\,dr.
$$
Observe that  
 $P(\Omega_{s})=1$ so that $\Omega_{s}\in\cF_{0}$.
Also $B_{\infty}(s)<\infty$ for any $s$ and~$\omega$.

By Lemma 2.6 of \cite{Kr_11} there   exists  
a function $m_{s,t}$ on $[0,\infty)^{2}\times\Omega$ with the
properties described in the statement of the lemma
but satisfying \eqref{6.17.4} with $\hat g$ in place of $g$.
Since $P(\Omega_{s})=1$ the integrals of $\hat g$ and $g$
coincide with probability one and the lemma is proved. \qed

\begin{remark}
                  \label{remark 6.17.1}
In light of \eqref{7,26.4}, for 
$k_{1}=1,...,d_{1}$ and almost
any $t_{1}\in(0,t_{0})$ we have
$$
m_{t_{1}}^{k_{1}}:= T_{0,t_{1}}\Big[\Big( 
\sigma^{i_{1}k_{1}}D_{i_{1}}   T_{t_{1},t_{0}}f\Big)^{2}\Big](x_{0})<\infty.
$$
Furthermore, $\sigma^{i_{1}k_{1}}D_{i_{1}}   T_{t_{1},t_{0}}f
\in E_{p,\beta}$ by Theorem \ref{theorem 11.29,1}
since $f\in E_{p,\beta}$.
It follows from Theorem \ref{theorem 7,27.1} that for those $t_{1}$ (a.s.)
$$
\sigma^{i_{1}k_{1}}D_{i_{1}}
   T_{t_{1},t_{0}}f(x_{t_{1}})=T_{0,t_{1}}(\sigma^{i_{1}k_{1}}D_{i_{1}}
   T_{t_{1},t_{0}}f)(x_{0})
$$
\begin{equation}
                                                     \label{6.17.20}
+\int_{0}^{t_{1}}\sigma^{i_{2}k_{2}}D_{i_{2}}T_{t_{2},t_{1}}\big(
\sigma^{i_{1}k_{1}}D_{i_{1}}
   T_{t_{1},t_{0}}f\big)(x_{t_{2}})\,dw^{k_{2}}_{t_{2}}.
\end{equation}
After that we want to substitute the result into \eqref{7,26.3}
to get
$$
f(x_{t_{0}})=T_{0,t_{0}}f(x_{0}) 
+\int_{0}^{t_{0}}T_{0,t_{1}}(\sigma^{i_{1}k_{1}}D_{i_{1}}
   T_{t_{1},t_{0}}f)(x_{0})\,dw^{k_{1}}_{t_{1}}
$$
\begin{equation}
                                                     \label{6.17.3}
+\int_{0}^{t_{0}}\Big(\int_{0}^{t_{1}}\sigma^{i_{2}k_{2}}D_{i_{2}}T_{t_{2},t_{1}}\big(
\sigma^{i_{1}k_{1}}D_{i_{1}}
   T_{t_{1},t_{0} }f\big)(x_{t_{2}})\,dw^{k_{2}}_{t_{2}}\Big)dw^{k_{1}}_{t_{1}}.
\end{equation}
The formal objection to do that is that we should know
that the stochastic integral in \eqref{6.17.20} is, for instance,
predictable  as a function of $(\omega,t_{1})$ and this may not
happen if we allow any version of the stochastic integral
to be taken for each $t_{1}$. However,  set 
$$
h^{k}(s,x)=
I_{s<t_{0},m^{k}_{s}<\infty}\sigma^{ik}D_{i}
   T_{s,t_{0}}f(x).
$$
It is not hard to see that $h^{k}(s,x)$
is a Borel function on $[0,t_{0}]
\times\bR^{d}$. 
Then observe that 
\begin{equation}
                                                     \label{6.17.5}
I^{k}(s,u):=\int_{0}^{u}I_{r<s }\sigma^{jm}D_{j}T_{r,s}
h^{k}(s,\cdot)(x_{r})\,dw^{m}_{r}
\end{equation}
  is the sum over $ m$ of stochastic integrals and 
$$
E \int_{0}^{\infty}I_{r<s }
\big| \sigma^{jm} D_{j}T_{r,s }h^{k}(s,\cdot)(x_{r})\big|^{2}\,dr
$$
$$
=E \int_{0}^{s} 
\big| \sigma^{jm} D_{j}T_{r,s }h^{k}(s,\cdot)(x_{r})\big|^{2}\,dr
\leq T_{0,s}\Big( \big(
h^{k}(s,\cdot)\big)^{2}\Big)(x_{0}),
$$
where the inequality is due to \eqref{7,26.4}. Also the last term is finite
(is zero if $m^{k}_{s}=\infty$).

It follows from Lemma \ref{lemma 6.17.1} that $I^{k}(s,u)=I^{k}(s,u,\omega)$
has a version which we denote again $I^{k}(s,u)$,
that is continuous in $u$ for each $s,\omega$ and
measurable with respect to $\cB (0,\infty) \otimes \cP $.
Then $I^{k}(s,s)$ is predictable and we take this modification
of the right-hand side of \eqref{6.17.20} in the right-hand side
of \eqref{6.17.3} thus justifying \eqref{6.17.3}. To be quite consistent
with this argument we should have 
inserted $I_{m^{k_{1}}_{t_{1}}<\infty}$
inside the stochastic integrals in
\eqref{6.17.3}, but this indicator 
equals one for almost all $t_{1}$
and we dropped it because changing
the integrands in It\^o's integrals
on sets of measure zero does not affect
the integral.

Similar argument justifies further
iterations of \eqref{6.17.3}.

\begin{remark}
            \label{remark 2.24.1}
A different, much more general, way of justification is presented in \cite{VK_76}. 
\end{remark}

  Introduce  
\begin{equation}
                                                         \label{6.26.5}
Q^{k}_{t,s }f(x)=\sigma^{ik}(t,x)D_{i}T_{t,s}f(x).
\end{equation}
In this 
\index{$C$@Operators!$Q^{k}_{t }$}%
notation \eqref{7,26.3} and \eqref{6.17.3} become, respectively,
$$
 f(x_{ t_{0}})=  
T_{0, t_{0}}f(x_{0})+\int_{0}^{t_{0}}Q^{k_{1}}_{t_{1},t_{0} }f(x_{t_{1}})\,dw^{k_{1}}_{t_{1}};
$$
$$
f(x_{t_{0}})=T_{0,t_{0}}f(x_{0}) 
+\int_{0}^{t_{0}}T_{0,t_{1}}Q^{k_{1}}_{t_{1},t_{0} }(x_{0})\,dw^{k_{1}}_{t_{1}}
$$
$$
+\int_{0}^{t_{0}}\Big(\int_{0}^{t_{1}}Q^{k_{2}}_{t_{2},t_{1} } 
Q^{k_{1}}_{t_{1},t_{0}}f (x_{t_{2}})\,dw^{k_{2}}_{t_{2}}\Big)dw^{k_{1}}_{t_{1}}.
$$
By induction we obtain that for any $n\geq1$ for any $t_{0}>0$  (a.s.) 
$$
f(x_{t_{0}})=T_{0,t_{0}}f(x_{0}) 
+\sum_{m=1}^{n}\int_{\Gamma_{m}(t_{0})}T_{0,t_{m}}
Q^{k_{m}}_{t_{m},t_{m-1} }\cdot...\cdot
Q^{k_{1}}_{t_{1},t_{0} }f(x_{0})\,dw^{k_{m}}_{t_{m}}\cdot...\cdot dw^{k_{1}}_{t_{1}}
$$
\begin{equation}
                                                                \label{6.18.1}
+\int_{\Gamma_{n+1}(t_{0})}
Q^{k_{n+1}}_{t_{n+1},t_{n}}\cdot...\cdot
Q^{k_{1}}_{t_{1},t_{0} }f(x_{t_{n+1}})\,dw^{k_{n+1}}_{t_{n+1}}\cdot...\cdot dw^{k_{1}}_{t_{1}},
\end{equation}
where 
$$
\Gamma_{m}(t_{0})=\{(t_{1},...,t_{m})
:t_{0}>t_{1}>...>t_{m}>0\}
$$
and by 
\index{$B$@Sets!$\Gamma_{m}(t_{0})$}%
the expressions like 
$$
\int_{\Gamma_{m}(t_{0})}:::\,dw^{k_{m}}_{t_{m}}\cdot...\cdot dw^{k_{1}}_{t_{1}}
$$
we mean
$$
\int_{ 0}^{t_{0}}\,dw^{k_{1}}_{t_{1}}\int_{  0}^{t_{1}}\,dw^{k_{2}}_{t_{2}}
...\int_{ 0}^{t_{m-1}}:::\,dw^{k_{m}}_{t_{m}}.
$$
By taking the expectations of the squares of the sides in
\eqref{6.18.1} we conclude that
$$
T_{0,t_{0}}f^{2}(x_{0})=\big(T_{0,t_{0}}f(x_{0})\big)^{2}
$$
$$
+
\sum_{m=1}^{n}\int_{\Gamma_{m}(t_{0})}
\sum_{k_{1},...,k_{m }}\big[T_{0,t_{m}}
Q^{k_{m}}_{t_{m},t_{m-1}}\cdot...\cdot
Q^{k_{1}}_{t_{1},t_{0}}f(x_{0})\big]^{2}\,d t_{m} \cdot...\cdot d t_{1} 
$$
\begin{equation}
                                                                \label{6.27.7}
+\int_{\Gamma_{n+1}(t_{0}) }\sum_{k_{1},...,k_{n +1}}T_{t_{n+1}}
\big[
Q^{k_{n+1}}_{t_{n +1},t_{n } }\cdot...\cdot
Q^{k_{1}}_{t_{1},t_{0} }f\big]^{2}(x_{0} )\,d t_{n+1 } \cdot...\cdot d t_{1}.
\end{equation}

In particular,
 the sequence of
$$
\int_{\Gamma_{n }(t_{0})}\sum_{k_{1},...,k_{n }}T_{t_{n}}
\big[
Q^{k_{n }}_{t_{n},t_{n-1} }\cdot...\cdot
Q^{k_{1}}_{t_{1},t_{0} }f\big]^{2}(x_{0})\,d t_{n } \cdot...\cdot d t_{1}
$$
is decreasing.

\end{remark}

\begin{remark}
                                                          \label{remark 6.27.1}
It turns out that proving {\em directly\/} that
each term on the right-hand side of \eqref{6.27.7}
is finite presents significant difficulties. However,
observe that, due to Theorem
\ref{theorem 11.29,1}, for   $f\in E_{p,\beta} $  and $q>2$ (recall $q$ from
$\beta\leq d/p+2/q$)   
we have
$$
 \big|T_{t_{m}}
Q^{k_{m}}_{t_{m},t_{m-1}}\cdot...\cdot
Q^{k_{1}}_{t_{1},t_{0}}f(x)\big|  
$$
$$
\leq\frac{N}{t_{m}^{\beta/2}(t_{m-1}-t_{m})^{\nu}
\cdot...\cdot(t_{0}-t_{1})^{\nu}}\|f\|_{E_{p,\beta}}, 
$$
where $\nu=(q+2)/(2q)<1$ and $N$  depends only on $t_{0},d$, $\delta,p,q,p_{0},\rho_{a}$, and $\rho_{b}$.
Furthermore,  
$$
\int_{\Gamma_{m}(t_{0})}\frac{1}{t_{m}^{\beta/2}(t_{m-1}-t_{m})^{\nu}
\cdot...\cdot(t_{0}-t_{1})^{\nu}}\,d t_{m} \cdot...\cdot d t_{1} <\infty.
$$
 
\end{remark}

Recall that
 $\cF^{w}_{t}$ is the completion of $\sigma(w_{s}:s\leq t)$. 
By a result of It\^o (\cite{It_51}),  
any $\xi$ with $E|\xi|^{2}<\infty$
is $\cF^{w}_{t_{0}}$-measurable iff,
for any $n=1,2,...$, $k_{1},...,k_{n}=1,...,d_{1}$, there exists
(nonrandom)
$f^{k_{1},...,k_{n}}_{n}(t_{1},...,t_{n})$ of class
$L_{2}(\Gamma_{n}(t_{0}))$ such that
$$
\xi=E\xi+\sum_{n=1}^{\infty}
\int_{\Gamma_{n}(t_{0})}f^{k_{1},...,k_{n}}_{n}(t_{1},...,t_{n})
\,dw^{k_{n}}_{t_{n}}\cdot...\cdot dw^{k_{1}}_{t_{1}}.
$$
Furthermore, $f^{k_{1},...,k_{n}}_{n}(t_{1},...,t_{n})$ are defined uniquely (as elements of the space $L_{2}(\Gamma_{n}(t_{0}))$) and
$$
E\xi^{2}=(E\xi)^{2} +
\sum_{n=1}^{\infty}\sum_{k_{1},...,k_{n}}
\int_{\Gamma_{n}(t_{0})}|f^{k_{1},...,k_{n}}_{n}(t_{1},...,t_{n})|^{2}\,dt_{n}\cdot...\cdot dt_{1}.
$$

Denote by $\gW_{n}(t_{0})$, $n=0,1,...$, the subspace
on $L_{2}(\Omega,\cF^{w}_{t_{0}},P)$ generated by constants and if, $n\geq1$, by constants and
$$
\int_{\Gamma_{m}(t_{0})}f (t_{1},...,t_{m})
\,dw^{k_{m}}_{t_{m}}\cdot...\cdot dw^{k_{1}}_{t_{1}}
$$
as $f$ runs over $L_{2}(\Gamma^{m}_{t_{0}})$,
$m=1,...,n$, and $k_{1},...,k_{m}=1,...,d_{1}$. Let
$\Pi_{n}(t_{0})$ be the projection 
\index{$C$@Operators!$\Pi_{n}(t_{0})$}%
operator 
in $L_{2}(\Omega,\cF ,P)$ on $\gW_{n}(t_{0})$. We know that, if $\xi\in L_{2}(\Omega,\cF,P)$,
then $E(\xi\mid \cF^{w}_{t_{0}})$ is the orthogonal projection of $\xi$ on 
$L_{2}(\Omega,\cF^{w}_{t_{0}},P)$. It follows that
  $\Pi_{n}(t_{0})\xi=
\Pi_{n}(t_{0})E(\xi\mid \cF^{w}_{t_{0}})$
and $\Pi_{n}(t_{0})\xi\to
 E(\xi\mid \cF^{w}_{t_{0}})$ in $L_{2}(\Omega,\cF,P)$
 as $n\to\infty$.

Since the last term in \eqref{6.18.1}
is orthogonal to $\gW_{n}(t)$,
$$
\Pi_{n}(t_{0})\xi=T_{0,t_{0}}f(x_{0}) 
+\sum_{m=1}^{n}\int_{\Gamma_{m}(t_{0})}T_{0,t_{m}}
Q^{k_{m}}_{t_{m},t_{m-1} }\cdot...\cdot
Q^{k_{1}}_{t_{1},t_{0} }f(x_{0})
\,dw^{k_{m}}_{t_{m}}\cdot...\cdot dw^{k_{1}}_{t_{1}},
$$
where $\xi=f(x_{t_{0}})$, and
 we come to the following conclusions, in which
 for $n\geq 1$   and $t_{0}>...>t_{n}>0$  we
define
\index{$C$@Operators!$Q_{s_{n},...,s_{1}}$}%
\begin{equation}
                                                                \label{6.26.1}
Q_{t_{n},...,t_{0}}f(x)=\sum_{k_{1},...,k_{n }}
\big[
Q^{k_{n}}_{t_{n},t_{n-1}}\cdot...\cdot
Q^{k_{1}}_{t_{1},t_{0} }f\big]^{2}(x ).
\end{equation}

\begin{theorem}
                 \label{theorem 6.18.1} 
Let   $t_{0}>0$,
$f\in E_{p,\beta }$, and assume
that $Ef^{2}(x_{t_{0}})<\infty$. Then
for $\xi:=f(x_{t_{0}})$
$$
f(x_{t_{0}})-\Pi_{n}(t_{0})\xi
=\int_{\Gamma_{n+1}(t_{0})}
Q^{k_{n+1}}_{t_{n+1},t_{n}}\cdot...\cdot
Q^{k_{1}}_{t_{1},t_{0} }f(x_{t_{n+1}})\,dw^{k_{n+1}}_{t_{n+1}}\cdot...\cdot dw^{k_{1}}_{t_{1}},
$$
$$
E\big|f(x_{t_{0}})-\Pi_{n}(t_{0})\xi
\big|^{2}
$$
$$
= \int_{\Gamma_{n+1}(t_{0})}T_{0,t_{n+1}}Q_{t_{n+1},...,t_{0}}f(x_{0})\,
dt_{n+1}\cdot...\cdot dt_{1},
$$
$$
E \big(f(x_{t_{0}})\mid\cF^{w}_{t_{0}}\big)=T_{0,t_{0}}f(x_{0})
$$
$$
+\sum_{m=1}^{\infty}\int_{\Gamma_{m}(t_{0})}T_{0,t_{m}}
Q^{k_{m}}_{t_{m},t_{m-1}}\cdot...\cdot
Q^{k_{1}}_{t_{1},t_{0} }f(x_{0})\,dw^{k_{m}}_{t_{m}}\cdot...\cdot dw^{k_{1}}_{t_{1}},
$$
where the series converges in the mean square sense.
\end{theorem}

\begin{theorem}
                                             \label{theorem 6.18.2}
Let   $t_{0}>0$,
$f\in E_{p,\beta}$, and assume
that $Ef^{2}(x_{t_{0}})<\infty$. Then
$f(x_{t_{0}})$ is $\cF^{w}_{t_{0}}$-measurable iff
\begin{equation}
                                                                \label{6.18.10}
\lim_{n\to\infty}  \int_{\Gamma_{n}(t_{0})}T_{0,t_{n }}
Q_{t_{n},...,t_{0}}f(x_{0})
\,d t_{n } \cdot...\cdot d t_{1}=0.
\end{equation}
Furthermore, under either of the above equivalent conditions
$$
 f(x_{t_{0}}) =T_{t_{0}}f(x_{0})
$$
\begin{equation}
                   \label{7.9.1}
+\sum_{m=1}^{\infty}\int_{\Gamma_{m}( t_{0})}T_{0,t_{m}}
Q^{k_{m}}_{t_{m-1}-t_{m}}\cdot...\cdot
Q^{k_{1}}_{t_{0}-t_{1}}f(x_{0})\,dw^{k_{m}}_{t_{m}}\cdot...\cdot dw^{k_{1}}_{t_{1}}.
\end{equation}

\end{theorem}

\begin{theorem}
                                             \label{theorem 6.18.3}
If equation \eqref{3.15,1} has two $E_{q,p,\beta}$-admissible
 solutions  which
are not indistinguishable, then it does not have any
strong $E_{q,p,\beta}$-admissible  solution. In particular, if \eqref{3.15,1} has at least
one strong $E_{q,p,\beta}$-admissible  solution, then any  $E_{q,p,\beta}$-admissible   solution  is strong and unique.

\end{theorem}

Indeed, if one of solutions is strong,  then
\eqref{6.18.10} holds, but then
for both solutions \eqref{7.9.1}
holds, so that $f(x_{t_{0}})$ is independent of which 
solution we take. The arbitrariness of $f$ and $t_{0}$ shows that the solutions are indistinguishable. Of course, this fact we have
already obtained under less restrictive conditions in Theorem \ref{theorem 2.28.1}.

\begin{remark}
                                               \label{remark 7.10.1}  
The criterion \eqref{6.18.10} is 
proved under the assumptions,
 which involve   $\sigma $ 
  and it turns out that for some choice
of  $\sigma $  \eqref{6.18.10} may hold and for another, with the same $a$, fail
to hold. Something even more peculiar
things may happen.  

To illustrate this we take $b\equiv 0$.   Then we take $d_{1}=d= 2$
and following \cite{KZ_75} set $\sigma^{1}(x)=x/|x|$, $\sigma^{2}(x)=x^{*}/|x|$,
where $x^{*}=(-x^{2},x^{1} )$
for $x\ne 0$, $\sigma^{ik}(0)=\delta^{ik}$. Then $a^{ij}(x)=\delta^{ij}$,
equation \eqref{3.15,1} has a solution for any $x_{0}$ (on a probability space), and each solution
is a Wiener process starting from $x_{0}$, thus admitting estimate \eqref{3.15.3}. For $x_{0}\ne0$
the solutions are strong because they never reach the origin, the only point where $\sigma$ is not smooth, and, hence, \eqref{6.18.10}
holds. However, for $x_{0}=0$
there are no strong solutions, because, as is easy to see,
rotation in $x^{1}x^{2}$ coordinates  by any angle brings 
any solution
 to another solution
of the same equation. Therefore, for $x_{0}=0$ equation 
\eqref{6.18.10} does not hold. 
Also observe for the future that in this example $D\sigma \in E_{p,1}$ for any $p\in(1,d)$ and $\not
\in L_{d,\loc}\cup E_{d,1 }$.

One can construct similar examples for $d\geq3$ starting
from the following   with $d=3$, $d_{1}=9$, and $\sigma^{k}$'s
that are the $k$th columns of the matrix
(cf.~\eqref{6.3.4})
$$
\frac{1}{|x|}
\begin{pmatrix}
x^{1} &  x^{2} & x^{3} & 0 &  0   & 0   & 0  & 0  & 0 \\
0 & 0 & 0 & x^{1}  & x^{2}   & x^{3}   & 0  & 0  & 0 \\
0 & 0 & 0 & 0  & 0   & 0  & x^{1}& x^{2} & x^{3} 
\end{pmatrix}
$$
with an appropriate definition of $0/0$ so that
 $a^{ij}=\delta^{ij}$. Again any solution
of \eqref{3.15,1} is a Wiener process
starting at $x_{0}$, so that no question
concerning \eqref{3.15.3} arises, and  if $x_{t}$ is a solution of \eqref{3.15,1} with $x_{0}=0$,
then $-x_{t}$ is also a solution of \eqref{3.15,1} with $x_{0}=0$, and thus
there is no strong solutions.
However, again if $x_{0}\ne0$, the solution is strong because it never reaches the origin. This have the same
implications as above concerning
equation 
\eqref{6.18.10}, that is purely analytical statement. Even in these simplest situations the problem of finding
an {\em analytical\/} proof that 
\eqref{6.18.10} holds iff $x_{0}\ne0$
(and $T_{t}$ is the heat semigroup)
seems to be   very challenging.

We are going to prove in the future
that \eqref{6.18.10} holds
under additional assumptions on $D\sigma$ by showing that
the series of what is under the limit sign converges.

\end{remark}

It is also worth noting an immediate
consequence of having strong solutions of any equation.

\begin{theorem}
                                             \label{theorem 7.1.1}
If equation \eqref{3.15,1} has a strong
solution on one probability space then it has a strong solution
on any other probability space carrying a $d_{1}$-dimensional
Wiener process.

\end{theorem}  

We also have a stability result. 
\begin{theorem}
              \label{theorem 2.28.2}
Let   $\sigma(n)$, $n=1,2,...$, be $d\times d_{1}$-matrix valued functions
on $\bR^{d+1}$ such that for any ball $B\in\bB$ 
and almost every $t$
$$
 \int_{B}|\sigma(n,t,x)-\sigma(t,x)|\,dx\to0
$$ 
 as $n\to \infty$. Suppose that  
 $a(n):=\sigma(n)\sigma^{*}(n),b(n)$ satisfy the assumptions
stated   in the introduction to the
chapter and at the beginning of Section \ref{section 4.16,1}
(with the same constants)
and $b(n)\to b$ in $L_{p}(C)$ for any  
$C\in \bC$.
 Finally, assume that we are given nonrandom $x_{0}(n)\to x_{0}$
  and on a probability space equation
 \eqref{3.15,1} has an $E_{q,p,\beta}$-admissible 
strong solution $x_{t}$ 
and  equation \eqref{3.15,1} with $\sigma(n),b(n), x_{0}(n)$ 
in place
of $\sigma,b,x_{0}$ has also an $E_{q,p,\beta}$-admissible 
strong solution $x_{t}(n)$ for each $n$ 
(on the same probability space with the same Wiener process). Then for any $T
\in(0,\infty)$, $m\geq 1$,
$$
E\sup_{t\leq T}|x_{t}-x_{t}(n)|^{m}
\to 0
$$
as $n\to \infty$.

\end{theorem}

Proof. By writing the equation for
$y_{t}(n)=x_{t}(n)-x_{0}(n)+x_{0}$ we reduce the situation to the one where
$x_{0}(n)=x_{0}$. In that case,
in light of Corollary \ref{corollary 3.14.6}, it suffices to prove that $E|x_{t}(n)- x_{t}|\to 0$   for any $t$, or that 
\begin{equation}
                      \label{2.28.3}
E|f(x_{t}(n))-f(x_{t})|^{2}\to 0
\end{equation}
for any $f\in C^{\infty}_{0}$. 

By Theorem \ref{theorem 6.18.2} for any $t$
$$
 f(x_{t}(n,x_{0})) =T_{0,t}(n)f(x_{0})
$$
\begin{equation}
                                                                \label{7.2.2}
+\sum_{m=1}^{\infty}\int_{\Gamma_{m}(t)}T_{0,t_{m}}(n)
Q^{k_{m}}_{t_{m},t_{m-1}}(n)\cdot...\cdot
Q^{k_{1}}_{t_{1},t}(n)f(x_{0})\,dw^{k_{m}}_{t_{m}}\cdot...\cdot dw^{k_{1}}_{t_{1}},
\end{equation}
where $T_{t,s}(n)$ and $Q^{k}_{t,s}(n)$ are the operators corresponding
to $\sigma (n)$, $b(n)$.   Since $Ef^{2}(x_{t}(n,x(0)))
\to Ef^{2}(x_{t} )$ (see Theorem \ref{theorem 12.4,2}), to prove \eqref{2.28.3}, it suffices to prove  
that $f(x_{t}(n,x(0)))\to f(x_{t} )$ weakly in $L_{2}(\Omega,\cF^{w}_{t},P)$.
Furthermore, according to \cite{It_51} the linear combinations
of constants and the multiple It\^o integrals of the type
$$
\int_{\Gamma_{m}(t)}\phi(t_{1},...,t_{m})\,dw^{k_{m}}_{t_{m}}\cdot...\cdot dw^{k_{1}}_{t_{1}},
$$
where $m$ is arbitrary and $\phi$ is an arbitrary bounded (nonrandom) Borel function,
are dense in $L_{2}(\Omega,\cF^{w}_{t},P)$.
Therefore, it suffices to prove that for all such $m$ and $\phi$
 $$
Ef\big(x_{t}(n,x(0))\big)\int_{\Gamma_{m}(t)}
\phi(t_{1},...,t_{m})\,dw^{k_{m}}_{t_{m}}\cdot...\cdot dw^{k_{1}}_{t_{1}}
$$ 
 $$
\to Ef(x_{t} )\int_{\Gamma_{m}(t)}
\phi(t_{1},...,t_{m})\,dw^{k_{m}}_{t_{m}}\cdot...\cdot dw^{k_{1}}_{t_{1}}.
$$
In light of Theorem \ref{theorem 6.18.1} this is equivalent to proving that
 $$
\int_{\Gamma_{m}(t)}\phi(t_{1},...,t_{m})T_{0,t_{m}}(n)
Q^{k_{m}}_{t_{m},t_{m-1} }(n)\cdot...\cdot
Q^{k_{1}}_{t_{1},t }(n)f(x_{0})\,dt_{m}\cdot...\cdot dt_{1}
$$ 
 $$
\to \int_{\Gamma_{m}(t)}\phi(t_{1},...,t_{m})T_{0,t_{m}} 
Q^{k_{m}}_{t_{m},t_{m-1}} \cdot...\cdot
Q^{k_{1}}_{t_{1},t} f(x_{0})\,dt_{m}\cdot...\cdot dt_{1}.
$$
This relation is indeed true, which follows by the dominated
convergence theorem from Theorems
\ref{theorem 7,24.1}, \ref{theorem 12.4,2}, 
and \ref{theorem 10,21,2}
and Remark \ref{remark 6.27.1}.
The theorem is proved. \qed   

 \mysection[Estimates for coefficients   
 smooth in $x$  and $\scB=0$]{Some estimates for 
 $B^{0,\infty}$ coefficients. Case  $\scB=0$}
                                   \label{section 7.3.2}

Introduce $B^{0,\infty}$ as the set
\index{$A$@Sets of functions!$B^{0,\infty}$}%
 of functions $f(t,x)$ on $\bR^{d+1}$ such that
they are Borel in $t$ and for each $t$ are infinitely differentiable with respect to $x$
with each derivative being a bounded function
on $\bR^{d+1}$.
In addition to the assumptions on $\sigma,b$
stated in the introduction to this chapter
and at the beginning of Section \ref{section 4.16,1}
here we suppose that   $\sigma,b\in B^{0,\infty}$.

Let $(\Omega,\cF,P)$ be a complete probability space,
let $\{\cF_{t}\}$ be an increasing filtration of 
$\sigma$-fields $\cF_{t}\subset \cF$, that are complete.
Let 
$w_{t}$ be a $d_{1}$-dimensional Wiener process relative to
$\{\cF_{t}\}$.  Take $x,\eta\in \bR^{d}$ ,
$t\in\bR$, use 
\index{$C$@Operators!$u_{(\eta)}$}%
the notation
$$
u_{(\eta)}(t,x)=\eta^{i}D_{i}u(t,x),
$$
and consider the following system   
\begin{equation}
                                                        \label{6.20.3}
x_{s}=x+\int_{0}^{s}\sigma (t+r,x_{r})\,dw _{r}+
\int_{0}^{s}b(t+r,x_{r})\,dr,
\end{equation}
\begin{equation}
                                                        \label{6.20.4}
\eta_{s}=\eta+\int_{0}^{s}\sigma _{(\eta_{r})}(t+r,x_{r})\,dw _{r}
+\int_{0}^{s}b_{(\eta_{r})}(t+r,x_{r})\,dr .
\end{equation}
As is well known, \eqref{6.20.3} 
 has a unique solution which we denote by $x_{s}(t,x)$.
By substituting it into \eqref{6.20.4} we see that the coefficients
of \eqref{6.20.4} grow linearly in $\eta$ and hence 
\eqref{6.20.4} also has a unique solution which we denote by
$\eta_{s}(t,x,\eta)$. By the way, observe that equation \eqref{6.20.4}
is linear with respect to $\eta_{r}$. Therefore
$\eta_{t}(x,\eta)$ is an affine function of $\eta$.
For the uniformity of notation we set $x_{s}(t,x,\eta)=x_{s}(t,x)$.
It is also well known 
  (see, for instance, Sections 2.7 and 2.8 of
\cite{Kr_77}) that, as a function of $x$ and $(x,\eta)$, the processes
$x_{s}(t,x)$ and $\eta_{s}(t,x,\eta)$
 are infinitely differentiable in an appropriate sense
(specified below),
their derivatives satisfy the equations which are obtained by formal
differentiation of \eqref{6.20.3} and \eqref{6.20.4},
respectively, and, for any $n\geq0,T\in(0,\infty)$,  $l_{k},\xi_{k}\in\bR^{d}$,
$k=1,...,n$ (if $n\geq1$),
$x,\eta\in \bR^{d}$, $t\in\bR$, and $q\geq 1$,
\begin{equation}
                                                        \label{6.21.1}
 E\sup_{s\leq T}\Big|\Big(\prod_{k=1}^{n}(lb)D _{(l_{k},\xi_{k})}
\Big)(x_{s},\eta_{s})(t,x,\eta)
\Big|^{q}\leq N(1+|\eta|^{m}),
\end{equation}
where $N $ is a certain constant  independent of $(x,\eta)$,
$m=m(n,q)$, and,
for instance,
by $(lb)D _{(l,\xi)} \eta_{s} (t,x,\eta)$ we mean a process $\zeta_{s}$
such that, for any $q\geq1$ and $S\in(0,\infty)$
$$
\lim_{\varepsilon\downarrow0}
E\sup_{s\leq S}\big|\zeta_{s}-\varepsilon^{-1}
\big(\eta_{s} (t,x+\varepsilon l,\eta+\varepsilon\xi)
-\eta_{s} (t,x,\eta)\big)\big|^{q}=0.
$$

\begin{lemma} 
                                                     \label{lemma 6.21.1}
 Let $\eta\in\bR^{d}$   and
 $\xi_{s}(t,x,\eta)=(lb)D_{\eta}x_{s}(t,x)$. Then

(i) $\xi_{s}(t,x,\eta)$ satisfies \eqref{6.20.4}, hence, coincides with $\eta_{s}(t,x,\eta)$ for every $(t,x,\eta)$
with probability one for all $s$.

(ii) If $f(x)$ is infinitely differentiable with bounded derivatives,
then
\begin{equation}
                                                        \label{6.21.5}
Ef_{(\xi_{s}(t,x,\eta))}(x_{s}(t,x))\Big(=
E\big(f_{(\xi_{s}(t,x,\eta))}\big)(x_{s}(t,x))\Big)
=\big(Ef(x_{s}(t,x))\big)_{(\eta)}.
\end{equation}
\end{lemma}

Proof. Assertion (i) is alluded to above and is well known (see, for instance, \cite{Kr_77}).
  Assertion (ii) follows from (i)
and the fact that (see, for instance, \cite{Kr_77})
$$
\big(Ef(x_{s}(t,x))\big)_{(\eta)}=Ef_{(\xi_{s}(t,x,\eta))}(x_{s}(t,x)).
$$
The lemma is proved.\qed

Here is a more general result.
\begin{lemma}
                                                     \label{lemma 6.21.10}
Let $f(x,\eta)$ be infinitely differentiable and such that
each of its derivatives grows  as $|x|+|\eta|\to\infty$
not faster than polynomially. Let $T\in\bR$. Then

(i) for $t\leq T$, the function
$u(t,x,\eta):=Ef\big((x_{T-t},\eta_{T-t})(t,x,\eta)\big)$
is infinitely differentiable in $(x,\eta)$ and 
each of its derivatives   by absolute
value is bounded on each finite   interval in
$(-\infty,T]$
by a constant times $(1+|x|+|\eta|)^{m}$
for some $m$,

(ii) for each $x,\eta$ the function $u(t,x,\eta)$ is Lipschitz continuous with respect to $t\in[0,T]$, 

(iii) in $(0,T)\times \bR^{2d}$  (a.e.)
$\partial_{t}u(t,x,\eta)$ exists and
$$
0= \partial_{t}u(t,x,\eta)+ (1/2)\sigma^{ik}\sigma^{jk}(t,x)u_{x^{i}x^{j}} (t,x,\eta)
+\sigma^{ik}\sigma_{(\eta)}^{jk}(t,x)u_{x^{i}\eta^{j}} (t,x,\eta)
$$
$$
+(1/2)\sigma_{(\eta)}^{ik} \sigma_{(\eta)}^{jk}(t,x)u_{\eta^{i}\eta^{j}}(t,x,\eta)
+b^{i}(t,x)u_{x^{i}} (t,x,\eta)+b^{i}_{(\eta)}(t,x)u_{\eta^{i}} (t,x,\eta)
$$
\begin{equation}
                          \label{6.21.3}
=:\partial_{t}u(t,x,\eta)+\check \cL(t,x,\eta)u(t,x,\eta).
\end{equation}
\end{lemma}

Proof. Assertion (i) is known from above.
To prove the rest, first suppose that $\sigma,b $
are infinitely differentiable in both $t$ and
$x$ with each derivative being bounded. In that case the result follows directly from Theorem
2.9.10 of \cite{Kr_77}. In the general case
take a $\zeta\in C^{\infty}_{0}(\bR)$
with unit integral
and for $\varepsilon>0$ introduce
$\zeta_{\varepsilon}(t)=\varepsilon^{-1}\zeta
(t/\varepsilon)$, $(\sigma^{\varepsilon},b
^{\varepsilon})(t,x)=(\sigma,b)(t,x)*\zeta_{\varepsilon}(t)$, where the convolution is performed with respect to $t$. Denote
by $x^{\varepsilon}_{t},\eta^{\varepsilon}_{t}$
the corresponding processes and set
$$
u^{\varepsilon}(t,x,\eta):=Ef\big((x^{\varepsilon}_{T-t},\eta^{\varepsilon}_{T-t})(t,x,\eta)\big).
$$
Since the assertions of the lemma are true for
$u^{\varepsilon}$ its derivative in $x,\eta$
admit the stated estimates (independent of $\varepsilon$) and then equation \eqref{6.21.3}
provides uniform in $\varepsilon$ estimates
of $\partial_{t}
u^{\varepsilon}(t,x,\eta)$.
 By Theorem 2.8.1
of \cite{Kr_77} $x^{\varepsilon}_{t},\eta^{\varepsilon}_{t}
\to
x _{t},\eta _{t}$, as $\varepsilon\downarrow0$, in such a sense that
$u^{\varepsilon}(t,x,\eta):=Ef\big((x^{\varepsilon}_{T-t},\eta^{\varepsilon}_{T-t})(t,x,\eta)\big)\to u^{\varepsilon}(t,x,\eta)$
at any point in $(-\infty,T]\times \bR^{2d}$.
By the results in \cite{Kr_77} also the derivatives in $\eta,x$ of $u^{\varepsilon}(t,x,\eta)$ converge to the corresponding
derivatives of $u (t,x,\eta)$.
By adding to this that, as is well known
$\sigma^{\varepsilon},b^{\varepsilon}$ and
their derivatives in $x$ converge to 
$\sigma ,b $ and their corresponding derivatives for every $x$ and almost any $t$,
we find in $(0,T)\times \bR^{2d}$  (a.e.) that
$$
\lim_{\varepsilon\downarrow0}\partial_{t}
u^{\varepsilon}(t,x,\eta)
=\check \cL(t,x,\eta)u(t,x,\eta).
$$
This easily proves (iii) and the lemma. \qed

Now comes one of the most important computations. The idea behind it
is the following. If we formally differentiate both parts of \eqref{7.9.1}
in the direction of $\eta$ and then take the expectations of the squares
of both sides, then we obtain an equality in \eqref{6.21.6} below, but the inequality
is more easily achieved and this is the only thing
we need. Naturally, by $T_{t,r}$ we mean
the operator acting by the formula
$$
T_{t,r}f(x)=Ef(x_{r-t}(t,x)).
$$

\begin{lemma}
                    \label{lemma 6.21.01}
 Let $x,\eta\in\bR^{d}$, $r\in\bR$,  and let $f\in C^{\infty}_{0}$.
 Then for any $t<r $
$(t_{0}=r)$ 
$$
E\big[f_{(\eta_{r-t}(t,x,\eta))}(x_{r-t}(t,x))\big]^{2}
\geq\Big[(T_{t, r}f(x))_{(\eta)}\Big]^{2}
$$
\begin{equation} 
                         \label{6.21.6}
+\sum_{n=1}^{\infty}\sum_{k_{1},...,k_{n}}
\int_{\Gamma_{n}(r-t)}\Big[\big(T_{t,t+t_{n}}Q^{k_{n}}_{t+t_{n},t+t_{n-1}}
\cdot...\cdot  Q^{k_{1}}_{t+t_{1}, r}f(x)\big)_{(\eta)}\Big]^{2}\,dt_{n}
\cdot...\cdot dt_{1}.
\end{equation}
\end{lemma}

Proof. For $t\leq r$ introduce   
$$
\check T_{t,r}u(x,\eta)=Eu\big((x_{r-t},\eta_{r-t}) (t,x,\eta)\big).
$$ 
Then,  
by using Lemma \ref{lemma 6.21.10} and
applying It\^o's formula to 
$$
\big(\check T_{ t+s,  r}u\big)\big((x_{s},\eta_{s})(t,x,\eta)\big)
$$
for smooth bounded
$u(x,\eta)$  by dropping for simplicity the arguments $t,x$ and $\eta$  
in $x_{\cdot}(t,x)$ and $\eta_{\cdot}(t,x,\eta)$, we get
$$
u(x_{r-t},\eta_{r-t})=\check T_{t, r}u(x,\eta)
$$
$$
+\int_{0}^{r-t}\Big[\sigma^{ik}(t+t_{1},x_{t_{1}})D_{x^{i}}
\check T_{t+t_{1}, r}u(x_{t_{1}},\eta_{t_{1}})$$
$$
+\sigma^{ik}_{(\eta_{t_{1}})}(t+t_{1},x_{t_{1}})D_{\eta^{i}}
\check T_{t+t_{1}, r}u(x _{t_{1} },\eta_{t_{1}} )\Big]\,dw^{k}_{t_{1}}.
$$
It follows that
$$
Eu^{2}(x_{r-t},\eta_{r-t})= \big(\check T_{t, r}u(x,\eta)\big)^{2}
$$
$$
+\sum_{k}\int_{0}^{r-t}E\Big[\sigma^{ik}(t+t_{1},x_{t_{1}})D_{x^{i}}
\check T_{t+t_{1}, r}u(x_{t_{1}},\eta_{t_{1}})
$$
\begin{equation}
                                                        \label{6.21.7}+
\sigma^{ik}_{(\eta_{t_{1}})}(t+t_{1},x_{t_{1}})D_{\eta^{i}}
\check T_{t+t_{1}, r}u(x _{t_{1}},\eta_{t_{1}})\Big]^{2}\,dt_{1}.
\end{equation}
By using Fatou's lemma, formulas like \eqref{6.21.5}, and well-known estimates of
the derivatives of solutions of It\^o's
equations with respect to initial data,
one easily carries \eqref{6.21.7},
with = replaced by $\geq$, over to smooth $u(x,\eta)$
whose derivatives have no more than polynomial growth
as $|x|+|\eta|\to\infty$. In particular, one can apply thus modified
\eqref{6.21.7} to $u(x,\eta)=f_{(\eta)}(x)$. Then,
after noting that, in light of \eqref{6.21.5}, in that case 
$$
\sigma^{ik}(t+t_{1},x)D_{x^{i}}
\check T_{t+t_{1}, r}u(x,\eta )+
 \sigma^{ik}_{(\eta )}(t+t_{1},x)D_{\eta^{i}}
\check T_{t+t_{1}, r}u(x,\eta ) 
$$
$$
=\sigma^{ik}(t+t_{1},x )D_{x^{i}}(T_{t+t_{1}, r}f(x))_{(\eta)}
+\sigma^{ik}_{(\eta )}(x )D_{\eta^{i}}(T_{t+t_{1}, r}f(x))_{(\eta)}
$$
$$
=\big(\sigma^{ik}(t+t_{1},x )D_{x^{i}} T_{t+t_{1}, r}f(x)\big)_{(\eta)}
=\big(Q^{k}_{t+t_{1}, r}f(x)\big)_{(\eta)},
$$
we obtain
$$
E\big[f_{(\eta_{r-t} )}(x_{r-t} )\big]^{2}
\geq\Big[(T_{t, r}f(x))_{(\eta)}\Big]^{2}
+\sum_{k_{1}}
\int_{0}^{r-t}E\big[(Q^{k_{1}}_{t+t_{1}, r}f)_{(\eta_{t_{1}} )}
(x_{t_{1}} )\big]^{2}\,dt_{1}.
$$
By applying this formula to $Q^{k_{1}}_{t+t_{1}, r}f$ in place of $f$
we get   
$$
E\big[f_{(\eta_{r-t})}(x_{r-t})\big]^{2}
\geq\Big[(T_{t, r}f(x))_{(\eta)}\Big]^{2}+\sum_{k_{1}}
\int_{0}^{r-t} \big[(T_{t,t+t_{1}}Q ^{k_{1}}_{t+t_{1}, r}f (x  ))_{(\eta)}\big]^{2}\,dt_{1} 
$$
$$
+\sum_{k_{1},k_{2}}\int_{0}^{r-t}dt_{1}
\int_{0}^{t_{1}}E
\big[(Q^{k_{2}}_{t+t_{2},t+t_{1}}Q^{k_{1}}_{t+t_{1}, r}f)_{(\eta_{t_{2})}}(x_{t_{2}})
\big]^{2}\,dt_{2}.
$$

 Using the induction shows that for any $n\geq1$ 
$$
E\big[f_{(\eta_{r-t} )}(x_{r-t} )\big]^{2}
\geq\Big[(T_{t, r}f(x))_{(\eta)}\Big]^{2}
$$
$$
+\sum_{m=1}^{n}\sum_{k_{1},...,k_{m}}
\int_{\Gamma^{m}_{r-t}}\Big[
I^{k_{1},...,k_{m}}(t_{1},...,t_{m})\Big]^{2}\,dt_{m}
\cdot...\cdot dt_{1}
$$
$$
+\sum_{k_{1},...,k_{n+1}}
\int_{\Gamma^{n+1}_{r-t}}E\Big[J^{k_{1},...,k_{n+1}}(t_{1},...,t_{n+1})\Big]^{2}\,dt_{n+1}
\cdot...\cdot dt_{1},
$$
where ($t_{0}=r$) 
$$
I^{k_{1},...,k_{m}}(t_{1},...,t_{m})=
\big(T_{t,t+t_{m}}Q^{k_{m}}_{t+t_{m},t+t_{m-1}}
\cdot...\cdot  Q^{k_{1}}_{t+t_{1}, r}f(x)\big)_{(\eta)},
$$
$$
J^{k_{1},...,k_{n+1}}(t_{1},...,t_{n+1})=\big( Q^{k_{n+1}}_{t+t_{n+1 },t+t_{n}}
\cdot...\cdot  Q^{k_{1}}_{t+t_{1}, r}f\big)_{(\eta_{t_{n+1}})}
(x_{t_{n+1}})
$$
This yields \eqref{6.21.6} and proves the lemma. \qed

Next, we want to estimate the left-hand side of \eqref{6.21.6}
which according to Lemma \ref{lemma 6.21.10} satisfies
\eqref{6.21.3}. 

In the future we will be interested in estimating not only
 the left-hand side of \eqref{6.21.6} but a slightly more
general quantity. Therefore, we take a  
nonnegative function $f(x,\eta)$, which
is a {\em polynomial\/} with respect to $\eta$
with coefficients depending on $x$
such that $f(\cdot,\eta)\in C^{\infty}_{0}$
for any $\eta$.
Then for $t_{0}\in(0,\infty)$ and $t\leq t_{0}$ denote 
$$
u(t,x,\eta)=\check T_{t,t_{0}}f(x,\eta).  
$$ 
According to Lemma \ref{lemma 6.21.10}
the function $u(t,x,\eta)$ salsifies
\eqref{6.21.3} and, since $\eta_{T-t}(t,x,\eta)$
is affine in $\eta$, $u(t,x,\eta)$ is a polynomial
in $\eta$.
 
 Introduce
$$
 \widehat {D\sigma }_{s,\rho}:=
\sup_{r\leq \rho}r\sup_{C\in \bC_{r}}
\dashnorm D\sigma\|_{L_{s}(C)}.
$$
\begin{theorem}
                   \label{theorem 6.21.1}
 
Let $n\in\{ 1,2,...\}, \lambda\geq 0$.
Then there
are constants $\widehat{D\sigma},\hat b\in(0,1)$,
depending only on $d,\delta,p_{0}$,  
$n $,  and the power of the
\index{$S$@Miscelenea!$\widehat{D\sigma}_{s,\rho}$}%
\index{$S$@Miscelenea!$\widehat{D\sigma}$}%
\index{$S$@Miscelenea!$\bar b_{R}$@$\hat b$}%
 polynomial $f(x,\eta)$,
such that if 
\begin{equation}   
                             \label{12.10,1}
 \widehat {D\sigma }_{p_{0},\rho_{0}}\leq 
e^{- \lambda\rho_{0}}\widehat{D\sigma},\quad\hat b_{p_{0} , \rho_{0} }
\leq e^{- \lambda\rho_{0}}\hat b,
\end{equation}
 then
\begin{equation}
                              \label{12.20,5}
\int_{\bR^{d}}\sup_{|\eta|\leq 1}|u(0,x,\eta)|^{2n}e^{-\lambda|x|}\,dx
\leq Ne^{ \alpha t_{0}}\int_{\bR^{d}}\sup_{|\eta|\leq 1}|f(x,\eta)|^{2n}e^{-\lambda|x|}\,dx,
\end{equation}
where 
$$
\alpha=N\rho_{0}^{-2}e^{ \lambda\rho_{0}}
$$
and the constants called $N$ depend  only on $d,\delta,p_{0}$,  
$n $ and the power of the polynomial $f(x,\eta)$.

\end{theorem}

By taking $\lambda=0$ and using the arbitrariness of $\rho_{0}$
we come to the following.
\begin{corollary}
                             \label{corollary 2.6.1}
 If $\|D\sigma\|_{\dot E_{p_{0},1}}\leq \widehat{D\sigma}$
 and $\| b\|_{\dot E_{p_{0},1}}\leq \hat b$, where
 $\widehat{D\sigma},\hat b$ are taken from Theorem \ref{theorem 6.21.1},
 then
 $$
 \int_{\bR^{d}}\sup_{|\eta|\leq 1}|u(0,x,\eta)|^{2n} \,dx
\leq Ne^{N t_{0}}\int_{\bR^{d}}\sup_{|\eta|\leq 1}|f(x,\eta)|^{2n} \,dx,
$$
where the constants called $N$ depend  only on $d,\delta,p_{0}$,  
$n $, $\rho_{0}$, and the power of the polynomial $f(x,\eta)$.

\end{corollary}  

The proof of Theorem \ref{theorem 6.21.1}  is rather long
and we present  it in a separate section.
One of the main objectives of Theorem 
\ref{theorem 6.21.1} is to provide the possibility
to derive from it the following 
fundamental result.

\begin{theorem}
                      \label{theorem 12.10,1}
Suppose that \eqref{12.10,1} is satisfied
with an integer $n >d/4$ and $\lambda=0 $. Let $f\in C^{\infty}_{0}$. Then 
\begin{equation}
                                \label{12.10,5}
 \sum_{n=1}^{\infty}\int_{\Gamma_{n}(t_{0})}T_{0,t_{n }}
Q_{t_{n},...,t_{0}}f(x_{0})\,d t_{n } \cdot...\cdot d t_{1}\leq
N\Big(\int_{\bR^{d}}|Df|^{4n}\,dx\Big)^{1/(2n)},
\end{equation}
where $N$ depends only on $d,\delta,p_{0},\rho_{0},n,t_{0}$.
\end{theorem}

Proof. Introduce $\bR^{d}$-valued functions $\sigma^{k}(t,x)=(\sigma^{ik}(t,x))$, $k=1,...,d_{1}$, and observe that  for $t_{n+1}<t_{n}$
\begin{equation}
                        \label{12.10,2}
 \big(T_{t_{n+1},t_{n+1}+t_{n}}Q^{k_{n}}_{t_{n+1}+t_{n},t_{n+1}+t_{n-1}}
\cdot...\cdot  Q^{k_{1}}_{t_{n+1}+t_{1}, t_{0}}f(x)\big)_{(\sigma^{k }(t_{n+1},x))} 
\end{equation}
$$
=Q^{k}_{t_{n+1},t_{n+1}+t_{n}}Q^{k_{n}}_{t_{n+1}+t_{n},t_{n+1}+t_{n-1}}
\cdot...\cdot  Q^{k_{1}}_{t_{n+1}+t_{1}, t_{0}}f(x) .
$$
Therefore for    
$$
u(t,x,\eta):=E|f_{(\eta_{t_{0}-t}(t,x,\eta)}(
x_{t_{0}-t}(t,x))|^{2}
$$
in light of \eqref{6.21.6} we find that
$$
\sum_{k}\int_{0}^{t_{0}}T_{0,t_{n+1}}u(t_{n+1},\cdot,\sigma^{k}(t_{n+1},\cdot)) (x_{0})\,dt_{n+1}
$$
\begin{equation}
                                \label{12.10,4}
\geq \sum_{n=1}^{\infty}\int_{\Gamma_{n+1}(t_{0})}T_{0,t_{n+1}}Q_{t_{n+1},...,t_{0}}f(x_{0})
\,dt_{n+1}\cdot...\cdot dt_{1}.  
\end{equation}

By Theorem \ref{theorem 6.21.1}
$$
\int_{\bR^{d}}\sup_{|\eta|\leq 1}|u(t_{n+1},x,\eta)|^{2n} \,dx
\leq N \int_{\bR^{d}}\sup_{|\eta|\leq 1}|f_{(\eta)}(x)|^{4n} \,dx.
$$
On account of our choice of $n$ and
 Remark \ref{remark 12.4,1} the left-hand side of \eqref{12.10,4} is estimated from above
by a constant times the right-hand side of
\eqref{12.10,5}. The theorem is proved. \qed

\mysection[Proof of Theorem \protect\ref{theorem 6.21.1}]{Proof of Theorem \protect\ref{theorem 6.21.1}}
                       \label{section 1.14.1}

We need the following,  which is similar to
Lemma 5.8 of \cite{Kr_25_1}.

\begin{lemma}
                         \label{lemma 12.8.1}

Let $n\geq1$ and suppose that for 
$i=1,...,n$ we are given $p_{i} >0$,
integers $k_{i}\geq1$, and   polynomials $A_{i}(\eta)$ of degree $k_{i}$
on $\bR^{d}$. Then there exists a constant
$N=N(d, n,p_{i},\kappa_{i})$ such that
\begin{equation}
                        \label{12.18.1}
|A_{1}|^{p_{1}}\cdot...\cdot|A_{n}|^{p_{n}}\leq 
N\int_{B_{1}} | A_{1}(\eta))|^{p_{1}}\cdot...\cdot | A_{n}(\eta)|^{p_{n}}
\,d\eta,
\end{equation}
where $|A_{i}|$ is the maximum of absolute
values of the coefficients of $A$, written
without similar terms.
\end{lemma}

Proof. As it is not hard to see it suffices
to prove that for any polynomial $A(\eta)$
of degree $k$ with $|A|=1$ and any $\gamma>0$ there exists $\varepsilon>0$,
depending only on $d,k,\gamma$,
such that
$$
|B_{1}\cap\{|A(\eta)|\leq\varepsilon \}
\leq  \gamma |B_{1}|.
$$
We are going to treat $A(\eta)$ as a random
variable on the probability space $(B_{1},dx/|B_{1}|)$. Observe that the set $\frA$
of the $A(\eta)$'s is compact in $C(\bar B_{1})$,
and, since for any polynomial its any level set  has Lebesgue measure zero, the distribution functions $F_{A}$ of the $A(\eta)$'s
form a compact set $\frF$ in $C[0,1]$. It follows
that for given $\gamma$ we can find a finite $\gamma/2$-net $F_{A_{1}},...,F_{A_{m}}$ in $\frF$ and $\varepsilon>0$ such that $F_{A_{i}}(\varepsilon)\leq\gamma/2$
for any $i=1,...,m$, and then for any $A\in\frA$
we can find $F_{A_{i}}$ such that
$$
F_{A}(|A(\eta)|\leq\varepsilon)
\leq F_{A_{i}}(|A(\eta)|\leq\varepsilon)
+\gamma/2\leq\gamma.
$$
\qed
 
Now we start proving the theorem.

{\em Step 1\/}.
As in the proof of Lemma \ref{lemma 11.26,3}
we assume that $\widehat { D\sigma }_{p_{0},\rho_{0}}
\leq 1$, $\hat b_{p_{0} , \rho_{0} }
\leq 1$. Then take a $C\in\bC_{\rho_{0}}$ and 
a nonnegative $\zeta
\in C^{\infty}_{0}(C)$  with the integral
of its square equal to one,
   multiply  \eqref{6.21.3}
by $\zeta^{2} u^{2n-1}$ and integrate by parts
with respect to $(t,x)$ regarding $\eta$
as a parameter in $B_{1}$.

Then as in the proof of Lemma \ref{lemma 11.26,3} (cf. \eqref{11.28,4}) for $s\leq t_{0}$ we find   
$$
\int_{\bR^{d}}\zeta^{2}(s,x)u^{2n}(s,x,\eta)\,dx+(\delta/8)
\int_{[s,t_{0}]\times \bR^{d}} \zeta^{2} u  ^{2n-2}
|Du |^{2}\,dxdt
$$
$$
\leq N\int_{\bR^{d}}\zeta^{2}(t_{0},\cdot)f ^{2n }   \,dx
+N \int_{[s,t_{0}]\times \bR^{d}}
(|\partial_{t}\zeta^{2}|+|D\zeta |^{2}) u^{2n}\,dxdt
$$
\begin{equation}
                               \label{12.7,2}
+N\int_{[s,t_{0}]\times \bR^{d}} \zeta^{2}(|D\sigma|+|b|)^{2}u^{2n} \,dxdt+\int_{[s,t_{0}]\times \bR^{d}}\zeta^{2}u^{2n-1}F\,dxdt,
\end{equation}
where
$$
F=\sigma^{ik}\sigma_{(\eta)}^{jk}(t,x)u_{x^{i}\eta^{j}} (t,x,\eta)
$$
$$
+(1/2)\sigma_{(\eta)}^{ik} \sigma_{(\eta)}^{jk}(t,x)u_{\eta^{i}\eta^{j}}(t,x,\eta)
 +b^{i}_{(\eta)}(t,x)u_{\eta^{i}} (t,x,\eta).
$$

Observe that ($a^{n}(a^{n-1}b)\leq\varepsilon^{-1}a^{2n}+\varepsilon a^{2n-2}b^{2}$)
for any $\varepsilon>0$
$$
\int_{[s,t_{0}]\times \bR^{d}}\zeta^{2}u^{2n-1}\sigma^{ik}\sigma_{(\eta)}^{jk} u_{x^{i}\eta^{j}} \,dxdt\leq
N \int_{[s,t_{0}]\times \bR^{d}} \zeta^{2} |D\sigma| ^{2}u^{2n} \,dxdt
$$
$$
+\varepsilon\int_{[s,t_{0}]\times \bR^{d}} \zeta^{2}u^{2n-2}|u_{x\eta}|^{2} \,dxdt,
$$
where and below we allow the constants $N$
to also depend on $\varepsilon$. Also
  $$
\int_{[s,t_{0}]\times \bR^{d}}\zeta^{2}u^{2n-1}
\sigma_{(\eta)}^{ik} \sigma_{(\eta)}^{jk} u_{\eta^{i}\eta^{j}} \,dxdt
$$
$$
\leq N\int_{[s,t_{0}]\times \bR^{d}}\zeta^{2}  |D\sigma|^{2}u^{2n-1}|u_{\eta\eta}|\,dxdt,
$$
 $$
\int_{[s,t_{0}]\times \bR^{d}}\zeta^{2}u^{2n-1}
b^{i}_{(\eta)} u_{\eta^{i}}\,dxdt
$$
$$
=-(2n-1)\int_{[s,t_{0}]\times \bR^{d}}\zeta^{2}
\big(b^{i} u^{n-1}u_{\eta^{i}}\big)\big(u^{n-1}u_{x^{j}}\eta^{j}\big)\,dxdt
$$
$$
-2\int_{[s,t_{0}]\times \bR^{d}}\zeta\zeta_{x^{j}}\eta^{j}\big(u^{n}
b^{i}\big)\big(u^{n-1}  u_{\eta^{i}}\big)\,dxdt
-\int_{[s,T]\times \bR^{d}}\zeta^{2}u^{2n-1}b^{i} u_{x^{j}\eta^{i}}\eta^{j}\,dxdt
$$
$$
\leq\varepsilon\int_{[s,t_{0}]\times \bR^{d}}\zeta^{2}
u^{2n-2}(|Du|^{2}+|u_{x\eta}|^{2})\,dxdt+N\int_{[s,t_{0}]\times \bR^{d}}\zeta^{2}|b|^{2}u^{2n-2}|u_{\eta}|^{2}\,dxdt
$$
$$
+N\int_{[s,t_{0}]\times \bR^{d}}\zeta^{2}|b|^{2}u^{2n } \,dxdt+N\int_{[s,t_{0}]\times \bR^{d}}|D\zeta|^{2}u^{2n-2}|u_{\eta}|^{2}\,dxdt.
$$
We substitute these estimate into \eqref{12.7,2} and get
$$
\int_{\bR^{d}}\zeta^{2}(s,x)u^{2n}(s,x,\eta)\,dx+
\int_{[s,t_{0}]\times \bR^{d}} \zeta^{2} u  ^{2n-2}
(|Du |^{2}-\varepsilon|u_{x\eta}|^{2})\,dxdt
$$
$$
\leq N\int_{\bR^{d}}\zeta^{2}(t_{0},\cdot)f ^{2n }   \,dx
+N \int_{[s,t_{0}]\times \bR^{d}}|D\zeta |^{2} u^{2n-2}\big(u^{2}+|u_{\eta}|^{2}\big)\,dxdt
$$
$$
+N\int_{[s,t_{0}]\times \bR^{d}} \zeta^{2}(|D\sigma|+|b|)^{2}u^{2n-2}\big(u^{2 }+u |u_{\eta\eta}|+|u_{\eta}|^{2}|\big) \,dxdt . 
$$
By integrating through this inequality,
using Lemma \ref{lemma 12.8.1} and \eqref{1.21.3}
and choosing
$\varepsilon$ appropriately we finally find
$$
\int_{\bR^{d}\times B_{1}}\zeta^{2}(s,x)u^{2n}(s,x,\eta)\,dxd\eta+
\int_{[s,t_{0}]\times \bR^{d}\times B_{1}} \zeta^{2} u  ^{2n-2}
 |Du |^{2} \,dxdtd\eta
$$
$$
\leq N\int_{\bR^{d}\times B_{1}}\zeta^{2}(t_{0},\cdot)f ^{2n }   \,dxd\eta
+N \rho_{0}^{-d-4}\int_{[s,t_{0}]\times \bR^{d}\times B_{1}}I_{C} u^{2n } \,dxdtd\eta
$$
\begin{equation}
                           \label{12.7,4}
+N\int_{[s,t_{0}]\times \bR^{d}\times B_{1}} \zeta^{2}(|D\sigma|+|b|)^{2}u^{2n } \,dxdtd\eta . 
\end{equation}

{\em Step 2\/}. Here we are dealing with the last term in \eqref{12.7,4}. Introduce $w=u^{n}$ and observe that for the function
$\zeta w$ as in the proof of Lemma \ref{lemma 11.26,3} we have 
$$
\partial_{t}(\zeta w)+\Delta(\zeta w)
 +2\zeta\Big(\frac{1}{n}-1\Big)a^{ij}(D_{i}(u^{n/2}))D_{j}(u^{n/2})
+G=0,
$$
where
$$
G=-w\partial_{t}\zeta-\Delta(\zeta w) 
+(1/2)\zeta a^{ij}D_{ij}w  
+\zeta b^{i}D_{i}w 
$$
$$ 
+n\zeta u^{n-1}\big(\sigma^{ik}\sigma^{jk}_{(\eta)}u_{x^{i}\eta^{j}}+(1/2)\sigma^{ik}_{(\eta)}\sigma^{jk}_{(\eta)}u_{\eta^{i}\eta^{j}}
 +b^{i}_{(\eta)}u_{\eta^{i}}\big).
$$
Then again as in the proof of Lemma \ref{lemma 11.26,3},  defining $w$ for $t>t_{0}$ as zero, we conclude that for $t<t_{0}$
$$
0\leq\zeta w(t,x,\eta)\leq P_{2,4}G
=h(t,x,\eta) 
$$
$$
+\big(J_{1}+J_{2}+(1/2)J_{3}+J_{4}+n[J_{5}+
(1/2)J_{6}+J_{7}+J_{8}]\big)(t,x,\eta),
$$
where   
$$
h(t,x,\eta)=\hat T _{t_{0}-t}[\zeta(t_{0},\cdot)f^{n}(\cdot,\eta)](x),
$$
$$
J_{1}=-P_{2,4}(w\partial_{t}\zeta),\quad
J_{2}=-P_{2,4}(\Delta(\zeta w))=-\big(
P_{2,4}( (\zeta w)_{x^{i}})\big)_{x^{i}},
$$
$$
J_{3}=P_{2,4}\big(\zeta a^{ij}w_{x^{i}x^{j}}  \big),\quad J_{4}= P_{2,4}\big(\zeta b^{i}D_{i}w\big)
$$  
$$
J_{5}=P_{2,4}\big(\zeta u^{n-1}\sigma^{ik}\sigma^{jk}_{(\eta)}u_{x^{i}\eta^{j}}\big),
\quad
J_{6}=P_{2,4}\big(\zeta u^{n-1}\sigma^{ik}_{(\eta)}\sigma^{jk}_{(\eta)}u_{\eta^{i}\eta^{j}}\big),
$$
$$
J_{7}=P_{2,4}\big(\zeta u^{n-1}b^{i}_{(\eta)}u_{\eta^{i}}\big)=(1/n)
\eta^{k}\big(P_{2,4}\big(\zeta  b^{i} w_{\eta^{i}})\big)_{x^{k}}
$$
$$
-(1/n)P_{2,4}\big(\zeta_{(\eta)}  b^{i} w_{\eta^{i}})-(1/n)P_{2,4}\big(\zeta   b^{i} w_{\eta^{i}(\eta)}).
$$

First, by Lemma \ref{lemma 11.29,1}
$$
\int_{[s,t_{0}]\times \bR^{d}}I_{C}b^{2}h^{2}
\,dxdt\leq N\hat b^{2}_{p_{0},\rho_{0}}
 \int_{\bR^{d}}\zeta^{2}(t_{0},x)f^{2n}(x,\eta)\,dx.
$$

Next,
since $P_{2,4}=NP_{1,4}P_{1,4}$
and  $I_{C}\in\dot E_{p_{0},p_{0},1}$ 
by Remark \ref{remark 11.30,1},  
by Theorem \ref{theorem 5.25,1} and
Corollary \ref{corollary 10.5,1}
(this combination will be used repeatedly
below)
$$
\int_{\bR^{d+1}_{s}}I_{C}|b|^{2}J_{1}^{2}\,dxdt
$$
\begin{equation}
                           \label{10.11,5}
\leq N\hat b^{2}_{p_{0},\rho_{0}}\int_{\bR^{d+1}_{s}} P^{2}_{1,4}(I_{C}w|\partial_{t}\zeta|)  \,dxdt
\leq N\hat b^{2}_{p_{0},\rho_{0}}\rho_{0}^{-d-4}\int_{\bR^{d+1}_{s}} I_{C} u^{2n}  \,dxdt,
\end{equation}

Then,
$
|J_{2}|
\leq NP_{1,8}(|D(\zeta w) |)
$,
 so that
\begin{equation}
                           \label{10.11,6}
   \int_{\bR^{d+1}_{s}}I_{C}|b|^{2}J_{2}^{2} \,dxdt
\leq N\hat b^{2}_{p_{0},\rho_{0}}\int_{\bR^{d+1}_{s}} |D(\zeta w) |^{2}\,dxdt.
\end{equation}

Dealing with $J_{3}$ observe that 
$$
P_{2,4}\big(\zeta  a^{ij} w_{x^{i}x^{j}}\big)=
\big[P_{2,4}\big(\zeta  a^{ij}w_{ x^{j}}\big)\big]_{x^{i}}
$$
$$
-P_{2,4}\big(\zeta_{x^{i}} a^{ij} w_{ x^{j}}\big)
-P_{2,4}\big(\zeta  [\sigma^{ik}_{x^{i}}\sigma^{jk}+\sigma^{ik}\sigma^{jk}_{x^{i}}] w_{ x^{j}}\big) .
$$
It follows that 
$$
 \int_{\bR^{d+1}_{s}}I_{C}|b|^{2}J_{3}^{2} \,dxdt\leq N \int_{\bR^{d+1}_{s}}I_{C}|b|^{2}P_{1,8}^{2}(\zeta|Dw |) \,dxdt
$$
$$
+N\hat b^{2}_{p_{0},\rho_{0}}\int_{\bR^{d+1}_{s}}
P_{1,8}^{2}\big(I_{C}|D\zeta |\,|Dw |
+I_{C}|D\sigma |\zeta|Dw |\big) \,dxdt
$$
\begin{equation}
                           \label{10.11,7}
 \leq N\hat b^{2}_{p_{0},\rho_{0}}
\rho_{0}^{-d-2}\int_{\bR^{d+1}_{s}} I_{C} |Dw|^{2}\,dxdt,
\end{equation}
where we used that $\widehat{D\sigma}_{p_{0},\rho_{0}}\leq 1$.

Next,
$$
\int_{\bR^{d+1}_{s}}I_{C}|b|^{2}J_{4}^{2} \,dxdt
$$ 
$$
\leq N\hat b^{2}_{p_{0},\rho_{0}}
\int_{\bR^{d+1}_{s}}P_{1,8}^{2}(|b|\zeta|Dw|)\,dxdt\leq N \hat b^{2}_{p_{0},\rho_{0}}
\int_{\bR^{d+1}_{s}} \zeta^{2} u^{2n-2}|Du|^{2} \,dxdt
$$
\begin{equation}
                           \label{10.13,2}
\leq N \hat b^{2}_{p_{0},\rho_{0}}\rho_{0}^{-d-2}
\int_{\bR^{d+1}_{s}} I_{C}  u^{2n-2}|Du|^{2} \,dxdt,
\end{equation}

$$
 \int_{\bR^{d+1}_{s}}I_{C}|b|^{2}J_{5}^{2} \,dxdt\leq N\hat b^{2}_{p_{0},\rho_{0}}\int_{\bR^{d+1}_{s}}P^{2}_{1,8}\big(I_{C}|D\sigma|\zeta u^{n-1}  |u_{x\eta}|\big)\,dxdt
$$
\begin{equation}
                           \label{10.11,8}
\leq N\hat b^{2}_{p_{0},\rho_{0}}\rho_{0}^{-d-2}\int_{\bR^{d+1}_{s}}   I_{C}  u^{2n-2}  |u_{x\eta}|^{2}\,dxdt.
\end{equation}
Similarly,
$$
 \int_{\bR^{d+1}_{s}}I_{C}|b|^{2}J_{6}^{2} \,dxdt\leq N\hat b^{2}_{p_{0},\rho_{0}}\int_{\bR^{d+1}_{s}}P^{2}_{1,8}\big(I_{C}|D\sigma|(\zeta|D\sigma| | u^{n-1}  |u_{\eta\eta}|)\big)\,dxdt
$$
\begin{equation}
                           \label{10.13,1}
\leq N\hat b^{2}_{p_{0},\rho_{0}}\int_{\bR^{d+1}_{s}} \zeta^{2} |D\sigma|^{2} | u^{2n -2}   |u_{\eta\eta}|^{2} \,dxdt.
\end{equation}

Finally,
$$
\int_{\bR^{d+1}_{s}}I_{C}|b|^{2}J_{7}^{2} \,dxdt\leq N\hat b^{2}_{p_{0},\rho_{0}}\int_{\bR^{d+1}_{s}} 
\zeta^{2}|b|^{2} |w_{\eta}|^{2}\,dxdt
$$
$$
+N\hat b^{2}_{p_{0},\rho_{0}} \int_{\bR^{d+1}_{s}} \big( 
|D\zeta|^{2}  |w_{\eta}|^{2}+\zeta^{2}|w_{\eta (\eta)}|^{2}\big)\,dxdt.
$$
  
Summing up the above estimates 
  yields  
$$
 \int_{\bR^{d+1}_{s}} \zeta^{2}|b|^{2}u^{2n}\,dxdt\leq N\hat b^{2}_{p_{0},\rho_{0}}
 \int_{\bR^{d}}\zeta^{2}(T,x)f^{2n}(x,\eta)\,dx
$$
$$
+N\hat b^{2}_{p_{0},\rho_{0}}
\int_{\bR^{d+1}_{s}}\Big(\rho_{0}^{-d-4}I_{C}
(u^{2n} +|w_{\eta}|^{2})
$$
$$
+\rho_{0}^{-d-2}I_{C}\big[u^{2n-2}\big(|Du|^{2} + |u_{x\eta}|^{2}\big)+|w_{\eta(\eta)}|^{2}\big]
$$
$$
+\zeta^{2}|D\sigma|^{2}u^{2n-2}|u_{\eta\eta}|^{2}
+\zeta^{2}  |b|^{2}|w_{\eta}|^{2}  \Big)\,dxdt.
$$

{\em Step 3\/}.
Obviously, similar  estimate is valid
for
$$
\int_{\bR^{d+1}_{s}}I_{C}|D\sigma|^{2}w^{2}\,dxdt,
$$
which after adding it to the above one,
integrating over $B_{1}$ with respect to $\eta$ and using Lemma \ref{lemma 12.8.1}
yields
 $$
 \int_{\bR^{d+1}_{s}\times B_{1}}\zeta^{2}(|b|^{2}+|D\sigma|^{2})u^{2n}\,dxdtd\eta
\leq N 
 \int_{\bR^{d}\times B_{1}}\zeta^{2}(t_{0},\cdot)f^{2n} \,dxd\eta
$$
$$
+ N\Big(\hat b^{2}_{p_{0},\rho_{0}}+\widehat{D\sigma}_{p_{0},\rho_{0}}^{2}\Big)\rho_{0}^{-d-4}
\int_{\bR^{d+1}_{s}\times B_{1}}   
I_{C}u^{2n} 
\,dxdtd\eta
$$
$$
+N\Big(\hat b^{2}_{p_{0},\rho_{0}}+\widehat{D\sigma}_{p_{0},\rho_{0}}^{2}\Big)\rho_{0}^{-d-2}
\int_{\bR^{d+1}_{s}\times B_{1}} I_{C}u^{2n-2}|Du|^{2}\,dxdtd\eta
$$
$$
+N_{1}\Big(\hat b^{2}_{p_{0},\rho_{0}}+\widehat{D\sigma}_{p_{0},\rho_{0}}^{2}\Big)
\int_{\bR^{d+1}_{s}\times B_{1}}
\zeta^{2}|D\sigma |^{2}u^{2n}\,dxdtd\eta.
$$

For
\begin{equation}
                                 \label{12.9,6}
N_{1}\Big(\hat b^{2}_{p_{0},\rho_{0}}+\widehat{D\sigma}_{p_{0},\rho_{0}}^{2}\Big)
\leq 1/2
\end{equation}
this implies that
 $$
 \int_{\bR^{d+1}_{s}\times B_{1}}\zeta^{2}(|b|^{2}+|D\sigma|^{2})u^{2n}\,dxdtd\eta
\leq N 
 \int_{\bR^{d}\times B_{1}}\zeta^{2}(T,\cdot)f^{2n} \,dxd\eta
$$
$$
+ N \rho_{0}^{-d-4}
\int_{\bR^{d+1}_{s}\times B_{1}}   
I_{C}u^{2n} 
\,dxdtd\eta
$$
$$
+N\Big(\hat b^{2}_{p_{0},\rho_{0}}+\widehat{D\sigma}_{p_{0},\rho_{0}}^{2}\Big)
\rho_{0}^{-d-2}
\int_{\bR^{d+1}_{s}\times B_{1}} I_{C}u^{2n-2}|Du|^{2}\,dxdtd\eta.
$$
Coming back to \eqref{12.7,4} we get
$$
\int_{\bR^{d}\times B_{1}}\zeta^{2}(s,x)u^{2n}(s,x,\eta)\,dxd\eta+
\int_{[s,t_{0}]\times \bR^{d}\times B_{1}} \zeta^{2} u  ^{2n-2}
 |Du |^{2} \,dxdtd\eta
$$
$$
\leq N\int_{\bR^{d}\times B_{1}}\zeta^{2}(t_{0},\cdot)f ^{2n }   \,dxd\eta
+N \rho_{0}^{-d-4}\int_{[s,t_{0}]\times \bR^{d}\times B_{1}}I_{C}u^{2n } \,dxdtd\eta
$$
$$
+N\Big(\hat b^{2}_{p_{0},\rho_{0}}+\widehat{D\sigma}_{p_{0},\rho_{0}}^{2}\Big)\rho_{0}^{-d-2}
\int_{[s,t_{0}]\times\bR^{d} \times B_{1}} I_{C}u^{2n-2}|Du|^{2}\,dxdtd\eta.
$$

After that we repeat the same manipulations
as at the end of the proof of Lemma \ref{lemma 11.26,3}
and similarly to \eqref{12.9,4} find
$$
e^{ \lambda\rho_{0}}\int_{\bR^{d}\times B_{1}} u^{2n}(s,x,\eta)e^{-\lambda|x|}\,dxd\eta
+
\int_{[s,t_{0}]\times \bR^{d}\times B_{1}}e^{-\lambda|x|} u  ^{2n-2}
 |Du |^{2} \,dxdtd\eta
$$  
$$
\leq Ne^{ \lambda\rho_{0}}\int_{\bR^{d}\times B_{1}}e^{-\lambda|x|}f ^{2n }   \,dxd\eta
+Ne^{2\lambda\rho_{0}}\rho_{0}^{-2} \int_{[s,t_{0}]\times \bR^{d}\times B_{1}}e^{-\lambda|x|} u^{2n } \,dxdtd\eta
$$
$$
+N_{2}e^{2\lambda\rho_{0}}\Big(\hat b^{2}_{p_{0},\rho_{0}}+\widehat{D\sigma}_{p_{0},\rho_{0}}^{2}\Big)
\int_{\bR^{d+1}_{s}\times B_{1}} e^{-\lambda|x|}u^{2n-2}|Du|^{2}\,dxdtd\eta.
$$

Now along with \eqref{12.9,6} we require
$$
 N_{2}e^{2\lambda\rho_{0}}\Big(\hat b^{2}_{p_{0},\rho_{0}}+\widehat{D\sigma}_{p_{0},\rho_{0}}^{2}\Big)
\leq 1.
$$
Then
$$
\int_{\bR^{d}\times B_{1}} u^{2n}(s,x,\eta)e^{-\lambda|x|}\,dxd\eta
\leq N \int_{\bR^{d}\times B_{1}}e^{-\lambda|x|}f ^{2n }   \,dxd\eta
$$
$$
+Ne^{ \lambda\rho_{0}}\rho_{0}^{-2} \int_{[s,t_{0}]\times \bR^{d}\times B_{1}}e^{-\lambda|x|} u^{2n } \,dxdtd\eta
$$ 
and \eqref{12.20,5} follows.   
\qed

\mysection[Existence of strong solutions. $\scB=0$]
{Existence of strong solutions. Case $\scB=0$}

                      \label{section 12.18,1}

The general set up of this chapter is that
there are $q_{0},p_{0},q,p, \beta$ such that
$$
p_{0}\in(2,d+2],\quad q_{0}\in(2,p_{0}],\quad
\beta\in (1,2),\quad p= p_{0}/\beta,q=q_{0}/\beta>2,
$$
$$
\beta\leq \frac{d}{p}+\frac{2}{q}.
$$
Also  $\rho_{a},\rho_{b}\in[\rho_{0},\infty)$.

  In addition to this in the present section
we impose additional assumptions coming after
short discussion showing that we need to slightly modify the assumptions in  Section
\ref{section 7.3.1}.

The reason for the modification is that
it is more appropriate (cf. Theorem \ref{theorem 6.21.1}) to work with $\sigma$
rather than with $a$. Of course, we are going
to use mollifications of $\sigma$, which
leads to   stochastic equations
close to the original one, and
 we need to express
our previous conditions on $a$ in terms of $\sigma$ to better understand what the mollifications do for conditions on $a$.   However, the
mollifications can easily destroy the nondegeneracy. For instance, in the one-dimensional case and $\sigma(x)=\sign x$
any mollified $\sigma$ will vanish at a point.
Therefore, we need an additional assumption
preventing this from happening. We need
the mollified $\sigma$   produce $a$ probably not 
of class $\bS_{\delta}$ but lying in a wider
class $\bS_{\delta'}$ with a fixed $\delta'
\in (0,\delta]$. In this connection,
once applying the previous results
proved for $\bS_{\delta}$-valued $a$, we need
them to be true for $\bS_{\delta'}$-valued $a$
and, hence, we need to change the assumptions
in these results accordingly.  
These are the reasons for the modifications.

Namely, $a$ is still
$\bS_{\delta}$-valued and $b$ is $\bR^{d}$-valued but in other   conditions
on $a,b$ we replace $\delta$ with $\delta'$.
To be more precise, recall that $\rho_{a},\rho_{b}\in[\rho_{0},\infty)$ and note
the following simple fact. Fix a $\delta'
\in(0,\delta]$

\begin{lemma}
                        \label{lemma 12.19,3}
Let $\usigma $ be a Borel $d\times d_{1}$-valued 
  and $\ub$ be a Borel $\bR^{d}$-valued functions on $\bR^{d+1}$. Assume that $\ua:=\usigma\usigma^*$ is $\bS_{\delta'}$-valued.
Then there exist   
\begin{equation}
                               \label{12.18,10}
\widehat{D\sigma}=\widehat{D\sigma}(d, \delta',q_{0},p_{0},\beta)>0,\quad
\hat{b}=\hat{b}(d, \delta',q_{0},p_{0},\beta,
\rho_{a})>0
\end{equation}
such that if
\begin{equation}
                               \label{12.18,2}
D\usigma\in L_{1,\loc}(\bR^{d+1}),\quad
\widehat { D\usigma}_{p_{0},\rho_{a}}\leq 
 \widehat{D\sigma},
\quad\hat{\ub}_{p_{0} , \rho_{b} }
< \hat b,
\end{equation}
then

 (i) we have
 \begin{equation}
                                \label{12,19,5}
e\widehat { D\ua }_{p_{0} ,\rho_{0}}\leq 
 \widehat{Da},
\quad e\hat {\ub}_{p_{0} , \rho_{0} }
\leq  \hat b,
\end{equation}
where $(\widehat{Da},\hat b)=
(\widehat{Da},\hat b)(d,\delta',p_{0})$, are taken from Lemma 
\ref{lemma 11.26,3} when $n=1$ there;

(ii) we have
\begin{equation}
                           \label{12.19,6}
N_{1}\hat {\ub}_{p_{0},q_{0}, \rho_{b}}<\sfb_{0}(d,\delta'),
\end{equation}
where $N_{1}=N_{1}(d,\delta',q ,p,\beta ,\rho_{a} )$ is taken from \eqref{11.15,4};

(iii) we have
\begin{equation}
                                \label{12.19,7}
\ua^{\shharp}_{\rho_{a}}\leq 
\hat  a,\quad\hat {\ub}_{q_{0},p_{0}, \rho_{b}}\leq
\hat  b ,
\end{equation}
where
$$
\hat a=\hat  a(d,\delta',q,p,\beta)>0,\quad\hat b=\hat  b(d,\delta',q,p,\beta,\rho_{a} )>0,
$$ 
are taken from Theorem \ref{theorem 5.8,20};

(iv) we have
\begin{equation}
                          \label{12.19,3}
e\widehat { D\usigma }_{p_{0},\rho_{0}}\leq 
 \widehat{D\sigma},\quad e\hat {\ub}_{p_{0} , \rho_{0} }
\leq  \hat b,
\end{equation}
where $(\widehat{D\sigma},\hat b )=(\widehat{D\sigma},\hat b )(d,\delta' ,p_{0},n)$ are taken from
Theorem \ref{theorem 6.21.1} with integer $n>d/4$ there and the degree of the polynomial equal to two.

\end{lemma} 

\begin{assumption}
                     \label{assumption 12.20,1}
Condition \eqref{12.18,2} holds with 
$\widehat{D\sigma},\hat b$ from \eqref{12.18,10}
and $\usigma=\sigma,\ub=b$.
\end{assumption}

This assumption is more restrictive than
the assumptions in Section \ref{section 7.3.1}  
because, for instance, $\sfb_{0}(d,\delta)$
is  an increasing function of $\delta$.
The reader may wonder why in \eqref{12.18,2},
 we have a strict inequality.
This is done for convenience allowing us
in Section \ref{section 2.24,1} to add
to $b$ another drift with small characteristic
and be able to use the results we are going to derive under Assumption \ref{assumption 12.20,1}. Actually, this feature does not
restrict generality because the right-hand sides in such inequalities are defined  very loosely
 and could be replaced with any close
quantities. Still, to keep our arguments straight we need to state \eqref{12.18,2}
with strict inequality concerning $b$ .
 
Finally, we impose the following,
in which we fix a nonnegative $\zeta\in C^{\infty}_{0}
(B_{1})$ with unit integral and for
$\varepsilon\in(0,1]$ set 
$$
\zeta_{\varepsilon}
(x)=\varepsilon^{-d}\zeta(x/\varepsilon),\quad  \sigma^{(\varepsilon)} (t,x)
= \sigma  (t,x)*\zeta_{\varepsilon}(x),
$$ 
where the convolution is performed with respect
to $x$.  

\begin{assumption}
                   \label{assumption 12.13,1}
Either $\sigma\in B^{0,\infty}$, or
there exists an $\varepsilon_{0}\in(0,1]$
such that $a^{\varepsilon}=\sigma^{(\varepsilon)}
\sigma^{(\varepsilon)*}$ is $\bS_{\delta'}$-valued for any $\varepsilon\in(0,\varepsilon_{0}]$.
\end{assumption}

\begin{remark}
                          \label{remark 2.18,1}
Having in mind mollifying $b$ as well as $\sigma$,
but in $(t,x)$, and using well-known properties
of mollifiers, it is easy to see that under Assumptions \ref{assumption 12.20,1} and \ref{assumption 12.13,1} there exists a sequence
$\sigma(n)$, $n=1,2,...$, of $d\times d_{1}$-matrix valued functions $\sigma(n)$ and a sequence
 $b(n)$ of $\bR^{d}$-valued functions
on $\bR^{d+1}$ of class $B^{0,\infty}$ such that $a(n):=\sigma(n)\sigma^{*}(n)$ is $\bS_{\delta'}$-valued, $\sigma(n)\to\sigma$ (a.e.), $b(n)\to b$ in $L_{p_{0},q_{0}}(C)$ for any cylinder $C$,
and, for each $n$, the conditions \eqref{12.18,2} through \eqref{12.19,3} are satisfied if $\sigma(n),a(n),b(n)$
are taken in place of $\usigma,\ua,\ub$.
\end{remark}

We discuss Assumption \ref{assumption 12.13,1}
further at the end of this section only noting that it is automatically
satisfied with $\delta'=\delta$ 
if $d_{1}=d$ and $\sigma=\sqrt{a}$.  

Recall that $(p_{0},q_{0},p,q,\beta)\in \sfA_{0} $,
\index{$B$@Sets!$\sfA_{0}$}%
 where
$$
\sfA_{0}=\Big\{(q_{0},p_{0},q,p,\beta):q_{0}=\beta q\leq p_{0}=\beta p
\leq d+2,\quad p,q>2, \quad \beta<2,
$$
$$
 1<\beta\leq \frac{d}{p}+\frac{2}{q}\Big\}.
$$

\begin{theorem}
             \label{theorem 2.27.2}
Take $x_{0}\in \bR^{d}$. Then under the above assumptions equation
\eqref{3.15,1} has  an $E_{q,p,\beta}$-admissible strong solution
 and, if the above assumptions are also satisfied
with $(q'_{0},p'_{0},q',p',\beta')\in \sfA_{0}$ in place of
$(q_{0},p_{0},q,p,\beta)$, then any  $E_{q',p',\beta'}$-admissible solution  of \eqref{3.15,1}
coincides with the above strong one.

\end{theorem}

Proof. Due to Theorem \ref{theorem 6.18.3}
or Theorem \ref{theorem 2.28.1}
it suffices to prove the existence of at least
one strong $E_{q,p,\beta}$-admissible solution. Since, as it is pointed out at the beginning of Section \ref{section 7.3.1} owing to Theorem \ref{theorem 12.16,6},  
there exists an $E_{q,p,\beta}$-admissible solution, it suffices to prove that
this solution is strong. In turn
in light of Theorem \ref{theorem 6.18.2}, for that to happen, 
it suffices to show that for any $t_{0}>0$ and 
$f\in C^{\infty}_{0}$
estimate \eqref{12.10,5} holds. 

The latter estimate is proved if $\sigma,b\in B^{0,\infty}$ 
and we are going to use
\eqref{12.10,5} for our $\sigma(n),b(n)$ from
Remark \ref{remark 2.18,1}  and send  $n\to\infty$, relying on the fact that $N$
in \eqref{12.10,5} is independent of $n$.
Denote by $T_{t,s}(n),Q^{k}_{t,s}
(n)$ the operators
$T_{t,s} ,Q^{k}_{t,s}
 $  constructed on the basis
of  $\sigma(n),b(n)$.
Then by using Fatou's lemma one easily sees
that to prove \eqref{12.10,5} it suffices to show that for any $n\geq 1$, $k_{1},...,k_{n}\in\{1,...,d_{1}\}$
$$
\int_{\Gamma_{n}(t_{0})}T_{0,t_{n }}\Big[
Q^{k_{n}}_{t_{n},t_{n-1}}...Q^{k_{1}}_{t_{1},t_{0}}f\Big]^{2}(x_{0})\,d t_{n } \cdot...\cdot d t_{1}
$$
\begin{equation}
                            \label{12.13,4}
\leq\nliminf_{m\to\infty}
\int_{\Gamma_{n}(t_{0})}T_{0,t_{n }}(m)\Big[
Q^{k_{n}}_{t_{n},t_{n-1}}(m)...Q^{k_{1}}_{t_{1},t_{0}}(m)f\Big]^{2}(x_{0})\,d t_{n } \cdot...\cdot d t_{1}=:J.
\end{equation}
We will prove \eqref{12.13,4} by showing that
it is true with ($t_{n+1}=0$)
$$
\Gamma_{n,\kappa}(t_{0})=\Gamma_{n}(t_{0})\cap
\{(t_{1},...,t_{n}): t_{i+1}+\kappa\leq t_{i}
,i=0,...,n \}
$$ in place
of $\Gamma_{n}(t_{0})$ for any $\kappa>0$
and then sending $\kappa\downarrow 0$
on the left and replacing $\Gamma_{n,\kappa}(t_{0})$ back with $\Gamma_{n}(t_{0})$ on the right.

The next simplification comes from the claim
that to prove the modified \eqref{12.13,4},
it suffices to show that for any 
Borel bounded
function
$\phi_{t_{1},...,t_{n}}(x)$ vanishing for large $|x|$ we have
$$
\int_{\Gamma_{n,\kappa}(t_{0})}T_{0,t_{n }}\Big[\phi_{t_{1},...,t_{n}}
Q^{k_{n}}_{t_{n},t_{n-1}}...Q^{k_{1}}_{t_{1},t_{0}}f\Big] (x_{0})\,d t_{n } \cdot...\cdot d t_{1}
$$
\begin{equation}
                            \label{12.13,5}
=\lim_{m\to\infty}
\int_{\Gamma_{n,\kappa}(t_{0})}T_{0,t_{n }}(m)\Big[\phi_{t_{1},...,t_{n}}
Q^{k_{n}}_{t_{n},t_{n-1}}(m)...Q^{k_{1}}_{t_{1},t_{0}}(m)f\Big] (x_{0})\,d t_{n } \cdot...\cdot d t_{1}.
\end{equation}

Indeed, if \eqref{12.13,5} holds then observe
that
by H\"older's inequality the right-hand side
is dominated by
$$
J^{1/2}\Big(\lim_{m\to\infty}
\int_{\Gamma_{n,\kappa}(t_{0})}T_{0,t_{n }}(m )[\phi_{t_{1},...,t_{n}}]^{2}
(x_{0})\,d t_{n } \cdot...\cdot d t_{1}\Big)^{1/2}
$$
$$
=J^{1/2}\Big( 
\int_{\Gamma_{n,\kappa}(t_{0})}T_{0,t_{n }} [\phi_{t_{1},...,t_{n}}]^{2}
(x_{0})\,d t_{n } \cdot...\cdot d t_{1}\Big)^{1/2},
$$
where the equality follows from Theorem
\ref{theorem 12.4,2}. Hence
$$
\int_{\Gamma_{n,\kappa}(t_{0})}T_{0,t_{n }}\Big[\phi_{t_{1},...,t_{n}}
Q^{k_{n}}_{t_{n},t_{n-1}}...Q^{k_{1}}_{t_{1},t_{0}}f\Big] (x_{0})\,d t_{n } \cdot...\cdot d t_{1}
$$
\begin{equation}
                            \label{12.13,6}
\leq J^{1/2}\Big( 
\int_{\Gamma_{n,\kappa}(t_{0})}T_{0,t_{n }} [\phi_{t_{1},...,t_{n}}]^{2}
(x_{0})\,d t_{n } \cdot...\cdot d t_{1}\Big)^{1/2}
\end{equation}
Here the left-hand side is, actually,
the scalar product of 
$$
\phi_{t_{1},...,t_{n}}(x),\quad
Q^{k_{n}}_{t_{n},t_{n-1}}...Q^{k_{1}}_{t_{1},t_{0}}f(x)
$$
in an $L_{2}$-space with measure
$$
p(0,x_{0},t_{n},y)\,dyd t_{n } \cdot...\cdot d t_{1}
$$
and the factor of $J^{1/2}$ is the norm
of $\phi_{t_{1},...,t_{n}}(x)$ in this space.
Therefore, the arbitrariness of $\phi$
immediately yields \eqref{12.13,4}.

Now observe that according to Remark \ref{remark 6.27.1} the integrands in \eqref{12.13,5} are uniformly bounded on
$\Gamma_{n,\kappa}(t_{0})$. Therefore,
to prove \eqref{12.13,5} it suffices to prove
that for each $(t_{1},...,t_{n})\in
\Gamma_{n,\kappa}(t_{0})$
$$
 T_{0,t_{n }}\Big[\phi_{t_{1},...,t_{n}}
Q^{k_{n}}_{t_{n},t_{n-1}}...Q^{k_{1}}_{t_{1},t_{0}}f\Big] (x_{0})   
$$
\begin{equation}
                            \label{12.14,1}
=\lim_{m\to\infty}
 T_{0,t_{n }}(m)\Big[\phi_{t_{1},...,t_{n}}
Q^{k_{n}}_{t_{n},t_{n-1}}(m)...Q^{k_{1}}_{t_{1},t_{0}}(m)f\Big] (x_{0}) .
\end{equation}
Of course, after that we may assume that
$\phi$ is independent of $(t_{1},...,t_{n})$.
Having in mind the telescoping argument we first note that
$$
\lim_{m\to\infty}
 T_{0,t_{n }}(m)\Big[\phi 
Q^{k_{n}}_{t_{n},t_{n-1}} ...Q^{k_{1}}_{t_{1},t_{0}} f\Big] (x_{0})
= T_{0,t_{n }}\Big[\phi 
Q^{k_{n}}_{t_{n},t_{n-1}}...Q^{k_{1}}_{t_{1},t_{0}}f\Big] (x_{0}) 
$$
in light of Theorem \ref{theorem 12.4,2}.
By adding a reference to Theorem \ref{theorem 10,21,2}, we get that
$$
\lim_{m\to\infty}
 T_{0,t_{n }}(m)\Big[\phi 
Q^{k_{n}}_{t_{n},t_{n-1}}(m)Q^{k_{n-1}}_{t_{n-1},t_{n-2}} ...Q^{k_{1}}_{t_{1},t_{0}} f\Big] (x_{0})
$$
$$
=\lim_{m\to\infty}
 T_{0,t_{n }}(m)\Big[\phi 
Q^{k_{n}}_{t_{n},t_{n-1}} ...Q^{k_{1}}_{t_{1},t_{0}} f\Big] (x_{0})
= T_{0,t_{n }}\Big[\phi 
Q^{k_{n}}_{t_{n},t_{n-1}}...Q^{k_{1}}_{t_{1},t_{0}}f\Big] (x_{0}) .
$$
Keeping going in this way, we obviously come
to \eqref{12.14,1} and this brings the proof of
the theorem to an end. \qed

We have proved that \eqref{12.10,5} holds.
Along with 
Theorem  \ref{theorem 6.18.1} 
 this implies the following fact to be used
 in case $\scB\not\equiv0$.

\begin{corollary}
                \label{corollary 7.26.1}
For any $f\in C^{\infty}_{0}$, $t_{0}>0$,
and $\xi:=f(x_{t_{0}})$ we have
\begin{equation}
                            \label{12.19,5}
\sum_{m=1}^{\infty}E|\xi -\Pi^{m}_{t_{0}}\xi |^{2}
\leq N\|Df\|^{2}_{L_{4n}},
\end{equation}
 where $N$ depends only on $d,\delta,p_{0},\rho_{0},t_{0}$, 
 and an integer $n>d/4$.
\end{corollary}

\begin{remark}
                       \label{remark 12.17,3}
As is pointed out in Theorem \ref{theorem 12.12.3} under condition $
d/p+1/q\leq 1$ (say, $b$ is bounded) any solution of  
\eqref{3.15,1} is $E_{q,p,\beta}$-admissible
with the meaning of $E_{q,p,\beta}$ changing
according to $p\geq q$ or $q\geq p$.
\end{remark}

Apart from the case that $\sigma=\sqrt a$,
there are a few others when
Assumption \ref{assumption 12.13,1}
is satisfied. For instance, if the matrix $\sigma$ consists of two blocks one of which is $I$, that
is the $d\times d$-identity matrix, multiplied by
$\gamma\ne0$. In that case $\sigma^{(\varepsilon)}\sigma^{(\varepsilon)*}$ dominates $\gamma^{2}I$. In the following lemma we
single out one more case used later if $\scB\not\equiv0$.

\begin{lemma}  
                        \label{lemma 7.2.1}  
Take $\delta'\in(0,\delta)$  and suppose that
for $\rho=\rho_{a}$ (note $\bB_{r}$)
\begin{equation}
                               \label{12.18,4}
\widetilde{ D\sigma }_{p_{0},\rho}:=\sup_{t}\sup_{r\leq \rho }r\sup_{B\in \bB_{r}}\dashnorm D\sigma(t,\cdot)\|_{L_{p_{0}}(B)}\leq \widetilde{D\sigma}(d,d_{1},\delta,\delta') ,
\end{equation}
where $\widetilde{D\sigma}(d,d_{1},\delta,\delta')$
is easily found from
\index{$S$@Miscelenea!$\widetilde{D\sigma}_{p_{0},\rho}$}%
\index{$S$@Miscelenea!$\widetilde{D\sigma}$}%
 the proof that follows.
Then  
for all $\varepsilon \in(0,\rho_{a}]$   the functions  
$\sigma^{(\varepsilon)},b^{(\varepsilon)}$ satisfy 
Assumptions \ref{assumption 12.20,1} and~\ref{assumption 12.13,1}.
 
\end{lemma}

Proof. The well-known properties of convolutions
show that we only need to check that
$\sigma^{(\varepsilon)}\sigma^{(\varepsilon)*}$
is $\bS_{\delta'}$-valued.
Note that 
$$
|\sigma^{(\varepsilon)*}( t,x)\lambda|\leq \zeta_{ \varepsilon}
(  x) *|\sigma^{*}(t,x)\lambda|\leq\delta^{-1/2}|\lambda|.
$$
Therefore, we need only prove that 
\begin{equation}
                                                    \label{7.2.1}
|\sigma^{(\varepsilon)*}( t,x)\lambda|\geq |\lambda|(\delta')^{-1/2} .
\end{equation}
Without losing generality we may and will
assume that $t=0,x=0$. Then observe that
for  any $y$ we have
$$
|\sigma^{(\varepsilon)*}( 0,0)\lambda|\geq |\sigma^{(\varepsilon)*}(0,y)\lambda|-
 |(\sigma^{(\varepsilon)*}( 0,0)-\sigma^{(\varepsilon)*}(0,y))\lambda|
$$
$$
\geq |\lambda|\delta^{1/2}
-|(\sigma^{(\varepsilon)*}( 0,0)-\sigma^{(\varepsilon)*}(0,y))\lambda| 
$$
$$
\geq |\lambda|\big(\delta^{1/2}
-| \sigma^{(\varepsilon)*}(0,0)-\sigma^{(\varepsilon)*}(0,y) |\big)
$$
Furthermore, by Poincar\'e's inequality
for $\varepsilon\leq \rho_{a}$
$$
\int_{\bR^{d}}| \sigma^{ (\varepsilon)*}( 0,0)-\sigma^{ *}(0,-y)  |\zeta_{(\varepsilon)}(y)\,dy
$$
$$
\leq \int_{\bR^{d}}\int_{\bR^{d}}
| \sigma^{*}( 0,-z)-\sigma^{*}(0-y)  |\zeta_{\varepsilon}(y)\zeta_{\varepsilon}(z)\,dydz
$$
$$
\leq N(d,d_{1})\varepsilon\dashint_{B_{\varepsilon}}
|D\sigma(0,x)|\,dx\leq N(d,d_{1})
\widetilde{D\sigma}_{p_{0},\rho_{a}}.
$$
 This certainly proves the lemma. \qed

\begin{remark}  
                      \label{remark 6.2.1}
{\em In this remark the $L_{q,p}$-norm is understood as
in \eqref{3.27.3} and as in \cite{RZ_25}\/}.
It is worth comparing Theorem \ref{theorem 2.27.2} with the corresponding results about existence and uniqueness
of strong solutions for It\^o equations
belonging to
  R\"ockner and Zhao   
\cite{RZ_25}. In this paper they consider the case
that $\sigma=(\delta^{ij})$   and 
there exist $q,p$ such that
\begin{equation}
                                   \label{6.2.2}
b\in L_{q,p},\quad
q,p\in (2,\infty),\quad \frac{d}{p}+\frac{2}{q}=1
\end{equation}
or $b\in C([0,T],L_{d})$, and they prove that, for any initial
data, equation \eqref{3.15,1} has a strong solution
on $[0,T]$
possessing the property

(b) given any $p,q$ satisfying
\begin{equation}
                                \label{7.29.5}
p,q\in(1,\infty),\quad \frac{d}{p}+\frac{2}{q}<2
\end{equation}
it holds that
\begin{equation}
                                \label{7.29.4}
E  \int_{0}^{T}|f( s,x_{s})|\,ds \leq 
N \sup_{C\in \bC_{1}}\|f \|_{L_{q,p}(C)},
\end{equation}
where $N$ is independent of $f$.

They also prove
that  {\em strong\/} solutions, possessing 
  property (b),  are   unique.

For simplicity we fix $T\in(0,\infty)$
and suppose that $b(t,x) =0$ for $t\not \in [0,T]$.

{\em Case of \eqref{6.2.2} ($p>d$ and $q>2$)\/}.  
Set $p_{0}=q_{0}= q\wedge p $  and take $\beta'\in(1,2)$ 
so that $p'=q'=p_{0}/\beta'>2$. 
Then $2<p_{0}\leq 2+d$, $d/p'+2/q'=\beta'(d+2)/p_{0}\geq
\beta'$
and for any $r>0$ and $C\in\bC_{r}$
\begin{equation}
                             \label{3.16,5}
\dashnorm b\|_{L_{p_{0},p_{0}}(C)}\leq
\dashnorm b\|_{L_{q,p }(C)}=Nr^{-1}\|b\|_{L_{q,p }(C)}
\end{equation}
with $N$ independent of $r,C$. The last norm
can be made arbitrarily small on account of
taking $r$ small enough. Since $\sigma$ is the unit matrix,
 it follows that Assumption
\ref{assumption 12.20,1} is satisfied for
small $\rho_{b}$ and $\delta'=1$. Assumption
\ref{assumption 12.13,1} is trivially satisfied.
 Now by Theorem \ref{theorem 2.27.2} the 
equation 
\begin{equation}
                                 \label{3.16,4}
x_{t}=w_{t}+\int_{0}^{t}b(s,x_{s})\,ds
\end{equation}
has an $E_{q',p',\beta'}$-admissible strong
solution $x_{\cdot}$.  By H\"older's inequality
it is also an $E_{q'',p'',\beta''}$-admissible strong
solution as long as $q''\geq q',p''\geq p',\beta''\geq
\beta'$ and if $1<\beta''\leq d/p''+2/q''$ (which
is really possible after choosing $\beta'$ close
to $1$), then 
Theorem \ref{theorem 2.28.1} says that $x_{\cdot}$ is a unique
$E_{q'',p'',\beta''}$-admissible solution of \eqref{3.16,4}.

Also, for any $T\in(0,\infty)$ 
there is a constant $N$ such that for any nonnegative Borel $f$ on $\bR^{d+1}$
\begin{equation}
                            \label{3.18,2}
E\int_{0}^{T}f(s,x_{s})\,ds\leq N\|f\|_{E_{q',p',\beta'}}
\leq N\|f\|_{E_{q,p,\beta'}}
=N\sup_{C\in \bC_{1}}
\|f\|_{L_{q,p}(C)}.
\end{equation}
 
It turns out that   there is a unique
(unconditional) strong solution if $p\geq d+1$. Indeed, 
in Remark \ref{remark 10.31.1} we have seen 
that all solutions have the same finite-dimensional
distributions. In particular, they are all
$E_{q',p',\beta'}$-admissible. Since
one of them is strong, they all are strong and coincide by Theorem \ref{theorem 2.28.1}.

{\em Case $p=d$ and $b\in C([0,T],L_{d})$\/}. In that case
\begin{equation}
                                  \label{6.3.1}
\lim_{r\downarrow 0}\sup_{t\in[0,T]}
\sup_{B\in \bB_{r}}\|b(t,\cdot) \|_{L_{d}(B)}=0.
\end{equation} 

Here  one can take $p_{0}=q_{0} =d$, choose any
$  \beta\in(1,(d+2)/d),\beta<d/2 $  and set
$(q,p)=(q_{0},p_{0})/\beta$. Then, in light of 
\eqref{6.3.1} similarly to
\eqref{3.16,5} we see that
 Assumption
\ref{assumption 12.20,1} is satisfied for
small $\rho_{b}$ and $\delta'=1$. Assumption
\ref{assumption 12.13,1} is trivially satisfied. Now by Theorem \ref{theorem 2.27.2}
equation \eqref{3.16,4} has a unique $E_{q,p ,\beta }$-admissible strong
solution.

 We see that, actually, condition that
$b\in C([0,T],L_{d})$ can be replaced with $\leq
\varepsilon$ in
\eqref{6.3.1} in place of $=0$, for $\varepsilon>0$ small enough, which holds, for instance,
if the norms $\|b(t,\cdot)\|_{L_{d}}$ are uniformly
sufficiently small, that is imposed as one of alternative conditions
in \cite{RZ_25}.

As in the case of weak solutions, we see
that we have a wider class of $b$
than in \cite{RZ_25} for  which
we prove strong solvability and conditional strong
uniqueness.

Of course, one has to say that
apart from strong solvability and conditional strong
uniqueness statements \cite{RZ_25}
  contains much more highly nontrivial
information about the solutions. It is also worth 
noting that the PDE version of
assertion (b) under the Ladyzhenskaya-Prodi-Serrin
condition is derived from a more general fact
in \cite{Kr_Heat}.
\end{remark}

\begin{remark}
                         \label{remark 3.21,3}

 To show that our class of $b$ is indeed     
wider than in \cite{RZ_25} recall that,
as we have seen in Remark \ref{remark 1.28.1},  {}
the function 
$$
f(t,x)= 
\frac{1}{|x|^{\gamma}(|x|+\sqrt t)^{1-\gamma}}
I_{t>0},\quad \gamma\in\Big(\frac{d}{d+1},\frac{2d}{2d+1}\Big),
$$
does not satisfy \eqref{6.2.2}, no matter
what $p,q$ are, yet equation \eqref{3.16,4}
with $|b|=cf$ and small enough $c$ has a  (weak) solution,
and all other 
solutions have the same finite-dimensional distributions.
Actually, the
computations   in Remark \ref{remark 1.28.1}
show that if $c>0$ is small enough and $|b|=cf$,
  Theorem \ref{theorem 2.27.2} is
available   and we see that
there exists a strong solution. By Remark 
\ref{remark 1.28.1}
and Theorem \ref{theorem 2.28.1} any other
solution coincides with the strong one.

  This example shows that involving somewhat
unnatural norms as in \eqref{4.3.2} could
be quite fruitful.

We also see that the Ladyzhenskaya-Prody-Serrin condition
\eqref{6.2.2} is rather rough in what concerns the
existence of strong solutions for equations with singular drift.
 Also note that this example is not covered by
the results of \cite{KM_24} because there the condition on $b$ is imposed for each $t$
uniformly in $t$. However, our results
do not cover the results of \cite{KM_24}
either because there the condition on $b$ are stated
in terms of form-boudedness.

\end{remark}

\mysection[Dependence on the starting point]{Dependence on the starting point}  
                                                \label{section 7.3.9}

We work in the framework of Section \ref{section 12.18,1}. However, we add in Lemma \ref{lemma 12.19,3} the following statement
in which $\kappa$ is an integer $>(d+2)/2$:

{\em we have
\begin{equation}
                              \label{2.17,3}
\widehat{D\usigma}\leq e^{-1} \widehat{D\sigma}(d,\delta' ,p_{0},1),\quad\hat {\ub}\leq e^{-1}\hat b (d,\delta' ,p_{0},1),
\end{equation}
where
$ (\widehat{D\sigma},\hat b)(d,\delta' ,p_{0},1)$ are from
Theorem \ref{theorem 6.21.1} with  the degree of the polynomial equal to $2\kappa$.}

This will reduce $\widehat{D\sigma},\hat b$
in the statement of the lemma, but we still require now stronger Assumptions 
 \ref{assumption 12.20,1} and \ref{assumption 12.13,1} to hold.  
Then, of course, \eqref{2.17,3} will hold
is we replace $\usigma,\ub$ with $\sigma(n),b(n)$
from Remark \ref{remark 2.18,1}.  

The results of this section are comparable
to those in Theorem 1.1 of \cite{RZ_25} but not so elaborated
as there
(where $\sigma$ is the unit matrix).

\begin{theorem}
                \label{theorem 6.29.2}
Under the above assumption, 
there is a function $x_{s}( x)=x_{s}(\omega, x)$ which, for each   $x=x_{0}$,
is an $E_{q,p,\beta}$-admissible 
strong solution of \eqref{3.15,1} with $t=0$,   and for each $\alpha<1-(d+2)/(2\kappa)$  and $\omega$
is $\alpha$-H\"older continuous with respect to $x$ and $(\alpha/2)$-H\"older 
continuous with respect to $s$ on each set $[0,T]\times \bar B_{R}$,
$T,R\in(0,\infty)$. Furthermore, for each
$s$ with probability one $x_{s}(\cdot)\in 
W^{1}_{2\kappa,\loc}(\bR^{d})$ and
\begin{equation}
                          \label{2.17,1}
E\int_{\bR^{d}}e^{-|x|/\rho_{0}}|Dx_{s}(x)|^{2\kappa}\,dx\leq N
\end{equation}
 for any $s\leq T\in(0,\infty)$, where $N$ depends only on $T,d,\delta,p_{0},\rho_{0},\kappa$.
\end{theorem}

To prove the theorem,
we need a version of one of Kolmogorov's
results.  
Let $\bZ^{d} _{n}(2)$ be
 the subset
of $[0,1]^{d} $   of  points
$z=( z^{1}2^{-n},...,z^{d}2 ^{-n}) $, where $z^{i}=0,1,2,...,2 ^{n} $. Define
$$
\bZ^{d} _{\infty}(2)=\bigcup_{n}\bZ^{d} _{n}(2).
$$
Also let $\bZ^{1} _{n}(4)$ be the lattice
in $[0,1] $ consisting of  points
$ z 4^{-n} $, where $z =0,1,2,...,4 ^{n} $. Define
$$
\bZ^{1} _{\infty}(4)=\bigcup_{n}\bZ^{1} _{n}(4).
$$
Here is Lemma 7.3.3 of \cite{Kr_25}.
 
\begin{lemma}
                      \label{lemma 3.5.2}
Let a   random field $u(t,x)$ be defined on $\bZ^{1}_{\infty}(4)\times \bZ^{d}_{\infty}(2)$.
Assume that there exist constants $ \gamma\geq2\kappa $, $K<\infty$
such that  for $t,s\in
\bZ^{1}_{\infty}(4),x\in\bZ^{d}_{\infty}(2)$
$$
E|u(t,x)-u(s,x)|^{\gamma}\leq 
K^{\gamma} |t-s |^{\gamma/2}, 
$$
$$
E\sup_{x,y\in \bZ^{d}_{\infty}(2)}\frac{|u(t,x)-u(t,y)|^{2\kappa}}
{|x-y|^{2\kappa-d}}\leq K^{2\kappa}.
$$
 Then, for every $0<\alpha<1-(d+2)/(2\kappa)$  
with probability one there exists
a continuous extension of $u$ on
$[0,1]^{d+1}$, called again $u$, and
an integer-valued $n=n(\omega,\alpha,\gamma,\kappa,d)$
such that for any $(t,x),(s,y)\in
[0,1]^{d+1}$ satisfying $|t-s|\leq 2^{-n}$ and $ |x-y | \leq 2^{-n}$  we have
\begin{equation}
                   \label{3.3.5}
|u(t,x)-u(s,y)|\leq N(\alpha, d)K
(|t-s|^{\alpha/2}+|x-y|^{\alpha}).
\end{equation}
 
\end{lemma}
 
 {\bf Proof of Theorem \ref{theorem 6.29.2}}. First we assume that $\sigma,b\in B^{0,\infty}$. 
In that case, as  it is known since \cite{BF_61}   (1961)
(see also \cite{Ku_90} 1990),  one can define $x_{s}( x)$ in such a way  
that it becomes differentiable in $x$ for all $(\omega,s)$
and the derivative  
$ \eta^{i}D_{i}x_{s}( x)$
of $x_{s}( x)$ in the direction of $\eta$
satisfies the same equation as $\eta_{s}(0,x,\eta)$,
for which \eqref{6.20.4} holds.  
Hence, for any $x$ with probability one 
$\eta_{s}(0,x,\eta)=\eta^{i}D_{i}x_{s}( x)$
for all $s\geq0$.

Take a smooth  $f(x)$ with compact support, $t_{0}\in(0,T]$,  and set 
$$
 u( x,\eta):=
E | f_{(\eta_{t_{0} }(0,x,\eta))}(x_{t_{0} }( x)) |^{2\kappa}.
$$   
By Theorem \ref{theorem 6.21.1}, with $n=1$ there, 
 for   $\lambda=1/\rho_{0}$   
\begin{equation}
                                  \label{7.2.4}
\int_{\bR^{ d}}e^{-\lambda|x|}\sup_{|\eta|\leq 1}u^{2}( x,\eta) \,dx 
\leq N \int_{\bR^{ d}}e^{-\lambda|x|}|Df (x)|^{4\kappa }\,dx ,
\end{equation}
where (and below) the  constants $N$  
depend  only on  $d$, $\delta,\delta'$, $p_{0},\rho_{0},\kappa,T$.
Next,  
$$
 E\int_{\bR^{d}}e^{-\lambda|x|}|D\big(f(x_{t_{0}}( x))\big)|^{2\kappa }\,dx 
\leq N\int_{\bR^{d}} e^{-\lambda|x|}\sup_{|\eta|\leq 1} E\big|\big(f(x_{t_{0}}( x))_{(\eta)}\big| ^{2\kappa } \, 
dx
$$
$$   
=N \int_{\bR^{d}} e^{-\lambda|x|}\sup_{|\eta|\leq 1} u( x,\eta) \,  dx.
$$
By using \eqref{7.2.4} and H\"older's inequality we obtain that
$$
 E\int_{\bR^{d}}e^{-\lambda|x|}|D\big(f(x_{t_{0}}( x))\big)|^{2\kappa }\,dx \leq
N  \Big(\int_{\bR^{ d}}e^{-\lambda|x|}|Df (x)|^{4\kappa }\,dx\Big)^{1/2} .
$$

We obtained this estimate for smooth $f$
with compact support. By using Fatou's lemma it is extended to all smooth functions. Clearly, one can also take
$\bR^{d}$-valued $f$'s. For $f(x)\equiv x$ we get
\begin{equation}
                       \label{2.8.5}
 E\int_{\bR^{d}}e^{-\lambda|x|}|D (x_{t_{0}}( x)\big)|^{2\kappa }\,dx \leq
N  \Big(\int_{\bR^{ d}}e^{-\lambda|x|} \,dx\Big)^{1/2}=:N_{0}.
\end{equation}

By Morrey's theorem (see, for instance, Theorem 
10.2.1 of \cite{Kr_08}) this implies that ($\kappa >d/2$)
\begin{equation}
                                                        \label{3.7.20}
E \sup_{x,y\in [0,1]^{d}}\frac{| x_{s}(x)  - x_{s}(y)  |^{2\kappa }
}{|x-y|^{2\kappa -d}} \leq
N  N_{0}  .
\end{equation} 
 Furthermore, owing to Corollary \ref{corollary 3.14.6} (here we need \eqref{12.19,6}, $\beta<2$, \eqref{12.19,7})
for any $\gamma>0$, $t,s\leq 1$, 
\begin{equation}
                      \label{3.7.10}
E| x_{t}(x)  - x_{s}(x)  |^{\gamma}
\leq  N(d,\delta,\gamma)
 |t-s |^{\gamma/2}.
\end{equation}

We proved \eqref{3.7.20} and \eqref{3.7.10} assuming that
$\sigma,b$ are smooth. In the case of general $\sigma,b$, by using their
smooth approximations,
Theorem \ref{theorem 2.28.2} (need \eqref{12,19,5}), and
Fatou's lemma we conclude that these
estimates also hold if we replace
$[0,1]^{d }$ with   $\bZ^{1}_{\infty}(4)\times \bZ^{d}_{\infty}(2)$ from Lemma \ref{lemma 3.5.2}. Then, by that lemma
with probability one $ x_{t}(x) $
extends by continuity from
$\bZ^{1}_{\infty}(4)\times \bZ^{d}_{\infty}(2)$ onto $[0,1]^{d+1}$. 

Next, Theorem \ref{theorem 2.28.2}  allows us to conclude that even
if $x_{0}\in [0,1]^{d}$ and $x_{0}\not\in \bZ^{d}_{\infty}(2)$, the extension 
$ x_{t}(x_{0}) $ of
$x_{t}(x)$
 is a strong $E_{q,p,\beta}$-admissible
solution of \eqref{3.15,1}  for $s\in[0,1]$. Therefore, with probability one
we have a continuous random field
of $E_{q,p,\beta}$-admissible strong solutions of \eqref{3.15,1}  defined
on $[0,1]^{d+1}$. Furthermore, by Lemma
 \ref{lemma 3.5.2},
 $x_{s}(x)$  
is $\alpha$-H\"older continuous with respect to $x$ and $(\alpha/2)$-H\"older 
continuous with respect to $s$ on
$[0,1]^{d+1}$.

Clearly, we can extend this result from
 $[0,1]^{d+1}$ to any $[0,T]\times\bar B_{R}$, $T,R<\infty$.
This proves the first part of Theorem
\ref{theorem 6.29.2}.

 Next, take $\sigma(n),b(n)$ from Remark
\ref{remark 2.18,1} and 
let $x_{s}(n,x)$ be the functions defined by using the proved above part of Theorem 
\ref{theorem 6.29.2} applied to equation
\eqref{3.15,1}  with $\sigma(n),b(n)$
in place of $\sigma,b$. 
  Then, in light of Corollary \ref{corollary 3.14.6} and Theorem \ref{theorem 2.28.2}, for any $r\geq 1$
and $R<\infty$  
$$ 
\sup_{n}E\int_{B_{R}}|x_{s}(x,n)|^{ r}\,dx<\infty,\quad 
\lim_{n\to\infty}E\int_{B_{R}}|x_{s}(x,n)-x_{s}(x )|^{r}\,dx=0.
$$
Furthermore, we know from the proof
of Theorem \ref{theorem 6.29.2} that
\eqref{2.17,1} holds with $x_{s}(x,n)$
in place of $x_{s}(x )$ (see \eqref{2.8.5}).
It follows that there is a subsequence
of $x_{s}(x,n)$, for simplicity denoted again by $x_{s}(x,n)$, such that, for any $R$,
$x_{s}(\cdot,n)\to x_{s}(\cdot)$ in $L_{2\kappa}(\Omega\times B_{R})$ and $Dx_{s}(x,n)\to v$ weakly in
$L_{2\kappa}\big(\Omega\times \bR^{d},P(d\omega)e^{-\lambda|x|}\,dx\big)$, where $v$
is   certain function such that
\begin{equation}
                   \label{3.9.4}
E\int_{\bR^{d}}e^{-\lambda|x|}|v( x)|^{2\kappa}
\,dx\leq\nliminf_{R\to\infty}
E\int_{\bR^{d}}e^{-\lambda|x|}|Dx_{s}(x,n)|^{2\kappa}
\,dx\leq N,
\end{equation}
where $N$ in the constant in \eqref{2.17,1}.

Now it only remains to prove that (a.s.)
$v=Dx_{s}(\cdot)$, that is (a.s.)
for any $\zeta\in C^{\infty}_{0}$
we have
\begin{equation}
                   \label{3.9.30}
\int_{\bR^{d}}x_{s}(x)D\zeta(x)\,dx=
-\int_{\bR^{d}}\zeta(x)v( x) \,dx.
\end{equation}

It follows from the above convergences
that for any $A\in\cF$
$$
EI_{A}\int_{\bR^{d}}x_{s}(x)D\zeta(x)\,dx=
\lim_{n\to\infty}EI_{A}\int_{\bR^{d}}x_{s}(x,n)D\zeta(x)\,dx
$$
$$
=-\lim_{n\to\infty}EI_{A}\int_{\bR^{d}}\zeta(x)Dx_{s}(x,n)\,dx=-EI_{A}
\int_{\bR^{d}}\zeta(x)v( x) \,dx.
$$
The arbitrariness of $A$ implies that,
for any $\zeta$, \eqref{3.9.30}
holds almost surely. Then \eqref{3.9.30}
holds for any $\omega\in\Omega'$, with some $\Omega'$ of full measure, for any $\zeta$ from
a countable family which is everywhere dense in
$L_{2\kappa/(2\kappa-1)}$. We can further restrict $\Omega'$ by requiring that on it
$$
\int_{\bR^{d}}e^{-\lambda|x|}|v( x)|^{2\kappa}
\,dx<\infty
$$
 (see \eqref{3.9.4}). After that, obviously,
\eqref{3.9.30} holds on the new $\Omega'$
for any $\zeta\in C^{\infty}_{0}$.
The theorem is proved. \qed 

\begin{remark}
                    \label{remark 3.21,5}
One can find additional information on the derivatives
of solutions in  Ladyzhenskaya-Prody-Serrin
case of $b$ and unit diffusion  in \cite{RZ_25}
and the references therein.
\end{remark}

\mysection[Strong solutions, $\scB\ne0$]{Strong solutions, $\scB\ne0$}

                         \label{section 2.24,1}

Here we return to the main setting of this chapter and in contrast to the
 previous sections we do not suppose that $\scB=0$.
It may look strange that adding drift $\scB$,
which was easily absorbed by Girsanov's theorem
  in the case of weak solutions, in the case of strong solutions  forces us to basically
restart treating strong solutions  under
much heavier assumptions than when $\scB\equiv0$. A partial consolation could be
that the case $\scB\ne0$ is not covered
in \cite{RZ_25} and \cite{KM_24}
even when $\sigma=(\delta^{ij})$.

We suppose that Assumption 
\ref{assumption 12.20,1} is satisfied with
$\delta'=\delta/2$ and consider
the equation   
\begin{equation}
                    \label{12.14,2}
x _{s}=x_{0}  +\int_{0}^{s}\sigma(r,x_{r})\,dw_{r}
+\int_{0}^{s}[b (r, x_{r})+\scB(r, x_{r})] \,dr,
\end{equation}
where $\scB (t,x)$ is a Borel $\bR^{d}$-valued
function such
\index{$S$@Miscelenea!$\bar b_{R}$@$[\scB]^{2}_{\infty}$}%
\index{$S$@Miscelenea!$\bar b_{R}$@$[\scB]^{2}_{s}$}%
 that  
\begin{equation}
                    \label{12.14,3}
[\scB]^{2}_{\infty}<
\infty,\quad [\scB]^{2}_{s}:=\int_{0}^{s}\sup_{\bR^{d}}|\scB(t,x)|^{2}\,dt.
\end{equation}

To state one more assumption take $\widetilde{ D\sigma }_{p_{0},\rho }$ from \eqref{12.18,4},
 similarly 
\index{$S$@Miscelenea!$\bar b_{R}$@$\tilde b_{p_{0},\rho }$}%
define
$$
\tilde b_{p_{0},\rho }=
\sup_{t}\sup_{r\leq \rho }r\sup_{B\in \bB_{r}}\dashnorm  b(t,\cdot)\|_{L_{p_{0}}(B)}
$$
and suppose that
\begin{equation}
                             \label{12.18,1}
\widetilde{ D\sigma }_{p_{0},\rho_{a}}\leq \widehat{D\sigma}\wedge \widetilde{ D\sigma },\quad
\tilde b_{p_{0},\rho_{b}}< \hat b,
\end{equation}
where $\widehat{D\sigma}$ and $\hat b $ are taken from \eqref{12.18,10}
and $\widetilde{ D\sigma }$ is from \eqref{12.18,4} 
with $\delta'=\delta/2$ (note strict inequality).  

Since $\widehat{ D\sigma }_{p_{0},\rho_{a}}
\leq\widetilde{ D\sigma }_{p_{0},\rho_{a}}$, 
the condition on $\sigma$ is now stronger than
in Section \ref{section 12.18,1}. After Lemma 
\ref{lemma 12.19,1}
we also need to additionally assume that
\eqref{12.20,1}  (coming
after some computations with $N_{1}$
depending only on $d ,\delta,p_{0},n$,
where $n$ is a fixed integer such that $n>d/4$) holds.

Our plan of proving strong solvability of \eqref{12.14,2}
is roughly the following. First we prove that if 
$\sigma,b,\scB$
are of class $B^{0,\infty}$, then for the evolution
family $T_{t,s}$ associated with $\sigma,b,\scB$,
\eqref{12.10,5} (=\eqref{12.10,50}) holds with $N$ depending only
on  $d,\delta,p_{0}$, $\rho_{0},n,t_{0},[\scB]_{t_{0}}$.

Then note that the condition on 
$\tilde b_{p_{0},\rho_{b}}$
(and not on $\widehat {(b+\scB)}_{p_{0},\rho_{b}}$)
still allows us to use the results of   
Sections  \ref{section 4.16,1} and \ref{section 7.3.1}
in case $\scB$ is bounded.
Indeed,
observe that for any $\varepsilon\in(0,1]$ and $r\leq \rho_{0}\wedge(\varepsilon/\sup|\scB|)=\rho_{0}'  $ and $C\in \bC_{r}$
$$
\dashnorm b+\scB \|_{L_{p_{0} }(C)}
\leq r^{-1}\tilde b_{p_{0},\rho_{0}}+ \sup|B| 
\leq r^{-1}(\tilde b_{p_{0},\rho_{0}}+\varepsilon).
$$
For $\varepsilon$ small enough the latter
quantity is dominated by $r^{-1}\hat b$
in light of the strict inequality in \eqref{12.18,1}.
This shows that there is an evolution family $T_{t,s}$
corresponding to $\sigma,b,\scB$ in the general case
provided $\scB$ is bounded.

Therefore, by repeating word for word
the proof of Theorem \ref{theorem 2.27.2} we can use
approximations and get 
that \eqref{12.10,5} holds in the general case
with bounded $\scB$ and $N$ depending only
on  $d,\delta,p_{0}$, $\rho_{0},n,t_{0},[\scB]_{t_{0}}$.
This and Theorem \ref{theorem 6.18.1}
leads to the crucial conclusion that
(using the notation from Section
\ref{section 12.18,1}) for any $E_{q,p,\beta}$-admissible 
solution of \eqref{12.14,2}
for any $f\in C^{\infty}_{0}$, $t_{0}>0$,
and $\xi:=f(x_{t_{0}})$ we have
\begin{equation}
                                           \label{4.17,1}
\sum_{m=1}^{\infty}E|\xi -\Pi^{m}_{t_{0}}\xi |^{2}
\leq N\|f\|_{L_{2n}},
\end{equation}
 where $N$ depends only on
 $d,\delta,p_{0},\rho_{0},t_{0}$,  $n $,
and  $[\scB]_{\infty}$.

Of course, \eqref{4.17,1} implies that the solution
is strong and our last step consists of proving
\eqref{4.17,1} for general $\scB$ (preserving $N$) and
any $E_{q,p,\beta}$-admissible 
solution of \eqref{12.14,2}, which does exist
due to Theorem \ref{theorem 12.16,6}.

Recall that $q_{0},p_{0},\beta,q,p$ are introduced
in the introduction to the chapter.
 
\begin{lemma}   
                     \label{lemma 12.19,1}
Let $\sigma,b,\scB\in B^{0,\infty}$.
  Let $f\in C^{\infty}_{0}$. Then for any integer $n>d/4$
there exist constants $\widetilde{D\sigma'}$
and $\tilde b'$, depending only
\index{$S$@Miscelenea!$\widetilde{D\sigma'}$}%
\index{$S$@Miscelenea!$\bar b_{R}$@$\tilde b'$}%
on $d ,\delta,p_{0},n$ such that if  
\begin{equation}
                              \label{12.20,1}
\widetilde{D\sigma}_{ p_{0},\rho_{0}}\leq \widetilde{D\sigma'},\quad \tilde b_{ p_{0},\rho_{0}}
\leq\tilde b' ,
\end{equation}
then
\begin{equation}
                                \label{12.10,50}
 \sum_{n=1}^{\infty}\int_{\Gamma_{n}(t_{0})}T_{0,t_{n }}
Q_{t_{n},...,t_{0}}f(x_{0})\,d t_{n } \cdot...\cdot d t_{1}\leq
N\Big(\int_{\bR^{d}}|Df|^{4n}\,dx\Big)^{1/(2n)},
\end{equation}
where $N$ depends only on $d,\delta,p_{0},\rho_{0},n,t_{0}$
and $T_{t,s}, Q_{...}$ are taken from 
Section~\ref{section 7.3.1}.
\end{lemma}

Proof. 
First we repeat what is done in Section
\ref{section 1.14.1} (containing the proof of  Theorem 
\ref{theorem 6.21.1})
with $b+\scB$ in place of $b$  
by taking the same functions $f,u$ as in Theorem 
\ref{theorem 6.21.1},
a $C\in\bC_{\rho_{0}}$   and 
a nonnegative $\zeta
\in C^{\infty}_{0}(C)$  with the integral
of its square equal to one. After we come
to \eqref{12.7,4} we use Corollary 2.1
of \cite{Kr_26_1} which implies that for each
$t\leq t_{0},\eta$ (we drop the arguments $t,\eta$)   
$$
\int_{\bR^{d}}|b+\scB|^{2}\zeta^{2}u^{2n}(t,x,\eta) \,dx
\leq 2\int_{\bR^{d}}I_{C}|b |^{2}\zeta^{2}u^{2n} \,dx
$$
$$
+2\sup_{x}|\scB(t,\cdot)|^{2}\int_{\bR^{d}}
\zeta^{2}u^{2n} \,dx
\leq N\tilde b_{ p_{0},\rho_{0}}^{2}\int_{\bR^{d}} |D(\zeta u^{n}) |^{2}\,dx
$$
$$
+2\sup_{x}|\scB(t,\cdot)|^{2}\int_{\bR^{d}}
\zeta^{2}u^{2n} \,dx\leq N
\tilde b_{ p_{0},\rho_{0}}^{2}\int_{\bR^{d}}  \zeta^{2}u^{2n-2}|Du|^{2}  \,dx
$$
$$
+2\sup_{x}|\scB(t,\cdot)|^{2}\int_{\bR^{d}}
\zeta^{2}u^{2n} \,dx+
N
\tilde b_{ p_{0},\rho_{0}}^{2}\int_{\bR^{d}}  |D\zeta|^{2}u^{2n}  \,dx.
$$
Similarly,
$$
\int_{\bR^{d}}|D\sigma|^{2}\zeta^{2}u^{2n} \,dx
\leq N
\widetilde{ D\sigma}_{ p_{0},\rho_{0}}^{2}\int_{\bR^{d}}  \zeta^{2}u^{2n-2}|Du|^{2}  \,dx
$$
$$
+
N
\widetilde{ D\sigma}_{ p_{0},\rho_{0}}^{2}\int_{\bR^{d}}  |D\zeta|^{2}u^{2n}  \,dx.
$$

After plugging in these estimates into \eqref{12.7,4} we get (assuming $\tilde b_{ p_{0},\rho_{0}}\leq1,\widetilde{ D\sigma}_{ p_{0},\rho_{0}}\leq1$)
$$
\int_{\bR^{d}\times B_{1}}\zeta^{2}(s,x)u^{2n}(s,x,\eta)\,dxd\eta+
\int_{[s,t_{0}]\times \bR^{d}\times B_{1}} \zeta^{2} u  ^{2n-2}
 |Du |^{2} \,dxdtd\eta
$$
$$
\leq N\int_{\bR^{d}\times B_{1}}\zeta^{2}(t_{0},\cdot)f ^{2n }   \,dxd\eta
+N \int_{[s,t_{0}]\times \bR^{d}\times B_{1}}|D\zeta |^{2} u^{2n } \,dxdtd\eta
$$
$$
+N\int_{[s,t_{0}]}\sup_{x}|\scB(t,\cdot)|^{2}\int_{\bR^{d}\times B_{1	}}
\zeta^{2}u^{2n}(t,x,\eta)\,\,dxd\eta dt
$$
$$
+N \int_{[s,t_{0}]\times \bR^{d}\times B_{1}}  |D\zeta|^{2}u^{2n}\,dxdtd\eta
$$
$$
+N(\widetilde{ D\sigma}_{ p_{0},\rho_{0}}^{2}+\tilde b_{ p_{0},\rho_{0}}^{2})\int_{[s,t_{0}]\times \bR^{d}\times B_{1}}   \zeta ^{2}u^{2n-2}|Du|^{2} \,dxdtd\eta. 
$$

Now we repeat the same manipulations
as at the end of the proof of Lemma \ref{lemma 11.26,3} taking there $\lambda=0$
and similarly to \eqref{12.9,4} find
$$
\int_{\bR^{d}\times B_{1}} u^{2n}(s,x,\eta)\,dxd\eta+
\int_{[s,t_{0}]\times \bR^{d}\times B_{1}}  u  ^{2n-2}
 |Du |^{2} \,dxdtd\eta
$$
$$
\leq N\int_{\bR^{d}\times B_{1}} f ^{2n }   \,dxd\eta
+N\int_{[s,t_{0}]}(\sup_{x}|\scB(t,\cdot)|^{2}+\rho_{0}^{-2}) \int_{\bR^{d}\times B_{1	}}
 u^{2n}(t,x,\eta)\,\,dxd\eta dt
$$
$$
+N_{1}(\widetilde{ D\sigma}_{ p_{0},\rho_{0}}^{2}+\tilde b_{ p_{0},\rho_{0}}^{2})\int_{[s,t_{0}]\times \bR^{d}\times B_{1}}  u^{2n-2}|Du|^{2} \,dxdtd\eta. 
$$
The last term will be absorbed by the left-hand side if  
\begin{equation}
                             \label{12.21,1}
N_{1}(\widetilde{ D\sigma}_{ p_{0},\rho_{0}}^{2}
+\tilde b_{ p_{0},\rho_{0}}^{2})\leq 1.
\end{equation}
This and Gronwall's inequality yield 
\begin{equation}
                           \label{12.20,7}
\int_{\bR^{d}}\sup_{\eta\in B_{1}} u^{2n}(0,x,\eta)\,dx 
\leq Ne^{N  t_{0}} 
\int_{\bR^{d}\times B_{1}} 
\sup_{\eta\in B_{1}}|f| ^{2n }(x,\eta) \,dx ,                       
\end{equation}
where $N$ depends only on $d,\delta,p_{0}, n,[\scB]_{t_{0}}$,
and the power of the polynomial $f(x,\eta)$.

After that we repeat the proof of
 Theorem \ref{theorem 12.10,1} and conclude
 that to finish the proof of the lemma, it suffices
 to have the following estimate (cf. \eqref{4.18,1})
 \begin{equation}
                                    \label{4.18,3}
 T_{0,s}f(x)\leq N(1\wedge\sqrt{s})^{-d/(r\eta)}\sup_{B\in\bB_{1}}
 \|f\|_{L_{r\eta}(B)}
\end{equation}
for $f\geq0$, $x\in\bR$, $r\geq p=p_{0}/\beta$, $\eta>1$
with $N$ depending only on $d,\delta,q_{0},p_{0},\beta,\eta,
[\scB]_{t_{0}}$,
and $r$.

Let $y_{s}$ be the solution of equation 
\eqref{12.14,2} in which
$\scB$ is dropped and let $\tilde T_{t,s}$ be the evolution
family corresponding to such equation. Then by Girsanov's
theorem and H\"older's inequality
$$
T_{0,s}f(x)=Ef(x_{s})=E\phi_{s}f(y_{s})\leq
(E\phi^{\xi}_{s})^{1/\xi}\big(\tilde 
T_{0,s}(f^{\eta})\big)^{1/\eta},
$$
where $1/\xi+1/\eta=1$ and
$$
\ln \phi_{s}=\int_{0}^{s}\gamma(u,y_{u})\,dw_{u}
-(1/2)\int_{0}^{s}|\gamma(u,y_{u})|^{2}\,du,\quad
\gamma=\sigma^{*}a^{-1}\scB.
$$
As we know from \eqref{4.18,1}
$$
\tilde T_{0,s}(f^{\eta})\leq N(1\wedge\sqrt{s})^{-d/r}\sup_{B\in\bB_{1}}
 \|f^{\eta}\|_{L_{r}(B)}.
$$
Also, since $[\scB]_{t_{0}}<\infty$, for any $\xi$
we have $E\phi^{\xi}_{s}\leq N$, where $N$
depends only on $d,\delta,\xi, [\scB]_{t_{0}]}$.
This yields \eqref{4.18,3} and proves the lemma. \qed

 \begin{theorem}
             \label{theorem 12.18,1}
 
Take $x_{0}\in \bR^{d}$. Then under the   assumptions stated at the beginning of the section 
(including \eqref{12.20,1}) equation
\eqref{12.14,2} has  an $E_{q,p,\beta}$-admissible strong solution
 and any other $E_{q,p,\beta}$-admissible solution  coincides with the above
strong one.
 
\end{theorem}

We prove this theorem after a long preparations.
As in the proof of Theorem \ref{theorem 2.27.2}
it suffices to show that any solution
of \eqref{12.14,2} from Theorem \ref{theorem 12.16,6}
($E_{q,p,\beta}$-admissible solution)
is strong. Actually, $d_{1}=d$ in Theorem \ref{theorem 12.16,6} but 
as it is explained at the beginning of Section
\ref{section 7.3.1} this is irrelevant. Thus, 
let $(\Omega,\cF,P)$ be a complete probability space 
carrying a $d_{1}$-dimensional process $w_{s}$
such that equation \eqref{12.14,2} has an
$E_{q,p,\beta}$-admissible solution $x_{s}$.
 We know that
\eqref{12.21,3} holds.

For $m=1,2,...$, define  
$$
\scB_{m}= \scB I_{|\scB|\leq m},\quad
\gamma_{n}=\sigma^{*}a^{-1}(\scB-\scB_{m}),
$$
$$
 \phi_{m}=-\int_{0}^{\infty}\gamma_{m}(s,x_{s})\,dw_{s}-(1/2)\int_{0}^{\infty}|\gamma_{m}(s,x_{s})
|^{2}\,ds.   
$$
Observe that $\gamma_{m}(t,x)$ is bounded by
a function of $t$ which is square integrable
over $(0,\infty)$ (see \eqref{12.14,3}). It follows
that for any $\alpha\in\bR$ we have
$E\exp(\alpha\phi_{m})<\infty$. Another useful fact
following from \eqref{12.14,3}
is that
\begin{equation}
                        \label{7.28.3}
\int_{\bR}\sup_{\bR^{d}}|\gamma_{m}(t,x)|^{2}\,dt \to 0
\end{equation}
as $n\to\infty$. 

Next, introduce $P^{m}(d\omega)=e^{\phi_{m}}P(d\omega)$,
$$
w^{(m)}_{t}=w_{t}+\int_{0}^{t}\gamma_{m}(s,x_{s})\,ds.
$$
By Girsanov's theorem $P^{m}$ is a probability measure, $w^{(m)}_{t}$ is a Wiener process
on $(\Omega,\cF,P^{m})$, and
\begin{equation}
                                \label{7.26.2}
x_{t}=\int_{0}^{t}\sigma(s,x_{s})\,dw^{(m)}_{s}
+\int_{0}^{t}(b(s,x_{s})+\scB_{m}(s,x_{s}))\,ds.
\end{equation}

Furthermore, owing to \eqref{12.21,3}, for any integer $k\geq 1$, $T\in
(0,\infty)$, and Borel $f\geq0$
$$
E^{m}\Big(\int_{0}^{T}f(s,x_{s})\,ds\Big)^{k}
=Ee^{\phi_{m}}\Big(\int_{0}^{T}f(s,x_{s})\,ds\Big)^{k}
$$
$$
\leq \Big(Ee^{2\phi_{m}}\Big)^{1/2}
\Big(Ee^{\phi_{m}}\Big(\int_{0}^{T}f(s,x_{s})\,ds\Big)^{2k}\Big)^{1/2}\leq N\|f\|^{k}_{E_{q,p,\beta}},
$$
where $N$ is independent of $f$. In particular,
$x_{s}$ is an $E_{q,p,\beta}$-admissible solution
of \eqref{7.26.2} relative to $(\Omega,\cF,P^{m})$.  

According to what was explained before 
Lemma \ref{lemma 12.19,1} about \eqref{4.17,1}, the process
$x_{t}$   is a strong ($\{\cF^{w^{(m)}}_{s}\}$-adapted) solution of \eqref{7.26.2} and
for $f\in C^{\infty}_{0}$, $t_{0}>0$ on $\Gamma_{r}(t_{0})$, $r=1,2,...$,
for $k_{i}=1,...,d_{1}$, $i=1,...,r$, 
 there exist deterministic functions 
$f^{m,k_{1},...,k_{r}} (t_{1},...,t_{r})$ 
square integrable over $\Gamma_{r}(t_{0})$  such that
$$
\sum_{r=1}^{\infty}
E^{m}\Big|f(x_{t_{0}})-c_{m}
$$
\begin{equation}
                             \label{7.26.3}
-
\sum_{i=1}^{r}\sum_{k_{1},...,k_{i}}
\int_{\Gamma_{i}(t_{0})}
f^{m,k_{1},...,k_{i}} (t_{1},...,t_{i})
\,dw^{(m)k_{i}}_{t_{i}}\cdot...\cdot dw^{(m)k_{1}}_{t_{1}}
\Big|^{2}\leq N\|Df\|_{L_{2n}}^{2},
\end{equation}
where $c_{m}=E^{m}f(x_{t_{0}})$ and $N$ is independent of $m$.

Since
$$
c_{m}^{2}+\sum_{i=1}^{\infty}\sum_{k_{1},...,k_{i}}
\|f^{m,k_{1},...,k_{i}}\|^{2}_{L_{2}(\Gamma_{i}(t_{0}))}=E^{m}f^{2}(x_{t_{0}}),
$$
and the right-hand side is bounded by a constant
independent of $m$, there is a subsequence $m'\to\infty$ such that $f^{m',k_{i},...,k_{1}}$
converge weakly in $L_{2}(\Gamma_{i}(t_{0}))$
to certain functions $f^{ k_{i},...,k_{1}}$.
Of course, $c_{m}\to Ef(x_{t_{0}})$.

This is the first step. 

Next, for completeness, we prove the following fact
 which can be extracted from \cite{Nu_06}
or \cite{Ja_97}.

\begin{lemma}
                          \label{lemma 7.27.1}
For any $i,s=1,2,...$, $t>0$, and 
$$
f(t_{1},...,t_{i})=\{
f^{k_{1},...,k_{i}}(
t_{1},...,t_{i}),k_{j}=1,...,d_{1}\},  
$$
 given
on $\Gamma_{i}(t_{0})$ and square integrable there
we have
\begin{equation}
                         \label{7.27.1}
E\Big(\int_{\Gamma_{i}(t_{0})}f(
t_{1},...,t_{i})\,dw_{t_{i}}\cdot...\cdot dw_{t_{1}}\Big)^{2s}
\leq N \|f\|_{L_{2}(\Gamma_{i}(t_{0}))}^{2s},
\end{equation}
where $N\, (<\infty)$ depends only $i,s,d_{1}$, and by the repeated stochastic
integral above we mean
\begin{equation}
                         \label{7.27.3}
\sum_{k_{1},...,k_{i}}
\int_{\Gamma_{i}(t_{0})}f^{k_{1},...,k_{i}}(
t_{1},...,t_{i})\,dw^{k_{i}}_{t_{i}}\cdot...\cdot dw^{k_{1}}_{t_{1}}.
\end{equation}
\end{lemma}

Proof. Clearly, it suffices to prove
\eqref{7.27.1} for each particular term in
\eqref{7.27.3}. Introduce $A_{s,i}$
as the supremum of
$$
E\Big(\int_{\Gamma_{i}(t_{0})}f(
t_{1},...,t_{i})\,dB^{i}_{t_{i}}\cdot...\cdot dB^{1}_{t_{1}}\Big)^{2s}
$$
taken over all sets of
$\{B^{1}_{\cdot},...,
B^{i}_{\cdot}\}\subset \{w^{1}_{\cdot},...,w^{d_{1}}_{\cdot}\}$ and functions $f(
t_{1},...,t_{i})$ on $\Gamma_{i}(t_{0})$ having the $L_{2}$-norm equal to one. To prove the lemma, we 
only need
to show that $A_{s,i}<\infty$ for all $s,i$.

We are going to use the induction on $i$.
If $i=1$, the stochastic integral is normally distributed and \eqref{7.27.1} is obvious,
so $A_{s,1}<\infty$

Suppose that for some $i\geq1$ and any $j=1,2,...,i$
we have $A_{s,j}<\infty$.
Then take $f (t_{1},...,t_{i+1})$ such that
$\|f\|_{L_{2}(\Gamma_{i+1}(t_{0}))}=1$, and observe
that by Burknolder-Davis-Gundy inequality
$$
I:=E\Big(\int_{\Gamma_{i+1}(t_{0})}f(
t_{1},...,t_{i+1})\,dB^{i+1}_{t_{i+1}}\cdot...\cdot dB^{1}_{t_{1}}\Big)^{2 s}\leq N(s)E\Big( 
\int_{0}^{t_{0}}I^{2}( t_{1})\,dt_{1}\Big)^{s}
$$
$$
=
N(s)\int_{(0,t_{0})^{s}}EI^{2}( t^{1}_{1})\cdot...\cdot
I^{2}( t^{1}_{s})\,dt^{1}_{1}\cdot...\cdot dt^{1}_{s},
$$
where
$$
I (t_{1})= 
\int_{\Gamma_{i}(t_{1})}f 
(t_{1},t_{2},....,t_{i+1})\,dB^{i+1 }_{t_{i+1}}\cdot...\cdot 
dB^{2}_{t_{2}}.
$$
Introduce $J(t)$ by $I(t)=J(t)\|f(t,\cdot)\|_{L_{2}(\Gamma_{i}(t))}$ and observe that
$$
EJ^{2}(t^{1}_{1})\cdot...\cdot J^{2}(t^{1}_{s})
\leq \Big(\prod_{k=1}^{s} EJ^{2s}(t^{1}_{k})\Big)^{1/s}
\leq A_{s,i},
$$
where the last inequality holds by assumption. It follows that
$$
I\leq N(s)A_{s,i}
\int_{(0,t_{0})^{s}}\|f(t_{1},\cdot)\|^{2}_{L_{2}(\Gamma_{i}(t_{1}))}\cdot...\cdot
\|f(t_{s},\cdot)\|^{2}_{L_{2}(\Gamma_{i}(t_{s}))}\,dt_{1}\cdot...\cdot dt_{s}.
$$
Since the last integral, obviously, equals
$\|f\|_{L_{2}(\Gamma_{m+1}(t_{0}))}^{2n}=1$, we have
$I\leq N(s) A_{s,i}$ and the arbitrariness
of $f$ and $B^{i}_{\cdot}$ implies that $A_{s,i+1}\leq
N(s) A_{s,i}$. This proves the lemma. \qed
 
Next, we prove three more auxiliary facts.

\begin{lemma}
                        \label{lemma 7.28.1}
If a real-valued $f\in L_{2}(\Gamma_{i}(t_{0}))$, then
for any $m,s=1,2,...$ and $k_{1},...,k_{i}\in\{1,...,d_{1}\}$
$$
I:=E\Big(\int_{\Gamma_{i}(t_{0})}
f(t_{1},...,t_{i})\,dw^{(m)k_{i}}_{t_{i}}\cdot...
\cdot dw^{(m)k_{1}}_{t_{1}}\Big)^{2s}\leq 
N\|f\|_{L_{2}(\Gamma_{i}(t_{0}))}^{2s},
$$
where (note $E$ not $E^{m}$) $N$ depends only on $s,i,d,d_{1},\delta$,
and $[\scB]_{\infty}$.

\end{lemma}

The proof of the lemma is obtained by observing that
owing to Lemma \ref{lemma 7.27.1} and Girsanov's theorem
$$
I=E^{m}e^{-\phi_{m}}\Big(\int_{\Gamma_{i}(t_{0})}
f(t_{1},...,t_{i})\,dw^{ k_{i}}_{t_{i}}\cdot...
\cdot dw^{ k_{1}}_{t_{1}}\Big)^{2k}
$$
$$
\leq \Big(E e^{- \phi_{m}}\Big)^{1/2}
\Big(E^{m}\Big(\int_{\Gamma_{i}(t_{0})}
f(t_{1},...,t_{i})\,dw^{ k_{i}}_{t_{i}}\cdot...
\cdot dw^{ k_{1}}_{t_{1}}\Big)^{4k}\Big)^{1/2}.
$$ \qed

\begin{lemma}
                        \label{lemma 7.28.2}
If a real-valued $f\in L_{2}(\Gamma_{i}(t_{0}))$, then
for any $s=1,2,...$ and $k_{1},...,k_{i}\in\{1,...,d_{1}\}$
$$
E\Big(\int_{\Gamma_{i}(t_{0})}
f  (t_{1},...,t_{i})
\,dw^{(m) k_{i}}_{t_{i}}\cdot...\cdot dw^{ (m)k_{1}}_{t_{1}}-
\int_{\Gamma_{i}(t_{0})}
f  (t_{1},...,t_{i})
\,dw^{ k_{i}}_{t_{i}}\cdot...\cdot dw^{ k_{1}}_{t_{1}}\Big)^{2}
$$
\begin{equation}
                                 \label{7.28.4}
\leq \varepsilon_{m}\|f\|_{L_{2}(\Gamma_{i}(t_{0}))}^{2},
\end{equation}
where $\varepsilon_{m}$ is independent of $f$
and $\varepsilon_{m}\to0$ as $m\to\infty$.
\end{lemma}

Proof. Having in mind a usual telescoping procedure
we see that it suffices to prove that for $j=1,...,i$,
with obvious agreements in the extreme cases $j=1$ or $i$,
(keep in mind that for $j\leq i$ we set $\prod_{r=i}^{j}dw_{t_{r}}=dw_{t_{i}}\cdot...\cdot dw_{t_{j}}$)
$$
K_{j,m}:=E\Big(\int_{\Gamma_{i}(t_{0})}
f (t_{1},...,t_{i})
\prod_{r=i}^{j}  
dw^{(m)k_{ r}}_{t_{r}}dw_{t_{j-1}}^{k_{j-1}}
\cdot...\cdot dw^{k_{1}}_{t_{1}}
$$
\begin{equation}
                                         \label{7.28.04}
-\int_{\Gamma_{i}(t_{0})}
f (t_{1},...,t_{i})
\prod_{r=i}^{j+1}  
dw^{(m)k_{ r}}_{t_{ r}}dw_{t_{j }}^{k_{j }}
\cdot...\cdot dw^{k_{1}}_{t_{1}}\Big)^{2}\leq \varepsilon_{m}\|f\|_{L_{2}(\Gamma^{i}_{t_{0}})}^{2}.
\end{equation}

 Note that
with $\gamma_{m}(  t_{j} ):=\gamma_{m}( t_{j} ,x_{t_{j}})$ we have
$$
K_{j,m}
=E\Big(\int_{\Gamma_{i}(t_{0})}
f (t_{1},...,t_{i})
\prod_{r=i}^{j+1}  
dw^{(m)k_{r}}_{t_{r}}  \gamma_{m}^{k_{j}} (  t_{j} )dt_{ j }dw_{t_{ j-1}}^{k_{j-1}}
\cdot...\cdot dw^{k_{1}}_{t_{1}}\Big)^{2} 
$$
 $$
=\int_{\Gamma_{j-1}(t_{0})}
J _{m}(t_{j-1},...,t_{1})
dt_{j-1}\cdot...\cdot dt_{1},
$$
where
$$
J _{m}(t_{j-1},...,t_{1})=E\Big(\int_{\Gamma_{i-j+1}( t_{j-1})}
f (t_{1},...,t_{i})
\prod_{r=i}^{j+1}  
dw^{(m)k_{r}}_{t_{r}} \gamma_{m}^{k_{j}} (  t_{j} )dt_{ j }\Big)^{2}
$$
$$
\leq E\int_{0}^{t_{j-1}}|\gamma_{m}^{k_{j}}|^{2} (  t_{j} )dt_{j}\int_{0}^{t_{j-1}}
\Big(\int_{\Gamma_{i-j }(t_{j})}
f (t_{1},...,t_{i})
\prod_{r=i}^{j+1}  
dw^{(m)k_{r}}_{t_{r}}\Big)^{2}dt_{j}.
$$
Here the first integral under the expectation sign
tends to zero as $n\to\infty$ uniformly
with respect to $t_{j-1},\omega$ (see \eqref{7.28.3}) and 
$$
E \int_{0}^{t_{j-1}}
\Big(\int_{\Gamma_{i-j }(t_{j})}
f (t_{1},...,t_{i})
\prod_{r=i}^{j+1}  
dw^{(m)k_{r}}_{t_{r}}\Big)^{2}dt_{j}
$$
$$
=\int_{\Gamma_{i-j+1}( t_{j-1})}| 
f (t_{1},...,t_{i})|^{2}dt_{j}\cdot...\cdot
dt_{i}.
$$
This easily implies \eqref{7.28.04} and the lemma is proved. \qed

\begin{lemma}
                                \label{lemma 7.29.1}
Let $f^{m}\to f$ weakly in $L_{2}(\Gamma_{i}(t_{0}))$ 
as $m\to \infty$ and $k_{1},...,k_{i}\in\{1,...,d_{1}\}$. Then
$$
 \int_{\Gamma_{i}(t_{0})}
f^{m}  (t_{1},...,t_{i})
\,dw^{ k_{i}}_{t_{i}}\cdot...\cdot dw^{ k_{1}}_{t_{1}}
\to \int_{\Gamma_{i}(t_{0})}
f  (t_{1},...,t_{i})
\,dw^{ k_{i}}_{t_{i}}\cdot...\cdot dw^{ k_{1}}_{t_{1}}
$$
weakly in $L_{2}(\Omega)$ as $n\to \infty$.
\end{lemma}

To prove the lemma, it suffices to observe that
for any $\eta\in L_{2}(\Omega)$ the
functional
$$
E\eta\int_{\Gamma_{i}(t_{0})}
f  (t_{1},...,t_{i})
\,dw^{ k_{i}}_{t_{i}}\cdot...\cdot dw^{ k_{1}}_{t_{1}}
$$
is bounded in $ L_{2}(\Gamma_{i}(t_{0}))$, hence
continuous and weakly continuous. \qed

Now note that

$$
M^{n'}:=\int_{\Gamma_{m}(t_{0})}
f^{n',k_{1},...,k_{m}} (t_{1},...,t_{m})
\,dw^{(n')k_{m}}_{t_{m}}\cdot...\cdot dw^{(n')k_{1}}_{t_{1}} 
$$
$$
=:\int_{\Gamma_{m}(t_{0})}
f^{n',k_{1},...,k_{m}} (t_{1},...,t_{m})
\,dw^{ k_{m}}_{t_{m}}\cdot...\cdot dw^{ k_{1}}_{t_{1}}
+J^{n'}=:I^{n'}+J^{n'},
$$
where $J^{n'}\to0$ in $L_{2}(\Omega)$ as $n'\to
\infty$ by Lemma \ref{lemma 7.28.2} and by Lemma \ref{lemma 7.29.1}
$$
I^{n'}\to
\int_{\Gamma_{m}(t_{0})}
f^{ k_{1},...,k_{m}} (t_{1},...,t_{m})
\,dw^{ k_{m}}_{t_{m}}\cdot...\cdot dw^{ k_{1}}_{t_{1}}
=:M
$$
weakly in $L_{2}(\Omega)$. Since $e^{\phi_{n}/2}\to1$ strongly
in $L_{2}(\Omega)$, we also have
that $e^{\phi_{n'}/2}M^{n'}\to M$
weakly in $L_{2}(\Omega)$.
Also, obviously, $e^{\phi_{n'}/2}
  f(x_{t_{0}})\to f(x_{t_{0}})$ and
$e^{\phi_{n'}/2} c_{n'}\to Ef(x_{t_{0}})$ weakly (strongly) in  $L_{2}(\Omega)$.

By Fatou's lemma the sum of the $\nliminf$'s of the
terms on the left-hand side of \eqref{7.26.3} with $n'$ in place of $n$
is less than the $\nliminf$ of the left-hand side
of \eqref{7.26.3} and, hence, is finite. Taking into account
that ``the norm of the weak limit is less than
the $\nliminf$ of the norms'' and taking into account
the above results we conclude that
$$
\sum_{m=1}^{\infty}
E \Big|f(x_{t_{0}})-Ef(x_{t_{0}})
$$
$$
-
\sum_{i=1}^{m}\sum_{k_{1},...,k_{i}}
\int_{\Gamma_{i}(t_{0})}
f^{ k_{1},...,k_{i}} (t_{1},...,t_{i})
\,dw^{ k_{i}}_{t_{i}}\cdot...\cdot dw^{ k_{1}}_{t_{1}}
\Big|^{2}\ <\infty,
$$
which implies that $f(x_{t_{0}})$ is $\cF^{w}_{t_{0}}$-measurable and the arbitrariness of $f$ and $t_{0}$,
finally, 
bring the proof of the theorem to an end.\qed

\chapter*{Appendix: A version of Gehring's lemma}

\setcounter{equation}{0}

Here we follow \cite{Kr_23}
and give a proof of the
 the parabolic
 version of the famous Gehring's lemma
stated without proof  as Proposition 1.3  in \cite{GS_82}
with the only hint that the proof
is similar to the one given
in the elliptic case in \cite{GM_79}. 
  The author found it quite hard
to make constructions in parabolic case
``similar'' to the elliptic ones given in
\cite{GM_79} and decided to give a complete
proof having a strong probabilistic flavor.
One might think that the only difference
between elliptic and parabolic cases is
different scaling. However, in the elliptic case the doubled cubes  strictly contain 
the original ones and in the parabolic case this is not so.   
Our proof is based on the ideas from \cite{GM_79}   but the organization
of the proof is different. In particular, this allows
us to easily track down the dependence of the constants
on $A$ and shows that $q$ is a decreasing function of $A$.
If $C=C_{R}(t,x)$ and $\mu>0$ by $\mu C$
\index{$S$@Miscelenea!$\mu C$}%
 we mean
$C_{\mu R}(t,x)$.

\begin{theorem*}\em
                                           \label{theorem 3.9.1}
Let in $C_{R}$ be given a measurable $f(t,x)\geq0$
such that, for some fixed $p,A,B,\mu\in(1,\infty)$ 
satisfying $A\leq B$ and   for all 
 $C\in\bC$  such that $\mu C\subset C_{R}$ we have
$$
\Big(\dashint_{C}f^{p}\,dz\Big)^{1/p}
\leq A\dashint_{\mu C}f\,dz.
$$
Then there exists $q=q(d,p,B)>p$ such that      
$$
\Big(\dashint_{C_{R/4}}f^{q}\,dz\Big)^{1/q}
\leq N(d,p,\mu)A \dashint_{  C_{R/2}}f\,dz.
$$
\end{theorem*}

Proof. It is convenient to work with parabolic 
boxes rather than cylinders. For $n= 0,1,...$
and $ k_{0}=0,1,...,2^{2(n+1)  }-1,k_{i}=-2^{n},-2^{n}+1,
...,2^{n}-1$,
 for $i\geq1$, introduce   $D_{k_{0},...,k_{d}}(n)$
 as
$$
[k_{0}2^{-2n},(k_{0}+1)2^{-2n})
\times[k_{1}2^{-n},(k_{1}+1)2^{-n}) 
\times...\times[k_{d}2^{-n},(k_{d}+1)2^{-n}).
$$
We call $2^{-n}$ the size of $D_{k_{0},...,k_{d}}(n)$.
These are dyadic parabolic boxes,
 subsets of $D_{0} :=[0,4)\times[-1,1)^{d}$.
Set $D_{1}=[0,1)\times[-1/2,1/2)^{d}$ and for any
box $D=[S,S+T)\times Q$, where $Q$ is a cube in $\bR^{d}$,
denote   $2D=[S,S+4T)\times 2Q$, where $2Q$
is the concentric cube with twice the side length of $Q$.

Routine arguments show that to prove the theorem,
it suffices to show that there exists 
$q=q(d,p,B)>p$ such that
\begin{equation}
                                       \label{3.9.1}
\Big(\dashint_{D_{1}}f^{q}\,dz\Big)^{1/q}
\leq N(d,p )A \dashint_{ 2D_{0}}f\,dz,
\end{equation}
provided that a nonnegative $f$ is defined in $2D_{0}$ and
\begin{equation}
                                       \label{3.21.5}
\Big(\dashint_{D}f^{p}\,dz\Big)^{1/p}
\leq A\dashint_{2D}f\,dz,
\end{equation}
for any $D=D_{k_{0},...,k_{d}}(n)$ such that
$ D\subset D_{0}$.
 
To proceed in so modified setting,
  for $n\geq 0$ introduce $\Sigma_{n}$ as 
the collection
of $D_{k_{0},...,k_{d}}(n)$. To be consistent with probability
language we add to $\Sigma_{n}$ the empty set.
 Then in the terminology
from \cite{Kr_08}
the family $\{\Sigma_{n}\}$
is a filtration of partitions of $D_{0}$. Observe that
for each $n\geq 0$ and $(t,x)\in D_{0}$ there is only one
element of $\Sigma_{n}$ containing $(t,x)$. We denote
it by $\Gamma_{n}(t,x)$.
Then for each $(t,x)\in D_{0}$ define $\gamma(t,x)$
as the least $n\geq 0$ such that $ 3\Gamma_{n}(t,x)
\subset D_{0}$. Clearly, if $\gamma(t,x)=n$ and $(s,y)
\in \Gamma_{n}(t,x)$, then $\gamma(s,y)=n$. Therefore,
the set $\{(t,x): \gamma(t,x)=n\}$ is the union of
some disjoint elements of $\Sigma_{n}$.
In the terminology
from \cite{Kr_08} this means that $\gamma$ is 
a stopping time relative to the filtration  $\{\Sigma_{n}\}$.

For
each $n\geq 0$ and measurable function $g\geq0$ on $D_{0}$
 one defines the function $g_{|n}$
which on each $\Gamma\in \Sigma_{n}$ equals its average 
over $\Gamma$.

Then for a fixed $\lambda>0$ and $(t,x)\in D_{0}$ we define 
$$
\tau_{\lambda}(t,x)=\inf\{m\geq \gamma(t,x):g_{|m}(t,x)>\lambda\},
\quad(\inf\emptyset:=\infty).
$$
The set $\{\tau_{\lambda}<\infty\}$ is similar to
what one usually
gets by applying the Riesz-Calder\'on-Zygmund decomposition.
However, we are following the averages of $g$ only on dyadic boxes
where $\gamma$ is constant. Otherwise
we continue in the usual way.

Observe that $D_{0}\cap\{g>\lambda\}\subset  
D_{0}\cap \{\tau_{\lambda}<\infty\}$ (a.e.) because of the Lebesgue
differentiation theorem. 

Next, assume that, for a constant $\bar g$, we have $g_{|\gamma}
\leq \bar g$ and take $\lambda>\bar g$ so that $\tau_{\lambda}>\gamma$. Then note that
the set $D_{0}\cap \{ \tau_{\lambda}<\infty\}$
 is either empty
or is the disjoint union
of some nonempty $\Gamma_{i}\in \Sigma_{m_{i}}$, $i=1,2,...$,
on each of which $\tau_{\lambda}=m_{i}$.
Trivially,
$$
\int_{\Gamma_{i}}g\,dz=\int_{\Gamma_{i}}g_{|m_{i}}\,dz
=\int_{\Gamma_{i}}g_{|\tau_{\lambda}}\,dz,
$$
which implies that
$$
\int_{D_{0}}gI_{ \tau_{\lambda}<\infty}\,dz
= \int_{D_{0}}g_{|\tau_{\lambda}}
I_{ \tau_{\lambda}<\infty}\,dz.
$$
Furthermore, on the set $D_{0}\cap \{ \tau_{\lambda}<\infty\}$
we have $g_{|\tau_{\lambda}}>\lambda$, $g_{|\tau_{\lambda}-1}
\leq \lambda$ and, since $g_{|m}\leq 2^{d+2}g_{|m-1}$,
we have $g_{|\tau_{\lambda}}\leq\nu^{-1}\lambda$,
where $\nu=2^{-d-2}$. It follows that
$$
\nu\lambda^{-1}\int_{D_{0}}gI_{g>\lambda}\,dz\leq
\nu\lambda^{-1}\int_{D_{0}}gI_{\tau_{\lambda}<\infty}\,dz
=\nu\lambda^{-1}\int_{D_{0}}g_{|\tau_{\lambda}}
I_{\tau_{\lambda}<\infty}\,dz
$$
\begin{equation}
                                               \label{3.21.4}
\leq
|D_{0}\cap\{ \tau_{\lambda}<\infty\}|.
\end{equation}
  We apply this to $g=\phi f^{p}$,
where $\phi(t,x)=[(4-t)^{1/2}\wedge\min_{i}(1-|x^{i}|)]^{d+2}$. As is easy to see
on $D_{0}$ we have
\begin{equation}
                                            \label{3.21.10}
(\phi f^{p})_{|\gamma}\leq N(d)\int_{D_{0}}f^{p}\,dz=:\bar g.
\end{equation}

Next, define $\tilde\Gamma_{1}$ as the largest
(by size)
of the above $\Gamma_{i}$'s and by induction set
$\tilde\Gamma_{i+1}$ to be one of the largest
of $\{\Gamma_{k},k=1,2,...\}\setminus
\{\tilde \Gamma_{k},k=1,2,...,i\}$ such that
its double has no intersection with the doubles
of $\{\tilde \Gamma_{k},k=1,2,...,i\}$. There could be many $\tilde\Gamma_{i }$'s of the same size.  Let
$s_{i}$ denote the size of $\tilde \Gamma_{i}$. We claim that
\begin{equation}
                         \label{4.11.1}
|D_{0}\cap\{\tau_{\lambda}<\infty\}|\leq N(d)
\sum_{i}|\tilde \Gamma_{i}|.
\end{equation}

To prove \eqref{4.11.1} define $\hat\Gamma_{i}$
to be the union of $5\tilde \Gamma_{i}$ and its reflection in its lower base. It turns out that
\begin{equation}
                         \label{4.11.2}
D_{0}\cap\{\tau_{\lambda}<\infty\}\subset
\bigcup_{i}\hat\Gamma_{i}.
\end{equation}
Indeed, if it is not true, then there is a $\Gamma_{i}$,
which is not completely covered by the right-hand side of \eqref{4.11.2}. Let $s$ be the size of $\Gamma_{i}$.
Then there is the largest $k$ such that $s_{k}\geq s $ and $2\Gamma_{i}$ has a nonempty intersection
with at least one of $2\tilde \Gamma_{r}$, $r\geq k$
(because otherwise $\Gamma_{i}\in\{\tilde\Gamma_{r},r\leq k+1\}$). Then, since 
$s_{k}\geq s $, as is easy to see, $\Gamma_{i}
\subset \hat\Gamma_{k}$. This proves \eqref{4.11.2},
which owing to $|\hat\Gamma_{i}|\leq 2\cdot 5^{d+2}
|\tilde\Gamma_{i}|$, implies \eqref{4.11.1}.

Also note that, since $\tau>\gamma$, each of  $\tilde\Gamma^{i}$ 
is a parabolic dyadic box of size $2^{-m_{i}}$
which is the subset of a parabolic dyadic box, say
$\check\Gamma^{j}$,
of size $2^{-k}$, where  $k\leq m_{i}$ is the value
of $\gamma$ on $\check\Gamma^{j}$. It follows by construction
  that
$3\check\Gamma^{j}\subset D_{0}$. In particular, $3\tilde\Gamma^{i}
\subset D_{0}$. Also   the ratio
$\phi(z_{1})/\phi(z_{2})$ is bounded by a constant $N$
as long as $z_{1},z_{2}\in \tilde\Gamma^{i}$.
Therefore,
$$
\lambda|\tilde\Gamma^{i}|^{p}
\leq  
|\tilde\Gamma^{i}|^{p}\dashint_{\tilde \Gamma^{i }}\phi f^{p}\,dz\leq N
|\tilde \Gamma^{i}|^{p}\max_{\tilde \Gamma^{i}}\phi
\dashint_{\tilde \Gamma^{i}}  f^{p}\,dz
$$
$$
\leq
NA^{p} \min_{\tilde \Gamma^{i}}\phi
\Big(\int_{2\tilde\Gamma^{i}}f\,dz\Big)^{p}\leq
NA^{p}  
\Big(\int_{2\tilde\Gamma^{i}}\phi^{1/p}f\,dz\Big)^{p},
$$
$$
|\tilde\Gamma^{i}|\leq N_{1}\frac{A}{ \lambda^{1/p}}
\int_{2\tilde\Gamma^{i}}\phi^{1/p}f\,dz.
$$
One of inconveniences of the last estimate is that we do not have
control of $f$ on $2\tilde\Gamma^{i}$. 
In a similar situation Gehring suggested to sacrifice
some part of what is on the right to be absorbed by the left-hand side
but restrict values of $f$. So following him 
we dominate the right-hand side
by
$$
N_{1}\frac{A}{ \lambda^{1/p}}
\int_{2\tilde\Gamma^{i}}I_{\phi f^{p}>s}\phi^{1/p}f\,dz+
N_{1}\frac{As^{1/p}}{ \lambda^{1/p}}|2\tilde\Gamma^{i}|,
$$
where $s>0$ is arbitrary. For $s= N^{-p}_{2}A ^{-p}\lambda$,
where $N_{2}= N_{1}2^{d+2} $, we get
$$
|\Gamma^{i}|\leq N\frac{ A}{ \lambda^{1/p}}
\int_{2 \tilde\Gamma^{i} }I_{\phi f^{p}>s}\phi^{1/p}f\,dz
$$
and hence, coming back to \eqref{3.21.4}
(and recalling that $2 \tilde\Gamma^{i}$'s are disjoint and $3\tilde\Gamma^{i}
\subset D_{0}$),
for any $\lambda>\bar g$, we obtain
$$
\nu\lambda^{-1}\int_{D_{0}}\phi
f^{p}I_{\phi f^{p}>\lambda}\,dz
\leq  N A\lambda^{-1/p}\int_{  D_{0} }
\phi^{1/p}f I_{\phi f^{p}>N^{-p}_{2}A ^{-p}\lambda}\,dz.
$$
Multiply both sides by $\lambda^{\alpha}$, $\alpha\in(0,1]$, 
and integrate between $\bar g$
and an arbitrary finite $\Lambda>\bar g$ to get
$$
\alpha^{-1}\int_{D_{0}}\phi f^{p} ((\phi f^{p})\wedge \Lambda)^{\alpha}\,dz
-\alpha^{-1}\int_{D_{0}}\phi f^{p} ((\phi f^{p})\wedge \bar g)^{\alpha}\,dz
$$
$$
\leq N (\alpha+1-1/p)^{-1}A \int_{D_{0}}\phi^{1/p}f
  \Big((N_{2}A\phi^{1/p}f)^{p}\wedge \Lambda\Big)^{\alpha+1-1/p}\,dz .
$$

Here
$$
\int_{D_{0}}f^{p} ((\phi f^{p})\wedge \bar g)^{\alpha}\,dz
\leq \bar g^{\alpha}\int_{D_{0}}\phi f^{p}\,dz
\leq N \big(\dashint_{D_{0}} f^{p}\,dz\Big)^{1+\alpha}.
$$
Also 
$$
\phi^{1/p}f \Big((N_{2}A\phi^{1/p}f)^{p}\wedge \Lambda\Big)^{\alpha+1-1/p}
\leq (N_{2}A)^{p(\alpha+1)-1}\phi^{1/p}f((\phi f^{p})\wedge \Lambda)^{\alpha+1-1/p}
$$
$$
\leq (N_{2}A)^{p(\alpha+1)-1}\phi f^{p}((\phi f^{p})\wedge \Lambda)^{\alpha}.
$$
We   conclude that
$$
\int_{D_{0}}\phi f^{p} ((\phi f^{p})\wedge \Lambda)^{\alpha}\,dz
\leq N\Big(\int_{D_{0}}f^{p}\,dz\Big)^{1+\alpha}
$$
$$
+
N_{3}\alpha(\alpha+1-1/p)^{-1}A^{p(\alpha+1)}
\int_{D_{0}}\phi f^{p} ((\phi f^{p})\wedge \Lambda)^{\alpha}\,dz.
$$
Now choose $\alpha\leq1$ so that
$$
N_{3}\alpha(\alpha+1-1/p)^{-1}B^{2p }\leq 1/2.
$$
Then we obtain
$$
\int_{D_{0}}\phi f^{p} ((\phi f^{p})\wedge \Lambda)^{\alpha}\,dz 
\leq N\Big(\int_{D_{0}}f^{p}\,dz\Big)^{1+\alpha},
$$
which after sending $\Lambda\to\infty$ and using
\eqref{3.21.5} yields the result with $q=p(1+\alpha)$.
 The theorem is proved.\qed

 \printindex

\end{document}